\documentclass[12pt, titlepage, oneside,a4paper]{book}

\usepackage[english]{babel}
\usepackage{makeidx} 
\usepackage{hyperref}   

\usepackage{amssymb,amsmath,amsthm}
\usepackage{mathrsfs}
\usepackage{dsfont}

\makeindex       

\numberwithin{equation}{section}

\def\N{\mathbb N}
\everymath{\displaystyle}
\newenvironment{pf}{\noindent{\sc Proof}.\enspace}{\rule{2mm}{2mm}\smallskip}
\newenvironment{pfn}{\noindent{\sc Proof} \enspace}{\rule{2mm}{2mm}\smallskip}

\newtheorem{theorem}{Theorem}[section]
\newtheorem{proposition}[theorem]{Proposition}
\newtheorem{lemma}[theorem]{Lemma}
\newtheorem{corollary}[theorem]{Corollary}

 \newtheorem{remark}[theorem]{Remark}
\newtheorem{remarks}{Remark}[section]
\newtheorem{definition}[theorem]{Definition}

\font\teneufm=eufm10
\font\seveneufm=eufm7
\font\fiveeufm=eufm5
\newfam\eufmfam
\textfont\eufmfam=\teneufm
\scriptfont\eufmfam=\seveneufm
\scriptscriptfont\eufmfam=\fiveeufm

\providecommand{\norm}[1]{\lVert#1\rVert}

\newcommand{\Lip}{{\rm Lip}}
\newcommand{\lip}{{\rm lip}}
\newcommand{\Ab}{\mathscr A }
\newcommand{\Bb}{\mathscr B}
\newcommand{\mF}{\mathcal{F}}
\newcommand{\ph}{\varphi}
\newcommand{\St}{\mathbb{S}}

\newcommand{\ompaph}{\om \cdot \partial_\ph}
 
\newcommand{\dnew}{\d_1} 
\newcommand{\co}{\mathfrak c } 

\newcommand{\bt}{\underline{\tau}} 
\newcommand{\ac}{\zeta}

\newcommand{\genset}{{\cal G}}
\newcommand{\dcopen}{ $C^{\infty}$-dense open }

\newcommand{\barK}{K_{\bar V}}

\newcommand{\bh}{{\breve h}}
\newcommand{\bg}{{\breve g}}
\newcommand{\br}{{\breve r}}

\newcommand{\rf}{{\mathfrak r}}

\newcommand{\Dom}{ \om \! \cdot \! \pa_\vphi}
\newcommand{\fracchi}{{\mathfrak I} }
\newcommand{\Om}{\Omega}
\newcommand{\es}{{|{\mathbb S}|}}
\newcommand{\gap}{\mathfrak g}
\newcommand{\la}{\langle}
\newcommand{\ra}{\rangle}
\newcommand{\pa}{\partial}
\newcommand{\form}{\varkappa}
\newcommand{\dom}{\zeta}

\newcommand{\lin}{{\cal L}}
\newcommand{\matr}{{\cal M}}
\newcommand{\mapp}{{\cal H}_+}
\newcommand{\norso}[1]{| #1 |_{s_0}}
\newcommand{\nors}[1]{| #1  |_{s}}
\newcommand{\norsone}[1]{| #1  |_{s_1}}

\newcommand{\norsob}[1]{| #1  |_{s_0 + b}}

\newcommand{\norsb}[1]{| #1  |_{s+ b}}

\newcommand{\ind}{\mathfrak n}
\newcommand{\remain}{{\cal R}}
\newcommand{\cc}{\eta}

\newcommand{\be}{\begin{equation}}
\newcommand{\ee}{\end{equation}}
\newcommand{\teta}{\theta}
\newcommand{\om}{\omega}
\newcommand{\e}{\varepsilon}
\newcommand{\ep}{\epsilon}

\newcommand{\ui}{{\underline i}}
\newcommand{\uQ}{{\underline Q}}

\newcommand{\uF}{{\underline \fracchi}}
\newcommand{\uy}{{\underline y}}
\newcommand{\uth}{{\underline \theta}}
\newcommand{\uTh}{{\underline \vartheta}}
\newcommand{\uz}{{\underline z}}

\newcommand{\ov}{\overline}
\newcommand{\wtilde}{\widetilde}

\renewcommand{\a }{\alpha }
\renewcommand{\b }{\beta }
\newcommand{\s }{\sigma }
\newcommand{\ppavphi}{\bar \om_\e  \cdot  \pa_{\vphi}  } 
\newcommand{\ppamuvphi}{\bar \mu  \cdot \pa_{\vphi}  } 
\newcommand{\ii }{{\rm i} }
\renewcommand{\d }{\delta }
\newcommand{\loss}{\varsigma}
\newcommand{\D }{\Delta}
\newcommand{\g }{\gamma}

\renewcommand{\l }{\lambda }

\newcommand{\m }{\mu }

\newcommand{\vphi}{\varphi }

\renewcommand{\t }{\tau }
\renewcommand{\o }{\omega }
\renewcommand{\O }{\Omega }
\newcommand{\C}{\mathbb{C}}
\newcommand{\Z}{\mathbb{Z}}
\newcommand{\T}{\mathbb{T}}
\newcommand{\R}{\mathbb{R}}

\newenvironment{Remark}{\begin{remark} \rm}{\rule{2mm}{2mm}\end{remark}}

\newcommand{\dps}{\displaystyle}

\newcommand{\norma}{|}

\newcommand{\FA}{\mathfrak A}
\newcommand{\FR}{\mathfrak R}
\newcommand{\FD}{\mathfrak D}
\newcommand{\FW}{\mathfrak W}
\newcommand{\Fm}{\mathfrak m}

\newcommand{\ddI}{{\mathfrak I}_{\delta}}

\title{\bf Quasi-periodic solutions of \\
nonlinear wave equations \\
on the $ d $-dimensional torus 
}
\author{Massimiliano Berti\\SISSA, Via Bonomea 265\\ 34136, Trieste, Italy \\
{\tt berti@sissa.it}.
\\\&\\
Philippe  Bolle, \\
Avignon Universit\'e, LMA EA2151\\
F-84000 Avignon, France, \\
{\tt philippe.bolle@univ-avignon.fr}.
}
\date{} 

\begin{document}

\frontmatter
\maketitle

\chapter*{Preface}

\thispagestyle{empty}

Many Partial Differential Equations (PDEs) arising in
physics can be seen as infinite dimensional Hamiltonian systems
\be\label{HS-PDE}
\partial_t z = J (\nabla_z H)(z) \, , \quad z \in E \, , 
\ee
where the Hamiltonian function
$ H : E   \to \R  $
is defined on an \emph{infinite} dimensional Hilbert 
space $ E $ of functions $ z  := z(x) $, and
$ J $ is a non-degenerate antisymmetric operator.

Main examples 
are the {\it nonlinear wave} 
equation (NLW)
\be \label{NLWeq1}
u_{tt} - \Delta u  + V(x) u + g(x,u) = 0   \, , 
\ee
the nonlinear Schr\"odinger equation (NLS), 
the  beam equation and the higher dimensional membrane equation, 
the water waves equations, i.e. the Euler equations of Hydrodynamics 
describing the evolution of 
an incompressible  
irrotational fluid under the action of gravity and surface tension,
as well as its approximate models 
like the Korteweg de Vries (KdV) equation, 
the Boussinesq, Benjamin-Ono,  
Kadomtsev-Petviashvili  (KP) 
equations, \ldots, among many others.    
We refer to  \cite{Ku1} for a general introduction to  Hamiltonian PDEs. 

\smallskip

In this monograph we shall 
adopt a ``Dynamical Systems" 
point of view,   regarding the nonlinear wave
equation \eqref{NLWeq1}, equipped with periodic boundary conditions 
$ x \in \T^d := (\R/ 2 \pi \Z)^d $, as an infinite dimensional Hamiltonian system, and
we shall prove the existence of Cantor families of finite dimensional invariant tori, filled by 
quasi-periodic solutions of \eqref{NLWeq1}.  The first results in this direction are due to Bourgain \cite{B5}. 
The search of invariant sets for the flow
is an essential change of paradigm in the study of hyperbolic equations, 
with respect to the more traditional pursuit of  the initial value problem. This perspective 
has allowed to find many new results, inspired by finite dimensional Hamiltonian systems,
for Hamiltonian PDEs. 

\smallskip

When the space variable $ x $
belongs to a {\it compact} manifold, 
say $ x \in [0, \pi] $ with Dirichlet boundary conditions 
or $ x \in \T^d $ (periodic boundary conditions), 
the dynamics of a Hamiltonian PDE \eqref{HS-PDE}, like \eqref{NLWeq1},  is expected to  have a  ``recurrent" behaviour in time,   
with many   {\it periodic} and
 {\it quasi-periodic} solutions, i.e.
solutions
(defined for {\it all} times)  of the form 
\be \label{QP-sol}
 u (t) = U (\om t ) \in E  
  \qquad {\rm where} \qquad   \T^\nu \ni 
  \vphi \mapsto U( \vphi ) \in  E 
 \ee
is $ 2 \pi $-periodic in the angular variables 
$ \vphi := (\vphi_1, \ldots, \vphi_\nu)  $ and  the frequency  vector $ \om \in \R^\nu $   is nonresonant, namely 
$ \om \cdot \ell  \neq 0 $, $ \forall \ell  \in \Z^\nu \setminus \{0\} $.
When $ \nu = 1 $ the solution $ u(t)  $ is periodic in time, with period $ 2 \pi \slash \om $.
If $ U(\om t ) $ is a quasi-periodic solution then, since  the orbit
$ \{ \om t \}_{t \in \R } $ is {\it dense}  on $ \T^\nu $, the  
torus-manifold  $  U ( \T^\nu )  \subset E$
is invariant under the flow of \eqref{HS-PDE}.  

Note that 
the linear wave equation \eqref{NLWeq1} with $ g = 0 $,  
\be\label{eq:LIN}
u_{tt} - \Delta u + V(x) u = 0  \, , \quad x \in \T^d  \, ,
\ee
possesses many quasi-periodic solutions.
Indeed the self-adjoint operator
$ - \Delta + V(x) $ has a complete $ L^2 $-orthonormal  basis  of eigenfunctions $ \Psi_j(x) $, $ j \in \N $,  with eigenvalues 
$ \lambda_j \to + \infty $,  
\be\label{def:eigen}
( - \Delta + V(x)) \, \Psi_j(x) = \lambda_j \, \Psi_j (x) \, , \quad j \in \N \, .
\ee
Supposing for simplicity  that $ - \Delta + V(x) >  0 $, the eigenvalues 
$\lambda_j=\mu_j^2$, $\mu_j>0$, are positive, and 
{\it all} the 
solutions of \eqref{eq:LIN} are 
\be\label{linear:sol}
 \sum_{j \in \N}   \a_j  \cos ( \mu_j t + \theta_j )\Psi_j(x) \, , \quad  \a_j , \ \theta_j \in \R   \, ,
\ee
which, depending on the resonance properties of the linear frequencies $ \mu_j = \mu_j (V) $, 
are periodic,  quasi-periodic, or almost-periodic in time
(i.e. quasi-periodic with infinitely many frequencies).

What happens  to these solutions under the effect of the nonlinearity $g(x, u)$ ?

\smallskip

There exist special nonlinear  equations 
for which all the solutions are still periodic, quasi-periodic  or almost-periodic in time, for example 
the  
KdV and  Benjamin-Ono equations, $ 1d$ defocusing cubic-NLS,  \ldots . 
These are  completely integrable PDEs.
However, for generic nonlinearities, 
one expects, in analogy with the celebrated 
Poincar\'e non-existence theorem of prime integrals for nearly integrable Hamiltonian systems, 
that this is not the case. 

\smallskip

On the other hand,  
for sufficiently small Hamiltonian perturbations of a non degenerate integrable system
in $ \T^n \times \R^n $,
the classical KAM --Kolmogorov-Arnold-Moser--
theorem proves 
  the persistence of quasi-periodic solutions 
 with  {\it Diophantine}
frequency  vectors $ \omega \in \R^n $, i.e. vectors satisfying for some  $ \g >  0 $ and $ \t \geq  n - 1 $, 
the non-resonance condition   
\be \label{dioph:preface}
|\om \cdot \ell | \geq \frac{\gamma}{|\ell|^\tau} \, , \quad \forall \ell \in \Z^{n} \setminus \{0\} \, . 
\ee
Such frequencies form a Cantor set of $ \R^n $ of positive measure if $ \t >  n - 1 $. 
These quasi-periodic solutions (which densely fill invariant Lagrangian tori)  were constructed
by Kolmogorov \cite{K} and Arnold \cite{Ar} for analytic systems  
using an iterative Newton scheme, and
by Moser \cite{Mo62}-\cite{M67} for 
differentiable perturbations by introducing smoothing operators.
This scheme then gave rise to abstract Nash-Moser implicit function theorems like 
the ones due to Zehnder in \cite{Z1,Z2}, see also
\cite{N} and Section \ref{sec:NNM}. 

What happens for infinite dimensional systems like PDEs? 

\begin{itemize}
\item The central question of KAM theory for PDEs is: 
{\it do ``most" of the periodic, quasi-periodic, almost-periodic solutions of an integrable PDE (linear or nonlinear)
persist, just slightly deformed, under the effect of a nonlinear perturbation?}
\end{itemize}

KAM theory for Partial Differential Equations 
started a bit more than thirty years ago with the 
pioneering works of Kuksin \cite{Ku} and 
Wayne \cite{W1},  
about the existence of quasi-periodic solutions 
for semilinear perturbations of
1-dimensional linear wave and Schr\"odinger equations
in the interval $[0, \pi] $.  
These results are based 
on an extension of the KAM perturbative
approach
developed for the search of lower dimensional tori in 
finite dimensional systems, see \cite{M67}, \cite{El88}, \cite{P89}, 
and relies on the  verification of the so called second order Melnikov non-resonance conditions.

Nowadays KAM theory for 1-$d$ partial differential equations 
has reached a satisfactory 
level of comprehension concerning  quasi-periodic solutions, 
while questions concerning  almost periodic solutions remain quite open. The known results
include bifurcation of 
small amplitude solutions  \cite{KP2}, \cite{Po3}, \cite{BBi11}, 
perturbations of large finite gap solutions  \cite{K2-KdV}, \cite{Ku1}, \cite{BoK}, \cite{KaP}, \cite{BKM},  extension to periodic boundary conditions 
\cite{CW}, \cite{Bo1}, \cite{CY}, \cite{GY0},
use of weak non-degeneracy conditions \cite{BaBM}, 
nonlinearities with derivatives \cite{Liu-Yuan},  \cite{ZGY}, 
\cite{Berti-Biasco-Procesi-rev-DNLW} up to 
quasi-linear ones 
\cite{BBM-Airy}-\cite{BBM-auto}, \cite{FP}, including water-waves equations 
\cite{BM16}, \cite{BBHM}, applications to quantum harmonic oscillator \cite{GT}, \cite{Bam1}-\cite{Bam2},
and a  few examples of almost periodic solutions 
\cite{Po4}, \cite{B7}. 
We describe these developments more in detail in Chapter \ref{Ch1.5}. 

\smallskip

Also KAM  theory for multidimensional 
PDEs still contains few results
and a satisfactory  picture is  under construction. 
If  the space dimension $ d $ is two or more, major difficulties are the following ones: 
\begin{enumerate}
\item  \label{pre1} the eigenvalues $ \mu_j^2 $ of the 
Schr\"odinger operator 
$  - \Delta + V(x) $ 
in \eqref{def:eigen}
appear in huge clusters of increasing  size. 
For example, if $ V(x) = 0 $, and $x \in \T^d $, 
they are
$$
|k|^2 = k_1^2 + \ldots + k_d^2 \, , \quad k = (k_1, \ldots, k_d) \in \Z^d  .
$$
\item \label{pre2}
The eigenfunctions $ \Psi_j (x) $ may be 
 ``not localized" with respect to the exponentials. 
 Roughly speaking, this means that there is 
 no one-to-one correspondence  $h: \N \to \Z^d$ such that 
the entries  $ (\Psi_j,  e^{\ii k \cdot x })_{L^2}$ of the change-of-basis matrix
(between $ ( \Psi_j )_{j\in \N} $  and $ (e^{\ii k \cdot x})_{k\in \Z^d} $)
decay  rapidly   to zero as  $ |k-h(j)| \to  \infty $. 
\end{enumerate}

The first existence result of time periodic solutions for the 
nonlinear 
wave equation 
$$
u_{tt} - \Delta u + m u = u^3 + {\rm h.o.t.}  \, , \quad x \in \T^d \, , \ d \geq 2 \, , 
$$
has been proved by Bourgain in \cite{B-Gafa}, 
extending the Craig-Wayne  approach \cite{CW}, originally developed 
if  $ x \in \T $.  
Further existence results  of periodic solutions 
have been proved in Berti-Bolle \cite{BBARMA} for merely differentiable nonlinearities, 
Berti-Bolle-Procesi \cite{BBP10} for  Zoll manifolds, 
Gentile-Procesi \cite{GP} using  Lindstedt series techniques,  and  
 Delort \cite{D} for NLS using paradifferential calculus. 

\smallskip

The first breakthrough result about existence of quasi-periodic solutions 
for space multidimensional PDEs 
was due to Bourgain \cite{B3} for  analytic Hamiltonian NLS  equations of the form
\be \label{NLS:model}
\ii u_t =  \Delta u + M_\sigma u + \e \partial_{\bar u} H(u, \bar u)
\ee
with $ x \in \T^2 $, where $ M_\sigma = {\rm Op}( \sigma_k )$ is a Fourier multiplier 
supported on finitely many sites $ {\mathbb E} \subset \Z^2 $, i.e. 
$  \sigma_{k} = 0  $, $  \forall  k \in  \Z^2 \setminus {\mathbb E} $.
The $ \sigma_k   $, $ k \in {\mathbb E} $,  
play the role of external parameters
used to verify suitable non-resonance conditions.
Note that  the eigenfunctions of $ \Delta  + M_\sigma$
are the exponentials 
$ e^{\ii k \cdot x } $ and so the above mentioned problem \ref{pre2} is not present.

Later on,
using tools  of sub-harmonic analysis 
 previously developed  for  quasi-periodic Anderson localization theory,  
 in
  Bourgain-Goldstein-Schlag  \cite{BGS}, \cite{B6}, 
Bourgain \cite{B5}  was able to extend this result 
in any space dimension $ d $, and also for  nonlinear wave equations of the form 
\be\label{model-B}
u_{tt} - \Delta u + M_\sigma u + \e F'(u) = 0 \, , \quad x \in \T^d \, , 
\ee
where  $ F (u) $ is a polynomial in $ u $.
Here $ F' (u) $ denotes the derivative of $ F $. 
We also mention the existence results of quasi-periodic solutions 
of Bourgain-Wang \cite{BW1}-\cite{BW2} for NLS and NLW equations under a 
random perturbation. The stochastic case  is a priori easier than the deterministic one because 
it is simpler to  verify the non-resonance conditions with  a random variable.  

\smallskip

Quasi-periodic solutions $ u (t,x) = U(\om t, x ) $ of 
\eqref{model-B}
with a frequency vector $ \om \in \R^\nu $, 
namely solutions $ U(\vphi, x) $, $ \vphi \in \T^\nu $, of   
\be\label{model-B-te}
(\om \cdot \pa_{\vphi})^2 U - \Delta U + M_\sigma U + \e F'(U) = 0 \, , 
\ee
are constructed by  a Newton scheme.
The main analysis
 concerns 
 finite dimensional restrictions 
of the quasi-periodic operators  obtained by linearizing \eqref{model-B-te} at each step of the Newton iteration, 
\be\label{QP-lin:OP}
{\it \Pi}_N \big( (\om \cdot \pa_{\vphi})^2  - \Delta + M_\sigma  + 
\e b(\vphi, x)  \big)_{| {\cal H}_N}  \, ,  
\ee
 where $ b(\vphi, x) = F'' (U(\vphi,x)) $ and  
 $ {\it \Pi}_N $ denotes the projection on the finite dimensional subspace 
 $$ 
 {\cal H}_N := \Big\{ h = \sum_{|(\ell, k)| \leq N} h_{\ell, k} e^{\ii (\ell \cdot \vphi + k \cdot x)} \, , \, \ell \in \Z^\nu \, , k \in \Z^d \Big\} \, . 
 $$ 
The matrix which represents \eqref{QP-lin:OP} in the exponential basis is a 
perturbation  of the diagonal matrix
$ {\rm Diag} ( - (\om \cdot 	\ell)^2 + |k|^2 + \sigma_k ) $ 
with  off-diagonal entries 
$ \e (\widehat b_ {\ell - \ell', k - k'}) $
which decay exponentially to zero as $ |(\ell-\ell', k - k')| \to + \infty $, 
assuming that $b$ is analytic (or sub-exponentially, if $ b $ is Gevrey). 
The goal is to prove that such matrix 
is invertible, for most values of the external parameters, and that its inverse 
has an exponential (or Gevrey)  off-diagonal decay.
It is not difficult to impose lower bounds for the eigenvalues of the self-adjoint 
operator \eqref{QP-lin:OP} for most values of the parameters. 
These ``first order Melnikov" non-resonance conditions are essentially 
the minimal assumptions for proving the persistence of quasi-periodic solutions
of \eqref{model-B}, and 
provide
 estimates of the inverse of the operator  \eqref{QP-lin:OP} in $ L^2 $ norm.
In order to prove fast off-diagonal decay estimates for the inverse matrix,
Bourgain's technique  is a ``multiscale" inductive analysis based on the repeated use of the 
``resolvent identity".
An essential ingredient is that the  ``singular"  sites  
\be\label{sing:cone}
(\ell, k) \in \Z^\nu \times \Z^d  \quad {\rm such \ that } \quad 
| - ( \om \cdot \ell )^2 + |k|^2 + \sigma_k | \leq 1  
\ee
are separated into clusters which are sufficiently distant from one another.
However,  the information \eqref{sing:cone} about just the linear frequencies of
\eqref{model-B}  is not sufficient (unlike for time-periodic solutions \cite{B-Gafa}) 
in order to prove that the  inverse matrix has an exponential (or Gevrey) off-diagonal decay.
Also finer non-resonance conditions at each scale 
along the induction need to be  verified. 
We describe the multiscale approach in Section \ref{sec:MULTI} and
we prove novel multiscale  results in Chapter \ref{sec:multiscale}.

These techniques have been extended 
in   the recent work of Wang \cite{Wang1} 
for the nonlinear  Klein-Gordon equation
$$
u_{tt}  - \Delta u + u + u^{p+1} = 0 \, , \quad p \in \N \, , \ x \in \T^d \, , 
$$
that, 
unlike \eqref{model-B}, is  
parameter independent. A key step is 
to verify that suitable non-resonance conditions are fulfilled 
for most ``initial data". We refer to \cite{Wang} for a corresponding result for NLS.

\smallskip

Another stream of important results for multidimensional PDEs have been 
inaugurated in the breakthrough paper \cite{EK}
of Eliasson-Kuksin  for the NLS equation \eqref{NLS:model}. 
In this paper the authors are able to 
block diagonalize, and reduce to constant coefficients, the quasi-periodic Hamiltonian operator
obtained at each step of the iteration. 
This KAM reducibility approach extends the
perturbative 
theory developed  for $ 1d$-PDEs, 
by  verifying the so called second order Melnikov non-resonance conditions.
It allows to prove directly also the linear stability of the 
quasi-periodic solutions. 
Other results in this direction
have been proved for the $ 2d$-cubic NLS by Geng-Xu-You \cite{GXY},  
by Procesi-Procesi \cite{PP}, \cite{PP3} in any space dimension  and arbitrary polynomial nonlinearities, 
by Geng-You \cite{GY1.5} and Eliasson-Gr\'ebert-Kuksin \cite{EGK} for beam equations.   
Unfortunately,
the second order Melnikov conditions are 
violated  
for nonlinear wave equations for which 
an analogous reducibility result does not hold. 
We describe the KAM reducibility approach with PDEs applications in  Sections \ref{sec:red} and \ref{sec:RR}.

\medskip

We now present more in detail the goal of this research monograph. 
The main result 
is  the existence of  small amplitude
time quasi-periodic solutions  for autonomous nonlinear wave equations 
of the form 
\be\label{lanostra}
u_{tt} - \Delta u + V(x) u + g(x, u) = 0 \, , \quad x \in \T^d \, ,  
\quad g (x,u) = a(x) u^3 + O(u^4 ) \, , 
\ee
in {\it any} space dimension $ d \geq 1 $,
where $ V(x)  $ is a smooth {\it multiplicative} potential such that 
$ - \Delta + V(x) > 0 $, 
 and the nonlinearity is $ C^\infty $. 
Given a finite set $ {\mathbb S} \subset \N $ (tangential sites), we construct
quasi-periodic solutions $u ( \omega t ,x )  $ with frequency vector
$ (\om_j)_{j \in {\mathbb S}}$,  of the form
\be \label{sol:NLW} 
u ( \omega t ,x ) = 
\sum_{j \in {\mathbb S}}  
\a_j \cos ( \omega_j t ) \Psi_j (x) + r (\omega t , x)  \, , \quad 
\omega_j = \mu_j + O(|\a|)  \, ,
\ee
where 
the remainder $ r (\vphi, x )  $ is $ o(|\a| ) $-small   in some Sobolev space;
here $ \a := ( \a_j)_{j \in {\mathbb S} } $.
The solutions \eqref{sol:NLW} 
are thus a small deformation of linear solutions \eqref{linear:sol},  
supported on 
the ``tangential" space spanned by  
the eigenfunctions $ (\Psi_j (x))_{j \in {\mathbb S} }$,  
with a much smaller component in the normal subspace.
These quasi-periodic solutions 
 of \eqref{lanostra} exist for generic potentials $ V(x) $, coefficients 
$ a(x)$ and ``most" small values of the amplitudes 
$   (\a_j)_{j \in {\mathbb S} }  $. 
The precise statement is given in Theorems \ref{thm:main} and \ref{Prop:genericity}.

\smallskip

The proof of this result requires  various mathematical methods which this book aims to present in a systematic and self-contained way. A complete outline of the 
steps of proof is  presented in Section \ref{sec:ideas}. 
Here we just mention that we shall use a Nash-Moser iterative scheme in scales of Sobolev spaces for the search of an invariant  torus embedding supporting quasi-periodic solutions, with a 
frequency vector 
$ \omega $ to be determined.
One key step is to establish the existence of an 
approximate inverse for the  operators obtained by
linearizing the nonlinear wave equation at any approximate 
quasi-periodic solution $ u(\omega t,x) $,
and to prove that such an approximate inverse
satisfies tame estimates  in Sobolev spaces, with loss of derivatives due to the small 
divisors. These  linearized operators have the form
$$
h \mapsto (\om \cdot \pa_\vphi)^2 h - \Delta h+ V(x) h + (\pa_u g)(x, u(\omega t,x)) h 
$$
with coefficients depending on $ x \in \T^d $ and $ \vphi \in \T^\es $. 
The construction of an approximate inverse requires several steps. 
After writing the wave equation as a Hamiltonian system in infinite dimension, 
the first step is to use 
a symplectic change of variables to approximately decouple the tangential
and normal 
components of the linearized operator. This is a rather general procedure for autonomous 
Hamiltonian PDEs, which reduces the problem to the search of an approximate inverse for a  quasi-periodic 
Hamiltonian  linear operator acting in the subspace normal to the torus, see 
Chapter \ref{sezione almost approximate inverse} and Appendix \ref{sec:2}. 

\smallskip

In order to avoid the difficulty posed by the violation of 
the second order Melnikov non-resonance conditions required by 
a KAM reducibility scheme, we develop a 
multiscale inductive approach {\`a la} Bourgain,  which is particularly delicate since the eigenfunctions $ \Psi_j (x) $ of $ - \Delta + V(x) $ defined in \eqref{def:eigen} are not localized near the exponentials.
In particular the matrix elements $ (\Psi_j, a(x) \Psi_{j'})_{L^2} $ representing 
the multiplication operator with respect to  the basis of the 
eigenfunctions $ \Psi_j (x) $ 
do not decay, in general,   as $ j- j' \to  \infty $.    
In 
Chapter \ref{sec:multiscale} we provide the complete 
 proof of the multiscale proposition  (which is fully self-contained  together 
 with the Appendix \ref{App:mult})
which we shall use in Chapters \ref{sec:splitting}-\ref{sec:proof.Almost-inv}.   
These results extend the multiscale analysis developed for
forced NLW and NLS equations in  \cite{BB12}-\cite{BBo10}.  

The presence of a multiplicative potential $ V(x) $ in \eqref{lanostra} makes also difficult to control the variations of the 
tangential and normal frequencies due to the effect of the nonlinearity 
$ a(x) u^3 + O(u^4) $ with respect to parameters. 
In this monograph, after a careful bifurcation analysis of the quasi-periodic solutions, 
we are able to use 
just the  length $ |\omega| $ of the frequency vector
as an internal parameter to verify all the 
non-resonance conditions along the iteration. 
The frequency vector is constrained to a fixed direction, see \eqref{frequencycolinear}-\eqref{def omep}.
The measure estimates  rely on positivity arguments for the variation
of parameter dependent 
families of self-adjoint  matrices, see Section \ref{sec:mult5}. These properties,
see \eqref{Hyp1},  are verified
for the  linearized operators obtained along the iteration. 

The genericity of the 
non-resonance and non-degeneracy conditions that we require on the potential $ V(x) $ and the coefficient
$ a(x) $ in the nonlinearity $ a(x) u^3 + O(u^4) $, are finally verified in Chapter \ref{section:gener1}. 

\smallskip

The techniques developed above for the NLW equation \eqref{lanostra}
would certainly apply to prove a corresponding result for 
nonlinear Schr\"odinger equations. However we have decided to focus on the 
NLW equation because,
as explained above,  there are less results available for the latter one. 
This context seems to better illustrate
the advantages of the present approach in comparison to  
that of reducibility. 

\smallskip

A  feature of the monograph is to present  the proofs, techniques and ideas developed  
in a self-contained and expanded manner, 
with the  hope to enhance further developments. 
We also aim to describe the connections of this 
result  with previous works 
in the literature.  The techniques developed in this monograph  
have deep connections with those used in Anderson localization theory
and we hope that the detailed  presentation in this manuscript 
of all technical aspects of the proofs, will allow a deeper interchange between the 
scientific communities of Anderson-localization and of ``KAM for PDEs".

\bigskip
\bigskip
\bigskip

\noindent
Massimiliano Berti, Philippe Bolle,


\thispagestyle{empty}

\noindent
{\it Keywords:} Nonlinear wave equation, KAM for PDEs,  
quasi-periodic solutions, small divisors,  invariant tori, multiplicative potential, 
Nash-Moser theory, infinite dimensional Hamiltonian systems, 
 multiscale analysis.
\\[4mm]
{\it 2000AMS subject classification:} 35Q55, 37K55, 37K50.
\\[4mm]
This research was partially supported by 
PRIN 2015 ``Variational methods with applications to problems in mathematical physics and geometry"
and by  ANR-15-CE40-0001  ``BEKAM (Beyond KAM theory)".

\medskip

The authors thank Thomas Kappeler and the referees for several suggestions which improved the presentation of the manuscript. 
 
\tableofcontents

\mainmatter


\chapter{Introduction}

 In this  introductory chapter 
we present in detail the main results 
of this monograph (Theorems \ref{thm:main} and \ref{Prop:genericity}) concerning
the existence of  quasi-periodic solutions of multidimensional 
nonlinear wave equations with periodic boundary conditions, with a short description of 
the related literature. 
A comprehensive  introduction to 
KAM theory for PDEs is provided in Chapter \ref{Ch1.5}. 

\section{Main result and historical context}

We consider autonomous nonlinear wave equations\index{Nonlinear wave equation} (NLW)
\be\label{NLW1}
u_{tt} - \Delta u + V(x) u + g(x,u) = 0  \, , \quad x \in \T^d := \R^d \slash (2 \pi \Z)^d \, ,  
\ee
in any space dimension $ d \geq 1 $, 
where  $ V(x) \in C^{\infty} (\T^d, \R ) $ is a real valued {\it multiplicative} potential and 
the nonlinearity $ g \in C^\infty (\T^d \times \R, \R )$ has the form 
\be\label{nonlinearity}
g(x, u) = a(x) u^3 + O( u^4 ) 
\ee
with $ a(x) \in C^{\infty} (\T^d, \R ) $. 
We require that the elliptic operator\index{Multiplicative potential} $- \Delta + V(x) $ be positive definite, namely there exist $ \b >  0 $ such that 
\be\label{positive}
- \Delta + V(x)  > \b \,  {\rm Id}   \, . 
\ee
Condition \eqref{positive} is  satisfied, in particular, if the  potential $ V (x) \geq 0 $ and $ V (x) \not \equiv 0 $.
\smallskip

In this monograph we 
prove the existence of small amplitude time quasi-periodic solutions of  \eqref{NLW1}. 
We remind that a 
solution $ u(t,x) $ of \eqref{NLW1} is time quasi-periodic with frequency vector  
$ \omega \in \R^\nu $, $ \nu \in \N_+ = \{1, 2, \ldots \} $,  
if it has the form\index{Quasi-periodic solution}  
$$ 
u(t,x)  = U( \omega t , x) 
$$ 
where $ U : \T^\nu \times \T^d \to \R $ is a $ C^2 $-function and 
$ \omega \in \R^\nu $ is a nonresonant vector, namely 
$$ 
\omega \cdot \ell \neq 0 \, , \quad  \forall 
\ell \in \Z^\nu \setminus \{ 0 \} \, . 
$$ 
If $ \nu = 1 $ a solution of this form is time-periodic with period $ 2 \pi / \omega $.

\smallskip

Small amplitude solutions of the nonlinear wave equation 
\eqref{NLW1} are approximated, at a first degree of accuracy, by 
solutions of the 
linear wave equation\index{Linear wave equation} 
\be\label{Linear}
u_{tt} - \Delta u + V(x) u = 0 \, , \quad x \in \T^d \, . 
\ee
In the sequel, for any $f,g \in L^2 (\T^d, \C)$, we denote by $(f,g)_{L^2}$ the standard $L^2$-inner product
$$
(f,g)_{L^2} = \int_{\T^d} f(x) \overline{g(x)} \, dx \, . 
$$
There is a $ L^2$-orthonormal basis $ \{  \Psi_j \}_{j \in \N} $, $ \N = \{0, 1, \ldots \} $, of $L^2(\T^d)$ 
such that each $\Psi_j$ is an eigenfunction of the 
Schr\"odinger operator\index{Schrodinger operator}
$ - \D + V(x) $. More precisely, 
\be\label{auto-funzioni}
( - \Delta + V(x) ) \Psi_j (x) = \mu_j^2 \Psi_j (x)  \, , 
\ee
where the eigenvalues of $- \Delta + V(x)$
\be\label{eige-order}
0 < \beta \leq \mu_0^2 \leq \mu_1^2 \leq \ldots  \leq \mu_j^2 \leq  \ldots \, , \quad  \mu_j>0 \, , \quad  (\mu_j) \to + \infty \, ,   
\ee
are written in increasing order and with multiplicities.
The solutions of the linear wave equation \eqref{Linear} are given by the linear superpositions of {\it normal modes} oscillations, 
\be\label{soluzione-lineare-inf}
\sum_{j\in \N}  \a_j  \cos (  \mu_j t + \theta_j ) \Psi_j (x) \, , \qquad  \a_j  \, ,  \theta_j \in \R  \, . 
\ee
All the solutions \eqref{soluzione-lineare-inf} of \eqref{Linear} are periodic, 
or quasi-periodic, or almost periodic in time, 
with linear frequencies of oscillations $ \mu_j $, depending on the resonance properties of $ \mu_j $
(which depend on the potential $ V (x) $) and how many normal mode amplitudes 
$ \a_j $ are not zero. 
In particular,  if $ \a_j = 0 $ for any  index $ j $ 
except a finite set $ {\mathbb S} $ (tangential sites),
and the  frequency vector $ \bar \mu := ( \mu_j )_{j \in {\mathbb S}} $ is nonresonant, 
then the linear solutions \eqref{soluzione-lineare-inf} are 
quasi-periodic in time. 

The main question   we pose is the following:
\begin{itemize}
\item
Do small amplitude quasi-periodic solutions of the nonlinear wave equation \eqref{NLW1}
 exist?
\end{itemize}
The main results presented in this monograph, Theorems \ref{thm:main} and \ref{Prop:genericity},   state that: 
\begin{itemize} 
\item {\it small amplitude 
quasi-periodic  solutions  \eqref{soluzione-lineare-inf} 
of the linear wave equation  \eqref{Linear}, which are  supported on finitely many indices 
$ j \in {\mathbb S} $,  persist, slightly deformed, 
 as quasi-periodic solutions\index{Quasi-periodic solution}  of the 
nonlinear wave equation \eqref{NLW1},  with a frequency vector 
$ \omega $ close to $ \bar \mu $, for ``generic" potentials $ V(x) $ and coefficients $ a(x) $
 and ``most" amplitudes $ ( \a_j )_{j \in {\mathbb S}} $}. 
\end{itemize}

The potentials $ V(x) $ and the coefficients $ a(x) $ such that Theorem \ref{thm:main} 
holds are generic in a very strong sense; in particular  they are $ C^\infty $-dense, 
according to Definition \ref{Def:dense-open}, in the set 
$$ 
\big( {\cal P}  \cap C^\infty (\T^d) \big) \times C^\infty (\T^d)  \, ,
$$
where 
$ {\cal P} := \big\{ V (x) \in H^s (\T^d) \ : \  - \Delta + V(x) > 0  \big\} $,
see \eqref{C-infty-dense}.  

\medskip

Theorem  \ref{thm:main}  is a KAM (Kolmogorov-Arnold-Moser) 
type perturbative result. We construct recursively 
an embedded invariant torus $ \T^\nu \ni \vphi \mapsto i(\vphi) $ 
taking values in the phase space (that we describe carefully below),  supporting
quasi-periodic solutions of \eqref{NLW1} with frequency vector  $ \omega $ (to be determined). We employ
a modified Nash-Moser iterative scheme for the search of zeros 
$$
{\cal F} (\lambda; i ) = 0 
$$
of a nonlinear operator ${\cal F} $  
acting on scales of Sobolev spaces of maps $ i $, depending on a 
suitable parameter $ \lambda $, see Chapter \ref{sec:thmNM}.
 As in a Newton scheme, the core of the problem consists in the analysis of the 
linearized operators 
$$
d_i {\cal F} (\lambda; \underline{i}) 
$$
at any approximate solution $\underline{i} $ at each step of the iteration, see Section \ref{sec:NNM}.
The main task is to prove that  $ d_i {\cal F} (\lambda; \underline{i})  $
admits an  
approximate inverse, for most values of the parameters, 
which satisfies suitable 
quantitative tame estimates in high Sobolev norms. The approximate inverse will
be unbounded, i.e. it  
loses derivatives, due to  the presence of {\it small divisors}.  
As we shall describe in detail in Section \ref{sec:ideas}, 
the construction of an approximate inverse for the linearized operators 
obtained from \eqref{NLW1}
is a subtle problem. Major difficulties 
come from 
complicated resonance phenomena between the frequency vector $ \omega $ of the 
expected quasi-periodic solutions  and the  multiple normal mode frequencies of oscillations,
shifted by the nonlinearity, and the fact that the normal mode eigenfunctions 
$ \Psi_j (x) $ are not ``localized close to the exponentials".

\smallskip
 
 We now make a short historical introduction to 
KAM theory for partial differential equations, 
that we shall expand in Chapter \ref{Ch1.5}.
 
\smallskip 
 
As we already mentioned in the preface, 
in small divisors problems for PDEs, as  \eqref{NLW1}, 
the space dimension $ d = 1 $ or $ d \geq 2 $ 
makes a fundamental difference, due to the very different 
properties of the eigenvalues and eigenfunctions of the 
Schr\"odinger operator $ - \Delta + V(x) $
on $ \T^d $ for $ d = 1 $ and $ d \geq 2 $. If the space dimension
$ d = 1 $ we also call $ - \pa_{xx} + V(x) $ a Sturm-Liouville operator. 

\smallskip

The first KAM  existence results of quasi-periodic solutions\index{Quasi-periodic solution} 
were proved by Kuksin \cite{Ku}, see also \cite{Ku1}, and Wayne \cite{W1} 
 for $ 1$-$d$  wave and  Schr\"odinger (NLS) equations on the interval $ x \in [0, \pi] $
with  Dirichlet boundary conditions and analytic  nonlinearities, see \eqref{NLW1-L}-\eqref{NLS1}. These pioneering theorems were limited 
to Dirichlet boundary conditions because the eigenvalues $ \mu_j^2 $ of 
the Sturm-Liouville operator $ - \pa_{xx}  + V(x) $ had to be simple. 
Indeed the 
KAM scheme in \cite{Ku1}, \cite{W1}, see also \cite{Po2},
reduces the linearized equations along the 
iteration
to a diagonal form, with coefficients constant in  time.
This process requires 
``second-order Melnikov" non-resonance conditions, which concern  lower bounds
for differences among the linear frequencies. 
In these papers the potential 
$ V(x) $ is used as a parameter to impose 
non-resonance conditions. 
Once the linearized PDEs obtained along the iteration are reduced to diagonal,
constant in time, form, 
it is  easy to prove that the corresponding 
linear operators are invertible, for most values of the parameters, 
with good estimates of
their inverses  in high norms (with of course a loss of derivatives). 
We refer to Section \ref{sec:red} for a more detailed explanation of the KAM reducibility approach. 

Subsequently  these results have been extended by  
P\"oschel \cite{Po3} for  parameter independent nonlinear Klein-Gordon equations like \eqref{NLWKAM}, 
and by Kuksin-P\"oschel \cite{KP2} for NLS equations like \eqref{NLSKAM}.
A major novelty of these papers was the use of Birkhoff normal form techniques 
to verify (weak) non-resonance conditions among the perturbed 
frequencies, tuning the amplitudes of the solutions as parameters. 

\smallskip

In the case  $ x \in \T $, the eigenvalues of the Sturm-Liouville operator $  - \pa_{xx}  + V(x) $ are 
asymptotically double, and therefore the previous 
second order Melnikov non-resonance conditions are violated. In this case the first existence results were obtained
by Craig-Wayne \cite{CW}  for time periodic solutions of analytic 
nonlinear Klein-Gordon equations  
(see also \cite{C} and \cite{BB06} for completely resonant
wave equations), and then extended by Bourgain 
\cite{Bo1} for time quasi-periodic solutions. The proofs are based on 
a Lyapunov-Schmidt bifurcation approach and a Nash-Moser implicit function iterative scheme. 
The key point of these papers is that there is no 
diagonalization of
the linearized equations at each step of the Nash-Moser iteration.
The advantage is to require only minimal non-resonance conditions 
which are easily verified for  PDEs also in presence of multiple frequencies 
(the second order Melnikov non-resonance conditions are not used). 
On the other hand, 
a difficulty of this approach is that
 the linearized equations obtained along the iteration are variable coefficients PDEs. 
 As a consequence it is hard to prove that the corresponding 
linear operators are invertible with  estimates of
their inverses  in high norms, sufficient to imply the convergence of the iterative scheme. 
Relying on a ``resolvent" type analysis inspired by the work  of 
Fr\"olich-Spencer \cite{FS} in the context of Anderson localization,  Craig-Wayne \cite{CW} were able to solve 
this problem for time periodic solutions in $ d  = 1 $, 
and  Bourgain in \cite{Bo1}  also for quasi-periodic solutions. 
Key properties of this approach are:
\begin{description}
\item ($i$) ``separation properties" between  singular sites, namely the Fourier indices 
$ (\ell, j) $ of the small divisors $ | (\om \cdot \ell )^2 - j^2 | \leq C $ in the case of (NLW); 
\item 
($ii$) ``localization" of the eigenfunctions of the Sturm-Liouville operator  
$ - \pa_{xx} + V(x) $\index{Eigenfunctions of Sturm-Liouville operator}  with respect to the exponential basis $ (e^{\ii k x})_{k \in \Z} $, namely
that the Fourier coefficients $ (\hat \Psi_{j})_k $ 
converge  rapidly to zero when $ | |k| - j | \to  \infty $. 
This property is always true if $ d = 1 $, see e.g. \cite{CW}.   
\end{description}

Property ($ii$) implies that the matrix which represents, in the eigenfunction basis, 
the multiplication operator defined by an analytic (resp. Sobolev) 
function has an exponentially (resp. polynomially) fast decay  off the diagonal. 
Then the ``separation properties"  ($i$)  
imply a  very ``weak interaction" between the singular sites.
If the singular sites were ``too many", the inverse operator would be ``too unbounded" 
to prevent the convergence of the iterative scheme. 
This approach is particularly promising in the presence of multiple normal mode frequencies
and it constitutes  the basis of the present monograph. We describe it in more detail in Section 
\ref{sec:MULTI}.

\smallskip

Later on, Chierchia-You \cite{CY} were able to extend  the KAM reducibility approach
to prove existence and stability of small amplitude quasi-periodic solutions of  $1$-$d$ NLW on $ \T $ with an external potential. We also mention 
the KAM reducibility results  in Berti-Biasco-Procesi \cite{Berti-Biasco-Procesi-Ham-DNLW}-\cite{Berti-Biasco-Procesi-rev-DNLW} for $ 1$-$d$ derivative wave equations.

\smallskip

In the case  the space dimension $ d $ is greater or equal to two,  major difficulties are:
\begin{enumerate}
\item \label{dif10}
the eigenvalues\index{Eigenvalues of Schr\"odinger operator} $ \mu_j^2 $ of $ - \Delta + V(x) $\index{Schr\"odinger operator} in \eqref{auto-funzioni} may be of high multiplicity,
or not sufficiently separated from each other in a suitable quantitative sense, required 
by the 
perturbation theory developed for $1$-$d$-PDEs;
\item  \label{dif20}
the eigenfunctions \index{Eigenfunctions of Sturm-Liouville operator} $ \Psi_j (x) $ of $ - \Delta + V(x) $ may be not 
``localized" with respect to the exponentials, see \cite{FKT}. 
\end{enumerate}

As  discussed in the preface, 
if $ d \geq 2 $, the first KAM  existence result for nonlinear wave equations 
has been proved  for time periodic solutions by Bourgain \cite{B-Gafa}, 
see also the extensions in \cite{BBARMA}, \cite{BBP10},  \cite{GP}. 
Concerning quasi-periodic solutions in dimension $ d $ greater or equal to $2 $, the first existence
result was proved by Bourgain in Chapter 20 of \cite{B5}: it deals with 
wave type equations of the form\index{Nonlinear wave equation}
$$
u_{tt} - \Delta u + M_\sigma u + \e F' (u) = 0  \, ,
$$
where 
$ M_\sigma = {\rm Op}( \sigma_k )$ is a Fourier multiplier 
supported on finitely many sites $ {\mathbb E} \subset \Z^d  $, i.e. 
$  \sigma_{k} = 0  $, $  \forall  k \in  \Z^d \setminus {\mathbb E} $.
The $ \sigma_k $, $ k \in {\mathbb E}  $, are   
used as external parameters, and $ F $ is a polynomial nonlinearity, with 
$ F' $ denoting the derivative of $ F $. 
Note that the linear equation 
$$ 
u_{tt} - \Delta u + M_\sigma u = 0  
$$ 
is diagonal in 
the exponential basis $ e^{\ii k \cdot x }$, $ k \in \Z^d $, unlike the linear wave equation \eqref{Linear}. 
We also mention the paper by Wang \cite{Wang} for the completely resonant NLS 
equation \eqref{c-res:NLS} and 
the Anderson localization result of Bourgain-Wang \cite{BW1} for  time quasi-periodic
random linear Schr\"odinger and wave equations. 

As already mentioned, a major difficulty of this approach 
is that the linearized equations obtained along the iteration
are   PDEs with variable coefficients. 
A key property which plays a fundamental role in 
\cite{B5} (as well as in previous papers as 
 \cite{B3} for NLS) 
 for proving estimates for the inverse of  linear operators 
$$
{\it \Pi}_N \big( (\om \cdot \pa_{\vphi})^2  - \Delta + M_\sigma  + \e b(\vphi, x)  \big)_{|{\cal H}_N}  \, ,  
$$
(see  \eqref{QP-lin:OP})
is that the matrix which represents the multiplication operator for a smooth function 
$ b(x) $ in the exponential basis $ \{ e^{\ii k \cdot x } \} $, $ k \in \Z^d $, 
has a sufficiently fast off-diagonal decay. 
Indeed the multiplication operator is represented in Fourier space by 
the T\"oplitz matrix $ ({\hat b}_{k - k' })_{k, k' \in \Z^d} $, 
with entries given by the Fourier coefficients $ {\hat b}_{K} $ of the function
$ b(x) $, constant on the diagonal 
$ k - k' = K $. The smoother the function $ b(x) $ is, the faster is the decay of 
$  {\hat b}_{k-k'} $ as
$ |k - k | \to + \infty $.  
We refer to Section \ref{sec:MULTI} for more explanations on this approach.

Weaker forms of this property, as for example those required in 
Berti-Corsi-Procesi \cite{BCP}, \cite{BP}
may be sufficient  for dealing with the eigenfunctions of $ - \Delta $ on compact Lie groups.
However, 
in general, the  elements $ (\Psi_j, b(x) \Psi_{j'})_{L^2} $ 
of the matrix which represents  the multiplication operator  with respect to  the basis of the 
eigenfunctions $ \Psi_j (x) $ of $ - \Delta + V(x) $ on $ \T^d $ (see  \eqref{auto-funzioni})  
do not decay sufficiently fast, if $ d \geq 2 $. 
This was proved by
Feldman, Kn\"orrer, Trubowitz in \cite{FKT} and it is the difficulty mentioned above 
in item \ref{dif20}. 
We remark that 
weak properties of localization have been proved by Wang \cite{ Wang2} in $ d =2 $ for 
potentials $ V(x) $ which are trigonometric polynomials.   

\smallskip
In the present  monograph  
we shall not use 
 properties of localizations of the eigenfunctions $ \Psi_j (x) $. 
A major reason why we are able to avoid the use of such properties
is that our Nash-Moser iterative scheme 
requires only  very weak tame estimates
for the approximate inverse of the linearized operators  as stated in \eqref{LN-1tame-s}, 
see the end of Subsection \ref{subsec:Ch4}. 
Such conditions are
close to the optimal ones, as a
famous counterexample of Lojiaciewitz-Zehnder in \cite{LZ} shows. 

\smallskip

The properties of the 
exponential basis $ e^{\ii k \cdot x } $, $ k \in \Z^d $,  play a key role 
also
for developing the KAM perturbative diagonalization/reducibility techniques.
Indeed,  no reducibility results are available so far for multidimensional  PDEs
in presence of a multiplicative potential which is not small.  
Concerning 
 higher space dimensional PDEs we refer to the results in Eliasson-Kuksin 
 \cite{EK} for the NLS equation \eqref{EKAnn}
 with a convolution potential on $ \T^d $ used as a parameter,  
Geng-You \cite{GY1.5} and  Eliasson-Gr\'ebert-Kuksin \cite{EGK1} 
 for  beam equations with a constant mass potential,  Procesi-Procesi \cite{PP} for 
 the completely resonant NLS equation\index{Nonlinear Schr\"odinger equation} \eqref{c-res:NLS}, 
 Gr\'ebert-Paturel  \cite{GrP} for the Klein-Gordon equation \eqref{GPSd} on 
 $ {\mathbb S}^d $ and  Gr\'ebert-Paturel   \cite{GrP1} for multidimensional harmonic oscillators.

On the other hand, no  reducibility results for NLW  on $ \T^d $ are known so far.  Actually 
a serious difficulty which appears is the following: 
the infinitely many second order Melnikov non-resonance conditions required by the KAM-diagonalization approach 
are strongly violated by the linear unperturbed 
frequencies  of oscillations of the Klein-Gordon equation 
$$ 
u_{tt} - \Delta u + m u = 0 \, , 
$$ 
see \cite{EGK}. A key difference with respect to the Schr\"odinger equation is that the 
linear frequencies of  the wave equations are $ \sim | k | $, $ k \in \Z^d $, 
while for the NLS and the  beam equations, they are $ \sim |k|^2 $, respectively $ \sim |k |^4 $, 
and $ | k |^2 $, $ |k|^4 $ are integer. 
Also for the multidimensional harmonic oscillator, the linear frequencies are, up to a translation,  integer numbers. 
Although no reducibility results are known so far for the NLW equation, 
a result of ``almost" reducibility for linear quasi-periodically forced
Klein-Gordon equations has been presented in \cite{E17}, \cite{EGK}. 

\smallskip

Existence of quasi-periodic solutions\index{Quasi-periodic solution}  for  wave equations on $ \T^d $ 
with a time-quasi periodic  forcing nonlinearity\index{Nonlinear wave equation}  of class $C^\infty$ 
or $C^q$ with $q$ large enough,
\be\label{NLW:forced-in}
u_{tt} - \Delta u + V(x) u = \e f(\om t, x, u ) \, , \quad x \in \T^d \, , 
\ee
has been proved in Berti-Bolle \cite{BB12},  extending the multiscale approach of 
Bourgain \cite{B5}.
The forcing frequency vector $ \omega $, which in \cite{BB12} 
 is constrained to a fixed direction $ \omega = \l \bar \om $, with $ \l \in [1/2, 3/2] $,  plays the role of 
 an external parameter.
In \cite{BCP} a corresponding result has been extended
for NLW equations
 on compact Lie groups, 
in \cite{BK} for Zoll manifolds, in \cite{BMa} for general flat tori, 
and  in \cite{CM} for forced Kirkhoff equations.  

Existence of quasi-periodic solutions for autonomous non-linear Klein Gordon
\index{Klein Gordon} equations 
\be\label{KG-NON}
u_{tt} - \Delta u + u + u^{p+1} + h.o.t.  = 0 \, , \quad p \in \N \, , \quad x \in \T^d \, , 
\ee
have been recently presented  by Wang \cite{Wang1}, relying on 
a bifurcation analysis to study the modulation of the frequencies induced by the nonlinearity $ u^{p+1} $, 
and  multiscale methods of \cite{B5} for implementing a Nash-Moser iteration.  
The result proves the continuation of quasi-periodic solutions supported on ``good" tangential sites. 

\smallskip

Among all  the works discussed above, the papers \cite{BB12}-\cite{BBo10} on 
 forced NLW and NLS\index{Nonlinear Schr\"odinger equation} equations
are related most closely to  the present monograph.  
The passage to prove KAM results for autonomous nonlinear  wave equations 
with a multiplicative potential as \eqref{NLW1} is a non trivial task. 
 It  requires  a bifurcation analysis which distinguishes 
the tangential directions where the major part of the oscillation of the quasi-periodic solutions
takes place, and the normal ones, see the form \eqref{sol:uel} of the quasi-periodic solutions
proved in Theorem \ref{thm:main}. 
When the multiplicative potential\index{Multiplicative potential}   $ V(x) $ changes, both the tangential and the normal frequencies
vary simultaneously in an intricate way (unlike the case of the convolution potential). 
This makes difficult to verify
 the  non-resonance conditions 
 required by the Nash-Moser iteration.
In particular, the choice of  the 
parameters adopted in order to fulfill all these conditions is relevant.   
  In this  monograph
we choose any finite set $ {\mathbb S} \subset \N_+ $ of tangential sites, we  
fix the potential $ V(x) $ and the coefficient function $ a(x) $ appearing in the nonlinearity \eqref{nonlinearity}
(in such a way that generic non-resonance and non-degeneracy conditions hold, see Theorem \ref{Prop:genericity}) 
and then we prove, in Theorem
\ref{thm:main}, the existence of quasi-periodic solutions of \eqref{NLW1}
 for most values of the one dimensional internal parameter $ \lambda  $ introduced in \eqref{frequencycolinear}. 
The parameter $ \lambda $
 amounts just to a {\it time rescaling} of the frequency vector  $ \omega $. 
We also deduce a density result for the frequencies  of the quasi-periodic 
solutions close to the unperturbed vector $ \bar \mu $.
We shall explain more in detail the choice of this parameter in Section \ref{sec:ideas}. 

\section{Statement of the main results}

In this section we state in detail the main results of this monograph, which are Theorems \ref{thm:main} and \ref{Prop:genericity}.

Under the rescaling $ u \mapsto \e u $, $ \e >  0 $,  the equation \eqref{NLW1} is transformed into
the nonlinear wave equation\index{Nonlinear wave equation}
\be\label{NLW2}
u_{tt} - \Delta u + V(x) u +  \e^2 g(\e, x, u ) = 0
\ee
with the  $ C^\infty $ nonlinearity 
\be\label{nonlinearity:gep}
 g(\e, x, u ) := \e^{-3} g(x, \e u) = a(x) u^3 + O( \e u^4 ) \, . 
\ee
Recall that we list the eigenvalues $ (\mu_j^2)_{j \in \N} $ of $ - \Delta + V(x) $ 
in increasing order,  see \eqref{eige-order}, 
and that we choose a corresponding $L^2$-orthonormal sequence of 
eigenfunctions $ \{  \Psi_j \}_{j \in \N} $.

We choose arbitrarily a finite set of indices $ {\mathbb S} \subset \N  $, 
called the ``tangential sites". We denote by $ \es \in \N $ the cardinality of $ {\mathbb S}$
and we list the  tangential sites in increasing order, $ {j_1}< \ldots < {j_\es}  $.
We look for quasi-periodic solutions of \eqref{NLW2}
which are perturbations  of normal modes oscillations 
supported on $ j \in {\mathbb S} $.   
We denote by
\be\label{unp-tangential}
\bar \mu :=  
(\mu_j)_{j \in {\mathbb S}} = (\mu_{j_1}, \ldots, \mu_{j_\es} ) 
\in \R^\es \, , \quad  \mu_j >  0  \, , 
\ee
the 
frequency vector of the  quasi-periodic  solutions  
\be\label{soluzione-lineare}
\sum_{j\in {\mathbb S}}  \mu_j^{-1/2} \sqrt{2 \xi_j}  \cos (  \mu_j t 
) \Psi_j (x) \, , \qquad  \xi_j > 0  \, , 
\ee
of the linear wave equation \eqref{Linear}. The components of $\bar \mu$ are called the 
unperturbed tangential frequencies and $ (\xi_j)_{j \in {\mathbb S}} $ the unperturbed actions.
We shall call the indices   in the complementary set $  {\mathbb S}^c := \N \setminus {\mathbb S} $ 
the ``normal" sites,  and  the corresponding 
$ \mu_j $, $  j \in {\mathbb S}^c $, the unperturbed  ``normal" frequencies.

\smallskip

Since  \eqref{NLW2} is an {\it autonomous}  PDE,  the frequency vector 
$ \om \in \R^\es $ of its expected quasi-periodic solutions $ u(\om t, x ) $
is an unknown, that we introduce as an explicit parameter in the equation, looking 
for  solutions $ u(\vphi, x)$, $ \vphi = (\vphi_1, \ldots, \vphi_\es ) \in \T^\es $, of 
\be\label{NLW2QP}
(\om \cdot \partial_\vphi)^2 u - \Delta u + V(x) u  + \e^2 g(\e, x, u ) = 0 \, . 
\ee
The frequency vector $ \om \in \R^\es $ of the expected quasi-periodic solutions of 
\eqref{NLW2} will be  $ O(\e^2) $-close to 
the unperturbed tangential frequency vector $ \bar \mu $ in  \eqref{unp-tangential}, 
see  \eqref{frequencycolinear}-\eqref{def omep}. 

Since the nonlinear wave equation  \eqref{NLW2} is time-reversible\index{Reversible PDE} (see Appendix \ref{App0}), 
it makes sense to look  for solutions of \eqref{NLW2} which are even in $ t $.
Since \eqref{NLW2} is autonomous, additional solutions are obtained from 
these even ones by time translation.   
Thus we look for solutions $ u(\vphi, x) $ of \eqref{NLW2QP} which are 
 {\it even} in  $\vphi $. This 
induces a small simplification in the proof, 
see Remark \ref{rem:rev}.  

In order to prove, for $ \e $ small enough, the existence of solutions of 
\eqref{NLW2}
close to the solutions \eqref{soluzione-lineare} 
of the linear wave equation \eqref{Linear},   
we first require  non-resonance conditions for  the  unperturbed linear frequencies $ \mu_j $, $ j \in \N $, which will be 
 verified by generic  potentials $ V(x) $, see Theorem \ref{Prop:genericity}.
\\[2mm]
{\bf Diophantine and first order Melnikov non-resonance conditions.}
We  assume that  
\begin{itemize}
\item
the tangential frequency vector 
$ \bar \mu $ in \eqref{unp-tangential} is Diophantine\index{Diophantine vector}, i.e.  
for some constants $ \g_0 \in (0,1) $, $ \t_0 >  \es -1 $, 
\be\label{diop}
| \bar \mu \cdot \ell  | \geq \frac{\g_0}{ \langle \ell \rangle^{\tau_0}} \, , \quad \forall \ell  \in \Z^\es \setminus \{0\} \, , \qquad 
\langle \ell \rangle := \max\{1, | \ell | \} \, , 
\ee
where $ | \ell | := \max\{|\ell_1|, \ldots, | \ell_\es | \} $. Note that  \eqref{diop} implies,
in particular, that the $ \mu_j^2 $, $ j \in {\mathbb S} $, 
are simple eigenvalues of $ - \Delta + V(x) $.
\item the unperturbed ``first order Melnikov"\index{First  Melnikov non-resonance conditions} 
non-resonance conditions hold:
\be\label{1Mel}
| \bar \mu \cdot \ell + \mu_j | \geq \frac{\g_0}{ \langle \ell \rangle^{\tau_0}} \, , \quad \forall \ell \in \Z^\es \, , 
\ j \in {\mathbb S}^c  \, .
\ee
\end{itemize}
 The non-resonance conditions \eqref{diop}, \eqref{1Mel} imply, in particular,  that 
the linear equation \eqref{Linear}
has no other quasi-periodic solutions  with  frequency $ \bar \mu $, even in $ t $, 
except the trivial ones \eqref{soluzione-lineare}. 

\smallskip

In order to prove ``separation properties" of the small divisors
as required by the multiscale analysis that we perform in Chapter \ref{sec:multiscale}, we require,  
as in \cite{BB12}, that 
\begin{itemize}
\item
the tangential frequency vector 
$  \bar \mu  $ in \eqref{unp-tangential} satisfies the quadratic Diophantine  condition\index{Quadratic Diophantine condition} 
\be\label{NRgamma0}
\Big| n + \sum_{i, j \in {\mathbb S}, i \leq j}  p_{ij} \mu_i \mu_j \Big| \geq  \frac{\g_0}{ \langle p \rangle^{ \t_0}} \, , \quad 
\forall (n, p) \in \Z \times \Z^{\frac{\es(\es+1)}{2}} \setminus \{(0,0)\} \, . 
\ee
\end{itemize}
The non-resonance conditions  \eqref{diop}, \eqref{1Mel} and  \eqref{NRgamma0}  
are assumptions  on the potential $ V(x) $, which are ``generic" in the sense of Kolmogorov measure. In 
\cite{KaKu},   \eqref{diop}, \eqref{1Mel} are proved to hold for ``most" potentials.
Genericity results are stated in Theorem \ref{Prop:genericity}, proved in Chapter \ref{section:gener1}.

We emphasize that along the monograph the constant $ \g_0 \in (0,1) $  in \eqref{diop}, 
\eqref{1Mel}, \eqref{NRgamma0} is regarded as fixed, and we shall often omit to track 
its dependence  
in the estimates. 
\\[2mm]
{\bf Birkhoff matrices.} 
We are interested in quasi-periodic solutions of \eqref{NLW2} which bifurcate for small $\e>0$
from a solution of the form \eqref{soluzione-lineare} of the linear wave equation.
In order to prove their existence, it is important to know precisely how the 
tangential and the normal frequencies 
change with respect to the unperturbed actions 
$ (\xi_j)_{j \in {\mathbb S}} $, under the effect of the nonlinearity $ \e^2 a(x) u^3 + O( \e^3 u^4 ) $. 
This is described in terms of the 
``Birkhoff" matrices\index{Birkhoff matrices}
\be\label{def:AB}
\Ab := (\mu_k^{-1} G^j_k \mu_j^{-1})_{j, k \in {\mathbb S}} \, , \quad 
\Bb :=   (\mu_j^{-1} G^k_j  \mu_k^{-1} )_{j \in {\mathbb S}^c, k \in {\mathbb S}} \, , 
\ee
where, for any $ j, k \in \N $,  
\be\label{def G}
G^j_k := G^j_k (V, a) :=  \begin{cases}
 (3/2) ( \Psi_j^2, a(x) \Psi_k^2 )_{L^2}, \quad {j \neq k} \, , \cr
 (3/4) ( \Psi_j^2, a(x) \Psi_j^2 )_{L^2},  \quad {j = k} 
\end{cases}
\ee 
and $ \Psi_j (x) $ are the eigenfunctions of $ - \Delta + V(x) $ introduced in \eqref{auto-funzioni}. 
Note that the matrix $ (G^j_k ) $ 
depends on the coefficient $ a(x) $ and the eigenfunctions $ \Psi_j $,
thus on  the potential $ V(x) $.
The $\es \times \es $ symmetric matrix\index{Twist matrix} $ \Ab $ is  called the ``twist"-matrix. 
The matrices $ \Ab, \Bb $ describe the shift of the tangential and  normal frequencies 
induced by the nonlinearity $ a(x) u^3 $ as  
they appear in the fourth order Birkhoff normal form of
\eqref{NLW1}-\eqref{nonlinearity}. 
Actually, we prove in Section \ref{sec:shifted-tan} that, 
up to terms $ O(\e^4 )$,  the tangential frequency vector $ \om $
of a small amplitude quasi-periodic solution of \eqref{NLW1}-\eqref{nonlinearity} close to \eqref{soluzione-lineare} is 
given by the action-to-frequency map
\be\label{act-to-fre}
\xi \mapsto \bar \mu + \e^2 \Ab  \xi  \, , \quad  \xi \in \R_+^\es \, ,  \quad \R_+ := (0, + \infty) \, .
\ee
On the other hand the perturbed normal frequencies are shifted  by the matrix $ \Bb $ as described in Lemma 
\ref{lem:ns}. 
We assume that\index{Twist condition} 
\begin{itemize} 
\item 
{\bf (Twist condition)}
\be\label{A twist}
{\rm det} \, \Ab \neq  0 \, ,
\ee
\end{itemize}
and therefore the action-to-frequency map in \eqref{act-to-fre} is invertible. 
The non-degeneracy, or ``twist"-condition \eqref{A twist}, is satisfied 
for generic   $ V(x) $ and  
$ a(x) $, as stated in Theorem \ref{Prop:genericity} (see in particular Corollary \ref{twist2} and 
Remark \ref{cor:twist}). 
\\[2mm]
{\bf 
Second order Melnikov non-resonance conditions.}
We also assume second order  Melnikov non-resonance conditions which 
concern only {\it finitely} many  
unperturbed normal frequencies. 
We have first to introduce an important decomposition of the normal indices $ j \in {\mathbb S}^c $.  
Note that, since $ \mu_j \to + \infty $,  the indices $ j \in {\mathbb S}^c $ such that 
$  \mu_j -  (\Bb \Ab^{-1} \bar \mu)_j <  0 $ are finitely many. Let $  \gap \in \R $ be defined by
\be\label{def:gs}
- \gap := \min \big\{   \mu_j -  ( \Bb \Ab^{-1} \bar \mu)_j  \, , j \in {\mathbb S}^c \big\} \, . 
\ee
We decompose  the set of normal indices as
\be\label{taglio:pos-neg-0}
{\mathbb S}^c  = {\mathbb F}  \cup {\mathbb G} \, , \quad 
{\mathbb G} := {\mathbb S}^c \setminus {\mathbb F} \, , 
\ee
where
\be
\begin{aligned}\label{taglio:pos-neg}
& {\mathbb F} := \Big\{  j \in {\mathbb S}^c \, : \, | \mu_j -  ( \Bb \Ab^{-1} \bar \mu)_j| \leq  \gap \Big\} \, , 	\\
& {\mathbb G}  := \Big\{  j \in {\mathbb S}^c \, : \, \mu_j -  ( \Bb \Ab^{-1} \bar \mu)_j  >  \gap \Big\} \, .
\end{aligned}
\ee
The set ${\mathbb F} $ is always finite,  and  is empty if $ \gap < 0 $. 
The relevance of the decomposition \eqref{taglio:pos-neg-0}  of the normal sites, concerns 
the variation of the normal frequencies with respect to the length of the tangential frequency vector, as we describe 
in \eqref{derivata lam} below, see also Lemma \ref{choice:M}.
If all the numbers  
 $ \mu_j -  ( \Bb \Ab^{-1} \bar \mu)_j $, $j \in {\mathbb S}^c$, were positive, then, by  \eqref{derivata lam}, 
the growth rates 
of the eigenvalues of the linear operators which we need to invert in the Nash Moser
iteration would all have the same sign. This would allow to obtain measure estimates by simple arguments, as in
\cite{BB12}, \cite{BBo10}, where the forced case is considered.
In general $  \gap > 0 $ and  we shall be able to
decouple, for most values of the
parameter $ \l $, 
the linearized operators
obtained at each step of the nonlinear Nash-Moser iteration,
acting in the normal subspace $ H_{\mathbb S}^\bot $, 
 along $ H_{\mathbb F} $ and its orthogonal $ H_{\mathbb F}^\bot $ $ = H_{\mathbb G}$,
 defined in \eqref{spliFG}. 
 We discuss the relevance of this decomposition in Section \ref{sec:ideas}.

\smallskip

We assume the following\index{Second Melnikov non-resonance conditions} 
\begin{itemize}
\item 
unperturbed  ``second order Melnikov" non-resonance conditions (part $1$):
\begin{align}
& | \bar \mu \cdot \ell +  \mu_j -  \mu_k | \geq \frac{\g_0}{ \langle \ell \rangle^{\tau_0}} \, , \ 
\forall (\ell, j, k) \in \Z^\es \times {\mathbb F} \times {\mathbb S}^c \, , \ \  (\ell, j, k) \neq (0,j,j) \, ,  \label{2Mel+} \\
& | \bar \mu \cdot \ell +  \mu_j +  \mu_k | \geq \frac{\g_0}{ \langle \ell \rangle^{\tau_0}} \, , \  
\forall (\ell, j, k) \in \Z^\es \times {\mathbb F} \times {\mathbb S}^c \, .   \label{2Mel} 
\end{align}
\end{itemize}

Note that  \eqref{2Mel+} implies, in particular, that the finitely many 
eigenvalues $ \mu_j^2 $ of $ - \Delta + V(x) $, $ j \in {\mathbb F} $, are simple
(clearly all the other  eigenvalues $ \mu_j^2 $, $ j \in {\mathbb G} $, could be highly degenerate).  

In order to verify a key positivity property for the variations of the restricted linearized operator 
with respect to $ \lambda $ (Lemma \ref{pos:def-var}), 
we  assume further 
\begin{itemize}
\item 
unperturbed  ``second order Melnikov" non-resonance conditions (part $2$):
\begin{align}
& | \bar \mu \cdot \ell +  \mu_j -  \mu_k | \geq \frac{\g_0}{ \langle \ell \rangle^{\tau_0}} \, , \quad 
\forall (\ell, j, k) \in  \Z^\es \times ({\mathbb M} \setminus {\mathbb F}) \times {\mathbb S}^c  \, , \label{2Mel rafforzate+} \\ 
& \qquad \qquad \qquad  \qquad \qquad \  \qquad (\ell, j, k) \neq (0,j,j) \, ,  \nonumber \\
& | \bar \mu \cdot \ell +  \mu_j +  \mu_k | \geq \frac{\g_0}{ \langle \ell \rangle^{\tau_0}} \, , \quad 
\forall (\ell, j, k) \in  \Z^\es \times ({\mathbb M} \setminus {\mathbb F}) \times {\mathbb S}^c  \, ,   \label{2Mel rafforzate}
\end{align}
where 
\be\label{set:M}
{\mathbb M} := \big\{ j \in {\mathbb S}^c \, : \, |j| \leq C_1 \big\} 
\ee
and  the constant $ C_1 := C_1 (  V, a ) > 0 $ is taken large enough, so that 
$ {\mathbb F} \subset {\mathbb M} $ and  \eqref{cond:suM} holds. 
\end{itemize}
Clearly the conditions \eqref{2Mel+}-\eqref{2Mel}  and \eqref{2Mel rafforzate+}-\eqref{2Mel rafforzate}
 could have been written together, requiring such conditions for $ j \in {\mathbb M}$, without distinguishing 
 the cases $ j \in {\mathbb F} $ and $ j \in {\mathbb M} \setminus  {\mathbb F} $.
 However, for conceptual clarity, in view of their different role in the proof, we prefer to state them separately. 
The above conditions \eqref{2Mel+}-\eqref{2Mel rafforzate}
on the unperturbed frequencies allow to perform one step of averaging and so 
to diagonalize, up to $ O(\e^4) $, the normal frequencies supported on
$ {\mathbb M} $, see Proposition \ref{prop:op-averaged}.
This is the only step where conditions \eqref{2Mel rafforzate+}-\eqref{2Mel rafforzate}
play a role. Conditions \eqref{2Mel+}-\eqref{2Mel} are used also in the
splitting step of Chapter  \ref{sec:splitting}, see Lemma \ref{lemma:measure1}.

Conditions \eqref{2Mel+}-\eqref{2Mel rafforzate} depend on 
the potential $ V(x) $ and also on the coefficient $ a(x) $, 
because the constant $ C_1 $ in \eqref{set:M} (hence the set $\mathbb M$) depends 
on $\|a \|_{L^\infty}$ and $\|{\cal A}^{-1} \|$. 
Given $(V_0,a_0)$ such that the matrix $\Ab$ defined in \eqref{def:AB}  is invertible and $s>d/2$, the set 
$\mathbb M$ can be chosen constant in some open neighborhood $U$ of
$(V_0,a_0)$ in $H^s$. In $U$, conditions
  \eqref{2Mel+}-\eqref{2Mel rafforzate} are
generic for $ V (x) $, as it is 
 proved in Chapter \ref{section:gener1} (see  Theorem \ref{Prop:genericity}). 
\\[2mm]
{\bf Non-degeneracy conditions.}
We also require the following {\it finitely} many 
\begin{itemize}
\item
non-degeneracy conditions:  
\begin{align}\label{non-reso}
& (\mu_j - [ \Bb \Ab^{-1} \bar \mu]_j) - (\mu_k - [ \Bb \Ab^{-1} \bar \mu]_k) \neq 0 \, , \quad \forall j, k \in {\mathbb F} \, ,
 \ j \neq k \, , 
\\
& (\mu_j - [ \Bb \Ab^{-1} \bar \mu]_j) + (\mu_k - [ \Bb \Ab^{-1} \bar \mu]_k) \neq 0 \, , \quad \forall j, k \in {\mathbb F} \, , 
\label{non-reso1}
\end{align}
where $ \Ab $ and $ \Bb $ are the Birkhoff matrices\index{Birkhoff matrices} defined in \eqref{def:AB}. 
\end{itemize}
Such assumptions are  similar 
to the non-degeneracy conditions required for the continuation of elliptic tori for finite dimensional systems in \cite{El88}, 
\cite{P89} and for PDEs in \cite{KP2}, \cite{Po3}, \cite{BBi11}. 
Note that the finitely many non-degeneracy conditions \eqref{non-reso} depend on 
 $ V (x) $ and $ a(x) $ and  we  prove in Theorem \ref{Prop:genericity} 
that they are generic in  $ (V, a )  $. 
\\[2mm]
{\bf Parameter.}
We now introduce 
the $ 1$-dimensional parameter that we shall use to perform the measure estimates. 

In view of \eqref{act-to-fre}
the frequency vector $ \om $ has to belong to the cone of the ``admissible" frequencies\index{Admissible frequencies}
$ \bar \mu + \e^2 \Ab ( \R_+^\es ) $, more precisely we require that $ \om $ belongs  to the image 
\be\label{def:admissible}
 {\cal A}:=  \bar \mu + \e^2 
\Big\{ \Ab \xi \, : \, \frac12 \leq \xi_j \leq 4 \, , \forall j \in 
{\mathbb S}  \Big\}  \subset \R^\es \, 
\ee
of the compact set of actions $ \xi \in [1/2,4]^\es $ 
under the  approximate action-to-frequency map \eqref{act-to-fre}.   
Then, in view of the method that we shall use for the measure estimates for the linearized operator,  
we look for quasi-periodic solutions with frequency vector
\be\label{frequencycolinear}
\om = (1+\e^2 \lambda) \bar \om_\e \, , \quad \lambda \in \Lambda := [- \l_0, \l_0 ] \, ,  
\ee
constrained to a fixed admissible direction 
\be\label{def omep}
\bar \om_\e :=  \bar \mu + \e^2 \dom \, , \quad  \dom \in   \Ab ([1,2]^\es ) \, . 
\ee
Note that in general we can not choose $ \bar \om_\e = \bar \mu $, because 
$ \dom = 0 $ might not belong to $ \Ab ([1,2]^\es ) $. 
We fix $ \zeta $ below so that the Diophantine conditions \eqref{dioep}-\eqref{NRgt1} hold. 

In \eqref{frequencycolinear} there exists $ \l_0  > 0 $ small, 
independent of $ \e > 0 $ and of $ \dom \in   \Ab ([1,2]^\es )  $, such that, 
\be\label{Lambda-unif}
\forall  \l \in \Lambda := [- \l_0, \l_0] \  \, , \quad  \om = (1+\e^2 \l) \bar \om_\e \in {\cal A} 
\ee
are still admissible (see \eqref{def:admissible}) and, using \eqref{def omep},
\be\label{xil}
\begin{aligned}
& \ \ (1+\e^2 \l) \bar \om_\e = \bar \mu + \e^2  \Ab (\xi ) \quad
\Longleftrightarrow \\ 
& \xi := \xi(\l) = (1 + \e^2 \l )  \Ab^{-1} \dom + \lambda  \Ab^{-1}  \bar \mu \, . 
\end{aligned}
\ee
We shall use the $ 1 $-dimensional  ``parameter" 
$ \l \in \Lambda   := [-\l_0, \l_0]$ 
in order to verify all the non-resonance 
conditions required for the frequency vector  $ \om $ 
in the proof of Theorem \ref{thm:main}.

\smallskip
 
For $ \e $ small fixed, we take the vector $ \dom   $ 
such that the direction $ \bar \om_\e $ in \eqref{def omep} still verifies 
Diophantine conditions like  \eqref{diop}, \eqref{NRgamma0} with the different exponents
\be\label{def:tau1}
 \g_1 := \g_0/2 \, , \quad  \tau_1 := 3 \tau_0  + \es (\es + 1 )  + 5  > \tau_0 \, , 
\ee 
 namely\index{Diophantine vector}
 \begin{align}\label{dioep}
&  | \bar \om_\e \cdot \ell | \geq \frac{\g_1}{ \langle \ell \rangle^{ \t_1}} \, , \quad \forall \ell \in \Z^\es \setminus \{ 0 \} \, , \\
& 
\label{NRgt1}
\Big| n + \sum_{1\leq i\leq j \leq \es}  p_{ij} {({\bar \om}_\e)}_i {({\bar \om}_\e)}_j \Big| \geq  \frac{\g_1}{ \langle p \rangle^{ \t_1}} \, , \quad 
\forall (n, p) \in \Z \times \Z^{\frac{\es(\es+1)}{2}} \setminus \{(0,0)\} \, . 
\end{align}
This is possible by Lemma \ref{lemma:rho-dioph}. Actually 
the  vector $  \bar \om_\e = \bar {\mu} + \e^2 \dom $
satisfies \eqref{dioep}-\eqref{NRgt1} for all $ \dom \in   \Ab ([1,2]^\es ) $ except a small set  of measure $ O(\e) $. 
In \eqref{NRgt1},  
we denote, for $ i = 1, \ldots, \es $,  the $ i$-component 
$ ({\bar \om}_\e)_i = \mu_{j_i} + \e^2 \zeta_i $, where
$ j_1, \ldots, j_\es $ are the tangential sites ordered according to \eqref{unp-tangential}. 
\\[2mm]
{\bf Main result. } We may now  state in detail the main result of this monograph, 
concerning existence of quasi-periodic solutions of the nonlinear wave equation 
\eqref{NLW1}. 

Let $ {\mathcal H}^s  $ denote the standard Sobolev space 
$ {\mathcal H}^s ( \T^\es \times \T^d;  {\R}) $ of  functions 
 $ u : \T^\es \times \T^d \to \R $,  
 with norm $ \| \ \|_s $, see \eqref{Sobo:sp1}.

\begin{theorem}\label{thm:main}\index{Main result} {\bf (Quasi-periodic solutions for the
 nonlinear wave equation \eqref{NLW1})}. 
 Let $ a $ be a function in $ C^\infty (\T^d, \R ) $. Take 
$ V  $ in $ C^\infty (\T^d, \R ) $ such that 
$ - \Delta + V(x) \geq \beta {\rm Id} $ for some $ \b > 0 $
(see 	\eqref{positive}) and let 
$ \{  \Psi_j \}_{j \in \N} $  be
an $ L^2 $-basis of eigenfunctions  of 
$ - \Delta + V(x) $, with eigenvalues 
$ \mu_j^2 $ written in increasing order and with multiplicity as in \eqref{eige-order}. 
 Fix finitely many tangential sites $ {\mathbb S} \subset \N $. 

Assume the following 
conditions: 
\begin{description}
\item (i)
the unperturbed frequency vector 
$  \bar \mu = (\mu_j)_{j \in {\mathbb S}}\in \R^\es $ in \eqref{unp-tangential} satisfies the Diophantine conditions \eqref{diop} and \eqref{NRgamma0}; 
\item 
(ii) the unperturbed first and second order Melnikov non-resonance conditions 
\eqref{1Mel} and  \eqref{2Mel+}-\eqref{2Mel rafforzate} 
hold; 
\item 
(iii)  the Birkhoff matrix $ \Ab $ defined in \eqref{def:AB} satisfies the 
twist condition \eqref{A twist}; 
\item 
(iv) the finitely many non-degeneracy conditions \eqref{non-reso}-\eqref{non-reso1} hold.
\end{description}
Finally, fix  a vector 
$$ 
\bar \om_\e := \bar \mu + \e^2 \dom \, , \quad  \dom \in   \Ab ([1,2]^\es ) \subset \R^\es \, , 
$$
such that the Diophantine conditions \eqref{dioep}-\eqref{NRgt1} hold.  

Then the following holds: 
\begin{enumerate}
\item 
there exists a Cantor like set $ {\cal G}_{\e,\dom} \subset \Lambda $ 
(with $ \Lambda $ defined as in \eqref{Lambda-unif})
with asymptotically full measure,
i.e. 
$$ 
| \Lambda \setminus {\cal G}_{\e,\dom} | \to 0 \quad \hbox{as} \quad  \e \to 0 \, ; 
$$
\item 
for all $ \l \in {\cal G}_{\e,\dom} $,  there exists
a solution 
$  u_{\e, \lambda}  \in C^\infty  (\T^\es \times \T^d;  {\R}) $  of the equation \eqref{NLW2QP}, 
of the form
\be\label{sol:uel}
u_{\e, \lambda} ( \vphi,x ) = \sum_{j \in {\mathbb S}}   
\mu_j^{\, -1/2} \sqrt{2 \xi_j}  \cos ( \vphi_j ) \Psi_j (x) + r_\e (\vphi, x)  \, ,
\ee
which is  even in $ \vphi $,
where
$$ 
\om = (1+ \e^2 \l) \bar \om_\e \, , \quad \bar \om_\e = \bar \mu + \e^2 \dom \, , 
$$ 
and  $ \xi := \xi(\l) \in [1/2,4]^{\mathbb S} $ is given in \eqref{xil}.
For any $ s \geq s_0 > (\es + d) / 2 $,  the remainder term $r_\e $ satisfies  
$ \| r_\e \|_{s} \to 0  $  as $\e \to 0$.
\end{enumerate}
As a consequence $ \e u_{\e, \lambda} ( \om t , x ) $ is a quasi-periodic solution\index{Quasi-periodic solution}  of 
the nonlinear wave equation \eqref{NLW1} with frequency vector
$ \omega = (1+ \e^2 \l) \bar \om_\e $.
\end{theorem}

Theorem \ref{thm:main} is a direct consequence of Theorem \ref{thm:NM}.
Note that the non-resonance  conditions ($i$)-($ii$) in Theorem \ref{thm:main} 
depend on the tangential sites $ {\mathbb S} $ and the potential $ V(x) $,  
whereas the non-degeneracy  conditions ($iii$)-($iv$) depend on $ {\mathbb S} , V(x) $,
the choice of the basis $ \{ \Psi_j \}_{j \in \N}$ of eigenvectors of $ - \Delta + V(x) $,  and 
the coefficient $ a(x) $ in the nonlinear term $ g(x, u) = a(x) u^3 + O(u^4) $ in \eqref{nonlinearity}. 

In Section \ref{sec:ideas} we shall provide a detailed presentation of the strategy of proof 
of Theorem \ref{thm:main}. 
Let us now make some comments on this result.

\begin{enumerate}
\item {\bf Measure estimate of $ {\cal G}_{\e,\dom} $.} The speed of convergence of 
$| \Lambda \setminus {\cal G}_{\e,\dom} |$ to $0$ does not depend
on $\dom$. More precisely ($\gamma_0, \gamma_1, \tau_0, \tau_1$ being fixed) there is a map
$\e \mapsto b(\e)$, satisfying $ \lim_{\e \to 0} b(\e)=0 $,  
such that, for  all $\dom \in \Ab ([1,2]^\es )$ with 
$ \bar \om_\e = \bar \mu + \e^2 \dom $ satisfying the Diophantine conditions 
\eqref{dioep}-\eqref{NRgt1},  it follows that 
$ | \Lambda \setminus {\cal G}_{\e,\dom} | \leq b(\e) $.
\item {\bf Density.}
Integrating in $ \l $ along all possible 
admissible directions $ \bar \om_\e $ in \eqref{def omep}, we deduce the
existence of  quasi-periodic solutions of \eqref{NLW1}
for a set of frequency vectors $ \omega $ of positive measure.  
More precisely, defining the convex subsets of $\R^\es $,  
\be\label{defC1C2}
\begin{aligned}
& {\cal C}_2 := \bar \mu + \R_+ {\cal C}_1 \, ,  \\
& {\cal C}_1 := \Ab ([1,2]^\es ) + \Lambda \bar \mu :=  
\big\{ \dom + \l \bar \mu  \   ; \  \dom \in \Ab ([1,2]^\es ) \, , \, \l \in \Lambda \big\} \, , 
\end{aligned}
\ee
the set $\Omega$ of the frequency vectors $\om$ of the quasi-periodic solutions
of \eqref{NLW1} provided by Theorem \ref{thm:main} has Lebesgue density $1$ at $\bar \mu$ in ${\cal C}_2$, {\it i.e.}
\be\label{density-Leb}
\dps \lim_{r \to 0+} \frac{|\Omega \cap {\cal C}_2 \cap B(\bar \mu,r)|}{|{\cal C}_2 \cap B(\bar \mu,r)|} =1 \, ,
\ee
where $B(\bar \mu,r) $ denotes the ball in $ \R^\es $ of radius $ r $ centered at $ \bar \mu $
(see the proof below Theorem \ref{thm:NM}). 
Moreover, we restrict ourself to $\dom \in \Ab ([1,2]^\es )$ just to fix ideas, and 
we could replace this condition by
$\dom \in \Ab ([r,R]^\es )$, for any $0<r<R$ (at the cost of stronger smallness conditions for $\l_0$ and $\e$ if $r$ is small 
and $R$ is large). Therefore we could obtain a similar density result with ${\cal C}'_2:= \bar \mu +  \Ab (\R_+^\es )$
instead of ${\cal C}_2$.  
\item
{\bf Lipschitz dependence.} The solution $ u_{\e, \lambda}  $ is a Lipschitz function of 
$ \lambda \in {\cal G}_{\e, \zeta} $ with values in  $ {\cal H}^s $, for any 
$ s \geq s_0 $. 
\item {\bf Regularity.}
Theorem \ref{thm:main} also  holds  if 
the nonlinearity $ g (x, u) $   and the potential $ V(x) $ in \eqref{NLW1} are of class 
$ C^q $ for  $ q $ large enough, proving the existence of a solution 
$ u_{\e, \lambda} $ in the Sobolev space 
$ {\mathcal H}^{\bar s } $ for some finite $ \bar s > s_0 > (\es + d) \slash 2 $, see Remark \ref{rem:Cq}.  
\end{enumerate}

Theorem \ref{Prop:genericity} below proves that, for any choice of finitely many tangential sites 
$ {\mathbb S } \subset \N $, 
all the non-resonance and non-degeneracy assumptions  ($i$)-($iv$)
required in Theorem \ref{thm:main} are 
 verified for generic potentials $ V(x) $ 
and coefficients $ a(x) $  in the nonlinear term $ g(x, u) = a(x) u^3 + O(u^4) $ in \eqref{nonlinearity}.  
In order to state a precise result we introduce the following definition. 

\begin{definition}\label{Def:dense-open} {\bf ($C^\infty$-dense open)}
Given an open subset ${\cal U}$ of  $H^s (\T^d)$ (resp. $H^s (\T^d) \times H^s (\T^d)$) a subset ${\cal V} $ 
of ${\cal U}$  is said to be $C^\infty$-dense open in ${\cal U}$ if 
\begin{enumerate}
\item
${\cal V} $ is open for the
topology defined by the $H^s (\T^d)$-norm,  
\item 
${\cal V} $ is $C^\infty$-dense in $ {\cal U} $, in the sense that, for any $w \in {\cal U}$,
there is a sequence $(h_n) \in C^\infty(\T^d)$ (resp.  $C^\infty(\T^d) \times C^\infty(\T^d)$) such that 
$ w+h_n \in {\cal V}$,  for all $n \in \N$, and $ h_n \to 0 $ in 
$H^r $ for any $ r \geq 0 $.
\end{enumerate}
\end{definition}

Let $ s > d /2 $ and 
define the subset of potentials 
\be\label{pos:poten}
{\cal P} := \big\{ V \in H^s (\T^d) \ : \  - \Delta + V(x) > 0  \big\} \, .
\ee
The set ${\cal P} $ is open in $H^s(\T^d)$ and convex, thus connected. 

Given a finite subset ${\mathbb S} \subset \N $ of tangential sites,  
consider the set $ \wtilde{\cal G} $ of potentials  $ V(x) $ and coefficients  $ a(x) $ 
such that  the non-resonance and non-degeneracy 
conditions ($i$)-($iv$) required in Theorem \ref{thm:main} hold, namely 
\be\label{def:tildeG}
\begin{aligned}
\wtilde{\cal G} & := \Big\{ (V,a) \in {\cal P} \times H^s (\T^d)  \,  :  \,  
(i){\rm-}(iv) \ {\rm in  \ Theorem \ \ref{thm:main} \ }  \hbox{hold} \Big\} \, .
 \end{aligned} 
\ee
Note that the non-degeneracy properties $(iii)$-$(iv)$ may depend on the choice \eqref{auto-funzioni}
of the basis $ \{ \Psi_j \}_{j \in \N}$ of eigenfunctions of $ - \Delta + V(x) $ (if some eigenvalues are not simple).
In the above definition of $\wtilde{\cal G}$, it is understood that properties $(i)$-$(iv)$ hold 
for some choice of the basis $ \{ \Psi_j \}_{j \in \N}$.


Given  a subspace $ E $ of $ L^2(\T^d) $ we denote by $ E^{\bot}  := E^{\bot_{L^2}} $
its orthogonal complement with respect to the $ L^2 $ scalar product. 
 
\begin{theorem}\label{Prop:genericity}\index{Genericity result} {\bf (Genericity)}
Let $ s > d /2 $. 
The set
\be\label{C-infty-dense}
\wtilde{\cal G} \cap \big(C^\infty (\T^d) \times C^\infty (\T^d)\big) \quad
{ is}  \quad 
C^\infty{\rm -dense} \ in \  
\big({\cal P} \cap C^\infty (\T^d) \big) \times C^\infty (\T^d)  
\ee  
where  $ {\cal P} \subset H^s(\T^d)  $ is the open and connected set of potentials $ V(x)  $ defined in \eqref{pos:poten}.
 
More precisely, there is a  \dcopen subset 
$ \genset $ of $  {\cal P} \times H^s(\T^d)  $ and 
a  $|\mathbb S|$-dimensional linear
subspace $E$ of $C^\infty (\T^d)$ such that,  
 for all $ v_2 (x)  \in E^{\bot}  \cap H^s(\T^d) $, $ a (x) \in H^s (\T^d)  $, 
the Lebesgue measure (on the finite dimensional space $ E \simeq \R^{|\mathbb S|}$)
\be\label{Measure:gene}
\big| 
\big\{  v_1 \in E  \ : \  ( v_1+ v_2, a) 
\in {\cal G} \setminus \wtilde{\cal G} \big\} 
\big| = 0 \, .  
\ee
\end{theorem} 

Theorem \ref{Prop:genericity} is proved in Chapter \ref{section:gener1}.
 
 \smallskip

For the convenience of the reader,
we provide in the next chapter 
a non technical survey about  the main methods and results  in 
KAM theory for PDEs.


\section{Basic notation} 

We collect here some basic notation used throughout the monograph.

For $k=(k_1, \ldots , k_q) \in \Z^q$, we set 
$$
|k|:=\max \{ |k_1|, \ldots , |k_q| \} \, ,
$$
and, for $(\ell,j) \in \Z^\es \times \Z^d$,
$$
\langle \ell,j \rangle := \max(| \ell |,|j|,1) \, . 
$$
For $ s \in \R $ we denote the Sobolev spaces\index{Sobolev spaces}  
\be
\begin{aligned}\label{Sobo:sp1}
{\mathcal H}^s  & :=  {\mathcal H}^s ( \T^\es \times \T^d;  {\C}^r ) \\
& :=   \Big\{  u(\vphi,x)= \sum_{ (\ell,j) \in \Z^\es \times \Z^d}  u_{\ell,j} 
e^{\ii (\ell \cdot \vphi + j \cdot x)}  :  \, \| u \|_s^2 :=  \sum_{(\ell,j) \in \Z^{\es + d }}  
|u_{\ell,j}|^2 \langle \ell,j  \rangle^{2s}< \infty  \, \Big\} 
\end{aligned}
\ee
and we use the same notation $ {\mathcal H}^s  $ also for the subspace of  real valued functions. 
Moreover we denote by  $ H^s := H^s_x  $ the Sobolev space of functions $ u(x) $ in $ H^s (\T^d, \C )$ and   $ H^s_\vphi  $ the Sobolev space of functions $ u(\vphi ) $ 
in $ H^s (\T^\es, \C )$. 
We denote by $ b := \es + d $. 
For
\be\label{algebra-Sobolev}
s \geq s_0 > (\es + d) \slash 2 \, , 
\ee
we have the continuous embedding 
$$ 
{\mathcal H}^s (\T^\es \times \T^d ) \hookrightarrow  C^0 (\T^\es \times \T^d ) 
$$ 
and each ${\mathcal H}^s  $ is an algebra with respect to the product of functions.

\smallskip

Let $ E $ be a Banach space. Given a continuous  map $ u : \T^\es \to E $, $ \vphi \mapsto u(\vphi) $, 
we  denote by $\widehat{u}(\ell ) \in E $, $\ell \in \Z^\es$, 
its Fourier coefficients
$$
\widehat{u}(\ell ) := \frac{1}{(2 \pi)^{\es}} \int_{\T^\es} u(\vphi) e^{- \ii \ell \cdot \vphi } \, d \vphi 
\, , 
$$ 
and its average
$$
\langle u \rangle :=  \widehat{u}(0)  
:= \frac{1}{(2 \pi)^\es} \int_{\T^\es} u (\vphi) \, d \vphi \, . 
$$
Given an irrational vector $ \om \in \R^\es $, i.e. $ \om \cdot \ell \neq 0 $, $ \forall \ell \in \Z^\es \setminus \{0\} $, and 
a function 
$ g(\vphi) \in\R^\es $ with zero average, we define the solution $ h (\vphi) $ 
of $ \om \cdot \pa_\vphi h = g  $, with zero average,  
\be\label{op-inv-KAM}
h(\vphi) = (\om \cdot \pa_\vphi)^{-1} g := \sum_{\ell \in \Z^\es \setminus \{0\}}  \frac{\widehat{g}(\ell )}{\ii \om \cdot \ell} 
e^{\ii \ell \cdot \vphi  } \, . 
\ee
Let $ E $ be a Banach space 
with norm $ \| \ \|_E $. 
Given a function $ f : \Lambda := [-\l_0, \l_0] \subset \R  \to E $ 
we define its Lipschitz norm 
\be\label{def:Lip-norm}
\begin{aligned}
\| f \|_{\Lip} := \| f \|_{\Lip, \Lambda} := \| f \|_{\Lip,E} := 
\sup_{\l \in \Lambda} \| f \|_E +  |f|_\lip
\, , \\
|f|_\lip := |f|_{\lip, \Lambda} := |f|_{\lip,E} := \sup_{ 
 \l_1, \l_2 \in \Lambda,
 \l_1 \neq \l_2}  \frac{ \| f(\l_2) - f(\l_1) \|_E }{ |\l_2- \l_1|} \, .
\end{aligned}
\ee
If a function  $ f :  \wtilde \Lambda \subset \Lambda \to E $ is defined only 
on a subset $ \wtilde \Lambda $ of $ \Lambda $
we shall still denote by 
$ \| f \|_{\Lip} := \| f \|_{\Lip, \wtilde \Lambda} := \| f \|_{\Lip, E} $ the norm in \eqref{def:Lip-norm} where the sup-norm 
and the Lipschitz seminorm are 
intended in $\wtilde  \Lambda $, without  
specifying explicitly the domain $ \wtilde \Lambda $. 

If the Banach space $ E $ is the Sobolev space 
$ {\cal H}^s $ then  we denote more simply 
$ \| \ \|_{\Lip,{\cal H}^s} =
\| \ \|_{\Lip,s} $. If $ E = \R $ then $ \| \ \|_{\Lip,\R} =
\| \ \|_{\Lip} $. 

If $ A (\lambda) $ is a function, operator, \ldots, which depends on a
 parameter $ \lambda $,   we  shall use the following notation for  the partial quotient\index{Partial quotient}
\be\label{part-quo}
\frac{\Delta A}{\Delta \l} := \frac{A(\l_2) - A(\l_1)}{\l_2-\l_1} \, , \quad \forall  \l_1 \neq \l_2  \, . 
\ee
Given a family of functions, or linear self-adjoint operators  $ A(\lambda)  $ on a Hilbert space $ H $, 
defined for  all $ \lambda \in \wtilde \Lambda $, we shall use the  notation 
\be\label{partial-increase}
{\mathfrak d}_\l A(\l) \geq \b  {\rm Id} 
\qquad \Longleftrightarrow  \qquad \frac{\Delta A}{\Delta \l } \geq \b {\rm Id} \, , 
\quad \forall \l_1, \l_2 \in \wtilde \Lambda \, , \ \l_1 \neq  \l_2  \, , 
\ee
where, for a self-adjoint operator, 
$$ 
A \geq \b{\rm Id}  
$$ 
means as usual
$$ 
(Aw,w)_{H} \geq \b \| w \|_{H}^2 \, , \quad \forall  w \in H \, .
$$
Given linear operators $A, B  $ we denote their commutator\index{Commutator} by 
\be\label{Ad-op}
{\rm Ad}_A B := [A,B] := AB - B A \, . 
\ee
We define  $ D_V := \sqrt{- \Delta + V(x) } $ and $ D_m := 
\sqrt{- \Delta + m } $ for some $ m >  0 $. 
\begin{itemize}
\item 
Given $ x \in \R $ we denote by $\lceil x \rceil $  the smallest integer greater or equal to $ x $, and by 
$ [x] $ the integer part\index{Integer part} of $ x $, i.e. the greatest integer smaller or equal to $ x $;
\item $\N = \{0, 1,  \ldots \} $ denote the natural numbers, and $\N_+ = \{ 1,  \ldots \} $ the positive integers;  
\item 
Given $ L \in \N $, we denote by $ [\![0,L ]\!] $ the integers in the interval $[0,L] $;
\item 
We use the notation $ a \lesssim_s b  $ to mean $ a \leq C(s) b $ for some positive 
constant $ C(s) $, and 
$ a \sim_s b $ means that 
$ C_1(s) b \leq a \leq C_2(s) b $ for
positive constants $ C_1(s), C_2(s)  $;
\item
Given functions $ a, b : (0,\e_0) \to \R $ we write 
\be\label{def:allb} 
a (\e)  \ll b (\e)  \qquad  \Longleftrightarrow \qquad \lim_{\e \to 0} \frac{a(\e)}{b(\e )} = 0  \, . 
\ee
\end{itemize}
In the Monograph we denote by ${\mathbb S} $, $ {\mathbb F} $, 
$ {\mathbb G}   $, $ {\mathbb M}   $ subsets of the natural numbers $ \N $, with
$$
\N = {\mathbb S} \cup {\mathbb S}^c \, , \quad {\mathbb F} \cup {\mathbb G} = 
{\mathbb S}^c \, , \quad  {\mathbb F} \cap {\mathbb G} = \emptyset \, ,  
\quad {\mathbb F} \subset {\mathbb M } \, . 
$$
We refer to Chapter \ref{Ch:3} for the detailed notation of operators, matrices, 
decay norms, \ldots 

\smallskip

For simplicity of notation we may write either $ \xi \in \R^{\mathbb S} $
or $ \xi \in \R^{\es} $. 

Along the monograph we shall use the letter $ j $ to denote a space-index: it may be 
\begin{itemize}
\item  $ j $ in $ \N $, when we use the $ L^2 $-basis 
$ \{  \Psi_j \}_{j \in \N} $ defined in \eqref{auto-funzioni}; 
\item 
$ j = (j_1, \ldots, j_d ) $ in  $ \Z^d $, when we use the  exponential 
basis $ \{ e^{ \ii j \cdot x} \}_{j \in \Z^d } $. 
\end{itemize}

\chapter{KAM for PDEs and strategy of proof} \label{Ch1.5}


In this chapter we first describe the basic ideas of a Newton-Nash-Moser 
algorithm to prove existence of quasi-periodic solutions. 
The key step consists in the analysis of the 
linearized operators obtained at each step of the iteration, prove that they are 
approximately invertible, for most values of suitable parameters, with quantitative 
tame estimates for the inverse in high 
norms.  
Then we describe the  two main approaches that have been developed
for achieving this task for PDEs:  
\begin{enumerate}
\item
the ``reducibility" approach, that we describe in Section   \ref{sec:red}, with its main applications to KAM theory (Section \ref{sec:RR}); 
\item the ``multiscale" approach, presented in Section \ref{sec:MULTI}. 
\end{enumerate}
Finally, in Section \ref{sec:ideas}, we shall provide a detailed account
of the proof of   Theorem \ref{thm:main}.

\section{The Newton-Nash-Moser algorithm}\label{sec:NNM}

The classical implicit function theorem \index{Implicit function Theorem} is concerned
with the solvability of the equation
\be\label{Ffh}
 {\cal F} (x,y) = 0
\ee
where $ {\cal F}  : X \times Y \to Z $ is a smooth map, 
$ X $, $ Y $, $ Z $ are Banach spaces, 
and there exists
$(x_0, y_0) \in X \times Y $  
such that 
$$
 {\cal F}   ( x_0, y_0 ) = 0 \, .
$$
If $ x $ is close to $x_0 $ we want to solve (\ref{Ffh})
finding $ y= y(x)$.

The main assumption of the classical implicit function theorem 
is that the partial derivative $ (d_y  {\cal F} ) (x_0, y_0 ) : Y \to Z $
possesses a {\it bounded} inverse. Note that, in this case, by a simple perturbative argument, 
$ (d_y  {\cal F} ) (x, y )$ possesses a bounded inverse for any $(x,y)$ in some neighborhood 
of $(x_0,y_0)$.

However
there are many situations where 
$$ 
(d_y  {\cal F} ) (x_0, y_0 ) \ {\rm has \ an \ {\it unbounded} \  inverse} \, , 
$$
for example due to small divisors. 

\smallskip

An approach to these class of problems has been proposed by Nash \cite{Nash}
for the isometric embedding problem of a Riemannian manifold.
Subsequently, Moser \cite{Mo61} has  
proposed an abstract modified algorithm in scales of Banach spaces,
with new applications in small divisor problems 
\cite{Mo62}, \cite{Mo1}, \cite{M67}. 
Further extensions  and improvements have been obtained  
by Zehnder \cite{Z1}-\cite{Z2}, see also \cite{N},  
and 
H\"ormander \cite{Ho1}, \cite{Ho2}, see also 
Hamilton \cite{Ha} and the more recent papers \cite{BBP10},  \cite{BH}, \cite{ES}.

\smallskip

There are many variants to present the Nash-Moser theorems, 
according to the applications one has in mind. 
The main idea in order to find a 
solution  of the nonlinear equation $ {\cal F} ( x, y) = 0 $
is to replace the usual Picard iteration method with a modified
Newton iterative scheme. 
One starts with a sufficiently good  approximate solution 
$  (x_0, y_0)  $, i.e. 
$  {\cal F}(x_0, y_0) $ is small enough. 
Then,  
given an approximate solution $ y $, 
we look for a better approximate solution 
of the equation $ {\cal F} (x,y) = 0 $, 
$$ 
y' = y + h  \, , 
$$  
by a  Taylor  expansion (for simplicity we omit to write the dependence on $ x $) 
\be\label{Newton-s0}
{\cal F} ( y') = {\cal F} ( y + h )  =   {\cal F} (  y) +  d_y {\cal F} ( y)[h] + O( \| h \|^2 ) \, . 
\ee
The idea of the classical Newton iterative scheme is to define $ h $ as the solution of 
\be\label{Newton-s1}
{\cal F} (  y) +  d_y {\cal F} ( y)[h] = 0 \, . 
\ee
The main difficulty of this scheme is to prove that $d_y {\cal F} ( y) $
is invertible, for any $ y $ in a neighborhood of the unperturbed solution $ y_0 $, 
in order to define the new correction
\be\label{Newton-s}
h = -  (d_y {\cal F} ( y))^{-1} {\cal F} (  y)  \, . 
\ee
The invertibility (in a weak sense) of $d_y {\cal F} ( y) $ for $y$ close enough to $y_0$ is not a
simple consequence of the existence of an unbounded inverse of $d_y {\cal F} ( y_0) $.
In addition, working in scales of Banach spaces, one has to  prove that the inverse $  (d_y {\cal F} ( y))^{-1} $ 
satisfies suitable bounds in high norms, with a ``loss
of derivatives'', of course. The ``loss of derivatives'' roughly means that to estimate 
$  (d_y {\cal F} ( y))^{-1} [h] $ in some norm, one must use a higher norm for $h$ (and possibly for $y$).
On the other hand, the main advantage of 
the scheme \eqref{Newton-s0}-\eqref{Newton-s} is to be quadratic, namely
we have   
$$
{\cal F} ( y') =   O( \| h \|^2 ) = 
O( \|{\cal F} (  y) \|^2 ) \, .
$$
As a consequence, 
the iterates will converge
to the expected solution at a super-exponential rate,
compensating the divergences in the scheme due to 
the ``loss of derivatives''. 
For the analytical aspects of the convergence in 
scales of Banach spaces of analytic functions
we refer for example to 
Zehnder 
\cite{N},  
Section 6.1.

A Newton scheme of this kind
was also used  by Kolmogorov \cite{K} 
 and Arnold \cite{Ar}  for proving their celebrated persistence result 
 of quasi-periodic solutions in nearly integrable 
 analytic Hamiltonian systems.  

When one works   
in scales of Banach spaces of functions with finite differentiability,
like $ C^k $ or Sobolev spaces $ H^s $,  
the Newton scheme  \eqref{Newton-s} will not converge because 
the inverse $ (d_y {\cal F} (y))^{-1}  $ is unbounded, 
say  $ (d_y {\cal F} (y))^{-1}  : H^s \mapsto H^{s- \tau} $ for some $ \tau > 0 $
(it ``loses" $ \tau  $ derivatives) and therefore, 
after finitely many steps,  the approximate solutions are no longer regular.  

Then 
Moser \cite{Mo61} suggested 
 to perform 
 a smoothing 
 of the approximate solutions
obtained at each step of the Newton scheme, like a Fourier series truncation 
with an increasing number of  
harmonics, see the smoothing operators \eqref{def:projectorN}.
The convergence of the iterative scheme is then guaranteed 
requiring that  the inverse  $ (d_y {\cal F} (y))^{-1} $ 
satisfies  ``tame" estimates. 
We  state them precisely 
in the case that $ {\cal F} $ acts on a scale  
$ (H^s)_{s \geq s_0} $ of Sobolev spaces
with norms $ \| \ \|_s $:
there are constants 
$ p $, $ \rho >  0 $ (``loss of derivatives") such that, for any
$ s \in [s_0, S] $, for all $  g \in H^{s+\rho}$, 
 \be\label{AINV0SOB}
\| (d_y {\cal F} (y))^{-1} g \|_{s} \leq C(s, \|y\|_{s_0+p}) \big( \| g \|_{s +\rho} + \| g \|_{s_0} \| y \|_{s+ \rho}  \big)  \, . 
\ee
These tame estimates
are sufficient for the convergence 
of the Nash-Moser iterative scheme, 
requiring that the initial datum $ \| {\cal F} (x_0,y_0)\|_{s_1}  $
is small enough  for some
$ s_1 := s_1 (p, \rho ) > s_0 $. 
For the analytical aspects of the convergence 
we refer for example to 
Moser \cite{Mo1}, Zehnder 
\cite{N},  Chapter 6, and \cite{BB07}, \cite{BBP10}.

\begin{remark}
{\bf (Parameters)}
If the nonlinear operator $ {\cal F}(\l; x,y) $ depends on some parameter $ \l $,  
tame estimates of the form  \eqref{AINV0SOB} may hold 
for ``most" values of $ \l $. Thus the solution $ y(x; \lambda) $ of the implicit equation
$ {\cal F}(\l; x,y) = 0 $ 
exists for most values of the parameters $\lambda $. 
This is 
the typical situation that one encounters in small divisor problems, as in 
the present monograph, see the nonlinear operator 
\eqref{operatorF}. 
\end{remark}

\begin{remark}\label{WTE}
{\bf (Weaker tame estimates)} 
Weaker tame estimates of the form \eqref{LN-1tame}
or \eqref{LN-1tame-s} are sufficient 
for the convergence of the
Nash-Moser scheme, see \cite{BBP10},  \cite{BB12}, \cite{BBo10}. 
In this monograph we shall indeed verify  tame estimates of this weaker type.
\end{remark}

Before concluding this section
we mention a  variant of the above Newton-Nash-Moser scheme 
 \eqref{Newton-s0}-\eqref{Newton-s} (that we shall employ). 
Zehnder \cite{Z1} noted that  
it is sufficient to find only an approximate right inverse of  $ d_y {\cal F} ( y) $, namely a linear operator
$ T(y)   $ such that
\be\label{InvApp}
d_y {\cal F} ( y) \circ T(y) - {\rm Id}  = O( \| {\cal F} (  y)\| ) \, .
\ee
Remark that, at a solution $ y $ of the equation $ {\cal F} (  y) = 0 $, the operator $ T(y) $ is an exact right inverse of $ d_y {\cal F} ( y) $. 
Then we define the new approximate solution 
$$
y' = y + h \, , \quad h := - T(y) {\cal F}( y)   \, , 
$$
obtaining,  by \eqref{InvApp}, that 
$$
{\cal F} ( y') =   {\cal F} (  y) -  d_y {\cal F} ( y)[ T(y) {\cal F} ( y)] + O( \| h \|^2 ) = 
O( \|{\cal F} (  y) \|^2 ) \, . 
$$
This is still  a quadratic scheme. 

Naturally, for ensuring the convergence of this iterative scheme in 
a scale of Banach spaces with finite differentiability,  
the approximate inverse $ T(y)$ has to satisfy
tame estimates as in \eqref{AINV0SOB}.

The above modification of the Newton-Nash-Moser  iterative scheme is useful because 
the invertibility of $ d_y {\cal F} ( y) $
may be a  difficult task, but the search of an approximate right inverse 
$ T(y) $ may be much simpler. In Chapter \ref{sezione almost approximate inverse} we shall indeed implement this device for the nonlinear operator 
$ \mF $ 
defined in \eqref{operatorF}. 

\smallskip

As we outlined in this section, the core of a Nash-Moser scheme concerns  
the  invertibility of  the linearized 
operators which arise at each step of the iteration with 
tame estimates like \eqref{AINV0SOB}, or  \eqref{LN-1tame}, \eqref{LN-1tame-s} for its inverse.
In the next sections we present the two major 
set of ideas and techniques which have been 
employed so far to prove the invertibility of  the linearized 
operators which arise in the search of quasi-periodic solutions
for PDEs. These are:  
\begin{enumerate}
\item the ``reducibility" approach, that we present in Section \ref{sec:red} with its main applications to KAM theory for PDEs (Section \ref{sec:RR}); 
\item  the ``multiscale" approach,  described in Section \ref{sec:MULTI}, which is 
at the basis of the proof of  the main result of this  monograph, Theorem \ref{thm:main}.  
\end{enumerate}

\section{The   reducibility approach to KAM for PDEs}\label{sec:red}

The goal of this section is to present the 
perturbative reducibility\index{Reducibility} approach for a time 
quasi-periodic linear operator (subsection \ref{sec:PR}) and then   
describe its main applications to 
KAM theory for PDEs (Section \ref{sec:RR}). 

\subsection{Transformation laws and reducibility} 

Consider a quasi-periodically time dependent  linear system
\be\label{sis1}
u_t + A(\om t ) u = 0  \, ,  \quad u \in H \, . 
\ee
Here the phase space $H$ may be a finite or infinite dimensional Hilbert space with scalar product
$ \langle \, , \,  \rangle $; $A$ (smoothly) maps any $ \vphi \in \T^\nu $,  $ \nu \in \N_+ $, 
to a linear operator $ A(\vphi ) $ acting on   
$ H $;
 $ \om \in \R^\nu \backslash \{ 0 \}  $ is the  frequency vector. We suppose that $ \omega $ is a nonresonant vector,
i.e. $ \omega \cdot \ell \neq 0 $, $ \forall \ell \in \Z^\nu \setminus \{ 0 \} $, thus any trajectory
$(\vphi_0 + \omega t)_{t \in \R}$ of the linear flow 
spanned by $\omega$ is dense in the torus $ \T^\nu $. 

Assume that $\Phi $ (smoothly) maps any $\varphi \in \T^\nu$ to an invertible bounded linear operator
$\Phi(\vphi) : H \to H$. 
Under  the quasi-periodically time dependent transformation
\be\label{QP-floquet}
u = \Phi (\om t ) [v]  \, , 
\ee
system \eqref{sis1} transforms into  
\be\label{sis2}
v_t + B(\om t) v = 0
\ee
with the new linear operator 
\be\label{newB}
B(\vphi ) =  \Phi^{-1} (\vphi) (\omega \cdot \pa_\vphi \Phi) (\vphi) +  \Phi (\vphi)^{-1} A (\vphi) \Phi  (\vphi) \, .
\ee
\begin{remark}\label{rem:Ham}
Suppose that $ H $ is endowed with a symplectic form  $ \Omega $ defined by 
$ \Omega (u,v) := \langle J^{-1} u, v \rangle  $, $ \forall u, v \in H $, 
where $ J $ is an antisymmetric, non-degenerate operator. 
If $ A (\om t) $ is Hamiltonian, namely $ A (\om t) = J S (\omega t) $ 
where $ S( \omega t ) $ is a (possibly unbounded) self-adjoint operator, 
 and $ \Phi (\omega t ) $  is symplectic, then the new operator 
$ B (\om t) $ is Hamiltonian as well, see Lemma \ref{lem:sym}.
\end{remark}

\noindent
{\bf Reducibility.}
If the operator $ B $ in \eqref{sis2} is a diagonal,  time independent operator, i.e. 
\be\label{B-reduced}
B (\om t ) = B = {\rm Diag}_j ( b_j )
\ee
in a suitable basis of $ H $, 
then
\eqref{sis2} reduces to the 
decoupled  
scalar linear ordinary differential equations 
\be\label{bj-diago}
\dot v_j + b_j v_j = 0 
\ee
where $ (v_j) $ denote the coordinates 
of $ v $  in the basis of eigenvectors of   $ B $.
Then  \eqref{sis2} is  integrated in a straightforward way, 
$$
v_j (t) = e^{- b_j t } v_j (0) \, , 
$$
and all the solutions
of system \eqref{sis1} are obtained via the  
change of variable  \eqref{QP-floquet}.
We say that
\eqref{sis1} has been {\it reduced} to constant coefficients by the change of variable 
$ \Phi $. 

\begin{remark}
If all the $ b_j $ in \eqref{bj-diago} are purely imaginary, then the linear  system
\eqref{sis1} is stable (in the sense of Lyapunov),
otherwise, it is unstable. 
\end{remark}

Let $GL(H)$ denote the group of invertible bounded linear operators of $H$. We shall  say  that  system \eqref{sis1} is reducible if 
there exists a (smooth) map  $ \Phi : \T^\nu \to GL(H) $ such that the operator
$ B(\vphi) $ in \eqref{newB}
is  a constant coefficient block-diagonal operator $ B $, i.e. 
$ b_j $ in \eqref{B-reduced}-\eqref{bj-diago} are finite dimensional 
matrices, constant in time. The spectrum of each matrix $ b_j $ determines 
the stability/instability 
properties of the system \eqref{sis1}.

\smallskip

If $ \omega \in \R  $ (time-periodic forcing) and the phase space 
$ H $ is  finite dimensional, 
the classical Floquet theory proves that any  time periodic linear system \eqref{sis1} is reducible, see
e.g. \cite{Ek}, Chapter I. 
On the other hand, if $ \omega \in \R^\nu  $, $ \nu \geq 2 $, 
it is known that there exist pathological 
non reducible linear systems, see e.g. \cite{BHS}-Chapter 1.

\smallskip

If  $ A(\om t) $ is a small perturbation of a constant coefficient operator, 
perturbative algorithms for reducibility can be implemented. In the next subsection
we describe this strategy in the simplest setting.
This approach was systematically adopted 
by Moser \cite{M67}  for developing 
finite dimensional KAM theory  (in a much more general context). 

\subsection{Perturbative reducibility}\label{sec:PR}

In the framework of the study of quasi-periodic solutions  of a PDE, \eqref{sis1}
stands for the linearized PDE at a quasi-periodic (approximate) solution of frequency vector 
$\omega$. 
In a Nash-Moser scheme, one has to invert the operator 
\be\label{inizA}
\om \cdot \pa_\vphi  + A(\vphi)  
\,  , \quad \vphi \in \T^\nu 
\, , 
\ee
which acts on maps $U : \T^\nu \to H$. We now assume  that
$A(\vphi)$ is a perturbation of a diagonal operator $D$ acting on $H$ and  not depending on $\vphi $, namely
\be\label{diagDhi}
A(\vphi) =  D  + R (\vphi) \, , \qquad   D = {\rm Diag}(\ii d_j)_{j \in \Z} = {\rm Op} ( \ii d_j) \, , \quad d_j \in \R \, , 
\ee
where the eigenvalues $ \ii d_j $ are simple. 
The  family of operators $ R(\vphi ) $ is acting on 
the phase space $ H $,  and plays the role of a small perturbation.  

\begin{remark}
We suppose that  $ d_j $, $ j \in \Z $, are real, i.e.
$ u = 0 $ is an elliptic equilibrium for the linear system
$ u_t + D u = 0 $. 
If some $ {{\rm Im} \, d}_j \neq 0 $,  then 
there are hyperbolic directions. These eigenvalues 
do not create resonance phenomena and 
perturbative reducibility theory is easier. For nonlinear 
systems, this case  corresponds to the search of whiskered tori, see e.g. 
\cite{Gr74} for finite dimensional systems, and 
\cite{FDLS} for PDEs. 
\end{remark}

We  look for a  transformation $ \Phi (\vphi)  $, $ \vphi \in \T^\nu $,  of the phase space $ H $, 
as in \eqref{QP-floquet} which removes from $ R (\vphi) $ 
the angles $ \vphi $ up to terms of size $ \sim O( |R|^2 ) $. 
We present below only the algebraic 
aspect of the reducibility scheme, without specifying the norms in the  
phase space $ H $ or the   operator norms of $ \Phi (\vphi) , R $, etc...

For computational purposes, 
it is convenient to  transform the linear system \eqref{sis1} under the flow 
$ \Phi_F (\vphi, \tau) $ generated by an auxiliary linear equation
\be\label{aux-F}
\pa_\tau \Phi_F (\vphi, \tau)   = F (\vphi) \Phi_F (\vphi, \tau)  \, , \quad \Phi_F (\vphi, 0) = {\rm Id } \, ,   
\ee
generated by a linear operator $ F (\vphi) $ to be chosen (which  could also be $ \tau $-dependent).
This amounts 
to computing the Lie derivative of $ A (\vphi)  $ in the direction of  the vector field $  F (\vphi) $. 
Note that, if $ F(\vphi) $ is bounded, then the flow \eqref{aux-F} is well posed
 and given by $ e^{\tau F(\vphi)} $. 
 This is always the case for a finite dimensional system, but it may be an issue for infinite dimensional systems,  if the generator $ F(\vphi )$ is not bounded.  

Given a linear operator $ A_0 (\vphi) $, the conjugated operator under the flow
$ \Phi_F (\vphi, \tau) $ generated by  \eqref{aux-F}, 
$$
A(\vphi, \tau) := \Phi_F (\vphi, \tau) A_0(\vphi) \Phi_F (\vphi, \tau)^{-1} \, , 
$$
 satisfies the Heisenberg equation
\be\label{eq:Hei}
\begin{cases}
\pa_\tau A(\vphi, \tau) 
= [ F (\vphi)  , A(\vphi, \tau)]  \cr
A(\vphi, \tau)_{|\tau=0} = A_0 (\vphi)  
\end{cases}
\ee
where $[A, B] := A \circ B - B \circ A $ denotes the commutator between two linear operators $ A, B $. 
Then, by a Taylor expansion, using \eqref{eq:Hei},  we obtain 
the formal Lie\index{Lie expansion} expansion 
\be\label{Lie-ex}
A(\vphi, \tau)_{|\tau = 1 } = \sum_{k \geq 0} \frac1{k!} {\rm Ad}^k(A_0)= A_0 (\vphi) + {\rm Ad}_F A_0 + \frac{1}{2} {\rm Ad}_F^2 A_0 + \ldots
\ee
where ${\rm Ad}_F [ \, \cdot \, ] := [F, \, \cdot \, ] $. 
One may expect  this expansion to be  certainly 
convergent if $ F  $ and $ A_0 $ are bounded 
in suitable norms. 

\smallskip

Conjugating \eqref{inizA} under the flow generated by \eqref{aux-F} we then obtain an operator of the form 
\be\label{trasop1}
{\om \cdot \pa_\vphi} + D  -  {\om \cdot \pa_\vphi} F (\vphi) +    [ F(\vphi), D  ] + R(\vphi)  
+ \ {\rm smaller \ terms} \ldots \, .
\ee
We want to choose $ F ( \vphi ) $ in such a way to solve the ``homological" equation 
\be\label{primopez}
-{\om \cdot \pa_\vphi} F  (\vphi) +  R(\vphi)  + [F(\vphi), D]  = [R] 
\ee
where 
\be\label{RNF} 
[ R ] := {\rm Diag}_{j} ( \widehat R^j_j (0)) \, , \quad  \widehat R^j_j (0) := \frac{1}{(2\pi)^\nu} \int_{\T^\nu} R_j^j (\vphi) \, d \vphi \, , 
\ee
is the normal form part of the operator  $ R (\vphi) $, independent of $ \vphi $,  that we can not eliminate.  
Representing the linear operators 
$ F (\vphi) =  (F^j_k (\vphi))_{j,k \in \Z} $ and $ R (\vphi) = (R^j_k (\vphi))_{j,k \in \Z} $ 
as matrices,
and computing the commutator with the diagonal operator $ D$  in  \eqref{diagDhi}
we obtain that \eqref{primopez} is represented as
$$
-{\om \cdot \pa_\vphi} F^j_k (\vphi)   +  R^j_k  (\vphi)  + \ii (d_j- d_k) F^j_k (\vphi)   = [R ]^j_k   \, .
$$
Performing the Fourier expansion in $ \vphi $, 
$$
F^j_k (\vphi)  = \sum_{\ell \in \Z^\nu}  \widehat F^j_k (\ell) e^{\ii \ell \cdot \vphi} \, , \quad
R^j_k (\vphi)  = \sum_{\ell \in \Z^\nu}  \widehat R^j_k (\ell) e^{\ii \ell \cdot \vphi} \, ,
$$
the latter equation reduces to the infinitely many scalar equations 
\be\label{eq:homo}
- \ii \om \cdot \ell \,  \widehat F^j_k (\ell)   +   \widehat R^j_k  (\ell)  
+ \ii (d_j- d_k)  \widehat F^j_k (\ell)  = [R ]^j_k \delta_{\ell, 0} \, , 
\quad j, k \in \Z \, , \quad \ell \in \Z^\nu \, ,
\ee
where $ \delta_{\ell,0} := 1 $ if $ \ell = 0 $ and zero otherwise. 
Assuming the 
so called {\it second-order Melnikov non-resonance conditions} 
\be\label{sec-Me}
|\om \cdot \ell + d_j - d_k | \geq \frac{\g}{\langle \ell \rangle^\tau} \, , \quad \forall (\ell, j, k) \neq (0, j, j ) \, ,
\ee
for some $ \gamma, \tau  > 0 $, 
we can define the solution of the homological equations \eqref{eq:homo} (see \eqref{RNF})
\be\label{defF}
 \widehat F^j_k (\ell) := 
\begin{cases}
\displaystyle \frac{ - \widehat R^j_k (\ell)}{ \ii (- \om \cdot \ell + d_j - d_k)} \ \,  \quad \forall (\ell, j, k) \neq (0, j, j )  \cr
0 \qquad  \qquad \qquad \qquad \ \ \,  \ \forall  (\ell, j, k) = (0, j, j ) \, .
\end{cases}
\ee
Therefore the transformed operator \eqref{trasop1} becomes 
\be\label{new+}
{\om \cdot \pa_\vphi} + D_+ + {\rm smaller \  terms}
\ee
where
\be\label{newD+}
D_+ := D + [R] = ( \ii d_j + [R]_j^j)_{j \in \Z} 
\ee
is the new diagonal operator, constant in $ \vphi $. 
We can iterate this step 
to reduce also the small terms of order $ O( |R|^2 \g^{-1}) $ which are left 
in \eqref{new+}, and so on. 
Note that, if $ R (\vphi) $ is a bounded operator and depends smoothly enough on $\vphi$, then $ F (\vphi) $ defined in 
\eqref{defF}, with the denominators satisfying \eqref{sec-Me}, is bounded as well and thus \eqref{aux-F} certainly defines a flow by standard Banach space ODE techniques.  
On the other hand the loss of time derivatives induced on $ F (\vphi)  $ by the divisors in \eqref{sec-Me} 
can be recovered by a smoothing procedure in the angles $ \vphi $, like a truncation in 
Fourier space. 

\begin{remark}\label{rem:Ham1}
If the operator $ A (\vphi) $ in \eqref{inizA}-\eqref{diagDhi} is Hamiltonian, as defined in Remark \ref{rem:Ham}
(with a symplectic form $ J $ which commutes with  $ D $), 
then $ F(\vphi ) $ is Hamiltonian, its flow $ \Phi_F ( \vphi, \tau ) $  
is symplectic, and the new operator in \eqref{new+} 
 is Hamiltonian as well.
\end{remark}

In order to continue the iteration one also needs to impose non-resonance conditions as in \eqref{sec-Me}
at each step and therefore we need information about 
the perturbed normal form $ D_+ $ in \eqref{newD+}, in particular the
asymptotic of $ [R]_j^j $. If these steps work, then, 
after an infinite iteration, one could conjugate the quasi-periodic operator \eqref{inizA} to a diagonal, constant in $ \vphi $, operator of the form  
\be \label{Dinf}
\om \cdot \pa_\vphi + {\rm Diag}_j ( \ii d_j^\infty) \, , \quad \ii d_j^\infty = \ii d_j + [R]_j^j + \ldots \, . 
\ee
At this stage, 
imposing the first order Melnikov non-resonance conditions
$$
| \om \cdot \ell + d_j^\infty | \geq \frac{\gamma}{\langle \ell \rangle^\tau} \, , \quad \forall \ell,  j \, ,  
$$
the diagonal linear operator \eqref{Dinf} is invertible 
with an inverse which loses $ \tau $ time-derivatives.
One then verifies that all the changes of coordinates
-that have been constructed iteratively to conjugate \eqref{inizA} to \eqref{Dinf}-
map spaces of high regularity in itselves. In conclusion 
this approach  enables to prove
 the existence of an inverse of the initial quasi-periodic linear operator 
 \eqref{inizA} which 
satisfies  tame estimates in high norms (with loss of $ \tau $ derivatives). 

\smallskip

This is the essence of  the Newton-Nash-Moser-KAM perturbative 
reducibility\index{Reducibility} scheme, that has been used for proving KAM results for 
$ 1 $d NLW and NLS equations with Dirichlet boundary conditions in 
\cite{Ku}, \cite{W1}, \cite{Po2}, as we shall describe 
 in the next section. 

\smallskip
The following questions arise naturally: 

\begin{enumerate}
\item \label{dif1}
{\it What happens if the eigenvalues $ d_j $ are multiple?} This is the common 
situation for $ 1$-$d$ PDEs with periodic boundary conditions or in higher space dimensions. In such a case
it is conceivable to reduce $ \om \cdot \pa_\vphi + A(\vphi) $ 
to a block-diagonal normal form linear system of the form \eqref{Dinf} where
$ d_\infty^j $ are finite dimensional matrices.  
\item \label{dif2}
{\it What happens if the operator $ R (\vphi)  $ in \eqref{diagDhi}
is unbounded?} This is the common situation for PDEs with 
nonlinearities which contain derivatives. In such a case
also the operator $ F (\vphi) $  defined in \eqref{defF} is unbounded and therefore \eqref{aux-F} might not define a flow.
\item \label{dif3}
{\it What happens if, instead of the Melnikov non-resonance  conditions \eqref{sec-Me}, 
we assume only 
\be\label{sec-Me-loss}
|\om \cdot \ell + d_j - d_k | \geq \frac{\g}{\langle \ell \rangle^\tau 
\langle j \rangle^{{\mathtt d}} \langle k \rangle^{{\mathtt d}} } \, , \quad \forall (\ell, j, k) \neq (0, j, j ) \, , 
\ee
for some  $ {\mathtt d} > 0  $, which induce a  loss of  space derivatives?} This is 
the common situation when the dispersion relation $ d_j \sim j^\alpha $, $ \a < 1 $, has a sublinear growth. 
Also in this case the operator $ F (\vphi) $  defined in \eqref{defF} would be unbounded. 
This situation appears for example for pure gravity water waves equations. 
\end{enumerate}

We describe  below some answers to the above questions. 

\begin{remark}
Another natural question is what happens if the number of frequencies 
is $ + \infty $, namely for {\it almost-periodic solutions}. 
Only a few results are available so far, all for PDEs in $1d$. 
The first result  is  due to P\"oschel \cite{Po4} for a regularizing  Schr\"odinger equation
and almost-periodic solutions
with a very fast decay in Fourier space.  
Then 
Bourgain \cite{B6} proved existence of almost-periodic 
solutions with  exponential Fourier decay
for a  semilinear Schr\"odinger equation with a convolution potential.  
In both papers the potential is regarded as an ``infinite dimensional" parameter. 
\end{remark}

\section{Reducibility results} \label{sec:RR}

We present in this separate section the main results about KAM theory for PDEs 
based on the reducibility\index{Reducibility} scheme described in the previous section.

\subsection{KAM for $1 d$ NLW and NLS with Dirichlet boundary conditions}

The iterative reducibility scheme outlined in
Section \ref{sec:PR} 
 has been effectively implemented 
by Kuksin \cite{Ku} and Wayne \cite{W1} for proving existence of quasi-periodic solutions\index{Quasi-periodic solution}  of 
1-$d$ semilinear wave\index{Nonlinear wave equation} 
\be\label{NLW1-L}
y_{tt} - y_{xx} + V(x) y + \e f (x, y) = 0 \, , \quad y(0) = y(\pi) = 0  \, , 
\ee
and Schr\"odinger equations\index{Nonlinear Schr\"odinger equation}
\be\label{NLS1}
\ii u_t - u_{xx} + V(x) u + \e f(|u|^2) u = 0 \, , \quad u (0) = u(\pi) = 0   \, , 
\ee
with Dirichlet boundary conditions.  
These equations are regarded as a perturbation of the linear PDEs\index{Linear Schr\"odinger equation} \index{Linear wave equation} 
\be\label{lin:PDE}
y_{tt} - y_{xx} + V(x) y   = 0 \, , \quad 
\ii u_t - u_{xx} + V(x) u  = 0  \, , 
\ee
which depend on the potential $ V (x) $, used as a parameter. 

The linearized operators obtained at an approximate quasi-periodic solution are, 
for NLS,  
\be\label{lin:NLS}
h \mapsto 
\ii \om \cdot \pa_\vphi h - h_{xx} + V(x) h + \e q(\vphi, x) h  + \e p(\vphi, x) \bar h  
\ee
with $ q(\vphi, x ) \in \R $, $ p(\vphi, x ) \in \C $, and, for NLW, 
\be\label{lin:NLW}
y \mapsto  ({\om \cdot \pa_\vphi})^2 y - y_{xx} + V(x) y + \e a (\vphi, x) y  \, ,
\ee
with $ a(\vphi, x ) \in \R $,
that, in the complex variable 
$$ 
h = D_V^{\frac12} y + \ii D_V^{- \frac12} y_t \, , \quad 
D_V := \sqrt{- \Delta + V(x) } \, ,
$$ 
assumes the form 
\be\label{lin:NLS-bis}
h \mapsto  \om \cdot \pa_\vphi h + \ii D_V h + 
\ii  \frac{\e}{2} D_V^{-\frac12} a(\vphi, x)  D_V^{-\frac12}  ( h + \bar h ) \, .
\ee
Coupling these equations with their complex conjugated component,
we have to invert the quasi-periodic operators, acting on  $ \begin{pmatrix}
h  \\
\bar h       
\end{pmatrix} $, given, for NLS,  by
\be\label{LeqNLS}
\om \cdot \pa_\vphi + 
\begin{pmatrix}
 \ii  \pa_{xx} -  \ii  V(x)  & 0 \\
0   &     - \ii \pa_{xx} +  \ii V(x)        
\end{pmatrix} - \e 
\begin{pmatrix}
 \ii  q (\vphi, x) &  \ii  p(\vphi, x) \\
 - \ii \bar p(\vphi, x)  &  - \ii q (\vphi, x)    
\end{pmatrix} 
\ee
and, for NLW, 
\be\label{LeqNLW}
\om \cdot \pa_\vphi +  \begin{pmatrix}
\ii D_V  & 0 \\
0   &   - \ii D_V       
\end{pmatrix} + \ii \frac{\e}{2} 
D_V^{-\frac12}
\begin{pmatrix}
a (\vphi, x) & a (\vphi, x) \\
-a(\vphi, x)  &  -a(\vphi, x)     
\end{pmatrix} D_V^{-\frac12} \, .
\ee
These linear operators 
 have the form \eqref{inizA}-\eqref{diagDhi} with an operator $ R (\vphi) $ 
which is  bounded. Actually note that for NLW the perturbative term in \eqref{LeqNLW}
 is also $ 1 $-smoothing.
Moreover 
the eigenvalues\index{Eigenvalues of Sturm-Liouville operator}  $ \mu_j^2 $, $j \in \N \backslash \{ 0 \}$, of the Sturm-Liouville operator 
$ - \pa_{xx} + V(x) $\index{Sturm-Liouville operator} with Dirichlet boundary conditions are simple and the quasi-periodic operators 
\eqref{LeqNLS} and \eqref{LeqNLW} take the form \eqref{inizA}-\eqref{diagDhi}, where 
the eigenvalues of $D$ are 
$$ 
\pm \ii \mu_j^2 \, , \ \mu_j^2 \sim j^2 ,  
\ {\rm for \  NLS} \, , \quad \pm \ii \mu_j \, , \ \mu_j \sim j \, , 
\ {\rm  for \ NLW} \, \, .  
$$ 
Then  it is not hard  to impose second order Melnikov  
non-resonance  conditions
 as in \eqref{sec-Me}. 
In view of these observations, it is possible to implement 
the KAM reducibility scheme presented above  
to prove the existence of quasi-periodic solutions for \eqref{NLW1-L}-\eqref{NLS1}
(actually the KAM iteration in \cite{Ku}, \cite{W1} is a bit different but the previous argument catches  its essence).   

Later on these results  have been extended in Kuksin-P\"oschel 
\cite{KP2} to parameter independent Schr\"odinger equations 
\be\label{NLSKAM}
\begin{cases}
\ii u_t = u_{xx} +  f(|u|^2) u \, , \\
u (0) = u(\pi) = 0   \, ,
\end{cases}
\quad {\rm where} \quad f(0)=0 \, , \  f'(0) \neq 0  \, ,  
\ee
and in P\"oschel \cite{Po3} to nonlinear Klein-Gordon equations\index{Klein Gordon} 
\be\label{NLWKAM}
y_{tt} - y_{xx} + m y   = y^3 + {\rm h.o.t}  \, , \quad y (0) = y(\pi) = 0   \, . 
\ee
The main new difficulty of these equations is that the linear equations 
$$
\ii u_t = u_{xx} \, , \quad  y_{tt} - y_{xx} + m y  =  0   \, ,
$$
have resonant invariant tori. 
Actually all the solutions of the first equation, 
\be\label{NLS:linsol}
u(t,x) = \sum_{j \in \Z} u_j(0) e^{\ii j^2 t } e^{\ii j x }\, ,  
\ee
are  $ 2 \pi $-periodic in time (for this reason \eqref{NLSKAM} is called a 
completely resonant PDE) and the Klein-Gordon linear frequencies
$ \sqrt{j^2 + m }$   may be resonant for several values of the  
mass $ m $.  
The new key idea in \cite{Po3}, \cite{KP2} is to compute 
 precisely how the nonlinearity in \eqref{NLSKAM}-\eqref{NLWKAM}  modulates 
the tangential and normal 
frequencies  of the expected quasi-periodic solutions.
In particular, a Birkhoff normal form analysis enables to prove that
the tangential frequencies 
vary  diffeomorphically with the ``amplitudes" of the solutions. This non-degeneracy property 
allows then  to prove that 
the Melnikov non-resonance conditions are satisfied for most amplitudes.
We note however the following difficulty for the equations \eqref{NLSKAM}-\eqref{NLWKAM}
which is not present for \eqref{NLW1-L}-\eqref{NLS1}: 
the  frequency vector $ \om $ may satisfy only a Diophantine condition \eqref{diop} 
with a constant $ \gamma_0 $ which tends to $  0 $  as the solution tends to $ 0 $ 
(the linear frequencies in \eqref{NLS:linsol} are integers), 
and similarly for the second order Melnikov conditions \eqref{sec-Me}.
Note that the remainders in \eqref{new+} 
have size $ O( |R|^2 \g^{-1}) $ and, 
as a consequence, careful estimates have to be performed  to overcome  this  
``singular" perturbation issue.  

\subsection{KAM for $1 d$ NLW and NLS with periodic boundary conditions}

The above results do not apply for periodic boundary conditions $ x \in \T $
because 
two eigenvalues of $ - \pa_{xx} + V(x) $ coincide (or are too close), and thus the 
second order Melnikov non resonance conditions \eqref{sec-Me} are violated.  
This is the first instance where the difficulty mentioned in item \ref{dif1} 
appears.

Historically  this difficulty was first solved  by Craig-Wayne \cite{CW} and Bourgain \cite{Bo1}
developing a multiscale approach, based on repeated use of the ``resolvent identity" in the spirit of 
the work \cite{FS} by Fr\"olich-Spencer for Anderson localization. 
This approach  does not require the second order Melnikov non-resonance conditions.
We describe it in Section \ref{sec:MULTI}. 
Developments of this multiscale approach are  the basis of the present monograph. 

\smallskip
 
 The KAM reducibility\index{Reducibility} approach was  extended  later by Chierchia-You \cite{CY}
for semi-linear wave equations like \eqref{NLW1-L} with periodic boundary conditions. 
Because of the near resonance between pairs of frequencies, 
the linearized operators \eqref{lin:NLW}-\eqref{lin:NLS-bis} are reduced to 
 a diagonal system of $ 2 \times 2 $  self-adjoint matrices, namely of the form \eqref{Dinf}
 with $ d_j^\infty \in {\rm Mat}(2 \times 2; \C) $, by 
 requiring at each step second-order Melnikov non-resonance conditions of the form
\be\label{sec-Me-doppie}
|\om \cdot \ell + d_j - d_k | \geq \frac{\g}{\langle \ell \rangle^\tau} \, , 
\quad \forall (\ell, j, k) \neq (0, j, \pm j ) \, .
\ee 
Note that we do not require in \eqref{sec-Me-doppie} 
non-resonance conditions for $ \ell = 0 $ and $ k = \pm j $. 
 Since NLW is a second order equation, the nonlinear perturbative part of its Hamiltonian vector field 
 is regularizing of order $ 1$ (it gains one space derivative), see \eqref{LeqNLW}, 
 and this is 
 sufficient to prove that 
the perturbed frequencies of \eqref{LeqNLW} satisfy an asymptotic estimate like
$$ 
\mu_j (\e) = \mu_j + O(\e |j|^{-1})  =  |j| + O( |j|^{-1})
$$
as $ |j| \to + \infty $,  where  $ \mu_j^2 $ denote the eigenvalues\index{Eigenvalues of Sturm-Liouville operator}  of the Sturm-Liouville operator $ - \pa_{xx} + V(x) $\index{Sturm-Liouville operator}. 
Thanks to this asymptotic expansion it is sufficient to impose, for each $ \ell \in \Z^\nu $, 
only finitely many second order Melnikov non-resonance conditions as \eqref{sec-Me-doppie}, 
by requiring first order Melnikov conditions like
\be\label{primi-concr}
|\omega \cdot \ell + h | \geq  \gamma  \langle \ell \rangle^{-\tau} \, ,  
\ee
for all $ (\ell, h) \in (\Z^\nu \times \Z) \setminus (0,0) $. 
Indeed, if
 $ |j| , |k| >   C | \ell |^\t \g^{-1}  $, for an appropriate constant $C>0$, we get 
\begin{align}
   |\omega \cdot \ell +\mu_j (\e)  -\mu_k (\e)  | &  
   \geq |\omega \cdot \ell + |j| - |k| |  - O(1/\min(|k|,|j|)) \nonumber 
\\
& \stackrel{\eqref{primi-concr}} \geq   \g  \langle \ell \rangle^{-\tau}  -O(1/\min(|k|,|j|))  \geq \frac{\g}{2} \langle \ell \rangle^{-\tau} \label{CY:wave}
\end{align}
noting that $ |j| - |k| $ is an integer.
Moreover if $||k|-|j|| \geq C| \ell |$ for another appropriate constant $C>0$ then 
$|\om \cdot \ell + \mu_j(\e) -\m_k(\e)| \geq | \ell |$.  Hence, under 
\eqref{primi-concr}, for $ \ell $ given, 
the second order Melnikov conditions with time-index $ \ell $ 
are automatically satisfied for all $(j,k)$ 
except a finite  number. 

\begin{remark} 
If  the Hamiltonian nonlinearity does not depend on the space variable $ x $, the equations
\eqref{NLSKAM}-\eqref{NLWKAM} 
are  invariant under space translations and therefore possess a prime integral by Noether Theorem.
Geng-You \cite{GY0}, \cite{GY} were the first to 
exploit such  conservation  law, which is preserved  along the 
KAM iteration, to  fulfill the non-resonance conditions. The main observation is that 
such symmetry 
enables to prove that many monomials are a-priori never present 
along the KAM iteration.  
In particular, this symmetry removes the degeneracy produced 
by  the multiple normal frequencies. 
\end{remark}

For  semilinear Schr\"odinger equations like \eqref{NLS1}, \eqref{NLSKAM}, 
the nonlinear vector field is not smoothing. Correspondingly, note that in the linearized operator 
\eqref{LeqNLS} the remainder $ R(\vphi) $ is a matrix of multiplication operators. 
In such a case the basic perturbative estimate for the eigenvalues gives
$$
\mu_j (\e) = \mu_j^2 + O(\e) \, , 
$$
which is not sufficient to verify second order 
Melnikov non resonance conditions like \eqref{sec-Me-doppie}, 
in particular
$$ 
| \om \cdot \ell + \mu_j (\e) -  \mu_{-j} (\e) | \geq \frac{\gamma}{\langle \ell \rangle^\tau}  \, , \quad \forall \ell \in \Z^\nu \setminus \{0 \} \, , \ j \in \N \, , 
$$ 
for most values of the parameters. 

The first KAM 
reducibility result for NLS with  $ x \in \T $ has been 
proved by Eliasson-Kuksin in \cite{EK} as a particular case of a much 
more general result valid for tori $ \T^d $ of  any space dimension $ d \geq 1 $, 
that we discuss below. The key point is to extract, using the notion of T\"oplitz-Lipschitz matrices,  
the first order asymptotic expansion of the perturbed eigenvalues. 
For perturbations of $ 1 $-dimensional Schr\"odinger equations,  another recent 
approach to obtain the improved asymptotics of the perturbed frequencies, i.e. 
$$
\mu_j (\e) = j^2  + c + O( \e / |j| ) \, , 
$$
for some constant $ c $ independent of $ j $, 
is developed in Berti-Kappeler-Montalto \cite{BKM}.
This expansion, which allows to verify the second order Melnikov 
non-resonance conditions \eqref{sec-Me-doppie}, is obtained
 via a regularization technique based on pseudo-differential ideas, that we  explain below. 
The approach in \cite{BKM}  applies to semilinear perturbations (also $ x $-dependent) 
of any large 
``finite gap" solution of 
$$
   \ii  u_t  = - \partial_{xx} u +  |u|^2 u +  \e   f  (x, u )  \, , \quad x \in \T  \, .
$$
Let us explain the term ``finite gap" solutions. The $ 1d$-cubic NLS
\be\label{cubicNLS}
   \ii  u_t  = - \partial_{xx} u +  |u|^2 u  \, , \quad x \in \T \, , 
\ee 
possesses global analytic action-angle variables, 
in the form of 
 Birkhoff coordinates,  see \cite{GrK}, and  
the whole infinite dimensional phase space 
is foliated by quasi-periodic -called ``finite gap" solutions- 
and almost-periodic solutions. 
The Birkhoff coordinates are a cartesian smooth 
version of the action-angle variables 
to  avoid the singularity  when one action component vanishes, i.e.  
close to the elliptic equilibrium $ u = 0 $. 
This situation generalizes what happens for a finite dimensional Hamiltonian  system in $ \R^{2n} $ 
which possesses $ n $-independent prime integrals in involution. According to the
celebrated  Liouville-Arnold theorem (see e.g. \cite{Ar-book}), 
in  suitable local symplectic angle-action variables 
 $ (\theta, I ) \in \T^n \times \R^n $,  the integrable Hamiltonian $ H(I) $ depends only on the actions  and the dynamics 
 is described by 
$$
\dot  \theta =  \pa_I H(I) \, , \quad \dot I = 0 \, . 
 $$ 
Thus the phase space is foliated by the invariant tori $ \T^n \times \{ \xi \} $, $ \xi \in \R^n $, filled by the quasi-periodic solutions 
$ \theta (t) = \theta_0 + \omega (\xi)  t $, $ I(t) = \xi  $,  with frequency vector  $ \omega (\xi) = (\pa_I H)(\xi ) $. The analogous construction close to an elliptic equilibrium, where the action-angle variables become singular,  is provided by the R\"ussmann-Vey-Ito theorem 
\cite{Ru}, \cite{Vey}, \cite{Ito}. We refer to \cite{KaP} for an introduction.    

Other integrable PDEs which possess Birkhoff coordinates are KdV \cite{KaP}  and mKdV \cite{KaT},  see Appendix \ref{KdV:int}. 

\begin{remark}
The Birkhoff normal form construction of \cite{KP2} discussed for
the NLS equation \eqref{NLSKAM} provides, close to $ u = 0 $, an 
approximation of the global Birkhoff coordinates of the $ 1d$-cubic NLS.
This is sufficient information to prove a KAM result for small
amplitude solutions. 
\end{remark}

\subsection{Space multidimensional PDEs}

For space multidimensional PDEs 
the reducibility\index{Reducibility} approach has been first worked out for 
semilinear Schr\"odinger equations
\be\label{EKAnn}
- \ii u_t = - \Delta u + V * u + \e \pa_{\bar u} F(x, u, \bar u ) \, , \quad x \in \T^d \, , 
\ee 
with a convolution potential by Eliasson-Kuksin \cite{EK}, \cite{EK1}. 
This is a much more difficult situation with respect to
the $ 1$-$d$-case because  the eigenvalues of $ - \Delta + V(x)  $ 
appear in clusters of unbounded size. This is the difficulty mentioned in item \ref{dif1}. 
In such a case the reducibility result that  one could look for 
is to block-diagonalize the quasi-periodic Schr\"odinger linear operator 
$$ 
h \mapsto 
\ii \om \cdot \pa_\vphi h - \Delta h + V * h + \e q(\vphi, x) h  + \e p(\vphi, x) \bar h  \, , 
$$
i.e. to obtain an operator as
\eqref{Dinf} with finite dimensional  blocks 
$ d_j^\infty $, which are self-adjoint matrices  of increasing dimension as $ j \to +\infty $. 
The convolution potential $ V $ plays the role of ``external parameters". 
Eliasson-Kuksin  introduced in \cite{EK}  the notion of T\"oplitz-Lipschitz matrices in order to 
extract asymptotic information on the eigenvalues, 
and so verify 
the second order Melnikov non resonance conditions. 
The quasi-periodic solutions of \eqref{EKAnn} obtained in  \cite{EK} are linearly stable. 

\begin{remark}
The reducibility techniques in  \cite{EK} enable to prove a stability result for
all the solutions of 
the  linear Schr\"odinger equation 
$$
\ii u_t =  \Delta u + \e V(\omega t, x )  \, , \quad x \in \T^d \, ,
$$
with a small quasi-periodic analytic potential $ V(\omega t, x ) $. 
For all frequencies $ \omega \in \R^\nu $, except a set of 
measure tending to $ 0 $ as  $ \e \to 0 $,  
the Sobolev norms of any solution $ u(t, \cdot ) $ satisfy
$$
\| u (t, \cdot) \|_{H^s(\T^d)} \sim \| u(0, \cdot)\|_{H^s(\T^d)} \, , \quad \forall t  \in \R \, . 
$$
\end{remark}

Subsequently for the cubic NLS equation 
\be\label{NLS-cubic}\index{Nonlinear Schr\"odinger equation}
\ii u_t = - \Delta u + |u|^2 u \, , \quad   \ x \in \T^2 \,     , 
\ee
which is parameter independent and  completely resonant,  Geng-Xu-You \cite{GXY} proved  
a KAM result, without reducibility,  using a Birkhoff normal form analysis.
We remark that the Birkhoff normal form  of \eqref{NLS-cubic} 
is not-integrable, unlike  in space dimension $ d = 1 $, 
causing additional difficulties with respect to \cite{KP2}.  

For  completely resonant NLS equations in any space dimension and a polynomial nonlinearity 
\be\label{c-res:NLS}
\ii u_t = - \Delta u + |u|^{2p} u \, , \quad p \in \N \, , \ x \in \T^d \,  , 
\ee
Procesi-Procesi \cite{PP1} realized a  systematic study of the resonant Birkhoff normal form 
and, using the notion  of quasi-T\"oplitz matrices developed in Procesi-Xu \cite{PX},  
proved in \cite{PP}, \cite{PP3}, the existence  of reducible quasi-periodic
solutions of  \eqref{c-res:NLS}.

\begin{remark}
The resonant Birkhoff normal form  of \eqref{NLS-cubic} is exploited in \cite{CKSTT}, 
\cite{GK}, to construct chaotic orbits with a growth of the Sobolev norm. 
This ``norm inflation"  phenomenon is reminiscent of the 
Arnold diffusion problem \cite{Ar1} for finite dimensional Hamiltonian systems.
Similar results for  the  NLS equation 
\eqref{c-res:NLS} have been proved in \cite{GHP}. 
\end{remark}

 KAM results have been  proved for parameter dependent beam equations by Geng-You 
 \cite{GY1.5}, Procesi \cite{Pbeam}, 
 and, more recently, in Eliasson-Gr\'ebert-Kuksin \cite{EGK1} for
multidimensional 
beam equations\index{Beam equation} like 
$$
u_{tt} + \Delta^2 u + m u + \pa_u G(x, u) = 0 \, , \quad x \in \T^d \, , \quad u \in \R \, . 
$$
We also mention  the work \cite{GrP} of Gr\'ebert-Paturel  concerning the existence of reducible quasi-periodic solutions of Klein-Gordon\index{Klein Gordon} equations on the sphere
 $ {\mathbb S}^d $, 
\be\label{GPSd}
 u_{tt} - \Delta u + m u+ \delta M_\rho u + \e g(x, u) = 0 \, , \quad x \in {\mathbb S}^d \, , \quad u \in \R \, , 
\ee
where $ \Delta $ is the Laplace-Beltrami operator and $ M_\rho $ is a Fourier multiplier.  

On the other hand,
 if $ x \in \T^d $,  the infinitely many 
 second order Melnikov conditions, required to  block-diagonalize the 
 quasi-periodic linear  wave operator 
$$ 
({\om \cdot \pa_\vphi})^2  - \Delta + V *   \, +  \, \e a (\vphi, x)   \, , \quad x \in \T^d \, , 
  $$ 
 are violated for $ d \geq 2 $, 
 and no  reducibility results are available so far. 
Nevertheless, results of ``almost" reducibility have been announced 
in  Eliasson \cite{E17}, Eliasson-Gr\'ebert-Kuksin \cite{EGK}. 

Before concluding this subsection, 
we also mention the KAM result by Gr\'ebert-Thomann \cite{GT}  for 
smoothing nonlinear perturbations of the $ 1d$ harmonic oscillator and 
Gr\'ebert-Paturel \cite{GrP1} in  higher space dimension. 

\subsection{$ 1$-d quasi and fully nonlinear PDEs, Water Waves}

Another situation where the reducibility  approach that we described
in subsection \ref{sec:PR}
encounters a serious difficulty is when the non-diagonal  
remainder $ R (\vphi) $ is unbounded.
This is the difficulty mentioned in item \ref{dif2}. 
In such a case, the auxiliary vector field $ F (\vphi) $ defined in \eqref{defF} is unbounded as well.
Therefore it may not define a flow, and the 
 iterative reducibility scheme described in subsection \ref{sec:PR}
 would formally produce
remainders which accumulate more and more derivatives.  

\subsubsection{KAM for semilinear PDEs with derivatives} 

The first KAM results for PDEs with an {\it unbounded} nonlinearity
have been proved by Kuksin \cite{K2-KdV} and, then, Kappeler-P\"oschel \cite{KaP}, 
for perturbations of finite-gap solutions of\index{Perturbed KdV equation}   
\begin{equation}\label{KdV-QL-semi}
u_t + u_{xxx} + \pa_x u^2 + \e \pa_x (\pa_u f)(x, u ) = 0 \, , \quad x \in \T \, .
\end{equation}
The corresponding quasi-periodic 
linearized operator at an approximate quasi-periodic solution $ u  $ 
has the form
$$
h \mapsto \om \cdot \pa_\vphi  h + \pa_{xxx} h + \pa_x ( 2 u h ) + \e \pa_x ( a h ) \, ,
\quad a :=  (\pa_{uu} f)(x, u) \, . 
$$
The key  idea in \cite{K2-KdV} is to exploit the fact that the frequencies  of KdV grow asymptotically as $  \sim j^3 $ 
as $ j \to + \infty $.
Therefore one can impose second order Melnikov non-resonance conditions like
$$ 
| \om \cdot \ell +  j^3  - i^3| \geq \gamma (j^2 + i^2)/ \langle \ell \rangle^{\tau}  \, , \quad i \neq j \, , 
$$ 
which gain $ 2 $ space derivatives
(outside the diagonal $ i = j $), sufficient to compensate the loss of one space derivative 
produced by the  vector field $ \e \partial_{x} (\partial_u f) (x, u) $.
On the diagonal $ \ell \neq 0,  i = j $,  one does not solve the homological equations,
with the consequence  that the KAM normal form is $ \vphi $-dependent.
In order to overcome this difficulty one can make  use of the so called 
Kuksin Lemma to invert the corresponding quasi-periodic scalar operator. 
Subsequently, developing an improved version of the Kuksin Lemma,  
Liu-Yuan  in \cite{Liu-Yuan} proved KAM results for 
semilinear perturbations of  Hamiltonian derivative NLS and  Benjiamin-Ono equations
and Zhang-Gao-Yuan \cite{ZGY} for the reversible\index{Reversible PDE} derivative NLS equation
$$ 
\ii u_t + u_{xx} = | u_x |^2 u  \, , \quad  u(0) = u(\pi) = 0  \, .
$$
These PDEs  are more difficult  than KdV  
because 
the linear frequencies 
grow  like $ \sim j^2 $ and not $ \sim j^3 $, and therefore one
gains only $ 1 $ space derivative when solving the homological equations.  

These methods do not apply for derivative wave equations where the dispersion relation is asymptotically linear.  
Such a case has been addressed more recently by  Berti-Biasco-Procesi \cite{Berti-Biasco-Procesi-Ham-DNLW}-\cite{Berti-Biasco-Procesi-rev-DNLW} who 
proved 
 the existence and the stability of quasi-periodic solutions of autonomous derivative Klein-Gordon equations 
 \be\label{DNLW}
 y_{tt} - y_{xx} + m y =
g(x, y, y_x, y_t) 
\ee
satisfying reversibility conditions. These  assumptions
 rule out  nonlinearities like $ y_t^3 $, $ y_x^3 $, for which no 
periodic nor quasi-periodic solutions exist (with these nonlinearities all the solutions dissipate to zero).
The key point in \cite{Berti-Biasco-Procesi-Ham-DNLW}-\cite{Berti-Biasco-Procesi-rev-DNLW} was to adapt the  notion of quasi-T\"oplitz vector field
introduced in \cite{PX}  
 to obtain the higher order  {\it asymptotic expansion} of the perturbed normal 
 frequencies 
$$
 \mu_j (\e) = \sqrt{j^2 + m } + a_\pm + O(1/j) \, , \quad {\rm as} \ j \to \pm \infty \, , 
$$
for suitable constants $ a_\pm $ (of the size $ a_\pm =  O(\e) $ of the solution $ y = O(\e) $). 
Thanks to this asymptotic expansion it is sufficient to verify, for each $ \ell \in \Z^\nu $, that  
only finitely many second order Melnikov non-resonance conditions hold.
Indeed, using an argument as in \eqref{CY:wave}, 
infinitely many conditions in \eqref{sec-Me-doppie} are already verified by imposing only
 first order Melnikov conditions like
$$
|\omega \cdot \ell + h | \geq  \gamma  \langle \ell \rangle^{-\tau} \, ,  \qquad 
|\omega \cdot \ell +  (a_+ - a_-) + h | \geq  \gamma  \langle \ell \rangle^{-\tau} \, , 
$$
for all $ (\ell, h) \in \Z^\nu \times \Z $ such that $ \omega \cdot \ell + h $
and $ \omega \cdot \ell +  (a_+ - a_-) + h  $
do not vanish identically. 

\subsubsection{KAM for quasi-linear and fully nonlinear PDEs}

All the above results still concern semi-linear perturbations, namely when the number of derivatives which are present in the nonlinearity is strictly lower than the order of the linear differential operator. 
The first existence results of small amplitude
quasi-periodic solutions for quasi-linear PDEs have been proved 
by Baldi-Berti-Montalto   in \cite{BBM-Airy} 
for fully nonlinear perturbations of the Airy equation\index{Airy equation}
\be\label{Ayry-P}
u_{t} + u_{xxx} + \e f(\omega t , x , u, u_{x}, u_{xx}, u_{xxx} ) = 0 \, , \quad 
x \in \T \, ,  
\ee
and in \cite{BBM-auto}  
for  quasi-linear autonomous 
perturbed KdV equations\index{Perturbed KdV equation} 
\begin{equation}\label{KdV-QL}
u_t + u_{xxx} + \pa_x u^2 + {\cal N} (x, u, u_x, u_{xx}, u_{xxx}) = 0   
\end{equation}
where the Hamiltonian nonlinearity
\be\label{qlpert}
{\cal N} (x, u, u_x, u_{xx}, u_{xxx}) := 
- \partial_x \big[ (\partial_u f)(x, u,u_x)  - \partial_{x} ((\partial_{u_x} f)(x, u,u_x)) \big]  
\ee
vanishes at the origin as $ O(u^4) $. See \cite{Giu} when $ {\cal N} $  vanishes only quadratically. 

\smallskip

The main new tool which has been introduced to solve this problem is a systematic use of
{\it pseudo-differential} calculus. 
The key point is to {\it reduce} to constant coefficients  the linear PDE 
\be\label{eq-KdV-L}
u_t  + (1 + a_3(\om t ,x)) u_{xxx} + a_2(\om t,x) u_{xx} + a_1(\om t,x) u_x + a_0 (\om t,x) u = 0 \, , 
\ee
which is obtained by linearizing \eqref{Ayry-P} at an approximate quasi-periodic solution $ u(\om t , x)  $.
The coefficients $ a_i (\om t, x ) $, $ i = 0, \ldots, 3 $ have size $  O(\e ) $.  
Instead of trying to diminish the size of the $ (t,x) $-dependent terms in 
 \eqref{eq-KdV-L}, as in the scheme  outlined in subsection 
\ref{sec:PR} --the big difficulty \ref{dif2} would appear--, the aim is to 
conjugate \eqref{eq-KdV-L}  to a system like 
\be \label{Kred2}
u_t + m_3 u_{xxx} + m_1 u_x + {\cal R}_0 (\om t ) u = 0 
\ee
where $ m_3 =  1 + O(\e)$,
$  m_1 = O(\e ) $ are constants and $ {\cal R}_0 (\om t ) $ is a zero order operator,
still time dependent.
To do this,  \eqref{eq-KdV-L} is conjugated with a time quasi-periodic 
change of variable  (as in \eqref{QP-floquet})
\be\label{comp-diff}
 u = \Phi (\om t )[v]  = v( t , x + \beta (\om t, x)) \, , 
\ee
 induced by the composition with a diffeomorphism $ x \mapsto x + \beta (\vphi, x) $ of
$ \T_x $ (requiring $ |\b_x | < 1 $).
The conjugated system \eqref{sis2}-\eqref{newB} is
$$
v_t + \Phi^{-1} (\om t )(  (1+ a_3 (\om t, x)) (1+ \b_x (\om t, x))^3 ) v_{xxx} ( t, x ) + {\rm lower \ order \ operators} = 0 
$$
and  therefore one chooses a periodic function $ \beta (\vphi, x) $  such that 
$$ 
(1+ a_3 (\vphi, x)) (1+ \b_x (\vphi, x))^3 = m_3 (\vphi) 
$$ 
is independent of $  x $. 
Since $ \beta_x (\vphi, x ) $ has zero space average, this is possible with
$$
m_3 (\vphi) = \Big(\frac{1}{2 \pi} \int_{\T} \frac{dx}{(1+ a_3 (\vphi, x))^{\frac13}} \Big)^{-3} \, . 
$$
The $ \vphi $ dependence of $ m_3 (\vphi) $ can also be eliminated  at the highest order 
using a  quasi-periodic reparametrization of time and, using other pseudo-differential transformations,
we can reduce also the lower order terms 
to constant coefficients  obtaining \eqref{Kred2}.
The reduction \eqref{Kred2} 
implies the accurate asymptotic expansion of the perturbed frequencies 
$$ 
\mu_j (\e) =  - \ii m_3 j^3 +  \ii m_1 j + O(\e) \, .
$$
Now it is possible to 
verify the second order Melnikov non-resonance conditions  required by a KAM  reducibility\index{Reducibility} scheme (as outlined in subsection \ref{sec:PR})
to diagonalize $ {\cal R}_0 (\om t ) $, 
completing the  reduction of  \eqref{Kred2}, thus \eqref{eq-KdV-L}.

\smallskip

These  techniques have been then employed by Feola-Procesi \cite{FP} 
for  quasi-linear forced perturbations of   Schr\"odinger equations and in \cite{CFP}, \cite{FP1} for 
the search of small amplitude analytic solutions of autonomous PDEs. 
These kind of ideas have been also successfully generalized
for unbounded perturbations of 
harmonic oscillators by Bambusi  \cite{Bam1}, \cite{Bam2} and
Bambusi-Montalto \cite{BaMo}.

\begin{remark}
The recent paper  \cite{BKM1} of  
Berti-Kappeler-Montalto 
proves   the persistence  of also large finite gap solutions
of the perturbed quasi-linear 
KdV equation \eqref{KdV-QL}. 
This answers positively a longstanding question 
concerning 
the persistence of 
 quasi-periodic solutions, with arbitrary size,  of integrable PDEs subject  
 to strongly nonlinear perturbations. 
\end{remark}

The KdV and the NLS equation are partial differential equations and the pseudo-differential tools required are essentially 
commutators of multiplication operators and Fourier multipliers. On the other hand, 
for the water waves equations, that we now present, 
the theory of pseudo-differential operators has to be used in full strength.

\subsubsection{Water waves equations}

The water waves equations\index{Water waves equations} for a
perfect, incompressible, inviscid, irrotational fluid
occupying the time dependent region
\be\label{domain-fluid}
{\cal D}_\eta := \big\{ (x, y) \in \T \times \R \, : \, - h < y < \eta (t,x) \big\} \, , \quad \T := \T_x := \R / 2 \pi \Z \, , 
\ee
under the action of gravity, and  possible capillary forces at the free surface, 
are the Euler equations of hydrodynamics  combined with conditions at the boundary 
of the fluid: 
\begin{equation}\label{water}
\begin{cases}
\partial_t \Phi + \frac12 | \nabla \Phi |^2 + g \eta = \kappa 
\pa_x \Big( \frac{ \eta_x }{\sqrt{1 +   \eta_x^2}} \Big)    
\qquad  \qquad \qquad  \, {\rm at} \ y = \eta (t, x)  \cr
\Delta \Phi =0 \qquad \qquad \qquad \  \  \,   \qquad \qquad \qquad  \qquad \qquad \ \qquad {\rm in} \ {\cal D}_{\eta}  \cr
\pa_y \Phi = 0  \qquad  \qquad \qquad \quad \qquad \qquad \qquad \qquad \qquad \qquad  {\rm at} \   y = - h   \cr
\partial_t \eta = \partial_y \Phi - \partial_x \eta 
\partial_x \Phi \qquad \qquad \qquad \quad \qquad \qquad \qquad \ \  {\rm at} \  y = \eta (t, x) 
\end{cases}
\end{equation}
where $ g $ is the acceleration of gravity and $ \kappa $
 is the surface tension coefficient. The unknowns of the problem \eqref{water} 
are the free surface $ y = \eta (t, x) $
and the  velocity potential $ \Phi : {\cal D}_\eta \to \R $, 
i.e. the irrotational  velocity field   of the fluid
$  v =\nabla_{x,y} \Phi  $.  
The first equation in \eqref{water}  is the Bernoulli condition according to which
 the jump of pressure across the free surface is proportional to the mean curvature. 
 The second equation  in \eqref{water} is the incompressibility property
$ {\rm div} \, v = 0 $. 
 The third equation expresses the impermeability of the bottom of the ocean. 
 The last condition in \eqref{water} 
means  that the fluid particles on the free surface $ y = \eta (x, t) $ remain forever on it along the fluid evolution. 

Following Zakharov \cite{Zakharov1968} and Craig-Sulem \cite{CrSu}, the evolution problem \eqref{water} may be written as an infinite-dimensional Hamiltonian system
in the  unknowns $ (\eta(t,x), \psi(t,x) ) $ where
$ \psi(t,x)=\Phi(t,x,\eta(t,x))   $
is, at each instant $ t $,   the trace at the free boundary of the velocity potential. 
Given 
$\eta(t,x)$ 
and 
$\psi(t,x) $  
there is a unique solution 
$\Phi(t,x,y)$ of the elliptic problem 
$$
\begin{cases}
\Delta \Phi = 0 & \text{in } 
\{-  h < y < \eta(t,x)\} \\
\partial_y \Phi = 0 & \text{on } y = - h \\
\Phi = \psi & \text{on } \{y = \eta(t,x)\} \, . 
\end{cases}
$$
System \eqref{water} is then equivalent to the Zakharov-Craig-Sulem system
\begin{equation}\label{WW}
\begin{cases}
\pa_t \eta =  G(\eta) \psi \\
\pa_t \psi = - g \eta - \dfrac{\psi_x^2}{2} + \dfrac{1}{2(1+\eta_x^2)} \big( G(\eta) \psi + \eta_x \psi_x \big)^2 +
\kappa \pa_x \Big( \frac{ \eta_x }{\sqrt{1 +   \eta_x^2}} \Big)  
\end{cases}
\end{equation}
where  $ G(\eta) := G(\eta; h )$ is the Dirichlet-Neumann operator defined by
\be\label{DN}
G(\eta) \psi := \big( \Phi_y - \eta_x \Phi_x \big)_{|y = \eta (t,x)}  \, .
\ee
Note that $ G(\eta ) $ maps the Dirichlet datum $\psi$ 
to the (normalized) normal derivative 
of $\Phi$ at the top boundary. 
The operator $ G(\eta) $ is linear in $ \psi $, 
self-adjoint with respect to the $ L^2 $ scalar product, positive-semidefinite, 
and its kernel contains only the constant functions. 
The Dirichlet-Neumann operator depends smoothly 
on the wave profile $ \eta $, and it is a {\it pseudo-differential} operator with principal symbol 
$ D \tanh (hD) $.

Furthermore the equations \eqref{WW} are the Hamiltonian system
\be\label{HS}
 \pa_t \eta = \nabla_\psi H (\eta, \psi) \, , \quad  \pa_t \psi = - \nabla_\eta H (\eta, \psi)  \, ,
\ee
where $ \nabla $ denotes the $ L^2 $-gradient,  
and  the Hamiltonian 
\be\label{Hamiltonian}
H(\eta, \psi) = 
 \frac12 \int_\T \psi \, G(\eta, h ) \psi \, dx + \frac{g}{2} \int_{\T} \eta^2  \, dx  +  \kappa \int_{\T}  \sqrt{1 + \eta_x^2} \, dx 
\ee
is the sum of the kinetic,   potential and capillary energies
expressed in terms of the variables  $ (\eta, \psi ) $.

The water waves system \eqref{WW}-\eqref{HS} exhibits several symmetries. 
First of all, the mass 
$$ 
\int_\T \eta ( x) \, dx  
$$ 
is a first integral of \eqref{WW}. 
Moreover \eqref{WW} is invariant under spatial translations
and Noether's theorem implies that the momentum 
$$
 \int_{\T} \eta_x (x) \psi (x)  \, dx 
$$
is a prime integral of \eqref{HS}. 
In addition, the subspace of functions that are even in $ x $, 
\be\label{even in x}
\eta (x) = \eta (-x) \, , \quad \psi (x) = \psi (- x )  \, , 
\ee
is invariant under \eqref{WW}. 
In this case also the velocity potential $ \Phi(x,y) $ is even  and $ 2 \pi $-periodic in $ x $ and so the 
$ x$-component  of the velocity field $ v = (\Phi_x, \Phi_y) $ vanishes  at $ x = k \pi $, for all $ k \in \Z $. 
Hence   there is no flow of fluid through the lines $ x = k \pi $, $ k \in \Z $, and 
a solution of \eqref{WW} satisfying 
 \eqref{even in x} describes 
the motion of a liquid confined  between two vertical walls. 

We also note that the water waves system  
 \eqref{WW}-\eqref{HS} is reversible\index{Reversible PDE} with respect to the involution  $ S : (\eta, \psi) \mapsto (\eta, - \psi) $,
i.e.  the Hamiltonian $H$ in \eqref{Hamiltonian} 
 is even in $ \psi $, see Appendix \ref{App0}.
As a consequence it is natural to look for solutions of \eqref{WW} satisfying 
\be\label{odd-even}
u(-t ) = S u(t) \, , \quad i.e. \quad \eta (-t, x) = \eta(t, x) \, , \ 
\psi(-t,x ) = - \psi (t,x ) \quad \forall t, x \in \R \, .
\ee
Solutions of the water waves equations \eqref{WW} satisfying 
\eqref{even in x} and \eqref{odd-even} are called  {\it standing water waves}.

The phase space of \eqref{WW} containing $(\eta, \psi)$ is (a dense subspace) of 
$H^1_0 (\T) \times {\dot H}^1 (\T)$; $H^1_0(\T)$ denotes the subspace of $H^1(\T)$ of zero average functions;
$$
{\dot H}^1 (\T) := H^1 (\T) \slash_{ \sim} 
$$
is the quotient space of $H^1 (\T)$  by the equivalence relation $\sim$ defined in this way :
$ \psi_1  \sim \psi_2  $ if and only if $ \psi_1 (x) - \psi_2 (x) = c $ is a constant.  
For simplicity of notation we denote the equivalence class $ [\psi] $ by $ \psi $. Note that the second equation in \eqref{WW}
is in $  {\dot H}^1 (\T) $, as it is natural because only the gradient of the velocity potential has a physical meaning.

 Linearizing \eqref{WW} at the equilibrium $(\eta, \psi) = (0,0)$ we get 
\begin{equation} \label{Lom}
\left\{
\begin{aligned}
&\partial_t \eta = G(0)\psi,\\
&\partial_t\psi = - g \eta  + \kappa  \eta_{xx}  
\end{aligned}
\right.
\end{equation}
where  $ G(0) = D \tanh (hD)   $ 
is the Dirichlet-Neumann operator at the flat surface $ \eta = 0  $. 
The linear frequencies of oscillations of \eqref{Lom} are
\be\label{Linear-Fre-WW}
\om_j = \sqrt{j \tanh (h j) (g  + \kappa j^2)} \, , \quad j \in \Z \, ,
\ee
which, in the phase space of even functions  \eqref{even in x},
 are {\it simple}.  Also note 
 that 
$$ 
\begin{aligned}
& {\rm if} \  \kappa > 0 \ {\rm (capillary-gravity \ waves)}  \quad \Longrightarrow \quad
\omega_j \sim  |j|^{3/2} \  {\rm as} \ |j| \to \infty \, , \\
& {\rm if} \  \kappa = 0  \  {\rm (pure \ gravity \ waves)} \qquad \ \ \quad \Longrightarrow \quad
 \omega_j \sim |j|^{1/2}  \  {\rm as} \ |j| \to \infty \, . 
 \end{aligned}
$$

\noindent
{\bf {KAM for water waves.}}
The first  existence results of small amplitude time-periodic gravity standing wave solutions 
for bi-dimensional fluids has been proved by Plotinkov and Toland \cite{PlTo}
in finite depth and by Iooss, Plotnikov and Toland in \cite{Ioo-Plo-Tol}  
in infinite depth. 
More recently, the existence of time periodic gravity-capillary standing wave solutions in infinite depth 
has been proved by Alazard-Baldi \cite{AB}. 

The main result in \cite{BM16}  proves that most of the standing wave solutions 
of the linear system \eqref{Lom}, 
which are  supported on finitely many space Fourier modes,  
can be continued to standing wave solutions of the nonlinear water-waves  system \eqref{WW} for most 
values of the surface tension parameter $ \kappa  \in [\kappa_1, \kappa_2] $.

A key step is the reducibility to constant coefficients of the quasi-periodic operator 
 $ {\cal L}_\om $ obtained 
linearizing \eqref{WW} at a quasi-periodic approximate solution. 
After the introduction of a linearized Alinhac good unknown, 
and using in a systematic way pseudo-differential calculus, it is possible to transform 
$ {\cal L}_\om $ into
a complex quasi-periodic linear operator of the form  
\be\label{forma-lineare-preliminary}
(h, \bar h) \mapsto  \big(  \om \cdot \pa_\vphi   + \ii \mathtt m_3   
|D|^{\frac12} (1 - \kappa \pa_{xx})^{\frac12} + \ii
{\mathtt m}_1 |D|^{\frac12} \big) h  + {\cal R} (\vphi) [h, \bar h]
\ee
where $ \mathtt m_3, {\mathtt m}_1 \in \R $ are constants  satisfying
$ \mathtt m_3 \approx 1 $, $ {\mathtt m}_1 \approx 0 $, 
and the remainder  $ {\cal R} (\vphi) $ is a small bounded operator. 
Then a KAM reducibility scheme
completes the diagonalization of the linearized operator $ {\cal L}_\om $.  
The required second order Melnikov non-resonance conditions 
are fulfilled for most values of the surface tension parameter $ \kappa $
generalizing   ideas of degenerate KAM theory
for PDEs  \cite{BaBM}.  One exploits that the linear frequencies $ \om_j  $ 
in \eqref{Linear-Fre-WW}
are analytic and non degenerate in $ \kappa $, and the sharp asymptotic expansion of
the perturbed frequencies obtained by the regularization procedure.

In the case of pure gravity water waves, i.e. $ \kappa = 0 $, 
the linear frequencies of oscillation are (see \eqref{Linear-Fre-WW})
\be\label{LIN:fre}
\om_j := \om_j ( h) :=  \sqrt{g j \tanh(h j)} \, , \quad j \geq 1 \, ,
\ee
and three major further difficulties 
in proving the existence of time quasi-periodic solutions\index{Quasi-periodic solution}  
are: 
\begin{itemize}
\item[$(i)$] The nonlinear water waves system \eqref{WW} (with $ \kappa = 0 $) is a singular perturbation
of \eqref{Lom} (with $ \kappa = 0 $) in the sense that
the quasi-periodic linearized operator assumes the form 
\[
\om \cdot \pa_\ph + \ii |D|^{\frac12} \tanh^{\frac12} (h |D|)+ V(\ph,x) \pa_x
\]
and  the term $V(\ph,x) \pa_x$ is now a  \emph{singular} perturbation of the linear dispersion relation operator
$  \ii |D|^{\frac12} \tanh^{\frac12} (h |D|) $
(on the contrary, for the gravity-capillary case 
the transport term $V (\ph,x) \pa_x $ is a lower order perturbation of $|D|^{\frac32} $, see \eqref{forma-lineare-preliminary}).

\item[$(ii)$] The dispersion relation \eqref{LIN:fre} is sublinear, i.e.
$ \om_j  \sim \sqrt{j} $ for $ j \to \infty $, and therefore 
it is only possible to impose second order Melnikov non-resonance conditions 
as in \eqref{sec-Me-loss} which lose space derivatives. 

\item[$(iii)$] 
The linear frequencies $ \omega_j ( h) $ in \eqref{LIN:fre} 
vary with $ h $  of just exponentially small quantities.
\end{itemize}

The main result in Baldi-Berti-Haus-Montalto  \cite{BBHM}  proves 
the existence  of pure gravity standing water waves solutions.

The difficulty (i) is solved proving a  straightening theorem for a quasi-periodic transport operator:
there is a  quasi-periodic change of variables of the form 
$ x \mapsto x + \b(\vphi ,x) $ which conjugates 
$$
\ompaph + V(\ph,x) \pa_x 
$$
 to the constant coefficient vector field $\ompaph $,  
 for $ V (\vphi, x) $ small. 
 This perturbative rectification result 
 is a classical small divisor problem, solved
for perturbations of a Diophantine vector field 
at the beginning of KAM theory, see e.g.  \cite{Z1, Z2}.
Note  that, despite the fact that $  \om \in \R^\nu $ is Diophantine\index{Diophantine vector},  the constant vector field $ \om \cdot \pa_\vphi  $
is resonant on the higher dimensional torus $ \T^{\nu}_{\vphi} \times \T_x $. 
We exploit in a crucial way the {\it symmetry} induced by the {\it reversible} structure\index{Reversible vector field} of the water waves equations, i.e.  $ V(\vphi, x ) $ is odd in $ \vphi $.
 This problem  amounts to proving that {\it all} the solutions 
of the quasi periodically time-dependent scalar characteristic 
equation $ \dot x = V(\om t, x )$, $ x \in \T $, are quasi-periodic in time with frequency $ \omega $.

The difficulty (ii) is overcome performing a regularizing procedure which conjugates the linearized operator, 
obtained along the Nash-Moser iteration, 
to a diagonal, constant coefficient linear one, up to a sufficiently {\it smoothing} operator. 
In this way the subsequent KAM reducibility\index{Reducibility} scheme converges also in presence of  
very  weak Melnikov non-resonance conditions as in \eqref{sec-Me-loss}
which lose space derivatives.  
This regularization strategy is in principle applicable to a broad class of PDEs
where the second order Melnikov non-resonance conditions lose space derivatives. 

The difficulty (iii)  is solved by an improvement of the degenerate KAM theory for PDEs in 
Bambusi-Berti-Magistrelli \cite{BaBM}, which allows to prove that all the Melnikov non-resonance 
are fulfilled for most values of $ h $. 

\begin{remark}
We can introduce 
the space wavelength $ 2  \pi \varsigma $ as an {\it internal} free parameter in the water waves equations
\eqref{WW}. 
Rescaling properly  time, space and amplitude of the solution $(\eta(t,x), \psi(t,x))$ 
we obtain   
system \eqref{WW} where the gravity constant $ g = 1 $ and the depth parameter $ h$ 
depends linearly on  $  \varsigma  $.
In this way  \cite{BBHM} proves existence results for a {\it fixed} 
equation, i.e. a fixed depth $ h $, for most values of the space wavelength $ 2 \pi \varsigma $. 
\end{remark}

\section{The multiscale approach to KAM for PDEs}\label{sec:MULTI}

We present now some key ideas about another approach
developed for analyzing linear quasi-periodic systems, 
in order to prove KAM results for PDEs, started with the seminal paper of 
Craig-Wayne \cite{CW} and
strongly extended by Bourgain \cite{Bo1}-\cite{B5}. 
This set of ideas and techniques  -referred as  ``multiscale analysis"- 
is at the basis of the present monograph. 
For this reason 
we find convenient to present it in some detail. 

For definiteness we consider 
a quasi-periodic linear wave operator 
\be \label{model-B-int} 
(\om \cdot \pa_{\vphi})^2  - \Delta + m  + \e b(\vphi, x) \, , \quad  \vphi \in \T^\nu  \, , \
x \in \T^d \, ,  
\ee
acting on a dense subspace of $ L^2 (\T^\nu \times \T^d ) $,  
where $ m > 0 $ and $ b(\vphi, x)  $ is a smooth function. A linear operator as
in \eqref{model-B-int}  is obtained by linearizing a nonlinear wave equation 
at a smooth approximate quasi-periodic solution.

\begin{remark}\label{rem:ANA} The choice of the parameters in  \eqref{model-B-int}
makes a significant difference. 
If \eqref{model-B-int} arises linearizing a quasi-periodically forced NLW equation as
\eqref{NLW:forced-in}, 
the frequency vector $ \omega  $ can be regarded as a free parameter
belonging to a subset of $ \R^\nu $, independent of $ \e $. 
On the other hand, if \eqref{model-B-int} arises linearizing an 
autonomous NLW equation like \eqref{NLW1}, or \eqref{KG-NON}, the frequency $ \om  $ and $ \e $ are linked.
In particular $ \om $ has to be  
$ \e^2 $-close 
to some  frequency vector $ ( (|j_i|^2 + m)^{1/2} )_{i=1, \ldots, \nu} $ 
 (actually $ \om $ has to belong to the region  of 
 ``admissible" frequencies as in \eqref{def:admissible}). 
This difficulty is present in this monograph 
as well as in \cite{Wang1}.  
\end{remark}

In the exponential basis $ \{ e^{\ii (\ell \cdot \vphi + j \cdot x )} \}_{\ell \in \Z^\nu, j \in \Z^d } $ 
the linear operator \eqref{model-B-int}
 is represented by the self-adjoint matrix 
\be\label{DomT} 
\begin{aligned}
& A := D  +  \e T \, , \\
& D := {\rm Diag}_{(\ell, j) \in \Z^\nu \times \Z^d} \big( 
- (\om \cdot 	\ell)^2 + |j|^2 + m \big) \, , \quad 
T = (T_{\ell,j}^{\ell',j'}):= (\widehat b_ {\ell - \ell', j - j'}) \, ,  
\end{aligned}
\ee
where  
$ \widehat b_ {\ell - \ell', j - j'} $ are the Fourier coefficients of the function
$ b(\vphi, x ) $; these coefficients
 decay rapidly to zero as $ |(\ell-\ell', j - j')| \to + \infty $, exponentially fast,
if  $ b(\vphi, x )  $ is analytic, polynomially, if  $ b(\vphi, x )  $ is a Sobolev function. 
Note that the matrix $ T $ is T\"oplitz, namely it has constant entries 
on the diagonals $ (\ell-\ell', j-j') = (L, J )  \in \Z^\nu \times \Z^d $.  

\begin{remark}
The analytic/Gevrey setting has been considered in \cite{Bo1}-\cite{B5}
and the Sobolev case in \cite{BBARMA}-\cite{BCP}, as well as in Chapter \ref{sec:multiscale} of the present monograph. 
The Sobolev regularity of $  b(\vphi, x) $ has to be large enough, see 
Remark \ref{rem:SOBO}. 
In KAM for PDE applications this requires 
that the nonlinearity (thus the solution) is sufficiently many times differentiable.
For finite dimensional Hamiltonian systems, it has been rigorously proved that,  
otherwise, 
 all the invariant tori 
could be destroyed and only  discontinuous  Aubry-Mather invariant sets survive, see e.g. \cite{He1}.
\end{remark}

The infinite set of the eigenvalues of the diagonal operator $ D  $, 
$$
 - (\om \cdot 	\ell)^2 + |j|^2 + m   \, , \quad  \  \ell \in \Z^\nu \, , \ j \in \Z^d \, ,
$$
accumulate to zero (small divisors) and therefore the matrix $ T  $ in \eqref{DomT}, 
which represents in Fourier space the multiplication operator for the function $ b(\vphi, x ) $,  
is a ``singular" perturbation of $ D $. As a consequence it is not obvious 
at all that the self-adjoint operator $ D + \e T  $ has still pure point spectrum
with a basis of eigenfunctions with exponential/polynomial decay, 
 for $ \e $ small, 
for most  values of the frequency vector $ \om \in \R^\nu $.
This is the main problem addressed in Anderson localization theory. 


\smallskip

Actually, for the convergence of a Nash-Moser scheme in applications to
KAM for PDEs (cf. \eqref{AINV0SOB} and Remark \ref{WTE}),  it 
is sufficient to prove the following:
\begin{enumerate}
\item
for any $ N $ large, the finite dimensional restrictions of the operator in \eqref{model-B-int}, 
\be\label{QP-lin:OP-int}
{\cal L}_N := 
{\it \Pi}_N \big( (\om \cdot \pa_{\vphi})^2  - \Delta + m  + \e b(\vphi, x)  \big)_{| {\cal H}_N}  \, ,  
\ee
are invertible for most values of the parameters, 
where $ {\it \Pi}_N $ denotes the projection on the finite dimensional subspace 
\be\label{calHN-primo}
 {\cal H}_N := \Big\{ h (\vphi, x) 
 = \sum_{|(\ell, j)| \leq N} h_{\ell, j} e^{\ii (\ell \cdot \vphi + j \cdot x)}   
 \Big\} \, .
\ee 
\item 
 The inverse  ${\cal L}_N^{-1} $ satisfies, for some $ \mu > 0 $, $ s_1 > 0 $,  
 tame  estimates  as 
\be\label{LN-1tame}
\| {\cal L}_N^{-1} h \|_s \leq C(s) N^\mu \big( \| h \|_{s} + \| b \|_s \| h \|_{s_1} \big) \, , \
\forall h \in {\cal H}_N \,  , \  \forall s \geq s_1 \, , 
\ee 
where $ \| \ \|_s  $ denotes the Sobolev norm in \eqref{Sobo:sp1}. 
\end{enumerate}

\begin{remark}
Note that $ \mu > 0 $ represents in  \eqref{LN-1tame} 
a ``loss of derivatives" due to the small divisors. 
Since the multiplication operator $ h \mapsto b h $
satisfies a tame estimate like \eqref{LN-1tame} with $ \mu = 0 $, the estimate
\eqref{LN-1tame} means that
 ${\cal L}_N^{-1}$ acts, on the Sobolev scale $ {\cal H}^s $,  somehow
as an unbounded differential operator of order $ \mu $.
\end{remark}

Also weaker tame estimates as 
\be\label{LN-1tame-s}
\| {\cal L}_N^{-1} h \|_s \leq C(s) N^{\tau'} \big( (N^{ \varsigma s} + \|b\|_s)
\| h \|_{s_1} +  N^{ \varsigma s_1} \| h \|_{s}  \big) \, , \
\forall h \in {\cal H}_N \,  , \  \forall s \geq s_1 \, , 
\ee 
for $ \varsigma < 1  $,  
are sufficient for the convergence of the Nash-Moser scheme. 
Note that such conditions are much weaker than \eqref{LN-1tame} 
because the loss of derivatives $ \tau' + \varsigma s $ 
increases with $ s $. 
Conditions like \eqref{LN-1tame-s} are essentially optimal
for the convergence, compare with Lojasiewicz-Zehnder \cite{LZ}.

\begin{remark}
It would be also sufficient 
to prove the existence of a right inverse of the operator 
$ [{\cal L}]_N^{2N} := {\it \Pi}_N \big( (\om \cdot \pa_{\vphi})^2  - \Delta + m  + \e b(\vphi, x)  \big)_{| {\cal H}_{2N}}  $
satisfying \eqref{LN-1tame} or \eqref{LN-1tame-s}, as we do in 
Section \ref{sec:rightinv}. We remind that a linear operator 
acting between finite dimensional vector spaces 
admits a right inverse if it  is surjective\index{Right inverse}, see Definition \ref{def:right-inv}. 
\end{remark}

In the time periodic setting, i.e.  $ \nu = 1 $, we establish the stronger 
tame estimates \eqref{LN-1tame}
whereas, in the time quasi-periodic setting, i.e.  $ \nu \geq 2 $, we prove the weaker
tame estimates \eqref{LN-1tame-s}. 
Actually, 
 in the present monograph, as in \cite{BB12}-\cite{BBo10}, we shall prove that 
the approximate inverse satisfies tame estimates of the weaker 
form \eqref{LN-1tame-s}.

In order to achieve \eqref{LN-1tame} or \eqref{LN-1tame-s} there are two main steps:
\begin{enumerate}
\item ($L^2$-{\it estimates}) \label{item_L2-e}
Impose lower bounds for the eigenvalues of the self-adjoint 
operator \eqref{QP-lin:OP-int}, for most values of the parameters. 
These ``first order Melnikov" non-resonance conditions 
provide
 estimates of the inverse of $ {\cal L}_N^{-1} $ 
 in $ L^2 $-norm. 
\item ($ {\cal H}^s$-{\it estimates}) \label{item_Hs-e}
Prove the estimates \eqref{LN-1tame} or \eqref{LN-1tame-s}
 in high Sobolev norms.  In the language of Anderson localization theory,  
this amounts to proving
polynomially  fast off-diagonal decay estimates for the inverse matrix $ {\cal L}_N^{-1} $. 
\end{enumerate}

In the forced case when the frequency vector $ \om \in \R^\nu $ provides independent parameters,  item 
\ref{item_L2-e} is not a too difficult task, using results about the eigenvalues of 
self-adjoint matrices depending on parameters, as Lemmata \ref{lem:posE} and
\ref{variatione}. On the other hand, in the autonomous case, for the difficulties 
discussed in Remark \ref{rem:ANA}, this is a more delicate task (that we address in this monograph).  

In the sequel we concentrate 
on the  analysis  of item \ref{item_Hs-e}.  
An essential ingredient is the decomposition into   ``singular" and 
``regular" sites, for some $ \rho > 0 $,   
\be\label{sing:cone-int}
\begin{aligned}
S & := \Big\{ (\ell, j) \in \Z^\nu \times \Z^d  \quad {\rm such \ that } \quad 
| - ( \om \cdot \ell )^2 + |j|^2 + m | < \rho  \Big\} \\
R & := \Big\{ (\ell, j) \in \Z^\nu \times \Z^d  \quad {\rm such \ that } \quad 
| - ( \om \cdot \ell )^2 + |j|^2 + m | \geq \rho  \Big\} \, .
\end{aligned}
\ee
It is clear, indeed, that, in order to achieve \eqref{LN-1tame}, conditions that 
limit the quantity of singular sites have to be fulfilled,  
otherwise the  inverse operator $ {\cal L}_N^{-1} $
would be  ``too big"  in any sense and
\eqref{LN-1tame-s} would be violated.

\begin{remark}
If in \eqref{model-B-int} the 
constant $ m $  is replaced by a (not small) multiplicative potential 
$ V(x) $ (this is the case considered in the present monograph), 
it is natural 
to define the singular sites as
$ | - ( \om \cdot \ell )^2 + |j|^2 + m | <  \Theta $
for some $ \Theta $ large depending on $ V(x) $, where
$ m $ is the mean value of $ V(x) $.
\end{remark}
 
We first consider the easier case $ \nu = 1 $ (time-periodic solutions). 

\subsection{Time periodic case}\label{sec:MULTI-P}

 Existence of time periodic solutions for NLW equations on $ \T^d $ have been 
first obtained in \cite{B-Gafa} in an analytic setting. In the exposition below we 
 follow \cite{BBARMA}, which works in the context of Sobolev spaces.  

The following ``separation properties" 
of the singular sites \eqref{sing:cone-int}
are sufficient  for proving tame estimates as \eqref{LN-1tame}: the singular sites 
$ S $ in
the box $[ - N, N]^{1+d} $ are partitioned into disjoint clusters  $ \Omega_\alpha $, 
\be\label{S-deco}
S \cap [-N,N]^{1+d} = \bigcup_{\a  } \Omega_\alpha \, , 
\ee
satisfying 
\begin{itemize}
\item 
{\bf (H1)} {\bf (Dyadic)} $ M_\a \leq 2 m_\a $, $ \forall \a $, where 
$ M_\a := \max_{(\ell,j) \in \Omega_\a} |(\ell,j)|$ and 
$ m_\a := \min_{(\ell,j) \in \Omega_\a} |(\ell,j)| $;  
\item
{\bf (H2)} {\bf (Separation at infinity)} There is $ \d = \d(d) > 0 $ (independent of $ N $) 
such that 
\be\label{pr:H2}
{\rm d} (\Omega_\a, \Omega_\b) := \min_{(\ell,j) \in \Om_\a, 
(\ell',j') \in \Om_\s} |(\ell,j) - (\ell',j')| \geq ( M_\a + M_\s)^{\d} \, , \quad \forall \a \neq \s \, . 
\ee
\end{itemize}
Note that in (H2)  the clusters $ \Omega_\a  $ of 
singular sites are ``separated at infinity", namely the distance between 
distinct  clusters 
increases when the Fourier indices tend to infinity. 
A partition of the singular sites as \eqref{S-deco}, satisfying (H1)-(H2),  has been 
proved in \cite{BBARMA} assuming that $ \om^2 $ is Diophantine\index{Diophantine vector}. 

\begin{remark}\label{rem:SOBO}
We require that the function $ b(\vphi, x) $ in \eqref{model-B-int}
has the same (high) regularity in $ (\vphi, x) $, i.e. $ \| b \|_s < \infty $ for some
large Sobolev index $ s $,  because
the clusters $ \Omega_\a $  in \eqref{pr:H2} are separated  in 
time-space Fourier indices only. 
 This implies, in  KAM applications,
 that  the 
 solutions that we obtain will have the same high Sobolev regularity in time and space. 
\end{remark}

We have to solve the linear system
\be\label{LNhb}
{\cal L}_N g = h \, , \quad g, h \in {\cal H}_{N} \, . 
\ee
Given a subset $ \Omega $ of $ [-N, N]^{1+d} \subset \Z \times \Z^d $ we denote by 
$ {\cal H}_{\Omega} $ the vector space spanned by $ \{ e^{\ii (\ell \vphi + j \cdot x )}, 
(\ell, j) \in \Omega \} $ and by $ \Pi_{\Omega}  $ the corresponding orthogonal projector. 
With this notation the subspace $ {\cal H}_N $ in \eqref{calHN-primo} 
coincides with $ {\cal H}_{[-N,N]^{1+d}} $. 
Given a linear operator $ L $ of $ {\cal H}_N $ and 
another subset $ \Omega' \subset [-N, N]^{1+d} $ we denote
$ L_{\Omega'}^{\Omega} := \Pi_{{\cal H}_{\Omega'}} L_{| {\cal H}_{\Omega}} $.
For simplicity we also set $  L_{\Omega'}^{\Omega} := 
\Pi_{{\cal H}_{\Omega'}} ({\cal L}_N)_{| {\cal H}_{\Omega}}  $. 

Let us  for simplicity of notation  still denote by $ S, R $ the sets of the 
singular and regular sites \eqref{sing:cone-int} intersected with the box   $[-N,N]^{1+d}  $.
According to the splitting  $ [-N,N]^{1+d} = S \cup R $, 
  we introduce the orthogonal decomposition 
$$
{\cal H}_N  = {\cal H}_S \oplus {\cal H}_R \, .
$$
Writing the unique decomposition $ g = g_S + g_R $,
$ g_R \in {\cal H}_R $, $ g_S \in {\cal H}_S $, and similarly for $ h $, 
the linear system  \eqref{LNhb} then amounts to 
$$
\begin{cases}
L_R^R g_R + L_R^S g_S = h_R  \\
L_S^R g_R + L_S^S g_S = h_S \, . 
\end{cases}
$$
Note that, by \eqref{DomT}, the coupling terms
$ L_R^S =  
\e T_R^S $, 
$ L^R_S =  
\e T^R_S  $,  
have polynomial off-diagonal decay. 

By standard perturbative arguments, the operator $ L_R^R $, which is 
the restriction of $ {\cal L}_N $ to the  regular sites,  is invertible
for $ \e \rho^{-1} \ll 1 $, and therefore, 
solving the first equation in $ g_R $, and inserting the result in the second one, we are reduced to
solve 
$$
\big( L_S^S - L_S^R (L_R^R)^{-1} L^S_R  \big) g_S = h_S -  
L_S^R (L_R^R)^{-1}  h_R \, .  
$$
Thus the main task is now to invert the self-adjoint matrix 
$$
U := L_S^S - L_S^R (L_R^R)^{-1} L^S_R \, , \quad U : {\cal H}_S \to {\cal H}_S \, .  
$$
This  reduction procedure is sometimes referred as a resolvent identity.  

According to \eqref{S-deco} we have the orthogonal decomposition 
$ {\cal H}_S = \oplus_{\a} {\cal H}_{\Omega_\a} $
which induces a block decomposition for the operator
$$
U = ( U_{\Omega_\s}^{\Omega_\a} )_{\a, \s } \, , \quad 
  U_{\Omega_\s}^{\Omega_\a} := 
  L_{\Omega_\s }^{\Omega_\a} - 
  L_{\Omega_\s}^R (L_R^R)^{-1} L^{\Omega_\a}_R \, . 
$$
Then we decompose $ U $ in block-diagonal and off-diagonal parts: 
\be\label{defDR}
U = {\cal D} + {\cal R } \, , \quad {\cal D} := {\rm Diag}_{\a } 
U_{\Omega_\a}^{\Omega_\a} \, , \quad 
{\cal R } := ( U_{\Omega_\s}^{\Omega_\a} )_{\a \neq \s} \, .
\ee
Since the matrix $ T $ has off-diagonal decay 
(see \eqref{DomT} and recall that the function $ b (\vphi, x ) $ is smooth enough) and 
the matrices with off-diagonal decay form an algebra (with interpolation estimates) it 
 easily results  an off-diagonal decay of the matrix $ {\cal R } $ like 
\be\label{off-diag-R}
\| U_{\Omega_\s}^{\Omega_\a} \|_{{\cal L} (L^2)} \leq \frac{\e \, C(s) }{ {\rm d}(\Om_\a, \Om_\s)^{s- \frac{1+d}{2}}} \, , \quad \a \neq \s \, .
\ee
In the simplest case that the $ \om $ are independent parameters (forced case), 
it is not too difficult to impose that 
each self-adjoint  operator $ U_{\Omega_\a}^{\Omega_\a} $
is invertible for most parameters and that
$$
\| (U_{\Omega_\a}^{\Omega_\a})^{-1}\|_{{\cal L} (L^2)} \leq M_\a^{\tau} \, , \quad \forall \a  \, , 
$$
for some $\tau$ large enough. Since, by (H1),  each cluster $\Omega_\a $ is {\it dyadic}, this $ L^2 $-estimate implies also a $ {\cal H}^s $-Sobolev bound: for all $ h \in {\cal H}_{\Omega_\a} $, 
$$
\begin{aligned}
\| (U_{\Omega_\a}^{\Omega_\a})^{-1} h \|_{s} \leq M_\a^s 
\| (U_{\Omega_\a}^{\Omega_\a})^{-1} h  \|_{L^2} \leq 
M_\a^s M_\a^\tau \| h \|_{L^2} 
& \leq \frac{M_\a^s}{m_\a^s} M_\a^\tau \| h \|_s \\
& \leq 2^s N^\tau \| h \|_s 
\end{aligned}
$$
by (H1) and  $ M_\a \leq N $. 
As a consequence
the whole operator $ {\cal D}$ defined in \eqref{defDR}
is invertible with a Sobolev  estimate 
$$
\| {\cal D}^{-1} h \|_s \leq C(s) N^\tau \| h \|_s \, , \quad \forall h \in {\cal H}_N \, .  
$$
Finally, using the off-diagonal decay \eqref{off-diag-R} and 
the ``separation at infinity" property \eqref{pr:H2}, the operator
$ {\cal D}^{-1}{\cal R} $ is bounded in $ L^2 $, and 
it is easy to reproduce a Neumann-series argument to
prove the invertibility of  $ U = {\cal D} ( I + {\cal D}^{-1}{\cal R}) $ 
with  tame estimates for the inverse $ U^{-1}$ in Sobolev norms, implying 
\eqref{LN-1tame}.

\subsection{Quasi-periodic case}\label{sec:MULTI-QP}

In the quasi-periodic setting,  i.e. $ \nu \geq 2 $, 
the proof that  the operator $ {\cal L}_N $ defined in \eqref{QP-lin:OP-int} is invertible 
and its inverse satisfies the tame estimates  \eqref{LN-1tame-s}  is more difficult. 
Indeed the `separation at infinity" property  (H2) never holds in the quasi-periodic case, 
not even  for finite dimensional systems. For example,  the operator 
$ \om \cdot \pa_\vphi $, is represented in the  Fourier basis 
as the diagonal matrix $ {\rm Diag}_{\ell \in \Z^\nu } ( \ii \omega \cdot \ell  ) $.
If the frequency vector $ \om  \in \R^\nu $  is  Diophantine\index{Diophantine vector}, then  the singular sites  $ \ell \in \Z^\nu  $ such that 
$$ 
| \omega \cdot \ell | \leq \rho 
$$ 
are ``uniformly distributed" in a neighborhood of the hyperplane
$ \omega \cdot \ell = 0 $, 
with nearby indices at distance $ O( \rho^{- \a}) $ for some $ \a >  0 $. 
Therefore,  
unlike  in the time periodic case,  the decomposition \eqref{sing:cone-int} 
into singular and regular sites of the unperturbed linear operator is not sufficient, 
and a finer analysis has to be performed. 

\smallskip
In the sequel we follow the exposition of 
\cite{BB12}. Let 
 $$ 
 A = D + \e T  \, , \quad 
 D := {\rm Diag}_{(\ell, j) \in \Z^\nu \times \Z^d} ( - (\om \cdot 	\ell)^2 + |j|^2 + m ) \, , 
 \quad T := (\widehat b_ {\ell - \ell', j - j'}) \, , 
 $$ 
 be the self-adjoint matrix in \eqref{DomT} that represents the second order operator \eqref{model-B-int}.
We suppose (forced case) that $ \omega $ and $ \e $ are unrelated and 
 that  $ \omega $ is constrained to a fixed Diophantine direction\index{Diophantine vector} 
\be\label{omega-su-linea}
\omega = \lambda \bar \omega \, , \quad \lambda \in \Big[ \frac12, \frac32 \Big] \, , 
\quad | \bar \om \cdot \ell | \geq \frac{\g_0}{|\ell|^{\tau_0}} \, ,
\ \forall \ell \in \Z^\nu \setminus \{0\} \, . 
\ee 
A way to deal with  
the lack of ``separation properties
at infinity" of the singular sites,  is to implement an
 inductive procedure. Introduced a rapidly increasing sequence $(N_n)$ of scales, 
defined by 
\be\label{def:Nk-multi-int}
N_n = \big[N_0^{\chi^{n}}\big] \, , \quad n \geq 0  \, , 
\ee
the aim is to obtain, for most parameters,  
off-diagonal
decay estimates for the inverses $ A_{N_n}^{-1} $ of the restrictions 
\be\label{def:ANn-i}
A_{N_n} := {\it \Pi}_{N_n} A_{| {\cal H}_{N_n}} \, . 
\ee
In \eqref{def:Nk-multi-int} 
the constants $ N_0 $ and $  \chi $ are chosen large enough.
The estimates at step $n$ (for matrices of size $N_n$) 
partially rely  on information about the invertibility and the 
off-diagonal decay of ``most" inverses $ A_{N_{n-1}, \ell_0, j_0}^{-1}  $
of  submatrices 
$$ 
A_{N_{n-1}, \ell_0, j_0} 
:= A_{| \ell - \ell_0| \leq N_{n-1}, |j - j_0| \leq N_{n-1}}  
$$ of size $ N_{n-1} $.   
This program  suggests the name ``multiscale analysis". 

In order to give a precise statement 
we first introduce decay norms.
Given a matrix
$ M = (M^{i'}_i)_{i' \in B, i \in C} $,   
where  $ B, C $ are  subsets of $ \Z^{\nu + d} $, we define its  $ s $-norm
$$
\norma M \norma_s^2 :=  \sum_{n \in \Z^{\nu + d}} [M(n)]^2 \langle n \rangle^{2s}  
\quad {\rm where}  \quad  \langle n \rangle := \max (1,|n|) \, , 
$$
and 
$$
[M(n)] := \begin{cases}
\max_{i-i'=n, i \in C, i' \in B} |M^{i'}_i|  \ \ \, \, \quad   {\rm if}  \ \   n \in C- B \\
0 \quad \ \qquad \qquad \qquad \qquad  {\rm if}  \ \  n \notin C - B \, . 
\end{cases}
$$

\begin{remark}
The  norm $ |T|_s $ of the matrix which represents the 
multiplication operator for the function $ b(\vphi, x ) $ is equal to 
$ | T |_s = \| b \|_s $. 
Product of  matrices (when it makes sense) with finite $ s $-norm  
satisfy algebra and interpolation inequalities, see Appendix \ref{sec:off}. 
\end{remark}

We now outline how to 
prove that the finite dimensional matrices $  A_{N_n} $ in \eqref{def:ANn-i}
are invertible for ``most"  
parameters $ \l \in [1/2,3/2]  $ and satisfy, for all $ n \geq 0 $, 
\be\label{Ln1}
\nors{ A_{N_n}^{-1}} \leq C(s) N_n^{\tau'} \big( N_n^{\varsigma s}   + \| b \|_s \big) 
 \, , \   \ \varsigma 
\in (0,1/2) \, ,  \ \t' >  0 \, , \  \forall s >  s_0 \, .
\ee
Such estimates imply  the off-diagonal decay of  the entries of the inverse matrix
$$
| ( A_{N_n}^{-1})_{i'}^{i}| \leq C(s) 
\frac{N_n^{\tau'} (N_n^{\varsigma s} + \| b \|_s )}{ \langle 
i - i' \rangle^s} \, ,  \quad |i|, |i'| \leq N_n \, ,
$$
and Sobolev tame  estimates as \eqref{LN-1tame-s} (by Lemma \ref{sobonorm}), 
assuming that $ \| b \|_{s_1} \leq C $. 

\smallskip

In order to prove \eqref{Ln1} at the initial scale $ N_0 $ we impose that, for most parameters, the eigenvalues of the diagonal matrix $ D $ satisfy 
$$  
|- (\om \cdot 	\ell)^2 + |j|^2 + m| \geq N_0^{-\tau }  \, , \quad \forall 
|(\ell, j)| \leq N_0  \, . 
$$ 
Then, for $ \e $ small enough, we deduce, by a direct Neumann series argument,
the invertibility of $ A_{N_0} = D_{N_0} + \e T_{N_0} $ and the decay of the inverse
$ A_{N_0}^{-1} $.  

In order to proceed at the higher scales, we use 
a multiscale analysis. 

\smallskip

\noindent
{\bf $ L^2 $-bounds.}
The first step is to  show that,
for ``most" parameters, the eigenvalues of 
$ A_{N_n}  $ 
are in modulus bounded from below by  $ O(N_n^{-\tau}) $ and so the  $ L^2 $-norm of the inverse  satisfies 
\be\label{Ln-1}
\| A_{N_n}^{-1} \|_0= O(N_n^\tau) \qquad {\rm where} \qquad 
\| \ \|_0 := \| \ \|_{{\cal L} (L^2)}.
\ee 
The proof is based on  eigenvalue variation arguments. Recalling 
\eqref{omega-su-linea}, 
dividing $ A_{N_n} $ by $ \l^2 $, 
 and setting $ \xi := 1 / \l^2 $, 
we observe that the derivative with respect to $ \xi $ satisfies 
\be\label{der-xi-pos}
\begin{aligned}
\partial_\xi  (\xi A_{N_n})  & =  {\rm Diag}_{|(\ell,j)| \leq N_n } (| j|^2 + m )
+ 
O \big( \e \| T \|_{0} + \e \| \partial_\l T \|_{0}\big)  \geq c > 0  \, , 
\end{aligned}
\ee
for $ \e $ small, i.e. it is positive definite. 
So, the eigenvalues  $ \mu_{\ell,j} (\xi ) $ 
of the self-adjoint matrix $ \xi A_{N_n} $ satisfy (see Lemma \ref{lem:posE})
$$
\partial_\xi \mu_{\ell,j} (\xi ) \geq c > 0   \, ,  \quad \forall |(\ell,j)| \leq N_n \, ,
$$
which  implies  
\eqref{Ln-1}  except in  a set of $ \l $'s of 
measure $ O(   N_n^{-\t+ d + \nu}) $. 

\begin{remark}
Monotonicity arguments for proving lower bounds for the moduli of
the eigenvalues of (huge) self-adjoint matrices  have been also used in
\cite{EK}, \cite{BBo10}, \cite{BB12}. 
 Note that the eigenvalues could be degenerate for some values of the parameters. 
This approach to verify ``large deviation estimates" is very robust  
and  the measure estimates that we perform at each step
of the iteration are not inductive, as those in \cite{B5}. 
\end{remark}

\noindent
{\bf  Multiscale Step.} 
The bounds \eqref{Ln1} for the decay norms of $ A_{N_n}^{-1} $ 
follow by an inductive application of the
 multiscale step proposition \ref{propinvA}, that we now describe. 
 A matrix $ A \in {\cal M}_E^E $, $ E \subset \Z^{\nu + d } $, with 
$ {\rm diam}(E) \leq  4 N $ is  called {\sc $ N $-good} if it is invertible and 
$$ 
 \nors{A^{-1}} \leq N^{\tau'+\varsigma s} \, , \quad \forall s \in [s_0, s_1] \, ,
$$
for some $ s_1 := s_1 (d,\nu) $ large.  Otherwise we say that $ A $ is $ N $-bad.

The  aim of the multiscale step is to deduce that 
a matrix  $ A \in \matr_E^E $ with 
\be\label{newN'-N-i}
{\rm diam}(E) \leq N' := N^\chi \quad {\rm with} \quad \chi \gg 1 \, , 
\ee 
is $ N' $-good, knowing 
\begin{itemize}
\item 
{\bf (H1)} {\bf (Off-diagonal decay)}
$\norsone{A- {\rm Diag}(A)} \leq \Upsilon $ where $ {\rm Diag}(A) := 
( \d_{i,i'} A_i^{i'})_{i, i' \in E} $;
\item 
{\bf (H2)}  {\bf ($ L^2 $-bound)} $ \| A^{-1} \|_0 \leq (N')^{\tau} $ where we set
$ \| \ \|_0 := \| \ \|_{ {\cal L}(L^2)} $;  
\end{itemize}
\noindent and suitable 
assumptions concerning the   $ N $-dimensional submatrices 
along the diagonal of $ A $. 

\begin{definition}  {\bf ($(A,N)$-bad sites)}
An index $ i \in E   $ is called
\begin{enumerate}
\item 
$(A,N)$-{\sc regular} if there is $ F \subset E$ containing $i$ such that
${\rm diam}(F) \leq 4N$,  ${\rm d}(i, E\backslash F) \geq N/2$ and
$A_F^F$ is $N$-good;
\item 
$(A,N)$-{\sc bad} if it is singular (i.e. $ |A_i^i| < \rho $) and $(A,N)$-{\sc regular}. 
\end{enumerate}
\end{definition}

We suppose the following ``mild" separation properties for the $ (A,N) $-bad sites. 
\begin{itemize}
\item 
{\bf (H3)} {\bf (Separation properties)}
There is a partition of the  $(A,N)$-bad sites $ B = \cup_{\alpha} \Omega_\alpha $ with
\be\label{sepabad-int}
{\rm diam}(\Omega_\alpha) \leq N^{C_1} \, , \quad {\rm d}(\Omega_\alpha , \Omega_\beta) \geq N^2 \ , \ \forall \alpha \neq \beta \, ,
\ee
for some $ C_1 := C_1(d,\nu) \geq 2 $.
\end{itemize}
The multiscale step proposition \ref{propinvA} deduces
that $ A $ is $ N' $-good and \eqref{Ln1} holds at the new scale $ N' $,  by (H1)-(H2)-(H3),
with suitable relations between the constants  $ \chi  $, $ C_1 $, $ \varsigma $, $ s_1 $, 
$ \tau $. 
Roughly, the main conditions on the exponents are 
$$ 
C_1 < \varsigma \chi \qquad  {\rm and } \qquad 
  s_1 \gg  \chi \, \t \, .
 $$
  The first  means that the size $ N^{C_1} $
of any bad cluster $ \O_\a $  is small with respect to the size 
$ N' := N^\chi $ of the matrix $ A $.
The second means that the Sobolev regularity 
$ s_1 $ is large enough to ``separate" the resonance effects of two nearby bad clusters $ \O_\a $, $ \O_\b $.

\smallskip

\noindent
{\bf Separation properties.} We apply the 
multiscale step Proposition \ref{propinvA} to the matrix $ A_{N_{n+1}} $.
The key property to verify is (H3). 
A first key ingredient is the following co-variance property:
consider the  family of infinite dimensional matrices 
$$
A(\theta) := D(\teta) +  \e T  \, , \quad 
D (\teta) := {\rm diag}_{(\ell,j) \in \Z^\nu \times \Z^d} 
\Big( - (\om \cdot \ell + \theta )^2 + |j|^2 + m \Big) \, ,  
$$
and its $ (2N+1)^{\nu + d } $-restrictions
$ A_{N,\ell_0, j_0}(\theta) := A_{| \ell - \ell_0| \leq N, |j - j_0| \leq N}(\theta) $
centered at any $ (\ell_0, j_0 ) \in \Z^\nu \times \Z^d $. 
Since the matrix $ T $ in \eqref{DomT} is T\"oplitz we have 
the co-variance property
\be\label{shifted-i}
A_{N, \ell_0,  j_0} (\theta) = A_{N, j_0,0} ( \teta +  \om \cdot \ell_0 ) \, . 
\ee
In order to deduce (H3), it is sufficient to prove  the ``separation properties" \eqref{sepabad-int} for 
the  {\sc $ N_n $-bad} sites of 
$ A $, namely the indices $ ( \ell_0, j_0) $ 
which are singular 
\be\label{singul-i}
| - (\om \cdot \ell_0)^2 + |j_0|^2 + m | \leq \rho \, , 
\ee  
and for which there exists a site $ (\ell,j) $, with $ |(\ell,j) - (\ell_0,j_0)| \leq N $, such that 
$ A_{N_n, \ell, j} $ is $ N_n $-bad. 
Such separation properties are obtained for all the parameters $ \lambda $ 
which are $ N_n $-good, namely such that
\be
\begin{aligned}
\forall \, j_0 \in \Z^d \, , \quad
& B_{N_n} (j_0;  \l)  :=  
\Big\{ \teta \in \R \, : \,   A_{N_n, 0,j_0} (\theta)  \ {\rm is \ } N_n-{\rm bad}  \Big\} 
\label{BNcomponents-i}   \\
& \subset \bigcup_{q = 1, \ldots, N_n^{C(d,\nu)}} I_q  
\ \  {\rm where} \ I_q \ {\rm are \ intervals \ with } \ | I_q| \leq N_n^{-\t}  \, .  
\end{aligned}
\ee
We first use the covariance property \eqref{shifted-i}
and the ``complexity" information  \eqref{BNcomponents-i}  
to bound the number of ``bad" time-Fourier components. Indeed
$$ 
A_{N_n, \ell_0,j_0}  \ {\rm is \ } N_n\text{-}{\rm bad  }  \
 \Longleftrightarrow \  
A_{N_n, 0,j_0}(\om \cdot \ell_0)  \ {\rm is \ } N_n\text{-}{\rm bad  } 
\  \Longleftrightarrow  \ \om \cdot \ell_0 \in B_{N_n}(j_0;\l)  \, . 
$$
Then, using  that $ \om $ is Diophantine\index{Diophantine vector},  the complexity  bound \eqref{BNcomponents-i} implies that, 
for each fixed $ j_0 $, there are at most $ O( N_n^{C(d,\nu)}) $ sites 
$ (\ell_0,j_0) $ in the larger box $ | \ell_0| \leq N_{n+1} $, 
which are $ N_n $-bad.

Next,  we prove that a $ N_n^2 $-``chain" of singular sites, 
i.e. a sequence of distinct integer vectors 
$ (\ell_1, j_1),  \ldots , (\ell_L, j_L) $ satisfying 
 \eqref{singul-i} and 
$$
| \ell_{q+1}-  \ell_q |+ | j_{q+1}-  j_q | \leq N_n^2 \, , \quad
 q = 1, \ldots, L \, ,
 $$ 
which are also  $ N_n $-bad, has a ``length"  $ L $ bounded by 
$ L \leq N_n^{C_1 (d, \nu)} $.  
This implies a partition of the $( A_{N_{n+1}}, N_n) $-bad sites 
as in \eqref{sepabad-int} at order $ N_n $. 
In this step we require that  
$ \om \in \R^\nu $ satisfies 
the quadratic Diophantine condition\index{Quadratic Diophantine condition}
\be\label{irra-slo}
\Big|n+ \sum_{1 \leq i\leq j \leq \nu} { \om}_i { \om}_j p_{ij} \Big|  \geq \frac{\g_2}{ | p |^{\t_2}} \, , 
\ \  \forall p := ( p_{ij}) \in  \Z^{\frac{\nu ( \nu + 1 )}{2}} \, , \  \forall n \in \Z \, , 
\, (p, n ) \neq (0,0) \, , 
\ee
for some positive $ \gamma_2, \tau_2 $. 

\begin{remark}
The singular sites \eqref{singul-i}
are integer vectors close to  a ``cone" and
 \eqref{irra-slo} can be seen as an irrationality condition
on its slopes. 
 For  NLS equations,  \eqref{irra-slo} is not required, because the 
 corresponding singular sites $ (\ell, j) $ satisfy 
 $ |  \om \cdot \ell  + |j|^2   | \leq C $, i.e. they are close to a paraboloid. 
 We refer to  Berti-Maspero \cite{BMa} to avoid the use of the quadratic Diophantine 
 condition
 \eqref{irra-slo} for  NLW equations. 
\end{remark}

\noindent
{\bf Measure and ``complexity" estimates.}
In order to conclude the inductive proof  we have to verify  that ``most" parameters 
$ \l $ are $ N_n $-good, according to 
\eqref{BNcomponents-i}. 
We prove first  that,  except a set 
of measure $ O( N_n^{-1}) $, all parameters $ \l \in [1/2, 3/2 ] $
 are $ N_n $-good in a $ L^2$-sense, namely
\be 
\begin{aligned}
\forall \, j_0 \in \Z^d \, , \quad
& B_{N_n}^0(j_0; \l)  :=  
\Big\{ \teta \in \R \, : \,   \| A_{N_n, 0, j_0}^{-1} (\theta)\|_0 > N_n^{\t}  \Big\} \label{tetabadweak0}   \\
& \subset  \bigcup_{q = 1, \ldots, N_n^{C(d, \nu)}} I_q  
\ \  {\rm where} \ I_q \ {\rm are \  intervals \ with } \ | I_q| \leq N_n^{-\t}  \, .  
\end{aligned}
\ee
The proof is again based on  eigenvalue variation\index{Eigenvalues of self-adjoint matrices} arguments as in 
\eqref{der-xi-pos}, using
that $ -\Delta + m  $ is positive definite.

Finally, the multiscale Proposition step \ref{propinvA}, and the fact that 
the separation properties  of the $ N_n $-bad sites of $ A(\teta) $ hold uniformly in
 $ \theta \in \R $, 
imply  inductively that 
the parameters $ \l $ which are  $ N_n $-good in  $ L^2$-sense, 
 are actually $ N_n $-good, i.e. \eqref{BNcomponents-i} holds, 
 concluding the inductive argument.

\subsection{The multiscale analysis of Chapter  \ref{sec:multiscale} }
\label{subsec:Ch4}

In this monograph we consider {\it autonomous} nonlinear wave equations
\eqref{NLW1}, whose small amplitude quasi-periodic solutions of frequency vector $\om$ 
correspond to solutions
of \eqref{NLW2QP}. 
The  quasi-periodic linear operator 
\be\label{ome-LO}
(\om \cdot \pa_{\vphi})^2  - \Delta + V(x)  + \e^2 
(\pa_u g)(x, u(\vphi,x))  \, , \quad  \vphi \in \T^\nu  \, , \
x \in \T^d \, ,  
\ee
is obtained by linearizing \eqref{NLW2QP} at an approximate 
solution $u$. The analysis of this operator, 
acting in  a subspace orthogonal  to the  unperturbed linear 
`tangential" solutions (see \eqref{sol:uel})
$$ 
\Big\{ \sum_{j \in {\mathbb S}}  
\a_j \cos ( \vphi_j ) \Psi_j (x)  \, , \ \a_j \in \R \Big\} \subset
{\rm Ker} \big( (\bar \mu  \cdot \pa_{\vphi})^2  - \Delta + V(x) \big) \, , 
$$
turns out to be  much more difficult  than in the forced case; it requires considerable
extra-work with respect to the previous section. 
Major difficulties are the following:
\begin{itemize}
\item
(i) 
the frequency vector 
$ \omega $ and the small parameter $ \e $ are linked in \eqref{ome-LO}, see the discussion 
in Remark \ref{rem:ANA}. In particular 
$ \omega $ has to be $ \e^2 $-close to the unperturbed frequency vector 
$ \bar \mu = ( \mu_j )_{j \in {\mathbb S}} $ in \eqref{unp-tangential}. More precisely 
$ \omega $ varies approximately according to 
the frequency-to-action map \eqref{act-to-fre}, and 
thus
it has to belong to the region of admissible frequencies\index{Admissible frequencies} \eqref{def:admissible}.
Moreover  the multiplicative function 
$$ 
(\pa_u g)(x, u(\vphi,x)) = 3 a(x) (u(\vphi, x))^2 + \ldots 
$$ 
in \eqref{ome-LO}
(recall the form of the nonlinearity \eqref{nonlinearity:gep})
depends itself on $ \omega $, 
via the approximate quasi-periodic solution $ u(\vphi, x ) $,  
and, as the tangential frequency vector $ \omega $
changes, also the normal frequencies undergo a significant modification. 
At least for finitely many modes, the shift of the normal 
frequencies  due to the effect of $ 3 a(x) (u(\vphi, x))^2 $ can be 
approximately described 
in terms of the Birkhoff matrix $ \Bb $  in \eqref{def:AB} as
$  \mu_j +  \e^2  (\Bb \xi)_j $ (see \eqref{Omega-normali-0}). Because of all these 
constraints,  
positivity properties like  \eqref{der-xi-pos}
in general fail, see \eqref{derivata lam} and Remark \ref{rem:per-split}.
This implies 
difficulties for imposing non-resonance conditions and for proving suitable 
complexity bounds. 
\item
(ii) An additional difficulty in \eqref{ome-LO} is the presence
of the  multiplicative   potential $ V(x) $, which is not diagonal in the Fourier basis.
\end{itemize}
In this monograph we shall still be able to use 
positivity arguments for imposing non-resonance conditions and proving
complexity bounds along the multiscale analysis. This requires  
to write \eqref{NLW2}  as a first order Hamiltonian system and 
to perform a block-decomposition  of the corresponding first order 
quasi-periodic linear  operator acting in the normal subspace.
We refer to the next section
\ref{sec:ideas} for an explanation about this procedure (which is necessarily technical) 
and here we limit ourselves to describing the  quasi-periodic operators 
that we shall analyze in Chapter \ref{sec:multiscale} with  multiscale techniques.

\smallskip

The multiscale Proposition\index{Multiscale proposition} \ref{propmultiscale} of 
Chapter \ref{sec:multiscale} provides the existence of 
a  right   inverse of finite dimensional restrictions 
of  self-adjoint linear operators like $ {\cal L}_r $ in \eqref{fo2-p}
and $ {\cal L}_{r,\mu} $ in \eqref{fo1-p}, which have the following form:
 
\begin{enumerate}
\item  for $ \omega = (1+ \e^2 \l) \bar \om_\e $, 
$ \l \in \widetilde \Lambda \subset \Lambda $,  
\be\label{fo2-p}
\begin{aligned}
& {\cal L}_r  =  J \om \cdot \partial_\vphi + D_V + r (\e, \l, \vphi) \\
&   {\rm acting \ on } \ 
h \in (L^2 (\T^\es \times \T^d, \R))^2   \, , 
\end{aligned}
\ee
where 
\begin{itemize}
\item 
 $ D_V = \sqrt{- \Delta + V(x) } $;
\item   $ J $ is the symplectic matrix 
$ 
J = \begin{pmatrix}
0  &   {\rm Id}   \\
 - {\rm Id} & 0     \\   
\end{pmatrix}, 
$ acting in $(L^2 (\T^\es \times \T^d, \R))^2$ by $ J(h_1,h_2)=(h_2,-h_1)$;
\item  For each $\vphi \in \T^\es$, $r(\e,\l, \vphi)$ is a (self-adjoint) linear operator
of $(L^2 (\T^d , \R))^2$. As an operator acting in $(L^2 (\T^\es \times \T^d, \R))^2$, 
the family $(r(\e,\l, \vphi))_{\vphi \in \T^\es}$ maps $h$ to $g$ such that
$$
g(\vphi, \cdot )=r(\e,\l, \vphi) h( \vphi , \cdot) \, ;
$$
it is self-adjoint and the entries of its matrix representation (in the Fourier basis) are $2\times 2$ matrices 
and have a polynomial off-diagonal decay.
\end{itemize}
A quasi-periodic linear  operator of the form $ {\cal L}_r $  arises 
by writing the linear wave operator \eqref{ome-LO} 
as a first order system as in  \eqref{lin:NLS-bis}, \eqref{LeqNLW} (using
real variables instead of complex ones).  
We constrain the frequency vector  
$ \omega = (1+ \e^2 \l) \bar \om_\e $ to a straight line, because, in this way, we 
will be able to prove that, for ``most" values of $ \lambda $,  these 
self-adjoint operators have 
eigenvalues different from zero. 
\item For $ \omega = (1+ \e^2 \l) \bar \om_\e $, $ \l \in \widetilde \Lambda \subset \Lambda $, 
\be\label{fo1-p} 
\begin{aligned}
& {\cal L}_{r,\mu} =  J \om \cdot \partial_\vphi + D_V + \mu (\e, \l) 
{\cal J} \Pi_{\mathbb G} + r (\e, \l, \vphi)  \\
&  \quad \qquad {\rm acting \ on } \ h \in (L^2 (\T^\es \times \T^d, \R))^4 \, , 
\end{aligned}
\ee
where 
\begin{itemize}
\item 
$ {\cal J} h := J h J $ and the left/right action of $ J $ 
on $ \R^4 $ are defined in \eqref{leftJA}, \eqref{rightJA};
\item $ \mu (\e, \l) $ is a scalar;  
\item   $ \Pi_{\mathbb G} $ is the $L^2$-orthogonal projector 
{\it in the space variable} on the infinite
dimensional subspace 
$$
H_{\mathbb G}  :=  H_{ {\mathbb S} \cup {\mathbb F} }^\bot \, , \quad  
H_{ {\mathbb S} \cup {\mathbb F} }  := 
\Big\{ \sum_{j \in {\mathbb S} \cup {\mathbb F}} (Q_j, P_j) \Psi_j (x) \, , \ 
(Q_j, P_j)  \in \R^2 \Big\} \subset (L^2(\T^d))^2
$$
(the set $ {\mathbb G} \subset \N $ is defined in 
\eqref{taglio:pos-neg}) :  the action of $ \Pi_{\mathbb G} $ in $(L^2 (\T^\es \times \T^d , \R))^4$ maps $h=(h^{(1)},h^{(2)}) \in 
(L^2 (\T^\es \times \T^d , \R))^2 \times (L^2 (\T^\es \times \T^d , \R))^2$ to
$g=(g^{(1)},g^{(2)})$ such that
$$
g^{(i)} (\vphi , \cdot)=\Pi_{\mathbb G} h^{(i)} (\vphi , \cdot) \ , \qquad i=1,2 \, ;
$$
\item  $r(\e, \l, \vphi)  $, acting in $(L^2 (\T^\es \times \T^d , \R))^4$, is a self-adjoint operator; 
the entries of its matrix representation (in the Fourier basis) are $4\times 4$ matrices 
and have a polynomial off-diagonal decay. 
\end{itemize}

 Quasi-periodic linear  operators of the form $ {\cal L}_{r,\mu} $  will arise
 in Chapter \ref{sec:splitting} in the process of block diagonalizing
a quasi-periodic linear operator as $ {\cal L}_r $  with respect to the 
splitting $  H_{{\mathbb S} \cup {\mathbb F}} \oplus H_{{\mathbb G}} $
(low-high frequencies).
The scalar $ \mu  (\e, \l)  $ will be (approximately) one 
of the finitely many eigenvalues 
of the restriction of $ {\cal L}_r $ to the subspace $ H_{ {\mathbb F}}  $.
We explain in Section \ref{sec:ideas} why (and how) we perform this block-diagonalization. 
\end{enumerate}

Notice that 
\begin{itemize}
\item the operators $ {\cal L}_{r} $ and $  {\cal L}_{r,\mu} $
defined in \eqref{fo2-p}, \eqref{fo1-p} 
are defined for all $ \lambda $ belonging 
to  a subset 
$ \widetilde \Lambda \subset \Lambda $ (which may shrink during the Nash-Moser iteration);
\item
 $ {\cal L}_{r} $ and $  {\cal L}_{r,\mu} $ 
are  first order vector valued quasi-periodic operators, unlike 
\eqref{ome-LO} 
which is a second order scalar quasi-periodic operator (it acts in ``configuration space").

The operator \eqref{fo2-p} arises (essentially) 
by writing the linear wave operator 
\eqref{ome-LO} as a first order system as in  \eqref{lin:NLS-bis}, 
where 
the non-diagonal vector field is
$ 1 $-smoothing. Therefore it is natural to require that the operators 
$ r = r(\e, \l, \vphi ) $ in \eqref{fo2-p}, \eqref{fo1-p} satisfy 
the 
decay condition
\be\label{rs1+}
| r |_{+,s_1} := | D_m^{\frac12} r D_m^{\frac12} |_{s_1}  = O( \e^2 )  
\quad {\rm where} \quad D_m := \sqrt{- \Delta + m }  
\ee
for some $ m >  0 $.
\end{itemize}
A key property of the operators $ {\cal L}_r $ and $ {\cal L}_{r, \mu} $,
that we shall be able to verify at each step of the Nash-Moser iteration, 
is the monotonicity property with respect to the parameter $\lambda$,
in item \ref{assu:pos} of Definition \ref{definition:Xr}. 

We first introduce the following notation: given a family $ (A(\lambda))_{\lambda \in \widetilde \Lambda}  $ of linear self-adjoint operators,
\be\label{nota-frak}
{\mathfrak d}_\l A(\l) \leq - c \,  {\rm Id} 
\qquad \Longleftrightarrow  \qquad \frac{A(\l_2) - A(\l_1)}{\l_2 - \l_1 }  \leq - 
c \, {\rm Id} \, , 
\quad \forall 
\l_1 \neq  \l_2  \, , 
\ee
where
$ A \leq - c \, {\rm Id}  $ means as usual 
$ (Aw,w)_{L^2} \leq - c  \| w \|_{L^2}^2 $.   
\begin{itemize}
\item
For $ {\cal L}_r $ the condition of monotonicity takes the form 
\be\label{vare2}
{\mathfrak d}_\l 
\Big( \frac{D_V  
+ r   }{1+ \e^2 \l} \Big) \leq - 
c \, \e^2  {\rm Id} \, , \quad c > 0 \, .  
\ee
 The condition  for
$ {\cal L}_{r, \mu} $ is similar. 
Such property  allows to prove the measure estimates stated in 
Properties  \ref{list1-m}-\ref{list3-m} of Proposition \ref{propmultiscale}. 
\end{itemize}

All the precise assumptions  on the operators $ {\cal L}_r $ in \eqref{fo2-p}
and ${\cal L}_{r, \mu} $ in \eqref{fo1-p} are  
stated   in Definition \ref{definition:Xr}
(we neglect in \eqref{fo2-p}, \eqref{fo1-p} the projector $ {\mathfrak c} \Pi_{\mathbb S} $ which has a purely technical role, see Remark \ref{rem:co}).

\begin{remark}
We shall use the results of the multiscale Proposition \ref{propmultiscale} 
about the approximate invertibility of the operator 
$ {\cal L}_r $
 in \eqref{fo2-p} in Chapter \ref{sec:proof.Almost-inv},
 and, for the operator $ {\cal L}_{r,\mu} $ defined in \eqref{fo1-p},  in Chapter \ref{sec:splitting}. 
 In the next Section \ref{sec:ideas} we shall describe their role in the proof of Theorem \ref{thm:main}. 
\end{remark}

The application of multiscale techniques to the operators \eqref{fo2-p}, \eqref{fo1-p}
is more involved than for the second order quasi-periodic operator \eqref{model-B-int} 
presented in the previous section. We finally 
describe the main steps for the proof of the multiscale Proposition
\ref{propmultiscale}. This may serve as a road-map for the technical  proof given in 
Chapter \ref{sec:multiscale}. 
For definiteness we consider 
$ {\cal L}_{r,\mu} $. 

In the Fourier exponential basis,  the operator $ {\cal L}_{r,\mu} $ 
 is represented by a self-adjoint matrix 
\be\label{ADomT}
{\mathtt A}= {\mathtt D}_\om  +  {\mathtt T} 
\ee
with a diagonal part 
 (see \eqref{diag-secondo-caso}) 
\be\label{eigeDom0}
{\mathtt D}_\om = {\rm Diag}_{(\ell,j) \in \Z^{\es + d }}
\begin{pmatrix}
  \langle j \rangle_m - \mu &  \ii \om \cdot \ell  & 0 & 0  \\
 -  \ii \om \cdot \ell  &   \langle j \rangle_m - \mu   &  0 & 0  \\
 0 & 0  &    \langle j \rangle_m + \mu &  \ii \om \cdot \ell \\
 0 & 0  &   - \ii \om \cdot \ell &  \langle j \rangle_m + \mu
\end{pmatrix} \, , 
\ee
where, for some $ m > 0 $, 
$$ 
\langle j \rangle_m := \sqrt{|j|^2 + m } 
$$  
and  $ \mu := \mu (\e, \l ) \in \R $. 
The 
matrix $ {\mathtt T} $ in \eqref{ADomT} is   
$$
\begin{aligned}
& {\mathtt T} := ( {\mathtt T}_{\ell,j}^{\ell', j'})_{(\ell,j) \in \Z^{\es+d}, 
{(\ell', j')} \in \Z^{\es+d}} \, , \\
&  {\mathtt T}_{\ell,j}^{\ell',j'} :=   (D_V-D_m)_j^{j'}  - \mu [{\cal J} \Pi_{\mathbb S \cup \mathbb F}]_j^{j'} + r_{\ell,j}^{\ell', j'} \in {\rm Mat}(4 \times 4) \, .
\end{aligned}
$$
The matrix $ {\mathtt T} $ is  {\it T\"oplitz} in $ \ell $, 
namely $ {\mathtt T}_{\ell,j}^{\ell',j'} $ depends only on the indices $ \ell - \ell', j, j'  $.

The infinitely many eigenvalues of the matrix $ { \mathtt D}_\om  $ are 
$$
\sqrt{|j|^2+m}  \pm  \mu \pm  \om \cdot \ell  \, , \quad j \in \Z^d \, , \  \ell \in \Z^\es \, .
$$
By Proposition \ref{lemma-DVvsDm}
and  Lemma \ref{pisig},  the matrix $ {\mathtt T}  $ satisfies the 
off-diagonal decay 
\be\label{Ts1+-i}
| {\mathtt T} |_{+,s_1} :=  | D_m^{1/2} {\mathtt T}D_m^{1/2}  |_{s_1}   < + \infty \, , 
\ee
cf. with \eqref{rs1+} and \eqref{off-diago:T0}.
  
We introduce the index $ \mathfrak a \in {\mathfrak I} := \{ 1, 2,3,4\} $ to distinguish each component in
the $ 4 \times 4 $ matrix \eqref{eigeDom0}. 
Then the singular sites $  (\ell,j, \mathfrak a ) \in \Z^\es \times \Z^d \times {\mathfrak I} $ 
of $ {\cal L}_{r, \mu} $ are those integer
vectors such that, for some
choice of the signs $ \s_1 (\mathfrak a), \s_2 (\mathfrak a)  \in \{ -1, 1 \} $,  
\be\label{sing-sit-0}
\big| \sqrt{|j|^2+m}  + \s_1 (\mathfrak a) \mu + \s_2  (\mathfrak a)  \om \cdot \ell  \big| < 
\frac{\Theta}{ \sqrt{|j|^2+m} } \, .
\ee
In \eqref{sing-sit-0} 
 the constant $ \Theta := \Theta (V) $ is chosen large enough depending on the 
 multiplicative potential\index{Multiplicative potential}  $ V(x) $ 
(which is not small). 
Note that, if $  (\ell, j, \mathfrak a) $ is singular, then we recover the 
second order bound
$$
\big| |j|^2+m -  (\s_1( \mathfrak a ) \mu + \s_2 ( \mathfrak a )  \om \cdot \ell)^2 \big| 
\leq C \Theta \, ,
$$
which is similar  to \eqref{singul-i}. 

The invertibility of the restricted operator 
$ {\cal L}_{r,\mu,N} := {\it \Pi}_N ( {\cal L}_{r,\mu})_{| {\cal H}_N} $
and the proof of the off-diagonal decay estimates \eqref{Lip-sDM}
is  obtained  in Section \ref{sec:Green} 
by an inductive
application of the multiscale step Proposition \ref{propinv},  
which is  deduced from the corresponding Proposition \ref{propinvA}.
On the other hand,  Section \ref{sec:rightinv} contains the proof of the existence
 of a right inverse of 
 $ [{\cal L}_{r,\mu}]_N^{2N} := {\it \Pi}_N ( {\cal L}_{r,\mu})_{| {\cal H}_{2N}} $ 
 satisfying  \eqref{est:Right+1} 
 at the small scales $ N < N(\e) $ (note  
that the size of $ N $ depends on
$ \e $). 

In view of the multiscale proof, we 
consider the family of operators 
\be\label{L-teta-tr}
{\cal L}_{r, \mu} (\theta) := 
J \om \cdot \partial_\vphi + D_V + \ii \theta J + \mu {\cal J} \Pi_{\mathbb G} + r \, ,  
\quad \theta \in \R \, , 
\ee
which is represented, in Fourier basis, by a self-adjoint 
matrix denoted $ {\mathtt A}(\theta ) $.

The monotonicity assumption in item \ref{assu:pos} of Definition \ref{definition:Xr},
see \eqref{vare2}, 
allows to  obtain effective measure estimates and
complexity bounds similar to \eqref{tetabadweak0}.  
Note however that 
\begin{itemize}
\item by \eqref{vare2} the eigenvalues 
of $ {\cal L}_{r,\mu,N} $ vary in $ \lambda $ with only a $ O(\e^2) $-speed,  
creating  further difficulties with respect to the previous section. 
\end{itemize}
The verification of assumption (H3) of the multiscale step Proposition
\ref{propinv},
concerning separation properties of the $ N $-bad sites, is  a key part of the analysis.
The new notion of $ N $-bad site is introduced in Definition \ref{GBsite}, 
according to  the new definition 
of $ N $-good matrix  introduced in Definition \ref{goodmatrix}
(the difference is due to the fact that 
the singular sites are defined as in \eqref{sing-sit-0}, and
\eqref{Ts1+-i} holds).

Then, in Section \ref{sec:sepabad}, we prove that a $ \Gamma $-chain of 
singular sites is not too long, under the bound \eqref{cardinality} for its time components, see 
Lemma \ref{Bourgain}. Finally in Lemma \ref{Ntime} we are able to 
conclude, for the parameters $ \l $ which are $N$-good,  
i.e. which satisfy
\be\label{BNcomponents0-i}
\begin{aligned} 
 \forall \, j_0 \in \Z^d  \, , \  
B_{N}(j_0;  \l)  & :=    
\Big\{ \teta \in \R \, : \,   {\mathtt A}_{N, 0, j_0}(\theta) \ {\rm is} \ \hbox{$N$-bad}  \Big\} 
 \subset \bigcup_{q = 1, \ldots, N^{C(d,\es, \tau_0)}} I_q \, , 
\\
& \qquad  {\rm where} \  I_q  \ {\rm are \ intervals \ with \ measure } \  | I_q| \leq N^{-\t} \,  
\end{aligned}
\ee
(Definition \ref{def:freqgood}),
an  upper bound 
$ L \leq  N^{C_3}  $ (see \eqref{finale})
for the length $ L $ of a $ N^2 $-chain of $ N $-bad sites. This result 
implies the partition of the 
 $ N $-bad sites into clusters separated by $ N^2 $, according to the 
 assumption (H3) of the multiscale step Proposition \ref{propinv}.  

\begin{remark}
 Lemma \ref{Ntime}  requires a significant improvement with respect to the arguments in \cite{BB12}, described in the previous section:
the exponent  $ \tau $  in \eqref{BNcomponents0-i} is large, but independent of 
$ \chi $,  which defines the new scale $ N' = N^\chi $
in the multiscale step, 
see Remark \ref{tau-large}. This improvement  is required by the fact that 
by \eqref{vare2} the eigenvalues 
of $ {\cal L}_{r,\mu,N} $ vary in $ \lambda $ with only a $ O(\e^2) $-speed. 
\end{remark}

We finally comment 
why the multiscale analysis works also 
for  operators  $ {\cal L}_{r} $, $ {\cal L}_{r,\mu} $
of the form \eqref{fo2-p}, \eqref{fo1-p}, where $ D_V $ is not a Fourier multiplier. 
The reason is similar to \cite{BBo10}, \cite{BB12}. 
In Lemma \ref{Inizioind} we have to prove that {\it all} the parameters $ \l \in \Lambda $
are $ N $-good (Definition \ref{def:freqgood}) at the small scales $ N \leq N_0 $. 
We proceed as follows. We regard the operator $ {\cal L}_{r, \mu} (\theta) $ in 
\eqref{L-teta-tr} 
as a small perturbation of the 
operator 
$$ 
{\cal L}_{0, \mu} (\theta) = 
J \om \cdot \partial_\vphi + D_V + \ii \theta J + \mu {\cal J} \Pi_{\mathbb G}
$$
which is $ \vphi $-independent. 
Thus 
a lower bound on the modulus of the 
eigenvalues  of its restriction to $ {\cal H}_N $ 
implies an estimate in $ L^2 $ norm of its inverse, 
and, 
thanks to  {\it separation properties of the eigenvalues $ |j|^2 $ 
of $ - \Delta  $},  also a bound of its $ s$-decay norm as
$ O(N^{\tau'+ \varsigma s }) $. 
By a Neumann series perturbative argument this bound also persists for 
$ {\it \Pi}_N {\cal L}_{r, \mu} (\theta)_{| {\cal H}_N } $, taking 
$ \e $ small enough, up to the scales $ N \leq N_0 $. The proof is done precisely in Lemmata \ref{lemma1}-\ref{Inizioind}. The set of $ \theta $ such that the spectrum of
$ {\it \Pi}_{N} {\cal L}_{0, \mu} (\theta)_{| {\cal H}_N } $ is at a distance
$ O(N^{-\tau}) $ from $ 0 $ is contained into a union of 
intervals like \eqref{BNcomponents0-i},
implying the claimed complexity bounds. 
 The proof at  higher scales follows by the induction multiscale process.

\section{Outline of proof of Theorem \ref{thm:main}} \label{sec:ideas}

In this section  we present in detail the plan of proof of Theorem\index{Nonlinear Schr\"odinger equation} \ref{thm:main}, 
which occupies Chapters \ref{sec:Ham}-\ref{sec:NM}.
This section is a road map 
through 
the technical aspects of the proof.

\smallskip

In Chapter \ref{sec:Ham} we first write the second order wave equation \eqref{NLW2} 
as the first order Hamiltonian system 
\eqref{HSinK}-\eqref{K},
\be\label{HSqp-sym}
\begin{cases}
q_t = D_V p \cr
p_t = - D_V q - \e^2  D_V^{-\frac12} g(\e, x, D_V^{-\frac12} q ) \, , 
\end{cases}
\ee
where $ D_V := \sqrt{- \Delta + V(x)} $ is defined spectrally in \eqref{def:DV}, and 
the variables $ (q, p ) $ belong to a dense subspace of 
$$
H = L^2 (\T^d, \R) \times L^2 (\T^d, \R)  \, . 
$$
We prove polynomial off-diagonal decay of 
$ D_V $ in Section \ref{sec:offDV}. 

Having fixed finitely many tangential sites $ {\mathbb S} \subset \N $, 
we decompose (see \eqref{tang-norm}) the canonical variables $(q,p)$ 
into ``tangential"\index{Tangential variables} and ``normal"\index{Normal variables} components,
as
$$
(q,p) = \sum_{j \in {\mathbb S}} (q_j, p_j ) \Psi_j (x) + (Q,P) \, , \quad (q_j, p_j) \in \R^2 \, . 
$$
We recall that the functions  $ \Psi_j (x) $, $j \in \N$,  are the eigenfunctions of
the\index{Eigenfunctions of Sturm-Liouville operator}  Schr\"odinger operator $ - \Delta + V(x) $ that were introduced 
 in \eqref{auto-funzioni};
$ (Q,P) $ belongs to the  subspace $  H_{\mathbb S}^\bot $ of $H$, 
which is the $ L^2$-orthogonal   of 
\be\label{H-defS}
H_{\mathbb S} := \Big\{ {\mathop \sum}_{j \in {\mathbb S}} (q_j, p_j ) \Psi_j (x) \, , \,
(q_j, p_j) \in \R^2 \Big\} \, 
\ee
ans is called the normal subspace.
The dynamics of \eqref{HSqp-sym} on the symplectic subspaces 
$ H_{\mathbb S} $ and $ H_{\mathbb S}^\bot $ 
is handled quite differently.  
 On the 
 subspace $ H_{\mathbb S} $ we introduce action-angle variables\index{Action-Angle variables}
  $ (\theta, I) $,  by setting (see \eqref{AA})
$$
(q_j,p_j) := \sqrt{2 I_j} \big( \cos \teta_j, - \sin \teta_j \big) \, , \quad \forall j \in {\mathbb S } \, .   
$$
In these new coordinates
 the solutions \eqref{soluzione-lineare} of the linear wave equation  \eqref{Linear}
are described as 
the continuous family of quasi-periodic solutions\index{Quasi-periodic solution}  
\be\label{family-QP0}
\teta (t) =  
\bar \mu t \, , \ I (t) = \xi \, ,   \ Q (t) = P (t) =  0 \, , \quad \xi \in \R^\es_+ \, ,
\ee
with a frequency vector $ \bar \mu = ( \mu_j )_{j \in {\mathbb S}} $  
which is independent of the unperturbed actions $ \xi \in \R^\es_+ $.
Note that, by the assumption 
\eqref{diop}, the unperturbed tangential frequency vector $\bar \mu $ is Diophantine\index{Diophantine vector}. 
Introducing the translated action variable $ y  $ by setting 
$$ 
I = \xi + y \, , \quad y \in \R^\es \, , 
$$
the quasi-periodic solutions \eqref{family-QP0}
densely fill  the 
invariant torus 
$$ 
\T^\es_\theta \times \{0\}_y \times \{(0,0)\}_{(Q,P)} \, .
$$ 
In the variables $ (\theta, y, Q, P ) $ 
the Hamiltonian system \eqref{HSqp-sym} assumes the form 
\be \label{HS00expI}
\begin{cases}
 \dot \teta - \bar \mu - \e^2 \partial_y R (\theta, y, Q, \xi )  = 0   \cr 
\dot y   +  \e^2 \partial_\teta R (\theta, y, Q, \xi ) = 0 \cr
\big( \pa_t  -  J D_V \big) (Q, P) + \e^2 \big(0, \nabla_Q R (\theta, y, Q, \xi ) \big) = 0 \, , 
\end{cases}
\ee
see \eqref{HS00exp},  where $ J $ is the symplectic matrix 
$$ 
J = \begin{pmatrix}
0  &   {\rm Id}   \\
 - {\rm Id} & 0     \\   
\end{pmatrix} \, , 
$$ 
and  
 $ R(\theta, y, Q, \xi ) $ is the Hamiltonian in \eqref{restoNN}.

The goal is now to  look, for $ \e $ sufficiently small, for 
quasi-periodic solutions 
\be\label{qp-pre}
t \mapsto \big( \om t + \vartheta (\om t ), y(\om t ), Q (\om t ), P (\om t ) \big)
\ee
of the nonlinear Hamiltonian system 
\eqref{HS00expI}, 
with a frequency vector $ \omega \in \R^\es $, $ O(\e^2 ) $-close to $ \bar \mu $, to be determined, and where the function 
$$ 
\vphi \mapsto 
(\vartheta (\vphi), y(\vphi ), Q (\vphi), P (\vphi)) \in  
\R^\es \times  \R^\es \times H_{\mathbb S}^\bot 
$$
is periodic in the variable  $ \vphi = (\vphi_j)_{j \in {\mathbb S} } \in \T^\es $, and close to
$ (0, 0, 0, 0 ) $.  

We shall be able to constrain 
the frequency vector  $ \omega $ in \eqref{qp-pre} to 
a fixed ``admissible" 
 direction $ \bar \om_\e $,  as stated in \eqref{frequencycolinear}-\eqref{def omep},
 namely 
 \be\label{omega-con} 
  \omega = (1 + \e^2 \l) \bar \om_\e \, , \quad \l \in \Lambda = [- \l_0, \l_0]  \, ,
 \ee
  where  $ \bar \om_\e \in \R^{\es} $ 
 satisfies the Diophantine\index{Diophantine vector} conditions \eqref{dioep} and \eqref{NRgt1}.
In Lemma \ref{lemma:rho-dioph} we prove that these conditions are actually satisfied by 
``most" vectors $ \bar \om_\e =  \bar \mu +\e^2 \zeta $ close to $ \bar \mu $.
We shall use the $ 1$-dimensional parameter $ \lambda $, which corresponds to a time-rescaling, in order to impose all the non-resonance conditions required by our KAM construction, in particular along the multiscale analysis of the linearized operator. 
Note that  
$ \lambda $ has to be considered as an ``internal" parameter of the wave equation \eqref{NLW1}. 

\begin{remark}
The existence of quasi-periodic solutions with tangential frequencies constrained along a fixed direction
had been proved, for finite dimensional autonomous Hamiltonian systems,  
by Eliasson \cite{El88} and Bourgain \cite{B1}, and 
for $ 1$-$d$ nonlinear autonomous wave and Sch\"rodinger equations in \cite{BBi11}.
Results in the easier case of quasi-periodically forced PDEs, 
where $ \om $ is an external parameter constrained to a fixed direction, 
have been obtained in \cite{BBo10} for NLS, \cite{BB12} for NLW, 
\cite{BBM-Airy} for KdV.
\end{remark}

The search of an embedded invariant torus 
 $$ 
\vphi \mapsto  i(\vphi) = (\vphi, 0,0) + (\vartheta (\vphi), y(\vphi ), (Q,P)(\vphi)) \, , 
 $$ 
 for  the Hamiltonian system \eqref{HS00expI}, supporting quasi-periodic solutions with 
 frequency $   \omega = (1 + \e^2 \l) \bar \om_\e $ as in \eqref{omega-con}, 
 amounts to solving the functional equation $ {\cal F} (\l; i ) = 0 $ where ${\cal F} $
 is the nonlinear operator defined in \eqref{operatorF},
\be
\begin{aligned}\label{operatorFint}
 {\cal F} ( \lambda; i  ) 
 = \left(
\begin{array}{c}
\Dom \vartheta (\vphi) + \om - \bar \mu    - \e^2 (\partial_y R) ( i(\vphi), \xi   )   \\
\Dom y (\vphi)  +  \e^2 (\partial_\teta R) ( i(\vphi), \xi  )  \\
\Dom (Q,P) (\vphi) -  J D_V  (Q,P) (\vphi) + \e^2 \big( 0, (\nabla_Q R) (i(\vphi), \xi) \big)
\end{array}
\right) 
\end{aligned} \, . 
\ee
In \eqref{operatorFint}, $(\vartheta , y , (Q,P) )$  belongs, for some $s\geq s_0$, to 
$$ 
 H^s_\vphi (\T^\es, \R^\es ) 
\times H^s_\vphi (\T^\es, \R^\es)  
\times ( {\mathcal H}^s (\T^\es \times \T^d, \R^2) \cap L^2( \T^\es, H_{\mathbb S}^\bot) ) \, . 
$$
The action of the operator $ {\cal F} $ on the above  scale of Sobolev spaces makes lose 
one $ (\vphi, x) $-derivative due to  the unbounded 
operators $ \Dom $ and $ D_V $. 

We have that 
 $$ 
 {\cal F} (\l;  \vphi, 0, 0, 0 ) = O(\e^2 ) \, ,  \quad  \forall  \l  \in \Lambda \, . 
 $$ 
Then,  in section
 \ref{NM:first-step},  using  the unperturbed first order  Melnikov non-resonance condition \eqref{1Mel},
we are able to construct the first approximate\index{First approximate solution} solution $ i_1 (\vphi) $ such that 
$$ 
{\cal F} (\l; i_1 (\vphi) ) = O(\e^4)  \, , \quad  \forall  \l  \in \Lambda\, . 
$$
Next, in Theorem \ref{thm:NM}, which contains the core of the result, we construct, by means of 
a Nash-Moser implicit function iterative scheme,  for ``most" values of $ \l \in \Lambda $, 
a solution $ i_\infty (\l; \vphi) $ of the equation 
$$ 
{\cal F} (\l; i_\infty (\l; \vphi) ) = 0 \, , 
$$
obtaining an invariant torus of 
 \eqref{HS00expI}  with frequency $ \omega = (1 + \e^2 \l) \bar \om_\e $.    
The iteration is performed in Chapter \ref{sec:NM}. In particular, in Theorem \ref{thm:NM-ite}
we construct a sequence of approximate solutions $ i_n (\l; \vphi) $ which converges 
to $ i_\infty  (\l; \vphi) $ for $ \l $ in a set of large measure.  
The key point is to prove the approximate invertibility of the 
linearized operators 
$$ 
d_{i} {\cal F}(\l; i_n (\l; \vphi) ) \, ,
$$ 
obtained at any approximate quasi-periodic solution $ i_n :=  i_n  (\l; \vphi) $ along the Nash-Moser iteration, 
for a large set of $ \l $'s, 
together with suitable  tame estimates for the approximate inverse in high Sobolev norms
(clearly, due to  small divisors, such approximate  inverse loses  derivatives). 
This is achieved by the analysis performed in Chapters   
\ref{sezione almost approximate inverse}-\ref{sec:proof.Almost-inv}.

\smallskip

In Chapter \ref{sezione almost approximate inverse}
 we  implement the general strategy proposed in \cite{BBField}.
 Instead of (approximately) inverting $ d_i {\cal F}(\l; i_n ) $,  
where all the $ (\vartheta, y, Q, P  ) $  components are coupled 
by the differential of the nonlinear term
$ ( - \pa_y R, \pa_\theta R, 0, - \nabla_Q R)    $ 
in \eqref{operatorFint}, 
we introduce suitable symplectic coordinates 
$$ 
(\phi, \zeta, w) \in \T^\es \times \R^\es \times H_{\mathbb S}^\bot 
$$ 
in which it is sufficient to  
(approximately) invert 
the linear operator $ {\mathbb D } (i_n) $ defined in  \eqref{operatore inverso approssimato}.
The advantage is that the components  of the operator 
$ {\mathbb D } (i_n) $ can be inverted in  a  triangular way. We mean that  
first one inverts the operator in the tangential action component $ \widehat \zeta $,
then the normal one for $ \widehat w \in H_{\mathbb S}^\bot  $, 
and finally  the operator for the tangential angles 
$ \widehat \phi $. 
This construction is implemented in detail 
in Chapter \ref{sezione almost approximate inverse}, completed with the results 
reported in Appendix
\ref{sec:2}. 

\begin{remark}
The above 
decomposition is deeply related  to the Nash-Moser approach 
for isotropic  tori of finite dimensional Hamiltonian systems in 
Herman-Fejoz \cite{HF}. A related construction for reversible \index{Reversible PDE}
PDEs is performed in  \cite{CFP}. 
\end{remark}

After this transformation, 
the main issue is reduced to proving 
the approximate invertibility of the  
linear operator ${\cal L}_\omega (i_n)  $ defined in \eqref{Lomega def} 
which 
acts on functions $  h : \T^\es  \to  H_{\mathbb S}^\bot $ with values
in the normal subspace $ H_{\mathbb S}^\bot  $. As proved in Lemma \ref{thm:Lin+FBR}, this operator 
has the form
\be\label{lin:restr}
{\cal L}_\om  (i_n) = \om \cdot \partial_\vphi  - J (D_V + \e^2 {\mathtt B}(\vphi)  + {\mathtt r}_\e (\vphi))  \, , 
\quad D_V = \sqrt{- \Delta + V(x)} \, , 
\ee
where  $ {\mathtt B}(\vphi) $ is the self-adjoint operator defined in  \eqref{form-of-B-ge}
and  $  {\mathtt r}_\e (\vphi)   $ is a self-adjoint remainder of size $ O(\e^4) $, more precisely 
it satisfies 
the quantitative bounds \eqref{estimate:rep-ge}-\eqref{estimate:rep-ge0}.
Note that 
$ J (D_V + \e^2 {\mathtt B}(\vphi)  + {\mathtt r}_\e (\vphi)) w $
is a linear Hamiltonian  vector field and that 
the corresponding quasi-periodic 
operator \eqref{lin:restr} is a small deformation of the operator 
obtained by linearizing the normal component of 
${\cal F}$, defined in \eqref{operatorFint},  at an approximate solution.  
A quasi-periodic operator as \eqref{lin:restr} is called Hamiltonian, see
Definition \ref{def:HAMS}.  

\smallskip

Chapters \ref{sec:6n}-\ref{sec:proof.Almost-inv} are devoted to prove the existence of an approximate right inverse 
of $ {\cal L}_\om  (i_n) $, as stated in Proposition \ref{prop:inv-ap-vero}, for ``most"
values of $ \lambda \in \Lambda $. This is obtained in several steps.

\begin{remark}\label{rem:per-split}
We cannot directly apply to the operator 
$ \frac{1}{1 + \e^2 \lambda }{\cal L}_\om (i_n) $
the approach developed for the quasi-periodically forced NLW and NLS equations
in \cite{BBo10}, \cite{BB12}, \cite{BCP}, described in Section 
\ref{sec:MULTI-QP}. Actually, 
since $ \pa_\lambda {\mathtt B}(\vphi) = O(1) $ the operator 
$ \pa_\l (\frac{1}{1 + \e^2 \lambda }{\cal L}_\om (i_n))  $ is not positive or negative definite, 
posing a  serious difficulty for verifying that the eigenvalues 
of its  finite dimensional restrictions are  in modulus bounded away from zero 
for most values of the parameter $ \lambda $. 
In order to overcome this problem we perform the splitting of the normal subspace 
that we describe below. 
\end{remark}

In Chapter \ref{sec:6n} 
 we first apply an averaging procedure to 
the operator 
\be\label{A-intro}
\begin{aligned}
 \frac{{\cal L}_\om}{1 + \e^2 \l } = & 
\ \bar \om_\e \cdot \partial_\vphi - J \Big( {\mathtt A} +  \frac{{\mathtt r}_\e }{1+\e^2 \lambda} \Big)  \, , \quad 
{\mathtt A}  := \frac{D_V}{1+ \e^2 \l} + \frac{\e^2 {\mathtt B}}{1+ \e^2 \l}   \, .
\end{aligned} 
\ee
Consider  the splitting  of the normal subspace  
$$
\begin{aligned}
& \qquad \qquad \qquad 
H_{\mathbb S}^\bot = H_{\mathbb M} \oplus H_{\mathbb M}^\bot \, , 
 \\
& H_{\mathbb M} := \Big\{ {\mathop \sum}_{j \in {\mathbb M}} (q_j, p_j ) \Psi_j (x) \, , \,
(q_j, p_j) \in \R^2 \Big\} \, , 
\end{aligned}
$$
where $ {\mathbb M} $ is the finite set of Lemma \ref{choice:M}, which  contains 
$ {\mathbb F} $ (we explain below its role). 
Taking in $ H_{\mathbb M} $ the basis 
$ \{ ( \Psi_j(x) ,0) \, ,  (0, \Psi_j(x) ) \}_{j \in {\mathbb M}} $, the 
linear operator ${\mathtt A} $ is represented by the matrix 
\be\label{def:primo-op-bis-in}
{\mathtt A} = 
 \begin{pmatrix}
 {\rm Diag}_{j \in {\mathbb M} } \frac{{\mu}_j}{1+ \e^2 \l }  {\rm Id}_2 &  0   \\
  0 &  \frac{D_V}{1+  \e^2 \l}  \\
\end{pmatrix} + \frac{\e^2}{1+ \e^2 \lambda} 
 \begin{pmatrix}
 {\mathtt B}_{{\mathbb M}}^{{\mathbb M}} &  {\mathtt B}_{{\mathbb M}}^{{\mathbb M}^c}    \\
 {\mathtt B}_{{\mathbb M}^c}^{{\mathbb M}}   &  {\mathtt B}_{{\mathbb M}^c}^{{\mathbb M}^c}   \\
\end{pmatrix} \, .
\ee
The $ \mu_j  $ are the eigenvalues of $ D_V = \sqrt{- \Delta + V(x)} $,  see  \eqref{auto-funzioni}.  
We perform an averaging transformation 
 to remove the off-diagonal terms of order $ O(\e^2 )$
in \eqref{def:primo-op-bis-in}, 
obtaining the quasi-periodic Hamiltonian operator \eqref{transf-op-1}-\eqref{newA+}, which 
has the form  
\be\label{Lnormal0}
\ppavphi - J ({\mathtt A}_0 + \varrho^+) 
\ee
where  $ {\mathtt A}_0 $ is split admissible\index{Split admissible operator} 
according to Definition \ref{def:calC} (see Lemma \ref{A-0-splitted})
and  the coupling operator 
$ \varrho^+ $ has size $ O(\e^4) $ in the $ | \ |_{+, s_1} $-norm defined in 
\eqref{rs1+}. More precisely it satisfies 
 $ |  \varrho^+  |_{\Lip, +, s_1} =
 O( \e^4)  $ where the norm $ | \ |_{\Lip,+,s_1}$ is introduced in Definition \ref{def:decay-sub}. 
In this step we use the unperturbed second order Melnikov conditions
 \eqref{2Mel+}-\eqref{2Mel rafforzate}.  
 
The main features of a split-admissible operator are the following ones:
\begin{enumerate}
\item A {\it split-admissible} operator is  self-adjoint and block diagonal 
with respect to the orthogonal splitting 
\be\label{spliFG}
\begin{aligned}
& \qquad \qquad \qquad 
H_{\mathbb S}^\bot  = H_{\mathbb F} \oplus H_{{\mathbb G}}  \, , 
 \\
& H_{\mathbb F} := \Big\{ {\mathop \sum}_{j \in {\mathbb F}} (q_j, p_j ) \Psi_j (x) \, , \,
(q_j, p_j) \in \R^2 \Big\} \, , \\
& H_{\mathbb G} := \Big\{ {\mathop \sum}_{j \in {\mathbb G}} (q_j, p_j ) \Psi_j (x) \, , \,
(q_j, p_j) \in \R^2 \Big\} \, \, , 
\end{aligned}
\ee
i.e. of  the form 
\be\label{form:A0-intro}
{\mathtt A}_0 = 
\begin{pmatrix}
D_0(\e , \lambda) & 0  \\
0  &  V_0(\e, \l, \varphi)  
\end{pmatrix}  \, .
\ee
We remind that the  subsets $ {\mathbb F} \subset {\mathbb M} $ and $ {\mathbb G} $  
 of the normal sites $ {\mathbb S}^c  $  in \eqref{spliFG} are defined in 
\eqref{taglio:pos-neg}, and that 
$$ 
{\mathbb F} \cup {\mathbb G} = {\mathbb S}^c  \, , \quad 
 {\mathbb F} \cap {\mathbb G} = \emptyset  \, . 
 $$
\item The operator $D_0 (\e, \lambda)  $ is,  
in the basis of the eigenfunctions $ \{ (\Psi_j,0) $, $(0,\Psi_j)\}_{ j \in {\mathbb F}} $, 
diagonal,  
\be\label{form:D0-intro}
D_0 (\e, \lambda) = {\rm Diag}_{j \in {\mathbb F}} \, \mu_j (\e,\l) {\rm Id}_2 \, , \quad \mu_j (\e,\l) \in \R \, ,
\ee
with eigenvalues $ \mu_j ( \e, \lambda ) = \mu_j + O (\e^2) $
(the $ \mu_j  $ are defined in  \eqref{auto-funzioni})  satisfying the non-degeneracy conditions
\eqref{Hyp2}-\eqref{Hyp4}. 
\item {\bf (Monotonicity)}
For all $ j \in {\mathbb F} $,  
\be\label{Hyp1-intro}
\begin{cases}
{\mathfrak d}_\l  \big( V_0 (\e , \l) + \mu_j (\e, \l) {\rm Id}  \big) \leq - c_1 \e^2 \  \cr
{\mathfrak d}_\l \big( V_0 (\e , \l) - \mu_j (\e, \l) {\rm Id} \big) \leq - c_1 \e^2  
 \end{cases}
\ee
where we use the notation \eqref{nota-frak}.
\end{enumerate}

\noindent
In order to prove Lemma  \ref{A-0-splitted}, saying that $ {\mathtt A}_0 $ is a 
split-admissible operator, 
the following is used:  
\begin{enumerate}
\item the precise knowledge of the shift of the 
tangential and normal frequencies\index{Shifted normal frequencies} by the nonlinearity $ a(x) u^3 +  O( u^4 )$ 
(via the Birkhoff matrices $ \Ab, \Bb $ in \eqref{def:AB});
\item 
 the additional second order Melnikov non-resonance 
conditions  \eqref{2Mel rafforzate+}-\eqref{2Mel rafforzate};
\item  
the choice \eqref{taglio:pos-neg} of the subsets $ {\mathbb F}$ and 
$ {\mathbb G} $;
\item  
the non-degeneracy conditions \eqref{non-reso}-\eqref{non-reso1}.
\end{enumerate}
More precisely we use
   \eqref{2Mel rafforzate+}-\eqref{2Mel rafforzate} and \eqref{taglio:pos-neg}  
to prove the monotonicity property \eqref{Hyp1-intro}, and the non-degeneracy conditions 
\eqref{non-reso}-\eqref{non-reso1} to prove \eqref{Hyp2}-\eqref{Hyp4}. 
The properties \eqref{Hyp2}-\eqref{Hyp1} allow to prove that the non-resonance conditions
required along the multiscale analysis of the linearized operator are fulfilled for a large set of values of the parameter $ \l $.
We comment about this issue below, around  \eqref{der-app}. 

\smallskip

The quasi-periodic Hamiltonian operator  $ \ppavphi - J ({\mathtt A}_0 + \varrho^+)  $ 
in \eqref{Lnormal0} is 
in a suitable form to apply Proposition \ref{prop-cruciale}, in order to prove that it admits 
an approximate right inverse. 
Proposition \ref{prop-cruciale}  is proved in Chapters \ref{sec:splitting} and \ref{sec:proof.Almost-inv}.

\smallskip

In Chapter \ref{sec:splitting} 
we block-diagonalize  $ \ppavphi - J ({\mathtt A}_0 + \varrho^+) $, 
according to the splitting $ H_{\mathbb S}^\bot = H_{\mathbb F} \oplus H_{\mathbb G} $,  
up to a very small coupling term, see Corollary  \ref{cor:split}. 
More precisely, we conjugate,  via symplectic transformations,  
the quasi-periodic Hamiltonian operator 
$$ 
 \ppavphi  - J ({\mathtt A}_0 + \varrho^+)  
 $$ 
 in \eqref{Lnormal0}, where 
the  block non-diagonal term  $ \varrho^+ $ has size $ O(\e^4) $ in $ | \ |_{+,s_1} $ norm, 
to the quasi-periodic Hamiltonian operator  
 in \eqref{goal-trasf-iterata}, 
\be\label{Lnormal-ind}
\bar \om_\e \cdot \partial_\vphi - J (A_\ind + \rho_\ind) \, , 
\ee
where the operator $A_\ind $ is split admissible\index{Split admissible operator}.
More precisely  $A_\ind $ has the form \eqref{form:A0-intro}, 
with
\be\label{Aind-intro}
\begin{aligned}
& \ \qquad A_\ind  = 
\begin{pmatrix}
D_\ind (\e , \lambda) & 0  \\
0  &  V_\ind (\e, \l, \varphi)  
\end{pmatrix}  \, ,  \\
& \| D_\ind - D_0 \|_0 = o(\e^2) \, , \ \| V_\ind - V_0 \|_0 = o(\e^2)  \, , 
\end{aligned}
\ee
(the norm $ \| \ \|_0 = \| \ \|_{{\cal L}(L^2)}$ is the operatorial $ L^2 $ norm)
and the  block non-diagonal self-adjoint term $ \rho_\ind $ is 
super-exponentially small, i.e. satisfies \eqref{form-A'-m} 
\be\label{re-ind}
|\rho_\ind|_{+,s_1} = O \big( (\e^3)^{ (\frac{3}{2})^{\ind-1}} \big) \, . 
\ee
The proof of the splitting Corollary  \ref{cor:split} is based on an iterative application of the `splitting step" Proposition \ref{prophomeq}. The goal of this proposition is 
to  block diagonalize a
Hamiltonian operator   of the form 
$$ 
\bar \om_\e \cdot \partial_\vphi - J (A_0 + \rho) 
\quad {\rm where} \quad 
\begin{aligned}
& \qquad \rho(\vphi) = \begin{pmatrix} \rho_1 (\vphi)& \rho_2 (\vphi)^* \\
\rho_2 (\vphi) & 0    \end{pmatrix} \in {\cal L} (H_{\mathbb S}^\bot ) \, , \\
& \rho_1 (\vphi) = \rho_1^* (\vphi) \in {\cal L} (H_{\mathbb F}) \, ,  \   \rho_2 (\vphi) \in {\cal L} (H_{\mathbb F}, H_{\mathbb G}) \, ,
\end{aligned}
$$
into a new Hamiltonian
operator 
$$
 \bar \om_\e \cdot \partial_\vphi - J (A_0^+ + \rho^+)  
 $$ 
where $ A_0^+ $ is block-diagonal with respect to $ H_{\mathbb F} \oplus 
H_{\mathbb G} $ and 
the coupling remainder
 $\rho^+ $ is 
much smaller than the previous one, i.e. $ \rho^+ = O(\| \rho \|^{3/2}) $
(compare the size of 
$ \rho $ in \eqref{rho-R0:small-0} and that of $ \rho^+ $ in \eqref{propo1}). 
The iteration is based on a  super-convergent Nash-Moser 
scheme, to compensate the loss of derivatives due to the small divisors. 

\smallskip

In this decoupling procedure (Proposition \ref{prophomeq}) a central role is played by 
the possibility of solving approximately  the homological equations, see
\eqref{lineqd}-\eqref{eqhoma},
\begin{align} 
& J \ppavphi d + D_0 d + J d J D_0= J \rho_1  \label{lineqd-int} \\
& J \ppavphi a -J V_0 J a + J a J D_0 =J \rho_2 \, ,  \label{eqhoma-int} 
\end{align}
where $ d(\vphi)  \in {\cal L} (H_{\mathbb F}) $ is self-adjoint, and 
$ a(\vphi) \in {\cal L} (H_{\mathbb F}, H_{\mathbb G}) $, for all $ \vphi \in \T^\es $. 

For solving the homological equation 
\eqref{lineqd-int} (Lemma  \ref{homdiag}), 
we need  the second order Melnikov non-resonance conditions 
\eqref{def:Lambda1}, which concern only the finitely many normal modes in 
${\mathbb F} $.
 Then  we use the non-degeneracy properties \eqref{Hyp2}-\eqref{Hyp4}
 to prove that  they are fulfilled for most values of $ \l $'s
(see  the measure estimate of  Lemma  
 \ref{lemma:measure1}). 
 
 On the other hand we use the monotonicity property \eqref{Hyp1-intro} to solve 
(approximately) the homological equation \eqref{eqhoma-int} for most values of $ \l $'s. This is a difficult step where we use the multiscale techniques of Chapter 
\ref{sec:multiscale}. Let us try to specify some key aspects of 
our approach. We have to solve, approximately,  each equation
\be\label{eqhomai1-int}
T_j (a^j)  = J \rho_2^j \, , \quad \forall j \in {\mathbb F } \, , 
\ee
where 
$$ 
\begin{aligned}
 a^j (\vphi) & := (a (\vphi))_{|H_{j}} \, , \quad 
 \rho_2^j (\vphi) := (\rho_2 (\vphi) )_{|H_j} \in {\cal L}(H_j, H_{\mathbb G}) \, , \\  
& 
 \quad \forall  \vphi \in \T^\es  \, , \qquad 
 H_j := \big\{ (q_j, p_j ) \Psi_j (x) \, , \, (q_j, p_j) \in \R^2 \big\} \, , 
\end{aligned}
$$  
and 
$ T_j  $ is the linear operator (see \eqref{eqhomai1})
$$
T_j (a^j) := J \ppavphi a^j -J V_0 J a^j + \mu_j (\e, \lambda) J a^j J  \, .
$$
We solve approximately \eqref{eqhomai1-int} by applying the multiscale Proposition \ref{propmultiscale}
 to the extended operator (see \eqref{TB1B2}), 
  \be \label{TB1B2-int}
 (1+ \e^2 \l) T^\sharp_j \, , \quad
 T^\sharp_j := 
\begin{pmatrix}
 J \ppavphi + \frac{D_V}{1+ \e^2 \l} + \frac{\co}{1+ \e^2 \l} \Pi_{\mathbb S}
 & 0    \\
 0 & T_j       
\end{pmatrix} 
\, ,
\ee
 which acts on $ \vphi $-dependent  functions with values in 
 $ {\cal L}(H_j, H ) = {\cal L}( H_j,  H_{\mathbb S \cup  \mathbb F} ) \oplus 
{\cal L}( H_j,  H_{\mathbb G} )  $, see \eqref{deco:HjFSG}.
The operator $ \Pi_{\mathbb S} $ is the $ L^2 $ projector on the 
subspace $ H_{\mathbb S} $ in \eqref{H-defS} and $\co $ is a positive constant. 
Identifying  (see \eqref{ident-H-4-me})
$$ 
{\cal L} (H_j, H ) \sim H \times H \sim ( L^2(\T^d, \R))^4 \, , 
$$
the operator $ T^\sharp_j $
can be regarded to act on (dense subspaces of) 
the whole $ ( L^2(\T^\es \times \T^d, \R))^4 $
 and not only in subspaces of $ \vphi $-dependent functions with values in the normal subspace
 $ H_{\mathbb S}^\bot \times H_{\mathbb S}^\bot $.  
The self-adjoint operator $ T^\sharp_j $ has the form \eqref{Ext1}, thus as  in Definition 
\ref{definition:Xr}-($ii$).  
We can apply the multiscale Proposition \ref{propmultiscale} to $ T^\sharp_j $. 
The monotonicity property 
\eqref{Hyp1-intro} implies, 
in Lemma \ref{pos:def-var}, 
 the sign condition  $ {\mathfrak d}_\lambda T_j^\sharp \leq - c \e^2 $ required in Definition 
 \ref{definition:Xr}-item \ref{assu:pos}.
  
An important property of the extension is the following:
\begin{itemize}
\item
By \eqref{TB1B2-int}, 
 the subspaces $ H_{\mathbb S} $ and $ H_{\mathbb S}^\bot $ 
are {\it invariant} for the action of the extended operator $ T^\sharp_j $, in the sense that it maps  a function of variable $\vphi$ taking its values in
$ H_{\mathbb S} $ (resp.  $ H_{\mathbb S}^\bot $) to a function taking its values in the same subspace. Therefore 
estimates for the approximate 
inverse of $ T^\sharp_j $ (obtained by the multiscale analysis with the
exponential basis) provide also estimates for the approximate 
inverse of $ T_j $ and thus for the 
approximate solution of \eqref{eqhomai1-int}. 
\end{itemize}

We now explain the relevance of the decomposition 
 $ {\mathbb S}^c = {\mathbb F} \cup {\mathbb G} $ of the normal sites introduced in \eqref{taglio:pos-neg-0}-\eqref{taglio:pos-neg}, and why we are able to obtain,  
in the splitting Corollary  \ref{cor:split},  
 a  block-diagonal operator as 
 in \eqref{Lnormal-ind}-\eqref{Aind-intro}, with an error
 like  \eqref{re-ind}, with 
 respect to the splitting \eqref{spliFG}. 

The decomposition  $ {\mathbb S}^c = {\mathbb F} \cup {\mathbb G} $ of the normal sites 
is important for  the proof of the `splitting step" Proposition \ref{prophomeq}. 
 As  proved in Section \ref{sec:shifted-tan},  
 under the effect of the nonlinearity $ a(x) u^3 + O(u^4) $, 
 the  tangential frequency vector  of the expected 
 quasi-periodic solutions of \eqref{HS00expI} is, 
up to terms of $ O(\e^4 ) $,  
\be\label{shiftedTF0}
\bar \mu + \e^2 \Ab  \xi  
\ee
where $ \Ab  $ is the twist matrix\index{Twist matrix} defined in \eqref{def:AB}. Moreover, 
as proved in Lemma \ref{lem:ns}, 
 the perturbed normal frequencies of the Hamiltonian linear 
 operator \eqref{Lnormal0},  
 with normal indices\index{Shifted normal frequencies} in a large finite  set $ {\mathbb M}  \subset {\mathbb S}^c $,  
 admit the expansion,  up to $ O(\e^4 ) $, 
\be\label{Omega-normali-0}
 \mu_j +  \e^2  (\Bb \xi)_j \, , \quad  \forall j \in {\mathbb M} \, ,  
\ee
 where $ \Bb $ is the Birkhoff matrix  defined in \eqref{def:AB}.\index{Birkhoff matrices}
 Thus the Hamiltonian linear operator \eqref{Lnormal0}, 
 restricted to $ H_{\mathbb M } $, is
 in diagonal form,  up to terms $ O(\e^4) $. 
 The set $ {\mathbb M} $ contains $ {\mathbb F} $
 and it is fixed in Lemma \ref{choice:M} large enough so that the sign condition
 \eqref{negative-on-up} holds. Note that,  in order 
  to get the expansion \eqref{Omega-normali-0} 
  for all the indices $ j \in  {\mathbb M} $, and not just in $ {\mathbb F} $,
we have assumed the further second order Melnikov conditions \eqref{2Mel rafforzate+}-\eqref{2Mel rafforzate}, and performed the averaging  Proposition \ref{prop:op-averaged}.
Expressing $ \xi $ in terms of  $ \om $ by inverting the relation
$$ 
\om =  \bar \mu + \e^2 \Ab  \xi  \, , \quad \om = ( 1 + \e^2 \l) \bar \om_\e \, , 
$$
(the twist matrix $ \Ab $ is invertible by \eqref{A twist}), 
 the shifted normal frequencies \eqref{Omega-normali-0} become
\be\label{Omega-normali}
\begin{aligned}
&  \mu_j -  [\Bb \Ab^{-1}  \bar \mu ]_j  +  [\Bb \Ab^{-1} \om ]_j = \\ 
& 
\mu_j -  [\Bb \Ab^{-1}  \bar \mu ]_j  +  (1+ \e^2 \l)[\Bb \Ab^{-1} \bar \om_\e ]_j =: \Omega_j (\e, \l) \, . 
\end{aligned}
\ee
Then we divide \eqref{Omega-normali} by $ 1+ \e^2 \l $ and we  consider the derivative 
\be\label{derivata lam}
\frac{d}{d \l} \frac{ { \Omega}_j (\e, \l )}{ 1+ \e^2 \l } =  
\frac{- \e^2 }{(1+ \e^2 \l)^2 } (  \mu_j -   [\Bb  \Ab^{-1} \bar \mu]_j ) \, , \quad
 \forall j \in {\mathbb M} \, . 
\ee
We now provide a heuristic argument which suggests
 a perturbative procedure to decouple the Hamiltonian linear operator \eqref{Lnormal0}
with respect to $ H_{\mathbb F} $ and $ H_{\mathbb G}  $, for most values of the parameter
$ \l $.
In order to diagonalize a linear operator, which is a small perturbation of a diagonal matrix 
 with eigenvalues  $ \frac{ { \Omega}_j (\e, \l )}{ 1+ \e^2 \l }  $, 
 by means of a  reducibility iterative scheme, one would like to impose 
second order Melnikov non-resonance conditions of the form 
\be\label{sm:naive}
| \omega \cdot \ell +  \Omega_j (\e, \l) \pm  \Omega_k (\e, \l)| \geq \frac{\g}{\langle \ell \rangle^\tau} \, , 
\ee
for all $ (\ell, j, k) \neq (0,j,j) $.  
Analogously, in order to decouple the Hamiltonian linear operator \eqref{Lnormal0}
with respect to $ H_{\mathbb F} $ and $ H_{\mathbb G}  $, 
the natural non-resonance conditions are (at least at the first step)
\be\label{NR:naive}
\begin{aligned}
& |f_{\ell, j, k} (\l) | \geq \frac{\gamma}{ \langle \ell \rangle^\tau}  \qquad {\rm where} 
\qquad 
f_{\ell, j, k} (\l) := \bar \om_\e \cdot \ell +  \frac{\Omega_j (\e, \l)}{1 + \e^2 \l} \pm \frac{\Omega_k (\e, \l)}{1+ \e^2 \l}  \, \, , \\
& \forall \ell \in \Z^\es \, ,  \ j \in {\mathbb G} \, , \ k \in {\mathbb F}  \, . 
\end{aligned}
\ee
Thanks to \eqref{derivata lam}, 
for $ j \in {\mathbb G} $, $ k \in {\mathbb F} $, 
by the definition of $ {\mathbb F}  $ and $ {\mathbb G} $ in \eqref{taglio:pos-neg},
\begin{align}\label{der-app}
\pa_\lambda f_{\ell, j, k} (\l) & = 
\frac{- \e^2 }{(1+ \e^2 \l)^2 } (  \mu_j -   [\Bb  \Ab^{-1} \bar \mu]_j ) \mp 
\frac{ \e^2 }{(1+ \e^2 \l)^2 } (  \mu_k -   [\Bb  \Ab^{-1} \bar \mu]_k ) \nonumber  \\
& \leq - \frac{\e^2}{(1 + \e^2 \l)^2} c  \, ,
\end{align}
for some $ c>0$. So $\pa_\lambda f_{\ell, j, k} (\l)$
is negative for all $ \l $, allowing to prove that \eqref{NR:naive} is fulfilled for most values of $ \l $. 
This explains the  relevance of the splitting \eqref{taglio:pos-neg-0}-\eqref{taglio:pos-neg} on the normal indices $ {\mathbb S}^c $. 

Actually $ {\mathbb G} $ is an infinite set and  we can not impose \eqref{sm:naive}
for all the $ j \in {\mathbb G} $, since we have only the expansion  \eqref{Omega-normali-0} 
for finitely many $ j \in {\mathbb M } $. 
The homological equation 
to be solved, in order to decouple $ H_{\mathbb F} $ and $ H_{\mathbb G} $,  
 is indeed \eqref{eqhoma-int}, or equivalently \eqref{eqhomai1-int}.  
 
We now roughly explain the role of the finite set  $ {\mathbb M} $ in the normal sites
$  {\mathbb S}^c $. 
We take $ {\mathbb M}$ large enough as in Lemma \ref{choice:M} 
in order to prove that the operator $ {\mathtt A_0} $ in \eqref{Lnormal0} is split-admissible, see Lemma \ref{A-0-splitted}, and so  is the operator
 $ A_\ind $ in \eqref{Lnormal-ind}-\eqref{Aind-intro}. 
Indeed, for $ {\mathbb M}$ large enough,  the infinite dimensional operator 
$$
\pa_\lambda {\mathtt A_0}^{{\mathbb M}^c}_{{\mathbb M}^c}  = 
\pa_\lambda {\mathtt A}^{{\mathbb M}^c}_{{\mathbb M}^c}  \, ,  \quad {\rm where} \quad
{\mathtt A} := \frac{D_V + \e^2 {\mathtt B}(\vphi)}{1+ \e^2 \l} \ \  {\rm (see \ \eqref{A-intro})} \, , 
$$
is strongly negative definite  (see \eqref{negative-on-up}). Jointly with
\eqref{derivata lam} and \eqref{der-app}, this allows to prove the sign conditions \eqref{Hyp1-intro}.
As already mentioned,  this property allows to 
use  monotonicity arguments for families of self-adjoint matrices
to verify that the non-resonance conditions 
required to solve the  homological equation \eqref{eqhomai1-int}
are fulfilled for most values of the parameter $ \l $. 

\smallskip

After all this procedure, the linear operator \eqref{Lnormal0} 
has been approximately block-diagonalized, 
according to the decomposition $ H_{\mathbb F} \oplus H_{\mathbb G} $, 
obtaining,  as in \eqref{Lnormal-ind}, \eqref{Aind-intro}, the quasi-periodic Hamiltonian 
operator 
\be\label{ome-Arho}
\bar \om_\e \cdot \partial_\vphi - J (A_\ind + \rho_\ind) \, , 
\ee
where $\rho_\ind $ is a very small coupling term according to \eqref{re-ind}.
It remains to prove that it admits 
an approximate inverse,  for most values of $ \l $'s, satisfying tame estimates.
This is achieved in Chapter \ref{sec:proof.Almost-inv}. For this task we  
apply once more the multiscale
Proposition \ref{propmultiscale}. 
More precisely we first find  in Section \ref{sec:9.1}
an approximate inverse of 
the block-diagonal operator 
$$
\bar \om_\e \cdot \partial_\vphi - J A_\ind    \, , \quad 
A_\ind = \begin{pmatrix}
D_\ind (\e , \lambda) & 0  \\
0  &  V_\ind (\e, \l, \varphi)  
\end{pmatrix}   \, , 
$$
applying the multiscale Proposition \ref{propmultiscale} 
to an extension of the self-adjoint 
operator 
$$ 
J \bar \om_\e \cdot \partial_\vphi  
+ V_\ind (\e, \l, \varphi)  
$$ 
acting on a dense subspace of the whole $ H $, and not just 
$ H_{\mathbb S}^\bot $, see \eqref{def:T-tilde}-\eqref{wT}. 
The extended operator has the form $ {\cal L}_r $ in \eqref{fo2-p}, see Definition 
\ref{definition:Xr}-(i),    
and satisfies the 
sign condition 
$$
{\mathfrak d}_\l ({\cal L}_r   (1 + \e^2 \l )^{-1}) \leq - c \e^2 \, . 
$$
This allows to verify 
lower bounds for the moduli of the eigenvalues 
of finite dimensional restrictions of ${\cal L}_r  $, for most values of
$ \l $'s.
These lower bounds  amount to  
first order Melnikov type non-resonance conditions.

Finally, once we have constructed an approximate inverse of 
the block-diagonal operator 
$ \bar \om_\e \cdot \partial_\vphi - J A_\ind $, 
we obtain  an approximate inverse of
$$ 
\bar \om_\e \cdot \partial_\vphi - J (A_\ind + \rho_\ind) 
$$ 
in  \eqref{ome-Arho}, 
taking into account  the small residual coupling term $ \rho_\ind $.
Since $ \rho_\ind $ satisfies \eqref{re-ind}, 
we conclude by using a Neumann series perturbative argument,  see Section \ref{sec:9.2}. 

In conclusion, 
after all this analysis, 
going back to the original coordinates, we  finally prove,
in Proposition \ref{prop:inv-ap-vero}, 
 the existence of an approximate right inverse of the quasi-periodic 
 Hamiltonian operator $ {\cal L}_\om  (i_n) $ in \eqref{lin:restr}, and thus of
 $ d_{i} {\cal F}(\l; i_n) $, for most values of $ \l $'s, satisfying tame estimates.   
 This enables to implement a Nash-Moser iterative scheme (Chapter \ref{sec:NM})
 which proves Theorem \ref{thm:NM}  and therefore Theorem  \ref{thm:main}.

\medskip

\noindent
{\bf Sobolev regularity thresholds.} Along the monograph we shall use four 
Sobolev indices
$$
s_0 \ll s_1 \ll s_2 \ll s_3 \, . 
$$
The first index 
$$ 
s_0 > (\es + d )/2 
$$ 
is fixed in \eqref{algebra-Sobolev}  to have the algebra and interpolation properties for all the Sobolev spaces 
$ {\mathcal H}^s $, $ s \geq s_0 $, defined in \eqref{Sobo:sp1}. We prove 
these properties
 in Section \ref{sec:35}. 

The index 
$$ 
s_1 \gg s_0 
$$ 
is the one required 
 for the multiscale Proposition \ref{propmultiscale} (see hypothesis 1 in  Definition \ref{definition:Xr})
 to have the sufficient off-diagonal decay
 to apply the multiscale step Proposition \ref{propinv}, 
see in particular  \eqref{s1} and assumption $(H1)$. It is the Sobolev threshold which 
defines the ``good" matrices in Definition \ref{goodmatrix}. 
Another condition of the type $ s_1 \gg s_0 $ appears  in the proof of Lemma \ref{loss-regularity-low}.

 Then the Sobolev index 
 $$ 
 s_2 \gg s_1 
 $$ 
 is used in the splitting step Proposition \ref{prophomeq},
 see \eqref{choice s2 s3}. This proposition is based on a Nash-Moser
 iterative scheme where $ s_2 $ represents a `high-norm", 
 see  \eqref{rho-R0:small}.

The largest index 
$$ 
s_3 \gg s_2 
$$ 
is finally used for the convergence of the Nash-Moser nonlinear iteration
in Chapter  \ref{sec:NM}. Note that in Theorem \ref{thm:NM-ite}
the divergence of the approximate solutions $ i_n $ 
in the high norm $ \| \ \|_{\Lip, s_3} $ is under control. We require in particular \eqref{choice s2 s3} and,
in Section \ref{sec:NM-ite}, also stronger largeness conditions for $ s_3 - s_2 $. 
Along the iteration  
 we shall verify (impose) that the Sobolev norms of the approximate quasi-periodic solutions remain bounded in 
$\| \ \|_{\Lip, s_1} $ and $ \| \ \|_{\Lip, s_2} $ norms (so that the assumptions of Propositions \ref{propmultiscale},
\ref{prop-cruciale} and  \ref{prophomeq} are fulfilled).

\chapter{Hamiltonian formulation}\label{sec:Ham}

In this chapter we write the nonlinear wave equation \eqref{NLW2} as a
first order Hamiltonian system in action-angle and normal variables. These 
coordinates are convenient for the proof of Theorem \ref{thm:main}. 
In Section \ref{sec:fre-ad} we  provide estimates about the 
measure of the admissible Diophantine directions $ \bar \om_\e $ of the frequency vector $ \omega = (1 + \e^2 \l) \bar \om_\e  $ of the quasi-periodic solutions
of Theorem \ref{thm:main}. 

\section{Hamiltonian form of NLW}

We  write the second order nonlinear wave equation \eqref{NLW2} as the first order system
\be\label{NLWS}
\begin{cases}
u_t = v  \cr 
v_t =  \Delta u - V(x) u  - \e^2 g(\e, x, u )
\end{cases} \, ,
\ee
which is the Hamiltonian system
\be\label{calJ}
\partial_t (u,v) = X_H (u,v) \, , \qquad
X_H (u, v) := J \nabla_{L^2} H(u,v) \, , 
\ee
where\index{Symplectic matrix} 
\be\label{sympl-J}
J := \begin{pmatrix}
0  &   {\rm Id}   \\
 - {\rm Id} & 0     \\   
\end{pmatrix} \, , 
\ee
is the symplectic matrix, and 
$ H $ is the Hamiltonian 
$$
H (u,v) := \int_{\T^d} \frac{v^2}{2} + \frac{ (\nabla u)^2}{2} + V(x)  \frac{u^2}{2}  + \e^2 G(\e, x, u)\, dx  \, .
$$
In the above definition of $H$, $G$ denotes a primitive of $g$ in the variable $ u $:
\be \label{svilG}
G(\e, x, u ) := \int_0^u g(\e, x, s) \, ds \, , \quad (\pa_u G)(\e, x, u ) = g(\e, x, u)  \, .
\ee 
For the sequel it is convenient to highlight also the fourth order term  of the nonlinearity $ g(x, u)$ in \eqref{nonlinearity}, i.e. 
writing
\be\label{nonlinearity1}
g(x, u) = a(x) u^3 + a_4 (x) u^4 + g_{\geq 5} ( x, u  )  \, , \quad g_{\geq 5}  ( x, u  ) = O( u^5 )  \, ,
\ee
so that  the rescaled nonlinearity $ g(\e, x, u ) $ in \eqref{nonlinearity:gep} has the expansion  
\be
\begin{aligned}\label{nonlinearity:gep1}
 g(\e, x, u ) := \e^{-3} g(x, \e u) & = a(x) u^3 + \e a_4 (x) u^4 + \e^2 {\mathfrak r} ( \e, x, u  )  \, , \\
& \qquad \qquad  {\mathfrak r} ( \e, x, u  ) := \e^{-5} g_{\geq 5}  ( x, \e u  ) \, , 
 \end{aligned}
\ee
and its primitive in \eqref{svilG} has the form 
\be \label{svilG1}
G(\e, x, u ) = 
 \frac{1}{4 } a(x)  u^4 + \frac{\e}{5}  a_4 (x) u^5 + 
 \e^2 {\mathfrak R} ( \e, x, u  ) 
\ee
where  $  ( \pa_u {\mathfrak R}) ( \e, x, u  ) = {\mathfrak r} ( \e, x, u  ) $.  

The phase space of \eqref{calJ} is a dense subspace of the real Hilbert space
$$
L^2 (\T^d) \times L^2 (\T^d) \, , \quad  L^2 (\T^d) := L^2 (\T^d, \R) \, ,  
$$
endowed  with 
the standard constant symplectic $ 2 $-form\index{Symplectic $ 2$-form}
\be\label{2fo}
\Omega ( (u,v) , (u',v') ) :=  (J (u,v) , (u',v') )_{L^2 \times L^2} =  
(v , u' )_{L^2} -  ( u , v' )_{L^2} \, . 
\ee
Notice that the Hamiltonian vector field  $ X_H $  is  characterized by the relation 
$$ 
\Omega ( X_H, \cdot) = - d H \, .
$$ 
System \eqref{NLWS} is reversible\index{Reversible vector field} with respect to the involution\index{Involution}
\be\label{invol}
S (u,v) := (u, - v ) \, , 	\quad S^2 = {\rm Id} \, , 
\ee
namely (see Appendix \ref{App0-1})
$$
X_H \circ S = - S \circ X_H  \, , \qquad {\rm equivalently } 
\qquad H \circ S = H \, . 
$$
Notice that the above equivalence is due to the fact that 
the involution $ S$ is antisymplectic, namely the pull-back 
$$ 
S^* \Omega = - \Omega \,    .
$$ 
Let us write \eqref{NLWS} in a more symmetric form. 
We  first introduce 
$ D_V := \sqrt{- \Delta + V(x)} $ as the closed (unbounded) 
 linear  operator of $L^2(\T^d)$ defined  by 
\be\label{def:DV}
D_V \Psi_j := \sqrt{- \Delta + V(x)} \Psi_j := \mu_j \Psi_j \, , \quad \forall j \in \N \, , 
\ee
where  $   \{ \Psi_j (x),  j \in \N  \} $  is the 
 orthonormal basis of  $ L^2 (\T^d) $ 
  formed by the eigenfunctions of $ - \Delta +  V(x) $ in \eqref{auto-funzioni}.
 
Under the symplectic transformation 
\be\label{D14}
q = D_V^{\frac12} u \, , \quad p = D_V^{-\frac12} v \, , 
\ee
the Hamiltonian system \eqref{NLWS} becomes 
\be\label{new-HS-qp}
\begin{cases}
q_t = D_V p \cr
p_t = - D_V q - \e^2  D_V^{-\frac12} g(\e, x, D_V^{-\frac12} q ) \, , 
\end{cases}
\ee
which is the Hamiltonian system 
\be\label{HSinK}
\partial_t (q,p) = J \nabla_{L^2} {\mathtt K}(q,p)  
\ee
with transformed Hamiltonian
\be\label{K}
{\mathtt K} (q,p) :=  \int_{\T^d} \frac{(D_V^{\frac12} p)^2}{2} 
+  \frac{(D_V^{\frac12} q)^2}{2} 
+  \e^2  G(\e, x,   D_V^{-\frac12} q   ) \, dx \, .
\ee
This Hamiltonian system is still reversible with respect to  $ S $ defined in \eqref{invol}, i.e.
$ {\mathtt K} \circ S = {\mathtt K} $, i.e. $ {\mathtt K} $ is even in $ p $.

\section{Action-angle and ``normal" variables}\label{sec:AA}

We look for quasi-periodic solutions of the Hamiltonian system \eqref{HSinK} which are mainly Fourier supported on 
the tangential sites $ {\mathbb S } \subset \N $. Thus 
we decompose 
\be
\begin{aligned}
\label{orthogonal-deco-SS}
& \qquad \qquad  L^2 (\T^d, \R) \times L^2 (\T^d, \R)
 = H_{\mathbb S} \oplus H_{\mathbb S}^\bot \,  ,   \\  
&  H_{\mathbb S} := \Big\{   (q(x), p(x)) =  \sum_{j \in {\mathbb S}} (q_j, p_j) \Psi_j (x)  \, , \quad (q_j , p_j) \in \R^2 \Big\} 
\end{aligned}
\ee
splitting the canonical variables $ (q,p) $ into  {\it tangential}\index{Tangential variables}  and  {\it normal} components\index{Normal variables} 
\be\label{tang-norm}
(q,p) = \sum_{j \in {\mathbb S}} (q_j, p_j ) \Psi_j (x) + (Q,P) 
\ee
where  $ (q_j,p_j) \in \R^2 $ and  $ (Q,P)  \in H_{\mathbb S}^\bot  $. 
Then  we introduce usual action-angle\index{Action-Angle variables} variables on the tangential sites by setting
\be\label{AA}
(q_j,p_j) := \sqrt{2 I_j} \big( \cos \teta_j, - \sin \teta_j \big) \, , \quad \forall j \in {\mathbb S } \, .   
\ee
The symplectic form\index{Symplectic $ 2$-form} \eqref{2fo} then becomes (recall that $ \Psi_j (x) $ are  $ L^2 $-orthonormal)
\be\label{2form}
{\cal W} := (d I \wedge d \theta) \, \oplus \,  \Omega  
\ee
where, for simplicity, we still denote by $ \Omega := \Omega_{|H_{\mathbb S}^\bot} $  the restriction of the $ 2$-form $ \Omega $, defined in \eqref{2fo}, to the symplectic subspace $ H_{\mathbb S}^\bot $.

By 
\eqref{tang-norm}, \eqref{AA},  the Hamiltonian in \eqref{K} 
then becomes, recalling that $ \bar \mu \in \R^\es $ is  the unperturbed tangential frequency vector 
defined in \eqref{unp-tangential},  
\begin{align}
K(\teta, I, Q, P) & = 
\bar \mu \cdot I +  
 \frac12 \int_{\T^d}  (D_V^{\frac12} Q)^2 + (D_V^{\frac12} P)^2  \, dx  \label{primalinea0} \\
  & \quad + \e^2
 \int_{\T^d} 
 G \Big(\e, x,    \sum_{j \in \mathbb S}  \mu_j^{-\frac12} \sqrt{2 I_j} \cos \teta_j \Psi_j (x)  + D_V^{-\frac12} Q   \Big) dx 
\,  . \nonumber 
\end{align}
For $ \e = 0 $, the Hamiltonian  system generated by \eqref{primalinea0} 
admits the continuous family of quasi-periodic solutions  
$$
\teta (t) = \teta_0 + \bar \mu t \, , \ I (t) = \xi \, ,   \ Q (t) = P (t) =  0 \, ,
$$
parametrized by the unperturbed tangential ``actions" $ \xi := (\xi_j )_{j \in \mathbb S} $, $ \xi_j >  0 $.  
The aim is to prove their persistence, being just slightly deformed, 
for $ \e $ small enough, for ``most" values of $ \xi $, and with a frequency close to
$ \bar \mu $.
 
Then we  introduce nearby coordinates 
by the 
symplectic transformation 
\be\label{tra-xi-y}
I = \xi + y \, ,
\ee 
and, substituting in 
\eqref{primalinea0}, 
we are reduced to study the parameter dependent family of Hamiltonians 
(that for simplicity we denote with the same letter $ K $) 
\begin{align}
K(\teta, y, Q, P, \xi)  & = c  +
\bar \mu \cdot y  +  \frac12  (D_V  Q, Q )_{L^2} + \frac12  (D_V  P, P)_{L^2}   + \e^2 R(\teta, y, Q, \xi)  \label{primalinea}
\end{align}
where 
\begin{align}\label{restoNN}
& R(\teta, y, Q, \xi)  :=   \int_{\T^d} G \big(\e, x,  v(\teta, y, \xi)  + D_V^{-\frac12} Q   \big) \, dx  
\end{align}
and 
\be\label{defv}
v(\teta, y, \xi) :=  \sum_{j\in \mathbb S}  \mu_j^{-\frac12} \sqrt{2 (\xi_j + y_j)} \cos \teta_j \Psi_j (x)  \, .
\ee
The  phase space of \eqref{primalinea} is now 
$$
\T^\es \times \R^\es \times H_{\mathbb S}^\bot  \ni (\teta, y, z)\, , \quad z := (Q,P) \in H_{\mathbb S}^\bot \, , 
$$
endowed with the symplectic structure (see \eqref{2form})
\be\label{2form-y}
{\cal W} := (d y \wedge d \theta) \, \oplus \,  \Omega  \, , 
\ee
 so that 
the Hamilton equations generated by \eqref{primalinea} 
have  the  form  
\be\label{HS0}
\begin{cases}
\dot \teta = \partial_y K (\teta, y, z, \xi)   \cr 
\dot y  = - \partial_\teta K (\teta, y,z, \xi) \cr
\dot z = J \nabla_z K (\teta, y,z, \xi) 
\end{cases}
\ee
or, in expanded form, 
\be\label{HS00exp}
\begin{cases}
 \dot \teta - \bar \mu - \e^2 \partial_y R (\theta, y, Q, \xi )  = 0   \cr 
\dot y   +  \e^2 \partial_\teta R (\theta, y, Q, \xi ) = 0 \cr
\big( \pa_t  -  J D_V \big) (Q, P) + \e^2 \big(0, \nabla_Q R (\theta, y, Q, \xi ) \big) = 0 \, .
\end{cases}
\ee
By \eqref{restoNN} and \eqref{svilG}, we have that $ \pa_y R := ( \pa_{y_m} R)_{m=1, \ldots, \es} \in \R^\es $ has components 
\be\label{gradyR}
\pa_{y_m} R ( \theta, y, Q, \xi)  = \int_{\T^d} g(\e, x, v(\theta, y, \xi) +  D_V^{- \frac12} Q ) 
 \frac{\mu_m^{-\frac12}}{\sqrt{2(\xi_m+ y_m)}} \cos (\theta_m) \Psi_m (x) \, dx 
\ee
and 
\be\label{gradRQ}
\nabla_Q R ( \theta, y, Q, \xi)  = D_V^{- \frac12} g \big( \e, x, v(\theta, y, \xi) +  D_V^{- \frac12} Q \big) \, .
\ee
We note that the $ 2 $-form $ {\cal W} $ in  \eqref{2form-y}  is exact, i.e.
$ {\cal W} = d \form  $
where $ \form $ is the Liouville $1$-form\index{Liouville $ 1$-form}
\be  \label{Lambda 1 form}
\form_{(\theta, y, z)} [ \widehat \theta, \widehat y, \widehat z] =  \sum_{j=1, \ldots, \es}  
y_j \widehat \theta_j +  \frac12 (J z, \widehat z)_{L^2_x} \, . 
\ee
The Hamiltonian system \eqref{HS0} is reversible with respect to the involution\index{Involution}
\be \label{involuzione tilde rho}
\tilde S ( \teta, y, Q,P ) := ( - \teta, y, Q, -P) 
\ee
which is nothing but \eqref{invol} in the variables \eqref{tang-norm}-\eqref{AA}. This means that 
$$
K \circ \tilde S = K \, , \quad K ( - \teta, y, Q, - P ) = K ( \teta, y, Q,P ) \, ,
$$
and the Hamiltonian $ R $ in \eqref{restoNN} satisfies 
\be\label{revers:R}
R ( - \teta, y, Q, \xi ) = R ( \teta, y, Q, \xi ) \, . 
\ee

\section{Admissible Diophantine directions $\bar \om_\e  $}\label{sec:fre-ad}

As\index{Admissible frequencies} explained in the introduction, we look for quasi-periodic solutions of \eqref{HS0}-\eqref{HS00exp} with  frequency vector
$$ 
\om = (1 + \e^2 \l ) \bar \om_\e 
$$
restricted to a  fixed line, spanned by $ \bar \om_\e = \bar \mu + \e^2 \dom $ (see \eqref{frequencycolinear}-\eqref{def omep}), 
which has to  satisfy the Diophantine conditions\index{Diophantine vector}  \eqref{dioep}-\eqref{NRgt1}. 
We now prove that for ``most" vector  $ \zeta \in \R^\es $ these conditions are satisfied. 
 
\begin{lemma}\label{lemma:rho-dioph}
Assume \eqref{diop} and \eqref{NRgamma0}. Then there exists 
a subset $ B_\e \subset  \Ab ([1,2]^\es )  $ (where $\Ab $ is the invertible twist matrix in \eqref{def:AB}) with measure $| B_\e | \leq \e $, 
such that all the  vectors 
$$ 
\bar \om_\e = \bar {\mu} + \e^2 \dom \, , \quad 
\dom \in \Ab ([1,2]^\es )  \setminus B_\e \, ,
$$  
satisfy the Diophantine conditions
\eqref{dioep}-\eqref{NRgt1} with $\gamma_1, \tau_1$  defined in \eqref{def:tau1}.  
\end{lemma}

\begin{pf}
We first verify that most vectors $ \bar \om_\e = \bar {\mu} + \e^2 \dom $, $ \dom \in \Ab ([1,2]^\es )  $, satisfy \eqref{dioep}. 
Since $ \bar \mu $ is $ (\g_0, \t_0 )$-Diophantine, i.e. \eqref{diop} holds, then, for all
$ \ell \in \Z^{\es} \setminus \{ 0 \}  $,  
$$
| (\bar \mu + \e^2 \dom) \cdot \ell  | \geq  
\frac{\g_0}{ \langle \ell \rangle^{\tau_0}} - \e^2 C | \ell | \geq \frac{\g_0}{2 \langle \ell \rangle^{\tau_0}} \, , 
\quad \forall |\ell | \leq \Big(\frac{\g_0}{2 C \e^2}\Big)^{1/(\t_0+1)} \, . 
$$
Then it remains to estimate the measure of
\be\label{setBad}
{\cal B}_\e := {\mathop \bigcup}_{|\ell | > \Big(\frac{\g_0}{2 C \e^2}\Big)^{1/(\t_0+1)}} {\cal R}_\ell 
\ee
where
$$
{\cal R}_\ell := \Big\{ \dom \in \Ab ([1,2]^\es )  \ : 
\ |(\bar \mu + \e^2 \dom) \cdot \ell | \leq \frac{\gamma_1}{  \langle \ell \rangle^{\tau_1}} \, , 
\ \gamma_1 = \frac{\gamma_0}{2} \Big\} \, .
$$
Since the derivative
$ \frac{\ell}{|\ell|} \cdot \pa_{\dom} \big( (\bar \mu + \e^2 \dom) \cdot \ell \big) = \e^2 |\ell | $, then 
the measure $ | {\cal R}_\ell | \lesssim \frac{ \gamma_0}{\e^2  \langle \ell \rangle^{\tau_1+1}} $. 
Therefore
$$
| {\cal B}_\e | \lesssim  \frac{ \gamma_0}{\e^2 } {\mathop \sum}_{|\ell | > 
\Big(\frac{\g_0}{2 C \e^2}\Big)^{1/(\t_0+1)}}  \frac{1}{ \langle \ell \rangle^{\tau_1+1}} 
\leq C(\g_0) \e^{2\frac{(\t_1 +1- \es)}{\t_0+1} - 2} \leq \e 
$$
by \eqref{def:tau1}.

We now consider  the quadratic Diophantine condition \eqref{NRgt1}.
Let $ M := M_p $ be the $ (\es \times \es)-${\it symmetric} matrix such that 
$$
\sum_{1 \leq i \leq j \leq \es} 
\omega_i \omega_j p_{ij}  = M \om \cdot \om \, , \quad \forall \om \in \R^\es \, .
$$
The symmetric matrix $ M $ has coefficients   
\be\label{defMS}
M_{ij} := \frac{p_{ij}}{2} (1 + \d_{ij})  \, , \  \forall   1 \leq i \leq j \leq \es \, , \quad {\rm and} \quad M_{ij} = M_{ji} \, .
\ee
We want to prove that for most $ \dom \in \Ab ([1,2]^\es )  $ the vector 
$ \bar \om_\e  = \bar \mu + \e^2 \dom $ 
satisfies the non-resonance condition
$$ 
\big| n + M \bar \om_\e \cdot \bar \om_\e   \big| \geq \g_1 \langle p \rangle^{- \t_1} \, ,
\quad 
\forall (n, p) \in \Z \times \Z^{\frac{\es(\es+1)}{2}} \setminus \{(0,0)\}  \, .  
$$ 
Then we write  
\begin{align}
 n + M \bar \om_\e \cdot \bar \om_\e 
= 
 n + M \bar \mu  \cdot \bar \mu  + 2 \e^2 M \dom \cdot \bar \mu + \e^4 M  \dom \cdot \dom   \, . 
 \label{nr-voluta}
\end{align}
We first note that, by \eqref{NRgamma0} and
$ | M_{ij} |\lesssim \langle p \rangle $ by \eqref{defMS}, then, for all $ | \dom | \leq 1 $,  
$$
\big| n + M \bar \mu \cdot \bar \mu +  2 \e^2 M \dom \cdot \bar \mu + \e^4 M  \dom \cdot \dom  \big|
\geq 
| n + M \bar \mu \cdot \bar \mu | -  \e^2 C \langle p \rangle   \geq  \frac{\g_0/2}{\langle p \rangle^{\t_0}} 
$$
if $
\Big( \frac{\g_0 }{ 2 \e^2 C } \Big)^{1/ (\tau_0 + 1)}\geq  | p | $. 
Thus we erase values of $ \dom \in \Ab ([1,2]^\es )  $ only when 
\be\label{p-large0} 
|p| > \Big( \frac{\g_0 }{ 2 \e^2 C } \Big)^{1/ (\tau_0 + 1)} \, .
\ee
Since $ M $ is symmetric there is an orthonormal basis $ V := (v_1, \ldots, v_k) $ of eigenvectors 
 of  $M$  with real eigenvalues
$ \l_k :=  \lambda_k (p) $, i.e. $ M v_k = \l_k v_k $. 
Under the isometric change of variables $ \dom = V y $ we have to estimate 
\be\label{Rpy}
\begin{aligned}
|{\cal R}_{n,p} | = \Big|  \Big\{ &  y \in \R^\es \, , \, | y | \leq 1  \ :  \\
& \Big| \,
n + M \bar \mu  \cdot \bar \mu  + 2 \e^2  V y \cdot M \bar \mu + \e^4  \sum_{1 \leq i  \leq \es} \l_k y_k^2
 \Big| < \frac{\gamma_0}{2 \langle p \rangle^{\tau_1}} \Big\}\Big| \, . 
 \end{aligned}
\ee
Since
$ M^2 v_k = \l_k^2 v_k $, $ \forall k =1, \ldots, \es $,  we get 
$$
\sum_{k=1}^\es \l_k^2 = {\rm Tr} (M^2) = \sum_{i, j =1}^\es M_{ij}^2 
\stackrel{\eqref{defMS}} \geq \frac{| p |^2 }{ 2}  \, . 
$$
Hence there is an index $ k_0 \in \{1, \ldots , \es \} $ 
such that $ | \l_{k_0} | \geq  |p| / \sqrt{2\es} $ and  the derivative
\be\label{Rpy1}
\begin{aligned}
\Big| \partial^2_{y_{k_0} } \Big( 
n + M \bar \mu  \cdot \bar \mu  + 2 \e^2  V y \cdot M \bar \mu + \e^4  \sum_{1 \leq i  \leq \es} \l_k y_k^2 
\Big) \Big| & =  \e^4 |2 \l_{k_0}| \\ 
& \geq \e^4
\sqrt{2} \, |p| \slash \sqrt{\es}  \, . 
\end{aligned}
\ee
As a consequence of \eqref{Rpy} and \eqref{Rpy1} we deduce  the measure estimate
$$ 
|{\cal R}_{n,p} | \lesssim \e^{-2}  \sqrt{\frac{\gamma_0}{ \langle p \rangle^{\tau_1 +1}} }  \, .  
$$
Recalling \eqref{p-large0}, and since $ {\cal R}_{n,p} = \emptyset $ if $ | n | \geq C \langle p \rangle $,
 we have 
$$
\begin{aligned}
\Big| \bigcup_{n, p \in \Z^{ \frac{\es (\es+1)}{2}} 
\setminus \{0 \}}  \!\! \!\! {\cal R}_{n,p}  \Big| 
\lesssim  \sum_{|p| > \Big( \frac{\g_0 }{ 2 \e^2 C } \Big)^{1/ (\tau_0 + 1)}}  \!\! \!\!  
  \e^{-2} \langle p \rangle \sqrt{\frac{\gamma_0}{ \langle p \rangle^{\tau_1 +1}} } 
 & \lesssim_{\g_0} 
 \e^{ \big[ \frac{\tau_1-1}{\tau_0+1} - \frac{\es(\es+1)}{\tau_0+1} - 2\big] } \\ 
 & \leq \e 
 \end{aligned}
$$
for $ \e $ small, by \eqref{def:tau1}.
\end{pf}

\begin{remark}
The measure  of the set  $ | B_\e | \leq \e $ is 
smaller than $ \e^p $, for any $ p $, at the expense of taking a larger Diophantine exponent $ \tau_1 $. 
We have written $ | B_\e | \leq \e $  for definiteness. 
\end{remark}

We finally notice that, for $ \omega = (1 + \e^2 \l ) \bar \om_\e  $ with  
a Diophantine vector $\bar \om_\e $ satisfying  \eqref{dioep}, for any  
zero average function $ g (\vphi) $ we have that  the function 
$  (\om \cdot \pa_\vphi)^{-1} g $, 
defined in \eqref{op-inv-KAM}, satisfies 
\be\label{diof-est}
\| (\om \cdot \pa_\vphi)^{-1} g \|_{\Lip,s} \leq C \g_1^{-1}  \| g \|_{\Lip,s + \tau_1}   \, . 
\ee

\chapter{Functional setting}\label{Ch:3}

In this chapter we collect 
all the  properties of the phase spaces,
linear operators, norms, interpolation inequalities used through the monograph.
Of particular importance for proving Theorem \ref{thm:main} 
is the result of Section \ref{sec:offDV}.

\section{Phase space and basis}

The phase space of the nonlinear wave equation
 \eqref{new-HS-qp} is a  dense subspace of 
the real Hilbert space\index{Phase space} 
\be\label{def:H}
H := L^2 (\T^d) \times L^2 (\T^d) \, , \quad  L^2 (\T^d) := L^2 (\T^d, \R) \, . 
\ee
We shall denote for convenience an element of $ H $ either as
$ h = ( h^{(1)},  h^{(2)}) $ a row, either as a column
$
h = \begin{pmatrix}
h^{(1)}    \\
h^{(2)}     \\   
\end{pmatrix}  
$ 
with components $ h^{(l)} \in L^2 (\T^d) $, $ l = 1, 2 $. 

In the exponential basis any function of $ H $ can be decomposed as
\be\label{basis-exp}
(q(x), p(x)) = \sum_{j \in \Z^d} (q_j, p_j) e^{\ii j \cdot x} \, , 
\quad q_{-j} = \ov{q_j}  \, , \ p_{-j} = \ov{ p_j}   \, . 
\ee
We will also use the orthonormal basis 
$$  
 \{ \Psi_j (x),  j \in \N  \} 
 $$
 of  $ L^2 (\T^d) $ formed by the eigenfunctions of $ - \Delta +  V(x) $ defined in \eqref{auto-funzioni} with eigenvalues 
$ \mu_j^2 $.
We then consider the Hilbert spaces
\be\label{spaces-cal-Hs}
{\mathtt H}^s_x := \Big\{  u := \sum_{j \in \N} u_j \Psi_j \, : \, \| u \|_{{\mathtt H}^s_x}^2 := 
\big( (- \Delta + V(x))^s u, u \big)_{L^2} = 
\sum_{j \in \N} \mu_j^{2s} |u_j|^2  < \infty \Big\} \, . 
\ee
Clearly  $ {\mathtt H}^0_x = L^2 $. 
Actually, for any $ s \geq 0 $,  the ``spectral" norm  $ \| u \|_{{\mathtt H}^s_x} $ 
is equivalent to the usual Sobolev norm   
\be\label{equiv-norms-s}
\| u \|_{{\mathtt H}^s_x} \simeq_s \| u \|_{H^s_x} :=  \big( (- \Delta)^s u, u \big)_{L^2} + (u, u)_{L^2} \big)^{\frac12} \, . 
\ee
For $ s \in \N $, the equivalence \eqref{equiv-norms-s} can be directly proved noting that $ V(x) $ is a lower order perturbation of the Laplacian
$ - \Delta $. Then, for $ s \in \R \setminus \N $, $ [s] < s <  [s] + 1 $, the equivalence
\eqref{equiv-norms-s}  follows by the classical  interpolation result 
stating that the Hilbert space $ {\mathtt H}^s_x  $ in \eqref{spaces-cal-Hs}, respectively the Sobolev space $ H^s_x $, 
 is the interpolation space between 
$ {\mathtt H}^{[s]}_x $ 
and $ {\mathtt H}^{[s]+1}_x $, respectively between $ H^{[s]}_x $  and $H^{[s]+1}_x $, and the
equivalence \eqref{equiv-norms-s}  for integers $ s $. 
\\[1mm]
{\bf Tangential and normal subspaces.  } 
Given the finite set  $ {\mathbb S} \subset \N $, we 
consider the $ L^2 $-orthogonal 
decomposition of the phase space in tangential and normal subspaces
as in  \eqref{orthogonal-deco-SS},  
$$
H = H_{\mathbb S} \oplus H_{\mathbb S}^\bot \, , 
$$
where
$$
H_{\mathbb S} = \Big\{   (q(x), p(x)) =  \sum_{j \in {\mathbb S}} (q_j, p_j) \Psi_j (x)  \, , \quad (q_j , p_j) \in \R^2 \Big\} \, .
$$
In addition, recalling the disjoint splitting
$$ 
\N = {\mathbb S} \cup {\mathbb F} \cup {\mathbb G} 
$$ 
where 
$ {\mathbb F} , {\mathbb G} \subset \N $ are  defined in \eqref{taglio:pos-neg-0}-\eqref{taglio:pos-neg},  
we further decompose the normal subspace $H_{\mathbb S}^\bot $ as
 \be\label{orthogonal-deco-FS}
H_{\mathbb S}^\bot = H_{\mathbb F} \oplus H_{{\mathbb G}} 
\ee
where 
\be\label{subspaceG}
\begin{aligned}
& 
H_{\mathbb G} := \Big\{   (Q(x), P(x)) \in H \ : \ (Q,P) \bot H_{\mathbb S} \, , \ (Q,P) \bot H_{\mathbb F} \Big\} \, , \\
& H_{\mathbb F} := \Big\{   (Q(x), P(x)) \in H \ : \ (Q,P) \bot H_{\mathbb S} \, , \ (Q,P) \bot H_{\mathbb G} \Big\}  \, . 
\end{aligned}
\ee
Thus 
\be\label{H:SFG}
H = H_{\mathbb S} \oplus H_{\mathbb S}^\bot 
= H_{\mathbb S} \oplus H_{\mathbb F} \oplus H_{{\mathbb G}} 
= H_{\mathbb S \cup \mathbb F} \oplus H_{{\mathbb G}}   \, . 
\ee
Accordingly  we denote by $ \Pi_{\mathbb S} $, $ \Pi_{\mathbb F} $, $ \Pi_{{\mathbb G}} $,  
$ \Pi_{\mathbb S \cup \mathbb F} $, the orthogonal $ L^2 $-projectors
on $ H_{\mathbb S} $, $ H_{\mathbb F} $, $ H_{\mathbb G} $, $ H_{\mathbb S \cup \mathbb F} $. We define
$$ 
\Pi_{\mathbb S}^\bot := {\rm Id} - \Pi_{\mathbb S} = \Pi_{\N \setminus \mathbb S} 
$$
and similarly for the other subspaces.

Decomposing  the finite dimensional space
\be\label{deco-F-Fj}
H_{\mathbb F} = \oplus_{j \in {\mathbb F}} H_j   \, , 
\ee
 we denote by $ \Pi_j $ the  $ L^2 $-projectors onto $ H_j $. 
In each real subspace $ H_j $, $ j \in {\mathbb F}$,  we take the basis  
$ ( \Psi_j(x) ,0) \, ,  (0, \Psi_j(x) )  $,  namely we represent 
\be\label{Fj-eigenfunctions}
H_j = \Big\{  
q (\Psi_j(x) ,0) + p ( 0, \Psi_j(x)) \, , \  q , p \in \R \Big\} \, . 
\ee
Thus $ H_j $ is isometrically isomorphic to $ \R^2 $. 
\\[1mm]
{\bf Symplectic operator $ J $.}
We define the linear {\it symplectic} operator\index{Symplectic matrix} $ J \in {\cal L}( H ) $  as 
\be\label{symplec-op}
J 
 \begin{pmatrix}
Q (x)    \\
P(x)     \\   
\end{pmatrix} :=
\begin{pmatrix}
P(x)    \\
- Q(x)     \\   
\end{pmatrix} \, , 
\quad 
J = \begin{pmatrix}
0  &   {\rm Id}   \\
 - {\rm Id} & 0     \\   
\end{pmatrix}  \, .
\ee
The symplectic operator $ J $ leaves invariant the symplectic subspaces $  H_{\mathbb S} $,
$  H_{\mathbb S}^\bot  $,  $ H_{\mathbb F} $, $  H_{{\mathbb G}} $. 
For simplicity of notation we  shall still denote by $ J $ the restriction  of the symplectic operator 
$$
J := J_{|H_{\mathbb S}}  \in {\cal L}(H_{\mathbb S}) \, , \  
J := J_{|H_{\mathbb S}^\bot}  \in {\cal L}(H_{\mathbb S}^\bot) \, , \  
J := J_{|H_{\mathbb F}}  \in {\cal L}(H_{\mathbb F}) \, , \  
J := J_{|H_{{\mathbb G}}}  \in {\cal L}(H_{{\mathbb G}}) \, .
$$
The symplectic operator $ J $ 
is represented, in the basis of the exponentials (see \eqref{basis-exp})
$$ 
\big\{ (e^{\ii j \cdot  x}, 0), (0, e^{\ii j \cdot x}) \, ,  j \in \Z^d \big\}  \, , 
$$ 
by the matrix 
$$
J = {\rm Diag}_{j \in \Z^d}
\begin{pmatrix}
0  &   1   \\
 -1 & 0     \\   
\end{pmatrix}  \, . 
$$
Note also that the symplectic operator $ J $ leaves invariant each  subspace  $ H_j $, $ j \in {\mathbb F} $, 
and it is represented,  in the eigenfunction basis \eqref{Fj-eigenfunctions},  by the same symplectic matrix 
$
\begin{pmatrix}
0  &   1   \\
 -1 & 0     \\   
\end{pmatrix} 
$.

\section{Linear operators and matrix representation}

According to the decomposition $ H_{\mathbb S}^\bot = H_{\mathbb F} \oplus H_{\mathbb G} $ in \eqref{orthogonal-deco-FS}
a  linear operator $ A $ of $ H_{\mathbb S}^\bot $ can be represented\index{Matrix representation of linear operators}  by a matrix of operators as
\be\label{decoFGsotto}
\begin{aligned}
&
\begin{pmatrix}
 A^{\mathbb F}_{\mathbb F} & A_{\mathbb F}^{\mathbb G} \\
A^{\mathbb F}_{\mathbb G} &  A^{\mathbb G}_{\mathbb G}
\end{pmatrix} \, , \\  
& A^{\mathbb F}_{\mathbb F} := \Pi_{\mathbb F} A_{|H_{\mathbb F}} \, , \
A_{\mathbb F}^{\mathbb G} := \Pi_{\mathbb F} A_{|H_{\mathbb G}} \, , \
A^{\mathbb F}_{\mathbb G} := \Pi_{\mathbb G} A_{|H_{\mathbb F}} \, , \
A^{\mathbb G}_{\mathbb G} := \Pi_{\mathbb G} A_{|H_{\mathbb G}} \, .
\end{aligned}
\ee
Moreover the decomposition \eqref{deco-F-Fj} induces the splitting
\be\label{op-a-deco-sez}
{\cal L}( H_{\mathbb F}, H ) = \oplus_{j \in {\mathbb F} }  {\cal L}( H_j, H ) \, ,
\ee
namely  a linear operator   $ A \in {\cal L}( H_{\mathbb F}, H ) $ can be written as 
\be\label{a-deco-sez}
A = (a^j)_{j \in {\mathbb F}}  \, , \quad a^j := A_{| H_j }  \in {\cal L}( H_j, H ) \, .
\ee
In each space $ {\cal L}( H_j, H ) $ we define the scalar product 
\be\label{sc-prod-L}
\langle a, b \rangle_0 := {\rm Tr} ( b^* a )  \, , \quad a, b \in {\cal L}( H_j, H )  \, ,  
\ee
where $ b^* \in  {\cal L}(  H , H_j ) $ denotes the adjoint of $ b $ with respect to the scalar product in $ H $.
Note that $ b^* a \in {\cal L} ( H_j, H_j ) $ 
is represented, in the basis \eqref{Fj-eigenfunctions},  by the 
$ 2 \times 2 $ real matrix whose elements are $ L^2 $ scalar products
\be\label{matrix-adjoint}
\left(
\begin{array}{cc}
(  b(\Psi_j,0), a(\Psi_j,0))_H  & (  b(\Psi_j,0), a(0, \Psi_j))_H     \\
 (  b(0,\Psi_j), a(\Psi_j,0))_H   &    (  b(0,\Psi_j), a(0,\Psi_j))_H  
\end{array}
\right) \, . 
\ee
Using the basis \eqref{Fj-eigenfunctions}, the space of linear operators  $ {\cal L}( H_j, H)  $ 
can be identified with  $ H \times H $, 
\be\label{ident-H-4-me}
{\cal L}( H_j, H ) \simeq H \times H  = (L^2 (\T^d, \R))^4 \, , 
\ee
identifying  $ a \in {\cal L}(H_j, H) $  with the vector
\be\label{ident-H-4}
\begin{aligned}
& \qquad \qquad \quad ( a^{(1)}, a^{(2)}, a^{(3)}, a^{(4)} ) \in H \times H   \\
& a(\Psi_j,0) =: (a^{(1)}, a^{(2)}) \in H \, , \quad  
a(0,\Psi_j) =: (a^{(3)}, a^{(4)}) \in H \, , 
\end{aligned}
\ee
so that 
\be \label{identific-matrix}
a \big( q  (\Psi_j , 0) + p (0, \Psi_j  )\big) = q (a^{(1)}, a^{(2)}) + p (a^{(3)}, a^{(4)}) 
\, , \quad \forall (q, p) \in \R^2 \, . 
\ee
With this identification and \eqref{matrix-adjoint} the scalar product \eqref{sc-prod-L} takes the form 
\be\label{scalar-pr-in-H0}
\langle a, b \rangle_0 := {\rm Tr} ( b^* a ) = \sum_{l=1}^4  (a^{(l)}, b^{(l)})_{L^2} 
\ee
and the induced norm 
\be\label{Trastara}
\langle a, a \rangle_0 = {\rm Tr}( a^* a) =  \sum_{l=1}^4 \| a^{(l)} \|^2_{L^2(\T^d)} \, . 
\ee
By the identification \eqref{ident-H-4} and taking in  $ H $  the exponential basis \eqref{basis-exp}, 
a linear operator $ a \in  {\cal  L}( H_j , H ) $ 
can be also identified with the sequence of $ 2 \times 2 $ matrices
$ (a_k^j)_{k \in \Z^d} $,  
$$
 a_k^j  := 
\begin{pmatrix}
 \widehat{a}^{(1)}_k &   \widehat{a}^{(3)}_k   \\
 \widehat{a}^{(2)}_k &  \widehat{a}^{(4)}_k   \\
\end{pmatrix} \in {\rm Mat}_2 (\C)  \, , \quad a^{(l)} = \sum_{k \in \Z^d}  \widehat{a}^{(l)}_k e^{\ii k \cdot x}  \, , \quad
l = 1, 2, 3 , 4 \, , 
$$
so that, by \eqref{identific-matrix}, 
$$
a \big( q  (\Psi_j, 0) + p (0, \Psi_j  )\big) = 
 \sum_{k \in \Z^d} a_k^j (q,p) e^{\ii k x } \, , \quad \forall (q, p) \in \R^2 \, . 
$$
Similarly,
using the basis $ \{ \Psi_j (x) \}_{j \in \N} $ a linear operator $ a \in  {\cal  L}( H_j , H ) $ 
can be also identified with the sequence of $ 2 \times 2 $ matrices
$ ({\mathtt a}_k^j)_{k \in \N} $,  
\be\label{a1-a4}
 {\mathtt a}_k^j  := 
\begin{pmatrix}
{{\mathtt a}}^{(1)}_k &  {{\mathtt a}}^{(3)}_k   \\
{{\mathtt a}}^{(2)}_k &  {{\mathtt a}}^{(4)}_k   \\
\end{pmatrix}  \, , \quad {{\mathtt a}}^{(l)}_k :=  (\Psi_k, a^{(l)})_{L^2} \, , \ \  
l = 1, 2, 3 , 4 \, , 
\ee
so that, by \eqref{identific-matrix}, 
\be\label{per-lemma7.2}
a \big( q  (\Psi_j, 0) + p (0, \Psi_j  )\big) = \sum_{k \in \N} 
{\mathtt a}_k^j (q,p) \Psi_k (x)  \, , \quad \forall (q, p) \in \R^2 \, . 
\ee
In addition, if $ a \in {\cal L}(H_j, H)  $ is identified with the vector 
$(a^{(1)}, a^{(2)} , a^{(3)} , a^{(4)}) \in H \times H$ as in \eqref{ident-H-4}, 
then, $ J a  \in {\cal L}(H_j, H) $ where $ J $ is the symplectic operator 
in \eqref{symplec-op}, can be identified with
\be\label{actionJ}
J a = ( a^{(2)}, - a^{(1)}, a^{(4)}, - a^{(3)} ) \, , 
\ee
and 
$ a J   \in {\cal L}(H_j, H) $  with
\be\label{actionJ-right}
a J = (a^{(1)}, a^{(2)} , a^{(3)} , a^{(4)}) J  = ( - a^{(3)},  - a^{(4)},  a^{(1)}, a^{(2)} ) \, .
\ee
Each space $ {\cal L}( H_j, H )  $ 
admits  the orthogonal decomposition
\be\label{ortho-operator}
{\cal L}( H_j, H ) = {\cal L}( H_j, H_{\mathbb S \cup  \mathbb F} ) \oplus 
{\cal L}( H_j, H_{\mathbb G} ) \, , 
\ee
defined, for any $ a^j \in {\cal L}( H_j, H ) $,  by 
\be\label{ort-deco-op}
a^j = \Pi_{\mathbb S \cup  \mathbb F} a^j + \Pi_{{\mathbb G}} a^j \, , \quad 
\Pi_{\mathbb S \cup  \mathbb F} a^j \in {\cal L}( H_j, H_{\mathbb S \cup  \mathbb F} ), \ \Pi_{{\mathbb G}} a^j \in  
{\cal L}( H_j, H_{{\mathbb G}})
\ee
where $ \Pi_{{\mathbb S \cup  \mathbb F }} $ and $\Pi_{{\mathbb G}} $ are the 
$ L^2 $-orthogonal projectors respectively on the subspaces $ H_{{\mathbb S \cup  \mathbb F}} = 
H_{\mathbb F} \oplus  H_{\mathbb S} $ and  $ H_{{\mathbb G}} $,  see \eqref{orthogonal-deco-FS}. 
 
A (possibly unbounded) linear operator $A$ acting on the Hilbert space
\be\label{phase-space-bfH}
(i) \ {\bf H } := H  = (L^2(\T^d))^2 \, , \qquad (ii) \ {\bf H } := H \times H = (L^2(\T^d))^4 
\ee
can be  represented in the Fourier basis of $L^2(\T^d)$ by 
a matrix $ (A^j_{j'})_{j,j' \in \Z^d} $ with  $A^j_{j'} \in {\rm Mat}_{2 \times 2}(\C)  $ 
in case ($i$), respectively  $ {\rm Mat}_{4\times 4} (\C) $ in case ($ii$),
by the relation
$$
A \Big( \sum_{j \in \Z^d} h_j e^{\ii j \cdot x} \Big)= \sum_{j' \in \Z^d} \Big( \sum_{j \in \Z^d} A^j_{j'} h_j \Big) e^{\ii j' \cdot x} \,
$$
with $ h_j  \in \C^{2}$ in case (i), respectively $ h_j  \in \C^4$ in case (ii).  

\smallskip

We decompose  the space of $ 2 \times 2 $-real matrices 
\be\label{deco:M+M-}
{\rm Mat}_2( \R ) = M_+ \oplus M_- 
\ee
where  $ M_+ $, respectively $ M_- $,  is the subspace of the $ 2 \times 2 $-matrices which commute, 
respectively anti-commute, with the symplectic matrix $ J $. 
A basis of $ M_+ $  is formed by 
\be\label{def:E1E2}
M_1 := J = 
\begin{pmatrix}
 0 &  1   \\
 -1 & 0   \\
\end{pmatrix}   \, ,\quad 
M_2 := {\rm Id}_2 = 
\begin{pmatrix}
 1 &  0   \\
  0 & 1   \\
\end{pmatrix}  \, ,  
\ee
and a basis of $ M_- $ is formed by 
\be\label{def:E3E4}
M_3 := 
\begin{pmatrix}
 1 &  0   \\
 0 & -1   \\
\end{pmatrix}  \, \quad  
M_4 := 
\begin{pmatrix}
 0 &  1   \\
 1 & 0   \\
\end{pmatrix} \, .
\ee
Notice that the matrices $ \{ M_1, M_2, M_3, M_4 \} $ form also a basis for the $ 2 \times 2 $-complex matrices  
$ {\rm Mat}_2( \C ) $.

We shall denote by $ \pi^+ $, $ \pi^- $ the projectors  on $ M^+ $, respectively $ M^- $. We shall use that, for a $ 2 \times 2 $
real symmetric matrix
\be\label{proiettore:sim}
M = \begin{pmatrix}
 a &  b   \\
 b & c   \\
\end{pmatrix} \, , \qquad 
\pi_+ (M)  = \frac{a+c}{2} \, {\rm Id}_2 = \frac{{\rm Tr}(M)}{2} {\rm Id}_2 \, .  
\ee
{\bf $ \vphi $-dependent families of functions and operators.} 
In this monograph  we often identify a function 
$$ 
h \in L^2 (\T^\es, L^2(\T^d)) \, , \quad h : \vphi \mapsto h (\vphi) \in L^2(\T^d) \, ,
$$
with the function $ h(\vphi, \cdot ) (x) = h(\vphi, x) $ of time-space,  i.e.  
$$ 
L^2 (\T^\es, L^2(\T^d)) \equiv L^2(\T^\es \times \T^d) \, . 
$$  
Correspondingly, we regard   a $ \vphi $-dependent family of  (possibly unbounded) 
operators 
$$ 
A : \T^\es \to {\cal L}(H_1, H_2) \, , \quad  \vphi \mapsto A(\vphi) \in {\cal L}(H_1, H_2)  \, , 
$$
acting  between Hilbert spaces $ H_1, H_2 $,  
as an operator $ A $
 which acts on functions $ h(\varphi , x) \in L^2 ( \T^\es, H_1 ) $ of space-time. When 
 $ H_1  = H_2 = L^2 (\T^d) $  we regard $ A $  as the linear operator 
$ A : L^2(\T^\es \times \T^d) \to L^2(\T^\es \times \T^d) $  defined by
\be  \label{AphiA}
(A h) (\varphi , x) := (A(\varphi) h(\varphi, \cdot ))(x) \, .  
\ee
For simplicity of notation we still denote such operator by $ A $. 

Given 
$ a, b : \T^\es \to {\cal L}( H_j, H ) $
we define the scalar product  
\be\label{sc-pr-operators-vphi}
\langle a, b \rangle_0 := \int_{\T^\es} {\rm Tr} \big( b^* (\vphi) a(\vphi)  \big) \, d \vphi \, , 
\ee 
that, with a slight abuse of notation,  we denote with the same symbol \eqref{sc-prod-L}. 
Using  \eqref{ident-H-4} we may identify 
$ a, b : \T^\es \to {\cal L}( H_j, H )  $ with $ a, b : \T^\es \to H \times H $ defined by 
$$ 
a(\vphi) =  (  a^{(1)}, a^{(2)}, a^{(3)}, a^{(4)})(\vphi) \, , \ \ 
b(\vphi) =  (  b^{(1)}, b^{(2)}, b^{(3)}, b^{(4)})(\vphi) \, , 
$$ 
and, 
by \eqref{sc-pr-operators-vphi} and  \eqref{scalar-pr-in-H0},  \eqref{Trastara}, 
\begin{align}\label{equiv-n-pr-vphi}
& \langle  a, b \rangle_0 = \int_{\T^\es} {\rm Tr}(a(\vphi)^* b(\vphi)) \, d\vphi 
= \sum_{l=1}^4 ( a^{(l)}, b^{(l)})_{L^2(\T^\es  \times \T^d) }  \, , \\
& \label{equiv-n-4}
\| a \|_0^2 = \langle a, a \rangle_0 = \int_{\T^\es} {\rm Tr}(a(\vphi)^* a(\vphi)) \, d\vphi 
= \sum_{l=1}^4 \| a^{(l)} \|^2_{L^2(\T^\es  \times \T^d) }  \, .
\end{align}
Consider a $ \vphi $-dependent family of 
linear operators $ A (\vphi) $, acting through \eqref{AphiA} on $ \vphi $-dependent family of functions
with values in the Hilbert space 
$ {\bf H} $  defined in \eqref{phase-space-bfH}. We   identify the family $(A(\vphi))$  
with an  infinite dimensional matrix 
$ (A^{\ell,j}_{\ell',j'})_{(\ell, j), (\ell',j') \in \Z^{\es+d}} $, 
whose entries $ A^{\ell,j}_{\ell',j'}$ are in $ {\rm Mat}_{2 \times 2} (\C) $ 
in case ($i$), and in  $ {\rm Mat}_{4 \times 4} (\C) $  in case ($ii$); this matrix represents the operator $A$ in the
Fourier exponential basis of $   (L^2(\T^\es \times \T^d))^2 $ and is defined by the relation
\be\label{matrix-A-space-time}
A \Big( \sum_{(\ell, j) \in \Z^{\es+d}} h_{\ell, j} e^{\ii ( \ell \cdot \vphi + j \cdot x )} \Big) 
= \sum_{(\ell', j') \in \Z^{\es+d}} \Big( \sum_{(\ell, j) \in \Z^{\es+d}} A^{\ell, j}_{\ell',j'} h_{\ell,j} \Big) e^{\ii (\ell' \cdot \vphi + j' \cdot x)} 
\ee
 where $ h_{\ell, j}  \in \C^{2}$ in case \eqref{phase-space-bfH}-($i$),  
 respectively $ h_{\ell, j}  \in \C^4$ in case \eqref{phase-space-bfH}-($ii$).
 
 Taking in $   (L^2(\T^\es \times \T^d))^2 $ 
 the basis 
 $$ 
 \big\{ e^{\ii \ell \cdot \vphi } (\Psi_j (x), 0), e^{\ii \ell \cdot \vphi } (0, \Psi_j (x)) \big\}_{j \in \N }
 $$
 we identify $ A  $ with another matrix $ (A^{\ell,j}_{\ell',j'})_{(\ell, j), (\ell',j') \in \Z^{\es} \times \N} $ by the relation
 \be\label{matrix-A-space-time-psi}
A \Big( \sum_{(\ell, j) \in \Z^{\es} \times \N} h_{\ell, j} e^{\ii  \ell \cdot \vphi  } \Psi_j (x) \Big) 
= \sum_{(\ell', j') \in \Z^{\es} \times \N} \Big( \sum_{(\ell, j) \in \Z^{\es} \times \N} 
A^{\ell, j}_{\ell',j'} h_{\ell,j} \Big) e^{\ii \ell' \cdot \vphi} \Psi_{j'}(x)  
\ee
 where $ h_{\ell, j}  \in \C^{2}$. 
\\[1mm]
{\bf Normal form.} 
We consider 
a $ \vphi $-dependent family  of linear operators $ A(\vphi) $ acting in $ H_{\mathbb S}^\bot $.
According to the decomposition $ H_{\mathbb S}^\bot = H_{\mathbb F} \oplus H_{\mathbb G} $, 
each $ A(\vphi) $ can be represented  by a matrix as in \eqref{decoFGsotto},
\be\label{decoFGsotto-vphi}
A(\vphi)  = 
\begin{pmatrix}
[A(\vphi)]^{\mathbb F}_{\mathbb F} & [A(\vphi)]_{\mathbb F}^{\mathbb G} \\
[A(\vphi)]^{\mathbb F}_{\mathbb G} &  [A(\vphi)]^{\mathbb G}_{\mathbb G}
\end{pmatrix} \, .
\ee
We denote by 
$ \Pi_{\mathtt D} A (\vphi)$ the operator of  $ H_{\mathbb S}^\bot $ represented by the matrix
\be\label{def:pro+}
\Pi_{\mathtt D} A (\vphi)  =  \begin{pmatrix}
D_+ (A^{\mathbb F}_{\mathbb F}) & 0 \\
0 &   [A(\vphi)]^{\mathbb G}_{\mathbb G}
\end{pmatrix}
\ee
where, in the basis  $ \{ (\Psi_j,0),(0,\Psi_j) \}_{j \in \mathbb F}$ of $H_{\mathbb F}$, 
\be\label{new normal form D+}
D_+ (A^{\mathbb F}_{\mathbb F}) := {\rm Diag}_{j \in \mathbb F} \big( 
\pi_+ [{ \widehat{A}^j_j}(0)] \big) \, ;
\ee
 $ \pi_+ : {\rm Mat}_2(\R) \to M_+ $ denotes the projector on $ M_+ $ associated to the decomposition \eqref{deco:M+M-};
  $ \widehat{A}^j_j (0) $ is the $ \vphi $-average
$$
 \widehat{A}^j_j (0) := \frac{1}{(2\pi)^\es} \int_{\T^\es} A^j_j (\vphi) \, d \vphi \, . 
$$ 
We define also 
\be\label{lem:off-dia}
\Pi_{\mathtt O} A := A - \Pi_{\mathtt D} A \, . 
\ee
We remark that, along the monograph, the functions and 
the operators may depend on a one-dimensional parameter $ \l \in \wtilde \Lambda \subset 
\Lambda $ in a Lipschitz way,  with  norm  defined as in \eqref{def:Lip-norm}.  
\\[1mm]
{\bf Hamiltonian and symplectic operators.}
 Along the paper we shall preserve the Hamiltonian structure of the vector fields.

\begin{definition}\label{def:HAMS}
A $ \vphi $-dependent family of linear operators $ X(\vphi) : {\mathcal D}(X) \subset H \to H $, 
defined on a dense subspace $ {\mathcal D} (X) $ of $ H $ independent of 
$ \vphi \in \T^{\es} $, is \textsc{Hamiltonian}\index{Hamiltonian operator} if 
$$ 
X(\vphi ) = J A( \vphi ) 
$$ 
for some  
real linear operator $ A (\vphi ) $ which is self-adjoint with respect to the $ L^2 $ scalar product.  
We also say that $ \om \cdot \partial_{\vphi} - J A ( \vphi ) $ is Hamiltonian. 
\end{definition}

We mean that  $ A $ is self-adjoint if its domain of definition 
$ {\mathcal D}(A) $ is dense in $ H $, and 
$ ( Ah, k) = (h, Ak ) $, for all $ h , k \in {\mathcal D}(A) $.  

\begin{definition}
A $\vphi $-dependent family of linear operators 
$ \Phi (\vphi ) : H \to H $, $ \forall \vphi \in \T^{\es} $, is {\sc symplectic}\index{Symplectic map} if   
 \be\label{mappa simplettica}
 \Omega( \Phi  (\vphi) u, \Phi  (\vphi)  v) = \Omega (u, v) \, , \quad \forall u,v \in H \, , 
\ee 
where the symplectic 2-form $ \Omega $ is defined in \eqref{2fo}. 
Equivalently 
$$ 
\Phi^*  (\vphi) J \Phi (\vphi)  = J  \, , \quad \forall \vphi \in \T^\es \, .
$$
\end{definition}

A Hamiltonian operator transforms into an Hamiltonian one under a symplectic transformation.

\begin{lemma}\label{lem:sym}
Let  $ \Phi(\vphi)  $, $ \vphi \in \T^{\es} $,  be a family of linear symplectic transformations and  
$ A^* (\vphi)  = A (\vphi)  $, for all $ \vphi \in \T^\es $. 
Then 
$$
 \Phi^{-1} (\vphi) \big( \om \cdot \pa_{\vphi} -  J A(\vphi)  \big) \Phi (\vphi) = 
\om \cdot \pa_\vphi - J A_+ (\vphi) 
$$
where $A_+ (\vphi)  $ is self-adjoint. 
Thus $ \om \cdot \pa_{\vphi} -  J A_+ (\vphi) $ is Hamiltonian.
\end{lemma}

\begin{pf}
We have that
\begin{align}
\Phi^{-1} (\vphi) \big( \om \cdot \pa_{\vphi} -  J A(\vphi)  \big) \Phi (\vphi) 
& = 
\om \cdot \pa_\vphi + \Phi^{-1} (\vphi) (\om \cdot \pa_{\vphi}  \Phi) (\vphi)  - 
\Phi^{-1} (\vphi)   J A(\vphi)  \Phi (\vphi) \nonumber \\
& = \om \cdot \pa_\vphi - J  A_+ (\vphi) \nonumber
\end{align}
with
\begin{align}
A_+ (\vphi) & 
=  J \Phi^{-1} (\vphi) \, (\om \cdot \pa_{\vphi} \Phi) (\vphi)  -  
J \Phi^{-1} (\vphi)   J A(\vphi)  \Phi (\vphi) \nonumber \\
& = J \Phi^{-1} (\vphi) \, (\om \cdot \pa_{\vphi}  \Phi) (\vphi)  + 
\Phi^* (\vphi)  A(\vphi)  \Phi (\vphi) \label{simp2}
\end{align}
using that $ \Phi( \vphi) $ is symplectic. 
Since $ A(\vphi) $ is self-adjoint,  the last operator in \eqref{simp2} is clearly self-adjoint.
In order to prove that also the first operator in \eqref{simp2} 
is self-adjoint we note that, since  $ \Phi( \vphi) $ is symplectic, 
\be\label{sinew}
\Phi^*(\vphi) J (\om \cdot \pa_\vphi  \Phi)(\vphi ) + (\om \cdot \pa_\vphi  \Phi)^*(\vphi ) J 
\Phi (\vphi)  = 0 \, .
\ee
Thus, using that $ \Phi (\vphi) $ is symplectic, 
\begin{align*}
\big( J \Phi^{-1} (\vphi) (\om \cdot \pa_\vphi  \Phi)(\vphi ) \big)^* & = 
-    (\om \cdot \pa_\vphi  \Phi)^*(\vphi ) \, (\Phi^{*})^{-1} (\vphi) J \\
& = -    (\om \cdot \pa_\vphi  \Phi)^*(\vphi ) \, J  \Phi (\vphi) \\
& \stackrel{\eqref{sinew}} = 
\Phi^{*} (\vphi) J  (\om \cdot \pa_\vphi  \Phi) (\vphi )  =  J \Phi^{-1} (\vphi) (\om \cdot \pa_\vphi  \Phi)(\vphi ) \, .
\end{align*}
We have proved that $ A^+ (\vphi) $ is self-adjoint. 
\end{pf}

\section{Decay norms}\label{sec:DN}

Let $ b :=  \es +  d $. For  $ B \subset \Z^b  $ we introduce the subspace
\be\label{H:subspaces}
{\mathcal H}^s_B := \Big\{ u = \sum_{i \in \Z^b}  u_{i} e_i  
\in {\mathcal H}^s  \ : \ u_i \in \C^r \, , \ u_{i} = 0 \ {\rm if} 
\ i \notin B \Big\} 
\ee
where  $ e_i  := e^{\ii (\ell \cdot \vphi + j \cdot x)}  $, $ i = (\ell, j ) \in \Z^\es \times \Z^d $,  and
$ {\mathcal H}^s $ is the Sobolev space \index{Sobolev spaces} 
\be
\label{def:Hs}
{\mathcal H}^s   :=  {\mathcal H}^s ( \T^\es \times \T^d;  \C^r) :=   \Big\{  u = \sum_{i \in \Z^b}  u_{i} e_i  \   :  \, \| u \|_s^2 :=  \sum_{i \in \Z^{b }}  |u_{i}|^2 \langle i \rangle^{2s}  
\Big\} \, . 
\ee
Clearly $ {\mathcal H}^0  =  L^2 ( \T^\es \times \T^d;  \C^r) $, and,
for $ s > b / 2 $, we have the continuous embedding 
$
{\mathcal H}^s  \subset C^0 (  \T^\es \times \T^d;  \C^r) $. 

For a Lipschitz family of functions $ f : \Lambda \mapsto {\mathcal H}^s $, 
$ \lambda \mapsto f(\lambda)  $, 
we define, as in \eqref{def:Lip-norm},
\be\label{Sobo-Lip}
\| f \|_{\Lip,s} := 
\sup_{\l \in \Lambda} \| f \|_s + \sup_{ 
 \l_1, \l_2 \in \Lambda,
 \l_1 \neq \l_2}  \frac{ \| f(\l_2) - f(\l_1) \|_s }{ |\l_2- \l_1|}  \, . 
\ee

\begin{remark}\label{rem:compo}
In Chapter \ref{sec:multiscale} we shall distinguish the components of 
the vector $ u_i = (u_{i, {\mathfrak a}})_{\mathfrak a \in {\mathfrak I}  } \in \C^r $ where
$ {\mathfrak I} = \{1, 2\}$ if $ r = 2 $, and  $ {\mathfrak I} = \{1, 2,3,4 \} $  if $ r = 4 $.   
In this case we also write an element of $ {\mathcal H}^s $ as 
$$
u = 
\sum_{i, \mathfrak a \in \Z^b \times {\mathfrak I}} 
u_{i, \mathfrak a} e_{i, \mathfrak a} \, ,  \ u_{i, \mathfrak a} \in \C \, , 
\quad e_{i, \mathfrak a}  := e_{\mathfrak a} e_i \, ,
$$
where $ e_{\mathfrak a} := (0, \ldots, \underbrace{1}_{{\mathfrak a}-th}, \ldots,0 )  
$, $ \mathfrak a \in {\mathfrak I }$, denotes the canonical basis of $ \C^r $. 
\end{remark}

 When $ B $ is finite, the space $ {\mathcal H}^s_B $ does not 
depend on $ s $ and will be denoted $ {\mathcal H}_B $.   
For $ B, C \subset \Z^b  $ finite, 
we identify the space $ \lin^B_C $ of the linear maps $ L : {\mathcal H}_B \to {\mathcal H}_C $  with the space of matrices
\be\label{sub-matrices-MBC}
\matr^B_C := \Big\{ M = (M^{i'}_i)_{i' \in B, i \in C} \, , \ M^{i'}_{i} \in {\rm Mat}(r \times r; \C) \Big\} 
\ee 
identifying $ L $ with the matrix $ M $ with entries
$$ 
M_{i}^{i'} =  (M_{i, \mathfrak a}^{i', \mathfrak a'} )_{\mathfrak a, \mathfrak a' \in {\mathfrak I}} 
\in {\rm Mat}(r \times r, \C ) \, , \quad
M_{i, {\mathfrak a}}^{i', {\mathfrak a'}}  
:=  (L e_{i', \mathfrak a'}, e_{i, \mathfrak a})_0 \, , 
$$ 
where $ ( \, , \, )_0 := (2\pi)^{-b} ( \, , \, )_{L^2}$ denotes the normalized $ L^2 $-scalar product.

Following \cite{BBo10} we shall use $ s $-decay norms which quantify the polynomial 
decay off the diagonal of the matrix entries.

\begin{definition} \label{defnormatr} {\bf ($s$-norm)}
The $ s $-norm\index{Decay norm} of a matrix $ M \in \matr^B_C $ is defined by
$$
| M |_s^2 :=  \sum_{n \in \Z^b} [M(n)]^2 \langle n \rangle^{2s}  
$$
where $ \langle n \rangle := \max (|n|,1)$ (see \eqref{def:Hs}), 
$$
[M(n)] := \begin{cases}
\sup_{i-i'=n} |M^{i'}_i|  \ \ \quad   {\rm if}  \ \   n \in  C- B\\
0  \qquad \qquad \qquad  {\rm if}  \ \  n \notin C - B \, , 
\end{cases}
$$
where $  | \ | $ denotes a norm of the matrices $ {\rm Mat}(r \times r, \C )$. 
\end{definition}

We shall use the above definition also if $B$ or $C$ are not finite (with the difference that $| M|_s$ may be infinite).

The $ s $-norm is modeled on matrices  which represent the multiplication operator.
The (T\"oplitz) matrix $ T $ which represents the multiplication 
operator 
$$ 
M_g : {\mathcal H}^s (\T^\es \times \T^d;  \C) \to {\mathcal H}^s (\T^\es \times \T^d;  \C) \, , 
\quad h \mapsto g h \, , 
$$ 
by a function  $ g  \in {\mathcal H}^s (\T^\es \times \T^d;  \C) $, 
$ s \geq s_0 $,  satisfies 
\be\label{multi-op-mat} 
|T|_s \sim \| g \|_s \, , \quad 
|T|_{\Lip,s} \sim \| g \|_{\Lip,s} \, . 
\ee 
The $ s$-norm  satisfies algebra and interpolation inequalities 
and  control the higher Sobolev norms as in \eqref{opernorm} below:   
as proved in \cite{BBo10}, for all $ s \geq s_0 > b / 2 $
\be\label{inter-AB-tutto}
| A B |_s \lesssim_s  | A |_{s_0} | B |_s +  | A |_s | B |_{s_0} \, , 
\ee
and 
for any  subset $ B, C \subset \Z^b $,  $ \forall M \in \matr^B_C $, $  w \in {\cal H}_B $ we have 
\begin{align}\label{opernorm}
\| Mw \|_s & \lesssim_s | M |_{s_0} \|w\|_s +  | M |_s \|w\|_{s_0} \, ,\\
\label{opernorm-Lip}
\| Mw \|_{\Lip,s} &  \lesssim_s | M |_{\Lip, s_0} \|w\|_{\Lip,s} +  | M |_{\Lip,s}
 \|w\|_{\Lip, s_0} \, .
\end{align}
The above inequalities can be easily obtained from the definition of the norms $| \ |_s$ and
the functional interpolation inequality\index{Interpolation inequalities}
\be  \label{intbasic}
\forall u,v \in {\cal H}^s \, , \ 
\|uv\|_s \lesssim_s \|u\|_{s_0}  \|v\|_{s} + \|u\|_{s} \|v\|_{s_0} \, .
\ee
Actually, \eqref{intbasic} can be slightly improved to obtain (see Lemma \ref{lem:int-impro})
\be  \label{intbasic+}
\|uv\|_s \leq C_0  \|u\|_{s_0}  \|v\|_{s} +  C(s) \|u\|_{s} \|v\|_{s_0} 
\ee
where only the second constant may depend on $s$ and $C_0 $ depends only on $ s_0 $
(we recall that $s_0$ is fixed once for all). 
From \eqref{intbasic+}
can be derived the following slight improvements of \eqref{inter-AB-tutto} and 
\eqref{opernorm-Lip}, which will be used in Chapter \ref{sec:proof.Almost-inv}:
\be\label{inter-AB-tutto+}
| A B |_s \leq C_0 | A |_{s_0} | B |_s +  C(s)  | A |_s | B |_{s_0} 
\ee
and 
\be  \label{opernorm-Lip+}
\| Mw \|_{\Lip,s}   \leq C_0  | M |_{\Lip, s_0} \|w\|_{\Lip,s} +  C(s) | M |_{\Lip,s}
 \|w\|_{\Lip, s_0} \, .
\ee
We also note that, denoting by $ A^* $ the adjoint matrix of $ A $, we have 
\be\label{A-s-adj}
|A|_s = |A^*|_s \, , \quad |A|_{\Lip,s} = |A^*|_{\Lip,s}  \, .   
\ee
The following lemma is the analogue of the smoothing properties  of the projection operators.
See Lemma 3.6 of \cite{BBo10}, and Lemma \ref{norcompA}. 

\begin{lemma} \label{norcomp}
{\bf (Smoothing)} 
Let $ M  \in \matr^B_C $. Then, $ \forall s' \geq  s \geq 0 $,
\be \label{Sm1}
M_i^{i'}=0 \, , \ \forall |i - i' | < N \quad \Longrightarrow  \quad  |M|_s \leq N^{-(s'-s)} |M|_{s'} \, , 
\ee
and similarly for the Lipschitz norm $ | \ |_{\Lip, s} $. 
\end{lemma}

We now define the decay norm for an operator $ A : E \to F $ defined on a 
closed subspace $ E \subset L^2 $ with range in a closed subset of $ L^2 $. 

\begin{definition}\label{def:decay-sub}
Let  $ E, F $ be closed subspaces of $ L^2 \equiv L^2( \T^\es \times \T^d, \C^r)$.
Given a linear operator $A : E \to F $, we extend it to a linear operator 
$ \wtilde{A} :  L^2  \to L^2 $  acting on the whole  $ L^2 = E \oplus E^\bot $, with 
 image  in $ F  $,  by defining\index{Decay norm}  
\be\label{op:extended}
\wtilde{A}_{|E^\bot} := 0 \, . 
\ee
Then, for $ s \geq  0 $ we define the (possibly infinite)
$ s$-decay norm 
\be\label{def:s-decay}
|A|_s := |\wtilde {A}|_s  \, , 
\ee
and 
\be\label{def:n+s}
|A|_{+,s} := | D_m^{\frac12} {\wtilde A} D_m^{\frac12}|_s 
\ee
where $ D_m $ is  the Fourier multiplier operator 
\be\label{def:Dm}
D_m := \sqrt{- \Delta + m} \, , \quad D_m (e^{\ii j \cdot x}) := \sqrt{|j|^2 + m} \, e^{\ii j \cdot x} \, , \quad j \in \Z^d \, , 
\ee
and $ m > 0 $ is a positive constant.

For a Lipschitz family of operators $ A(\l) : E \to E $, $ \lambda \in \Lambda $, 
we associate the Lipschitz norms $ |A|_{\Lip, s}  $, $ |A|_{\Lip, +,s}  $, accordingly to \eqref{def:Lip-norm}. 
\end{definition}

The norm \eqref{def:n+s} is stronger than \eqref{def:s-decay}, actually, since
$ 1 \lesssim_m (|j|^2 + m)^{\frac14}  (|j'|^2 + m)^{\frac14} $, we have  
\be\label{simple-embed}
  | A |_{s} \lesssim_m | A |_{+,s} \, , \quad 
| A |_{\Lip, s} \lesssim_m | A |_{\Lip, +,s}     \, . 
\ee

\begin{lemma}\label{normprod} 
{\bf (Tame estimates\index{Tame estimates} for composition)}
Let $ A, B : E \to E $ be linear operators acting on a closed subspace $ E \subset L^2 $. 
Then, the following tame estimates hold: for all $s \geq s_0 > (\es +d)/2$, 
\begin{align}\label{tame-s-decay}
& |AB|_s \lesssim_s |A|_{s_0} |B|_s + |A|_s |B|_{s_0} \\
& |AB|_{\Lip, s} \lesssim_s |A|_{\Lip, s_0} |B|_{\Lip, s} + |A|_{\Lip, s} |B|_{\Lip, s_0} \, , \label{tame-s-decayLip}
\end{align}
more precisely
\be \label{tame-s-decayLip+}
 |AB|_{\Lip, s} \leq C_0 |A|_{\Lip, s_0} |B|_{\Lip, s} +  C(s) |A|_{\Lip, s} |B|_{\Lip, s_0}
\ee
and 
\begin{align}
& |AB|_{+,s} + |BA|_{+,s} \lesssim_s  |A|_{s_0+ \frac12} |B|_{+,s} + |A|_{s+ \frac12} |B|_{+,s_0}   \label{normprod1} \\
& |AB|_{\Lip, +,s} + |BA|_{\Lip, +,s}  \lesssim_s
 |A|_{\Lip, s_0+ \frac12} |B|_{\Lip, +,s} + |A|_{\Lip, s+ \frac12}  
|B|_{\Lip, +,s_0}  \, . \label{inter-norma+s}
\end{align}
Notice that in \eqref{normprod1}, resp. \eqref{inter-norma+s}, 
the operator $ A  $ is  estimated in $ | \ |_{s+ \frac12}$ norm, resp. $ | \ |_{\Lip, s+ \frac12}$, and not in 
$ | \ |_{+, s} $, resp. $ | \ |_{\Lip, +, s} $. 
\end{lemma}

\begin{pf}
Notice first 
that the operation introduced in \eqref{op:extended} of extension of an operator commutes with the composition:
 if $A, B : E \to E $ are linear operators acting in $ E $, then 
$\wtilde{AB}= \wtilde{A}  \, \wtilde{B}$.
Thus \eqref{tame-s-decay}-\eqref{tame-s-decayLip} and \eqref{tame-s-decayLip+} follow by the interpolation inequalities 
\eqref{inter-AB-tutto} and \eqref{inter-AB-tutto+}. 

In order to prove \eqref{normprod1}-\eqref{inter-norma+s} we first show that, given a linear 
operator acting on  the whole $ L^2 $, 
\be\label{lemma-primo-comm}
| D_m^{-\frac12} \, A \, D_m^{\frac12}  |_s  \, , \   | D_m^{\frac12} \, A \, D_m^{- \frac12}  |_s \lesssim_s | A |_{s + \frac12 } \, .
\ee
We prove \eqref{lemma-primo-comm} for 
$$
D_m^{\frac12} \, A \, D_m^{- \frac12} = A + [D_m^{\frac12}, A] D_m^{-\frac12}  \, . 
$$
Since $  |D_m^{-\frac12}|_{s} \leq C(m) $, $ \forall s $, it is sufficient to prove that
$ | [D_m^{\frac12}, A] |_s \lesssim_s | A |_{s + \frac12 } $. Since, $ \forall i , j \in \Z^d $, 
\begin{align*}
& \big| (|j|^2+m)^{1/4} - (|i|^2+m)^{1/4}\big| = \nonumber \\ 
& \frac{\big| |j| - |i| \big| (|j|+|i|)}{\big((|j|^2+m)^{1/4} 
 + (|i|^2+m)^{1/4}\big)((|j|^2+m)^{\frac12} + (|i|^2+m)^{\frac12})} \\
& \leq \frac{| j - i | }{\big((|j|^2+m)^{1/4} 
+ (|i|^2+m)^{1/4}\big)} \leq  \sqrt{| j - i |} \, , 
\end{align*}
the matrix elements of the commutator $ [D_m^{\frac12}, A]^i_j $ satisfy 
$$
|[D_m^{\frac12}, A]^i_j| = 
|A^i_j | \big| (|j|^2+m)^{1/4} - (|i|^2+m)^{1/4}\big| \leq  |A^i_j | \sqrt{| j - i |} 
$$
and  therefore \eqref{lemma-primo-comm} follows. 

We now prove the estimates \eqref{normprod1}-\eqref{inter-norma+s}. 
We consider the extended operator   $ \wtilde {A B} =  \wtilde {A} \wtilde {B} $ 
and  write
$$
D_m^{\frac12} \wtilde {A B} D_m^{\frac12}   = \big(D_m^{\frac12} \wtilde{A} D_m^{-\frac12}\big)   \big(D_m^{\frac12} 
\wtilde{B} D_m^{\frac12} \big) \, . 
$$
Hence \eqref{normprod1}-\eqref{inter-norma+s}
follow, recalling  \eqref{def:n+s}, by  \eqref{inter-AB-tutto}
and  \eqref{lemma-primo-comm}.
\end{pf}

By iterating \eqref{tame-s-decayLip} we deduce that, there exists 
$ C(s) \geq 1 $, non decreasing in $ s \geq s_0 $,  such that 
\be\label{AkLip}
|A^k|_{\Lip,s} \leq (C(s))^{k} |A|_{\Lip,s_0}^{k-1} |A|_{\Lip,s} \, , \quad \forall k \geq 1 \, . 
\ee
\begin{lemma}
Let $ A: E \to E $ be a linear operator acting on a closed subspace $ E \subset L^2 $. 
Then its operatorial norm satisfies  
\be\label{incl-plus}
\| A \|_{0}  \lesssim_{s_0} | A |_{s_0} 
\lesssim_{s_0} | A |_{+, s_0} \, .
\ee
\end{lemma}

\begin{pf}
Let $ \wtilde A $ be  the extended operator in $ L^2 $ defined in \eqref{op:extended}.
Then, using Lemma 3.8 in \cite{BBo10} (see Lemma \ref{Lem:L2s0}), and
\eqref{def:s-decay}, we get  
$ \| A \|_0 \leq $ $ \| \wtilde A \|_0 \lesssim_{s_0}  $  $ | \wtilde A |_{s_0} 
= $
$ | A |_{s_0}
\lesssim_{s_0}  | A |_{+s_0}  $
by \eqref{simple-embed}. 
\end{pf}

We shall also use the following elementary inequality: given a matrix $ A \in {\cal M}^B_C $ where $ B $ and $ C $ 
are included in $ [-N, N]^b $, then 
\be\label{el:AD}
| A  |_s \lesssim N | D^{-\frac12}_m A D^{-\frac12}_m |_s \, . 
\ee
We shall use several times the following simple lemma.

\begin{lemma}  \label{gchi-r}
Let $ g (\l, \vphi, x), \chi (\l, \vphi, x ) $ be a Lipschitz family of functions 
in $ {\mathcal H}^s (\T^\es \times \T^d, \C)$ for $s\geq s_0$. 
Then the operator $L$ defined by
$$
L[h] (\vphi,x) := \big( h(\vphi, \cdot ) , g(\vphi, \cdot) \big)_{L^2_x} \chi (\vphi,x) \, , 
\quad \forall h \in {\mathcal H}^0 (\T^\es \times \T^d, \C) \, , 
$$
satisfies 
\be \label{boundsimpleL}
|L|_{\Lip, s}  \lesssim_s \| g\|_{\Lip, s_0} \| \chi\|_{\Lip,s} + 
\| g \|_{\Lip, s} \| \chi\|_{\Lip, s_0} \, . 
\ee
\end{lemma}

\begin{pf}
We write $ L = M_{\chi} P_0 M_{\ov{g}}$ as the composition 
of the multiplication operators $ M_{\chi}, M_{\ov{g}}$ for the functions 
$ \chi,   \ov{g} $ respectively, 
and the  mean value projector $ P_0 $ defined as
$$
P_0 h (\vphi) := \frac{1}{(2\pi)^d} \int_{\T^d} h(\vphi , x) \, dx \, , \quad \forall h \in {\mathcal H}^0 \,  . 
$$
For $i=(j,\ell)$, $i'=(j',\ell') \in \Z^d \times \Z^\es$,  its entries are 
$(P_0)^{i'}_i=\delta^0_j \delta^{j'}_0 \delta^{\ell '}_\ell $, and therefore 
\be\label{projector0} 
|P_0|_{\Lip,s} = |P_0|_s \leq 1 \, ,\quad  \forall s \, . 
\ee 
We derive \eqref{boundsimpleL} by 
\eqref{multi-op-mat}, \eqref{projector0}  and 
the tame estimates \eqref{inter-AB-tutto} for the composition of operators. 
\end{pf}

%

Now, given a finite set  $ {\mathbb M} \subset \N $, we estimate the  $ s $-decay norm of the $ L^2 $-orthogonal 
projector  
$ \Pi_{\mathbb M} $ on  
the
subspace of $ L^2 (\T^d, \R) \times L^2 (\T^d, \R) $ defined by
\be\label{def:HE}
H_{\mathbb M} := \Big\{   (q(x), p(x)) := \sum_{j \in {\mathbb M}} (q_j, p_j) \Psi_j (x)  \, , \quad q_j , p_j \in \R \Big\} \, .
\ee
In the next lemma we regard $ \Pi_{\mathbb M} $ as an operator acting on functions $ h(\vphi, x )$.  

\begin{lemma}\label{pisig} {\bf (Off-diagonal decay of $ \Pi_{\mathbb M} $)} 
Let $ {\mathbb M} $ be a finite subset of $ \N $. Then, for all $ s \geq 0 $, 
there is  a constant $ C(s) := C(s, {\mathbb M})  > 0  $ such that   
\be\label{prpr1}
| \Pi_{\mathbb M}|_{+,s} = |\Pi_{\mathbb M}|_{\Lip, +, s}  \leq C(s) \, .
\ee
In addition, setting $ \Pi_{\mathbb M}^\bot := {\rm Id} - \Pi_{\mathbb M} = 
\Pi_{{\mathbb M}^c}$, we have  
\be\label{sepr2}
|\Pi_{\mathbb M}|_{\Lip,s}  \leq C(s)  \, , \qquad
|\Pi_{{\mathbb M}^c}|_{\Lip,s}  
= |\Pi_{\mathbb M}^\bot|_{\Lip,s}  \leq C(s) \, . 
\ee
\end{lemma}

\begin{pf}
To prove \eqref{prpr1} we have to estimate $ | \Pi^\bot_{\mathbb M} |_{\Lip,+,s} = 
 | D_m^{1/2} \Pi^\bot_{\mathbb M} D_m^{1/2} |_{\Lip,s}  $. 
For any $ h = (h^{(1)} , h^{(2)}) \in \big( L^2(\T^\es \times \T^d) \big)^2 $ we have 
\begin{align} 
(D_m^{\frac12} \Pi^\bot_{\mathbb M} D_m^{\frac12} h) (\vphi, x) =
\sum_{j \in {\mathbb M}} 
\begin{pmatrix}
(h^{(1)} (\vphi, \cdot), D_m^{\frac12} \Psi_j)_{L^2_x} \\
(h^{(2)} (\vphi, \cdot), D_m^{\frac12} \Psi_j)_{L^2_x} 
\end{pmatrix}  D_m^{\frac12} \Psi_j  \, . \label{con-plus-M}
\end{align}
Now, using Lemma \ref{gchi-r}, the fact that  each $\Psi_j (x) $ is in $ C^\infty $
and $\mathbb M $ is finite, we deduce, by \eqref{con-plus-M}, the estimate \eqref{prpr1}.
The first estimate in \eqref{sepr2}  is trivial  because 
$ | \Pi_{\mathbb M}|_{\Lip, s} \lesssim |\Pi_{\mathbb M} |_{\Lip, +,s} $.
The second estimate in \eqref{sepr2}  follows by $ \Pi_{\mathbb M} + \Pi_{\mathbb M}^\bot = {\rm Id} $
and $ |{\rm Id}|_{\Lip,s} = 1  $. 
\end{pf}

\begin{lemma} \label{est:PiDO}
Given a $\vphi$-dependent family of linear operators $A(\vphi)$ acting in 
$ H_{\mathbb S}^\bot $, the operators 
$ \Pi_{\mathtt D} A $ and $ \Pi_{\mathtt O} A $  defined respectively 
in \eqref{def:pro+} and \eqref{lem:off-dia}  satisfy 
\begin{align}
& |\Pi_{\mathtt D} A|_{\Lip, +,s}+ |\Pi_{\mathtt O} A |_{\Lip, +,s} \leq C(s) | A |_{\Lip, +,s} \, .
\label{estimate-deco-resto}
\end{align}
\end{lemma}

\begin{pf}
We consider the extension of the operator  
$\Pi_{\mathtt D} A $ defined in  \eqref{def:pro+}-\eqref{new normal form D+}.
This extension,  
acting on $ \big(L^2(\T^\es \times \T^d) \big)^2$, is defined as
$  A_1 +A_2  $
where the operators $A_1$ and $A_2$ are 
\be \label{defA1}
A_1 := \Pi_{\mathbb G} A   \Pi_{\mathbb   G} \, , 
\ee
and, for any $h=(h^{(1)} , h^{(2)}) \in \big( L^2(\T^\es \times \T^d) \big)^2$,
\be \label{defA2}
(A_2 h) (\vphi, x) :=\sum_{j \in {\mathbb F}} \pi_+ \widehat{A}_j^j (0) 
\begin{pmatrix}
(h^{(1)} (\vphi, \cdot), \Psi_j)_{L^2_x} \\
(h^{(2)} (\vphi, \cdot), \Psi_j)_{L^2_x} 
\end{pmatrix}  \Psi_j  (x) \, . 
\ee
We remark that, for $j \in {\mathbb F}$, $\widehat{A}_j^j (0) $ is the $ 2 \times 2 $  matrix
\be\label{Ajj0}
\widehat{A}_j^j (0)  =
\begin{pmatrix}
( A (\Psi_j,0), (\Psi_j,0) )_0 & ( A (0,\Psi_j), (\Psi_j,0) )_0 \\
( A (\Psi_j,0), (0, \Psi_j) )_0 & ( A (0, \Psi_j), (0, \Psi_j) )_0
\end{pmatrix}  
\ee
where $(\ , \ )_0 $ is the (normalized) $L^2$ scalar product in 
$ \big(L^2(\T^\es \times \T^d) \big)^2$.

By  Definition \ref{def:decay-sub}, 
we have  $| \Pi_{\mathtt D} A|_{\Lip, +,s}=|A_1 +A_2|_{\Lip, +,s}$.
By \eqref{defA1}, \eqref{inter-norma+s} and \eqref{sepr2} (with 
 ${\mathbb M}^c = {\mathbb G} $) 
 we get 
\be\label{esA1d}
|A_1|_{\Lip, +,s} \lesssim_s   | \Pi_{\mathbb G} |^2_{\Lip, s+\frac12} |A|_{\Lip, +,s}  \lesssim_s  |A|_{\Lip, +,s}  \, . 
\ee
For $A_2$, we apply Lemma \ref{gchi-r}. Using that, by \eqref{defA2},
$$
(D_m^{\frac12}A_2 D_m^{\frac12} h) (\vphi, x)=\sum_{j \in {\mathbb F}} \pi_+ 
\widehat{A}_j^j (0) 
\begin{pmatrix}
(h^{(1)} (\vphi, \cdot), D_m^{\frac12} \Psi_j)_{L^2_x} \\
(h^{(2)} (\vphi, \cdot), D_m^{\frac12} \Psi_j)_{L^2_x} 
\end{pmatrix}  D_m^{\frac12} \Psi_j  
$$
with $\widehat{A}_j^j (0)  $ given in \eqref{Ajj0}, 
that $\Psi_j \in C^\infty(\T^d) $ and $\mathbb F$ is finite, we deduce by 
\eqref{boundsimpleL} that 
\be \label{est1A2}
|A_2|_{\Lip, +,s} = 
|D_m^{\frac12}A_2 D_m^{\frac12}|_{\Lip,s} 
\lesssim_s \max_{j \in {\mathbb F}} \| A  \|_{\Lip,0} \, . 
\ee
Finally \eqref{esA1d}, \eqref{est1A2} imply 
$$
|\Pi_{\mathtt D} A|_{\Lip, +,s} \leq |A_1|_{\Lip, +,s} + |A_2|_{\Lip, +,s} \leq C(s) |A|_{\Lip, +,s} 
$$
and
$$
|\Pi_{\mathtt O} A|_{\Lip, +,s}=|A- \Pi_{\mathtt D} A|_{\Lip, +,s} \leq |A|_{\Lip, +,s} +| \Pi_{\mathtt D} A|_{\Lip, +,s}  \leq C(s) |A|_{\Lip, +,s}  \, .
$$
Thus \eqref{estimate-deco-resto} is proved. 
\end{pf}

We now mention some norm equivalences and estimates that shall be used in the sequel. 

Let us consider a $ \vphi$-dependent family of operators 
$ \rho(\vphi) \in {\cal L} (H_{\mathbb S}^\bot ) $, that, 
according to the  splitting \eqref{decoFGsotto}, have the form
$$
\rho(\vphi)=\begin{pmatrix} \rho_1 (\vphi)& \rho_2 (\vphi)^* \\
\rho_2 (\vphi) & 0    \end{pmatrix} \in {\cal L} (H_{\mathbb S}^\bot ) \, , \quad  
\rho_1 (\vphi) \in {\cal L} (H_{\mathbb F}) \, ,  \   \rho_2 (\vphi) \in {\cal L} (H_{\mathbb F}, H_{\mathbb G}) \, . 
$$
Recalling  \eqref{a-deco-sez}, we can identify each  $ \rho_l(\vphi) $, $l=1,2 $,  with 
$$
\rho_l(\vphi)=(\rho_{l,j}(\vphi))_{j \in {\mathbb F}} \, ,  \   {\rm where} \ \  
\rho_{l,j}(\vphi) := (\rho_{l})_{| H_j} \in   
\begin{cases}  {\cal L}(H_j, H_{\mathbb F}) \  {\rm if}\  l=1 \\  
 {\cal L}(H_j, H_{\mathbb G})  \  {\rm if}\  l=2 \, . 
\end{cases}
$$
Moreover, according to \eqref{ident-H-4}, we can identify 
each operator $ \rho_{l,j}(\vphi) $ with the vector 
\be\label{l1l2FG}
 (\rho^{(1)}_{l,j}(\vphi), \rho^{(2)}_{l,j}(\vphi), \rho^{(3)}_{l,j}(\vphi) , \rho^{(4)}_{l,j}(\vphi)) \in 
\begin{cases}   H_{\mathbb F}  \times   H_{\mathbb F} \  \ \, {\rm if}  \  \  l=1 \\  
H_{\mathbb G}\times H_{\mathbb G} \  \ {\rm if}\  \ l=2 \, , 
\end{cases}
\ee
where
\be\label{ide-327}
\begin{pmatrix}  \rho^{(1)}_{l,j} \\ \rho^{(2)}_{l,j}  \end{pmatrix} := \rho_l    \begin{pmatrix}  \Psi_j  \\ 0 \end{pmatrix}   \     ,    \     
\begin{pmatrix}  \rho^{(3)}_{l,j} \\ \rho^{(4)}_{l,j}  \end{pmatrix} := \rho_l  \begin{pmatrix}  0 \\ \Psi_j  \end{pmatrix}  \, .
\ee 
We define the Sobolev norms of each $\rho_l (\vphi) $, $ l = 1, 2 $,  as
\be \label{opSob}
\|\rho_l \|^2_s = \max_{j \in {\mathbb F}} 
\Big( \sum_{k=1}^4 \|\rho_{l,j}^{(k)} \|_s^2  \Big) \, ,
\ee
and $ \|\rho_l \|_{\Lip, s} $ according to   \eqref{def:Lip-norm}. 

\begin{lemma}\label{Lemma:rho}
We have  
\begin{align}\label{special-form}
& | \rho_1 |_{\Lip, +,s}  \sim_s | \rho_1 |_{\Lip, s} \sim_s \| \rho_1 \|_{\Lip, s} \, , \\ 
& | \rho_2 |_{\Lip, +,s}  \lesssim_s | \rho_2 |_{\Lip, s + \frac12}  \, , \    
 | \rho_2 |_{\Lip, s} \sim_s \| \rho_2 \|_{\Lip, s}  \, . 
\label{special-form-Lip}
\end{align}
\end{lemma}
\begin{pf}
\\[1mm]
{\bf Step1.} We first justify that, for $l=1,2$,  
\be \label{rhoeq}
|\rho_l|_{\Lip,s} \sim_s \|\rho_l\|_{\Lip,s} \, . 
\ee  
By Definition \ref{def:decay-sub}, the $ s$-norm $|\rho_l|_{\Lip,s} = 
|\rho_l \Pi_{\mathbb F}|_{\Lip,s}$.
Now, recalling \eqref{ide-327},  
the extended operator $ \rho_l \Pi_{\mathbb F} $ has the form,   
for any $h(\vphi)=(h^{(1)} (\vphi) , h^{(2)} (\vphi)) \in H $, 
\be \label{rhoPi}
(\rho_l \Pi_{\mathbb F} h)(\vphi)= \sum_{j \in {\mathbb F}} ( h^{(1)} (\vphi)  , \Psi_j)_{L^2_x} 
\begin{pmatrix}  \rho^{(1)}_{l,j}(\vphi)  \\ \rho^{(2)}_{l,j}(\vphi)  \end{pmatrix}
+( h^{(2)} (\vphi)  , \Psi_j)_{L^2_x} 
\begin{pmatrix}  \rho^{(3)}_{l,j}(\vphi)  \\ \rho^{(4)}_{l,j}(\vphi)  \end{pmatrix}  \, . 
\ee
Using that $\Psi_j \in C^\infty (\T^d) $,  for all $j\in {\mathbb F}$,  
and that $\mathbb F$ is finite, we derive, from
\eqref{rhoPi} and Lemma \ref{gchi-r}, the estimate $ |\rho |_{\Lip,s} \lesssim_s \|\rho_l\|_{\Lip,s}$.

The reverse inequality $\|\rho_l\|_{\Lip,s} \lesssim_s |\rho_l|_{\Lip,s}$ follows by 
\eqref{ide-327}, \eqref{opernorm-Lip},  and the fact that 
$  \Psi_j \in C^\infty (\T^d) $.  
This concludes the proof of \eqref{rhoeq} and so of the second 
equivalences in  \eqref{special-form}-\eqref{special-form-Lip}.
\\[1mm]
{\bf Step 2.}  We now prove that 
\be  \label{rho1+est}
|\rho_1|_{\Lip , +,s} \lesssim_s \|\rho_1\|_{\Lip,s} \, .
\ee
By Definition   \ref{def:decay-sub}, the norm $|\rho_1|_{\Lip,+,s} = |D_m^{\frac12 } \rho_1 \Pi_{\mathbb F} D_m^{\frac12 }|_{\Lip,s}$ and, by 
\eqref{rhoPi}, 
\begin{align} \label{rhoPi+}
(D_m^{\frac12 } \rho_1 \Pi_{\mathbb F} D_m^{\frac12 } h)(\vphi)&=
\sum_{j \in {\mathbb F}} ( h^{(1)} (\vphi), D_m^{\frac12 } \Psi_j)_{L^2_x} 
\begin{pmatrix}  D_m^{\frac12 }  \rho^{(1)}_{1,j}(\vphi)  \\ 
D_m^{\frac12 }  \rho^{(2)}_{1,j}(\vphi)  \end{pmatrix} \nonumber \\
&+\sum_{j \in {\mathbb F}} ( h^{(2)} (\vphi), D_m^{\frac12 } \Psi_j)_{L^2_x}
\begin{pmatrix}  D_m^{\frac12 }  \rho^{(3)}_{1,j}(\vphi)  \\ 
D_m^{\frac12 }  \rho^{(4)}_{1,j}(\vphi)  \end{pmatrix}   \, . 
\end{align}
Now, applying Lemma \ref{gchi-r},  we deduce by \eqref{rhoPi+}
and the fact that $\Psi_j$  are $C^\infty$ (and independent on $\l$)
and $ \mathbb F$ is finite, that 
\be\label{interm-ste2}
|\rho_1|_{\Lip,+,s} = |D_m^{\frac12 } \rho_1 \Pi_{\mathbb F} D_m^{\frac12 }|_{\Lip,s} 
\lesssim_s \max_{j \in {\mathbb F}, k=1, \ldots, 4} \| D_m^{\frac12 } \rho^{(k)}_{1,j} \|_{\Lip,s }\, .
\ee
Now, by \eqref{l1l2FG}, 
each $ \rho^{(k)}_{1,j} $, $ k=1, \ldots, 4 $, is a function of the form 
$u(\vphi)=\sum_{j \in {\mathbb F} }u_j(\vphi) \Psi_j $. 
We claim that, for functions of this form,  
$ \| D_m^{\frac12} u \|_{\Lip,s} \sim_s \| u \|_{\Lip, s} $.  
Indeed 
$$
\| u \|_{\Lip ,s} \sim_s  \sum_{j \in {\mathbb F}} \| u_j \|_{\Lip, H^s(\T^\es)} \,  , 
$$
because $\mathbb F$ is finite and $\Psi_j$  is $C^\infty$ (and independent on $\l$). Similarly,
$$
\| D_m^{\frac12} u \|_{\Lip ,s} = \Big\| \sum_{j \in {\mathbb F} } u_j(\vphi) D_m^{\frac12} \Psi_j \Big\|_{\Lip , s} \sim_s  \sum_{j \in {\mathbb F}} \| u_j \|_{\Lip, H^s(\T^\es)}  \, ,
$$
and the claim follows. 
Applying this property to $ u =\rho^{(k)}_{1,j} $ we deduce, 
by \eqref{interm-ste2}, that 
$$
|\rho_1|_{\Lip,+,s}
\lesssim_s 
\max_{j \in {\mathbb F}, k=1, \ldots, 4} \| \rho^{(k)}_{1,j} \|_{\Lip,s }   
 $$
and thus  \eqref{rho1+est}. Finally, by \eqref{rhoeq} 
and \eqref{simple-embed} 
we have $ \| \rho_1 \|_{\Lip,s} \lesssim_{s} | \rho_1 |_{\Lip,+, s}  $
and thus we deduce the first equivalence in \eqref{special-form}. 
\\[1mm]
{\bf Step 3.} We finally check 
\be \label{rho2est}
|\rho_2|_{\Lip,+,s}  \lesssim_s |\rho_2|_{\Lip,s+ \frac12} \, , 
\ee
which is the first estimate  in \eqref{special-form}. 
We have 
$$
|\rho_2|_{\Lip,+,s}  = | \rho_2 \Pi_{\mathbb F}|_{\Lip,+,s} 
\stackrel{\eqref{inter-norma+s}} {\lesssim_s} 
|\rho_2|_{\Lip,s+\frac12}  | \Pi_{\mathbb F}|_{\Lip,+,s} 
\stackrel{ 
\eqref{prpr1}} {\lesssim_s} 
|\rho_2|_{\Lip,s+\frac12} 
$$
(applied with ${\mathbb M}={\mathbb F}$). 
\end{pf}

We also remind a standard perturbation lemma for operators which admit 
a right inverse. 

\begin{definition}\label{def:right-inv} {\bf (Right Inverse)}\index{Right inverse}
A matrix $ M \in \matr_C^B $ has a right inverse, that we denote by  $ M^{-1}
\in  \matr^C_B $,  if  $ M M^{-1} = {	\rm Id}_C $.  
\end{definition}

Note that $ M $ has a right inverse  if and only if $ M $ (considered as a linear map) is surjective.
The following lemma is proved as in  Lemma 3.9 of \cite{BBo10} (see also Lemma 
\ref{leftinvA}) by a Neumann series argument. 

\begin{lemma} \label{leftinv} 
There is a constant $c_0 >  0 $ such that, for any $C,B \subset \Z^b$,
for any $ M \in \matr_C^B $ having  a right  inverse $ M^{-1} \in \matr^C_B $,
 for any $ P $ in $ \matr_C^B $ with 
$ | M^{-1} |_{s_0} | P |_{s_0} \leq c_0 $ 
the matrix $ M + P $ has a right  inverse 
that  satisfies
\be \label{inv1}
\begin{aligned}
& |  (M + P)^{-1} |_{s_0} \leq 2 | M^{-1}  |_{s_0} \, \\
& | (M + P)^{-1} |_{s}  \
\lesssim_s  | M^{-1} |_s + | M^{-1} |^2_{s_0} | P |_s  \, , \ \forall s \geq s_0 \, . 
\end{aligned}
\ee
\end{lemma}

Finally, we report the following lemma (related to Lemma 2.1 of \cite{BBP10})
which will be used in Chapter \ref{sec:proof.Almost-inv}.

\begin{lemma} \label{pert+}
Let $ E $ be a closed subspace of ${\cal H}^0=L^2(\T^\es \times \T^d)$ and let
$E^s := E \cap {\cal H}^s $ 
for $s\geq s_0$. Let $R: E\to E$ be
a linear operator, satisfying, for some $s \geq s_1 \geq s_0$, $ \a \geq  0 $, 
\begin{align}
& \forall v \in E^{s_1} \, ,  \  \| Rv \|_{s_1} \leq \frac12 \| v \|_{s_1}  \label{Es1} \\
& \forall v \in E^s \, , \   \| Rv \|_{s} \leq \frac12 \| v \|_{s} + \alpha \| v \|_{s_1} \label{Ess} \, .
\end{align}
Then $ {\rm Id} +R $ is invertible as an operator of $E^{s_1}$, and 
\begin{align}
& \forall v \in E^{s_1} \, ,  \  \| ({\rm Id} +R)^{-1} v \|_{s_1} \leq 2 \| v \|_{s_1}  \label{IRs1} \\
& \forall v \in E^s \, , \  \   \| ({\rm Id} +R)^{-1} v \|_{s} \leq 2 \| v \|_{s} + 4 \alpha \| v \|_{s_1} \, .
\label{IRss} \end{align}
Moreover, assume that $R$ depends on the parameter $\l \in \wtilde{\Lambda}$ and satisfies also
\begin{align}
& \forall v(\l)  \in E^{s_1} \, ,  \  \| Rv \|_{\Lip, s_1} \leq \frac12 \| v \|_{\Lip, s_1}  \label{Es1lip} \\
& \forall v(\l) \in E^s \, , \  \ \| Rv \|_{\Lip, s} \leq \frac12 \| v \|_{\Lip, s} + \alpha \| v \|_{\Lip, s_1} \label{Esslip} \, .
\end{align}
Then
\begin{align}
& \forall v (\l) \in E^{s_1} \, ,  \  \| ({\rm Id} +R)^{-1} v \|_{\Lip, s_1} \leq 2 \| v \|_{\Lip , s_1}  \label{IRs1lip} \\
& \forall v(\l) \in E^s \, , \  \   \| ({\rm Id} +R)^{-1} v \|_{\Lip , s} \leq 2 \| v \|_{\Lip , s} + 4 \alpha \| v \|_{\Lip , s_1} \, .
\label{IRsslip} \end{align}
 
\end{lemma}

\begin{pf}
Since   $ E $ is a closed subspace of ${\cal H}^0 $,  each space $
 E^s = E \cap {\cal H}^s $, $ s \geq 0 $,  is complete.   
By \eqref{Es1} the operator $  {\rm Id} +R $ is invertible in $ E^{s_1} $ and \eqref{IRs1} holds. 
In order to prove \eqref{IRss}, let 
$ k = ( {\rm Id} + R)^{-1} v $ so that $ k = v - R k  $. Then
$$
\begin{aligned}
\| k \|_s \leq \|v \|_s + \| R k \|_s 
& \stackrel{\eqref{Ess}}  \leq \|v \|_s + \frac12 \| k \|_{s} + \alpha \| k \|_{s_1} \\
& \stackrel{\eqref{IRs1}}  \leq \|v \|_s + \frac12 \| k \|_{s} + \alpha 2 \| v \|_{s_1}
\end{aligned}
$$
and we deduce
$ \| k \|_s  \leq 2( \|v \|_s  + \alpha 2 \| v \|_{s_1} ) $ which is \eqref{IRss}.
The proof of \eqref{IRs1lip}  and \eqref{IRsslip} when we assume \eqref{Es1lip}  and \eqref{Esslip} follows the same line.
\end{pf}

\section{Off-diagonal decay of $ \sqrt{- \Delta + V(x)}$ }\label{sec:offDV} 

The goal of this section is to provide a direct proof of the fact  that 
the matrix which represents the operator\index{Multiplicative potential}  
$$ 
D_V = \sqrt{- \Delta + V(x)} \, , 
$$
defined in \eqref{def:DV}, 
in the exponential basis $ \{ e^{\ii j \cdot x} \}_{j \in \Z^d} $ has off-diagonal decay. 
This is required by the multiscale analysis performed in Chapter \ref{sec:multiscale}.
We compare $ D_V $ to the Fourier multiplier 
$ D_m = \sqrt{- \Delta + m} $ defined in  \eqref{def:Dm}.

\begin{proposition}\label{lemma-DVvsDm}
{\bf (Off-diagonal decay of $ D_V   - D_m $)}
There exists  a positive constant $ \Upsilon_s := \Upsilon_s (\| V \|_{C^{n_s}})  $, where 
$ n_s := \lceil s + \frac{d}{2}\rceil +1 \in \N $, such that 
\be\label{DeltaV2}
\big|  D_V   - D_m  \big|_{+,s} =
 \big| D_m^{\frac12} \big( D_V   - D_m \big) D_m^{\frac12}  \big|_s \leq \Upsilon_s \, . 
\ee
\end{proposition}

The rest of this section is devoted to the proof of Proposition \ref{lemma-DVvsDm}. 
In order to prove  that the matrix $ (A_i^j)_{i,j \in \Z^d} $ which represents, in the exponential basis 
$ \{ e^{\ii j \cdot x } \}_{j \in \Z^d}$,  
a linear operator $ A $ acting on a dense subspace of $ L^2 (\T^d) $,  has  polynomial off-diagonal decay, we shall use the following
criterium. Let  
$$
{\rm Ad}_{\pa_{x_k}}  := [ \partial_{x_k} , \,  \cdot \,  ] 
$$ 
denote the commutator with the partial derivative $ \partial_{x_k} $. 
Since the operator $ {\rm Ad}_{\pa_{x_k} }^{n} A $
 is  represented by the  matrix 
 $$  
 \big( \ii^n (i_k -j_k)^n A_i^j   \big)_{i,j \in \Z^d} \, ,
 $$  
it is sufficient that $ A $ and the operators $ {\rm Ad}_{\pa_{x_k} }^{n} A $, $ k = 1, \ldots, d $, for $ n $ large enough, 
 extend to  bounded operators in $ L^2 (\T^d) $.  

We shall use several times the following abstract lemma.

\begin{lemma} \label{norm-abs}
Let $(H, \langle \ , \  \rangle)$ be a separable 
Hilbert space with norm  $ \| u \| := \langle u , u  \rangle^{\frac12} $. Let 
$B : {\cal D}(B) \subset H \to H  $ be an unbounded symmetric operator 
with a  dense  domain of definition $ {\cal D}(B) \subset H $, 
satisfying:
\begin{itemize}
\item (H1)
There is $\beta>0$ such that  $  \la Bu, u\ra \geq \beta \| u \|^2$, $ \forall u \in {\cal D}(B) $.
\item (H2)  $B$ is invertible and $B^{-1} \in {\cal L}(H)$ is a compact operator. 
\end{itemize}
Let $A : {\cal D}(A) \subset H \to H $ be a symmetric linear operator, such that: 
\begin{itemize}
\item (H3)
${\cal D}(A)$, respectively in addition ${\cal D}(B^{\frac12} A B^{\frac12})$, contains ${\cal D}(B^p)$ for some $p \geq 1$.
\end{itemize}
Moreover we assume that: 
\begin{itemize}
\item (H4)
There is $\rho \geq 0$ such that
$  | \la Au , Bu \ra| \leq \rho \| u \|^2 $, $ \forall u \in {\cal D}(A) \cap {\cal D}(B) $.
\end{itemize}
Then $A$, respectively $ B^{\frac12} A B^{\frac12}$, can be extended to a bounded operator of $H$ 
(still denoted by $ A $, respectively $ B^{\frac12} A B^{\frac12} $) 
satisfying 
\be\label{A-A'}
\| A \|_{{\cal L} (H)} \leq \rho / \beta \, , \quad {\it respectively}  \ 
 \| B^{\frac12} A B^{\frac12} \|_{{\cal L} (H)} \leq \rho \, .
 \ee
\end{lemma}

\begin{pf}
The operator $B^{-1} \in {\cal L}(H) $ is compact and symmetric, and therefore 
there is an orthonormal basis 
$(\psi_k)_{k \geq 1} $ of $H$ of eigenfunctions of $B^{-1}$, i.e.    $B^{-1} \psi_k=\l_k \psi_k$, with eigenvalues 
$\l_k \in \R \backslash \{ 0 \}$, $(\l_k) \to 0$. 
Each  $ \psi_k $ 
is an eigenfunction of $ B $, i.e. 
$$ 
B \psi_k=\nu_k \psi_k \quad
{\rm  with \ eigenvalue }  \quad 
\nu_k :=\l_k^{-1} \, , \quad (\nu_k) \to \infty \, . 
$$
By assumption $(H1)$, each $\nu_k \geq \beta>0 $. 
Clearly  each eigenfunction $ \psi_k $  belongs to  the domain $ {\cal D}(B^p)$ of  $ B^p $ for any $ p \geq 1 $.

For any $N \geq 1$, we consider
the $N$-dimensional subspace 
$$
E_N:= {\rm Span}(\psi_1, \ldots , \psi_N)
$$ 
of $H$, 
and we denote by $\Pi_N$ the corresponding 
orthogonal projector on $E_N$. We have that $E_N \subset {\cal D}(B^p)$ for any $p \geq 1$, and therefore
assumption  $(H3)$ implies that  $E_N \subset {\cal D}(A)$.
\\[1mm]
{\sc Proof of $\| A \|_{{\cal L} (H)} \leq \rho / \beta $.} 
The operator $A_N := \Pi_N A_{|E_N}$ is a symmetric operator on the finite dimensional 
Hilbert space $(E_N, \la \ , \ \ra)$ and 
\be\label{normAN}
\| A_N \|_{{\cal L} (E_N)}= \max \big\{ |\l| \, ; \  \l \ \hbox{eigenvalue of} \ A_N \big\} \, .
\ee 
Let $\l$ be an eigenvalue of $A_N$  and $ u  \in E_N \backslash \{ 0 \}$ be an associated eigenvector,
i.e. $ A_N u = \lambda u $. 
Since $B(E_N) \subset E_N$, the vector $Bu $ is in $E_N$,
and we have
\be\label{steppino-int}
\l \la u , Bu \ra = \la  A_N u, Bu \ra = \la  \Pi_N A u, Bu \ra  = \la Au, Bu \ra \, .
\ee
Since $  \la u , Bu \ra $ is positive by ($H1$), by \eqref{steppino-int} and 
using assumption ($H4$) (note that $ u $ is in $ E_N \subset {\cal D}(A) \cap {\cal D}(B) $), we get
$$
|\l | \la u, Bu \ra = | \la Au, Bu \ra  | \leq \rho \| u \|^2 \, , 
$$
and, by assumption $(H1) $, we deduce  that $|\l| \leq  \rho / \beta$. 
By \eqref{normAN} we  conclude that, for any $ N \geq 1 $, 
\be\label{bound:ANL2}
\| A_N \|_{{\cal L} (E_N)} \leq \rho / \beta \, .
\ee
Defining  the subspace  $ E:=\bigcup_{N \geq 1} E_N$ of ${\cal D}(A)$, we deduce by \eqref{bound:ANL2} that 
\be \label{Auv}
  \forall (u,v) \in E \times E   \, , \  \la Au ,v \ra \leq \rho \beta^{-1} \|u \| \| v \| \, . 
\ee
Moreover, since $E$ is dense in $H $ (the 
$(\psi_k)_{k \geq 1} $ are an orthonormal basis of $ H$), the inequality \eqref{Auv} holds for all 
$(u,v) \in E \times H $, in particular  for all 
$(u,v) \in E \times {\cal D}(A)$. Therefore, since $A$ is symmetric, we obtain that 
\be\label{inclu22}
\forall u \in E \, , \forall v \in {\cal D}(A) \, , \  \la u, Av\ra \leq \rho \beta^{-1} \| u \| \| v \|  \, . 
\ee
By the density of $E$ in $H$, the inequality \eqref{inclu22} holds for all $ 
(u, v) \in H \times {\cal D}(A) $, and  we conclude that 
\be\label{concl-for-A}
\forall v \in {\cal D}(A)  \, , \quad \| Av \| \leq \rho \beta^{-1} \| v \|  \, . 
\ee
By continuity and \eqref{concl-for-A},  the operator 
$A$ can be extended to a bounded operator on the closure $ \overline{{\cal D} (A)} = H $ (that we still denote by $ A$) 
with operatorial norm $ \| A \|_{{\cal L} (H)} \leq \rho / \beta$, proving the
first estimate in \eqref{A-A'}. 
\\[1mm]
{\sc Proof of  $\| B^{\frac12} A B^{\frac12} \|_{{\cal L} (H)} \leq \rho $.}
The linear operator $ B^{\frac12}$ is defined on the basis $(\psi_k)_{k \geq 1}$ by setting 
$$
B^{\frac12} \psi_k := \sqrt{\nu_k} \psi_k \, . 
$$
Notice that, since in assumption $ (H3) $ we also require that, for some $p \geq 1$,  we have 
$ {\cal D}(B^p) \subseteq {\cal D}(B^{\frac12} A B^{\frac12}) $, we  have  
the inclusion $ E_N \subset {\cal D}(B^{\frac12} A B^{\frac12})$.
Since $B^{\frac12} (E_N) \subset E_N$, and the operators $B^{\frac12} $, $ A $ are symmetric on $ E_N \subset {\cal D} (A) $, 
then 
$A'_N:= B^{\frac12} \Pi_N A B^{\frac12}_{|E_N}$ is a symmetric operator of $E_N$. Let $\l'$ be an eigenvalue
of $A'_N$ and $ u'\in E_N \backslash \{ 0 \}$ be an associated eigenvector, i.e. $ A'_N u' = \l' u' $. 
Using that $y:= B^{\frac12} u'  \in E_N$ and $B y \in E_N$,
we obtain (recall that $ A_N = \Pi_N A_{|E_N} $)
\begin{align}\label{altro-caso-B12}
\l' \la u', Bu'\ra  = \la B^{\frac12}  A_N  B^{\frac12}u', Bu' \ra & =
\la B^{\frac12} A_N y , B^{\frac12} y \ra \nonumber \\
& = \la A_N y , By \ra =\la Ay, By\ra \, .
\end{align}
Since $  \la u' , Bu' \ra $ is positive by ($H1$), by \eqref{altro-caso-B12}
and assumption ($H4$) (note that $ u' $ is in $ E_N \subset {\cal D}(A) \cap {\cal D}(B) $), we get
\be\label{eccola:B12}
|\l'| \la u', Bu' \ra \leq \rho \| y \|^2 
= \rho  \la B^{\frac12} u' , B^{\frac12} u' \ra = \rho \la u', Bu' \ra \, .
\ee
Since $ \la u', Bu' \ra > 0 $ by ($H1$), we deduce by \eqref{eccola:B12} that 
 $|\l'| \leq \rho$ and thus $ \| A'_N \|_{{\cal L} (E_N)} \leq \rho $.  
Since $B^{\frac12}$ and $\Pi_N$ commute and   $E_N \subset {\cal D}(B^{\frac12}  A B^{\frac12})$, we have  
$$
A'_N= \Pi_N (A')_{|E_N} \, , \quad A' := B^{\frac12}  A B^{\frac12} \, , 
$$ 
and arguing by density as to prove \eqref{concl-for-A}  for $ A $, we deduce that 
for all $ v \in {\cal D}(A') $, $ \|A'v \| \leq \rho \| v \| $. 
Hence $A'$ can be extended to a bounded linear operator of $H$ with  norm 
$ \| A' \|_{{\cal L} (H)}  \leq \rho$.  This proves the second estimate in \eqref{A-A'}.
\end{pf}

In the sequel we shall apply Lemma \ref{norm-abs} with Hilbert space $H=L^2(\T^d)$ and an operator 
$ B \in \{ B_1, B_2, B_3 \} $ among 
\be\label{def:B1B2B3}
B_1 := D_m \, , \quad
B_2 := D_V \, , \quad 
B_3 := D_m + D_V \, ,
\ee
with  dense domain  
$$ 
{\cal D}(B_i) := H^1(\T^d) \, . 
$$ 
Notice that each operator $B_i $, $ i=1,2, 3 $,  satisfies assumption ($H1$) (recall \eqref{positive}) and $B_i^{-1}$  sends continuously $L^2(\T^d)$ 
into $H^1(\T^d)$ (note that $ \| D_V u \|_{L^2} \simeq \| u \|_{H^1_x}$ 
by \eqref{spaces-cal-Hs}-\eqref{equiv-norms-s}). Moreover, since $H^1(\T^d)$ is compactly embedded into $L^2(\T^d)$, 
each $B_i^{-1}$ is a compact operator of $H=L^2(\T^d)$, and hence also assumption $(H2)$ holds.

\begin{lemma} \label{L2norms} {\bf ($ L^2 $-bounds of $ D_V   - D_m $)}
Consider the linear operators 
\begin{align} 
& A_0 := 
D_V - D_m , \label{def:A0A0'L1} \\
& A'_0 := 
D_m^{\frac12} (D_V - D_m) D_m^{\frac12} \label{def:A0A0'L2}
\end{align} 
with domains $ {\cal D}(A_0) :=H^1(\T^d)$,  ${\cal D}(A'_0) :=H^2(\T^d)$. 
Then $ A_0 $ and  $ A'_0 $  can be extended to bounded linear operators of $L^2(\T^d)$ satisfying
\be\label{bounds:A0A0'}
\| A_0 \|_{{\cal L}(L^2)}, \| A'_0 \|_{{\cal L}(L^2)} \leq C(\| V \|_{L^\infty}) \, . 
\ee
\end{lemma}

\begin{pf}
Since $ D_m^2 = -\Delta +m. $ and $ D_V^2=-\Delta +V(x) $, we have 
$$
D_V(D_V-D_m) + (D_V-D_m)D_m=D_V^2-D_m^2= {\rm Op}(V(x) -m) \, ,
$$
where $ {\rm Op}(a)$ denotes the multiplication operator by the function $a(x) $. Hence, for all $ u \in H^2( \T^d) $,
$$
\begin{aligned}
\big| ( D_V(D_V-D_m) u, u )_{L^2} +  ( (D_V-D_m)D_m u ,u )_{L^2} \big| & = 
| ( (V(x)-m) u , u )_{L^2} | \\ 
& \leq (\|V \|_{L^\infty} + m) \| u \|^2_{L^2} \, ,
\end{aligned}
$$
which gives, by the symmetry of $D_V $ and $D_V-D_m$ 
that
\be\label{agg-int-L2}
\forall u \in H^2( \T^d) \, , \   
| ( (D_V-D_m) u, (D_V + D_m) u )_{L^2} | \leq
(\| V \|_{L^\infty} + m) \| u \|^2_{L^2} \, . 
\ee
Actually  \eqref{agg-int-L2}  holds  for all $u \in H^1(\T^d)$, 
by the density of $H^2(\T^d) $
in $H^1(\T^d)$ and the fact that $D_V$ and $D_m$ send continuously $H^1(\T^d)$ into $L^2(\T^d)$:
setting  $ A_0 = D_V - D_m $ as in \eqref{def:A0A0'L1} and $ B_3 = D_V + D_m $ as in \eqref{def:B1B2B3}, we have  
\be\label{ine-density}
\forall u \in H^1( \T^d) \, , \   |(A_0 u, B_3 u)_{L^2}|  \leq (\| V \|_{L^\infty} + m) \| u \|^2_{L^2} \, . 
\ee 
Thus the assumption $ (H4) $ of Lemma \ref{norm-abs} holds with $ A = A_0 $, $ B = B_3 $,  
Hilbert space 
$ H = L^2 (\T^d) $ and $ \rho = \| V \|_{L^\infty} + m $.  
Applying Lemma \ref{norm-abs} to $A_0$ and $B_3$ (also
assumption ($H3$)  holds since $ {\cal D}(A_0) = {\cal D} (B_3) =  H^1 (\T^d) $) 
we conclude that $A_0$ can be extended to a bounded 
operator of  $L^2(\T^d)$ with norm 
\be \label{A0boundL2}
 \| A_0 \|_{{\cal L}(L^2)} \leq C(\| V \|_{L^\infty}) 
 \ee 
(depending also on 
the constants $ m $ and $ \beta > 0 $ in \eqref{positive}).
This proves the first bound in \eqref{bounds:A0A0'}. 
Then, recalling the definitions of $ A_0, B_1, B_3 $ in \eqref{def:A0A0'L1}, \eqref{def:B1B2B3}, 
using  \eqref{ine-density} and 
\eqref{A0boundL2}, 
we also deduce that 
$$
\forall u \in H^1(\T^d) \, , \   | ( A_0 u , B_1 u )_{L^2} | = \frac{1}{2}| ( A_0 u , (B_3 - A_0) u )_{L^2} |
\leq C'(\| V \|_{L^\infty})  \| u \|_{L^2}^2  \, . 
$$
Thus the assumption $ (H4) $ of Lemma \ref{norm-abs} holds with $ A = A_0 $, $ B = B_1 $. 
Also assumption ($H3$)  holds since 
$$ 
{\cal D}(A_0) = {\cal D} (B_1) =  H^1 (\T^d) \quad 
{\rm and} 
\quad 
{\cal D} (B_1^{\frac12} A_0 B_1^{\frac12}) \supset  H^2 (\T^d) = {\cal D} (B_1^2)  \,.
$$
Therefore Lemma \ref{norm-abs} 
implies that $ A'_0 = B_1^{\frac12} A_0 B_1^{\frac12} $ can be extended to a bounded operator of  $L^2(\T^d)$,
 with operatorial norm   $  \| A_0' \|_{{\cal L}(L^2)} \leq C(\| V \|_{L^\infty}) $. 
 This proves the second bound in \eqref{bounds:A0A0'}. 
\end{pf}

\smallskip

Let $[A'_0]_i^j = (|i|^2+m)^{1/4} [D_V-D_m]_i^j (|j|^2+m)^{1/4} $, 
$i,j \in \Z^d$, denote the elements of the matrix representing the operator
$ A'_0 = D_m^{\frac12} (D_V - D_m) D_m^{\frac12} $ 
in the exponential basis. 
By Lemma \ref{L2norms}  we have that 
\be\label{diago-n=0}
|[A'_0]_i^j| \leq  \| A_0' \|_{{\cal L}(L^2)} \leq C(\| V \|_{L^\infty}) \, , \quad \forall i, j \in \Z^d \, . 
\ee
In order to prove also 
a polynomial off-diagonal decay for 
$ [A'_0]_i^j  $, $ i \neq j $,  we note that, for $ n \geq 1 $,  
\be\label{prop:st1}
{\rm Ad}_{\pa_{x_k} }^{n} A'_0 \quad 
{\rm is  \ represented \ by \ the \ matrix} \quad  
\Big( \ii^n (i_k -j_k)^n [A'_0]_i^j   \Big)_{i,j \in \Z^d} 
\ee
and then prove that $ {\rm Ad}_{\pa_{x_k} }^{n} A'_0 $ 
extends to a bounded operator in $ L^2 (\T^d) $ . 

\begin{lemma}  \label{estimwithad}
 {\bf ($ L^2 $-bounds of $ {\rm Ad}_{\pa_{x_k}}^n D_V    $)}
For any $ n \geq 1 $, $ k = 1, \ldots, d $
the operators  ${\rm Ad}_{\pa_{x_k}}^n D_V $ 
and $  D_m^{\frac12} ({\rm Ad}_{\pa_{x_k}}^n D_V) D_m^{\frac12}  $ can be 
extended to bounded operators of $ L^2 (\T^d) $: there exist  positive  constants $ C_n $ and $ C_n' $ 
depending on $\| V \|_{C^n}$  such that 
\begin{align}
& \big\| {\rm Ad}_{\pa_{x_k}}^n D_V \big\|_{{\cal L}(L^2)} \leq C_n \, , \label{claim} \\
& \big\| D_m^{\frac12} ({\rm Ad}_{\pa_{x_k}}^n D_V) D_m^{\frac12}  \big\|_{{\cal L}(L^2)} \leq C_n' \, .  \label{claim1}
\end{align}
\end{lemma}

\begin{pf}
We shall use 
 the following algebraic formulas: given linear operators $L_1$, $ L_2 $ we have
\begin{align} 
& \label{formula-ad1}
{\rm Ad}_{\partial_{x_k}} (L_1 L_2)= ( {\rm Ad}_{\partial_{x_k}} L_1) L_2 + L_1 ({\rm Ad}_{\partial_{x_k}} L_2) \,, \\
& \label{formula-adn}
{\rm Ad}_{\pa_{x_k}}^n ( L_1 L_2 ) = \sum_{n_1=0}^n  \begin{pmatrix} n  \\ n_1 \end{pmatrix} 
\big({\rm Ad}_{\pa_{x_k} }^{n_1} L_1 \big)
\big({\rm Ad}_{\pa_{x_k} }^{n-n_1} L_2 \big)  \, .  
\end{align}
We split the proof in two steps. 
\\[1mm]
{\bf 1st step.} We prove by iteration that, for all $ n \geq 1  $, there are constants $ C_n, C_n'' >  0 $ such that
\begin{align}
& \| ({\rm Ad}_{\pa_{x_k}}^n  D_V) D_V + D_V ({\rm Ad}_{\pa_{x_k}}^n D_V)\|_{{\cal L} (L^2)} \leq C''_n   \label{iter-ad1} \\
& \| {\rm Ad}_{\pa_{x_k}}^n  D_V\|_{{\cal L} (L^2)} \leq C_n  \, .  \label{iter-ad2} 
\end{align}
Clearly the estimates \eqref{iter-ad2} are \eqref{claim}. 
\\[1mm]
{\sc Initialization: proof of \eqref{iter-ad1}-\eqref{iter-ad2} for $n=1$.}  
Applying \eqref{formula-ad1} with $ L_1 = L_2 = D_V $ and since 
$ D_V^2 = -\Delta +V (x) $, we get 
\be\label{Stima:Ad}
({\rm Ad}_{\pa_{x_k}} D_V) D_V + D_V ({\rm Ad}_{\pa_{x_k}} D_V)   = {\rm Ad}_{\pa_{x_k}} (-\Delta +V(x))= 
{\rm Op}(V_{x_k}(x) ) 
\ee
where  $ {\rm Op}(V_{x_k}) $ is the multiplication operator by the function $ V_{x_k} (x) := (\pa_{x_k} V)(x) $. 
Hence \eqref{iter-ad1} for $ n = 1 $ 
holds with $C''_1=\|V\|_{C^1}$. 
In order to prove \eqref{iter-ad2} for $n=1$
we apply  Lemma \ref{norm-abs} to the operators 
$A_1 := {\rm Ad}_{\pa_{x_k}} D_V$ and $B_2=D_V$. Assumption ($H3$) holds because 
${\cal D}(A_1)=H^2(\T^d)={\cal D}(B_2^2) $.  
Note that, because of the $L^2$-antisymmetry of the operator $\pa_{x_k}$,  
if $L$ is  $L^2$-symmetric then so is $ {\rm Ad}_{\pa_{x_k}} L $. 
Hence $A_1 = {\rm Ad}_{\pa_{x_k}} D_V $ is symmetric. Also  assumption  ($H4$) holds because
\begin{align} \label{auto-dominio} 
\forall u \in H^2 (\T^d) \, , \ 
|( A_1 u , B_2 u )_{L^2}|  & = |( {\rm Ad}_{\pa_{x_k}} D_V u , D_V u )_{L^2}|  \\
& = \frac12 \big| \big(({\rm Ad}_{\pa_{x_k}} D_V) D_V + D_V ({\rm Ad}_{\pa_{x_k}} D_V) u,u \big)_{L^2} \big| \nonumber \\
&  \leq \frac{C''_1}{2} \|u\|^2_{L^2}  \, . \nonumber
\end{align}
by \eqref{iter-ad1} for $ n = 1 $.
Therefore Lemma \ref{norm-abs} implies that \eqref{iter-ad2} holds for $n=1$, for some constant 
$C_1$ depending on $\|V\|_{C^1}$. 
\\[1mm]
{\sc Iteration: proof of \eqref{iter-ad1}-\eqref{iter-ad2} for $n > 1$.}  
Assume by induction that \eqref{iter-ad1}-\eqref{iter-ad2} have been proved up to rank $ n-1  $.
We now prove them at rank $ n $. 
Applying \eqref{formula-adn} with $ L_1 = L_2 = D_V $ and since $ D_V^2 = -\Delta +V(x) $, we get 
\begin{align} 
& \big({\rm Ad}_{\pa_{x_k} }^{n} D_V \big) D_V + 
D_V   \big({\rm Ad}_{\pa_{x_k} }^{n} D_V \big) +
\sum_{n_1=1}^{n-1}  \begin{pmatrix} n  \\ n_1 \end{pmatrix} 
\big({\rm Ad}_{\pa_{x_k} }^{n_1} D_V \big)
\big({\rm Ad}_{\pa_{x_k} }^{n-n_1} D_V \big) \nonumber \\
& = {\rm Ad}_{\pa_{x_k} }^{n} (-\Delta + V(x))= {\rm Op}( (\pa_{x_k}^n V) ) \, . \label{commu-deca}
\end{align}
Let  
$$ 
T_n  := \big({\rm Ad}_{\pa_{x_k} }^{n} D_V \big) D_V + 
D_V   \big({\rm Ad}_{\pa_{x_k} }^{n} D_V \big) \, . 
$$ 
By \eqref{commu-deca}
 and using the inductive assumption \eqref{iter-ad2} at 
rank $ n-1  $, we obtain 
\be\label{Tn-bounds}
\|  T_n \|_{{\cal L}(L^2(\T^d))} \leq
 \| V \|_{C^n} +  \sum_{n_1=1}^{n-1} \begin{pmatrix} n  \\ n_1 \end{pmatrix} 
\| {\rm Ad}_{\pa_{x_k} }^{n_1} D_V \|_{ {\cal L}(L^2)} \| {\rm Ad}_{\pa_{x_k} }^{n-n_1} D_V \|_{{\cal L}(L^2)}
  \leq C''_n  
\ee
where the constant $C''_n$  depends on $\| V \|_{C^n}$. We have that
\be \label{nnn}
\forall u \in H^{n+1}(\T^d) \, , \  |({\rm Ad}_{\pa_{x_k} }^{n} D_V  u , D_V u)_{L^2}| = \frac{1}{2} |(T_n u ,u)_{L^2}| 
\stackrel{\eqref{Tn-bounds}} \leq \frac{C''_n}{2} \|u\|^2_{L^2} \, .  
\ee
We now apply Lemma \ref{norm-abs} to $A_n := {\rm Ad}_{\pa_{x_k} }^{n} D_V$ and $B_2=D_V$. Assumption
$ (H3)$ holds because 
${\cal D}(A_n)=H^{n+1}(\T^d)= {\cal D}(B_2^{n+1})$  
 by \eqref{equiv-norms-s}.
Assumption  ($H4$) holds by \eqref{nnn}. Therefore  the first inequality in 
\eqref{A-A'} implies that  $\|A_n\|_{{\cal L}(L^2(\T^d))} \leq  C_n (\| V \|_{C^n})$.
This proves \eqref{iter-ad2} at rank $ n $.
\\[2mm]
{\bf 2nd step.} {\sc Proof of \eqref{claim1} for any $ n \geq 1 $.}  
In order to prove that the operator  
$ D_m^{\frac12} \big({\rm Ad}_{\pa_{x_k} }^{n} D_V \big)   D_m^{\frac12}$ 
extends to  a bounded operator
of  $ L^2(\T^d)$ we  use Lemma \ref{norm-abs} with $ A = A_n = {\rm Ad}_{\pa_{x_k} }^{n} D_V  $ 
and $ B =  D_m $. 
Assumption ($H3$) holds because 
$$ 
{\cal D}(A_n) = H^{n+1}(\T^d)= {\cal D}(D_m^{n+1}) \quad  
{\rm and} \quad 
{\cal D} (D_m^{\frac12} A_n D_m^{\frac12}) \supset H^{n+2}(\T^d)= {\cal D}(D_m^{n+2})\,.
$$ 
Also Assumption (H4) holds because, by 
\eqref{nnn}, 
and the fact that $ A_n $ and 
$D_m-D_V$ are bounded operators of $L^2(\T^d)$ by  \eqref{iter-ad2} and Lemma \ref{L2norms}, we get 
$$
\begin{aligned}
\forall u \in H^{n+1}(\T^d) \, ,   \ |(A_n u, D_m u)_{L^2}| & \leq  |(A_n u , D_V u)_{L^2}| + |(A_n u , D_m u -D_V u)_{L^2}| \\ 
& \leq \rho_n \|u\|^2_{L^2} \, ,
\end{aligned}
$$
where $\rho_n$ depends on $\|V \|_{C^n}$. 
Then the second inequality in \eqref{A-A'} 
implies   \eqref{claim1}.
\end{pf}
\\[2mm]
{\bf Proof of Proposition \ref{lemma-DVvsDm} concluded.}
Recalling   \eqref{diago-n=0}  the  entries of 
the matrix  $ ([A'_0]_i^j)_{i,j \in \Z^d} $, which represents the operator
$A'_0 := D_m^{\frac12} (D_V - D_m) D_m^{\frac12}$ in the exponential basis, are bounded.  Furthermore,  
since each $\pa_{x_k}$, $ k =1, \ldots, d $,  commutes with $D_m$ and $D_m^{\frac12}$, we have that,
for any $ n \geq 1 $, 
$$
{\rm Ad}_{\pa_{x_k} }^{n} A'_0= D_m^{\frac12} ({\rm Ad}_{\pa_{x_k} }^{n} D_V) D_m^{\frac12} \, .
$$
By Lemma  \ref{estimwithad}, this operator can be extended to a bounded operator
on $ L^2 (\T^d) $ satisfying (see \eqref{claim1})
\be\label{prop:st2}
\|{\rm Ad}_{\pa_{x_k} }^{n} A'_0\|_{{\cal L} (L^2(\T^d))}
\leq C'_n (\|V\|_{C^n}) \, . 
\ee
By \eqref{diago-n=0}, \eqref{prop:st1}, \eqref{prop:st2} we deduce that 
$$
\begin{aligned}
& \forall n \geq 0 \, , \  \forall k =1, \ldots, d \, , \ 
\forall (i,j) \in \Z^d \times \Z^d \, , \\
& |i_k -j_k|^n |[A'_0]_i^j| \leq \|{\rm Ad}_{\pa_{x_k} }^{n} A'_0\|_{{\cal L} (L^2(\T^d))} \leq C'_n (\|V\|_{C^n}) 
\end{aligned}
$$
with  the convention ${\rm Ad}_{\pa_{x_k} }^{0} := {\rm Id}$. 
Hence, for $s \geq 0$ and 
$ n_s:=\lceil s + \frac{d}{2} \rceil +1 \geq s + \frac{d}{2} +1$, we deduce that the $s$-decay norm of $ A_0' $ satisfies 
$$
| A'_0 |_s^2 \leq 
K_0^2 + \sum_{\ell \in \Z^d \setminus \{0\}} \frac{K_{n_s}^2}{| \ell |^{2 n_s}} | \ell  |^{2s} \leq
K^2_{n_s} \sum_{\ell  \in \Z^d} \frac{1}{ \langle \ell  \rangle^{d+2}} =: \Upsilon_s^2
$$
for some constant $\Upsilon_s $ depending on $s$ and $\|V\|_{C^{n_s}}$ only. 
The proof of Proposition \ref{lemma-DVvsDm} is complete. \rule{2mm}{2mm}

\smallskip

With the same methods we also obtain the following propositions.
\begin{proposition}\label{prop:decay-1/2}
There exists  a positive constant $ C(s, \| V \|_{C^{n_s}})  $, where 
$ n_s := \lceil s + \frac{d}{2}\rceil +1 \in \N $, such that 
$$ 
|D_V^{\frac12} - D_m^{\frac12}|_s \leq C(s, \| V \|_{C^{n_s}}) \, .
$$ 
\end{proposition}

\begin{pf}
Since the proof closely the one of Proposition \ref{lemma-DVvsDm}
we shall indicate just the main steps.
One first proves that $ D_V^{\frac12} - D_m^{\frac12}  $ can be extended to a bounded operator of $ L^2 (\T^d )$, arguing 
as in Lemma \ref{L2norms}. Writing
$$
D_V^{\frac12} (D_V^{\frac12} - D_m^{\frac12} ) + (D_V^{\frac12} - D_m^{\frac12} ) D_m^{\frac12} = D_V - D_m \, ,
$$
using the symmetry of $D_V^{\frac12} $ and $ D_V^{\frac12} - D_m^{\frac12} $, 
the fact that
$ D_V - D_m $ is bounded on $ L^2 (\T^d )$ by Lemma \ref{L2norms} (see \eqref{bounds:A0A0'}), and the density of $ H^{1} $ in $H^{\frac12} $,
we deduce that 
$$
 \forall u \in H^{\frac12}( \T^d) \, , \quad   
| ( (D_V^{\frac12} - D_m^{\frac12}) u, (D_V^{\frac12} + D_m^{\frac12}) u )_{L^2} | \leq
C(\| V \|_{L^\infty}) \| u \|^2_{L^2} \, . 
$$
Then, arguing 
as in Lemma \ref{L2norms}, Lemma \ref{norm-abs} implies that $ D_V^{\frac12} - D_m^{\frac12} $ can be extended to a bounded operator of $ L^2 (\T^d)$ satisfying 
\be\label{DV12-DM12} 
\| D_V^{\frac12} - D_m^{\frac12} \|_{{\cal L}(L^2)} \leq C(\| V \|_{L^\infty}) \, .
\ee
Next one proves that, for all $ n \geq 1 $, $ k = 1, \ldots, d $,  the operators
${\rm Ad}_{\pa_{x_k}}^n D_V^{\frac12}$ can be extended to bounded operators on $ L^2 (\T^d)$ satisfying 
\be\label{claim:1/2}
\big\| {\rm Ad}_{\pa_{x_k}}^n D_V^{\frac12} \big\|_{{\cal L}(L^2)} \leq C_n(\| V \|_{C^n}) \, . 
\ee
It is sufficient to argue as in Lemma \ref{estimwithad}, applying \eqref{formula-ad1}-\eqref{formula-adn} to 
$ L_1 = L_2 = D_V^{\frac12} $ and $ D_V^{\frac12} D_V^{\frac12} = D_V $, and using that 
$ {\rm Ad}_{\pa_{x_k}}^n D_V $ are bounded operators of $ L^2 (\T^d) $ satisfying \eqref{claim}.
The proposition follows by \eqref{DV12-DM12} and \eqref{claim:1/2} as in the conclusion of the proof of
Proposition \ref{lemma-DVvsDm}.
\end{pf}

\begin{proposition}\label{prop:decay-prod}
There exists  a positive constant $ C(s, \| V \|_{C^{n_s}})  $, where 
$ n_s := \lceil s + \frac{d}{2}\rceil +1 \in \N $, such that 
\be\label{lafinala}
| D_m^{\frac12} D_V^{-\frac12}|_s \, , \quad  | D_V^{-\frac12}  D_m^{\frac12}|_s \leq  C(s, \| V \|_{C^{n_s}}) \, . 
\ee
\end{proposition}

\begin{pf} Writing
\begin{align*}
& D_V^{-\frac12}  D_m^{\frac12}  = {\rm Id} - D_V^{-\frac12} (D_V^{\frac12} - D_m^{\frac12}) \, , \\
& D_m^{\frac12} D_V^{-\frac12}   = {\rm Id} -   (D_V^{\frac12} - D_m^{\frac12})D_V^{-\frac12} \, , 
\end{align*}
the estimates \eqref{lafinala} follow by Proposition \ref{prop:decay-1/2}, \eqref{inter-AB-tutto}, 
and 
\be \label{ultima-ine}
|D_V^{-\frac12}|_s \leq C (s, \| V \|_{C^{n_s}}) \, .
\ee 
The bound  \eqref{ultima-ine} follows, as in the conclusion of 
the proof of Proposition \ref{lemma-DVvsDm}, by the fact that 
$ D_V^{-\frac12} $ is a bounded operator of $ L^2 (\T^d )$ and 
\be\label{claim:-1/2}
 \big\| {\rm Ad}_{\pa_{x_k}}^n D_V^{-\frac12} \big\|_{{\cal L}(L^2)} \leq C_n (\| V \|_{C^n}) \, , \quad
 \forall n \geq 1 \, , \ k = 1, \ldots, d  \, .  
\ee
The estimates \eqref{claim:-1/2} can be proved applying \eqref{formula-ad1}-\eqref{formula-adn} to 
$ L_1 = L_2 = D_V^{\frac12} $ and $ D_V^{\frac12} D_V^{-\frac12} = {\rm Id} $. For example, applying \eqref{formula-ad1},  
we get 
$$
\big({\rm Ad}_{\pa_{x_k}} D_V^{-\frac12}\big) = - D_V^{-\frac12} \big({\rm Ad}_{\pa_{x_k}} D_V^{\frac12} \big) D_V^{-\frac12}    
$$
which is a bounded operator on $ L^2 (\T^d)$ satisfying \eqref{claim:1/2}, and  $ D_V^{-\frac12} \in {\cal L}(L^2) $. 
By iteration, the estimate \eqref{claim:-1/2} follows for any $ n \geq 1 $. 
\end{pf}

\section{Interpolation inequalities}\label{sec:35}

We conclude this chapter with useful 
interpolation inequalities for the Sobolev spaces $  {\mathcal H}^s $ defined in 
\eqref{def:Hs}\index{Interpolation inequalities}. 

For any $ s \geq s_0 > (\es + d) / 2 $, we have the tame estimate\index{Tame estimates} for the product
\be\label{inter:pro-Lip}
\| f g  \|_{\Lip, s} \lesssim_s \| f \|_{\Lip, s} \| g \|_{\Lip, s_0} + \| f \|_{\Lip, s_0} \| g \|_{\Lip, s} \, . 
\ee
Actually we directly prove the improved tame estimate \eqref{intbasic+Lip} below, 
used in \cite{BBo10}, \cite{BBP10}.

\begin{lemma}\label{lem:int-impro}  
Let $ s \geq s_0 > (\es + d) / 2 $. Then 
\be  \label{intbasic+Lip}
\|uv\|_{\Lip,s} \leq C_0  \|u\|_{\Lip, s}  \|v\|_{\Lip,s_0} +  C(s) \|u\|_{\Lip, s_0} \|v\|_{\Lip, s} 
\ee
where the constant $ C_0 > 0 $ is independent of $ s \geq s_0  $.  
\end{lemma}

\begin{pf}
Denoting $ y := (\vphi, x) \in \T^b $, $ b = \es + d $, we  expand in Fourier series 
$$
 u (y) =\sum_{m \in \Z^b} u_m \, e^{\ii m \cdot y } \, ,  \qquad 
 v (y) =\sum_{m \in \Z^b} v_m \, e^{\ii m \cdot y } \, . 
$$
Thus 
\be\label{s1s2}
\| u \, v \|_s^2 = \sum_{m \in \Z^b} \Big| \sum_{k \in \Z^b} u_{k} v_{m-k} \Big|^2 
\langle m \rangle^{2 s} \leq
\sum_{m \in \Z^b}  \Big( \sum_{k \in \Z^b} |u_{k}|  |v_{m-k}| \Big)^2 \langle m \rangle^{2 s}  \leq 
S_1 + S_2 
\ee
where
\begin{align*}
& S_1 := 2\sum_{m \in \Z^b}  \Big( \sum_{|m| < |k| 2^{ 1/s} } | u_{k}|  |v_{m-k}| \Big)^2 
\langle m \rangle^{2 s} \\
& S_2 := 2\sum_{m \in \Z^b}  \Big( \sum_{|m| \geq |k| 2^{ 1/s}  } | u_{k}|  |v_{m-k}| \Big)^2 
\langle m \rangle^{2 s} \, . 
\end{align*}
The indices in the sum $ S_1 $ are restricted to $ |m| <   |k| 2^{ 1/s} $, thus
$ ( \langle m \rangle / \langle k \rangle )^s \leq 2 $, and, 
using Cauchy-Schwarz inequality, we deduce 
\begin{align}
S_1 & =  2\sum_{m \in \Z^b}  \sum_{|m| < |k| 2^{ 1/s}}
\Big(  | u_k | \langle k \rangle^{s}  
|v_{m-k}| \langle m-k \rangle^{s_0} \, 
 \frac{  \langle m \rangle^{ s} }{  \langle k \rangle^{ s} \langle m-k \rangle^{s_0} }  
 \Big)^2 \nonumber \\
& \leq  2\sum_{m \in \Z^b} \Big( \sum_{k \in \Z^b}
  | u_k |^2 \langle k \rangle^{2s}  
  |v_{m-k}|^2  \langle m-k \rangle^{2s_0} 4  \Big) \Big( \sum_{k \in \Z^b}
  \frac{1}{\langle m - k \rangle^{2 s_0}} \Big)  \nonumber \\ 
 & \leq  C(s_0) \| u \|_s^2 \| v \|_{s_0}^2 \, . \label{step1}
\end{align}   
On the other hand, the indices in the sum $ S_2 $  are restricted to
$ |k | \leq |m | 2^{- 1/s} $, and therefore  
$$
|m - k| \geq |m | - |k| \geq |m| ( 1 -   2^{-1/s}) 
$$
and $ ( \langle m \rangle / \langle m - k \rangle )^s \leq c(s)  $. As a consequence 
\begin{align}
S_2 & = 2 \sum_{m \in \Z^b}  \sum_{|k| \leq |m| 2^{- 1/s}}
\Big(  | u_k | \langle k \rangle^{s_0}  
|v_{m-k}| \langle m-k \rangle^{s} \, 
 \frac{  \langle m \rangle^{ s} }{  \langle k \rangle^{ s_0} \langle m-k \rangle^{s} }  
 \Big)^2 \nonumber \\
& \leq 2 \sum_{m \in \Z^b} \Big( \sum_{k \in \Z^b}
  | u_k |^2 \langle k \rangle^{2s_0}  
  |v_{m-k}|^2  \langle m-k \rangle^{2s} c(s)^2   \Big) \Big( \sum_{k \in \Z^b}
  \frac{1}{\langle k \rangle^{2 s_0}} \Big)  \nonumber \\ 
 & \leq  C(s) \| u \|_{s_0}^2 \| v \|_{s}^2 \, . \label{step2}
\end{align}   
By \eqref{s1s2} and the estimates \eqref{step1}-\eqref{step2} we deduce
$$
\| u \, v \|_s \leq C(s_0) \| u \|_{s} \| v \|_{s_0}  +  C(s) \| u \|_{s_0} \| v \|_{s}  \, .
$$
Recalling \eqref{Sobo-Lip}, the estimate \eqref{intbasic+Lip} follows.  
\end{pf}

In the case when $0\leq s \leq s_0$, the estimate of $ \| uv \|_s$ can be simplified.

\begin{lemma}\label{sleqs0}
Let $s_0 > (\es + d) / 2 $ and $0\leq s \leq s_0$. Then
\be  \label{tamesmalls}
\| u v  \|_{\Lip,s} \leq C(s_0) \| u \|_{\Lip,s_0} \| v \|_{\Lip,s} \, .
\ee
\end{lemma}

\begin{pf}
Denoting $b=(\es + d) / 2$, we have, as in the proof of Lemma \ref{lem:int-impro}, 
\begin{align} 
\| uv \|_s^2 
& \leq \sum_{m \in  \Z^b} \Big( \sum_{k \in \Z^b}\la  m \ra^{s} |u_k| |v_{m-k}|  \Big)^2
\nonumber 
\\ 
& \leq \sum_{m \in  \Z^b} \Big( \sum_{k \in \Z^b} 
 |u_k|\la k \ra^{s_0} |v_{m-k}|\la m -k \ra^{s} 
\frac{\la m \ra^{s}}{\la k \ra^{s_0} \la m -k \ra^{s} } \Big)^2 \nonumber \\
& \leq \sum_{m \in  \Z^b} 
\Big( \sum_{k \in \Z^b}   |u_k|^2 \la k \ra^{2s_0} |v_{m-k}|^2 \la m -k \ra^{2s} \Big)   
\Big( \sum_{k \in \Z^b} \frac{\la m \ra^{2s}}{\la k \ra^{2s_0} \la m -k \ra^{2s} } \Big)
\label{smalls-1}
\end{align}
by the Cauchy-Schwarz inequality.
We now use the following inequality: for $ 0 \leq s \leq s_0 $, 
\be \label{tamesmalls-2}
\forall m \in \Z^b \, , \ \forall k \in \Z^b \, , \
\frac{\la m \ra^s}{\la m-k \ra^s \la k \ra^{s_0}} \leq 2^s \Big(  \frac{1}{\la k \ra^{s_0}} + 
\frac{1}{\la m-k\ra^{s_0} }  \Big) \, . 
\ee
To prove \eqref{tamesmalls-2}, we distinguish two cases:
if $\la m-k \ra \geq \la m \ra / 2 $ then \eqref{tamesmalls-2} is trivial.
If $\la m-k \ra < \la m \ra / 2 $, then $\la k \ra >  \la m \ra / 2  > \la m-k \ra$ and, 
since $s_0-s \geq 0$,
$$
\la m \ra^s \leq 2^s \la k \ra^s = 2^s \frac{\la k \ra^{s_0}}{\la k \ra^{s_0-s}}
\leq 2^s \frac{\la k \ra^{s_0}}{\la m-k \ra^{s_0-s}} \, , 
$$
which implies \eqref{tamesmalls-2}.
By \eqref{tamesmalls-2} we deduce that, for any $ m \in \Z^b $, 
$$
\sum_{k \in \Z^b} \frac{\la m \ra^{2s}}{\la k \ra^{2s_0} \la m -k \ra^{2s} } 
\lesssim_s \sum_{k \in \Z^b}  
\frac{1}{\la k \ra^{2s_0}} + 
\sum_{k \in \Z^b} \frac{1}{\la m-k\ra^{2 s_0} }  \leq C(s_0)
$$
and therefore, by \eqref{smalls-1}, we obtain 
\be\label{uv-low}
\| uv \|_s \leq C(s_0) \| u \|_{s_0} \| v \|_s \, . 
\ee
Recalling \eqref{Sobo-Lip}, the estimate 
\eqref{tamesmalls} is a consequence of \eqref{uv-low}.
\end{pf}

As in any scale of Sobolev spaces with smoothing operators, the Sobolev norms 
$ \| \ \|_{s} $ defined  in \eqref{def:Hs} admit an
interpolation estimate. 

\begin{lemma}\label{Lem:loga} 
For any $ s_1 < s_2 $, $ s_1, s_2 \in \R $, and $ \theta \in [0,1] $, we have 
\be\label{inter-Sobo-Lip}
\| h \|_{\Lip, s} \leq 2 \| h \|_{\Lip, s_1}^\theta \| h \|_{\Lip, s_2}^{1-\theta} \, , \quad 
s := \theta s_1 + (1- \theta) s_2 \, . 
\ee
\end{lemma}

\begin{pf}
Recalling the definition of the Sobolev norm in \eqref{def:Hs}, we deduce,
by H\"older inequality, 
\begin{align*}
\| h \|_{s}^2  = \sum_{i \in \Z^b} 
| h_i |^2 \langle i \rangle^{2s} & = \sum_{i \in \Z^b} 
(| h_i |^2 \langle i \rangle^{2s_1})^{\theta} (| h_i |^2 \langle i \rangle^{2s_2})^{1- \theta} \\
& \leq 
\Big( \sum_{i \in \Z^b} | h_i |^2 \langle i \rangle^{2s_1} \Big)^{\theta} 
\Big( \sum_{i \in \Z^b} | h_i |^2 \langle i \rangle^{2s_2} \Big)^{1- \theta} \\
& = \| h \|_{s_1}^{2 \theta} \| h \|_{s_2}^{2 (1-\theta)} \, . 
\end{align*}
Thus $\| h \|_{s} \leq \| h \|_{s_1}^{ \theta} \| h \|_{s_2}^{ 1-\theta}  $, and, for a Lipschitz family of Sobolev functions, see \eqref{Sobo-Lip},  the inequality \eqref{inter-Sobo-Lip} follows. 
\end{pf}

As a corollary we deduce the following interpolation inequality\index{Interpolation inequalities}. 

\begin{lemma}\label{lem:double_interp}
	Let $ \alpha \leq a \leq b \leq \beta $ such that $ a + b = \alpha + \beta $. Then
\be\label{inte-ref}
	\| h \|_{\Lip, a} \| h \|_{\Lip, b} \leq 4 \| h \|_{\Lip, \alpha} \| h \|_{\Lip, \beta} \, .  
\ee
\end{lemma}

\begin{proof}
	Write $a= \theta  \alpha+(1- \theta )\beta$ and $b=\mu \alpha+(1-\mu)\beta$, where
	$$ 
	\theta =\frac{a-\beta}{\alpha-\beta}, \quad \mu =\frac{b-\beta}{\alpha-\beta}, \quad \theta +\mu = 1  \ .
	$$
	By the interpolation Lemma \ref{Lem:loga}, we get
	$$ 
	\norm{h}_{\Lip, a} \leq 2 \norm{h}_{\Lip, \alpha}^{\theta}
	\norm{h}_{\Lip, \beta}^{1-\theta}, 
\quad \norm{h}_{\Lip,b}\leq 2 \norm{h}_{\Lip, \alpha}^{\mu}\norm{h}_{\Lip, \beta}^{1-\mu} \, , 
	$$
	and \eqref{inte-ref} follows multiplying these inequalities.
\end{proof}

We finally recall the following  Moser tame estimates\index{Tame estimates}
for the  composition operator \index{Composition operator}
\be\label{comp:op}
(u_1, \ldots , u_p) \mapsto {\mathtt f}(u_1, \ldots , u_p)(\vphi, x) := f(\vphi, x, u_1(\vphi, x), \ldots , u_p(\vphi,x))  
\ee
induced by a smooth function $ f  $.  

\begin{lemma}\label{Moser norme pesate} {\bf (Composition operator)} \index{Moser estimates for composition operator}
Let $ f \in {\cal C}^{\infty}(\T^\es \times \T^d \times \R^p, \R )$. Fix $ s_0 > (d + \es)/2 $,
$ s_0 \in \N $.  
Given real valued functions $ u_i $, $1\leq i \leq p$, satisfying   $ \| u_i \|_{\Lip, s_0} \leq 1$,  
then, $ \forall  s \geq  s_0  $,  
\be\label{Moser-co}
  \| {\mathtt f}(u_1, \ldots , u_p) \|_{\Lip,s} \leq C(s, f ) \Big( 1 + \sum_{i=1}^p \| u_i \|_{\Lip,s} \Big) \, .
\ee
Assuming  also  that $ \| v_i \|_{\Lip, s_0} \leq 1$, $1\leq i \leq p$,  then
\begin{align}
 \|{\mathtt f}(v_1, \ldots , v_p) - {\mathtt f}(u_1, \ldots , u_p) \|_{\Lip,s} & \lesssim_{s,f} 
 \sum_{i=1}^p \| v_i-u_i \|_{\Lip ,s} + \label{Moser-Lip} \\
& \ \,  \quad  \Big(\sum_{i=1}^p \| u_i \|_{\Lip,s} + \| v_i \|_{\Lip,s} \Big) 
\sum_{i=1}^p \| v_i-u_i \|_{\Lip,s_0}   \, . \nonumber 
\end{align}
\end{lemma}

\begin{pf}
Let $ y := (\vphi, x) \in \T^{b} $, $ b = \es + d $. 
For simplicity of notation we consider only the case $ p = 1 $
	and denote $ (u_1, \ldots, u_p) = u_1 = u $.  
\\[1mm]	
{\sc Step 1.} For any $ s \geq 0 $,  
for any  function $ f \in C^\infty $, for all $   \| u \|_{s_0} \leq 1 $, 
\be\label{Moser-co-noLip}
\| {\mathtt f}(u) \|_{s} \leq C(s, f ) 
  \big( 1 + \| u \|_{s} \big) \, . 
\ee
This estimate was proved by Moser in \cite{Mo1}.  
 We propose here a 
different proof, following \cite{BB07}. Note that it is enough to prove 
\eqref{Moser-co-noLip} for $u \in C^\infty (\T^b)$. 
\\[1mm]
{\it Initialization:  \eqref{Moser-co-noLip} holds for  any $ s \in [0, s_0] $}. 
For a multi-index $\alpha=(\a_1, \ldots , \a_b) \in \N^b$, we denote
$ \pa_y^\a = \pa_{y_1}^{\a_1} \ldots \pa_{y_b}^{\a_b}  $ and we set 
$|\alpha|:= \a_1+ \ldots + \a_b$.
Recalling that $ s_0 $ is an integer, by the formula for the derivative
of a composition of functions,  
we estimate the Sobolev norm $  \| {\mathtt f}(u) \|_{s_0} $ by 
$$
C(s_0) \max \big\{ 
\|(\partial_y^\beta \partial_u^q f) (\cdot ,u) (\partial_y^{\alpha^{(1)}} u ) \ldots 
 (\partial_y^{\alpha^{(q)}} u) \|_0 \, ; \,  0 \leq q \leq s_0 \, ,\, |\alpha^{(1)}+ \ldots + \alpha^{(q)} + \beta| \leq s_0 
\big\} 
$$
where $ \alpha^{(1)},  \ldots,  \alpha^{(q)},  \beta \in \N^b $.
By the Sobolev embedding 
$ \| u \|_{L^\infty(\T^b)} \lesssim_{s_0} \| u \|_{s_0} $ we have 
$$
\| (\partial_y^\beta \partial_u^q f) (y,u(y)) \|_{L^\infty(\T^b)} \leq C(f, \| u \|_{s_0}) \, .
$$
Hence, in order to prove the bound 
$ \| {\mathtt f}(u) \|_{s_0} \lesssim_{f} 1 + \| u \|_{s_0} $, 
it is enough to check that,  for any $ 0 \leq q \leq s_0$, $|\alpha^{(1)} + \ldots + \alpha^{(q)}| \leq s_0$, the function 
\be \label{interms0}
v := (\partial_y^{\alpha^{(1)}} u)  \ldots (\partial_y^{\alpha^{(q)}} u ) \quad
{\rm satisfies} \quad \| v \|_0 \leq C(s_0) \| u \|_{s_0}^q \, . 
\ee 
Expanding in Fourier $ u(y) = \sum_{k \in \Z^b} u_k e^{\ii k \cdot y } $,  we have
\be\label{def:vtri}
\| v \|_0^2=\sum_{k \in \Z^b} \Big| \sum_{k_1+ \ldots + k_q=k} ( \ii k_1)^{\a^{(1)}} \ldots ( \ii k_q)^{\a^{(q)}} u_{k_1} \ldots u_{k_q} \Big|^2 \, ,
\ee
where, for $w\in \C^b $ and $\a \in \N^b$, we use the multi-index notation 
$w^\a:=\prod_{l=1}^b w_l^{\a_l}$.  Since $|\a^{(1)}| + \ldots +  |\a^{(1)}| \leq s_0$, the iterated Young inequality implies 
$$
| ( \ii k_1)^{\a^{(1)}} \ldots ( \ii k_q)^{\a^{(q)}} | \leq 
| k_1|^{|\a^{(1)}|} \ldots |k_q|^{|\a^{(q)}|} \leq \la  k_1 \ra^{s_0} + \ldots + \la  k_q \ra^{s_0} \, , 
$$
and therefore by \eqref{def:vtri} we have 
$$
\| v \|_0^2 \lesssim_q \sum_{i=1}^q S_i \quad {\rm with} \quad 
S_i :=\sum_{k \in \Z^b} \Big( \sum_{k_1+ \ldots + k_q=k} 
\la k_i \ra^{s_0} |u_{k_1}| \ldots |u_{k_q}| \Big)^2 \, . 
$$
An application of the Cauchy-Schwarz inequality, using  
$$
\sum_{k_2, \ldots, k_q\in \Z^b} \la k_2 \ra^{-2s_0} \ldots \la k_q \ra^{-2s_0}<\infty \, , 
$$
provides the bound  $S_1 \lesssim_{s_0} \| u \|_{s_0}^{2q} $. 
The bounds for the other $S_i$ are obtained similarly. This
proves \eqref{interms0} and therefore \eqref{Moser-co-noLip} for $ s = s_0 $.   
As a consequence,  \eqref{Moser-co-noLip} trivially holds also for any $ s \in [0, s_0 ] $: 
 in fact,  for all $ \| u \|_{s_0} \leq 1 $,  we have
$$
\| {\mathtt f}(u) \|_s \leq \| {\mathtt f}(u) \|_{s_0} \leq C(f) \leq C(f) (1 + \| u \|_s ) \, . 
$$	
{\it Induction: \eqref{Moser-co-noLip} holds for  $ s > s_0 $.}	
We  proceed by induction. Given some integer $k \geq s_0$, we  assume  
	that  \eqref{Moser-co-noLip}  holds for any $s \in [0,k]$, and for any  function $ f \in C^\infty $.  We are going to prove \eqref{Moser-co-noLip} for any $ s \in  ]k, k+1]$. Recall that
	\be\label{equivs} 
	\norm{u}_{s}^2  = 
         \sum_{k\in\Z^b} |u_k|^2 \langle k \rangle^{2s} \sim_s 
	\norm{u}_{L^2}^2 + \max_{i=1, \ldots, b} \norm{ \pa_{y_i }u }_{s-1}^2  \, . 
	\ee
Then    
	\begin{align}
	\| {\mathtt f}(u) \|_{s} &  \stackrel{\eqref{equivs}} {\sim_s} \| f(y, u(y)) \|_{L^2} + \max_{i=1, \dots, b} 
	\| \pa_{y_i}  (f(y, u(y))) \|_{s-1} \nonumber \\ 
	& \lesssim_{s,f} 
	1 + \max_{i=1, \dots, b }  
	\big( \|   ({\pa_{y_i} f})(y,  u(y))  \|_{s-1} + 
	\|   ({\pa_{u} f}) (y, u(y)) (\pa_{y_i} u)(y)  \|_{s-1}\big)
	\label{Herman-qf}  \, .  
	\end{align}
By the inductive assumption,
\be\label{Mo2}
 \|   ({\pa_{y_i} f})(y,  u(y))  \|_{s-1} \leq C(f) (1+ \| u \|_{s-1}) \, .
\ee
To estimate	$\|   ({\pa_{u} f}) (y, u(y)) (\pa_{y_i} u)(y)  \|_{s-1}$, we distinguish two cases:
\begin{itemize}
\item $1$st case: $k=s_0$. Thus $ s \in (s_0, s_0+1]$. Since  
$ s- 1  \in (s_0 -1, s_0 ] $ we apply Lemma \ref{sleqs0} obtaining  
\be\label{Mo3}
\|   ({\pa_{u} f}) (y, u(y)) (\pa_{y_i} u)(y)  \|_{s-1} \lesssim_{s_0} \|   ({\pa_{u} f}) (y, u(y))\|_{s_0} \| \pa_{y_i} u\|_{s-1}
\lesssim_{f} \|  u \|_{s} \, .  
\ee
\item $2$nd case: $k \geq s_0 + 1 $. For any $ s \in (k, k+1] $ 
we have $s-1 \geq s_0$, and,  by Lemma \ref{inter:pro-Lip} and the inductive assumption,
we obtain
$$
\|   ({\pa_{u} f}) (y, u(y)) (\pa_{y_i} u)(y)  \|_{s-1} \lesssim_{s,f} 
(1+\| u \|_{s-1}) \|  u\|_{s_0+1} + \| u  \|_{s} \, .
$$
By the interpolation inequality $ \| u \|_{s-1} \| u  \|_{s_0+1} \leq \| u \|_{s} \| u  \|_{s_0} $
(Lemma \ref{lem:double_interp}), and $s \geq s_0+1$,   we 
conclude that
\be\label{Mo4}
\|   ({\pa_{u} f}) (y, u(y)) (\pa_{y_i} u)(y)  \|_{s-1} \lesssim_{s,f} \| u  \|_{s} \, .
\ee
\end{itemize}	
Finally, by \eqref{Herman-qf}, \eqref{Mo2},  \eqref{Mo3},  \eqref{Mo4},   the estimate 
\eqref{Moser-co-noLip} holds for all $ s \in ]k,k+1] $. This concludes the iteration and the proof of \eqref{Moser-co-noLip}. 
\\[1mm]
{\sc Step 2. Proof of \eqref{Moser-co}. } 
In order to prove the Lipschitz estimate we write
\begin{align}
{\mathtt f}(v)(y) - {\mathtt f}(u)(y)
& = f(y, v(y) )  -  f(y, u(y)) \nonumber \\
& =  \int_0^1  (\pa_{u} f) (y, u (y) + \tau ( v -  u)(y)) ( v - u)(y) d \tau \, . \label{fv-fu}
\end{align}
Then 
\begin{align}
& \| {\mathtt f}(v) - {\mathtt f}(u) \|_{s}  \leq   \int_0^1  
\|(\pa_{u} f) (y, u  + \tau (v - u)) ( v - u) \|_{s }d \tau  \label{tame-diff} \\
& \stackrel{\eqref{inter:pro-Lip}} {\lesssim_s}  
\int_0^1   \|(\pa_u f) (y,  u  + \tau (v -  u)) \|_s \|  v - u \|_{s_0}   +
\| (\pa_u f) (y, u  + \tau ( v -  u)) \|_{s_0} \|  v - u \|_{s}  d \tau \, . \nonumber
\end{align}
Specializing \eqref{tame-diff} for $ v = u(\l_2)$ and $ u = u(\l_1) $,
using \eqref{Moser-co-noLip} and $  \| u \|_{\Lip, s_0} \leq 1$,  we deduce 
\be\label{Lip-co}
\| {\mathtt f}(u(\l_2)) - {\mathtt f} (u(\l_1) \|_s \lesssim_{s,f} 
\| u \|_{\Lip, s} |\l_2 - \l_1| \, , \quad  \forall \l_1, \l_2 \in \Lambda \, .
\ee
The estimates  \eqref{Moser-co-noLip} and \eqref{Lip-co} imply \eqref{Moser-co}.
\\[2mm]
{\sc Step 3. Proof of  \eqref{Moser-Lip}}.  
By \eqref{fv-fu},
\begin{align*}
\| {\mathtt f}( v) - {\mathtt f}( u) \|_{\Lip,s}
& \leq    \int_0^1  
\|(\pa_{u} f) (y,  u  + \tau ( v -  u)) ( v - u) \|_{\Lip,s }d \tau \\
& \stackrel{\eqref{inter:pro-Lip}} {\lesssim_s}    \int_0^1  
\|(\pa_{u} f) (y,  u  + \tau ( v -  u)) \|_{\Lip,s} \|  v - u \|_{\Lip,s_0 } \\ 
& \quad \quad \ \ \ +
\|(\pa_{u} f) (y,  u  + \tau ( v -  u)) \|_{\Lip,s_0} \|  v - u \|_{\Lip,s} 
\, d \tau 
\end{align*}
and, using \eqref{Moser-co} and $  \| u \|_{\Lip, s_0}, 
\| v \|_{\Lip, s_0} \leq 1 $,  
we deduce
\eqref{Moser-Lip} for $p=1$.

Estimates \eqref{Moser-co} and \eqref{Moser-Lip} for $p\geq 2$ can be obtained exactly in the same way.
\end{pf}

\chapter{Multiscale Analysis}\label{sec:multiscale}

The main result of this chapter is the abstract multiscale Proposition \ref{propmultiscale}, which provides  invertibility properties of finite dimensional restrictions of the Hamiltonian 
operator $ {\cal L}_{r, \mu} $ defined  in  \eqref{def:Lr}
for a large set of parameters  $ \l \in \Lambda $. This multiscale 
Proposition  \ref{propmultiscale} will be used
in Chapters \ref{sec:splitting} and \ref{sec:proof.Almost-inv}.

\section{Multiscale proposition}

Let  $ H := L^2(\T^d, \R) \times L^2(\T^d, \R) $ and consider the Hilbert spaces
\be \label{phase-space-multi}
(i) \ {\bf H} := H  \, , \qquad \quad  (ii)  \ { \bf H} :=H \times H  \, . 
\ee
Any vector $ h \in {\bf H} $ can be written as 
$$ 
(i) \ h =(h_1, h_2) \, , 
\quad \qquad (ii) \ h =(h_1,h_2,h_3,h_4) \, ,   \qquad h_i \in L^2 (\T^d, \R) \, .  
$$
On $ {\bf H} $, we define the linear operator   $ J $ as 
\begin{align}
J(h_1,h_2) & := ( h_2, - h_1) \qquad \qquad \qquad {\rm in \ case} \, (i) \label{simpl-rip} \\  
J(h_1,h_2,h_3,h_4) & := (h_2,-h_1, h_4, -h_3) \ \qquad {\rm in \ case} \, (ii) \, . \label{leftJA}
\end{align}
 Moreover, in case ($ii$), we also define the right action of $ J $ on $ {\bf H} $ as
\be\label{rightJA}
h J = (h_1,h_2,h_3,h_4)J := (-h_3,-h_4, h_1, h_2) \, . 
\ee
\begin{remark}
The definition \eqref{leftJA}, \eqref{rightJA} is motivated by the identification
of an operator $ a \in  {\cal L}(H_j, H) $ with $  (a^{(1)}, a^{(2)}, a^{(3)}, a^{(4)}) \in 
H \times H $ as in \eqref{ident-H-4-me}-\eqref{ident-H-4}.
Then the operators $ J a, a J \in  {\cal L}(H_j, H) $ are identified with 
the vectors of $ H \times H $ given in 
 \eqref{actionJ}, \eqref{actionJ-right}. This corresponds to the left and right 
 actions of $J$ in \eqref{leftJA}, \eqref{rightJA}. 
\end{remark}

We   denote by $\Pi_{{\mathbb S} \cup {\mathbb F}} $ 
the $L^2$-orthogonal projection
on the subspace $ H_{{\mathbb S} \cup {\mathbb F}} $ 
in $H$ defined in  \eqref{def:HE} or the analogous one in $H \times H$. Note that $\Pi_{{\mathbb S} \cup {\mathbb F}} $ commutes with the left (and right in  case \eqref{rightJA}) action 
of $ J $. We also denote 
$$ 
\Pi_{{\mathbb S} \cup {\mathbb F}}^\bot := 
{\rm Id} - \Pi_{{\mathbb S} \cup {\mathbb F}} = \Pi_{\mathbb G} \, ,
$$
see \eqref{taglio:pos-neg-0}. We fix  a constant $ \co \in (0,1) $ 
such that 
\be\label{diof:co}
| \bar \mu \cdot \ell + \mu_j + \co | \geq \frac{\gamma_0}{\langle \ell \rangle^{\tau_0}} \, , \quad \forall \ell \in \Z^\es \, , \  j \in {\mathbb S} \, .
\ee
By standard  arguments,  
condition \eqref{diof:co} is fulfilled by all $ \co  \in (0,1) $ except a set of  measure $ O(\gamma_0 )$. 
We explain  the purely technical 
role of the term $ \co \Pi_{\mathbb S} $ in  \eqref{op:i}-\eqref{def:Ar} below, 
in remarks \ref{rem:co} and \ref{rem:co-new} .

\begin{definition}\label{definition:Xr}
Given positive constants $ C_1,c_2 > 0 $, we define  the  class 
$  {\mathfrak C}(C_1,c_2) $ of $L^2$-self-adjoint operators acting on $ {\bf H} $, 
of the form, according to the cases ($i$)-($ii$) in \eqref{phase-space-multi},  
\begin{align}
& (i) \qquad X_r  :=   X_{r} (\e,\l,\vphi) = 
D_V +  \co \Pi_{\mathbb S}   + r (\e,  \l, \vphi) \label{op:i} \\
& (ii) \ \, \quad X_{r, \mu}  :=   X_{r, \mu} (\e,\l ,\vphi) = 
D_V +  \co \Pi_{\mathbb S} + \mu(\e,\l) {\cal J} \Pi_{{\mathbb S} \cup {\mathbb F}}^\bot  + r (\e,  \l, \vphi) \label{def:Ar}
\end{align}
defined for  $ \l \in \wtilde \Lambda \subset \Lambda $,  where $ D_V := \sqrt{ - \Delta + V(x) } $, $ \mu(\e,\l) \in \R $,
$ {\cal J} $ is the self-adjoint operator 
\be\label{def:opJ}
{\cal J} : H \times H  \to H \times H \, , \qquad 
h \mapsto {\cal J} h := J h J \, , 
\ee
(recall \eqref{leftJA}-\eqref{rightJA}) and such that  
\begin{enumerate}
\item 
\label{one-p} $ | r |_{\Lip, +,s_1}  \leq C_1 \e^2 $, for some $ s_1 > s_0 $, 
\item \label{mukinF}
$ |\mu (\e, \lambda) - \mu_k |_{\Lip} \leq C_1 \e^2 $ for some $ k \in {\mathbb F } $ (set defined in \eqref{taglio:pos-neg}), 
\item \label{assu:pos}
$ {\mathfrak d}_\l \Big(\frac{X_{r, \mu} }{1+\e^2 \l}\Big) \leq - c_2 \e^2 {\rm Id} $,  see the notation 
\eqref{partial-increase}.
\end{enumerate}
We assume that the non-resonance conditions 
\eqref{1Mel}, \eqref{2Mel+}-\eqref{2Mel}, \eqref{diof:co} hold. 

For simplicity of notation, in the  sequel  we shall denote by $ X_{r,\mu} $ also  the operator $ X_{r} $ in \eqref{op:i},
understanding that  $ X_r = X_{r, \mu} $ does not depend  on $ \mu $, i.e. $ {\cal J} = 0 $.  
\end{definition}

Note that the operator $ {\cal J} $ defined in \eqref{def:opJ} and 
$ \Pi_{{\mathbb S} \cup {\mathbb F} }^\bot = \Pi_{\mathbb G } $ commute.

\begin{remark}\label{rem:co}
The form of the operators $ X_r $, $ X_{r, \mu} $ in \eqref{op:i}, \eqref{def:Ar} is motivated  by the application
of the multiscale Proposition \ref{propmultiscale} to the operator 
$ {\cal L}_r $ in \eqref{applicazione-i} acting on $ H $, 
and  the operator $ {\cal L}_{r, \mu} $ in  \eqref{applicazione-ii} acting on $ H \times H $. 
We add  
the term  $  \co \Pi_{\mathbb S}  $ in \eqref{op:i}, \eqref{def:Ar} as a  purely  technical trick
to prove Lemma \ref{lem:Bel-inv}, see  Remark \ref{rem:co-new}.
\end{remark}

In the next proposition we prove invertibility properties of finite dimensional restrictions of the operator
\be\label{def:Lr}
 {\cal L}_{r, \mu} :=  J \om \cdot \partial_\vphi + X_{r, \mu} (\e, \l ) \, , \quad \om = (1+ \e^2 \l) \bar \om_\e   \, , 
\ee
for a large set of $ \l \in \Lambda $. 
 For $ N \in \N $ we define the subspace of trigonometric polynomials 
\be\label{def:EN-tr}
\begin{aligned}
& {\cal H}_N := \Big\{  u(\vphi, x) = \sum_{|(\ell, j)| \leq N} u_{\ell, j} e^{\ii (\ell \cdot \vphi + j \cdot x)} \, , \  u_{\ell,j}   \in \C^v  \Big\} \\ 
& {\rm where} \ \quad v :=\begin{cases}
2 \ \hbox{in case \eqref{phase-space-multi}-$(i)$} \\  4 \  \hbox{in case \eqref{phase-space-multi}-$(ii)$\, ,  } 
\end{cases} 
\end{aligned}
\ee
and we denote by  $ {\it \Pi}_N $ the corresponding $ L^2 $-projector: 
\be\label{def:projectorN}
u (\vphi, x) = \sum_{(\ell, j) \in \Z^\es \times \Z^d} u_{\ell, j} e^{\ii (\ell \cdot \vphi + j \cdot x)} \quad 
\mapsto \quad
{\it \Pi}_N u := \sum_{|(\ell, j)| \leq N} u_{\ell, j} e^{\ii (\ell \cdot \vphi + j \cdot x)} \, . 
\ee
The projectors $ {\it \Pi}_N $ satisfy the usual smoothing\index{Smoothing operators} estimates in Sobolev spaces: for any $ s, \b \geq 0 $,   
\begin{align}\label{smoothingS1S2}
 \| {\it \Pi}_N u \|_{s+\b} \leq N^\b \| u \|_{s} \, , 
\quad & \| {\it \Pi}_N^\bot u \|_s \leq N^{-\b} \| u \|_{s+\b} \, \\
\label{smoothingS1S2-Lip}
\| {\it \Pi}_N u \|_{\Lip, s+\b} \leq N^\b \| u \|_{\Lip, s} \, , 
\quad & \| {\it \Pi}_N^\bot u \|_{\Lip,s} \leq N^{-\b} \| u \|_{\Lip, s+\b} \, . 
\end{align}

We shall require that $ \omega $ satisfies the following  
quadratic Diophantine\index{Quadratic Diophantine condition} non-resonance condition. 

\begin{definition}\label{NRgamtau}
$ {\bf (NR)}_{\g, \t} $ 
Given  $ \g \in (0,1) $, $ \tau >  0 $, a vector $ \om \in \R^\es $
is $ {\bf (NR)}_{\g, \t} $ non-resonant, if, 
for any {\it non zero} polynomial $ P(X) \in \Z [X_1, \ldots, X_\es] $ of the form 
\be\label{NRom0}
P( X ) = n + \sum_{1\leq i\leq j \leq \es}  p_{ij} X_i X_j \, , \quad n, p_{ij} \in \Z \, , 
\ee  
we have
\be\label{NRom}
| P(\om)  | \geq \g \langle p \rangle^{- \t} \, , \quad \langle p \rangle := \max_{i,j=1, \ldots, \es}\{ 1, |p_{ij} | \} \,  . 
\ee
\end{definition}

The main result of this section is the following proposition. Set 
\be\label{def:varsigma}
\varsigma := 1 / 10 \, .
\ee
For simplicity of notation, 
 in the  next proposition $ {\cal L}_{r,\mu}$ also denotes  $ J \om \cdot \partial_\vphi + X_{r}  $,
understanding that  $ X_r = X_{r, \mu} $  in \eqref{op:i} does not depend  on $ \mu $.

\begin{proposition}\label{propmultiscale} {\bf (Multiscale)}\index{Multiscale proposition}
Let  $ \bar \om_\e \in \R^\es $ be  $(\gamma_1 , \tau_1 )$-Diophantine and satisfy 
property $ {\bf (NR)}_{\gamma_1, \tau_1 } $ in Definition \ref{NRgamtau} with $ \gamma_1, \tau_1 $
defined in \eqref{def:tau1}.
Then there are $ \e_0>0$, $ \tau' > 0 $,  $ \bar s_1 > s_0 $, 
$ \bar N \in \N $ (not depending on $X_{r,\mu}$ but possibly on the constants
$C_1 $, $c_2 $, $\gamma_0$, $\t_0$, $\gamma_1$,  $\t_1$) and $\t'_0>0$ (depending only on $\t_0$) such that the
following holds:

assume $ s_1 \geq \bar s_1 $ and take 
an operator $ X_{r, \mu} $ in $ {\mathfrak  C}(C_1, c_2) $ as in Definition \ref{definition:Xr},  
which is defined for all  $ \l \in \wtilde \Lambda $. 
Then for any $\e \in (0, \e_0)$, there are 
$ N (\e) \in \N $,  closed subsets $\Lambda ( \e; \cc , X_{r, \mu} ) \subset \wtilde \Lambda  $, $ \cc \in [1/2,1] $,
 satisfying 
\begin{enumerate}
\item  \label{list1-m}
$\Lambda (\e; \cc, X_{r, \mu} )\subseteq \Lambda (\e; \cc', X_{r, \mu} ) $, for all  $1/2 \leq \cc \leq \cc' \leq 1 $;
\item  \label{list2-m}
the complementary set 
$ \Lambda (\e; 1/2, X_{r, \mu})^c :=  \Lambda \setminus \Lambda (\e; 1/2, X_{r, \mu}) $ satisfies
\be\label{stima-Cantor-voluta}
| \Lambda (\e; 1/2, X_{r, \mu})^c \cap  \wtilde \Lambda | \lesssim \e \, ; 
\ee
\item  \label{list3-m}
if $ |r'  - r |_{+, s_1} +  |\mu' - \mu | \leq \d \leq \e^2 $,  then, for $ (1/2) + \sqrt{\d} \leq \cc \leq 1 $, 
\be\label{prop:spos}
|  \Lambda (\e; \cc, X_{r', \mu'} )^c \cap \Lambda (\e; \cc -\sqrt{\d}, X_{r, \mu}) \cap \wtilde \Lambda'| 
\lesssim \delta^\alpha \, ; 
\ee
\end{enumerate}
such that, 
\begin{enumerate}
\item  \label{item1-multiP} 
$ \forall  \bar N  \leq N < N(\e) $, $ \lambda \in \wtilde \Lambda   $,  the operator 
\be\label{to-have-RI}
[{\cal L}_{r, \mu}]_{N}^{2N} := {\it \Pi}_{N} ({\cal L}_{r, \mu})_{| {\mathcal H}_{2N}} 
\ee
has a right inverse $ \big([{\cal L}_{r, \mu}]_{N}^{2N}\big)^{-1}  : {\mathcal H}_N \to {\mathcal H}_{2N} $
satisfying, for all $ s \geq s_0 $, 
\be\label{est:Right+1}
\Big| \Big( \frac{ [{\cal L}_{r, \mu}]_{N}^{2N}}{1 + \e^2 \l}  \Big)^{-1}  \Big|_{\Lip, s} \,  \leq C(s)  
N^{\tau_0' +1}  \big(   N^{\loss s } + | r |_{\Lip, +,s} \big)  \, . 
\ee
Moreover, for all 
$ \l \in \wtilde \Lambda \cap \wtilde \Lambda' $, we have 
\be\label{MSLip1}
\big|  \big([{\cal L}_{r, \mu}]_{N}^{2N}\big)^{-1} -
\big([{\cal L}_{r', \mu'}]_{N}^{2N}\big)^{-1}  \big|_{s_1}
\lesssim_{s_1} N^{2(\tau_0'+ \varsigma s_1 ) + 1} \big(| \mu - \mu'  | N^2 + | r - r' |_{+,s_1} \big)
\ee
where $ \mu $, $ \mu' $, $ r $, $  r' $ are evaluated at fixed $ \lambda $. 

\item \label{item2-multiP} $ \forall N \geq N(\e)  $,  $ \lambda \in \Lambda (\e; 1, X_{r, \mu}) $,  
the operator 
\begin{equation}\label{def:Ln}
{\cal L}_{r, \mu, N}  := {\it \Pi}_{N} \big( J \om \cdot \partial_\vphi + X_{r, \mu} (\e, \l ) \big)_{| {\mathcal H}_{N}}  \, , 
\quad \om=(1+ \e^2 \l) \bar{\om}_\e \, , 
\end{equation}
is invertible and,  for all $ s \geq s_0 $, 
\be\label{Lip-sDM}
\Big|   \Big( \frac{{\cal L}_{r, \mu,N}}{1 + \e^2 \l}  \Big)^{-1} \Big|_{\Lip, s}  \leq C(s)  
N^{2(\tau' + \loss s_1) +3}  \big(   N^{\loss (s - s_1)} + | r |_{\Lip, +,s} \big) \, . 
\ee
Moreover for all $ \l \in \wtilde \Lambda \cap \wtilde \Lambda' $,
\be\label{MSLip2}
\big|  {\cal L}_{r, \mu, N}^{-1} - {\cal L}_{r', \mu', N}^{-1} \big|_{s_1}
\leq C(s_1) N^{2(\tau'+ \varsigma s_1 )+1} \big(| \mu - \mu'  | N^2 + | r - r' |_{+,s_1} \big)
\ee
where $ \mu $, $ \mu' $, $ r $, $  r' $ are evaluated at fixed $ \lambda $. 
\end{enumerate}
\end{proposition}

\begin{remark}
The measure 
of the set $\Lambda (\e; 1/2, X_{r, \mu} )^c $ is 
smaller than $ \e^p $, for any $ p $, at the expense of taking a larger constant $ \tau' $, see Remark \ref{rem:meas1}. 
We have written $ \e $ in \eqref{stima-Cantor-voluta}    for definiteness. 
\end{remark}

\begin{remark} \label{stablambda}
Properties \ref{list1-m}-\ref{list3-m} for the sets $\Lambda (\e; \cc, X_{r, \mu} )$ are stable under finite intersection:
if 
$$
(\Lambda^{(1)} (\e; \cc, X_{r, \mu} )), \ldots , (\Lambda^{(p)} (\e; \cc, X_{r, \mu} ))
$$ 
are families of closed
subsets of $\Lambda$ satisfying \ref{list1-m}-\ref{list3-m}, then the family $\dps (\bigcap_{1\leq k \leq p} 
\Lambda^{(k)} (\e; \cc, X_{r, \mu} ))$ still satisfies these properties. 
\end{remark}

The proof of Proposition \ref{propmultiscale} is based on the multiscale analysis of the papers 
\cite{BB12}, \cite{BBo10}, \cite{BCP} for quasi-periodically forced nonlinear wave and Schr\"odinger equations, but it
is more complicated  in the present autonomous setting. 
The proof of the multiscale Proposition \ref{propmultiscale} is given in the next sections \ref{sec:mult1}-\ref{sec:mult5}.

\section{Matrix representation}\label{sec:mult1}

We  decompose the operator $ {\cal L}_{r, \mu} $ 
in \eqref{def:Lr}, with $ X_{r, \mu} $ as in  \eqref{def:Ar} or \eqref{op:i} (in such a case
we mean that $ {\cal J} = 0 $), as 
\be
\begin{aligned}\label{split-Lom}
&  {\cal L}_{r, \mu} = {\cal D}_\om+ {\cal T} \, , \\
 & {\cal D}_\om := 
 J \om \cdot \partial_\vphi + D_m + \mu {\cal J}  \, , \\
&  {\cal T} := D_V-D_m - \mu {\cal J}  \Pi_{\mathbb S \cup \mathbb F} +  r + \co \Pi_{\mathbb S}  
\end{aligned}
\ee
where 
$ D_V := \sqrt{- \Delta + V(x) } $ is defined in  \eqref{def:DV},  
$ D_m := \sqrt{- \Delta + m} $  in  \eqref{def:Dm}, 
and $ \omega = (1+ \e^2 \l) \bar \om_\e $. 
By Proposition \ref{lemma-DVvsDm} the matrix which represents the operator $ D_V - D_m $
in the exponential basis has off-diagonal decay.

\smallskip

In what follows we   identify a linear operator $ A (\vphi) $, $ \vphi \in \T^\es $,  acting on functions
$ h(\vphi, x )$,  
with the  infinite dimensional matrix 
$ (A^{\ell,j}_{\ell',j'})_{\{ (\ell, j),  (\ell',j') \in \Z^{b} \} } $ of $ 2 \times 2 $ matrices $ A^{\ell,j}_{\ell',j'} $ 
in case \eqref{phase-space-multi}-($i$), respectively  $4\times 4$  in case \eqref{phase-space-multi}-($ii$),
 defined by the relation \eqref{matrix-A-space-time}.
In this way, the  operator  ${\cal L}_{r, \mu}  $ in \eqref{split-Lom} is represented by the infinite dimensional Hermitian matrix 
\be\label{Aomega}
{\mathtt A}(\e, \lambda) := {\mathtt A} (\e, \lambda; r) := {\mathtt D}_\om  +  {\mathtt T} \,  ,    
\ee
of $ 2 \times 2 $ matrices in case ($i$), resp. $ 4 \times 4 $ matrices in case ($ii$),  
where the diagonal part is, in case ($i$), 
\be
\begin{aligned}\label{diagD}
& \qquad \quad {\mathtt D}_\om = {\rm Diag}_{i \in \Z^b}
\begin{pmatrix}
 \langle j \rangle_m &  \ii \om \cdot \ell  \\ 
- \ii \om \cdot \ell  &
 \langle j \rangle_m  
\end{pmatrix} \, , \\
&  i := (\ell ,j) \in \Z^{b}  := \Z^\es \times \Z^d  \,, \quad \langle j \rangle_m := \sqrt{|j|^2 + m } \, .  
\end{aligned}
\ee
In case ($ii$), recalling the definition of $ {\cal J} $  in \eqref{def:opJ}, of the left and right action of 
$ J $ in \eqref{leftJA}-\eqref{rightJA}, 
and choosing the basis of $ \C^4 $
$$
\begin{aligned}
 \{ f_1, f_2, f_3, f_4 \} := & \,  \{ e_3- e_2, e_1 + e_4, e_1 - e_4, e_2 + e_3\} \, , \\
& \ e_{\mathfrak a} := (0, \ldots, \underbrace{1}_{{\mathfrak a}-th}, \ldots,0 ) \in \C^4 \, ,
\end{aligned} 
$$
we have 
\be\label{diag-secondo-caso}
{\mathtt D}_\om = {\rm Diag}_{(\ell,j) \in \Z^b}
\begin{pmatrix}
  \langle j \rangle_m - \mu &  \ii \om \cdot \ell  & 0 & 0  \\
 -  \ii \om \cdot \ell  &   \langle j \rangle_m - \mu   &  0 & 0  \\
 0 & 0  &    \langle j \rangle_m + \mu &  \ii \om \cdot \ell \\
 0 & 0  &   - \ii \om \cdot \ell &  \langle j \rangle_m + \mu
\end{pmatrix} \, .
\ee 
The off-diagonal matrix  
\be\label{Tmatrix}
{\mathtt T} := ( {\mathtt T}_i^{i'})_{i \in \Z^b , i'\in \Z^b} \, , \ \ b = \es + d \, , 
\quad  {\mathtt T}_i^{i'} :=   (D_V-D_m)_j^{j'}  - \mu [{\cal J} \Pi_{\mathbb S \cup \mathbb F}]_j^{j'} + r_i^{i'}\, , 
\ee
where $ {\mathtt T}_i^{i'} $ are $ 2 \times 2 $, resp. $ 4 \times 4 $ matrices, satisfies, by 
\eqref{DeltaV2}, Lemma \ref{pisig},  
property \ref{one-p} of Definition \ref{definition:Xr},  
\be\label{off-diago:T0}
| {\mathtt T} |_{+,s_1} \leq C(s_1) \, . 
\ee
Note that $ ( {\mathtt T}_i^{i'})^* = {\mathtt T}_{i'}^{i} $. 
Moreover, since the operator 
$ {\cal T} = {\cal T}(\vphi) $ in \eqref{split-Lom} 
 is a $ \vphi $-dependent family of operators acting on $ {\bf H }$ (as
$ r $ defined in \eqref{op:i}-\eqref{def:Ar} which is the only $ \vphi $-dependent operator),  
the matrix $ {\mathtt T} $ is  {\it T\"oplitz} in $ \ell $, 
namely $ {\mathtt T}_i^{i'} = {\mathtt T}_{\ell,j}^{\ell',j'} $ depends only on the indices $ \ell - \ell', j, j'  $.
We introduce a further index $ \mathfrak a \in {\mathfrak I } $ where
\be\label{def:indice a}
{\mathfrak I } :=  \{1,2\} \  \text{in case \eqref{phase-space-multi}-}(i) \, ,  \qquad {\mathfrak I } := \{1,2, 3, 4\} \ \text{in case 
\eqref{phase-space-multi}-}(ii) \, , 
\ee
to distinguish the matrix elements of each $ 2 \times 2 $, resp. $ 4 \times 4 $, matrix
$$
{\mathtt T}_i^{i'} := \big( {\mathtt T}_{i, \mathfrak a}^{i', \mathfrak a'} \big)_{\mathfrak a, \mathfrak a' \in {\mathfrak I}} \, . 
$$
Under the unitary change of variable  (basis of eigenvectors)
\be\label{def:U}
U := {\rm Diag}_{(\ell, j) \in \Z^b} 
\frac{1}{\sqrt{2}} 
\begin{pmatrix}
 1 & 1    \\
 \ii &  -\ii   \\   
\end{pmatrix}
\ee
the matrix $ {\mathtt D}_\om $ in \eqref{diagD} becomes completely diagonal 
\be\label{diag-caso1}
D_\om  := U^{-1} {\mathtt D}_\om  U  = 
{\rm Diag}_{(\ell, j) \in \Z^b}
\begin{pmatrix}
\langle j \rangle_m -  \om \cdot \ell   & 0    \\
 0 &   \langle j \rangle_m  +  \om \cdot \ell    \\   
\end{pmatrix}
\ee
and, under the unitary transformation 
\be\label{MATU2}
U := {\rm Diag}_{(\ell, j) \in \Z^b}\frac{1}{\sqrt{2}}
\begin{pmatrix}
1 & 1  & 0 & 0  \\
 \ii &  - \ii   &  0 & 0  \\
 0 & 0  &   1  &  1 \\
 0 & 0  &   \ii &  - \ii 
\end{pmatrix} \, ,  
\ee 
the matrix in \eqref{diag-secondo-caso} becomes completely diagonal 
\begin{align}
& D_\om  := U^{-1} {\mathtt D}_\om U = {\rm Diag}_{(\ell, j) \in \Z^b} \nonumber \\
& \begin{pmatrix}
 \langle j \rangle_m  - \mu - \om \cdot \ell  \! \! & \! \!  0 \! \!  & \! \! 0 \! \! & \! \! 0  \\
0 \!  \!&  \! \!  \langle j \rangle_m - \mu + \om \cdot \ell  \!  \! &  \! \! 0 \! \!& \! \! 0   \\
 0 \! \!& \! \! 0  \! \! &  \! \!  \langle j \rangle_m  + \mu -  \om \cdot \ell  \!\! & \!  \! 0  \\
0 \! \!& \! \! 0  \! \! & \! \!  0 \! \! & \! \!   \langle j \rangle_m + \mu +  \om \cdot \ell  
\end{pmatrix} \, . \label{diag-caso2}
\end{align}
Under the unitary change of variable $ U $ in \eqref{def:U} in case \eqref{phase-space-multi}-($i$), 
\eqref{MATU2} in case \eqref{phase-space-multi}-($ii$), 
the hermitian matrix $ \mathtt A $ in \eqref{Aomega} transforms in the hermitian matrix 
\be\label{unitary-A}
A(\e, \lambda) := A := U^{-1} {\mathtt A} U = D_\om + T  \, , \qquad T_i^{i'} := U^{-1} {\mathtt T}_i^{i'} U \, , 
\ee
where the off-diagonal term $ T  $ satisfies, by \eqref{off-diago:T0}, 
\be\label{off-diago:T}
 | T |_{+,s_1} \leq C(s_1) \, . 
\ee
We introduce the one-parameter family of infinite dimensional matrices 
\be\label{Atranslated}
A(\e, \lambda,\teta) := A(\e, \lambda) + \theta Y :=   D_\om  +  \theta \, Y +  T \, , \quad \theta \in \R \, , 
\ee
where
\be
\begin{aligned}\label{def:Y}
& Y := {\rm Diag}_{i \in \Z^b}
\begin{pmatrix}
- 1 & 0  \\ 
0  & 1 \end{pmatrix} \  \text{in case }\eqref{phase-space-multi}{\rm -}(i) \, , \\
&  Y  := {\rm Diag}_{i \in \Z^b} 
\begin{pmatrix}
- 1 &  0   & 0 & 0  \\
 0   &  1  &  0 & 0  \\
 0 & 0  & -   1 &  0  \\
 0 & 0  &    0  & 1 
\end{pmatrix} \  \text{in case }\eqref{phase-space-multi}{\rm-}(ii) \, .
\end{aligned}
\ee
The reason for adding  $ \theta Y $ is that,  translating  the time Fourier indices
$$ 
(\ell ,j) \mapsto ( \ell + \ell_0 ,  j )
$$
in $ A(\e, \lambda) $, gives $ A (\e, \lambda, \teta)$ with $ \teta = \om \cdot \ell_0 $. Note that 
the matrix $ T $  remains unchanged under translation because it is T\"oplitz with respect to $ \ell $.

\begin{remark}\label{def:r=0}
The matrix $ A(\e, \lambda,\teta)  := A(\e, \lambda,\teta; r )$ in  \eqref{Atranslated}
represents  the $L^2$-self-adjoint  operator 
\be\label{Lrmutheta}
{\cal L}_{r,\mu} (\theta):= J \om \cdot \partial_\varphi + \ii \theta J + D_V + 
\mu  {\cal J} \Pi^\bot_{{\mathbb S} \cup {\mathbb F}} + {\co} \Pi_{\mathbb S} + r  \, .
\ee
In Section \ref{sec:Green} we shall denote by $ {\breve A}(\e, \lambda,\teta) := A(\e, \lambda,\teta; 0 ) $
the matrix which represents 
\be\label{L0mu}
{\cal L}_{0,\mu} (\theta):= J \om \cdot \partial_\varphi + \ii \theta J + D_V + 
\mu  {\cal J} \Pi^\bot_{{\mathbb S} \cup {\mathbb F}} + {\co} \Pi_{\mathbb S} \, .
\ee
Note that 
 ${\cal L}_{0, \mu} (\theta)  $ is independent of $ \vphi $, thus $ {\breve A}(\e, \lambda,\teta)$ is diagonal in $ \ell \in \Z^\es $. 
\end{remark}

The eigenvalues 
of the $ 2 \times 2 $ matrix $ D_\om + \theta Y $, with  $ D_\om $
 in \eqref{diag-caso1} and $ Y $ defined in \eqref{def:Y}-case ($i$), 
 resp. $4 \times 4 $ matrix  with $ D_\om $  in \eqref{diag-caso2} and 
 $ Y $ defined in \eqref{def:Y}-case ($ii$), are 
\be\label{diao-1-caso}
d_{i, \mathfrak a} (\theta) = 
\begin{cases}
 \langle j \rangle_m  -  (\om \cdot \ell + \teta) \qquad \ {\rm if} \  \mathfrak a = 1 \cr
    \langle j \rangle_m  +  (\om \cdot \ell + \teta)  \ \qquad {\rm if} \ \mathfrak a = 2 \, , 
\end{cases}
\ee
and, in the second case, 
\be\label{diao-2-caso}
d_{i, \mathfrak a} (\theta) = 
\begin{cases}
\langle j \rangle_m  - \mu -  (\om \cdot \ell + \teta) \qquad {\rm if} \ \mathfrak a  = 1 \cr
 \langle j \rangle_m  - \mu +  (\om \cdot \ell + \teta) \qquad {\rm if} \ \mathfrak a = 2 \cr
 \langle j \rangle_m  + \mu -  (\om \cdot \ell + \teta) \qquad {\rm if} \ \mathfrak a = 3 \cr
   \langle j \rangle_m  + \mu +  (\om \cdot \ell + \teta) \qquad {\rm if} \ \mathfrak a = 4 \, .
\end{cases}
\ee

\smallskip

The main goal of the following sections is to prove polynomial off-diagonal decay for the inverse
of  the $ |{\mathfrak I}|(2N+1)^b $-dimensional (where $ |{\mathfrak I}| = 2 $, resp. $ 4 $,
in case \eqref{phase-space-multi}-($i$), resp. \eqref{phase-space-multi}-$ (ii)$) sub-matrices of $ A(\e, \l, \theta) $ centered at $ (\ell_0, j_0) $ denoted by
\be\label{ANl0}
A_{N,\ell_0, j_0}(\e, \l, \theta) := A_{|\ell - \ell_0| \leq N, |j - j_0| \leq N}(\e, \l, \theta) 
\ee  
where 
\be\label{supeuc}
|\ell | :=  \max \{| \ell_1|, \ldots, |\ell_\es|\} \, ,  \ \  |j| := \max \{|j_1|, \ldots, |j_d|\} \, .
\ee 
If $ \ell_0 = 0 $ we use the simpler notation
\be\label{ANj0}
A_{N, j_0}(\e, \l, \theta) :=  A_{N,0,j_0} (\e, \l, \theta) \, .
\ee
If also $ j_0 = 0 $,  we simply write
\be\label{notation:AN0}
A_{N}(\e, \l, \theta) :=  A_{N,0} (\e, \l, \theta) \, ,
\ee
and, for $ \teta = 0 $, we denote
\be\label{theta=0}
A_{N,j_0}(\e, \l) :=  A_{N,j_0} (\e, \l, 0) \, , \quad A_{N}(\e, \l) :=  A_{N,0} (\e, \l, 0) \, .   
\ee
\begin{remark}\label{rem:ANL}
The matrix $ A_{N}(\e, \l) $ in \eqref{theta=0} represents the truncated self-adjoint operator
$$
{\it \Pi}_N ( {\cal L}_{r, \mu} )_{| {\mathcal H}_N} = {\it \Pi}_N ( J \om \cdot \pa_\vphi + X_{r, \mu} )_{| {\mathcal H}_N}  = 
{\it \Pi}_N ( J \om \cdot \pa_\vphi + D_V + \mu  {\cal J} \Pi_{\mathbb S \cup \mathbb F}^\bot +  r )_{| {\mathcal H}_N}
$$
where $ {\mathcal H}_N $ is defined in \eqref{def:EN-tr} and $ {\it \Pi}_N $ in \eqref{def:projectorN}. 
\end{remark}

We have the following crucial {\it covariance} property 
\be\label{shifted}
A_{N, \ell_1,  j_1 } (\e, \l, \theta) = A_{N, j_1} (\e, \l , \teta +  \om \cdot \ell_1 ) \, . 
\ee

\section{Multiscale step}\label{sec:MS}

The main result of  this section is  the multiscale step Proposition \ref{propinv} which  is a variant of that
proved  in \cite{BBo10}. 
The constant $ \varsigma \in (0,1)$  is fixed and $\tau'>0$, $\Theta \geq 1$
are  real  parameters, on which we shall impose some condition in 
 Proposition \ref{propinv}. 

Given $ \Omega, \Omega' \subset  E \subset \Z^b \times {\mathfrak I} $, 
where $ {\mathfrak I } $ is defined in \eqref{def:indice a},  we define 
$$
{\rm diam}(E) := \sup_{k,k' \in E} |k-k'| \, , \qquad 
 {\rm d}(\Omega, \Omega') := \inf_{k \in \Omega, k' \in \Omega'} |k-k'|  \, , 
$$
where, for $ k= (i,\mathfrak a) $, $ k' := (i', \mathfrak a') \in \Z^b \times \fracchi $,  we set 
$$ 
|k - k' | := \begin{cases} 
1  \quad \qquad \qquad \, \text{if} \ i = i' \, , \mathfrak a \neq \mathfrak a' \, ,   \cr 
0  \quad \qquad \qquad \, \text{if} \ i = i' \, , \mathfrak a = \mathfrak a' \, ,   \cr 
|i-i'|  \qquad \quad \text{if} \ i \neq i'  \, . 
\end{cases}
$$ 
{\bf Notation:} Given a matrix $ A  \in {\cal M}_E^E $, when writing the matrix 
$ D_m^{1/2} A D_m^{1/2}  \in {\cal M}_E^E  $,   we understand that 
we apply the diagonal matrix $ D_m^{1/2} : E \to E $
  to the right/left of $ A $. 
  
\begin{definition}\label{goodmatrix}
{\bf ($N$-good/bad matrix)} The matrix $ A \in {\cal M}_E^E $, with $ E \subset \Z^b \times \mathfrak I $, 
$ {\rm diam}(E) \leq 4 N $,  is $ N $-good\index{N-good matrix} if $ A $ is invertible and
\be\label{Ngoodmat}
\forall s \in [s_0, s_1] \   , \ \ |D_m^{-1/2} A^{-1} D_m^{-1/2}|_s \leq N^{\tau'+\varsigma s} \, .
\ee
Otherwise $ A $ is $ N $-bad\index{N-bad matrix}.
\end{definition}

The above definition is different with respect to that of \cite{BBo10}: 
the matrix  $ A $ is $ N $-good according to 
Definition \ref{goodmatrix} if and only if  $ D_m^{1/2} A D_m^{1/2} $ is $ N $-good  
according to 
Definition \ref{goodmatrixA}.

\begin{definition}\label{regulars} {\bf (Regular/Singular sites)} 
The index $ k := (i,\mathfrak a) = (\ell, j, \mathfrak a) \in \Z^{b} \times  \mathfrak I $  (where $ \fracchi $ is defined in \eqref{def:indice a})
is  {\sc regular}\index{Regular site} for $ A $  if 
$$ 
|A_k^k| \geq \Theta \langle j \rangle^{-1} \, . 
$$ 
Otherwise $ k $ is\index{Singular site}  {\sc singular}.
\end{definition}

Also the above definition  is different with respect to that of \cite{BBo10}: 
the index 
$ k $ is regular for $ A $ according to 
Definition \ref{regulars} if and only if  $ k $ 
is regular  for $ D_m^{1/2} A D_m^{1/2} $ according to Definition \ref{regularsA}
with $ \Theta $ replaced by $ c(m) \Theta $
(because of the equivalence $ (|j|^2 + m)^{1/2} \sim_m \langle j \rangle $).

The constant $ \Theta := \Theta (V) $ will be chosen large enough depending on the potential $  V(x) $ 
in order to apply the multiscale proposition (as in \cite{BBo10}, \cite{BB12}). 

\begin{definition}\label{ANreg}
{\bf ($(A,N)$-good/bad site)}
For $ A \in \matr^E_E $, we say that $ k \in E \subset \Z^b \times \mathfrak I  $ is
\begin{itemize}
\item $(A,N)$-{\sc regular} if there is $ F \subset E$ such that
${\rm diam}(F) \leq 4N$,  ${\rm d}(k, E\backslash F) \geq N$ and
$ A_F^F $ is $N$-good. 
\item $(A,N)$-{\sc good}  if it is regular for $A$ or $(A,N)$-regular. 
Otherwise we say that $ k $ is $(A,N)$-{\sc bad}.
\end{itemize}
\end{definition}

Note that a site $ k $ is $(A,N)$-{\sc good} according to Definition \ref{ANreg} if and only if 
$ k $ is $(D_m^{1/2} A D_m^{1/2},N)$-{\sc good}   according to  \cite{BBo10}.  

\smallskip

Let us consider  the new larger scale 
\be\label{newscale}
N' = N^\chi 
\ee
with $ \chi > 1 $.  

\smallskip

For a matrix $ A \in \matr_E^E  $ we define 
$ {\rm Diag}(A) := ( \d_{kk'} A_k^{k'})_{k, k' \in E} $.

\begin{proposition} {\bf (Multiscale step)}\index{Multiscale step} \label{propinv}
Assume
\be\label{dtC}
\varsigma \in (0,1/2) \, ,  \ \tau' > 2 \tau + b + 1 \, , \ C_1 \geq 2 \, , 
\ee
and, setting $ \kappa := \tau' + b + s_0 $,
\begin{align}\label{chi1}
& \chi (\t' - 2 \t - b) >  3 (\kappa + (s_0+ b) C_1 ) \, , \\
& \chi \varsigma > C_1 \, \label{chil} \, , \\
& \label{s1} s_1 > 3 \kappa + \chi (\tau + b) + C_1 s_0 \, . 
\end{align}
For any given $ \Upsilon > 0 $, there exist $\Theta := \Theta (\Upsilon, s_1) > 0 $ large enough 
(appearing in  Definition \ref{regulars}),  and  $  N_0 (\Upsilon, \Theta ,  s_1) \in \N $  
 such that:
\\[1mm]
$ \forall N \geq N_0(\Upsilon, \Theta  , s_1) $, 
$ \forall E \subset \Z^b \times \mathfrak I  $ with 
${\rm diam}(E) \leq 4N'=4N^\chi $ (see \eqref{newscale}), if $ A \in \matr_E^E $ satisfies 
\begin{itemize}
\item 
{\bf (H1)} {\bf (Off-diagonal decay)} $ |A- {\rm Diag}(A)|_{+, s_1} \leq \Upsilon $ 
\item 
{\bf (H2)} {\bf ($ L^2 $-bound)}  $ \| D_m^{-1/2} A^{-1} D_m^{-1/2} \|_0 \leq (N')^{\tau}$
\item
{\bf (H3)} {\bf (Separation properties)}
There is a partition of the  $(A,N)$-bad sites $ B = \cup_{\alpha} \Omega_\alpha$ with
\be\label{sepabad}
{\rm diam}(\Omega_\alpha) \leq N^{C_1} \, , \quad {\rm d}(\Omega_\alpha , \Omega_\beta) \geq N^2 \ , \ \forall \alpha \neq \beta \, ,
\ee
\end{itemize}
then $ A $ is $ N' $-good. More precisely
\be\label{A-1alta}
\forall s \in [s_0,s_1]  \ , \  \  
|D_m^{-1/2}   A^{-1} D_m^{-1/2}|_s \leq \frac{1}{4} ({N'})^{\tau' } \big( ({N'})^{\varsigma s}+ | A- {\rm Diag}(A) |_{+,s} \big) \, ,
\ee
and, for all $ s \geq s_1 $, 
\be\label{multi:s}
|D_m^{-1/2}   A^{-1} D_m^{-1/2}|_s \leq C(s) ({N'})^{\tau' } \big( ({N'})^{\varsigma s}+ | A- {\rm Diag}(A) |_{+,s} \big) \, .
\ee
\end{proposition}

\begin{remark} The main difference with respect to the multiscale Proposition 4.1 in 
\cite{BBo10} is that, since 
the Definition \ref{regulars} of regular sites  is  weaker than that in 
\cite{BBo10}, 
we require the stronger assumption ($H1$)
concerning the off-diagonal decay  of $ A $ in $ | \ |_{+,s_1}$ norm defined in \eqref{def:n+s}, 
while  in   \cite{BBo10} we only require the off-diagonal decay of $ A $  in 
$ | \ |_{s_1} $ norm. Another difference is that we prove \eqref{A-1alta} (with the constant $ 1/ 4 $) for  
$ s \in [s_0,s_1] $, and not in a larger interval 
$  [s_0,S] $ for some $ S \geq s_1 $. For larger $ s \geq s_1 $
we prove \eqref{multi:s} with $ C(s) $. 
\end{remark}

\begin{pfn}{\sc of Proposition \ref{propinv}.} The multiscale step
Proposition \ref{propinv} follows by Proposition 4.1 in \cite{BBo10}, that we report
in the Appendix \ref{App:MU},  see Proposition \ref{propinvA}. 
Set 
\be \label{def:T}
{\cal T}  := A - {\rm Diag}(A) \, , 
\qquad |{\cal T}|_{+,s_1} \stackrel{(H1)} \leq \Upsilon \, , 
\ee
and consider the matrix
\be \label{def:B}
A_+ := D_m^{1/2} A D_m^{1/2} = {\rm Diag}(A_+) + T 
\ee
where
\be\label{def:BT}
{\rm Diag} (A_+) := D_m^{1/2}  {\rm Diag}(A) D_m^{1/2} \, , \quad 
T := D_m^{1/2}  {\cal T} D_m^{1/2} \, . 
\ee
We apply the multiscale Proposition 
\ref{propinvA}  to the matrix $ A_+ $. 
By (H2) the matrix $ A_+  $ is invertible and 
$$ 
\| A_+^{-1} \|_0 \stackrel{\eqref{def:B}} = 
\| D_m^{-1/2} A^{-1} D_m^{-1/2}\|_0  \stackrel{(H2)} \leq (N')^{\tau} \, .  
$$ 
Moreover, 
 the decay norm 
 $$ 
 | T |_{s_1}  \stackrel{\eqref{def:BT}}  
 =   | D_m^{1/2}  {\cal T} D_m^{1/2}  |_{s_1} 
 \stackrel{\eqref{def:n+s}}    =  | {\cal T} |_{+,s_1}  
 \stackrel{\eqref{def:T}} \leq \Upsilon \, . 
 $$ 
Finally the  
$(A_+,N)$-{\sc bad} sites according to Definition \ref{ANregA} coincide
with the $(A,N)$-{\sc bad} sites according to Definition \ref{ANreg} (with a $ \Theta $ replaced by  $ c(m) \Theta $). 
Hence by (H3) also the separation properties
required to  apply Proposition \ref{propinvA} hold, and we deduce that 
$$ 
\forall s \in [s_0, s_1] \, , \quad 
|A_+^{-1}|_s \leq \frac{1}{4} ({N'})^{\tau' } \big( ({N'})^{\varsigma s}+  |A_+- {\rm Diag}(A_+)|_s \big) \, , 
$$
that, recalling \eqref{def:B}, \eqref{def:T},  \eqref{def:BT}, 
implies \eqref{A-1alta}. The more general estimate \eqref{multi:s} follows by
\eqref{multi:s_APP}. 
\end{pfn}

\section{Separation properties of bad sites}\label{sec:sepabad}

The aim of this section is to verify the separation\index{Separation 
properties of bad sites} properties of the bad sites 
required in the multiscale step Proposition \ref{propinv}. 

Let  $ A := A(\e,\l,\teta ) $ be the infinite dimensional matrix defined in (\ref{Atranslated}). Given $ N \in \N $
and $ i = (\ell_0,j_0 ) $, recall that the submatrix  $ A_{N,i} $ is defined in (\ref{ANl0}). 

\begin{definition} \label{GBsiteRS} {\bf ($N$-regular/singular site)}
A site $ k := (i,\mathfrak a) \in \Z^b \times \mathfrak I $ is\index{N-regular site}\index{N-singular site}: 
\begin{itemize}
\item $ N$-{\sc regular } if $ A_{N,i} $ is $ N $-good (Definition \ref{goodmatrix}). 
\item $ N$-{\sc singular} if $ A_{N,i} $ is $ N $-bad (Definition \ref{goodmatrix}). 
\end{itemize}
\end{definition}
\noindent
We also define the $N$-good/bad sites of $ A $. 

\begin{definition} \label{GBsite} {\bf ($N$-good/bad site)} A site $ k := (i,\mathfrak a) \in \Z^b \times \mathfrak I $ is\index{N-good site}\index{N-bad site}: 
\begin{itemize}
\item  
$ N$-{\sc good} if 
\be\label{regNreg} 
k  \, {is \, regular} \, ({Def.} \, \ref{regulars})
\  {\rm or} \ 
{ \it all \,  the \, sites} \,  k' \,  { with} \,  {\rm d}(k',k) \leq N \, { are} \, N-{ regular} \, .
\ee
\item $ N $-{\sc bad} if
\be\label{regNreg-bad} 
k  \, {is \, singular} \, ({Def.} \, \ref{regulars})
\  {\rm and} \ 
 \exists  \, k'   \, { with} \,  {\rm d}(k',k) \leq N \, ,   \, k'  \, {\rm is } \, N{-singular}   \, .
\ee
\end{itemize}
\end{definition}

\begin{remark}\label{good2}
A site $ k $ which is $ N $-good  according to Definition \ref{GBsite}, is $(A_E^E,N)$-good according to Definition \ref{ANreg}, 
for any set $ E = E_0 \times \mathfrak I $ containing $ k $ where 
$ E_0  \subset \Z^b $ is a product of intervals of length $\geq N$.
\end{remark}

Let 
\be\label{tetabad}
B_{N}(j_0; \l )  := \Big\{ \teta \in \R \, : \,   A_{N, j_0}(\e, \l,\theta) \ {\rm is} \ N-bad  \Big\} \, . 
\ee

\begin{definition} {\bf ($ N$-good/bad parameters)} \label{def:freqgood}
A parameter $  \l \in \Lambda $ is $ N$-good\index{N-good parameter} for $ A $ if
\be\label{BNcomponents}
\forall \, j_0 \in \Z^d  \, , \quad 
B_{N}(j_0;  \l)  \subset \bigcup_{q = 1, \ldots, N^{\a-d-\es}} I_q \, , \quad \a := 3d+ 2 \es+4+ 3 \tau_0 \, , 
\ee
where $ I_q $ are  intervals with measure $ | I_q| \leq N^{-\t} $. 
Otherwise, we say that $  \l  $ is $ N$-bad\index{N-bad parameter}.
We define the set of $ N $-good parameters
\be\label{good}
{\cal G}_{N}  := \Big\{ \l \in \wtilde \Lambda 
\, : \,  \l  \ \  {\rm is \ } \ N-{\rm good \ for \ } A(\e, \lambda)  \Big\} \, . 
\ee
\end{definition}

The main result of this section is  Proposition \ref{prop:separation}
which enables to verify the assumption (H3) of Proposition \ref{propinv} 
for the submatrices $ A_{N',j_0}(\e, \l,\theta) $. 

\begin{proposition}\label{prop:separation}
{\bf (Separation properties of $ N $-bad sites)}
Let $ \t_1,\g_1,\t_2, \g_2  $ be fixed as in \eqref{def:tau1} and \eqref{ga2t2}, depending
on the parameters $\t_0, \g_0$ which appear in properties 
\eqref{diop}-\eqref{NRgamma0}, \eqref{2Mel+}-\eqref{2Mel rafforzate}.
Then there exist $ C_1 := C_1(d, \es, \t_0) \geq 2 $, $ \tau^\star := \tau^\star (d, \es, \t_0) $,  and  
$ \bar N := \bar N (\es,d, \g_0, \t_0, m,  \Theta)$   such that,  if $N \geq \bar N $ and 
\begin{itemize}
\item {\bf (i)} $ \l $ is $ N$-good for $ A $,
\item {\bf (ii)}  $ \tau \geq \tau^\star $, 
\item  {\bf (iii)} $\bar \om_\e$  satisfies \eqref{dioep}  and $ \om = (1 + \e^2 \l ) \bar \om_\e $ satisfies $ ({\bf NR})_{\gamma_2, \tau_2} $ 
(see Definition \ref{NRgamtau}),  
\end{itemize}
then, $  \forall \teta \in \R $,   
the $ N $-bad sites $ k =  (\ell,j, \mathfrak a) \in \Z^\es \times \Z^d \times {\mathfrak I } $  of $ A(\e, \l, \teta) $   
admit a partition $ \cup_\a \Omega_\a $ in disjoint clusters satisfying
\be\label{separ}
{\rm diam}(\Omega_\a) 
\leq N^{C_1(d,\es,\t_0)} \, ,   \quad {\rm d}(\Omega_\a, \Omega_\b) > N^2 \, , \ 
\forall \a \neq \b \, .
\ee
\end{proposition} 

We underline that  the estimates (\ref{separ}) are {\it uniform} in $ \teta $.

\begin{remark}\label{tau-large}
Hypothesis (ii) in Proposition \ref{prop:separation} just requires that the constant 
$ \tau $ is larger than some $ \tau^\star (d, \es, \t_0 )$.  This is important in the present 
autonomous setting for the choice of the constants in section
\ref{sec:constants}. 
On the contrary in  the corresponding propositions in the 
papers \cite{BB12}, \cite{BBo10}, the constant $ \tau $ was required to be large 
with the exponent $ \chi $ in \eqref{newscale}. 
\end{remark}

The rest of this section is devoted to the proof of Proposition \ref{prop:separation}. In some parts
of the proof, we may  point out the dependence of some constants on parameters such as $\t_2, \t_1$, which,
by \eqref{def:tau1} and \eqref{ga2t2}, amounts to a dependence on $\t_0$.

\begin{definition} \label{chain}
{\bf ($ \Gamma $-chain)} 
A sequence $ k_0, \ldots , k_L \in \Z^\es \times \Z^d \times {\mathfrak I}  $ of distinct integer vectors
satisfying 
$$ 
| k_{q+1} - k_q | \leq \Gamma \, , \quad \forall q = 0, \ldots, L - 1 \, ,
$$ 
 for some $ \Gamma  \geq 2 $,  is called a $ \Gamma $-chain\index{Chain} of length\index{Length of a $ \Gamma$-chain} $ L $.
\end{definition}

We want to prove an upper bound for the length of a $ \Gamma $-chain of 
$ \theta $-singular integer vectors, i.e.  satisfying  \eqref{sing-teta} below. 
In the next  lemma we obtain the upper bound \eqref{lunghez} 
under the assumption \eqref{cardinality}. Define the functions of signs
\be\label{def:signs}
\begin{aligned}
&  
\s_1, \s_2 : {\mathfrak I} \to \{ -1, + 1 \} \\ 
& \s_1(1) :=\s_1 (2):=\s_2 (1):=\s_2 (3) :=-1 \, , \\
& \s_1(3) :=\s_1(4) :=\s_2(2) :=\s_2(4) :=1 \, . 
\end{aligned}
\ee

\begin{lemma}\label{Bourgain} 
{\bf (Length of $ \Gamma$-chain of $ \theta $-singular sites)}
Assume that $ \om  $ satisfies the non-resonance condition $  {\bf (NR)}_{\g_2, \t_2} $ (Definition \ref{NRgamtau}). 
For $ \Gamma \geq \bar \Gamma (d, m,\Theta, \g_2, \t_2) $ large enough, 
consider a $ \Gamma $-chain  $( \ell_q, j_q, \mathfrak a_q)_{q=0, \ldots, L} 
\subset  \Z^\es \times \Z^d \times \fracchi $ of 
singular sites for the matrix $ A(\e, \lambda, \theta ) $, namely 
\be
\label{sing-teta}
 \forall q=0, \ldots, L, \quad 
\big|  \la j_q \ra_m \big(  \la j_q \ra_m + \s_1({\mathfrak a}_q) \mu + \s_2 ({\mathfrak a}_q) (\om \cdot \ell_q + \theta) \big) \big| < \Theta \, , 
\ee
with $ \s_1, \s_2 $  defined in  \eqref{def:signs},  
such that, $ \forall \tilde{\jmath} \in \Z^d $, the cardinality
\be\label{cardinality}
|\{ ( \ell_q,j_q, \mathfrak a_q)_{q = 0, \ldots, L} \, : \,  j_q = \tilde{\jmath} \}| \leq K \, .
\ee
Then there is $ C_2 := C_2(d,   \tau_2) >  0 $ such that its length is bounded by
\be\label{lunghez}
L \leq  ( \Gamma K)^{C_2} \, . 
\ee
Moreover, if $ \ell $ is fixed (i.e. $ \ell_q = \ell $, $\forall q = 0, \ldots, L $) the same result holds without assuming that
$ \om $ satisfies $  {\bf (NR)}_{\g_2, \t_2} $. 
\end{lemma}

The proof of Lemma  \ref{Bourgain}  is a  variant of Lemma 4.2 in \cite{BB12}. We split it in several steps. 

First  note that  it  is sufficient to 
bound the length of a $ \Gamma $-chain of singular sites when $ \theta = 0 $. 
Indeed, suppose first that $ \theta = \o \cdot \bar \ell $ for some $ \bar \ell \in \Z^\es $.
For a $ \Gamma $-chain of $ \theta $-singular sites $(\ell_q, j_q, {\mathfrak a}_q)_{q=0, \ldots, L} $, see \eqref{sing-teta},
the translated $ \Gamma $-chain
$(\ell_q + \bar \ell, j_q, {\mathfrak a}_q)_{q=0, \ldots, L} $, 
is formed by $ 0 $-singular sites, namely
$$
\big|  \la j_q \ra_m \big(  \la j_q \ra_m + \s_1({\mathfrak a}_q) \mu + \s_2 ({\mathfrak a}_q) (\om \cdot (\ell_q + \bar \ell)) \big) \big| < \Theta \, . 
$$
For any $ \teta \in \R $, we consider an approximating sequence 
$ \o \cdot {\bar \ell}_n \to \teta $, $ {\bar \ell}_n \in \Z^\es  $.  Indeed, by 
 Assumption \eqref{dioep}, $\o = (1 + \e^2 \l) \bar \om_\e $ is not colinear to an integer vector, and therefore the set
 $\{ \o \cdot \ell \ , \  \ell \in \Z^\es \}$ is dense in $\R$. 
A $ \Gamma $-chain of $ \theta $-singular sites (see \eqref{sing-teta}), is, for 
$ n $ large enough, also a $ \Gamma $-chain of $ \om \cdot {\bar \ell}_n $-singular sites. Then
 we bound its length arguing as in the above case.

We first prove Lemma \ref{Bourgain} in a particular case.  
\begin{lemma}\label{lem:step1}
Assume that $ \om  $ satisfies $  {\bf (NR)}_{\g_2, \t_2} $ (Definition \ref{NRgamtau}).
Let 
$ (\ell_q, j_q, {\mathfrak a}_q)_{q= 0, \ldots, L} $ be a $ \Gamma $-chain  
of 
integer vectors of  $ \Z^\es \times \Z^d \times {\mathfrak I} $ satisfying, 
$ \forall q = 0 , \ldots, L $,
\begin{align}\label{singusites:set1}
& \big| \sqrt{| j_q |^2+m}  + \s_{2}({\mathfrak a}_q) \, \om \cdot \ell_q  \big| \leq \frac{\Theta}{ \langle j_q  \rangle} \, , \qquad  \qquad
\quad  \
{ in \ case} \ (i) \ { of  \ Def. } \, \ref{definition:Xr} \, , \\
& \label{singusites}
\big| \sqrt{| j_q |^2+m} +  \s_{1}({\mathfrak a}_q) \mu + \s_{2}({\mathfrak a}_q) \, \om \cdot \ell_q  \big| \leq \frac{\Theta}{ \langle j_q  \rangle} \, , 
 \quad { in \ case} \ (ii) \ { of \ Def.}   \, \ref{definition:Xr} \, ,
\end{align}
where  $ \s_1, \s_2 $  are defined in  \eqref{def:signs}.  
Suppose, in case \eqref{singusites}, 
that the product of the signs $ \s_{1}({\mathfrak a}_q) \s_{2}({\mathfrak a}_q) $ 
is the same for any $ q \in [\![0,L ]\!] $. 
Then, for some constant $ C_1 := C_1 (d, \tau_2) $ and 
$C := C (m,\Theta,\g_2, d, \t_2)$, its length  $ L $ is bounded by 
\be\label{Lstep1}
L \leq C (\Gamma K)^{C_1} 
\ee
where $ K $ is defined in \eqref{cardinality}.

Moreover, if $ \ell $ is fixed (i.e. $ \ell_q = \ell $, $\forall q $), the lemma holds without assuming that $ \om $ satisfies $  {\bf (NR)}_{\g_2, \t_2} $. 
\end{lemma}

\begin{pf}
We make the proof when \eqref{singusites} holds, since \eqref{singusites:set1} is a particular case 
of \eqref{singusites} setting $ \mu = 0$ (note that, in the case \eqref{singusites:set1},  
the conclusion of the lemma follows without conditions on the signs 
$ \s_{2}({\mathfrak a}_q) $, see Remark \ref{rem:signs}).

We introduce  
the quadratic  form 
$ Q : \R \times \R^d \to \R $ defined by
\be\label{Qxy}
Q(x,y) := - x^2 + | y |^2  
\ee
and the  associated bilinear symmetric form $ \varPhi : (\R \times \R^d)^2 \to \R $ defined by 
\be\label{def:varphi}
\varPhi \big((x,y), (x',y') \big) := - x x' + y \cdot y' \, .  
\ee
Note that $ \varPhi $ is the sum of the bilinear forms
\be
\begin{aligned}
& \qquad \qquad \qquad \qquad \varPhi = - \varPhi_1 + \varPhi_2   \label{-v12} 
\\
& 
\varPhi_1 \big((x,y), (x',y') \big) :=  x x' \, , 
\quad \varPhi_2 \big((x,y), (x',y') \big) := y \cdot y' \, . 
\end{aligned}
\ee
\begin{lemma}\label{Phi-bi}
For all $  q, q_0 \in [\![0,L ]\!] $, 
\be\label{seconda} 
\big| \varPhi \big(( x_{q_0} , j_{q_0} ), ( x_q - x_{q_0}, j_q - j_{q_0}) \big) \big| \leq C\Gamma^2 |q - q_0 |^2 + C(\Theta) \, ,
\ee
where  $x_q := \o \cdot \ell_q $.
\end{lemma}

\begin{pf}
Set for brevity $ \s_{1,q} := \s_{1}({\mathfrak a}_q) $ and 
$ \s_{2,q} := \s_{2}({\mathfrak a}_q) $.
First note that  by \eqref{singusites} we have 
$$
\big| |j_q|^2+m - ( \s_{1,q} \mu + \s_{2,q} \om \cdot \ell_q  )^2 \big| 
\leq C \Theta \, . 
$$
Therefore
$$
\big| - (\om \cdot \ell_q)^2  + |j_q|^2  - 2  \s_{1,q} \, \s_{2,q} \, \mu \, \om \cdot \ell_q
 \big| \leq C' \, , \quad C' := C \Theta + m + \mu^2 \, , 
 $$
and so, recalling \eqref{Qxy}, for all $ q = 0 , \ldots, L $,  
\be\label{sing-h1}
| Q ( x_q, j_q ) - 2  \s_{1,q}  \s_{2,q} \, \mu \, x_q | \leq C' \qquad {\rm where} \qquad x_q := \o \cdot \ell_q \, . 
\ee
By the hypothesis of Lemma \ref{lem:step1}, 
$$ 
\s_{1,q} \s_{2,q}  =  \s_{1,q_0} \s_{2,q_0} \, , \quad \forall q, q_0 \in [\![0,L ]\!] \, , 
$$ 
and, by bilinearity, we get  
\begin{align}\label{seco-re-s2}
Q(x_q, j_q) - 2  \s_{1,q} \s_{2,q} \mu x_{q} & =  Q ( x_{q_0}, j_{q_0})  \nonumber \\
& \ + 2 \varPhi \big( ( x_{q_0}, j_{q_0}), (x_q - x_{q_0}, j_q - j_{q_0})\big)  +
Q (x_q - x_{q_0}, j_q - j_{q_0})  \nonumber \\
& \  - 2 \s_{1,q_0} \s_{2,q_0} \mu x_{q_0} - 2  \s_{1,q} \s_{2,q} \mu (x_{q}- x_{q_0}) \, . 
\end{align}
Recalling \eqref{Qxy}  and the Definition \ref{chain} of $ \Gamma $-chain we have, $ \forall q, q_0 \in [\![0,L ]\!] $, 
\be\label{upper-bound-Q}
\big| Q (x_q - x_{q_0}, j_q - j_{q_0}) - 2  \s_{1,q} \s_{2,q} \mu (x_{q}- x_{q_0}) \big| \leq C \Gamma^2 | q - q_0 |^2 \, . 
\ee
Hence \eqref{seco-re-s2}, \eqref{sing-h1}, \eqref{upper-bound-Q} imply 
\eqref{seconda}.
\end{pf}

\begin{remark}\label{rem:signs}
In the  case \eqref{singusites:set1},  
the conclusion of  Lemma \ref{Phi-bi} follows without conditions on the signs 
$ \s_{2}({\mathfrak a}_q) $. This is why Lemma \ref{lem:step1} holds
without condition on the signs 
$ \s_{2}({\mathfrak a}_q) $. 
\end{remark}

\noindent
{\bf Proof of Lemma \ref{lem:step1} continued.}  
We introduce the subspace of $ \R^{d+1} $
\be\label{G+}
\begin{aligned}
G & := {\rm Span}_\R \Big\{  ( x_q - x_{q'}, j_q - j_{q'}) \, : \, 0 \leq q, q' \leq L  \Big\} \\
& = 
{\rm Span}_\R \Big\{ ( x_q - x_{q_0}, j_q - j_{q_0})  \, : \, 0 \leq q \leq L  \Big\}
\end{aligned}
\ee
and we call $ g \leq d + 1 $ the dimension of $ G $. Introducing a small parameter
$\delta >0$, to be specified later (see \eqref{fixd}), we distinguish two cases. \\[1mm]
{\bf Case I}. 
$ \forall q_0 \in [\![0,L ]\!] $, 
\be\label{Case1} 
{\rm Span}_\R \big\{ ( x_q - x_{q_0}, j_q - j_{q_0}) \, : \,  | q - q_0 | \leq L^\d \, ,  \, 
q \in [\![0,L ]\!] \, \big\} = G \, . 
\ee
We select a basis of  $ G  \subset \R^{d+1}  $  
from  $ ( x_q - x_{q_0}, j_q - j_{q_0})  $ with $ | q - q_0 | \leq L^\d $, say 
\be\label{fs+}
f_s := ( x_{q_s} - x_{q_0}, j_{q_s} - j_{q_0}) = ( \om \cdot \Delta_s \ell , \Delta_s j ) \, , \ \ s = 1, \ldots, g \, ,
\ee
where $ (\Delta_s \ell, \Delta_s j) := ( \ell_{q_s} - \ell_{q_0}, j_{q_s} - j_{q_0} ) $
satisfies, by the Definition \ref{chain} of $ \Gamma $-chain, 
\be\label{lanuova}
\quad  |(\Delta_s \ell, \Delta_s j)|  \leq C 
\Gamma | q_s - q_0 | \leq C \Gamma L^\d \, . 
\ee
Hence
\be\label{boundfi}
| f_s | \leq C \, \Gamma L^\d \,  , \qquad \forall s =1, \ldots, g \, .
\ee
Then, in order to derive from \eqref{seconda} a bound on $(x_{q_0}, j_{q_0})$ or its projection onto $G$, 
we need a nondegeneracy property for $Q_{|G}$. The following lemma states it. 
\begin{lemma}\label{Omegainv}
Assume that $ \omega $ satisfies\index{Quadratic Diophantine condition}  $ {\bf (NR)}_{\g_2, \t_2} $. 
Then the matrix 
\be\label{Omega}
\Omega := ( \Omega_s^{s'} )_{s,s'=1}^g \, , \quad \Omega_s^{s'} :=  \varPhi ( f_{s'}, f_s)  \, ,
\ee
is invertible and
\be\label{Cramer}
|( \Omega^{-1})_s^{s'}| \leq C (\Gamma L^\d)^{C_3(d,\tau_2)} \, , \quad \forall s, s' = 1, \ldots , g \, ,
\ee
where the multiplicative constant $ C $ depends on $ \g_2$. 
\end{lemma}

\begin{pf}
According to the splitting \eqref{-v12} we write  $ \Omega $  as
\be\label{-1-2}
\Omega :=  \big(  - \varPhi_1 ( f_{s'}, f_s) + \varPhi_2 ( f_{s'}, f_s) \big)_{s,s'=1, \ldots, g} = - S + R 
\ee
where, by \eqref{fs+}, 
\be\label{RSss'}
\begin{aligned}
& S_{s}^{s'} := \varPhi_1 ( f_{s'}, f_s) =
 (\o \cdot \Delta_{s'} \ell ) (\o \cdot \Delta_s \ell)   \, , \\ 
&  R_{s}^{s'} := \varPhi_2 (f_{s'}, f_s) =  \Delta_{s'} j \cdot \Delta_s j  \, .
\end{aligned}
\ee
The matrix $ R = (R_1, \ldots, R_{g}) $ has integer entries (the $ R_i \in \Z^{g} $ denote the columns).
The matrix $ S := ( S_1, \ldots, S_g) $ has rank $ 1 $ since all its columns $  S_s \in \R^{g}  $ 
are colinear:
\be  \label{defSs}
S_s = (\om \cdot  \Delta_s \ell) ( \om \cdot \Delta_1 \ell , \ldots,   \om \cdot \Delta_g  \ell)^\top \, , \quad s = 1, \ldots g \, . 
\ee
We develop the determinant 
\begin{eqnarray}
P(\om) & := & {\rm det} \, \Omega  \stackrel{\eqref{-1-2}}  = {\rm det} (- S + R) \nonumber \\
&  = & {\rm det} (R ) - {\rm det} ( S_1, R_2, \ldots, R_{g}) - \ldots -
{\rm det} (R_1, \ldots, R_{g-1}, S_g)  \nonumber \\
&=&    {\rm det} (R ) - \sum_{1\leq s \leq g}  (-1)^{g+s} S_s\cdot  (R_1 \wedge \ldots \wedge R_{s-1} 
\wedge R_{s+1} \wedge \ldots \wedge R_g ) \label{polin}  
\end{eqnarray}
using that the determinant of matrices with $ 2 $ columns  $ S_{i} $, $ S_j $, $ i \neq j $,
is zero. By \eqref{defSs}, the expression  in \eqref{polin} is a polynomial in $ \om $ of degree $ 2 $ of the form \eqref{NRom0}
with coefficients 
\be\label{coeffb}
|(n,p)|  \stackrel{\eqref{RSss'}, \eqref{lanuova}} 
\leq C (\Gamma L^\d)^{2s} \leq  C (\Gamma L^\d)^{2(d+1)}   \, .
\ee 
If $ P \neq  0 $ then the non-resonance condition $ {\bf (NR)}_{\g_2, \t_2}  $ implies 
\be\label{determinant}
|{\rm det} \, \Omega | = |P(\om)| \stackrel{ \eqref{NRom} } \geq \frac{\g_2}{\langle p \rangle^{\t_2}} 
\stackrel{\eqref{coeffb} } \geq \frac{\gamma_2}{ C (\Gamma L^\d)^{2\t_2 (d+1)}}  \, . 
\ee
In order to conclude the proof of the lemma, we have to show that 
$ P (\om ) $ is not identically zero in $ \omega  $. 
We have that 
$$
P( \ii \omega)  = {\rm det} \big( \varPhi_1 ( f_{s'}, f_s) + \varPhi_2 ( f_{s'}, f_s)\big)_{s,s'=1,\ldots g} 
=  {\rm det} ( f_{s'} \cdot f_s)_{s,s'=1,\ldots g}   > 0
$$
because $ (f_s)_{1\leq s \leq g}  $ is a basis of $ G $. Thus  $ P $ is not the zero polynomial. 

By \eqref{determinant}, the Cramer rule, and \eqref{boundfi} we deduce \eqref{Cramer}.
\end{pf}

\begin{remark}
As recently proved in \cite{BMa} 
the same result holds also assuming just that $ \om $ is Diophantine, 
instead of the quadratic non-resonance condition $ {\bf (NR)}_{\g_2, \t_2}  $. 
\end{remark}

\noindent
{\bf Proof of Lemma \ref{lem:step1} continued.} 
We introduce 
$$
G^{\bot \varPhi} := \Big\{  z \in \R^{d+1}  \ : \  \varPhi (z, f) = 0 \, , \  \forall f \in G  \Big\} \, . 
$$
Since $ \Omega $ is invertible (Lemma \ref{Omegainv}), 
$ \varPhi_{|G}$ is nondegenerate, hence 
$$
\R^{d+1} = G \oplus G^{\bot \varPhi} 
$$
and we denote by $ P_G : \R^{d+1} \to G $ the corresponding projector onto $ G $.

\smallskip

We are going to estimate 
\be\label{PGf}
P_G ( x_{q_0}, j_{q_0} ) = \sum_{s' =1}^{g} a_{s'} f_{s'} \, . 
\ee
For all $ s =1, \ldots, g $, and since $ f_s \in G $,  we have
$$
\begin{aligned}
 \varPhi \big( ( x_{q_0}, j_{q_0} ), f_s \big) = 
\varPhi \big( P_G ( x_{q_0}, j_{q_0} ), f_s \big) 
& \stackrel{\eqref{PGf}} 
 =
 \sum_{s' =1}^{g} a_{s'} \varPhi ( f_{s'}, f_s)
 \end{aligned} 
$$
that we write as the linear system  
\be\label{def:ab}
\Omega a  = b \, , \qquad 
a := \left(\begin{array}{c} a_{1} \\ \vdots  \\ a_{g} \end{array}\right) \, , \quad 
b := \left(\begin{array}{c}  \varPhi \Big( ( x_{q_0}, j_{q_0} ), f_1 \Big)
\\ \vdots  \\  \varPhi \Big( ( x_{q_0}, j_{q_0} ), f_g \Big) \end{array}\right)
\ee
with $ \Omega $  defined in \eqref{Omega}.
\begin{lemma}\label{lem:12}
For all $ q_0 \in [\![0,L ]\!] $ we have 
\be\label{PGs}
| P_G ( x_{q_0}, j_{q_0} ) | \leq C( \Gamma L^\d )^{C_4(d,\tau_2)} \, ,
\ee
where $C$ depends on $\g_2, \Theta$.
\end{lemma}

\begin{pf}
We have
$$
|b| \stackrel{\eqref{def:ab}}{\lesssim}  \max_{1\leq s \leq g}  \Big|  \varPhi \Big( ( x_{q_0}, j_{q_0} ), f_s \Big)\Big|
\stackrel{\eqref{fs+}, \eqref{seconda} }{\leq}   C(\Theta)  \Gamma^2  \max_{1\leq s \leq g} |q_s-q_0|^2   \leq  C(\Theta) (\Gamma L^\d)^2 \, , 
$$
recalling that, by  \eqref{Case1}, the indices $(q_s)_{1\leq s \leq g}$ were
selected such that $|q_s-q_0| \leq L^\d$. 
Hence, by \eqref{def:ab} and \eqref{Cramer},   
\be\label{stima:a}
|a|  = | \Omega^{-1} b | \leq C(\g_2, \Theta) (\Gamma L^\d)^{C}  
\ee
for some constant $C := C(d,\t_2)$. 
We deduce \eqref{PGs} by   \eqref{PGf}, \eqref{stima:a} and \eqref{boundfi}.
\end{pf}

We now complete the proof of Lemma \ref{lem:step1} when case I holds.
As a consequence of Lemma \ref{lem:12}, for all $ q_1, q_2 \in [\![0,L ]\!] $, we have 
$$
| ( x_{q_1}, j_{q_1} ) -  ( x_{q_2}, j_{q_2} ) | =
\big| P_G \big( ( x_{q_1}, j_{q_1} ) -  ( x_{q_2}, j_{q_2} )\big) \big| 
\leq C ( \Gamma L^\d )^{C_4(d,\tau_2)} 
$$
where $ C $ depends on $ \g_2, \Theta $. 
Therefore, for all $ q_1, q_2 \in [\![0,L ]\!] $, we have 
$ | j_{q_1} - j_{q_2} | \leq C (\Gamma L^\d)^{C_4(d,\tau_2)} $, 
and so
$$
{\rm diam}  \{ j_{q} \ ; \ 0\leq q \leq L \}  \leq  C (\Gamma L^\d)^{C_4(d,\tau_2)} \, . 
$$
Since all the $ j_q $ are in  $ \Z^d $, 
their number (counted without multiplicity) 
does not exceed $ C (\Gamma L^\d)^{C_4(d,\tau_2) d} $,
for some other constant $ C $ which depends on $ \g_2, \tau_2, d $. 
Thus we have obtained the bound
$$
\sharp \{ j_{q} \ : \ 0\leq q \leq L \} \leq C (\Gamma L^\d)^{C_4(d,\tau_2) d}  \, . 
$$
By assumption \eqref{cardinality}, for each $q_0 \in [\![0,L ]\!]$, the number of $q\in [\![0,L ]\!]$ such that $ j_q = j_{q_0} $  is  at most $ K $, 
and so
$$
L \leq  C (\Gamma L^\d)^{C_5(d,\tau_2)}  K  \, , \qquad C_5(d,\tau_2) := C_4(d,\tau_2) d \, . 
$$
Choosing $ \d > 0 $ such that 
\be\label{fixd}
\d C_5(d,\tau_2)  < 1 \slash 2 \, , 
\ee
we get 
$ L \leq C (\Gamma^{C_5(d,\tau_2) } K )^2 $, for some 
multiplicative constant $C$ that may depend on $\g_2, \Theta, d$. This proves  \eqref{Lstep1}. 
\\[1mm] 
{\bf Case II.} There is $ q_0 \in [\![0,L ]\!] $ such that  
$$
{\rm dim} \, {\rm Span}_\R \big\{  ( x_q - x_{q_0},  j_q - j_{q_0} ) \, : \, 
| q - q_0 | \leq L^\d \, ,  \ q \in [\![0,L ]\!] \, 	\big\}  \leq g -1 \, ,
$$
namely all the vectors $ (x_q, j_q) $ stay in a affine subspace of dimension 
less than $ g - 1 $. 
Then we repeat on the sub-chain $ (\ell_q, j_q)   $, $ | q - q_0 | \leq L^\d $, the argument of case I, 
to obtain a bound for $ L^\d $ (and hence for $L$).

Applying  at most $ (d+1) $-times the above procedure, we obtain a bound for $L$ of the form 
$ L \leq  C( \Gamma K)^{C(d,\tau_2)} $. This concludes the proof of 
\eqref{Lstep1} of Lemma \ref{lem:step1}.

\smallskip

To prove the last statement of Lemma \ref{lem:step1}, note that, 
if $ \ell $ is fixed, then \eqref{seconda} reduces just to $ | j_q  \cdot ( j_q - j_{q_0}) | \leq C\Gamma^2 |q - q_0 |^2$ 
and the conclusion of the lemma follows as above, 
see also  Lemma 5.2 in \cite{BBo10} 
(actually it is the same argument for NLS in  \cite{B3}). 
\end{pf}

\smallskip

\noindent
{\bf Proof of Lemma \ref{Bourgain}. }
Consider a general  $ \Gamma $-chain 
$ (k_q)_{q= 0, \ldots, L} = (\ell_q, j_q, {\mathfrak a}_q)_{q= 0, \ldots, L} $ of singular sites satisfying
\eqref{singusites}. 
To fix ideas assume $ \s_{1} ( {\mathfrak a}_0) \s_{2} ({\mathfrak a}_0) >  0 $. 
Then the integer vectors $ k_q $ along the chain with the same sign 
$ \s_{1}({\mathfrak a}_q) \s_{2}({\mathfrak a}_q) > 0 $,  say $ k_{q_m} $,  
satisfy 
$$ 
q_{m+1} - q_m  \leq C (\Gamma K)^{C_1} \, , 
$$ 
by Lemma \ref{lem:step1} applied to each subchain of consecutive 
indices with $ \s_{1}({\mathfrak a}_q) \s_{2}({\mathfrak a}_q) < 0 $. 
Hence we deduce that all such  $ k_{q_m} $ form a $ \Gamma' := C \Gamma (\Gamma K)^{C_1} $-chain and 
all of them have the same sign  $ \s_{1} ({\mathfrak a}_{q_m}) \s_{2}
({\mathfrak a}_{q_m}) > 0 $. It follows, again 
by Lemma \ref{lem:step1}, that their length is bounded by  
$$
C (\Gamma' K )^{C_1}  =C \big( C \Gamma (\Gamma K)^{C_1} K \big)^{C_1} = C' (\Gamma K)^{C_1(C_1+1)}   \,  . 
$$
Hence the length of the original $ \Gamma $-chain $ (k_q)_{q= 0, \ldots, L} = 
(\ell_q, j_q, {\mathfrak a}_q)_{q= 0, \ldots, L} $ satisfies   
$$ 
L \leq C(\Gamma K)^{C_1} C' (\Gamma K)^{C_1(C_1+1)} \leq  (\Gamma K)^{C_2}  
$$ 
where $ C_2 := C_2 (d, \tau_2) =C_1(C_1+2)+1$, and provided $\Gamma$ is large enough,
depending on $m,\Theta, \g_2,\t_2, d $. This  proves \eqref{lunghez}.

Finally, the  last statement of Lemma \ref{Bourgain} follows as well
by the last statement of Lemma \ref{lem:step1}. 
\rule{2mm}{2mm}

\medskip

We fix
\be\label{tau-star}
\tau^\star := \tau_1 ( 2 + C_2 (3 + \alpha ) )  +1 
\ee
where $ \t_1 $ is the Diophantine exponent of $ \bar \om_\e $ in \eqref{dioep}, $ C_2 $ is the constant defined in Lemma
\ref{Bourgain}, 
and $ \a $ is defined in \eqref{BNcomponents}. 

The next lemma proves an upper bound for the length of a chain of $ N $-bad-sites.  

\begin{lemma}\label{Ntime}
{\bf (Length of $ \Gamma $-chain of $ N $-bad sites)} \index{Length of a $ \Gamma$-chain}
Assume (i)-(iii) of Proposition \ref{prop:separation} with $ \tau^\star 
$ defined in \eqref{tau-star}. Then, for $N$ large enough (depending on $m,\Theta,
\g_2, \t_2, \g_1, \t_1$), 
any $ N^2 $-chain $ (\ell_q, j_q, \mathfrak a_q)_{q= 0, \ldots, L} $  of $ N $-bad sites
 of $ A(\e, \l, \teta) $  (see Definition \ref{GBsite})
has length 
\be\label{finale}
L \leq  N^{(3+\alpha)C_2}  
\ee
where $ C_2 $ is  defined in Lemma \ref{Bourgain} and $ \a $  in \eqref{BNcomponents}. 
\end{lemma}

\begin{pf}
Arguing  by contradiction we assume that
\be\label{contra}
L > l := N^{(3+\alpha)C_2}  \, .
\ee
We consider the subchain $ (\ell_q, j_q, \mathfrak a_q)_{q= 0, \ldots, l} $ of $ N $-bad sites, which, 
recalling
\eqref{regNreg-bad}, are singular sites.
Then, for $N$ large enough, the assumption \eqref{cardinality}
of Lemma \ref{Bourgain} with $ \Gamma = N^2 $, and $ L $ replaced by $ l $,  
can not hold with $ K < N^{1+ \alpha} $, 
otherwise \eqref{lunghez} would imply $ l < (N^2 N^{1+\a} )^{C_2} = N^{(3+\alpha)C_2} $, contradicting
\eqref{contra}.  As a consequence there exists $ \tilde{\jmath} \in \Z^d $, 
and  distinct  indices  $ q_i  \in Ê[\![ 0,l]\!] $, 
$ i = 0, \ldots, M:=[N^{\a+1}/4]  $, such that 
\be
0 \leq  q_i  \leq l  
\qquad
{\rm and}
\qquad
j_{q_i}  = \tilde{\jmath} \, ,  \ \forall i = 0, \ldots, M:=\Big[\frac{N^{\a+1}}4 \Big]  \, .
\ee
Since the sites $ (\ell_{q_i}, \tilde{\jmath}, \mathfrak a_{q_i})_{i=0, \ldots, M} $ belong to 
a $ N^2 $-chain of length $ l $,  
the diameter of the set $ E := \{ \ell_{q_i},  \mathfrak a_{q_i}\}_{ i = 0, \ldots, M} \subset \Z^\es \times \mathfrak I $ satisfies 
$$
{\rm diam} (E) \leq N^2 l  \, . 
$$
Moreover, each site $ (\ell_{q_i},  \tilde{\jmath}, \mathfrak a_{q_i} ) $ is $ N $-bad, and therefore,
recalling \eqref{regNreg-bad}  it is in a 
$ N $-neighborhood of some $ N $-singular site. 
Let $ \mathfrak K $ be the number of $ N $-singular sites 
$ (\ell,  j, \mathfrak
 a ) $ such that $ | j - \tilde{\jmath} | \leq N $ and $ {\mathtt d}(E, (\ell, \mathfrak a)) \leq N $.
Then the cardinality $ | E | \leq C N^\es \mathfrak K $. 
Moreover there is $ \tilde{\jmath}' \in \Z^d $  with 
$ | \tilde{\jmath}' - \tilde{\jmath} | \leq N  $ such that 
there are at least $ \mathfrak K / (CN^{d})$ $ N $-singular sites $ (\ell,  \tilde{\jmath}', \mathfrak a ) $ 
with distance  $ {\mathtt d}(E, (\ell, \mathfrak a)) \leq N  $. 
Let 
$$ 
E'  := \Big\{  (\ell', \mathfrak a') \in   \Z^\es \times \mathfrak I \ : \  (\ell',  \tilde{\jmath}', \mathfrak a' )  \
{\rm is} \ N {\rm-singular} \ {\rm and} \    {\mathtt d}(E, (\ell', \mathfrak a')) \leq N 
  \Big\} \, . 
 $$
By what preceeds the cardinality
\be\label{cardE'}
\begin{aligned}
& |E'| \geq \frac{|E|}{C N^{d+\es}} \geq \frac{N^{\alpha + 1-d-\es}}{C} > N^{\alpha- d - \es} \, , \\
&  {\rm and} \quad {\rm diam} (E') \leq N^2 l + 2 N  \, . 
\end{aligned}
 \ee
Recalling Definition \ref{GBsiteRS} and the covariance property \eqref{shifted}, 
$$ 
\forall (\ell' , \mathfrak a') \in E' \, , \quad 
A_{N, \ell' , \tilde{\jmath}'} (\e, \lambda) = A_{N,  \tilde{\jmath}'} (\e, \lambda, \omega \cdot \ell'  ) \ \ 
{\rm is} \ \  N {\rm -bad} \, ,  
$$
and, recalling \eqref{tetabad}, 
\be\label{pingeon}
\forall (\ell' , \mathfrak a') \in E' \, , \quad
 \omega \cdot \ell' \in B_N ( \tilde{\jmath}'; \lambda ) \, . 
 \ee
 Since $ \lambda $ in $ N $-good (Definition \ref{def:freqgood}), 
\eqref{BNcomponents} holds and therefore 
\be\label{pengeon2}
 B_N (\tilde{\jmath}'; \lambda ) \subset  \bigcup_{q = 0, \ldots, N^{\a-d-\es}} \!\!\!\!\!\! I_q \quad
{\rm where} \ I_q  \  {\rm are  \ intervals \ with \ measure} \  | I_q| \leq N^{-\t} \, .
\ee
By \eqref{pingeon}, \eqref{pengeon2} and
since, by \eqref{cardE'},  the 
cardinality $ |E'| \geq   N^{\alpha- d - \es} $, 
 there are two distinct integer vectors $ \ell_1 $, $ \ell_2 \in E' $ such that
$ \om \cdot \ell_1 \, , \om \cdot \ell_2  $ belong to the same interval $ I_q $. 
Therefore
\be\label{up}
|\om \cdot ( \ell_1 - \ell_2)| = | \om \cdot \ell_1  -  \om \cdot \ell_2  | \leq |I_q| \leq N^{- \tau} \, .
\ee
Moreover, since by \eqref{dioep} the frequency vectors $ \om = (1+ \e^2 \l ) {\bar \om}_\e $, $ \forall \l \in \Lambda $, are Diophantine\index{Diophantine vector}, namely
$$
|\om \cdot \ell | \geq \frac{\g_1}{2 |\ell |^{\tau_1}} \, , \quad \forall \ell  \in \Z^{\es} \setminus \{ 0 \} \, , 
$$
we also deduce 
\begin{align} \label{lo}
|\om \cdot (\ell_1  - \ell_2)| \geq \frac{\g_1}{| \ell_1- \ell_2|^{\tau_1}}  
& \geq  \frac{\g_1}{ ({\rm diam}(E') )^{\tau_1}} \nonumber \\
& \stackrel{\eqref{cardE'}} 
\geq \frac{\g_1}{ (N^2 l + N)^{\tau_1}} \nonumber \\
& \stackrel{\eqref{contra}}{\geq}  \frac{\g_1}{ (2 N^{2 + (3+\alpha)C_2} )^{\tau_1}}  \, . 
\end{align}
The conditions \eqref{up}-\eqref{lo} contradict, for $N$ large enough,  the assumption that 
$ \tau \geq \tau^\star $ where  $ \tau^\star  $ is defined in \eqref{tau-star}. 
\end{pf}

\noindent
{\sc Proof of Proposition \ref{prop:separation} completed.} 
We introduce the following equivalence relation
in the set 
$$
{\cal S}_N := \Big\{ k = (\ell ,j, \mathfrak a) \in \Z^\es \times \Z^d \times {\mathfrak I} \,  :   \, k  \ {\rm is} \, 
N\text{-}{\rm \, bad \, } \ 
 {\rm for } \, A(\e, \l, \teta) 
\Big\}.
$$

\begin{definition}\label{equivalence} 
We say that $ x \equiv y $ if  there is a $N^2$-chain $ \{ k_q \}_{q = 0, \ldots, L} $
 in ${\cal S}_N$  
connecting $ x $ to $ y $, namely $ k_0 = x $, $ k_L = y $.
\end{definition}
This  equivalence relation induces a partition of the $ N$-bad sites
of $ A(\e, \l,\teta ) $, in disjoint equivalent  classes
$ \cup_\a \Omega_\a $,  satisfying, by Lemma \ref{Ntime}, 
\be\label{sepaquasi}
{\rm d}(\Omega_\a, \Omega_\b) >  N^2 \, , \quad {\rm diam}(\Omega_\a) 
\leq  N^2 N^{(3+\alpha)C_2}  = N^{C_1}
\ee
with $ C_1 := C_1(d, \es , \t_0) := 2 + (3+\alpha)C_2 $.  This proves \eqref{separ}.

\section{Definition of the sets $ \Lambda (\e; \cc, X_{r, \mu} ) $}\label{sec:constants}

In order to define the sets 
$ \Lambda (\e; \cc, X_{r, \mu} ) \subset \wtilde \Lambda  $ appearing in  the statement of Proposition \ref{propmultiscale}, 
we first fix the values of some constants. 
\begin{enumerate}
\item 
{\bf Choice of $ \tau $.} First we fix $ \tau  $ 
satisfying   
\be\label{def:tau}
\tau > \max \big\{ \tau^\star, 
9 d + 8 \es + 5 s_0 + 5  \big\}  \, ,
\ee
where the constant $ \tau^\star $ is defined in \eqref{tau-star}. 
Thus \eqref{def:tau} implies hypothesis $ (ii) $ of Proposition \ref{prop:separation} about
the separation properties of the bad sites. The second condition on $ \tau $
arises in the measure estimates of Section \ref{sec:mult5} (see Lemma \ref{measure:tot1}).
Moreover the condition \eqref{def:tau} on $ \tau $ is also used in the proof of Proposition 
\ref{prop:RI}, see \eqref{condpert}.   
\item 
{\bf Choice of $ \bar \chi  $.} 
Then we choose a constant $ \bar \chi $ 
such that 
\be\label{def:chi}
\bar \chi \varsigma > C_1 \, , \quad \bar \chi > \tau + s_0 + d  \, ,  
\ee
where $ \varsigma := 1 / 10  $ is fixed as in \eqref{def:varsigma} and  
the constant $ C_1 \geq 2 $ is defined in Proposition \ref{prop:separation}. 
The constant $ \bar \chi $ is the exponent which enters in the definition
in \eqref{def:Nk-multi} of the 
scales $ N_{k+1} = [N_k^{{\bar \chi}}] $
along the multiscale analysis. 
Note that the first inequality in \eqref{def:chi} is condition \eqref{chil}  in  the multiscale step Proposition \ref{propinv}. 
The second condition on $ \bar \chi  $ 
arises in the measure estimates of Section \ref{sec:mult5} (see Lemma \ref{lem:complexity1}).
\item
{\bf Choice of $ \tau' $.} 
Subsequently we choose $ \tau' $ large enough so that  the inequalities \eqref{dtC}-\eqref{chi1} hold for
all $ \chi \in [\bar \chi, \bar \chi^2 ] $. This is used in the multiscale argument in the proof of Proposition \ref{inv:N02}. 
We also take
\be\label{futau1}
 \tau' > \tilde \tau' + (\es/2 )
\ee
where $  \tilde \tau' >  \tau $ is the constant provided by Lemma \ref{lem:As} associated to 
$ \tilde \tau = \tau $. 
\item
{\bf Choice of $ s_1 $.}  Finally we choose the Sobolev index  $ s_1 $ large enough so that  \eqref{s1} 
holds for all $ \chi \in [\bar \chi, \bar \chi^2 ] $.  
\end{enumerate}

We define the set of $ L^2 $-$(N, \eta)$-good/bad parameters. 

\begin{definition} {\bf {($ L^2 $-$(N, \eta)$-good/bad parameters)}}
Given  $ N \in \N $, $ \eta \in (0,1] $,  let 
\begin{align}\label{tetabadweak}
B_{N}^0 (j_0; \l, \eta) & :=  \Big\{ \teta \in \R \, : \,   
\| D_m^{-1/2} A_{N, j_0}^{-1} (\e, \l,\theta) D_m^{-1/2} \|_0 > \eta N^{\t}  \Big\}  \\
&  = 
\Big\{ \teta \in \R \, : \,  \exists {\rm \ an \ eigenvalue \ of \ }  
D_m^{1/2} A_{N, j_0} (\e, \l,\theta) D_m^{1/2} \nonumber \\
&  \ \  \qquad \qquad \qquad \qquad {\rm with \ modulus \ less \ than} \ \eta^{-1} N^{-\t}  \Big\}
\nonumber 
\end{align}
where $ \|  \ \|_0 $ is the operatorial $ L^2 $-norm, and 
define the set of $ L^2 $-$(N, \eta)$-good parameters 
\begin{align}\label{weakgood}
{\cal G}_{N, \eta}^0    :=  & \Big\{   \l \in \wtilde \Lambda  \, : \,   
 \forall \, j_0 \in \Z^d \, , \quad
B_{N}^0(j_0;  \l, \eta)   \subset \bigcup_{q = 1, \ldots, N^{2d+\es+4+3\t_0}} I_q   \\
 & \ \  {\rm where} \ I_q \ {\rm are \  intervals \ with \ measure} \ | I_q| \leq N^{-\t}   \Big\}  \, . \nonumber
\end{align}
Otherwise we say that $ \lambda $ is  $ L^2 $-$(N, \eta)$-bad. 
\end{definition}
Given  $ N \in \N $, $ \eta \in (0,1] $, we also define 
\begin{align} \label{Binver}
{\mathtt G}_{N, \eta}^0 
& := \Big\{  \l \in \wtilde \Lambda  \, : \, \| D_m^{-1/2} A_{N}^{-1}( \e, \l)  D_m^{-1/2}  \|_0 \leq  \eta  N^{\t}   \Big\} \, .
\end{align}
Note that the sets $ {\cal G}_{N, \eta}^0  $, $  {\mathtt G}_{N, \eta}^0 $ 
are increasing in $ \eta $, namely 
\be\label{increasing}
\eta < \eta' \quad \Rightarrow \quad {\cal G}_{N, \eta}^0 \subset   {\cal G}_{N, \eta'}^0 \, , \quad  
{\mathtt G}_{N, \eta}^0  \subset {\mathtt G}_{N, \eta'}^0  \, . 
\ee
We also define the set 
\be\label{def:calG}
\tilde {\cal G} := \Big\{ \l \in \Lambda \, : \,  \om = (1+ \e^2 \l ) {\bar \om}_\e \quad \rm{satisfies} \quad
{({\bf NR})}_{\g_2, \t_2}  \Big\} 
\ee
(recall Definition \ref{NRgamtau}) with 
\be\label{ga2t2}
\gamma_2 := \frac{\gamma_1}{2} = \frac{\g_0}{4} \, , \quad  \tau_2  := \frac{\es (\es-1)}{2} + 2(\tau_1 + 2)  \, . 
\ee
and $ \g_1, \tau_1 $ defined in \eqref{def:tau1}. 

Fix 
$ N_0 := N_0 (\e)  $
such that 
\be\label{N_0:large}
1 \leq \e^2 N_0^{\tau + s_0 + d } \leq 2 
\ee
and define the increasing sequence of scales 
\be\label{def:Nk-multi}
N_k = \big[N_0^{{\bar \chi}^{k}}\big] \, , \quad k \geq 0 
 \, . 
\ee
\begin{remark}
Condition \eqref{N_0:large} is used in Lemma  \ref{lemma1}, see \eqref{lesss1}, 
and  in the proof of Proposition \ref{prop:RI}, see \eqref{condpert}.  
The first  inequality $ \e^{-2} \leq  N_0^{\tau + s_0 + d }  $ in \eqref{N_0:large} is also 
used in Lemma \ref{measure0}.  
\end{remark}

Finally we define, for $ \eta \in [1/2, 1] $,  the sets 
\be\label{def: Cantor like set-all-N}
\Lambda (\e; \cc, X_{r, \mu} )
:= \bigcap_{k \geq 1} {\cal G}_{N_k, \eta}^0 \bigcap_{N \geq N_0^2} {\mathtt G}_{N, \eta}^0 \bigcap \tilde {\cal G}
\ee
where $ {\cal G}_{N, \eta}^0  $ is defined in \eqref{weakgood}, the set $ {\mathtt G}_{N, \eta}^0 $
is defined in \eqref{Binver}, and   $ \tilde {\cal G}$ in \eqref{def:calG}. 
These are the sets 
$ \Lambda (\e; \cc, X_{r, \mu} ) \subset \wtilde \Lambda  $ appearing in  the statement of Proposition \ref{propmultiscale}. 
By \eqref{increasing} these sets clearly satisfy the property 
\ref{list1-m} listed in Proposition \ref{propmultiscale}. 
We shall prove the measure  properties \ref{list2-m} and \ref{list3-m}
in Section \ref{sec:mult5}. 

\begin{remark}\label{intersect-good}
The second intersection in \eqref{def: Cantor like set-all-N} is restricted to the indices $ N \geq N_0^2 $ for definiteness: 
we could have set $ N \geq  N_0^{\a(\tau)} $ for some exponent $ \alpha (\tau) $ which increases linearly with $ \tau $. 
Indeed, the  right invertibility properties  of $ {\it \Pi}_N [ {\cal L}_{r, \mu}]_{| {\mathcal H}_{2N}} $ at the 
scales $ N \leq N_0^2 $ are deduced in Section \ref{sec:rightinv} by the unperturbed Melnikov non-resonance conditions
\eqref{1Mel}, \eqref{2Mel+}-\eqref{2Mel}, \eqref{diof:co} and a perturbative argument
which  
holds for $ N \leq N_0^{\a(\tau)} $  with $ \alpha (\tau) $ linear in $ \tau $,
see 
\eqref{condpert}. 
\end{remark}

\section{Right inverse of $ [{\cal L}_{r, \mu}]_{N}^{2N}  $ for $ \bar N  \leq N < N_0^2 $}\label{sec:rightinv}

The goal of this section is to prove the following proposition, 
which implies item \ref{item1-multiP} of Proposition \ref{propmultiscale} with $ N (\e) := N_0^2  $
and  $ N_0 = N_0 (\e ) $ satisfying \eqref{N_0:large}. 

\begin{proposition}\label{prop:RI}
There are $ \bar N $ and $\e_0>0$ such that, for all $ \e \in (0, \e_0) $, for all $  \bar N  \leq N < N_0^2(\e) $, 
 the operator  
 $$
  [{\cal L}_{r, \mu}]_{N}^{2N} = {\it \Pi}_N [ {\cal L}_{r, \mu}]_{| {\mathcal H}_{2N}} \, ,
  $$ 
which is  defined for all $ \lambda \in \wtilde \Lambda $,  has a right inverse 
$ \big( [{\cal L}_{r, \mu}]_{N}^{2N} \big)^{-1}  : 
{\mathcal H}_N \to {\mathcal H}_{2N} $ satisfying, for all $ s \geq s_0 $,   
\be\label{rigth-inv+1}
\Big|  \Big( \frac{[{\cal L}_{r, \mu}]_{N}^{2N}}{1+ \e^2 \l}\Big)^{-1}   \Big|_{\Lip, s} 
\leq  C(s)  N^{\tau_0' +1 } \big( N^{ \loss s } + | r |_{\Lip, +,s} \big) 
\ee
where $\tau'_0$ is a constant depending only on $\tau_0$. Moreover \eqref{MSLip1} holds. 
\end{proposition}

The proof of Proposition \ref{prop:RI} is given in the rest of this section.

We decompose the operator $ {\cal L}_{r, \mu} $ 
in \eqref{def:Lr}, with $ X_{r, \mu} $ defined  in  \eqref{op:i}-\eqref{def:Ar}, as 
\be\label{Lr0l}
{\cal L}_{r, \mu} =  {\cal L}_{\mathtt d}  + \rho(\e, \l, \vphi) 	\\ 
\ee
with 
\begin{align} \label{def:calk} 
&  {\cal L}_{\mathtt d}   := J  \bar{\mu} \cdot \pa_{\vphi} + 
D_V + \mu_k {\cal J} \Pi_{{\mathbb S} \cup {\mathbb F}}^\bot +
\co \Pi_{\mathbb S} \\
& \rho(\e, \l, \vphi) := 
J (\om - \bar{\mu}) \cdot \pa_{\vphi} +  (\mu -\mu_k) {\cal J} \Pi_{{\mathbb S} \cup {\mathbb F}}^\bot + r(\e, \l, \vphi) \, . 
\label{def:rhopN} 
\end{align}
Note that the operator ${\cal L}_{\mathtt d}$ is independent of $ (\e,\l) $ and $ \vphi $, and $ \rho $ is small in $\e $
since 
\be\label{ombare2}
\om = (1+ \e^2 \lambda)  \bar \omega_\e =  
\bar \mu + \e^2 (\zeta +   \lambda \bar \mu + \e^2 \zeta  \lambda)=\bar \mu + O(\e^2) \, , \quad
| \om |_{\lip} = O(\e^2) \, ,
\ee
where $ \bar \mu $ is the vector defined in \eqref{unp-tangential} (see  \eqref{def omep}), 
\be\label{muke2}
\mu= \mu_k + O(\e^2) \quad \hbox{for some $k \in \mathbb F$} \, ,
\quad  | \mu |_{\lip} =O(\e^2)  \, , 
\ee
by item \ref{mukinF} of Definition \ref{definition:Xr}, and $ | r |_{\Lip, +,s_1}  = O (\e^2) $ by 
item \ref{one-p} of Definition \ref{definition:Xr}.

In order to prove Proposition \ref{prop:RI}
we shall first find a right inverse of $ [{\cal L}_{\mathtt d} ]_N^{2N} $, 
thanks to the non-resonance conditions  \eqref{1Mel}, \eqref{2Mel+}-\eqref{2Mel}, \eqref{diof:co},
and we shall prove that it has off-diagonal decay estimates, see \eqref{est:Lds}. Then we shall deduce the existence of 
a right inverse for  $ [{\cal L}_{r, \mu}]_N^{2N} $ by a perturbative Neumann series argument. 

Note that the space $  e^{\ii \ell \cdot \vphi} {\bf H}  $, where $ {\bf H}$ is defined in 
\eqref{phase-space-multi}, is invariant under $ {\cal L}_{\mathtt d}   $ and
\be\label{tempo-spazio}
D_m^{1/2} {\cal L}_{\mathtt d} D_m^{1/2} (e^{\ii \ell \cdot \vphi} h) = 
e^{\ii \ell \cdot \vphi} ( M_\ell h ) \, ,  \quad \forall h = h(x)  \in {\bf H} \, , 
\ee
where $ M_\ell $ is the operator, acting on functions $ h(x) \in {\bf H} $ of the space variable $ x $ only,  defined by 
\begin{align}\label{def:Bell}
M_\ell & := D_m^{1/2} \big( \ii \bar{\mu} \cdot \ell \,  J + D_V 
+ \mu_k {\cal J} \Pi_{{\mathbb S} \cup {\mathbb F}}^\bot + \co \Pi_{\mathbb S} \big)D_m^{1/2} \, . 
\end{align}

\begin{lemma}\label{lem:Bel-inv} {\bf (Invertibility of $ M_\ell $)}
For all $ \ell \in \Z^\es $, $| \ell | < N_0^2 $, the self-adjoint operator $ M_\ell $ of $ {\bf H} $ is invertible and 
\be\label{B-ell-inv}
\| M_\ell^{-1}\|_0 \lesssim \g_0^{-1} \langle \ell \rangle^{\tau_0} \, .  
\ee
\end{lemma}

\begin{pf}
We write $M_\ell =  D_m^{1/2} M_\ell' D_m^{1/2} $ where 
\be \label{defBl2}
M'_\ell := \ii \bar{\mu} \cdot \ell J + D_V + \mu_k {\cal J} \Pi_{{\mathbb S} \cup {\mathbb F}}^\bot + \co \Pi_{\mathbb S}  \, , 
\ee
and we remind that, in case ($i$) of \eqref{phase-space-multi}, $ M_\ell $ acts on  $ H $, and 
we have set  ${\cal J} = 0$ (see \eqref{op:i}), while
in  case ($ii$), $ M_\ell $ acts on  $ H \times H $ and $ {\cal J} ( h_1, h_2, h_3, h_4) = (-h_4, h_3, h_2, - h_1) $, 
by \eqref{def:opJ}, \eqref{leftJA}, \eqref{rightJA}.    

In case $(i)$ of \eqref{phase-space-multi}, $M'_\ell$ is represented, in the Hilbert basis 
$((\Psi_j,0) , (0, \Psi_j))_{j \in \N}$ of ${\bf H}=H$,  by the block-diagonal matrix 
${\rm Diag}_{j \in \N}  M'_{\ell , j} $ where
$$
M'_{\ell , j} := \left(
\begin{array}{cc}
\mu_j + \co \delta_{\mathbb S}^j & \ii \bar{\mu} \cdot \ell \\
- \ii \bar{\mu} \cdot \ell &  \mu_j + \co \delta_{\mathbb S}^j
\end{array}
\right)  \qquad 
{\rm and } \qquad  
\delta_{\mathbb S}^j := \begin{cases}
1 \quad  {\rm if} \  j \in {\mathbb S}  \cr
0 \quad  {\rm if} \  j \notin {\mathbb S}  \, .
\end{cases}  
$$ 
The eigenvalues of $M'_{\ell , j}$ are 
$$
\begin{cases}
\pm \bar{\mu} \cdot \ell + \m_j +\co  \quad \quad  {\rm if} \  j\in {\mathbb S} \cr
\pm \bar{\mu} \cdot \ell + \m_j  \  \qquad  \quad \   {\rm if} \  j\notin {\mathbb S} \, .
\end{cases}
$$
In case $(ii)$  of \eqref{phase-space-multi}, $M'_\ell$ is represented, 
in the Hilbert basis 
$$
\Big(\Big(0,- \frac{\Psi_j}{\sqrt{2}} \ , \ \frac{\Psi_j}{\sqrt{2}} ,0\Big) \, , \ 
\Big(\frac{\Psi_j}{\sqrt{2}} , 0, 0, \frac{\Psi_j}{\sqrt{2}}\Big) \, ,
\ \Big(\frac{\Psi_j}{\sqrt{2}} , 0,0, -\frac{\Psi_j}{\sqrt{2}}\Big) \, , \
\Big(0,\frac{\Psi_j}{\sqrt{2}} , \frac{\Psi_j}{\sqrt{2}} , 0\Big)\Big)_{j \in \N} 
$$ 
of ${\bf H}=H\times H$, 
by  the block-diagonal matrix ${\rm Diag}_{j \in \N} M'_{\ell , j}$, where
$$
\begin{aligned}
& M'_{\ell , j} := \\ 
& 
\begin{pmatrix}
\mu_j + \co \delta_{\mathbb S}^j - \mu_k \delta_{ (\mathbb S \cup \mathbb F)^c}^j  & \ii \bar{\mu} \cdot \ell  &0&0 \\
 - \ii  \bar{\mu} \cdot \ell&  \mu_j + \co \delta_{\mathbb S}^j - \mu_k \delta_{ (\mathbb S \cup \mathbb F)^c}^j & 0& 0 \\
 0&0&\mu_j + \co \delta_{\mathbb S}^j + \mu_k \delta_{ (\mathbb S \cup \mathbb F)^c}^j & \ii \bar{\mu} \cdot \ell \\
 0&0&- \ii \bar{\mu} \cdot \ell & \mu_j + \co \delta_{\mathbb S}^j + \mu_k \delta_{ (\mathbb S \cup \mathbb F)^c}^j
\end{pmatrix} 
\end{aligned}
$$  
and, in this case, the eigenvalues of $M'_{\ell,j}$ are
\be\label{nr}  
\begin{cases}
 \pm \bar{\mu} \cdot \ell + \mu_j  + \co   & \quad {\rm if} \quad   j \in {\mathbb S}  \, ,  \cr
 \pm  \bar{\mu} \cdot \ell + \mu_j   & \quad {\rm if} \quad    j \in {\mathbb F} \, ,  \cr
\pm   \bar{\mu} \cdot \ell + \mu_j \pm \mu_k   & \quad {\rm if} \quad  j \notin {\mathbb S} \cup {\mathbb F} \, . 
\end{cases}
\ee  
By \eqref{diof:co}, \eqref{1Mel}, \eqref{2Mel+}-\eqref{2Mel},
in both cases we have $  \|(M'_{\ell , j})^{-1} \| \leq  \gamma_0^{-1} \la \ell \ra^{\tau_0}$ for all $j \in \N$. Hence
$M'_\ell $ defined in  \eqref{defBl2} is invertible and  $ \|(M'_\ell)^{-1} \|_0 \leq  \gamma_0^{-1} \la \ell \ra^{\tau_0}$. 
In conclusion $M_\ell = D_m^{1/2} M'_\ell D_m^{1/2}$  is invertible and \eqref{B-ell-inv} holds, by the bound 
$ \| D_m^{-1/2} \|_0 \leq C(m) $.
\end{pf} 

\begin{remark}\label{rem:co-new}
The role of the term  $  \co \Pi_{\mathbb S}  $ in \eqref{op:i}, \eqref{def:Ar} is  
precisely to prove  Lemma \ref{lem:Bel-inv}. Otherwise, if $ \co = 0 $, then some of the eigenvalues 
$ \pm \bar \mu \cdot \ell + \mu_j  $, $ j \in {\mathbb S } $,  would vanish. 
We have some flexibility in defining  the 
extended operators  $ {\cal L}_r $, $ {\cal L}_{r, \mu} $  
in the complementary subspace $ H_{{\mathbb S} \cup {\mathbb F}} $: the important  point 
is  to define  a positive definite operator leaving $ H_{{\mathbb S} \cup {\mathbb F}} $ invariant. 
\end{remark}

Let ${\bf H}_N $ be the finite dimensional subspace  of  
$ {\bf H}$  in  \eqref{phase-space-multi}, 
\be\label{def:HN}
{\bf H}_N := \Big\{  h(x) = \sum_{|j| \leq N} u_{j} e^{\ii  j \cdot x}; \, u_j   \in \C^v \Big\} 
\subset {\bf H} \, ,  \
v :=\begin{cases}
2 \  \hbox{in case \eqref{phase-space-multi}-$(i)$} \\  4 \  \hbox{in case \eqref{phase-space-multi}-$(ii)$} 
\end{cases}
\ee
and denote by $ \Pi_{N} $  the corresponding $ L^2_x $-projector.

Note that, since the space ${\bf H}_N$ is not invariant under the operator $M_\ell$, the invertibility of 
the infinite dimensional operator 
$M_\ell$ does not imply the invertibility of the finite dimensional restriction 
$\Pi_N {M_\ell}_{|{\bf H}_N}$. However,  in the next lemma we prove that, for $N$ large enough,
the operator 
\be\label{def:BellN}
M_{\ell, N} :=  [M_{\ell}]_N^{2N}  := \Pi_{N} (M_\ell)_{|{\bf H}_{2N}} : {\bf H}_{2N} \to {\bf H}_N 
\ee
has a {\it right} inverse. 

We shall use 
the following decomposition of the operator $ M_\ell $ in \eqref{def:Bell}:  
\begin{align}\label{def:Bell1}
M_\ell  = D_m^{1/2} \big( 
 \ii \bar{\mu} \cdot \ell \,  J + D_m 
+ \mu_k {\cal J} + \bar R \big)D_m^{1/2}  \, , 
\end{align}
where
the operator  $ \ii \bar{\mu} \cdot \ell \,  J + D_m  + \mu_k {\cal J}  $
is diagonal in the exponential basis $ \{ e^{\ii  j \cdot x}  $, $ j \in \Z^d \} $,  
and the self-adjoint operator 
\be \label{calL-k}  
\bar R := D_V - D_m  - \mu_k {\cal J} \Pi_{{\mathbb S} \cup {\mathbb F}} + \co \Pi_{\mathbb S}  
\ee
satisfies, by \eqref{DeltaV2} and  Lemma \ref{pisig},  the off-diagonal estimate 
\be\label{def:R-decay}
\quad |  \bar R  |_{+, s} < C(s) < + \infty \, , \quad  \forall s \geq s_0  \, . 
\ee

\begin{lemma} {\bf (Right inverse of $M_{\ell, N}  $)}\label{lem:Bell-RI1}
Let $M_{\ell,N}^*   :  {\bf H}_{N} \to {\bf H}_{2N}$ denote the adjoint operator of $M_{\ell , N}$.
For all  $  N \geq \bar N   $ large enough, $ | \ell | \leq N <N_0^2 $, the operator 
$ M_{\ell, N} M^*_{\ell,N} : {\bf H}_N \to {\bf H}_N  $
is invertible and
\be\label{adjoint}
\| \big( M_{\ell, N}  M_{\ell, N}^* \big)^{-1}\|_0 \lesssim  N^{2 \t_0} \, . 
\ee
As a consequence, the operator $B_{\ell ,N}$ defined in \eqref{def:BellN} has the  right inverse 
 \be\label{il-right-inv}
M_{\ell, N}^{-1} :=  M_{\ell, N}^*  (M_{\ell, N}  M_{\ell, N}^* \big)^{-1} : {\bf H}_N \to {\bf H}_{2N} \, ,
\ee
which satisfies 
\be  \label{BlN1}
\| M_{\ell, N}^{-1} \|_0 \lesssim   N^{ \t_0} \, ,  \qquad
|M_{\ell, N}^{-1}|_{s_0} \lesssim  N^{ \t_0 + s_0 + \frac{d}{2}} \, . 
\ee
\end{lemma}

\begin{pf}
Since $ M_\ell = M_{\ell}^* $,
$ M_{\ell, N}^* =  \Pi_{2N} (M_\ell^*)_{|{\bf H}_{N}} =  \Pi_{2N} (M_\ell)_{|{\bf H}_{N}} $. 
Using \eqref{def:Bell1}-\eqref{calL-k}, we have 
\begin{align}\label{Bright-per}
\| M_{\ell, N}^* -  (M_\ell)_{|{\bf H}_N} \|_0 
= \| \Pi_{2N}^\bot (M_\ell)_{|{\bf H}_N} \|_0 
&  \stackrel{\eqref{def:Bell1}} =  
\| \Pi_{2N}^\bot (D_m^{1/2} {\bar R} D_m^{1/2})_{|{\bf H}_N}  \|_0 \nonumber \\
& \lesssim | \Pi_{2N}^\bot  (D_m^{1/2} \bar R D_m^{1/2} )_{|{\bf H}_N}  |_{s_0} \nonumber  \\
& \stackrel{ \eqref{Sm1}} {\lesssim_{s}} N^{-(s-s_0)} | \bar R |_{+, s} \nonumber \\
&  \stackrel{\eqref{def:R-decay}} \leq C(s) N^{-(s-s_0)} \, .   
\end{align}
For all $ | \ell | \leq N < N_0^2 $,   we get, by \eqref{B-ell-inv},  that, for any $  h \in {\bf H}_N  $,   
\be\label{Bell-basso}
\| M_\ell h \|_0 \gtrsim \g_0 \langle \ell \rangle^{- \tau_0} \| h \|_0  \gtrsim \g_0 N^{- \tau_0} \| h \|_0 \, .
\ee 
Choosing $s>s_0+\tau_0 + 2$ in  \eqref{Bright-per}, we deduce by  \eqref{Bell-basso}  that, 
for $ N \geq \bar N  $ large enough, $ \forall h \in {\bf H}_N $, 
\begin{align*}
\| M_{\ell,N}^* h \|_0 & \geq \| M_{\ell} h \|_0 - \big\| \big( M_{\ell, N}^* -  (M_\ell)_{|{\bf H}_N}\big) h  \big\|_0 \\
& \geq \g_0 c N^{ - \tau_0} \| h \|_0 -  C(s)N^{-\tau_0 - 2}  \| h \|_0 \\
&  \gtrsim \g_0  N^{- \tau_0}  \| h \|_0   
\end{align*}
and therefore 
$$
( M_{\ell, N}M_{\ell, N}^* h, h )_0 = \| M_{\ell, N}^* h \|_0^2 
 \gtrsim \g_0^2 N^{- 2 \tau_0}   \| h \|_0^2 \, .
$$
Since ${\bf H}_N $ is of finite dimension, we conclude that $M_{\ell, N}M_{\ell, N}^* $ is invertible 
and \eqref{adjoint} holds (we do not track  anymore the dependence with respect to the constant  $ \g_0 $).  

It is clear that $M_{\ell , N}^{-1}$ defined in \eqref{il-right-inv} is a right inverse of $M_{\ell,N}$. Moreover,
for all $  h \in {\bf H}_N $,   
\begin{align*}
\|M_{\ell, N}^{-1} h \|_0^2&=  
\big( M_{\ell, N} M_{\ell, N}^* (M_{\ell, N} M_{\ell, N}^*)^{-1} h , (M_{\ell, N} M_{\ell, N}^*)^{-1}h \big)_0 \\
& = \big( h, (M_{\ell, N} M_{\ell, N}^*)^{-1}h \big)_0 \\
&  \leq \|(M_{\ell, N} M_{\ell, N}^*)^{-1} \|_0 \| h \|_0^2   
\end{align*}
and \eqref{adjoint} implies that $ \|M_{\ell , N}^{-1} \|_0 \lesssim  N^{\tau_0}$.  The second inequality of \eqref{BlN1} is an 
obvious consequence of the latter.
\end{pf}

Our aim is now to obtain upper bounds of the $| \ |_s$-norms, for $ s \geq s_0 $,  of the right inverse operator 
$ M_{\ell, N}^{-1} $ 
defined in \eqref{il-right-inv}. 
We write
\be\label{decayR1}
\begin{aligned}
M_{\ell, N} = \Pi_{N} (M_\ell)_{|{\bf H}_{2N}}  & =
\big( \Pi_N M_\ell \Pi_N + \Pi_N M_\ell \Pi_N^\bot)_{ | {\bf H}_{2N}} \\
& =  {\cal D}_N + R_1 + R_2 
\end{aligned}
\ee
where, recalling \eqref{def:Bell1},  
\begin{align}
& {\cal D}_N := \Pi_N {\cal D}_{|{\bf H}_N}  \, , \quad 
{\cal D} :=  D_m^{1/2} \big( \ii \bar{\mu} \cdot \ell J  + D_m +  \mu_k {\cal J}  \big) D_m^{1/2} \, , \label{calD-N} \\
& R_1 := \Pi_N (D_m^{1/2} \bar R D_m^{1/2} \Pi_N )_{| {\bf H}_{2N}} 
 \nonumber \\
& 
R_2 := \Pi_N ( D_m^{1/2} \bar R D_m^{1/2}  )_{| {\bf H}_N^\bot \cap {\bf H}_{2N} } = 
\Pi_N ( D_m^{1/2} \bar R D_m^{1/2} \Pi_N^\bot )_{| {\bf H}_{2N}  }
 \, . \nonumber
\end{align}
By \eqref{def:R-decay} we have
\be\label{decayR1R2}
 | R_1 |_{s} \leq C(s)   \, , \quad | R_2 |_{ s} \leq C(s) \, , \quad \forall s \geq s_0 \,  . 
\ee
We decompose accordingly the adjoint operator as
\be\label{decayR2}
\begin{aligned}
M_{\ell, N}^* = \Pi_{2N} (M_\ell)_{| {\bf H}_N} & = 
 \Pi_{N} (M_\ell)_{| {\bf H}_N} + \Pi_{N}^\bot \Pi_{2N} (M_\ell)_{| {\bf H}_N} \\
 & =
 {\cal D}_N + R_1^* + R_2^* \, . 
 \end{aligned}
\ee
By \eqref{calD-N}, \eqref{A-s-adj}, \eqref{decayR1R2}, $ M_{\ell, N}^* $ satisfies the estimate 
\be\label{Belstar}
| M_{\ell, N}^*|_{s} \leq |{\cal D}_N|_s + |R_1^*|_s + |R_2^*|_s \leq CN^2 + C(s) \lesssim_s N^2 \, . 
\ee 
Thus the right inverse  $M_{\ell , N}^{-1} $ defined in \eqref{il-right-inv},
satisfies, by \eqref{inter-AB-tutto} and \eqref{Belstar},  
\begin{align} \label{BlNs}
|M_{\ell , N}^{-1}|_s  &  \lesssim_s   | M_{\ell, N}^*|_{s_0} | (M_{\ell,N} M^*_{\ell,N})^{-1}|_s + 
 | M_{\ell, N}^*|_{s} | (M_{\ell,N} M^*_{\ell,N})^{-1}|_{s_0} \nonumber  \\
 &  \lesssim_s  N^2 | (M_{\ell,N} M^*_{\ell,N})^{-1}|_s \, . 
\end{align}
In Lemma \ref{lem:BlBls} below we shall bound $| (M_{\ell,N} M^*_{\ell,N})^{-1}|_s$ by a multi-scale argument. Without any loss of generality, we
consider the case $(ii)$ of \eqref{phase-space-multi}.  We first give a general result which 
is a reformulation of the multiscale step 
Proposition 4.1 of \cite{BBo10}, with stronger assumptions. Here ${\mathfrak I} := \{ 1,2,3,4 \}$. 

\begin{lemma} \label{lem:As}
Given $\varsigma \in (0,1/2)$, $C_1 \geq 2$, $\tilde \tau>0$, there are $\tilde \tau'> \tilde \tau$ 
(depending only on $\tilde \tau$ and $d$), 
$s^*>s_0$, $\eta >0$, $\tilde N \geq 1$ with the following property. For any $N \geq \tilde N$, for any finite 
$\wtilde{E} \subset \Z^d \times {\mathfrak I}$ with ${\rm diam}(\wtilde{E}) \leq N$, for any 
$A=D+R \in {\cal M}_{\wtilde{E}}^{\wtilde{E}}$ with $D$ diagonal, assume that
\begin{itemize}
\item[i)] {\bf ($ L^2 $-bound)} $ \|A^{-1} \|_0 \leq N^{\tilde \tau} $, 
\item[ii)]  {\bf (Off-diagonal decay)} $|R|_{s^*} \leq \eta $,  
\item[iii)] {\bf (Separation properties of the singular sites)}  There is a partition of the singular sites 
$ \Omega : = \{ i \in \wtilde{E} \ : \  |D^i_i| < 1/4  \}
\subset \dps \cup_\alpha \Omega_\alpha $  with
$$
{\rm diam} (\Omega_\alpha) \leq N^{C_1 \varsigma / (C_1+1)} \, , \quad 
 {\rm d}(\Omega_\alpha , \Omega_\beta) \geq 
N^{2 \varsigma / (C_1+1)} \, , \   \forall \alpha \neq \beta \, . 
$$
\end{itemize}
Then 
\be \label{estAs}   
 |A^{-1}|_s \leq C(s) N^{\tilde \tau'} (N^{\varsigma s} + |R|_s) \, , \quad \forall s \geq s_0 \, . 
\ee  
\end{lemma}
As usual, what is interesting  in bound \eqref{estAs} is its ``tamed'' dependence with respect to
$N$ ($\varsigma <1$); the constant $C(s)$ may grow very strongly with $s$, but it does not depend on $N$.

\begin{lemma} \label{lem:BlBls}
There exist $t_0 $,  depending only on $\tau_0 $, and $\bar{N} $ such that for all $  \bar N  \leq N < N_0^2$, $| \ell | \leq N$, we have that, $ \forall s \geq s_0 $,
\be \label{BlBls}
 | (M_{\ell,N} M^*_{\ell,N})^{-1}|_s \leq C(s) N^{t_0 + \varsigma s} \, .
\ee
\end{lemma}

\begin{pf}
We identify the operator $ {\cal D}_N  $
in \eqref{calD-N} with the diagonal matrix  
$$ 
{\rm Diag}_{j\in [-N,N]^d}{\cal D}_j^j 
$$
where each ${\cal D}_j^j$ is in ${\cal L} (\C^4)$ (we are in case $(ii)$ of \eqref{phase-space-multi}).
Using the unitary basis of $\C^4$ defined in \eqref{MATU2}, we  identify each ${\cal D}_j^j$ with
 the $ 4 \times 4 $ diagonal matrix (see \eqref{diag-caso2})
\be\label{eige-diag}
{\cal D}_j^j=   {\rm  Diag} ({\cal D}_{(j, \mathfrak a)}^{(j , \mathfrak a)})_{{\mathfrak a} \in \{ 1,2,3,4 \}} \, , 
\quad {\cal D}_{(j, \mathfrak a)}^{(j , \mathfrak a)} := \la j \ra_m \big(  \la j \ra_m + \s_1({\mathfrak a}) \mu_k 
+ \s_2 ({\mathfrak a}) \bar{\mu} \cdot \ell  \big)
\ee
with signs $ \s_1(\mathfrak a), \s_2 (\mathfrak a ) $ defined as in \eqref{def:signs}.

Note that the singular sites of $ M_\ell $ are those 
in \eqref{sing-teta} with $\om = \bar{\mu}$, $\mu =\mu_k$, $ \theta = 0 $ and $ \ell $ fixed.

By Lemma \ref{Bourgain}, 
which we can apply 
with  $ K = 4 $ ($ \ell $ being fixed), for
$\Gamma \geq \bar{\Gamma}(\Theta)$, every $\Gamma$-chain
of singular sites for $M_\ell$ is of length smaller than $(4\Gamma)^{C_2}$. Let us take 
\be  \label{def:C1C2}  
C_1 = 2C_2+3 \, ,  \qquad \Gamma = N^{2/\chi}  
\quad  {\rm with }  \quad  \chi=\varsigma^{-1} (C_1+1) \, . 
\ee 
As a consequence, arguing as at the end of Section \ref{sec:sepabad}, we deduce that,  
for $N$ large enough (depending on $\Theta$), the set $\Omega$ of the singular sites for $M_{\ell,N}$ can be partitioned as  
\be \label{Bllpart}  
\Omega = \dps  \cup_\alpha \Omega_\alpha \  , \quad {\rm with}   \quad 
{\rm diam} (\Omega_\alpha)  \leq  (4\Gamma)^{C_2} \times \Gamma 
\leq N^{C_1/\chi} \, ,  \ {\rm d}(\Omega_\alpha , \Omega_\beta)  > N^{2/\chi} \, .
\ee 
We now multiply 
\be\label{proMlMstar}
M_{\ell, N} M^*_{\ell,N} \stackrel{ \eqref{decayR1}, \eqref{decayR2}} 
= {\cal D}_N^2 + {\cal D}_N R_1^* +  R_1 {\cal D}_N +  R_1 R_1^*  +   R_2 R_2^*
\ee
(notice that $ {\cal D}_N R_2^* $, $ R_2 {\cal D}_N  $, $ R_1 R_2^* $, $ R_2 R_1^* $ are zero)  
in both sides by the diagonal matrix 
\be\label{d-riscaN}
\begin{aligned} 
& d_{\Theta,N}^{-1} := (|{\cal D}_N| + \Theta {\rm Id}_N)^{-1} 
 \qquad  {\rm where } \\
 &    |{\cal D}_N| := 
{\rm Diag}_{|j| \leq N , {\mathfrak a}\in {\mathfrak I}}(|{\cal D}_{j, \mathfrak a}^{j, \mathfrak a}|)_{|j| \leq N} \, , \quad
{\rm Id}_N := \Pi_{N} {\rm Id}_{| {\bf H}_N} \, . 
\end{aligned} 
\ee 
We write
\be\label{def:PlN}
P_{\ell, N}  :=  d_{\Theta,N}^{-1} M_{\ell, N}  M_{\ell, N}^* d_{\Theta,N}^{-1} 
= \wtilde {\cal D}_N + \varrho_{\Theta,N}  
\ee
where, by \eqref{proMlMstar},  
\be\label{def:tildeDN}
\wtilde {\cal D}_N :=  d_{\Theta,N}^{-1} {\cal D}_N^2 d_{\Theta,N}^{-1} \, , 
\qquad 
 \varrho_{\Theta,N} := 
 d_{\Theta,N}^{-1}  \big( {\cal D}_N R_1^* +  R_1 {\cal D}_N +  R_1 R_1^*  + 
  R_2 R_2^* \big) d_{\Theta,N}^{-1} \, . 
\ee
We apply the multiscale Lemma \ref{lem:As} to $P_{\ell,N}$ with  $C_1$  defined in \eqref{def:C1C2} 
and $\tilde \tau : =2\tau_0+5 $.
Let us verify its assumptions. 
The operator $ P_{\ell, N} $ defined in \eqref{def:PlN}  
is invertible as $ M_{\ell, N}  M_{\ell, N}^* $ (Lemma \ref{lem:Bell-RI1}) and, by 
the definition of $ d_{\Theta,N} $ in  \eqref{d-riscaN} and \eqref{eige-diag}, 
for $ N_0^2 > N \geq \bar N $ large enough (depending on $ \Theta $), for all $ |\ell | \leq N $, and using 
\eqref{adjoint}, we obtain  
\begin{align}\label{PelN-inv}
\| P_{\ell, N}^{-1} \|_0 & = 
\big\|  d_{\Theta,N} \big( M_{\ell, N}  M_{\ell, N}^* \big)^{-1}  d_{\Theta,N} \big\|_0 \nonumber \\
& \lesssim N^4 \big\| \big( M_{\ell, N}  M_{\ell, N}^* \big)^{-1} \big\|_0 \nonumber \\ 
& \leq 
N^{2\tau_0+ 5} = N^{\tilde \tau}  
\end{align} 
by the definition of $\tilde \tau : =2\tau_0+5 $. Thus Assumption $i$) of Lemma \ref{lem:As} holds.    

Then we estimate the $| \ |_s$-decay norm of the operator $ \varrho_{\Theta,N} $ in \eqref{def:tildeDN}.
By \eqref{tame-s-decay},  \eqref{decayR1R2} and $  |d_{\Theta,N}^{-1}|_{s} \leq \Theta^{-1} $, 
$ | d_{\Theta,N}^{-1} {\cal D}_N|_{s} \leq 1 $ which directly follow by the definition
\eqref{d-riscaN}, we get, for all $s \geq s_0$,
\be  \label{estrhoN}
\begin{aligned}
|  \varrho_{\Theta,N} |_{s} & \lesssim_s  | d_{\Theta,N}^{-1} {\cal D}_N|_{s} |R_1|_{s} |d_{\Theta,N}^{-1}|_{s} 
+ |d_{\Theta,N}^{-1}|_{s}^2 \big( |R_1|_{s_0}  |R_1|_{s} + |R_2|_{s_0}  |R_2|_{s} \big) \\
& \lesssim_s   \Theta^{-1} \, . 
 \end{aligned}
\ee
In particular,
provided that $\Theta$ has been chosen large enough (depending on $\varsigma, C_1 , \tau_0$),
Assumption $ii)$ of Lemma \ref{lem:As} is satisfied.

To check Assumption $iii)$ 
it is enough to note that,  by the definition of $  \wtilde {\cal D}_N $  in \eqref{def:tildeDN}, and of 
$ d_{\Theta,N}^{-1} $ in \eqref{d-riscaN}, 
for all $i \in  [-N,N]^d \times {\mathfrak I}$, we have 
$$
\wtilde{\cal D}_i^i=\frac{|{\cal D}_i^i|^2}{(|{\cal D}_i^i| + \Theta)^2} \, , \qquad {\rm and} \quad |{\cal D}_i^i| \geq \Theta \ \  \Longleftrightarrow \ \  |\wtilde {\cal D}_i^i| \geq 1/4 \, . 
$$ 
As a consequence, the separation properties for the singular sites of $ M_{\ell,N} $ proved in
 \eqref{Bllpart}  (with $\chi=\varsigma^{-1} (C_1+1)$), imply that Assumption $iii)$ 
of Lemma \ref{lem:As} is satisfied, 
provided that $N$ is large enough (depending only on $\Theta$).
Lemma \ref{lem:As} implies that  there is $t_0$,  depending only on $\tau_0$, such that
(see \eqref{estAs})
$$
|P_{\ell,N}^{-1}|_s \lesssim_s N^{t_0} \big( N^{\varsigma s} + |\varrho_{\Theta,N}|_s \big) \, ,
$$
which gives, using \eqref{def:PlN}, \eqref{tame-s-decay},  $  |d_{\Theta,N}^{-1}|_{s} \leq 1$, 
\be \label{BlBlinterm}
\begin{aligned}
|(M_{\ell,N}   M_{\ell,N}^*)^{-1} |_s = | d_{\Theta ,N}^{-1} P_{\ell, N}^{-1} d_{\Theta ,N}^{-1} |_s
& \lesssim_s |d_{\Theta ,N}^{-1}|_s^2 |P_{\ell, N}^{-1}|_s \\ 
& \lesssim_s  
 N^{t_0} (N^{\varsigma s} + |\varrho_{\Theta,N}|_s) \, . 
 \end{aligned}
\ee
Finally, estimate \eqref{BlBls} follows by  \eqref{BlBlinterm} and \eqref{estrhoN}. 
\end{pf}

\noindent
{\sc Proof of Proposition \ref{prop:RI} concluded.} Recalling the decomposition \eqref{Lr0l},   
the first goal is to define a right inverse of the operator  
\be\label{def:LdN}
{\cal L}_{{\mathtt d},N}: {\mathcal H}_{2N}  \to {\mathcal H}_{N} \, , \quad 
{\cal L}_{{\mathtt d},N}:= {\it \Pi}_N D_m^{1/2} {\cal L}_{\mathtt d} D_m^{1/2}|_{{\mathcal H}_{2N}} \, , 
\ee
  where 
 ${\cal L}_{\mathtt d}$ is defined in \eqref{def:calk}. 
  Recalling \eqref{tempo-spazio}, \eqref{def:BellN} and Lemma \ref{lem:Bell-RI1},  
the linear operator
 ${\cal L}_{{\mathtt d},N}^{-1} : {\mathcal H}_N \to {\mathcal H}_{2N}$ defined by 
 $$
 {\cal L}_{{\mathtt d},N}^{-1} (e^{\ii \ell \cdot \varphi} g)=  e^{\ii \ell \cdot \varphi} M_{\ell,N}^{-1} (g) \, ,  \ \ 
  \forall \ell \in [-N,N]^\es  \, , \ \ \forall g \in {\bf H}_N \, , 
 $$
 is a right inverse of ${\cal L}_{{\mathtt d},N}$. Using that  $ {\cal L}_{{\mathtt d},N}^{-1}  $ is diagonal in time it results 
\be\label{stimaN2}
|{\cal L}_{{\mathtt d},N}^{-1}|_s \leq N^{\es /2} \max_{| \ell | \leq N} |M_{\ell,N}^{-1}|_s 	\, . 
\ee
By \eqref{BlNs} and Lemma \ref{lem:BlBls}, we have the bound
\be \label{BlNs2}
 |M_{\ell,N}^{-1}|_s  \lesssim_s N^{t_0 + \varsigma s +2} \, , \quad  \forall s \geq s_0 \, . 
\ee
Therefore, by  \eqref{stimaN2},  \eqref{BlNs2} and the second estimate in \eqref{BlN1}, we get 
\be \label{est:Lds}
  |{\cal L}_{{\mathtt d},N}^{-1}|_s \lesssim_s N^{\frac{\es}{2} + t_0 + \varsigma s + 2 } \ , \ 
   \forall s \geq s_0   \, , \qquad 
 |{\cal L}_{{\mathtt d},N}^{-1}|_{s_0} \lesssim   N^{\frac{b}{2} + \t_0 + s_0} \, , \quad b = d + \es \, . 
\ee
Finally we use a perturbative Neumann series argument to prove the existence of a
right inverse of 
\be
\begin{aligned}\label{calLrmuN}
{\cal L}_{r,\mu,N}^+ := D_m^{1/2} [{\cal L}_{r,\mu} ]^{2N}_N D_m^{1/2} 
& \stackrel{\eqref{to-have-RI}}  = {\it \Pi}_N D_m^{1/2} {\cal L}_{r,\mu} (\e,\l) {D_m^{1/2}}_{|{\cal H}_{2N}} \\ 
& \stackrel{ \eqref{Lr0l}, \eqref{def:LdN}}    = {\cal L}_{{\mathtt d},N} + \rho_N \, . 
\end{aligned}
\ee
Here
\be\label{def:rNr}
\rho_N:= {\it \Pi}_N D_m^{1/2} \rho (\e,\lambda, \vphi) {D_m^{1/2}}_{|{\cal H}_{2N}} 
\ee
satisfies, by  \eqref{def:rhopN}, \eqref{ombare2}, \eqref{muke2}, Lemma \ref{pisig},
and item \ref{one-p} of Definition \ref{definition:Xr}, 
\be \label{pert-s}
| \rho_N |_{\Lip, s} 
 \lesssim_s \e^2 N^2  + |r|_{\Lip, +,s} \, , \quad | \rho_N|_{\Lip, s_1}  \lesssim_{s_1} \e^2 N^2 \, . 
\ee
Using the second estimates in \eqref{est:Lds} and \eqref{pert-s},  
$ N \leq N_0^2 $, \eqref{N_0:large} and \eqref{def:tau}, we get
\be \label{condpert}
| {\cal L}_{{\mathtt d},N}^{-1} |_{s_0} |\rho_N|_{\Lip,s_0} \lesssim_{s_1} N^{\frac{b}{2}  + \t_0 + s_0+2} \e^2 
\leq N_0^{b  + 2\t_0 + 2s_0 +4} \e^2  \ll 1 \, .
\ee 
Hence $ {\rm Id}_{{\cal H}_N} + \rho_N {\cal L}_{{\mathtt d},N}^{-1} $ is invertible and 
the operator ${\cal L}_{r,\mu,N}^+ $  in \eqref{calLrmuN} has the right inverse 
\be\label{form-inv-SN}
({\cal L}_{r,\mu,N}^+)^{-1} = D_m^{-1/2} ([{\cal L}_{r,\mu} ]^{2N}_N)^{-1} D_m^{-1/2}  :={\cal L}_{{\mathtt d},N}^{-1} \big( 
{\rm Id}_{{\cal H}_N} + \rho_N {\cal L}_{{\mathtt d},N}^{-1} \big)^{-1} 
\ee
which satisfies the following tame estimates (see Lemma \ref{leftinv})
for all $ s \geq s_0 $, 
\be\label{lip-int-JEMS}
 |({\cal L}_{r,\mu,N}^+)^{-1}|_{\Lip,s} 
 \lesssim_s  |{\cal L}_{{\mathtt d},N}^{-1} |_{\Lip,s} + |{\cal L}_{{\mathtt d},N}^{-1}|_{\Lip,s_0}^2 |\rho_N|_{\Lip, s} \,.
\ee
Since $ {\cal L}_{{\mathtt d},N}^{-1} $ does not depend on $\lambda $, 
$ |{\cal L}_{{\mathtt d},N}^{-1}|_{\Lip,s} = |{\cal L}_{{\mathtt d},N}^{-1}|_{s} $, 
 and \eqref{form-inv-SN}, \eqref{lip-int-JEMS}, \eqref{est:Lds}, \eqref{pert-s} imply 
 \begin{align*}
 | D_m^{-1/2} ([{\cal L}_{r,\mu} ]^{2N}_N)^{-1} D_m^{-1/2} |_{\Lip,s}  
 & \lesssim_s N^{\t'_0} \big( N^{\varsigma s} + |r|_{\Lip, +,s} \big)  
\end{align*}
where
$$ 
\t'_0:=\max \big\{ (\es/2) + t_0+2 , b+2s_0+2\t_0 +4 \big\} \, . 
$$ 
Using \eqref{el:AD}, the estimate   \eqref{rigth-inv+1} follows. 

Let us now prove \eqref{MSLip1}. Calling $ \rho' $ the operator defined in \eqref{def:rhopN} 
associated to $ (\mu', r' ) $ we have
$ \rho - \rho' = (\mu - \mu') {\cal J} \Pi_{{\mathbb S} \cup {\mathbb F}}^\bot + r - r' $
and, by \eqref{form-inv-SN},
\begin{align}
& D_m^{-1/2} \big( ([{\cal L}_{r,\mu} ]^{2N}_N)^{-1} -  ([{\cal L}_{r',\mu'} ]^{2N}_N)^{-1} \big) D_m^{-1/2} \nonumber \\
&  = 
{\cal L}_{{\mathtt d},N}^{-1} \Big[ 
\big( {\rm Id}_{{\cal H}_N} + \rho_N {\cal L}_{{\mathtt d},N}^{-1} \big)^{-1} -
\big( {\rm Id}_{{\cal H}_N} + \rho_N' {\cal L}_{{\mathtt d},N}^{-1} \big)^{-1}  \Big] \nonumber \\
&  = 
{\cal L}_{{\mathtt d},N}^{-1} \big( {\rm Id}_{{\cal H}_N} + \rho_N {\cal L}_{{\mathtt d},N}^{-1} \big)^{-1} 
(\rho_N' - \rho_N) {\cal L}_{{\mathtt d},N}^{-1} 
\big( {\rm Id}_{{\cal H}_N} + \rho_N' {\cal L}_{{\mathtt d},N}^{-1} \big)^{-1}  \nonumber \\
& =  ({\cal L}_{r,\mu,N}^+)^{-1} (\rho_N' - \rho_N)  ({\cal L}_{r',\mu',N}^+)^{-1}  \label{pass-sem-ut} \, . 
\end{align}
In conclusion, by \eqref{pass-sem-ut} \eqref{tame-s-decay}, \eqref{lip-int-JEMS}, and 
$ | {\cal L}_{{\mathtt d},N}^{-1} |_{s_0} |\rho_N|_{\Lip,s_1}   \leq 1 $, we get
\begin{align}
\big| D_m^{-1/2} \big( ([{\cal L}_{r,\mu} ]^{2N}_N)^{-1} -  ([{\cal L}_{r',\mu'} ]^{2N}_N)^{-1} \big) D_m^{-1/2} \big|_{s_1} & 
\lesssim_{s_1}  | {\cal L}_{{\mathtt d},N}^{-1} |_{s_1}^2 |  \rho_N - \rho_N' |_{s_1} \nonumber \\
& {\lesssim_{s_1}}  
N^{2(\tau_0'+ \varsigma s_1 )} \big(| \mu - \mu'  | N^2 + | r - r' |_{+,s_1} \big)\nonumber
\end{align}
which, using \eqref{el:AD},  implies \eqref{MSLip1}. This completes the proof of Proposition \ref{prop:RI}.

\begin{remark}
In this section we have proved directly the Lipschitz estimate of $ ([{\cal L}_{r,\mu}]_N^{2N})^{-1} $ 
instead of arguing as in Proposition \ref{inv:N02}, 
because for a right inverse  we do not have the formula \eqref{formula-der-U}.
\end{remark}

\section{Inverse of $ {\cal L}_{r, \mu, N} $ for $ N \geq N_0^2 $}\label{sec:Green}

In Proposition \ref{prop:RI}
we proved item \ref{item1-multiP}  of Proposition \ref{propmultiscale}.
The aim of this section is prove item \ref{item2-multiP} of Proposition \ref{propmultiscale}, namely 
that for all $ N \geq N_0^2 $ and $ \lambda \in \Lambda (\e; 1, X_{r, \mu} ) $
the operator  $ {\cal L}_{r, \mu,N} $ defined in \eqref{def:Ln} is invertible and its inverse
$ {\cal L}_{r, \mu,N}^{-1} $  
has off-diagonal decay, see Proposition \ref{inv:N02}.
The proof is based on inductive applications of the multiscale step Proposition \ref{propinv},
thanks to Proposition \ref{prop:separation} about the separation properties of the bad sites.  

\subsubsection{The set $\Lambda (\e; 1, X_{r, \mu} ) $ is good at any scale}

We first prove the following proposition.

\begin{proposition}\label{prop:par-good} 
The $ \Lambda (\e; 1, X_{r, \mu} ) $ in \eqref{def: Cantor like set-all-N} satisfies 
\be\label{prima-inc} 
\Lambda (\e; 1, X_{r, \mu} ) \subset  
\bigcap_{N \leq N_0} {\cal G}_{N} \bigcap_{k \geq 1} {\cal G}_{N_{k}} 
\ee
where the sets $ {\cal G}_{N} $ are defined in \eqref{good}. 
\end{proposition}

We first consider the small scales $ N \leq N_0 $. 

\begin{lemma}\label{lemma1}
Let $ {\breve A}(\e, \lambda, \theta)  $
be the matrix  in \eqref{Atranslated} corresponding to $ r = 0 $, see Remark \ref{def:r=0}. 
There is $\bar N$ such that  for  all $\bar N \leq  N \leq N_0 \leq (2\e^{-2})^{ \frac{1}{\tau+s_0+d} } $
(see \eqref{N_0:large}), $ \forall \l \in \wtilde \Lambda  $, $ \forall j_0 \in \Z^d $,
$ \theta \in \R $, 
\be\label{N0good-1}
\begin{aligned}
& \quad  \|  D_m^{- 1/2} {\breve A}_{N,j_0}^{-1}(\e, \lambda, \teta) D_m^{-1/2}  \|_0 \leq {N^\tau}  \quad
\Longrightarrow  \\
& |D_m^{-1/2} A_{N, j_0}^{-1} (\e, \lambda, \theta) D_m^{-1/2}|_s \leq N^{\tau'+\varsigma s} \, , \  \forall s \in [s_0, s_1] \, , 
\end{aligned}
\ee
 namely the matrix $ A_{N, j_0}(\e, \lambda, \theta) $ is $ N $-good according to  Definition \ref{goodmatrix}.  
\end{lemma}

\begin{pf}
For brevity, the dependence of the operators 
with respect to $(\e, \lambda)$ is kept implicit.
\\[1mm]
{\bf Step $ 1 $.}  
We first  prove that there is  $\bar N$ 
such that $\forall N \geq \bar N$, $\forall j_0 \in \Z^d $, $\forall \theta \in \R$,
\be \label{rzerofirst}
\begin{aligned}
& \|  D_m^{- 1/2} {\breve A}_{N,j_0}^{-1}( \teta) D_m^{-1/2}  \|_0 \leq {N^\tau} \
\Longrightarrow \\
&  |D_m^{-1/2} \breve{A}_{N, j_0}^{-1} ( \theta) D_m^{-1/2}|_s \leq C(s) N^{\tau_1' +\varsigma s} \, , \  \forall s \geq s_0 \, ,
\end{aligned}
\ee  
where $ \tau'_1 := \tilde \tau' + (\es /2)$ and $\tilde \tau' > \tau $ is the constant 
provided by Lemma \ref{lem:As} associated to 
$ \tilde \tau = \tau $. 

Recall that the matrix $\breve{A} (\theta)$ represents the $L^2$-self-adjoint operator 
$ {\cal L}_{0,\mu} (\theta) $ defined in \eqref{L0mu}, 
which is independent of $ \vphi $. 
For all $\ell \in \Z^\es$,  $ h \in {\bf H}$, we have that 
$$
D_m^{1/2} {\cal L}_{0,\mu} (\theta) D_m^{1/2} (e^{\ii \ell \cdot \varphi} h)=
e^{\ii \ell \cdot \varphi} M_\ell (\theta) h \, , \quad M_\ell (\theta) = {\cal D}_\ell (\theta) + T_\mu \, , 
$$
where 
$$
\begin{aligned}
& {\cal D}_\ell (\theta) := \ii (\om \cdot \ell + \theta) J D_m + D_m^2 + \mu {\cal J} D_m \,  , \\
& T_\mu := D_m^{1/2} \big( D_V-D_m - \mu {\cal J}  \Pi_{{\mathbb S} \cup {\mathbb F}}
+\co \Pi_{\mathbb S}  \big) D_m^{1/2} \, . 
\end{aligned}
$$
Note that, by \eqref{DeltaV2}, Lemma \ref{pisig},  and since $ \mu = O(1) $
(item  \ref{mukinF} of Definition \ref{definition:Xr}), 
it results   
\be \label{est:Tspace}
|T_\mu |_s \leq C(s) \, , \quad \forall s \geq s_0 \, . 
\ee
In order to prove \eqref{rzerofirst},  since
\be\label{sMel}
|  D_m^{- 1/2} {\breve A}_{N,j_0}^{-1}( \teta) D_m^{-1/2}  |_s \leq N^{\es /2} \max_{|\ell| \leq N} 
| (M_\ell (\theta))_{N,j_0}^{-1} |_s  \, , 
\ee
it is sufficient to  bound the $| \ |_s $-norms of $ (M_\ell (\theta))_{N,j_0}^{-1} $.
We apply the multiscale Lemma \ref{lem:As}.

We identify as usual the operator $ {\cal D}_\ell (\theta)  $
with the diagonal matrix  $ {\rm Diag}_{j}(({\cal D}_\ell (\theta))_j^j ) $
where 
$$
\begin{aligned} 
& ({\cal D}_\ell (\theta))_j^j  = 
 {\rm  Diag} ([{\cal D}_\ell(\theta)]_{(j, \mathfrak a)}^{(j , \mathfrak a)})_{{\mathfrak a} \in \{ 1,2,3,4 \}} \, , \\
& [{\cal D}_\ell(\theta)]_{(j, \mathfrak a)}^{(j , \mathfrak a)}
:=  \la j \ra_m \big(  \la j \ra_m + \s_1({\mathfrak a}) \mu 
+ \s_2 ({\mathfrak a}) (\om \cdot \ell + \theta)\big)
\end{aligned}
$$
with signs $ \s_1(\mathfrak a), \s_2 (\mathfrak a ) $ defined as in \eqref{def:signs}.
Note that the singular sites of $ {\cal D}_\ell (\theta) $ are those 
in \eqref{sing-teta} with $ \ell $ fixed.
By Lemma \ref{Bourgain}, which we use with $K=4$, $\ell$ being fixed, there is $C_2$ (independent of $\Theta$) 
such that, for $\Gamma > \bar {\Gamma} (\Theta)$,  any $\Gamma$-chain of singular sites has length
$L \leq (4\Gamma )^{C_2}$.  As in \eqref{def:C1C2}, we can apply this result with $C_1=2C_2+3$,
$\chi=\varsigma^{-1} (C_1+1)$ and $\Gamma =N^{2/\chi}$ (for any $N \geq \bar N$ large enough), 
and we find that the operator $\Theta^{-1} (M_\ell (\theta))_{N,j_0}$ satisfies 
Assumption $iii$) of Lemma  \ref{lem:As}, where we take $\tilde \tau=\tau$. Assumption 
$ii$) is also satisfied, provided that $\Theta$ has been chosen large enough, more
precisely $\Theta \geq \eta^{-1} C(s^*)$, with the constant $C(s^*)$ of \eqref{est:Tspace}.

By Lemma \ref{lem:As}, there is $\tilde \tau' > \tau$ such that $\forall N \geq \bar N$, $\forall j_0 \in \Z^d$, 
$\forall \theta \in \R$, $\forall \ell \in \Z^\es$, 
\be\label{implicazb} 
\| (M_\ell (\theta))_{N,j_0}^{-1} \|_0 \leq N^\tau  \quad   \Longrightarrow  \quad 
\forall s \geq s_0 \, , \  | (M_\ell (\theta))_{N,j_0}^{-1} |_s \lesssim_s  N^{\tilde \tau '} \big( N^{\varsigma s} + |T_\mu|_s \big) \, . 
\ee
Since 
$$
 \|  D_m^{- 1/2} {\breve A}_{N,j_0}^{-1}( \teta) D_m^{-1/2}  \|_0=\max_{|\ell| \leq N} 
\| (M_\ell (\theta))_{N,j_0}^{-1} \|_0 \, , 
$$
the premise in \eqref{rzerofirst} implies the premise in \eqref{implicazb}
and therefore \eqref{sMel}, \eqref{implicazb}, \eqref{est:Tspace} imply \eqref{rzerofirst}
since  $ \tau'_1 := \tilde \tau' + \es /2$. 
\\[1mm]
{\bf Step $ 2$.} We now apply a perturbative argument to the operator 
\be\label{Dpp}
D_m^{ 1/2} { A}_{N,j_0}( \teta) D_m^{1/2} = D_m^{1/2} {\breve A}_{N,j_0}( \teta) D_m^{1/2}+
\rho_{N,j_0} \, ,
\ee
where $ \rho $ is the matrix which represents $ D_m^{1/2} r D_m^{1/2} $. By 
item \ref{one-p} of Definition \ref{definition:Xr},  
\be \label{rhos1}
| \rho_{N,j_0} |_{s_1} \leq |r|_{+,s_1} \leq C_1 \e^2 \, . 
\ee
If $\|  D_m^{- 1/2} {\breve A}_{N,j_0}^{-1}( \teta) D_m^{-1/2}  \|_0 \leq {N^\tau}$ then 
\be \label{lesss1}
|  D_m^{- 1/2} {\breve A}_{N,j_0}^{-1}( \teta) D_m^{-1/2}  |_{s_0} |\rho_{N, j_0}|_{s_1} \lesssim 
N^{\tau + s_0 + (d/2)} \e^2 \lesssim  N_0^{\tau + s_0 + (d/2)} \e^2 \ll  1 
\ee
since $ N_0 \leq (2\e^{-2})^{ \frac{1}{\tau+s_0+d} } $
(see \eqref{N_0:large}). 
Then Lemma \ref{leftinv} implies that
$ D_m^{ 1/2} { A}_{N,j_0}( \teta) D_m^{1/2} $ in \eqref{Dpp} 
is invertible, and 
$ \forall s \in [s_0,s_1] $ 
\begin{align*}
|  D_m^{- 1/2} { A}_{N,j_0}^{-1}( \teta) D_m^{-1/2}  |_{s}  & \lesssim_{s_1}  
|  D_m^{- 1/2} { \breve A}_{N,j_0}^{-1}( \teta) D_m^{-1/2}  |_{s} + 
|  D_m^{- 1/2} { \breve A}_{N,j_0}^{-1}( \teta) D_m^{-1/2}  |_{s_0}^2 |\rho_{N, j_0}|_{s_1}  \\
& \stackrel{\eqref{lesss1}} 
{\lesssim_{s_1}}   |  D_m^{- 1/2} { \breve A}_{N,j_0}^{-1}( \teta) D_m^{-1/2}  |_{s} 
\stackrel{\eqref{rzerofirst}} {\lesssim_{s_1}}
N^{\tau_1' +\varsigma s} \leq   N^{\tau' +\varsigma s} \, ,
\end{align*}
because  $\tau'>\tau_1' = \tilde \tau' + (\es/2 ) $ (see \eqref{futau1}) and for $\bar N$  large enough. 
\end{pf}

At small scales $ N \leq N_0 $, any $ \l \in \wtilde \Lambda $ is  $ N $-good. 

\begin{lemma}\label{Inizioind} {\bf (Initialization)}
For all $ N \leq N_0 $ and $ \e $ small, 
the set  $ {\cal G}_{N} $ defined in \eqref{good} is $  {\cal G}_{N}  = \wtilde \Lambda  $. 
\end{lemma}

\begin{pf}
Lemma \ref{lemma1} implies that 
 $ \forall \l \in \wtilde \Lambda  $, $ \forall j_0 \in \Z^d $, the set
$ B_{N} (j_0; \l ) $ defined in \eqref{tetabad} satisfies 
\begin{align}\label{smallBN0-1}
B_{N} (j_0; \l ) \subset {\breve B}_N^0
:= \Big\{ \theta \in \R   \, : \, &  \ 
 \|  D_m^{- 1/2} {\breve A}_{N,j_0}^{-1}(\e, \l, \teta) D_m^{-1/2}  \|_0 >  N^\tau / 2  \Big\} \, .
\end{align}
Thus, in order to prove that $  {\cal G}_{N}  = \wtilde \Lambda  $, it is sufficient to 
show that 
  the set $ {\breve B}_N^0 $ in \eqref{smallBN0-1} 
satisfies the complexity bound \eqref{BNcomponents}.  Note that,
since $ \| D_m^{-1/2} \|_0 \leq m^{- 1/4} $, we have 
\begin{align}
{\breve B}_N^0  
& \subset \big\{ \teta \in \R \, : \,    \|  {\breve A}_{N, j_0}^{-1} (\e, \l,\theta) \|_0 > C N^{\t}  \, , \ \  C := \sqrt{m} /2 \big\}   
\label{sec-line} \\ 
& = \big\{  \teta \in \R \, : \, \exists {\rm \ an \ eigenvalue \ of \ }  
{\breve A}_{N, j_0} (\e, \l,\theta)  \ {\rm with \ modulus \ less \ than} \ N^{-\t} / C  \big\} \, .  \nonumber 
\end{align}
Let $ \Pi_{N, j_0} $ denote the  $ L^2 $-projector on the subspace 
$$
H_{N, j_0} := \Big\{  (q(x), p(x)) = \sum_{|j - j_0| \leq N} (q_j, p_j) e^{\ii j \cdot x }   \Big\} \, . 
$$
Since $ {\breve A}_{N,j_0} (\e, \l, \teta) $ represents the operator ${\cal L}_{0, \mu} (\theta) $ 
in \eqref{L0mu} which 
does not depend on $ \vphi $ (see Remark \ref{def:r=0}), 
the spectrum of $ {\breve A}_{N, j_0} (\e, \l,\theta)  $ is formed by 
$$
\begin{aligned}
\pm (\om \cdot \ell + \teta)  - \beta_j \, ,   \  j = 1, \ldots , (2N+1)^d  \, , \ \ell \in \Z^\es, \\
\beta_j \ {\rm eigenvalue \ of \ }  \Pi_{N, j_0} \big( D_V + \mu {\cal J} \Pi_{\mathbb S \cup \mathbb F}^\bot + 
\co \Pi_{\mathbb S} \big) \Pi_{N, j_0} \, , 
\end{aligned}
$$
and, by \eqref{sec-line},  we have 
\be\label{incl-B-R}
\begin{aligned}
& \qquad {\breve B}_N^0 \subset 
\bigcup_{|\ell| \leq N, j=1, \ldots, (2N+1)^d, \s = \pm } {\cal R}_{\ell, j}^{ \sigma} \, , \\
& {\cal R}_{\ell, j}^{ \sigma}  := \big\{ \theta \in \R \, : \,  
| \sigma (\theta + \om \cdot \ell ) - \beta_j | \leq N^{-\tau} / C  \big\}  \, . 
\end{aligned}
\ee
It follows that  $ {\breve B}_N^0 $ is included in the union of $  N^{\es + d+ 1} $ intervals $ I_q $ of length $ 2 N^{-\tau}  / C  $. 
By eventually dividing the intervals $ I_q $ we deduce that 
$ {\breve B}_N^0 $  is included in the union of 
$  N^{d+2+ \es } $ intervals $ I_q $ of length $ N^{-\tau}  $.  
\end{pf}

\begin{lemma}\label{Inizioind-new}
For all $ k \geq 0 $ we have
\be\label{GN0N}
 {\cal G}_{N_k} \bigcap {\cal G}_{N_{k+1}, 1}^0 \bigcap \tilde{\cal G} \subset {\cal G}_{N_{k+1}}  
\ee
where the set $ {\cal G}_{N} $ is defined in \eqref{good}, ${\cal G}_{N, \eta}^0 $ in \eqref{weakgood} and
$ \tilde{\cal G} $ in \eqref{def:calG}. 
\end{lemma}

\begin{pf}
Let $ \l \in  {\cal G}_{N_k} \cap {\cal G}_{N_{k+1}, 1}^0 \cap \tilde{\cal G} $. 
In order to prove that $ \l \in {\cal G}_{N_{k+1}}   $ (Definition \ref{def:freqgood}), 
since  $ \lambda \in  {\cal G}_{N_{k+1}, 1}^0 $ (set defined in \eqref{weakgood}), 
it is sufficient to prove that the sets  $ B_{N_{k+1}} (j_0; \l) $  in \eqref{BNcomponents} and 
$ B_{N_{k+1}}^0 (j_0; \l, 1 ) $ in \eqref{tetabadweak} satisfy: 
$$
 \forall j_0 \in \Z^d \, ,  \quad  
 B_{N_{k+1}} (j_0; \l) \subset B_{N_{k+1}}^0 (j_0; \l, 1 ) \, ,
$$
or equivalently, that 
\be\label{inclu-L2-to-Good}
\begin{aligned}
& \| D_m^{-1/2} A_{N_{k+1},j_0}^{-1} ( \e, \l , \teta ) D_m^{-1/2} \|_0 \leq N_{k+1}^\t \quad \Longrightarrow  \\
&  
| D_m^{-1/2} A_{N_{k+1},j_0}^{-1} ( \e, \l , \teta ) D_m^{-1/2} |_s \leq N_{k+1}^{\t' + \varsigma s}  \, , \quad \forall s \in [s_0, s_1] \,.  
\end{aligned}
\ee
We prove \eqref{inclu-L2-to-Good} applying the multiscale step\index{Multiscale step} Proposition \ref{propinv}
to the matrix  $ A_{N_{k+1},j_0} $. By \eqref{off-diago:T} the assumption (H1)  holds.  The assumption
(H2) is the premise in \eqref{inclu-L2-to-Good}. 
Let us verify (H3). 
By Remark \ref{good2}, a site 
\be\label{defE}
k \in E := \Big((0,j_0) + [-N_{n+1}, N_{n+1}]^b\Big) \times {\mathfrak I}  \, ,  
\ee
which is $ N_k $-good for 
$ A(\e, \l, \teta ) := {\cal L}_{r, \mu} + \teta Y  $ 
(see Definition \ref{GBsite} with $ A = A(\e,\l,\teta) $) is  also 
$$
(A_{N_{n+1},j_0}(\e, \l, \teta), N_k)-{\rm good}
$$ 
(see Definition \ref{ANreg} with 
$ A = A_{N_{n+1},j_0}(\e, \l ,\teta) $). As a consequence we have the inclusion
\be\label{badincl}
\begin{aligned}
& \Big\{ \ (A_{N_{n+1},j_0}(\e, \l, \teta), N_k){\rm -bad \ sites} \ \Big\}  \subset \\ 
 & \Big\{ N_k{\rm -bad \ sites \ of} \  A(\e, \l ,\teta) \ {\rm with \ } | \ell | \leq N_{k+1} \Big\}
 \end{aligned}
\ee
and (H3) is proved if the latter $ N_k $-bad sites (in the right hand side of \eqref{badincl}) 
are contained in a disjoint union $ \cup_\a \Om_\a $ of clusters satisfying \eqref{sepabad} (with $ N = N_k $). 
This is a consequence of  Proposition \ref{prop:separation} applied to 
the infinite dimensional matrix  $ A(\e, \l ,\teta) $. 
Since $ \l \in {\cal G}_{N_k} $ 
 then
assumption (i) of Proposition \ref{prop:separation} holds with $ N = N_k $. 
Assumption (ii) holds by \eqref{def:tau}. 
Assumption (iii) of Proposition \ref{prop:separation} holds because $ \l \in \tilde {\cal G} $, see \eqref{def:calG}. 
Therefore the $ N_k $-bad sites of  $  A(\e, \l ,\teta) $ satisfy  \eqref{separ} with $ N = N_k $, and therefore
(H3) holds. 

Then the  multiscale step Proposition \ref{propinv} applied to the matrix  $ A_{N_{k+1},j_0} (\e, \l ,\teta) $ implies that if  
$$ 
\| D_m^{-1/2} A_{N_{k+1},j_0}^{-1} ( \e, \l , \teta ) D_m^{-1/2} \|_0 \leq N_{k+1}^\t 
$$
then 
\begin{align}
| D_m^{-1/2} A^{-1}_{N_{k+1},j_0}(\e, \l ,\teta ) D_m^{-1/2} |_s & \leq 
\frac{1}{4} N_{k+1}^{\tau' } \big( N_{k+1}^{\varsigma s}+    | T |_{+,s}  \big)  \nonumber \\
& \stackrel{\eqref{off-diago:T}} \leq N_{k+1}^{\t' + \varsigma s} \, , \
\forall s \in [s_0, s_1]  \, ,  \label{risultAN}
\end{align}
proving \eqref{inclu-L2-to-Good}. 
\end{pf}

\begin{corollary}\label{cor:ind}
For all $ n \geq 1 $ we have
\be\label{verso-la-conclu}
\bigcap_{k=1}^n {\cal G}_{N_{k}, 1}^0 \bigcap \tilde{\cal G} \subset {\cal G}_{N_{n}} \, . 
\ee
\end{corollary}

\begin{pf}
For $ n = 1 $ the inclusion \eqref{verso-la-conclu} follows by \eqref{GN0N} at $ k = 0 $ and the fact that
$ {\cal G}_{N_0} = \wtilde \Lambda  $ by Lemma \ref{Inizioind}. 
Then we argue by induction. Supposing that \eqref{verso-la-conclu} holds at the step $ n $  then 
$$
\bigcap_{k=1}^{n+1} {\cal G}_{N_{k}, 1}^0 \bigcap \tilde{\cal G} = 
 {\cal G}_{N_{n+1}, 1}^0 \bigcap \Big( \bigcap_{k=1}^{n} {\cal G}_{N_{k}, 1}^0 \bigcap \tilde{\cal G}\Big) 
 \stackrel{\eqref{verso-la-conclu}_n} \subset  {\cal G}_{N_{n+1}, 1}^0 \bigcap {\cal G}_{N_{n}} \bigcap  \tilde{\cal G}  
\stackrel{\eqref{GN0N}_n} \subset  {\cal G}_{N_{n+1}} 
$$
proving  \eqref{verso-la-conclu} at the step $ n+1$.
\end{pf}

\noindent
{\sc Proof of proposition \ref{prop:par-good} concluded.} Corollary \ref{cor:ind} implies that 
\be\label{verso-la-conclu-infty}
\bigcap_{k \geq 1} {\cal G}_{N_{k}, 1}^0 \bigcap \tilde{\cal G} \subset  \bigcap_{n \geq 1} {\cal G}_{N_n} \, . 
\ee
Then we conclude that the set $ \Lambda (\e; 1, X_{r, \mu} ) $
defined in \eqref{def: Cantor like set-all-N} satisfies 
$$
\Lambda (\e; 1, X_{r, \mu} )
= \bigcap_{k \geq 1} {\cal G}_{N_k, 1}^0 \bigcap_{N \geq N_0^2} {\mathtt G}_{N, 1}^0 \bigcap \tilde {\cal G}
\stackrel{\eqref{verso-la-conclu-infty}} \subset  \bigcap_{n \geq 1} {\cal G}_{N_n} 
$$
proving \eqref{prima-inc}, since $  {\cal G}_{N}  = \wtilde \Lambda  $, for all $ N \leq N_0 $,  by
 Lemma \ref{Inizioind}.  


\subsubsection{Inverse of $ {\cal L}_{r, \mu, N} $ for $ N \geq N_0^2 $}

We can finally prove the following proposition.  

\begin{proposition}\label{inv:N02}
For all $ N \geq N_0^2 $, $ \lambda \in \Lambda (\e; 1, X_{r, \mu} ) $,  the operator 
$ {\cal L}_{r, \mu,N} $ defined in \eqref{def:Ln} is invertible and satisfies \eqref{Lip-sDM}. Moreover 
\eqref{MSLip2} holds. 
\end{proposition}

\begin{pf}
Let $ \lambda \in \Lambda (\e; 1, X_{r, \mu} ) $. 
For all $ N \geq N_0^2 $ there is $ M \in \N $ such that $ N = M^\chi $, for some $ \chi \in [\bar \chi, \bar \chi^2 ] $
and $ \l \in {\cal G}_M $. 
In fact 
\begin{enumerate}
\item 
If $ N \geq N_1 $, then $ N \in [ N_{n+1}, N_{n+2 }]$ for some
$ n \in \N $, and we have $ N = N_n^\chi $ for some $ \chi \in [\bar \chi, \bar \chi^2 ] $. 
Moreover if $ \lambda \in \Lambda (\e; 1, X_{r, \mu} ) $ then $ \l \in  {\cal G}_{N_n } $ by \eqref{prima-inc}. 
\item 
If  $ N_0^2 \leq N < N_1 $, it is enough to write
$ N = M^{\bar \chi}$ for some integer $ M < N_0 $. Moreover if $ \lambda \in \Lambda (\e; 1, X_{r, \mu} ) $ then 
 $ \l \in  {\cal G}_{M} $ by \eqref{prima-inc}. 
\end{enumerate}

We  now apply the multiscale step Proposition \ref{propinv} to the matrix $ A_N (\e, \l) $ 
(which represents $ {\cal L}_{r, \mu,N} $ as stated in Remark \ref{rem:ANL}), for $ N \geq N_0^2 $, with 
$ E = [-N, N]^b  \times \mathfrak I  $ and $ N' \rightsquigarrow N $, $ N \rightsquigarrow  M $.  
The assumptions \eqref{dtC}-\eqref{s1} hold, for all $ \chi \in [ \bar \chi, \bar \chi^2 ] $,
by the choice of the constants $ \bar \chi $, $ \tau' $, $ s_1 $ at the beginning of Section \ref{sec:constants}.  
Assumption (H1) holds by \eqref{off-diago:T}.
Assumption (H2) holds because $ \l \in   \Lambda (\e; 1, X_{r, \mu} )   \subset {\mathtt G}_{N,1}^0 $ for 
$ N \geq N_0^2 $, see \eqref{def: Cantor like set-all-N}. 
Moreover, arguing as in Lemma \ref{Inizioind-new} -for the matrix $ A(\e, \theta, \l) $ with $ \teta = 0 $,  $ j_0 = 0 $-,
the hypothesis (H3) of Proposition \ref{propinv} holds. 
Then the multiscale step Proposition \ref{propinv} implies that, 
$ \forall \l \in  \Lambda (\e; 1, X_{r, \mu} )  \subset {\cal G}_{M} \cap {\mathtt G}_{N,1}^0 $, 
we have 
\be\label{boundLn}
| D_m^{-1/2} A_{N}^{-1}  D_m^{-1/2} |_s  \stackrel{(\ref{multi:s})} \leq 
C(s)  N^{\tau' } \big( N^{\varsigma s}  + \e^2  | r |_{+,s} \big) \leq C(s) 
N^{\tau' } \big( N^{\varsigma s} +   | r |_{+,s} \big) \, .
\ee
We claim the following direct consequence: on  the set $  \Lambda (\e; 1, X_{r, \mu}) $ we have 
\begin{align}
\Big|  D_m^{-1/2} \Big( \frac{{\cal L}_{r, \mu, N}}{1 + \e^2 \l}  \Big)^{-1} D_m^{-1/2} \Big|_{\Lip, s} \, & \leq C(s)  
N^{2(\tau' + \loss s_1 +1) }  
\big(   N^{\loss (s - s_1)} + | r |_{\Lip, +,s} \big) \, .  \label{cormul2}
\end{align}
For all $ \l $  in the set $ \Lambda (\e; 1, X_{r, \mu})  $,  
the operator 
$$ 
U(\e, \l) := D_m^{-1/2}  \Big(  \frac{{\cal L}_{r, \mu, N}}{1+ \e^2 \l } \Big)^{-1} D_m^{-1/2}  
$$
satisfies, by \eqref{boundLn} and $ | r |_{+,s_1} \leq C_1 \e^2 $ (see item \ref{one-p} of 
Definition \ref{definition:Xr}), the estimates  
\be\label{boundLnU}
 | U(\e, \l) |_s \leq  C(s)  N^{\tau'} (N^{ \loss s } + | r |_{+,s}) \, ,  
 \quad | U(\e, \l) |_{s_1} \leq  C(s_1)  N^{\tau' + \loss s_1 }  \, .
\ee
Moreover, for all $  \l_1, \l_2 \in \Lambda (\e; 1, X_{r, \mu})  $, we write (using that $\dps \frac{\om}{1+\e^2 \l}$
is independent of $\l$)
\begin{align}
\frac{U(\l_2) - U(\l_1)}{\l_2- \l_1} 
& = - U(\l_2) \frac{U^{-1}(\l_2) - U^{-1}(\l_1)}{\l_2- \l_1} U(\l_1) \nonumber \\
& = - U(\l_2) D_m^{1/2}
\frac{1}{\l_2-\l_1} \Big(  \frac{X_{r,\mu, N}(\l_2)}{1+ \e^2 \l_2} -  \frac{X_{r, \mu, N}(\l_1)}{1+ \e^2 \l_1} \Big) D_m^{1/2} U(\l_1) 
\label{formula-der-U}
\end{align}
where $ X_{r, \mu, N} := {\it \Pi}_N (X_{r, \mu})_{| {\mathcal H}_N} $. 
Decomposing
$$
\frac{1}{\l_2-\l_1} \Big(  \frac{X_{r, \mu, N}(\l_2)}{1+ \e^2 \l_2} -  \frac{X_{r, \mu, N}(\l_1)}{1+ \e^2 \l_1} \Big)
= \frac{X_{r, \mu, N}(\l_2)- X_{r, \mu, N}(\l_1)}{(\l_2-\l_1)(1+ \e^2 \l_2)}
-  \frac{\e^2 X_{r, \mu, N}(\l_1) }{(1+ \e^2 \l_2)(1+ \e^2 \l_1)}  
$$
we deduce 
by \eqref{formula-der-U}, \eqref{boundLnU} 
and  $ \mu = O(1) $, $  |\mu |_{\lip} = O(\e^2 ) $, 
$ | r |_{+,s_1} , | r |_{\lip, +,s_1} \leq C_1 \e^2  $,  the estimates
\be
\begin{aligned}  \label{der:promu}
&  | U |_{\lip, s} \lesssim_s N^{2(\t'+\varsigma s_1 +2 )} ( N^{\varsigma (s- s_1)} 
+ | r |_{+,s} +  | r |_{\lip, +,s}) \, , \\
&  | U |_{\lip, s_1} \lesssim_{s_1}  N^{2(\t'+\varsigma s_1 +1)}  \, . 
\end{aligned}
\ee
Finally \eqref{boundLnU} and \eqref{der:promu} imply \eqref{cormul2}. 
The inequalities \eqref{el:AD} and  \eqref{cormul2} imply \eqref{Lip-sDM}.

We finally prove \eqref{MSLip2}. 
Denoting  $ A_N' $ the matrix  which represents $ {\cal L}_{r', \mu',N} $
as  in Remark \ref{rem:ANL} we have that 
$ A_N  - A_N' $ represents  $ {\it \Pi}_N ( (\mu - \mu') {\cal J} \Pi_{{\mathbb S} \cup {\mathbb F}}^\bot + r - r' )_{| {\cal H}_N} $.
Then  it is enough to write 
\begin{align*}
& \big| D_m^{-1/2} \big( A_N^{-1} - (A_N')^{-1} \big) D_m^{-1/2} \big|_{s_1} \\
&  \lesssim_{s_1}  
| D_m^{-1/2} A_N^{-1} D_m^{-1/2}|_{s_1}  | A_N - A_N' |_{+,s_1} | D_m^{-1/2} (A_N')^{-1} D_m^{-1/2}|_{s_1} \\
& \lesssim_{s_1} N^{2(\tau'+ \varsigma s_1 )} \big(| \mu - \mu'  | N^2 + | r - r' |_{+,s_1} \big)
\end{align*} 
using \eqref{boundLn} at $ s = s_1 $ and the bounds $ | r |_{+, s_1},  | r' |_{+, s_1}  \lesssim \e^2 $. 
This estimate and \eqref{el:AD} imply \eqref{MSLip2}. 
\end{pf}

\section{Measure estimates} 
\label{sec:mult5}

The aim of this section is to prove the measure estimates \eqref{stima-Cantor-voluta}-\eqref{prop:spos} 
in Proposition \ref{propmultiscale}. 
 
\subsubsection{Preliminaries}

We first give  several lemmas on basic properties of eigenvalues
of self-adjoint matrices, which are a consequence of their 
variational characterization\index{Eigenvalues of self-adjoint matrices}.

\begin{lemma}\label{lem:posE}
 Let $ A (\xi) $  be a family of  self-adjoint  
 matrices in $ {\cal M}^E_E $, $ E $ finite, 
defined for $ \xi  \in \wtilde \Lambda \subseteq \R $,    
satisfying, for some $ \beta > 0 $,  
$$  
 {\mathfrak d}_\xi A (\xi)  \geq \beta {\rm Id}  
 $$ 
 (recall the notation \eqref{partial-increase}). 
We list  
the eigenvalues of $A(\xi)$ in non decreasing order 
$$ 
\mu_1 (\xi)  \leq \ldots \leq \mu_q (\xi)  \leq \ldots \leq \mu_{|E|} (\xi) 
$$
according to  their variational characterization 
\be\label{Varcar}
 \mu_q (\xi)  := \inf_{F \in {\cal F}_q } \max_{y \in F , \| y \|_0 = 1} \langle A (\xi)  y, y \rangle_0 
\ee
where $  { \cal F}_q $ is the set of all subspaces $ F $ of $ \C^{|E|} $ of dimension  $ q $.
Then 
\be\label{LBE}
{\mathfrak d}_\xi \mu_q (\xi) \geq \beta > 0 \, , \quad \forall q = 1, \ldots, | E |  \, . 
\ee
\end{lemma}

\begin{pf}
By the assumption $   {\mathfrak d}_\xi A (\xi)  \geq \beta {\rm Id}  $, we have, 
for all $ \xi_2 > \xi_1 $, $ \xi_1 , \xi_2 \in \wtilde \Lambda $, 
$  y \in F $, $ \| y \|_0 = 1 $, that 
$$ 
\langle A (\xi_2 ) y, y \rangle_0 > \langle A ( \xi_1 ) y, y \rangle_0 + \beta ( \xi_2 - \xi_1)  \, . 
$$
Therefore 
$$
\max_{y \in F , \| y \|_0 = 1}  \langle A (\xi_2 ) y, y \rangle_0 \geq \max_{y \in F , \| y \|_0 = 1}  
\langle A (\xi_1 ) y, y \rangle_0 + \beta ( \xi_2 - \xi_1) 
$$
and, by \eqref{Varcar}, for all $ q = 1, \ldots, |E| $,    
$$ 
\mu_q ( \xi_2 ) \geq   \mu_q (\xi_1) + \beta ( \xi_2 - \xi_1) \, . 
$$ 
Hence $  {\mathfrak d}_\xi \mu_q (\xi ) \geq  \beta $.
\end{pf}

\begin{lemma}\label{variatione}
i) Let $ A (\xi) $  be a family of  self-adjoint  matrices in ${\cal M}^E_E $, $ E $ finite, 
Lipschitz with respect to $ \xi   \in \wtilde \Lambda \subseteq  \R $, satisfying,   
for some $ \beta > 0 $, 
$$  
 {\mathfrak d}_\xi A (\xi)  \geq \beta {\rm Id} \, .
$$ 
Then there are intervals
$ (I_q)_{1 \leq q \leq |E|} $ in $ \R $ such that 
\be\label{me2-i}
\Big\{ \xi \in \wtilde \Lambda \, : \, \|A^{-1}(\xi)\|_0 \geq \alpha^{-1} \Big\} 
\subseteq \bigcup_{1 \leq q \leq |E|} I_q \, , 
\qquad  | I_q | \leq  2 \alpha \beta^{-1} \, ,  
\ee
and in particular the Lebesgue measure 
\be\label{meas-Lem1}
\Big| \Big\{ \xi \in \wtilde \Lambda \, : \, \|A^{-1}(\xi)\|_0 \geq \alpha^{-1} \Big\} \Big| \leq
2 |E| \alpha \beta^{-1}  \, . 
\ee
ii) Let $ A (\xi) := Z + \xi W $  be a family of  self-adjoint  matrices in ${\cal M}^E_E $, 
Lipschitz with respect to $ \xi  \in \wtilde \Lambda \subseteq \R $,  with 
 $ W $  invertible and 
$$ 
\beta_1 {\rm Id}  \leq Z \leq \beta_2 {\rm Id} 
$$ 
with $ \beta_1 > 0 $. Then there are intervals
$ (I_q)_{1 \leq q \leq |E|} $ such that 
\be\label{me2}
\Big\{ \xi \in \wtilde \Lambda \, : \, \|A^{-1}(\xi)\|_0 \geq \alpha^{-1} \Big\} 
\subseteq \bigcup_{1 \leq q \leq |E|} I_q \, , 
\qquad  | I_q | \leq  2 \alpha \beta_2 \beta_1^{-1} \|W^{-1}\|_0 \, . 
\ee
\end{lemma}

\begin{pf}   Proof of $i$). 
Let $ (\mu_q(\xi) )_{1\leq q \leq |E|} $ be the eigenvalues of $ A (\xi) $ listed as in Lemma \ref{lem:posE}. 
We have  
$$
\Big\{ \xi \in \wtilde \Lambda \, : \, \|A^{-1}(\xi)\|_0 \geq \alpha^{-1} \Big\} =
\bigcup_{1\leq q \leq |E|} \Big\{ \xi \in \wtilde \Lambda \, : \, \mu_q(\xi) \in [-\a , \a] \Big\}.
$$ 
By \eqref{LBE} each 
$$
I_q := \big\{ \xi \in \wtilde \Lambda \, : \, \mu_q(\xi) \in [-\a , \a] \big\}
$$ is included in an interval
of length less than  $ 2\a \b^{-1}$.
\\[1mm]
Proof of $ii$). 
Let $ U := W^{-1} Z $  and consider the family of self-adjoint matrices 
$$ 
\tilde A(\xi) := A(\xi)  U = (Z+ \xi W) U =  Z W^{-1} Z  + \xi Z  \, . 
$$
We have the inclusion
\begin{align}
\Big\{ \xi \in \wtilde \Lambda \, : \, \| A^{-1}(\xi) \|_0 \geq \a^{-1} \Big\} & \subset 
\Big\{ \xi \in \wtilde \Lambda \, : \, \| {\tilde A}^{-1}(\xi) \|_0 \geq (\a \| U \|_0) ^{-1} \Big\} \, . \label{princ1}
\end{align}
Since  $ {\mathfrak d}_\xi \tilde A(\xi)  \geq \beta_1 {\rm Id} $ we derive by item $i$)  that 
\begin{align}\label{princ15}
\Big\{ \xi \in \wtilde \Lambda \, : \, \| {\tilde A}^{-1}(\xi) \|_0 \geq (\a \| U \|_0) ^{-1} \Big\} 
 \subset  \bigcup_{1 \leq q \leq | E|} I_q 
\end{align}
where $ I_q $ are intervals with measure
\be\label{princ2}
| I_q | \leq 2 \a \| U \|_0 \beta_1^{-1} \leq 2 \a \| W^{-1} \|_0 \| Z \|_0  \b_1^{-1}\leq 2 \a    \| W^{-1} \|_0 \b_2 \b_1^{-1} \, .
\ee
Then \eqref{princ1}, \eqref{princ15}, \eqref{princ2} imply \eqref{me2}. 
\end{pf}

The variational characterization of the eigenvalues also implies  the following lemma. 

\begin{lemma}\label{Lips}
Let $ A $, $ A_1 $   
be self adjoint matrices  ${\cal M}^E_E $, $ E $ finite. Then their eigenvalues, ranked in nondecreasing order, 
satisfy the Lipschitz property
$$
|\mu_q (A) - \mu_q (A_1)| \leq \| A - A_1 \|_0  \, , \quad \forall q =1, \ldots, | E | \, . 
$$
\end{lemma}

We finish this section stating a simple perturbative lemma, proved by a Neumann series argument. 

\begin{lemma}\label{lem:pre1}
Let $ B, B' \in {\cal M}_E^E $ with $ E := [-N, N]^b \times \fracchi $, and assume that 
\be\label{lem:per1-L2}
\| D_m^{-1/2} B^{-1} D_m^{-1/2} \|_0 \leq K_1 \, , \quad \| B'- B \|_0 \leq \alpha \, . 
\ee
If $ 4 (N^2+m)^{1/2} \alpha K_1 \leq 1  $, then 
\be\label{concl:per1-L2}
\| D_m^{-1/2} (B')^{-1} D_m^{-1/2} \|_0 \leq K_1 + 4 (N^2+m)^{1/2} \alpha K_1^2  \, .
\ee
\end{lemma}

\subsubsection{Measure estimate of $ \Lambda \setminus  {\cal \tilde G} $}

We  estimate the complementary set $ \Lambda \setminus  {\cal \tilde G} $ 
where $ {\cal \tilde G} $  is the set defined in \eqref{def:calG}. 

\begin{lemma}\label{ultima stima misura tilde G}
$ | \Lambda \setminus {\cal \tilde G} | \leq \e^2 $.
\end{lemma}

\begin{pf}
Since 
 $ \bar \om_\e := \bar \omega $ satisfies 
 \eqref{NRgt1} 
we have  
\begin{align*}
\big| n +  (1+ \e^2 \lambda)^2 \sum_{1\leq i\leq j \leq \es} 
p_{ij} {\bar \om}_i {\bar \om}_j \big| &  \geq 
\big| n +  \sum_{1\leq i\leq j \leq \es} 
p_{ij} {\bar \om}_i {\bar \om}_j \big| - C | \lambda | \e^2 |p | \\
& \geq  \frac{\g_1}{\langle p \rangle^{\tau_1}} -  C' \e^2 |p | 
\geq \frac{\g_1/2}{\langle p \rangle^{\tau_1}} 
\end{align*}
for all 
$ |p| \leq \big( \frac{\g_1}{2 C' \e^2 } \big)^{\frac{1}{\tau_1 +1}} $. Moreover, 
if $ n = 0 $ then, for all $ \lambda \in \Lambda $, we have 
\begin{align*}
\big| (1+ \e^2 \lambda)^2 \sum_{1\leq i\leq j \leq \es} 
p_{ij} {\bar \om}_i {\bar \om}_j \big|  
\geq (1 - c \e^2 )^2 \frac{\g_1}{\langle p \rangle^{\tau_1}} \, . 
\end{align*}
As a consequence, recalling the definition of $ {\cal \tilde G} $ in \eqref{def:calG},  
Definition \ref{NRgamtau} and $ \gamma_2 := \gamma_1 / 2 $,  we have 
\be\label{p-large-lambda}
\Lambda \setminus {\cal \tilde G} \subset \bigcup_{(n,p) \in {\cal N}} {\cal R}_{n,p} 
\, , \quad {\cal N} := \Big\{ n \neq 0 \, ,  \ |p | > \Big( \frac{\g_1}{2 C' \e^2 } \Big)^{\frac{1}{\tau_1 +1}} \, , \ |n | \leq C |p| \Big\} \, , 
\ee
where
$$
{\cal R}_{n, p} := \Big\{ \lambda \in \Lambda  \, : \, 
| f_{n,p} (\lambda) | <    \frac{2 \g_2}{ \langle p  \rangle^{\t_2}}
\Big\}  \, , \quad  f_{n,p} (\lambda) :=  \frac{n}{(1+ \e^2 \l)^2}+  \sum_{1\leq i\leq j \leq \es} p_{ij} {\bar \om}_i {\bar \om}_j \, .
$$
Since 
$ \pa_\l f_{n,p} (\lambda) :=  - \frac{2n \e^2}{(1+ \e^2 \l)^3} $, if  
$ |n| \neq 0  $, we have $ | {\cal R}_{n, p}| \lesssim \e^{-2} \g_2 \langle p \rangle^{-\tau_2}  $ and 
by \eqref{p-large-lambda} 
\begin{align*}
|\Lambda \setminus {\cal \tilde G}|
\lesssim_{\g_0}  \sum_{|p | > \Big( \frac{\g_1}{2 C' \e^2 } \Big)^{\frac{1}{\tau_1 +1}}} |p|\frac{\e^{-2}}{|p|^{\tau_2}} 
& \lesssim_{\g_0} \e^{-2}  \int_{ \Big( \frac{\g_1}{2 C' \e^2 } \Big)^{\frac{1}{\tau_1 +1}} }^{+\infty} \rho^{-\t_2 + 
\frac12 \es (\es-1) } d \rho \\
& \lesssim_{\g_0} \e^{\frac{2(\t_2 - (\es (\es-1)/2)  - 2 - \tau_1 )}{\t_1 + 1}} \leq \e^2
\end{align*}
for $ \t_2 $ large with respect to $ \tau_1 $, i.e. $ \t_2 $ defined as in \eqref{ga2t2},  and $ \e $ small. 
\end{pf}

\subsubsection{Measure estimate of $ \wtilde \Lambda \setminus  {\mathtt G}^0_{N, \frac12} $ for $ N \geq N_0^2 $}

We  estimate the complementary set $ \wtilde \Lambda   \setminus {\mathtt G}^0_{N, \frac12} $ 
where $ {\mathtt G}^0_{N, \eta} $  is  defined in \eqref{Binver}. 

\begin{lemma}\label{measure0} 
If $ N \geq N_0^2 $
  then  
\be\label{measGN0} 
| \wtilde \Lambda \setminus {\mathtt G}_{N, \frac12}^0  | \leq 
\e^{-2}  N^{- \t +  d+  \es+1} \leq 
 N^{- \frac{\t}{2} + 2 d+  \es+ s_0} \, . 
\ee 
\end{lemma}

\begin{pf}
By \eqref{Binver} we have the inclusion 
$$
 \wtilde \Lambda \setminus {\mathtt G}_{N, \frac12}^0  \subset \big\{ \l \in \wtilde \Lambda \, : \, 
 \|  P_{N}^{-1} (\l)  \|_0 >  N^\tau / 4  \big\} 
$$
where (see Remark \ref{rem:ANL}) 
$$
P_{N} (\l) :=  D_m^{1/2} \frac{A_{N} (\e, \l)}{1 + \e^2 \l}  D_m^{1/2} = 
D_m^{1/2} {\it \Pi}_N \Big[ J \bar \om_\e  \cdot \pa_\vphi + \frac{X_{r, \mu}(\e, \l ) }{1 + \e^2 \l}  \Big]_{|{\cal H}_N} D_m^{1/2}   \, . 
$$
By assumption \ref{assu:pos} of Definition \ref{definition:Xr}, the matrix $ P_{N} (\lambda) $
satisfies  $ {\mathfrak d}_\l  P_{N} (\l)  \leq - c \e^2 $ with $ c := c_2 \sqrt{m} $. 
The first inequality in \eqref{measGN0} follows by Lemma \ref{variatione}-$i$)  (in particular \eqref{meas-Lem1}) applied 
to $ -  P_{N} (\l) $ with $ E = [-N, N]^{d+ \es}  \times \fracchi $, $ \alpha = 4 N^{- \tau} $, $ \beta = c \e^2 $, taking $ N $ large. The second inequality in \eqref{measGN0} follows because, by \eqref{N_0:large}, 
$$ 
\e^{-2} \leq N_0^{\tau+d+s_0}  \leq N^{\frac{\tau+d+s_0}{2}} 
$$ 
for $ N \geq N_0^2 $. 
\end{pf}

\subsubsection{Measure estimate of $ \wtilde \Lambda \setminus  {\cal G}^0_{N_k, \frac12} $ for $ k \geq 1 $}

We  estimate the complementary set $ \wtilde \Lambda \setminus  {\cal G}^0_{N_k, \frac12} $ 
where $  {\cal G}^0_{N, \eta}$  is  defined in \eqref{weakgood}.

\begin{proposition}\label{PNmeas}
For all $ N \geq N_1 := N_0^{\bar \chi} $ 
the  
set $ {\cal B}_{N, \frac12}^0 :=  \wtilde \Lambda \setminus {\cal G}_{N, \frac12}^0  $   has measure 
\be\label{measBN0}
|{\cal B}_{N, \frac12}^0 | \leq   N^{- 1} \, .
\ee
\end{proposition}

We first obtain complexity estimates for the set $ B_N^0 (j_0; \l, 1/2) $ defined in \eqref{tetabadweak}.
We argue differently for $ | j_0 | \geq 3 (1 + |\bar \mu |) N $ and $| j_0| < 3 (1 + |\bar \mu |) N $.

\begin{lemma}\label{lemma:1 inclusione}
For all $ \l \in \wtilde \Lambda $, 
for all $ | j_0 | \geq 3 (1 + |\bar \mu |) N $,  we have
\be\label{1-lemma-inclusione}
B_N^0 (j_0; \lambda, 1/2) \subset \bigcup_{q=1}^{N^{d+\es+1}} I_q
\ee
where $ | I_q| $ are intervals with length $ |I_q| \leq N^{-\tau} $.
\end{lemma}

\begin{pf}
Recalling  \eqref{ANj0}, \eqref{Atranslated}  and \eqref{def:Y}, we have
\begin{align}\label{svipi}
D_m^{1/2} A_{N,j_0}(\e, \l, \theta) D_m^{1/2} & 
= D_m^{1/2} A_{N, j_0}(\e, \l) D_m^{1/2} + \theta \, D_m^{1/2} Y_{N, j_0} D_m^{1/2}  \, .  
\end{align}
We claim that, if $ |j_0| \geq 3 (1 + |\bar \mu |) N  $ and $  N \geq \bar N(V,d,\es) $ is large,  then  
\be \label{positiv} 
\frac{|j_0|^2}{10}  {\rm Id}  \leq D_m^{1/2}  A_{N, j_0}(\e, \l) D_m^{1/2}  \leq 4 |j_0|^2 {\rm Id}   \, . 
\ee
Indeed by \eqref{unitary-A}, \eqref{diao-1-caso}-\eqref{diao-2-caso} 
the eigenvalues $ \nu_{\ell,j} $ of $ A_{N, j_0}(\e, \l) $ satisfy
\be\label{vpA}
\nu_{\ell,j} = \d_{\ell,j}^\pm  + O( \| D_m^{1/2} T_{N,j_0}(\e, \l) D_m^{1/2} \|_0 )  \quad {\rm where} \quad
\d_{\ell,j}^\pm :=   \langle j \rangle_m \big(  \langle j \rangle_m  \pm \om \cdot \ell  \pm \mu ) \, .
\ee
Since $ | \ell | \leq N $ and $ |j| \geq |j_0| -  N $ we see that, for $ |j_0| \geq 3 (1+ | \bar \mu |) N  $, 
for $ N \geq \bar N (V,d,\es)  $ large enough, 
\be\label{lo1}
\frac{ 2 | j_0|^2}{9} \leq \d_{\ell,j}^\pm  \leq  3 |j_0|^2  \, .
\ee 
Hence \eqref{vpA}, \eqref{lo1} and  \eqref{off-diago:T} imply (\ref{positiv}). 
As a consequence, Lemma \ref{variatione}-$ii$) applied to the matrix  in \eqref{svipi} 
with  $ Z = D_m^{1/2} A_{N, j_0}(\e, \l) D_m^{1/2} $, 
$ \a = 2 N^{-\tau} $, $ \b_1 = |j_0|^2 \slash 10 $, $ \b_2 \leq 4 |j_0|^2 $, $ W = D_m^{1/2} Y_{N, j_0} D_m^{1/2} $, $ \| W^{-1} \|_0 \leq C $, 
imply that 
$$
B_N^0 (j_0; \lambda, 1/2) \subset \bigcup_{q=1}^{4(2N+1)^{d+\es}} I_q \qquad {\rm with} \qquad | I_q | \leq  C N^{-\tau} \, . 
$$
Dividing further these intervals we  obtain \eqref{1-lemma-inclusione}.
\end{pf}

We now consider the case $ |j_0| \leq 3 (1 + |\bar \mu |) N $. We can no longer argue directly as in Lemma  
\ref{lemma:1 inclusione}. 
In this case the aim is to bound 
the measure of  
\be\label{B2N0}
B_{2,N}^0(j_0;  \l) :=
\Big\{   \teta \in \R \, : \,  \| D_m^{-1/2} A_{N, j_0}^{-1} (\e, \l,\theta) D_m^{-1/2}  \|_0 >  N^{\t} / 4  \Big\} \, .
\ee
The continuity property of the eigenvalues (Lemma \ref{Lips}) allows then to derive a complexity estimate for 
$ B_N^0 (j_0; \l, 1/2) $
in terms of the measure $ |B_{2,N}^0(j_0;  \l)| $ (Lemma \ref{complessita}). 
Lemma \ref{cor2} is devoted to the estimate of the  bi-dimensional Lebesgue measure 
$$
\Big| \Big\{ (\l, \theta) \in \wtilde \Lambda \times \R \, : \,  \theta \in B^0_{2,N} (j_0; \l) \Big\} \Big| \, . 
$$
Such an estimate is then used in Lemma \ref{intermed} to justify  that 
the measure of the section $ |B^0_{2,N} (j_0; \l)|$ has an appropriate bound for ``most'' $ \l $ 
(by  a Fubini type argument). 

We first show that, for $ |j_0| \leq 3 (1 + |\bar \mu |) N $, the set $ B^0_{2,N} (j_0; \l) $ is contained in an interval
of size $ O(N)  $ centered at the origin.
 
\begin{lemma}\label{corol1}
$ \forall | j_0 | < 3 (1 + |\bar \mu |) N $, $ \forall \l \in \wtilde \Lambda $, we have
$$ 
B_{2,N}^0(j_0;  \l) \subset I_N := \big( - 5 (1+ |\bar \mu |)  N, 5(1+ |\bar \mu|) N \big) \, .
$$
\end{lemma}

\begin{pf}
The eigenvalues $ \nu_{\ell,j}(\teta)$ of $ D_m^{1/2} A_{N,j_0}(\e,\l, \teta) D_m^{1/2} $ satisfy
\be \label{vpA1}
\begin{aligned}
& \nu_{\ell,j}(\theta) = \d_{\ell,j}^\pm (\theta)  + O( | T|_{+, s_1} ) \\
& {\rm where} \quad
\d_{\ell,j}^\pm :=   \langle j \rangle_m \big(  \langle j \rangle_m  \pm (\om \cdot \ell + \theta) \pm \mu )  \, .
\end{aligned}
\ee
If $ |\teta | \geq 5(1+ |\bar \mu|) N $ then, using also \eqref{off-diago:T}, 
each eigenvalue satisfies $ |\nu_{\ell,j}(\teta)| \geq 1 $, and
therefore $ \theta $ belongs to the complementary   of the set $  B_{2,N}^0(j_0;  \l) $
defined in \eqref{B2N0}.  
\end{pf}

\begin{lemma}\label{complessita}
$ \forall |j_0| \leq  3 (1 + |\bar \mu |) N $, $ \forall \l \in \wtilde \Lambda $, we have 
$$ 
B^0_N (j_0;  \l, 1/2) \subset \bigcup_{q=1,..., [\hat C \, {\mathtt M} N^{\tau+1}] } I_q 
$$
where $ I_q $ are  intervals with $ |I_q| \leq N^{-\tau } $ and  $ {\mathtt M} := | B_{2,N}^0(j_0;  \l) |$.
\end{lemma}

\begin{pf}
Suppose that $\theta \in B_N^0 (j_0; \l, 1/2)$ where  $ B_N^0 (j_0; \l, \eta) $ is defined in 
\eqref{tetabadweak}. Then there exists an eigenvalue of 
$ D_m^{1/2} A_{N,j_0}(\e, \l,\theta) D_m^{1/2}  $ with modulus less than $ 2 N^{-\t} $.  
Now, by \eqref{svipi}, and  since $ | j_0 | \leq 3 (1 + |\bar \mu |) N $, we have 
\begin{align*}
\| D_m^{1/2}  \big( A_{N,j_0}(\e, \l,\theta + \Delta  \theta) - A_{N,j_0}(\e, \l,\theta)\big)D_m^{1/2}  \|_0 
&  =
| \Delta \theta | \| D_m^{1/2}   Y_{N,j_0} D_m^{1/2}  \|_0 \\
& \leq  |\Delta \theta| \, 5 (1 + |\bar \mu | ) N   \, .
\end{align*}
Hence, by Lemma \ref{Lips},   if 
$ 5 (1 + |\bar \mu | ) N | \Delta \theta | \leq N^{-\t} $ then $ \theta + \Delta \theta \in B_{2,N}^0 (j_0; \l) $
because $ A_{N,j_0}(\e, \l,\theta + \Delta \theta)$
has an eigenvalue with modulus less than $ 4 N^{-\t}$. 
Hence
$$
[\theta - c N^{-(\t +1)} , \theta + c N^{-(\t +1)}] \subset  B_{2,N}^0 (j_0; \l) \, . 
$$ 
Therefore $ B_{N}^0(j_0; \l, 1/2)  $ is included in an union of  intervals $J_m$ with disjoint interiors,
\be\label{inequa}
B_{N}^0 (j_0; \l, 1/2)  \subset \bigcup_{m} J_m \subset B_{2,N}^0(j_0; \l), \quad {\rm with \ length} \quad  
|J_m|  \geq  2 c N^{-(\t+1)}
\ee
(if some of the intervals $ [\theta - cN^{-(\t+1)} , \theta +c N^{-(\t+1)} ]$ overlap, then we glue them together).  
We decompose each 
$ J_m $  as an union of (non overlapping) intervals $ I_q $ of length between 
$ c N^{-(\t+1)}/2 $ and  $ c N^{-(\tau +1)} $. Then, 
by (\ref{inequa}), we  get a new covering  
$$
\begin{aligned}
& B_{N}^0 (j_0; \l, 1/2)  \subset \bigcup_{q=1, \ldots, Q} I_q  \subset B_{2,N}^0(j_0; \l) \\
&   {\rm with} \ \   
c N^{-(\tau +1)} / 2 \leq |I_q | \leq c N^{-(\tau +1)} \leq N^{-\t}
\end{aligned}
$$
and, since the intervals $ I_q $ do not overlap, 
$$
Qc N^{- (\tau +1)}  \slash 2 \leq \sum_{q = 1}^Q | I_q |   
\leq |  B_{2,N}^0(j_0; \l) | =:   {\mathtt M} \, . 
$$
As a consequence $ Q \leq \hat C \,  {\mathtt M} \, N^{\t+1} $, proving the lemma. 
\end{pf}

In the next lemma we use the crucial sign condition assumption \ref{assu:pos} of Definition \ref{definition:Xr}.  

\begin{lemma}\label{cor2} $ \forall |j_0| < 3 (1+ |\bar \mu|) N $,  the set 
\be\label{B22N}
\begin{aligned}
& {\bf B}^{0}_{2,N}(j_0) := {\bf B}^{0}_{2,N}(j_0; \e) := \\
& \Big\{  (\l, \teta)  \in \wtilde \Lambda \times {\R} \, : \, 
\Big\| D_m^{-1/2} A_{N,  j_0}^{-1}(\e, \l,\teta) D_m^{-1/2} \Big\|_0 > N^{\t}/4  \Big\}
\end{aligned}
\ee
has measure 
\be\label{B2N} 
|{\bf B}^{0}_{2,N}(j_0)| \leq C \e^{-2} N^{-\t+d+\es+1}  \, . 
\ee 
\end{lemma}

\begin{pf} 
By Lemma \ref{corol1},  the set  $ {\bf B}^{0}_{2,N}(j_0)  \subset  \wtilde \Lambda \times I_N $.
In order to estimate the  ``bad" $ ( \l, \teta) $ for which at least one eigenvalue of 
$ D_m^{1/2} A_{N,j_0}(\e, \l, \teta ) D_m^{1/2} $ has modulus less than $ 4 N^{-\t} $, we introduce the variable
\be\label{changevaria}
\vartheta := \frac{\teta}{1 + \e^2 \l}  \quad {\rm where} \quad \vartheta \in 2 I_N \, , 
\ee
and we consider the self adjoint matrix (recall that $ \om = (1+ \e^2 \l) \bar \om_\e $)
\be
\begin{aligned}\label{etaxi}
P_{N, j_0} (\l) & :=  D_m^{1/2} \frac{A_{N,j_0} (\e, \l,\teta)}{1 + \e^2 \l}  D_m^{1/2}	\\
& = 
D_m^{1/2} 
\Big( J \bar \om_\e  \cdot \pa_\vphi + \frac{[X_{r, \mu}(\e, \l ) ]_{N,j_0} }{1 + \e^2 \l} +\vartheta Y_{N, j_0} \Big) D_m^{1/2}  \, . 
\end{aligned}
\ee
By the assumption \ref{assu:pos} of Definition \ref{definition:Xr} we get
$$
{\mathfrak d}_\l  P_{N, j_0} (\l)  \leq - c \e^2 \, , \quad c := c_2 \sqrt{m} \, . 
$$
By Lemma \ref{variatione}-$i$), for each fixed 
$ \vartheta $, the set of $ \l \in \wtilde \Lambda $ such that at least one eigenvalue is $ \leq  4 N^{-\t} $ has measure
at most $ O(  \e^{-2} N^{- \t + d + \es })$. Then, integrating on $ \vartheta \in I_N  $, whose length is $ |I_N| = O(N)$,  
we deduce \eqref{B2N}.
\end{pf}

As a consequence of Lemma \ref{cor2} for ``most" $  \l $ the measure  of  $ B_{2,N}^0 (j_0;  \l ) $ is ``small".

\begin{lemma}\label{intermed}
$ \forall |j_0| < 3  (1+ | \bar \mu|) N $,  the set
\be\label{def:FNj0}
{\cal F}_{N}(j_0) :=  \Big\{ \l \in \wtilde \Lambda \, : \, |B_{2,N}^0 (j_0;  \l)| \geq 
\e^{-2} \hat{C}^{-1} N^{- \t + 2d + \es + 3} \Big\} 
\ee
where $\hat{C}$ is the positive constant of Lemma \ref{complessita}, has measure 
\be\label{sectio}
| {\cal F}_{N}(j_0)| \leq C N^{-d - 2} \, .
\ee
\end{lemma}

\begin{pf}
By Fubini theorem, recalling \eqref{B22N} and \eqref{B2N0}, we have  
\be\label{Fubini}
|{\bf B}^{0}_{2,N}(j_0)| = \int_{\wtilde \Lambda} | B_{2,N}^0(j_0;  \l) | \, d \l  \, .
\ee
Let  $ \mu := \t - 2 d - \es - 3 $.  By (\ref{Fubini}) and (\ref{B2N}), 
\begin{eqnarray}
C \e^{-2} N^{-\t+d+ \es + 1} & \geq & 
\int_{\wtilde \Lambda} |B_{2,N}^0 (j_0;  \l) |  \, d \l \nonumber \\
& \geq & \e^{-2} \hat{C}^{-1} N^{-\mu} \Big| \Big\{ \l \in \wtilde \Lambda  \, : \, |B_{2,N}^0 (j_0; \l)| \geq 
\e^{-2} \hat{C}^{-1} N^{-\mu} \Big\} \Big| \nonumber \\
& := & \e^{-2} \hat{C}^{-1} N^{-\mu} |{\cal F}_{N}(j_0)| \nonumber
\end{eqnarray}
whence (\ref{sectio}).
\end{pf}

As a corollary we get

\begin{lemma}\label{lem:complexity1} 
 Let $ N \geq N_1 := [N_0^{\bar \chi}] $, see \eqref{def:Nk-multi}.
Then $ \forall |j_0| < 3 (1 + |\bar \mu |) N  $,  
$ \forall \l \notin {\cal F}_N(j_0) $, we have  
\be\label{new-inclusi} 
B_{N}^0 (j_0;  \l, 1/2) 
\subset \bigcup_{q=1, \ldots, N^{2d+\es+ 5}} I_q
\ee
with $ I_q $ intervals  satisfying $ |I_q | \leq N^{- \t} $. 
\end{lemma}

\begin{pf}
By the definition of $ {\cal F}_{N}(j_0) $ in \eqref{def:FNj0}, 
for all $ \l \notin {\cal F}_{N}(j_0) $, we have  
$$ 
|B_{2,N}^0(j_0; \l)| <  \e^{-2} \hat{C}^{-1} N^{- \t + 2d + \es + 3 }  \, . 
$$  
Then Lemma \ref{complessita}  implies
that  $ \forall |j_0| < 3 (1 + |\bar \mu |) N $, 
$$ 
B_{N}^0 (j_0;  \l, 1/2) \subset \bigcup_{q=1, \ldots, \e^{-2} N^{2d+\es+4}} I_q \, . 
$$ 
For all $ N \geq N_1 = [N_0^{\bar \chi}] $ 
we have  
$$
\e^{-2} N^{2d+\es+4} \stackrel{\eqref{N_0:large}}\leq N_0^{\tau + s_0 + d} N^{2d + \es + 4} \leq 
C N^{\frac{\tau + s_0+d}{\bar \chi} + 2d + \es + 4} \leq N^{2d + \es + 5} 
$$
by \eqref{def:chi}. This proves \eqref{new-inclusi}.
\end{pf}

\noindent
{\sc Proof of Proposition \ref{PNmeas} concluded.}  
By Lemmata \ref{lemma:1 inclusione}  and \ref{lem:complexity1}, 
for all $ N \geq N_1 $, $ \l \in \wtilde \Lambda $,  
$$ 
\l \notin \bigcup_{|j_0| < 3 (1 + |\bar \mu |) N} {\cal F}_{N}(j_0)   \quad \Longrightarrow \quad 
\l \in  {\cal G}_{N, \frac12}^0 
$$
(see the definition of $ {\cal G}_{N, \frac12}^0  $ in \eqref{weakgood})  and therefore 
\be\label{calBqf}
{\cal B}_{N, \frac12}^0 \subseteq \bigcup_{|j_0| < 3 (1 + |\bar \mu |) N} {\cal F}_{N}(j_0)  \, . 
\ee
In conclusion,  \eqref{calBqf} and \eqref{sectio} imply that, for $ N \geq N_1 $, 
$$
 | {\cal B}_{N, \frac12}^0 | \leq 
\sum_{|j_0| < 3 (1 + |\bar \mu |) N} |{\cal F}_{N}(j_0)|  
\lesssim  N^d  N^{-d - 2} \leq   N^{- 1} \, . 
$$

\subsubsection{Stability of the $ L^2 $-good parameters under variation of $ X_{r, \mu} $}

In order to prove \eqref{prop:spos} we prove 
the ``stability" of the sets 
$ {\mathtt G}_{N, \eta}^0  := {\mathtt G}_{N, \eta}^0 (X_{r, \mu}) $ and  
$ {\cal G}_{N, \eta}^0  := {\cal G}_{N, \eta}^0 (X_{r, \mu}) $ defined respectively 
in  \eqref{Binver} and \eqref{weakgood} with respect to small variations of the operator $ X_{r, \mu} $.

\begin{lemma}\label{lem:inc}
Assume $ | r' - r |_{+, s_1} + | \mu' - \mu | \leq \delta $. 
\\[1mm]
i) If $ N^{\tau+1} \sqrt{\delta} $ is small enough, then, for $ (1/2) + \sqrt{\d} \leq \eta \leq 1  $, 
\be\label{inclu1}
{\mathtt G}_{N, \eta- \sqrt{\delta}}^0 (X_{r, \mu}) \cap \wtilde \Lambda' 
\subset {\mathtt G}_{N, \eta}^0  (X_{r', \mu'}) \, . 
\ee
ii) If $ N_k^{\tau+1} \sqrt{\delta} $, $ k \geq 1 $,  is small enough, then, for $ (1/2) + \sqrt{\d} \leq \eta \leq 1  $,  
\be\label{inclu2}
{\cal G}_{N_k, \eta- \sqrt{\delta}}^0  (X_{r, \mu}) \cap  \wtilde \Lambda' 
\subset {\cal G}_{N_k, \eta}^0  (X_{r', \mu'}) \, . 
\ee
\end{lemma}

\begin{pf} Call $ A_N (\e, \l ) $, resp. $ A_N' (\e, \l ) $, the truncated operator associated to
$ {\cal L}_{r, \mu } $, resp. $ {\cal L}_{r', \mu' } $, see Remark \ref{rem:ANL}, defined for
$ \l \in \wtilde \Lambda $, resp. $ \l \in  \wtilde \Lambda'  $.  
By assumption, for any $ \l \in \wtilde \Lambda \cap \wtilde \Lambda' $
we have 
\be\label{ANA'N'-close}
\| A_N (\e, \l ) - A_N' (\e, \l ) \|_0 \lesssim \| r - r' \|_0 + | \mu - \mu' | \leq C \delta \, . 
\ee
Proof of $i)$.  Assume that $ \l \in {\mathtt G}_{N, \eta- \sqrt{\delta}}^0 (X_{r, \mu}) \cap \wtilde \Lambda' $ where
$ {\mathtt G}_{N, \eta}^0  := {\mathtt G}_{N, \eta}^0 (X_{r, \mu}) $ is defined 
in  \eqref{Binver}. 
Then
\be\label{upAN'}
\| D_m^{-1/2} A_N^{-1} (\e, \l ) D_m^{-1/2} \|_0 \leq (\eta - \sqrt{\delta}) N^\tau \, . 
\ee
Now we apply Lemma \ref{lem:pre1} to $ B = A_N (\e, \l ) $, $ B' = A_N' (\e, \l ) $.
By \eqref{upAN'}, \eqref{ANA'N'-close}, the assumption \eqref{lem:per1-L2} holds 
with $ K_1 = (\eta - \sqrt{\delta})N^\tau $ and $ \alpha = C \delta $. 
If $ \d N^{\tau+1}  $ is small enough then \eqref{concl:per1-L2} applies, and we deduce
$$
\| D_m^{-1/2} (A_N')^{-1} (\e, \l ) D_m^{-1/2} \|_0 \leq (\eta - \sqrt{\delta}) N^\tau +
4 C \d N^{2 \tau + 1 } \leq \eta N^\tau  
$$
provided that $ \sqrt{\d} N^{\tau+1} \leq 1 / (4C ) $. Hence 
$ \l \in {\mathtt G}_{N, \eta}^0 (X_{r', \mu'}) $, proving \eqref{inclu1}.
\\[1mm]
Proof of $ii)$.  
 Assume that $ \l \in {\cal G}_{N_k, \eta- \sqrt{\delta}}^0 (X_{r, \mu}) \cap \wtilde \Lambda' $ where
$ {\cal G}_{N, \eta}^0  := {\cal G}_{N, \eta}^0 (X_{r, \mu}) $
is defined in  \eqref{weakgood}.
Let  $ B_{N_k}^0 (j_0;\l, \eta ) $, resp. $ (B'_{N_k})^0 (j_0;\l, \eta ) $,  be the the set 
defined in \eqref{tetabadweak}  corresponding to  $ X_{r, \mu } $, resp. $ X_{r', \mu' } $, at $ N = N_k $. 
Applying the same perturbative argument of item $i$) 
to the matrices $ D_m^{-1/2} A_{N_k, j_0} (\e, \l, \theta ) D_m^{-1/2} $ and 
$ D_m^{-1/2} A_{N_k, j_0}' (\e, \l, \theta ) D_m^{-1/2} $, 
we prove that, if $ N_k^{\tau+1} \sqrt{\delta} $ is small enough, then, for all $ |j_0| \leq 3(1+ |\bar \mu|) N_k  $, 
we have the inclusion
$$
(B'_{N_{k}})^0 (j_0;\l, \eta ) \subset  B_{N_{k}}^0 (j_0;\l, \eta - \sqrt{\d}) \, . 
$$
Hence, by Lemma \ref{lemma:1 inclusione}, we have $ \l \in {\cal G}_{N_k, \eta}^0 (X_{r', \mu'}) $, proving \eqref{inclu2}. 
\end{pf}

\subsubsection{Conclusion: proof of \eqref{stima-Cantor-voluta}-\eqref{prop:spos}}

We finally prove that the sets
$ \Lambda (\e; \cc, X_{r, \mu} ) $, $ \eta \in [1/2, 1 ] $,  defined in \eqref{def: Cantor like set-all-N}
satisfy the measure estimates \eqref{stima-Cantor-voluta}-\eqref{prop:spos}. 
\\[1mm]
{\sc Proof of \eqref{stima-Cantor-voluta}. }
We have to estimate the measure of the complementary set  
\be\label{compl-CC}
\Lambda (\e; 1/2, X_{r, \mu} )^c \cap \wtilde \Lambda
= \bigcup_{k \geq 1}  {\cal B}_{N_k, \frac12}^0 
\bigcup_{N \geq N_0^2} 
\big( {\mathtt G}_{N, \frac12}^0 \big)^c  \bigcup \tilde {\cal G}^c 
\cap \wtilde \Lambda 
\ee
where  $ {\cal B}_{N_k, \frac12}^0 = \Lambda \setminus {\cal G}_{N_k, \frac12}^0  $
with $  {\cal G}^0_{N, \eta}$  defined in \eqref{weakgood};  the set $ {\mathtt G}_{N, \eta}^0 $
is defined in \eqref{Binver}; the set   $ \tilde {\cal G}$ is defined in \eqref{def:calG}. 

\begin{lemma}\label{measure:totG0} $
\Big| \bigcup_{k \geq 1} {\cal B}_{N_k, \frac12}^0  \Big| \leq \e^2 \, . 
$
\end{lemma}

\begin{pf} 
By  Proposition \ref{PNmeas} we have 
$$
\Big| \bigcup_{k \geq 1} {\cal B}_{N_k,  \frac12}^0  \Big|
\leq 
\sum_{k \geq 1}  N_k^{-1} \lesssim N_1^{-1} \stackrel{\eqref{def:Nk-multi}} 
\lesssim N_0^{- \bar \chi} \stackrel{\eqref{N_0:large}} \lesssim \e^{ \frac{2 \bar \chi}{\t + s_0 + d }} \leq \e^2 
$$
since $ \bar \chi $ is large 
according to the second inequality in \eqref{def:chi}. 
\end{pf}

\begin{lemma}\label{measure:tot1}
$
\Big|  \bigcup_{N \geq N_0^2} \big({\mathtt G}_{N, \frac12}^0 \big)^c  \cap \wtilde \Lambda  \Big| \leq \frac{\e}{3} \, . 
$
\end{lemma}

\begin{pf} 
By Lemma \ref{measure0} we have 
$$
\begin{aligned}
\Big| \bigcup_{N \geq N_0^2} \big({\mathtt G}_{N, \frac12}^0\big)^c \cap \wtilde \Lambda \Big| 
 \leq \sum_{N \geq N_0^2  }  N^{- \frac{\t}{2} + 2 d+ 2 \es+ s_0} 
& \lesssim  N_0^{- \t+ 4d + 4\es + 2s_0 + 2}   \\
& \stackrel{\eqref{N_0:large}}   \lesssim   \e^{2 \frac{ ( \t - 4 d - 4 \es - 2 s_0 - 2)}{\t + s_0+d}} 
\\
&  \leq \frac{\e}{3}  
\end{aligned}
$$
since $ \tau $ is large according to \eqref{def:tau}. 
\end{pf}

By Lemmata \ref{measure:totG0}, \ref{measure:tot1}, \ref{ultima stima misura tilde G} 
we deduce that the complementary set \eqref{compl-CC}
has measure
$$
\begin{aligned}
| \Lambda (\e; 1/2, X_{r, \mu} )^c \cap \wtilde \Lambda | 
& \leq \Big| \bigcup_{k \geq 1} {\cal B}^0_{N_k, \frac12}  \Big|  
+  \Big| \bigcup_{N \geq N_0^2} \big({\mathtt G}_{N, \frac12}^0\big)^c  \cap \wtilde{\Lambda}  \Big|  + | {\cal \tilde G}^c  \cap \wtilde{\Lambda}  | \\ 
& \leq   2 \e^2 + \frac{\e}{3} \leq  \e  
\end{aligned}
$$
proving \eqref{stima-Cantor-voluta}. 

\begin{remark}\label{rem:meas1}
We could prove that the measure of $\Lambda (\e; 1/2, X_{r, \mu} )^c $ is 
smaller than $ \e^p $, for any $ p $, optimizing the choice of the constants. 
Indeed  the measure of the set of Lemma  \ref{ultima stima misura tilde G}, respectively  \ref{measure:totG0},
decreases taking the constant $ \tau_2 $ in \eqref{ga2t2},
respectively $ \bar \chi $ in \eqref{def:chi},  larger.
The set in Lemma \ref{measure:tot1} has measure $ \e^{2-\a}$
as the constant $ \tau \ $ defined in \eqref{def:tau} increases. If we had intersected in 
\eqref{def: Cantor like set-all-N}
for  $ N \geq N_0^{\alpha(\tau)}$, as explained in Remark \ref{intersect-good}, the new set in Lemma \ref{measure:tot1} 
would have arbitrarily
small measure as well.  
\end{remark}

\noindent
{\sc Proof of \eqref{prop:spos}.} 
Recalling the definition of the sets 
$ \Lambda (\e; \cc, X_{r, \mu} ) $, $ \eta \in [1/2, 1 ] $,  in \eqref{def: Cantor like set-all-N}
and Lemma \ref{lem:inc},  we have that, for $ N \geq \bar N $, $ (1/2) + \sqrt{\d} \leq \eta \leq 1  $,  
\begin{align*}
& N^{\tau+2} \sqrt{\d} < 1  \qquad \stackrel{\eqref{inclu1}} \Longrightarrow \qquad 
 \Lambda (\e; \cc -\sqrt{\d}, X_{r, \mu} ) \cap {\mathtt G}_{N, \eta}^0 (X_{r', \mu'})^c 
 \cap \wtilde \Lambda'  = \emptyset 
  \\
& N_k^{\tau+2} \sqrt{\d} < 1 
\qquad \stackrel{\eqref{inclu2}} \Longrightarrow \qquad 
\Lambda (\e; \cc -\sqrt{\d}, X_{r, \mu} ) \cap {\cal G}_{N_k, \eta}^0 (X_{r', \mu'})^c  \cap \wtilde \Lambda'  = \emptyset  \, .  
\end{align*}
Hence 
\begin{align}
& \Lambda (\e; \cc -\sqrt{\d}, X_{r, \mu} ) \cap \Lambda (\e; \cc, X_{r', \mu'} )^c  \cap \wtilde \Lambda'  \nonumber \\
& \subset \bigcup_{ k \geq 1}  \Big( \Lambda (\e; \cc -\sqrt{\d}, X_{r, \mu} )  \cap {\cal G}_{N_k, \eta}^0 (X_{r', \mu'})^c 
 \Big)  \cap \wtilde \Lambda'   \label{ultima-inclu-M}  \\
& 
\ \   \bigcup_{ N \geq N_0^2 } \Big( \Lambda (\e; \cc -\sqrt{\d}, X_{r, \mu} )  \cap {\mathtt G}_{N, \eta}^0 (X_{r', \mu'})^c 
\Big)  \cap \wtilde \Lambda'  \nonumber \\
& 
\subset \bigcup_{ N_k^{\tau+2} \sqrt{\d} \geq 1}  {\cal G}_{N_k, \eta}^0 (X_{r', \mu'})^c 
\!\! \bigcup_{ N^{\tau+2} \sqrt{\d} \geq 1} {\mathtt G}_{N, \eta}^0 (X_{r', \mu'})^c   \cap \wtilde \Lambda'  \nonumber \\
& \subset \bigcup_{ N_k^{\tau+2} \sqrt{\d} \geq 1}  {\cal G}_{N_k, \frac12}^0 (X_{r', \mu'})^c 
\!\! \bigcup_{ N^{\tau+2} \sqrt{\d} \geq 1} {\mathtt G}_{N, \frac12}^0 (X_{r', \mu'})^c  \cap \wtilde \Lambda' \nonumber
\end{align}
by \eqref{increasing}.  
Finally, by \eqref{ultima-inclu-M},  Proposition \ref{PNmeas} 
and Lemma \ref{measure0}  we deduce the measure estimate 
$$
\begin{aligned}
|  \Lambda (\e; \cc -\sqrt{\d}, X_{r, \mu} ) \cap \Lambda (\e; \cc, X_{r', \mu'} )^c \cap \wtilde \Lambda'  | & \leq 
\d^{\frac{1}{2(\tau+2)}} + \d^{\frac{(\tau/2) - 2d + 2\es + 4}{2(\tau+2)}}  \\
& \leq 2 \d^{\frac{1}{2(\tau+2)}} 
\end{aligned}
$$
by 
\eqref{def:tau}.
This proves \eqref{prop:spos} for $ \a < 1 \slash (2 (\tau+2)) $.

\chapter{Nash-Moser theorem} \label{sec:thmNM}

The goal of this chapter is to state 
the Nash-Moser implicit function Theorem \ref{thm:NM}. This result 
proves  the existence of a torus embedding 
$ \vphi \mapsto i(\vphi ) $ of the form  \eqref{torus-embe}
which is a zero of the  nonlinear operator $ {\cal F} $ defined in \eqref{operatorF}.
Theorem \ref{thm:NM} implies, going back to the original coordinates,   
Theorem \ref{thm:main}. 

\section{Statement} 

In this section we state a Nash-Moser implicit function theorem (Theorem \ref{thm:NM}) which proves the 
existence of a solution\index{Torus embedding}
\be\label{torus-embe}
\vphi \mapsto i(\vphi ) =  ( \theta (\vphi) , y(\vphi) , z(\vphi) )  = (\vphi + \vartheta (\vphi) , y(\vphi) , z(\vphi) )   \, , 
\ee
with $ z (\vphi) = (Q(\vphi), P(\vphi) ) \in H_{\mathbb S}^\bot $,  $ \forall \vphi \in \T^\es $,
of the nonlinear operator
\be
\begin{aligned}\label{operatorF}
{\cal F} ( i  ) :=  {\cal F} ( \lambda; i  ) 
 & 
  := \om \cdot \partial_\vphi  i (\vphi ) - X_K ( i (\vphi ))  \\
& = \left(
\begin{array}{c}
\Dom \vartheta (\vphi)  + \om - \bar \mu - \e^2 (\partial_y R) ( i(\vphi), \xi   )   \\
\Dom y (\vphi)  +  \e^2 (\partial_\teta R) ( i(\vphi), \xi  )  \\
\Dom z (\vphi) -  J D_V  z (\vphi) - \e^2 \big( 0, (\nabla_Q R) (i(\vphi), \xi) \big)
\end{array}
\right) 
\end{aligned}
\ee
which depends on the  one dimensional parameter 
$$ 
\l \in \Lambda := [ - \l_0, \l_0 ] 
$$ 
(the set $ \Lambda$ is fixed in \eqref{Lambda-unif})
through the 
 frequency vector  
 $$
  \om = (1+ \e^2 \l ) \bar \om_\e 
  $$ 
  where  $ \bar \om_\e \in \R^\es $ is introduced in \eqref{def omep},
 and the amplitudes $ \xi := \xi (\lambda) $ are defined in \eqref{xil}.  
A solution $ i (\vphi ) $ of  
\eqref{operatorF} is an embedded  invariant torus
for  the Hamiltonian system \eqref{HS0}-\eqref{HS00exp},   filled by quasi-periodic 
solutions with frequency $ \omega $. 

We look for reversible solutions\index{Reversible solution} of $ {\cal F}(\l; i) = 0 $,  namely satisfying $ {\tilde S} i (\vphi ) = i (- \vphi) $ 
(the involution $ {\tilde S} $ is defined in \eqref{involuzione tilde rho}), i.e.  
\begin{equation}\label{parity solution}
\theta(-\vphi) = - \theta (\vphi) \, , \quad 
y(-\vphi) = y(\vphi) \, , \quad 
z (- \vphi ) = ( S z)(\vphi) \, . 
\end{equation}

\begin{remark}\label{rem:rev}
The reversibility property
slightly simplifies the argument in Proposition \ref{Prop:inversione}
because the right hand side in  \eqref{appro-eq-2} has zero average, and therefore the equation 
\eqref{appro-eq-2} is directly solvable. Otherwise 
we would have  to add a counterterm in the second component of the operator $ {\cal F} $ as in \cite{BBField}, \cite{BBM-auto}.
\end{remark}

The Sobolev norm of the periodic component of the embedded torus 
\begin{equation}\label{componente periodica}
{\mathfrak I}(\vphi)  := i (\vphi) - (\vphi,0,0) := ( {\vartheta} (\ph), y(\ph), z(\ph))\,, \quad \vartheta(\ph) := \teta (\vphi) - \vphi \, , 
\end{equation}
is 
\be\label{weighted-Sobo-three}
\|  {\mathfrak I}  \|_{\Lip, s} := \| \vartheta \|_{\Lip, H^s_\vphi} +  \| y  \|_{\Lip, H^s_\vphi} +  \| z \|_{\Lip, s} \, .
\ee
The solutions of  $ {\cal F} ( \lambda; i  ) = 0  $ will be found by a Nash-Moser iterative scheme.
Evaluating $ {\cal F} $ at the trivial embedding 
$$ 
i_0 (\vphi) := (\vphi,0,0) 
$$
we have 
\be
\begin{aligned}\label{operatorF0}
{\cal F} ( i_0 ) = 
& \left(
\begin{array}{c}
\om - \bar \mu   -  \e^2 (\partial_y R) (\vphi, 0, 0, \xi   )   \\
   		          \qquad  \qquad  \e^2 (\partial_\teta R) ( \vphi, 0, 0, \xi  )  \\
	  		 \quad \quad -  \e^2 \big(0, \nabla_Q R (\vphi, 0, 0, \xi) \big)
\end{array}
\right) 
\end{aligned}
\ee
which satisfies, since $ \om $ is $ O(\e^2)$-close to $ \bar \mu $, 
\be\label{0sol}
\| {\cal F} ( i_0 )\|_{\Lip,s} \leq C(s) \e^2  , \quad \forall s \geq s_0   \, .
\ee
In order to construct a better approximate solution we first compute in Section \ref{sec:shifted-tan} the 
shifted tangential frequency vector induced by the nonlinearity, 
up to $O(\e^4 ) $. 
Then in Section \ref{NM:first-step} we construct the first approximate solution $  i_1 (\vphi ) $, defined for all
$ \l \in \Lambda $,  by using the unperturbed Melnikov 
non-resonance conditions  \eqref{diop}-\eqref{1Mel} on the linear unperturbed frequencies, in such a way that
(see \eqref{approx1})
$$
\| {\cal F}(i_1) \|_{\Lip,s} \leq C(s) \e^4 \, , \quad \forall s \geq s_0 \, .
$$
Subsequently, given an approximate solution $ i_n (\vphi)  $, 
the main point 
 is to construct a much better approximate
solution $  i_{n+1}  (\vphi)  $.  We use an inductive  Nash-Moser iterative scheme. 
The key step concerns
the approximate right invertibility properties of  the linearized operators $ d_i {\cal F} (i_n) $
obtained along the iteration, that we  obtain restricting the values of 
$ \lambda $ to  subsets $ {\bf \Lambda}_n \subset \Lambda $ with large measure, see Theorem \ref{thm:NM-ite}. 
The following theorem will be proved in Chapter
\ref{sec:NM}, relying on the results of Chapters 
\ref{sezione almost approximate inverse}-\ref{sec:proof.Almost-inv} concerning the
invertibility properties of the linearized operator $ d_i {\cal F} (i_n) $. 

\begin{theorem}\label{thm:NM}
{\bf (Nash-Moser)}\index{Nash-Moser theorem}
Assume \eqref{positive} and the non-resonance conditions 
\eqref{diop}-\eqref{NRgamma0}, 
\eqref{2Mel+}-\eqref{2Mel rafforzate}. Assume also 
the twist condition
\eqref{A twist} and the  non-degeneracy conditions \eqref{non-reso}-\eqref{non-reso1}. 
Fix a direction 
$ \bar \om_\e := \bar \mu + \e^2 \dom $, $ \dom \in   \Ab ([1,2]^\es ) $,   as in \eqref{def omep}  
such that the Diophantine conditions \eqref{dioep}-\eqref{NRgt1} hold. 
Define $ \Lambda = [-\l_0, \l_0] $ as in \eqref{Lambda-unif}. 

Then there are Sobolev indices  $ s_2 > s_1 > s_0 $, a constant $ \e_0 >  0 $,   and,  
for all $  \e \in (0, \e_0) $  there exist 
\begin{enumerate}
\item \label{item1:TNM}
a Cantor-like set $ {\cal C}_\infty \subset   \Lambda  $ 
of asymptotically full measure as $ \e \to 0 $, i.e.
\be\label{C-infty-AM}
\lim_{\e \to 0} \frac{| {\cal C}_\infty|}{ | \Lambda | } = 1 \, , 
\ee
more precisely, 
there is a map $\e \mapsto b(\e)$, independent of $\dom \in \Ab ([1,2]^\es )$
such that the Diophantine conditions \eqref{dioep}-\eqref{NRgt1} hold, and
 satisfying   
$ | \Lambda \setminus {\cal C}_\infty | \leq b(\e) $, $ \lim_{\e \to 0} b(\e)=0 $; 
\item
a Lipschitz function 
$$
i_\infty (\vphi; \l)  - (\vphi, 0, 0) = ( \vartheta_\infty, y_\infty, z_\infty ) 
: {\cal C}_\infty \to H^{s_2}_\vphi \times H^{s_2}_\vphi \times 
({\mathcal H}^{s_2} \cap H_{\mathbb S}^\bot) 
$$
satisfying 
\be\label{bound-i-infty}
\| i_\infty - (\vphi, 0, 0) \|_{\Lip, s_1} \leq   C(s_1) \e^2  \, , \quad 
\| i_\infty - (\vphi, 0, 0) \|_{\Lip, s_2} \leq  \e  \, , 
\ee
such that the torus $ i_\infty (\vphi; \l)   $, $  \l \in {\cal C}_\infty $,  is a solution of 
$ {\cal F} ( \lambda; i_\infty (\lambda) ) = 0 $. 
\end{enumerate}

Moreover, for any $ \l \in {\cal C}_\infty $, 
the function $ i_\infty - (\vphi, 0, 0) $ is of class  $ C^\infty $ in $ (\vphi, x ) $,
and Lipschitz in $\lambda $ as a map valued in 
$ H^{s}_\vphi \times H^{s}_\vphi \times 
({\mathcal H}^{s} \cap H_{\mathbb S}^\bot)$, $ 	\forall s \geq s_2 $.

As a  consequence the embedded torus 
$ \vphi \mapsto i_\infty (\vphi; \l)  $ is invariant 
for  the Hamiltonian system \eqref{HS0}-\eqref{HS00exp},  
and it is filled by quasi-periodic solutions with frequency $ \om = (1+ \e^2 \l ) \bar \om_\e $.
\end{theorem}

Going back to the original coordinates
via \eqref{D14}, \eqref{tang-norm}, \eqref{AA}, \eqref{tra-xi-y},  
Theorem \ref{thm:NM} implies the existence, for all $ \l \in  {\cal C}_\infty $, 
of a quasi-periodic solution 
of the wave equation \eqref{NLW2} of the form
$$ 
u (t, x) = \sum_{j \in {\mathbb S}} \mu_j^{-\frac12} \sqrt{2(\xi_j 
+ (y_\infty)_j ( \omega t)) } \cos \big( \omega_j t + (\vartheta_\infty )_j( \omega t ) \big) \Psi_j (x) + D_V^{-\frac12} Q_\infty  ( \omega t ) 
$$
with frequency $ \om = (1+ \e^2 \l ) \bar \om_\e $. This proves  Theorem\index{Main result}  \ref{thm:main}
with $ {\mathcal G}_{\e, \zeta} := {\mathcal C}_\infty $.
By \eqref{bound-i-infty} and the fact that $ i_\infty - (\vphi, 0, 0) $ is in  $ C^\infty $,  we deduce, by using interpolation,  that, for any $ s \geq s_0 $,  
$ \| r_{\e} \|_s $ tends to $ 0 $ as $ \e \to 0 $.   

\smallskip

The proof of Theorem \ref{thm:NM} occupies the rest of the monograph 
from Section \ref{sec:shifted-tan} until Chapter
\ref{sec:NM}. 

We  first prove, as a corollary  of Theorem \ref{thm:NM},  
the result  \eqref{density-Leb} about the density, 
close to $ \bar \mu $,  of the 
frequency vectors $ \omega $ of the quasi-periodic solutions of 
\eqref{NLW2} obtained in Theorem \ref{thm:NM}. 

\subsection*{Proof of 
\eqref{density-Leb}} 

Let $\Omega$ be the set of the frequency vectors $\om$ of the quasi-periodic solutions
of \eqref{NLW1} provided by Theorem \ref{thm:main}.
Such frequency vectors 
have the form
\be\label{expr:omega}
\omega =  ( 1 + \e^2 \l) \bar \om_\e = 
\bar \mu  + \e^2 ( \dom + \lambda \bar \mu + \e^2   \lambda \zeta ) \, ,
\quad  \dom \in \Ab ([1,2]^\es ) \, , \  \lambda \in \Lambda \, , 
\ee
where 
\begin{itemize}
\item 
$ \dom \in \Ab ([1,2]^\es ) \setminus B_\e $ (the set $B_\e$ is defined in Lemma \ref{lemma:rho-dioph}),
so that $ \bar \om_\e = \bar \mu + \e^2   \zeta $ satisfies the Diophantine conditions
 \eqref{dioep}-\eqref{NRgt1};
 \item 
$ \lambda \in \Lambda \setminus  {\mathcal G}_{\e, \zeta} $,
$  {\mathcal G}_{\e, \zeta}  := {\cal C}_\infty $ being the set defined in Theorem 
\ref{thm:NM}.
\end{itemize}
We define 
$$
\begin{aligned}
{\frak B}_\e & := \Big\{ (\dom, \lambda) \in \Ab ([1,2]^\es ) \times \Lambda \, : 
\, \dom \in B_\e \ {\rm or} \ \lambda \notin {\mathcal G}_{\e, \zeta} \Big\} \\
{\frak G}_\e & := ( \Ab ([1,2]^\es ) \times \Lambda ) \setminus {\frak B}_\e \, .  
\end{aligned}
$$
 By Lemma \ref{lemma:rho-dioph} the Lebesgue measure of $B_\e$ satisfies
 $ | B_\e | \leq \e $, and, using also the measure estimate provided in 
  item \ref{item1:TNM} of Theorem \ref{thm:NM},  we 
 deduce that 
 \be\label{fraBep} 
 | {\frak B}_\e  | \leq 
 \e |\Lambda | + b(\e) | \Ab ([1,2]^\es )| 
 =:  b_1(\e ) \, ,
 \ee 
with $ \lim_{\e \to 0}  b_1 (\e )  = 0 $. 

In view of \eqref{expr:omega}, 
in order to prove \eqref{density-Leb}, we have to estimate the measure of the set  
\be\label{def:Bep}
B_\e' := \Big\{ \beta \in {\cal C}_1 := \Ab ([1,2]^\es ) + \Lambda \bar \mu \, : \, 
\not \exists (\dom, \l) \in {\frak G}_\e  \ {\rm such \ that} \  
\beta =  \dom + \l \bar \mu + \e^2 \l \dom  \Big\} \, .
\ee

\begin{lemma} \label{mesbe} 
$ | B_\e' | \to 0 $ as $ \e \to 0 $.  
\end{lemma}

\begin{pf}
Define the map 
 $$
 \Psi_\e : \Ab ([1,2]^\es ) \times \Lambda  \to \R^\es \times \Lambda \, , \quad
 \Psi_\e (\dom, \lambda) := \big( \dom + \l \bar \mu + \e^2 \lambda \dom, \lambda \big) \, ,  
 $$
 which is a diffeomorphism onto its image. Thus, recalling \eqref{fraBep},  
 \be\label{meaPsiB}
| \Psi_\e ( {\frak B}_\e )| \lesssim |{\frak B}_\e |  \lesssim b_1 (\e) \, . 
 \ee
For any $ \beta \in  {\cal C}_1 $, 
let
\be\label{def:Ubeta}
U_{\beta , \e} := \Big\{ \l \in \Lambda \, : \, \frac{\beta -  \l \bar \mu}{1 + \e^2 \l} \in  \Ab ([1,2]^\es ) \Big\} \, ,
\ee
i.e. $ \l \in U_{\beta , \e} $ if and only if 
$ (\beta, \lambda ) $ is in the image
$ \Psi_\e (  \Ab ([1,2]^\es )  \times \Lambda ) $.
Thus, recalling \eqref{def:Bep}, we deduce that 
$$
\beta \in B_\e'  \, , \ \l \in U_{\beta , \e} \quad 
\Longrightarrow \quad 
(\beta, \l) \in \Psi_\e ( {\frak B}_\e ) \, . 
$$
Therefore
\be\label{inte-sepa}
\int_{ B_\e'} | U_{\beta , \e} | d \beta \leq |\Psi_\e ( {\frak B}_\e )| \stackrel{\eqref{meaPsiB}} 
\lesssim b_1(\e) \, .
\ee
Our aim is now to  justify that the measure of $ U_{\beta , \e} $ satisfies 
$|U_{\beta , \e}| \geq \sqrt{b_1(\e)}$ for all 
$\beta \in \Ab ([1,2]^\es ) + \Lambda \bar \mu$ but a subset the measure of
which vanishes as $\e \to 0$.

First note that $\Ab ([1,2]^\es) $ is a convex subset of $\R^\es$, with interior
$\Ab ((1,2)^\es)$. For $ \e > 0 $, define 
\begin{align}
V_\e & := \Big\{ x \in \R^\es \ : \  \forall \l \in \Lambda \, , \ \frac{x}{1+\e^2 \l} \in 
\Ab ((1,2)^\es) \Big\} \label{defVep} \\
&=\Big\{   x\in  \R^\es  \ : \  \Big[\frac1{1+\e^2\l_0}, \frac1{1-\e^2\l_0}\Big] \cdot x \subset \Ab ((1,2)^\es) \Big\} \subset \Ab ((1,2)^\es) \, . \nonumber
\end{align}
Each set $ V_\e  $ 
is convex and open, 
\be\label{Veei}
V_{\e'} \subset V_\e \, ,  \ \forall \, 0<\e<\e' \qquad {\rm and} \qquad 
\bigcup_{\e>0} V_\e= \Ab ((1,2)^\es) \, .
\ee
Then, for any $\beta \in  {\cal C}_1 $ 
define
$$
U'_{\b,\e}:= \Big\{ \l \in \Lambda \  : \  \b-\l \bar \mu \in V_\e \Big\}
\stackrel{\eqref{defVep}, \eqref{def:Ubeta}} \subset  U_{\b,\e} \, .
$$
Recalling  \eqref{Veei} we have 
$ U'_{\b,\e'} \subset U'_{\b,\e} $, $ \forall 0<\e<\e' $. 
At last define, for $ \d > 0 $,
\be\label{def-e-b}
D_{\e, \d}:= \big\{  \beta \in {\cal C}_1 \ : \  |U'_{\b,\e}| \geq \d \big\} \, . 
\ee
The following properties holds: 
\begin{itemize}
\item[(i)] For  $0< \e < \e'$ and $0<\d<\d' $, we have $ D_{\e',\d'} \subset D_{\e,\d}$.
\item[(ii)] Since the sets $U'_{\beta,\e}$ are open (hence of strictly positive measure if nonempty),
$$
\bigcup_{\d>0} D_{\e , \d} = \big\{ 
\b \in {\cal C}_1 \ : \ U'_{\beta ,\e} \neq \emptyset \big\} =: D_\e \, .
$$
Moreover $ D_\e =V_\e + \Lambda \bar \mu $.
\item[(iii)]  By items (i)-(ii) and \eqref{Veei} we deduce
$$ 
\bigcup_{\e>0} \Big( \bigcup_{\d >0} D_{\e,\d} \Big)=\bigcup_{\e>0} D_\e = \Ab((1,2)^\es)+ \Lambda \bar \mu = {\cal C}_1 \, . 
$$
\end{itemize}
{\bf Claim:}
\be\label{claim-ch5}
\lim_{\e \to 0_+, \d \to 0_+} \big|{\cal C}_1\backslash D_{\e,\d}\big|=\big| {\cal C}_1 \backslash( \Ab((1,2)^\es)+ \Lambda \bar \mu ) \big| = 0 \, .
\ee
The first equality follows by items (i)-(iii) above. 
To justify this last equality, let us introduce the $(\es -1)$-dimensional linear 
subspace  of $\R^\es$, $E:=\bar \mu^\bot$. Let $K$ be the 
orthogonal projection of $\Ab ([1,2]^\es )$ onto $E$; note that $K$ is a convex compact subset of 
$E$ of nonempty interior in $E$. Moreover, since  $\Ab ([1,2]^\es )$ is convex, it can be decomposed as
$$
\Ab ([1,2]^\es ) = \bigcup_{x \in K} 
\big\{ x + [\a_-(x), \a_+(x)] \bar \mu  \big\} 
$$
where the functions $\alpha_+, \alpha_- : K \to \R$ are respectively concave and convex, with 
$\a_- (x)  \leq \a_+ (x) $, for all $ x \in K $,  and  $\a_-(x) < \a_+(x)$ for $x \in {\rm int}(K)$. Hence, since $\a_\pm$ are continuous
on $ {\rm int} (K)$, 
$$
\Ab ((1,2)^\es )= {\rm int}(\Ab ([1,2]^\es ))=
\big\{x + (\a_-(x), \a_+(x)) \bar \mu \, ; \  x \in {\rm int}(K) \big\} \, ,
$$
and 
\begin{align*}
{\cal C}_1 \backslash (\Ab ((1,2)^\es )+ \Lambda  \bar \mu) =
\Big\{ x+ \l \bar \mu \ :  & \ 
\big( x \in \partial K \, , \, \l \in [\a_-(x)-\l_0,\a_+(x)+\l_0] \big) \\
& \  {\rm or} \ 
\big( x \in K \, , \, \l =\a_{\pm}(x) \pm \l_0 \big) \Big\} \, ,
\end{align*}
which gives $\big|{\cal C}_1 \backslash (\Ab ((1,2)^\es )+ \Lambda  \bar \mu )\big|=0$.

\smallskip

Setting $D'_\e := D_{\e, \sqrt{b_1(\e)}}$, where $ D_{\e, \d} $ is defined in \eqref{def-e-b},
the estimate \eqref{claim-ch5} implies 
\be\label{limed}
\lim_{\e \to 0_+} \big|{\cal C}_1 \backslash  D'_\e \big| = 0 \, . 
\ee 
Moreover, by the definition of $D'_{\e}$ and the inclusion $U'_{\b,\e} \subset U_{\b,\e} $, 
we deduce   
$$
\forall \b \in D'_{\e} \, , \ \   | U_{\beta,\e} | \geq \sqrt{b_1(\e)} \, ,  
$$
and therefore 
\be\label{Be-dentro}
| B_\e' \cap D'_{\e} | \sqrt{b_1(\e)} \leq 
\int_{ B_\e' \cap D'_\e} | U_\beta | d \beta  
\stackrel{\eqref{inte-sepa}} 
\lesssim b_1(\e) \, . 
\ee
Finally, since  $ B_\e' \subset (B_\e' \cap D'_{\e}) \cup ({\cal C}_1 \backslash D'_{\e}) $, we deduce, by  \eqref{Be-dentro} and \eqref{limed}, that $ \lim_{\e \to 0 } |B_\e' | = 0 $.
Lemma \ref{mesbe} is proved. 
\end{pf}

Now, recalling \eqref{def:Bep} and \eqref{expr:omega},  we have 
\be \label{includedOm}
\bar \mu + \bigcup_{\e>0} \e^2 ({\cal C}_1 \backslash B_\e') \subset \Omega 
\ee
where $\Omega$ is the set of the frequency vectors $\om$ of the quasi-periodic solutions
of \eqref{NLW1} provided by Theorem \ref{thm:main}.
Note that,  
 by \eqref{Lambda-unif}, \eqref{def omep}, \eqref{def:admissible}, we get  that 
$$
  \zeta +
 \l \bar \mu + \l  \e^2 \zeta
 \in   \Ab \big( \big[\frac12,4 \big]^\es \big) \, , \quad 
 \forall \zeta \in \Ab ([1,2]^\es )  \, , \ \forall \l \in \Lambda \, , 
$$
and therefore ${\cal C}_1 = \Ab ([1,2]^\es ) + \Lambda \bar \mu $ 
does not contain $ 0 $. Moreover 
${\cal C}_1 $ is a compact convex subset of $\R^\es$, with nonempty interior, which implies that 
\be \label{propcone}
\big\{ x \in \R_+ {\cal C}_1 \ : \ \R_+ x \cap {\cal C}_1 \ \hbox{is a singleton} \, \big\} \subset \partial(\R_+ {\cal C}_1) \, .
\ee
Thus, given  $ y \in \R_+ {\cal C}_1 $, the measure
$ \big| \{ r > 0 \, : \,  y / r  \in {\cal C}_1 \} \big| > 0 $
except for $ y \in \partial( \R_+ {\cal C}_1) $, which is of zero measure. 
Using \eqref{propcone} and Lemma \ref{mesbe}, we can obtain
\be \label{dens1}
\lim_{r \to 0_+} \frac{\Big| \Big( \bigcup_{\e>0} \e^2 ({\cal C}_1 \backslash B_\e') \Big)
\cap B(0,r)  \Big|}{\Big| \R_+ {\cal C}_1 \cap B(0,r)  \Big|} =1 \, .
\ee
We omit the details. Recalling \eqref{includedOm}, 
\eqref{dens1} implies \eqref{density-Leb}.

\section{Shifted tangential frequencies up to $ O(\e^4 ) $}\label{sec:shifted-tan}

In this section\index{Shifted tangential frequencies} we evaluate the average of  the first component in \eqref{operatorF0}: 
\be\label{bif-eq}
\om - \bar \mu - \e^2 \langle \partial_y R (\vphi, 0, 0, \xi) \rangle   
\ee
 where
\be\label{bif-eqa}
\langle f \rangle := \frac{1}{(2 \pi)^\es} \int_{\T^\es} f (\vphi) \, d \vphi \, .
\ee
Evaluating \eqref{gradyR} at $ (\theta, y, Q) = (\vphi, 0,0) $, 
we get, inserting the expression of $ g (\e, x, u) $ in \eqref{nonlinearity:gep1}, that for each $ m =1, \ldots, \es $, 
\begin{align} 
\langle \partial_{y_m} R(\vphi, 0, 0, \xi) \rangle 
& =  \frac{1}{(2\pi)^\es}\int_{\T^\es } \int_{\T^d} g(\e, x, v(\vphi, 0, \xi) ) 
\frac{\mu_m^{-1/2}}{\sqrt{2 \xi_m }} \cos \vphi_m \, \Psi_m (x)\, dx \, d \vphi  \nonumber \\ 
& = {\mathtt r}_{m, 3} + \e {\mathtt r}_{m, 4} + \e^2 {\mathtt r}_{m, 5}   \label{shifted0}
\end{align}
where
\begin{align}
{\mathtt r}_{m, 3} & :=   \frac{1}{(2\pi)^\es}\int_{\T^\es } \int_{\T^d}  a(x)  
\big( v(\vphi, 0, \xi)  \big)^3  
 \frac{\mu_m^{-1/2}}{\sqrt{2 \xi_m }} \cos \vphi_m \, \Psi_m (x)\, dx \label{Rm3} \\
 {\mathtt r}_{m, 4} & :=   \frac{1}{(2\pi)^\es}\int_{\T^\es }  \int_{\T^d}  a_4 (x)  
\big( v(\vphi, 0, \xi)  \big)^4  
 \frac{\mu_m^{-1/2}}{\sqrt{2 \xi_m }} \cos \vphi_m \,  \Psi_m (x)\, dx \label{termine:succ} \\
 {\mathtt r}_{m, 5} & :=  
  \frac{1}{(2\pi)^\es}\int_{\T^\es }  \int_{\T^d}  {\mathfrak r} ( \e, x, v(\vphi, 0, \xi) ) \frac{\mu_m^{-1/2}}{\sqrt{2 \xi_m }} \cos \vphi_m \,  \Psi_m (x)\, dx  \, . \label{termine:r2-non-nullo}
\end{align}

We now compute the  terms in \eqref{shifted0}. 

\begin{lemma}\label{lem:rm3}
${\mathtt r}_{m,3} $ in \eqref{Rm3} is 
\be\label{eccola-first3}
{\mathtt r}_{m,3}  = [\Ab \xi]_m  
\ee 
where  $ \Ab  := ( \Ab^j_m)_{ j, m \in \mathbb S} $ is  the symmetric {\it twist} matrix\index{Twist matrix}  defined in \eqref{def:AB}.
\end{lemma}

\begin{pf}
Using  \eqref{defv}, we expand the integral  \eqref{Rm3} as 
 \begin{align}
{\mathtt r}_{m,3} & = 
 \sum_{j_1, j_2, j_3 \in {\mathbb S}}    \mu_{j_1}^{-1/2} \mu_{j_2}^{-1/2}  \mu_{j_3}^{-1/2} \mu_m^{-1/2} 
  \sqrt{\xi_{j_1}  \xi_{j_2} \xi_{j_3} } \frac{2}{\sqrt{ \xi_m }}  \nonumber \\
 & \quad \times 
 \frac{1}{(2\pi)^\es}\int_{\T^\es }   \cos \vphi_{j_1}\cos \vphi_{j_2}\cos \vphi_{j_3}  \cos \vphi_m d \vphi 
  \nonumber  \\ 
& \quad  \times \int_{\T^d}  a(x)   \Psi_{j_1}(x)   \Psi_{j_2}(x)   \Psi_{j_3}  (x) \Psi_m (x)\, dx \, . \label{passints}
\end{align}
The integral
$ \int_{\T^\es }   \cos \vphi_{j_1}\cos \vphi_{j_2}\cos \vphi_{j_3}  \cos \vphi_m d \vphi  $
does not vanish only if 
\begin{enumerate}
\item
$ j_1 = j_2 = j_3 = m $, 
\item
$ j_1 = j_2 \neq  j_3 = m $ and permutation of the indices ($ 3 $ times).
\end{enumerate}
Hence,  by \eqref{passints}, 
\begin{align*}
{\mathtt r}_{3,m} & =  2  \mu_{m}^{-2}   \xi_{m}    \frac{1}{(2 \pi)^\es}
 \int_{\T^\es }   \cos^4 \vphi_m d \vphi 
   \int_{\T^d}  a(x)   \Psi_{m}^4(x) \, dx  \nonumber  \\
& +  3  \sum_{j_1 \neq m}  \mu_{j_1}^{-1}  \mu_m^{-1} \xi_{j_1}  2  
  \frac{1}{(2 \pi)^\es}\int_{\T^\es }   \cos^2 \vphi_{j_1} \cos^2 \vphi_m d \vphi 
  \int_{\T^d}  a(x)   \Psi^2_{j_1}(x)   \Psi^2_m (x)\, dx \nonumber \\
 & =   \frac34   \mu_{m}^{-2}   T_m^m  \xi_{m}    +  \frac32  \sum_{j \neq m}  \mu_{j}^{-1}  \mu_m^{-1}  T^j_m \xi_{j}  
\end{align*}
having set 
$$ 
T^j_m :=  \int_{\T^d}  a(x)   \Psi^2_{j}(x)  \Psi^2_m (x)\, dx \, , \  j, m \in \mathbb S \, , 
$$
and noting that 
\begin{align*}
& \frac{2}{(2 \pi)^\es} \int_{\T^\es}   \cos^4 \vphi_m d \vphi =  
 \frac{1}{ \pi}  \int_\T \cos^4 \vphi_m d \vphi_m = \frac34   \\
&   \frac{6}{(2 \pi)^\es} \int_{\T^\es }   \cos^2 \vphi_{j_1} \cos^2 \vphi_m d \vphi =
 \frac{6}{(2 \pi)^2}  \Big( \int_{\T} \cos^2 \theta \, d \theta \Big)^2    = \frac32   \, . \nonumber 
\end{align*}
Recalling the definition of  the twist matrix 
$ \Ab  := ( \Ab^j_m)_{ j, m \in \mathbb S} $  in \eqref{def:AB}-\eqref{def G}
we deduce \eqref{eccola-first3}. 
\end{pf}

\begin{lemma}\label{lem:r1=0}
For all $ m \in {\mathbb S}$, each $ {\mathtt r}_{m,4} $ in \eqref{termine:succ} is $ {\mathtt r}_{m,4} = 0 $.  
\end{lemma}

\begin{pf}
Since the function $ v $ defined in \eqref{defv} satisfies the symmetry 
\be\label{simmetry:v}
v(\vphi + \vec \pi, 0, \xi ) = - v(\vphi, 0, \xi ) \, , \quad \vec \pi := (\pi, \ldots, \pi) \in \R^\es \, , 
\ee
the function 
$ g(\vphi) := (v(\vphi, 0, \xi ))^4 \cos ( \vphi_m) $ satisfies $ g(\vphi + \vec \pi ) = - g(\vphi) $
and therefore its integral
$$
\int_{\T^\es} g ( \vphi) d \vphi = \int_{\T^\es} g ( \vphi + \vec \pi) d \vphi = - \int_{\T^\es} g ( \vphi ) d \vphi 
$$
is equal to zero. Hence $ {\mathtt r}_{m,4} = 0 $. 
\end{pf}

By \eqref{shifted0}, Lemmata \ref{lem:rm3} and \ref{lem:r1=0}  we deduce that 
\be\label{eccola-first}
\langle \pa_y R (\vphi,0,0,\xi) \rangle = \Ab \xi +  \e^2 {\mathtt r}_5 (\e, \xi)    
\ee 
where   $ \Ab $ is  the symmetric {\it twist} matrix  defined in \eqref{def:AB}, and
$ {\mathtt r}_5 (\e, \xi)  := ({\mathtt r}_{5,m})_{m =1, \ldots, \es } \in \R^\es $.
As a consequence, the term \eqref{bif-eq} is $ O(\e^4) $, more precisely, 
since  $ \om = (1+\e^2 \l) \bar \om_\e  $, it results
\begin{align}\label{1comp-pro}
 \omega - \bar \mu - \e^2 \langle \partial_y R (\vphi, 0 , 0, \xi)\rangle  
  & = (1+\e^2 \l) \bar \om_\e  - \bar \mu - \e^2 \Ab \xi  - \e^4 {\mathtt r}_5 (\e, \xi ) \nonumber \\
&  \stackrel{\eqref{xil}} = - \e^4 {\mathtt r}_5 (\e, \xi ) \, .
\end{align}

\section{First approximate solution}\label{NM:first-step}

We now define the first approximate\index{First approximate solution} torus embedding solution
\be\label{def:i1}
i_1 (\vphi) = 
\big( \theta_1 (\ph), y_1(\ph), Q_1(\ph) , P_1 (\ph) \big)  \, , \qquad \theta_1 (\vphi ) = \vphi + \vartheta_1(\ph) \, , 
\ee
in such a way that $ {\cal F} ( i_1 ) = O(\e^4 ) $. 
Given a function $ f : \T^\es \to \R^\es $, we denote by $ \langle f \rangle \in \R^\es $
its average with respect to $ \vphi $ (as in \eqref{bif-eqa})  and by $ [f ]( \vphi)  $ its zero mean part, so that 
\be\label{def:quadra}
f (\vphi ) = \langle f \rangle + [f ](\vphi) \, . 
\ee

\begin{lemma}\label{NM:step1} {\bf (First approximate solution)}
Let $ \om = (1+\e^2 \l) \bar \om_\e  $ 
with $ \bar \om_\e   = \bar \mu + \e^2 \dom $ as in  \eqref{def omep} and 
 define $ \xi := \xi (\lambda) $ as in  \eqref{xil}.  
Then there exists a unique solution $ i_1 $, with the form  in \eqref{def:i1}, 
with average $ \langle y_1 \rangle = 0 $, independent of $ \l \in \Lambda $, 
of the system 
\be\label{HS00exp-1}
\begin{cases}
\bar \mu \cdot \partial_\vphi \vartheta_1  - \e^2 \big[ (\partial_y R) (\vphi,0  , 0, \xi) \big] = 0  \cr 
\bar \mu \cdot \partial_\vphi  y_1  +  \e^2 (\partial_\teta R) (\vphi, 0 , 0, \xi) = 0 \cr
\big(\bar \mu \cdot \partial_\vphi -  J D_V \big) (Q_1, P_1) 
- \e^2 \big( 0, (\nabla_Q R) (\vphi, 0 , 0, \xi) \big) = 0 
\end{cases}
\ee
satisfying, for all $ s \geq s_0 $, 
\be\label{estimate-i1}
\| i_1 - (\vphi, 0, 0) \|_{\Lip,s} = \| i_1 - (\vphi, 0, 0) \|_{s} \leq C(s) \e^2 \, . 
\ee
It results
\be\label{approx1} 
\| {\cal F}(i_1) \|_{\Lip,s} \leq C(s) \e^4 \, , \quad \forall s \geq s_0 \, . 
\ee
\end{lemma}

\begin{pf}
{\sc Solution of the first equation in \eqref{HS00exp-1}.} 
Since $ \bar \mu $ is a Diophantine vector by  \eqref{diop},
we solve the first equation in \eqref{HS00exp-1}, finding 
\be\label{def:theta1}
 \vartheta_1  =  \e^2 (\bar \mu \cdot \partial_\vphi)^{-1} \big[ (\partial_y R) (\vphi, 0 , 0, \xi) \big] 
 \ee
where $ (\bar \mu \cdot \partial_\vphi)^{-1}  $ is defined as in  \eqref{op-inv-KAM} (with $ \bar \mu $ instead of $ \om $). 
\\[1mm]
{\sc Solution of the second equation in \eqref{HS00exp-1}.} 
Since 
$$ 
(\partial_\teta R) (\vphi,  0 , 0, \xi) = \partial_\vphi \big( R (\vphi, 0 , 0, \xi)  \big) \, , 
$$  
it has zero average in $ \vphi $. Then, since $ \bar \mu $ is  Diophantine  by  \eqref{diop}, 
the second equation in \eqref{HS00exp-1} admits the unique solution 
with zero average
\be\label{def:y1}
 y_1  =  -  \e^2 (\bar \mu \cdot \partial_\vphi)^{-1}  \big[ (\partial_\teta R) (\vphi, 0 , 0, \xi) \big]   \, .
\ee
{\sc Solution of the third equation in \eqref{HS00exp-1}.}
The operator $ \bar \mu \cdot \partial_\vphi -  J D_V  $ is represented, 
in
 the basis $ \{ e^{\ii \ell \cdot \vphi } (\Psi_j (x), 0), e^{\ii \ell \cdot \vphi } (0, \Psi_j (x)) \}_{j \in \N }$
 (see  
\eqref{matrix-A-space-time-psi}) by the diagonal matrix
$$
{\rm Diag}_{\ell \in \Z^\nu, j \in \N} \begin{pmatrix}
 \ii \bar \mu \cdot \ell & - \mu_j  \\
 \mu_j &  \ii \bar \mu \cdot \ell    \\
\end{pmatrix} \, , 
$$
and therefore, by the unperturbed first Melnikov condition \eqref{1Mel}, it is invertible.
Moreover, arguing as in the end of Lemma \ref{homdiag-2-1}, it satisfies the estimate 
\be\label{stima-Hs-in}
\| \big(\bar \mu \cdot \partial_\vphi -  J D_V \big)^{-1} h \|_s \leq C(s) \| h \|_{s+ \tau_0} \, . 
\ee
Then the third equation in \eqref{HS00exp-1} admits the unique solution
\be\label{def:z1}
(Q_1, P_1) = \e^2 \big(\bar \mu \cdot \partial_\vphi -  J D_V \big)^{-1} \big(0, (\nabla_Q R) (\vphi, 0 , 0, \xi) \big) \, . 
\ee
 By the definition of 
$ i_1 (\vphi) = (\theta_1 (\vphi), y_1 (\vphi), z_1 (\vphi)) $ where 
$ z_1 (\vphi) := (Q_1 (\vphi) ,P_1 (\vphi)  ) $,  in 
\eqref{def:theta1}, \eqref{def:y1},  \eqref{def:z1}, the estimate \eqref{estimate-i1} follows from  \eqref{stima-Hs-in}  and 
the Diophantine condition \eqref{diop}. 
 Finally, 
 comparing \eqref{operatorF} with 
system \eqref{HS00exp-1}, we have 
$$
\begin{aligned}
& {\cal F}(i_1) = \\
&  
\left(
\begin{array}{c}
\!\! \om - \bar \mu - \e^2 \langle (\partial_y R) ( i_1(\vphi), \xi   )\rangle +  
(\om - \bar \mu) \cdot \pa_\vphi \vartheta_1 (\vphi)  - 
\e^2 \big( [\partial_y R ( i_1(\vphi), \xi   )] -  [\partial_y R (\vphi, 0, 0, \xi )] \big) \!\!  \\
(\om - \bar \mu) \cdot \pa_\vphi  y_1 (\vphi)  +  \e^2 
\big( (\partial_\teta R) ( i_1(\vphi), \xi  ) -   (\partial_\teta R) (\vphi, 0, 0, \xi  ) \big)  \\
(\om - \bar \mu) \cdot \pa_\vphi z_1 (\vphi) - \e^2 \Big( \big( 0, (\nabla_Q R) (i_1(\vphi), \xi) \big) - 
 \big( 0, (\nabla_Q R) (\vphi,0,0, \xi) \Big)
\end{array}
\right) 
\end{aligned}
$$
and \eqref{approx1} follows
using the estimate \eqref{estimate-i1}, the fact that 
$ \om = (1+\e^2 \l) \bar \om_\e  $ with $\bar \om_\e  = \bar \mu + \e^2 \zeta $, 
and \eqref{1comp-pro}.
 \end{pf}

The successive approximate solutions $ i_n $, $ n \geq 2 $,  
of the functional equation $ {\cal F}(i) = 0 $, 
are defined through a Nash-Moser iterative scheme.
The main point to define $ i_{n+1} $ is the construction of an approximate right inverse of the linearized operators 
$ d_{i} {\cal F}(i_n ) $ 
at the approximate torus $ i_n $, that we obtain in the next Chapters 
\ref{sezione almost approximate inverse}-\ref{sec:proof.Almost-inv}.

\chapter{Linearized operator at an approximate solution}\label{sezione almost approximate inverse}

In order to implement a convergent Nash-Moser scheme 
(Chapter \ref{sec:NM}) that leads to a solution of 
$ \mF(\lambda, i) = 0 $ where $ \mF(\lambda, i) $ is 
the nonlinear operator defined in \eqref{operatorF}, 
the key step is to 
prove the existence of an approximate right inverse of the linearized operator
$ d_{i} {\cal F}(\l; \ui) $ in \eqref{diF}. 
The first step  is  Proposition \ref{Prop:sec6} 
where we introduce suitable symplectic coordinates which reduce the problem to the  
search of an approximate inverse of the operator $ {\cal L}_\omega $ in 
\eqref{Lomega def} acting in the normal components only. 
This will be studied in Chapters \ref{sec:6n}-\ref{sec:proof.Almost-inv}. 

\section{Symplectic approximate decoupling}

We linearize $ \mF(\lambda, i)$   at an arbitrary  torus 
\be\label{torus-linearized}
 \ui (\vphi) = (\uth (\vphi) , \uy (\vphi), \uz (\vphi) ) \, , 
\ee
obtaining 
\be\label{diF}
d_{i} {\cal F}(\l; \ui)[\widehat \imath  ] =
\Dom \widehat \imath - d_i X_{K} ( \ui (\vphi) ) [\widehat \imath ]  \, .
\ee
We denote by 
\be\label{per-comp}
\uF (\vphi)  := \ui  (\vphi) - (\vphi,0,0) := 
( {\uTh} (\ph), \uy (\ph), \uz (\ph))\,, \quad \uTh (\ph) := \uth (\vphi) - \vphi \, ,  
\ee
the periodic component of the torus $ \vphi \mapsto \ui  (\vphi) $ with norm as in 
\eqref{weighted-Sobo-three}.

We assume the following 
condition for $ \ui $ which is satisfied by any approximate solution obtained along 
the Nash-Moser iteration performed in Chapter \ref{sec:NM} (see precisely \eqref{NM:indu-line2-trunc}):  
\begin{itemize}
\item 
The map $ \l \mapsto  \uF (\l) $
is Lipschitz  with respect  to  $ \l \in \Lambda_{\uF} \subset \Lambda $,  and   
\begin{equation}\label{ansatz 0}
\|  \uF  \|_{\Lip, {s_1 +2}} \leq C(s_1) \e^2 \, ,    \quad \| \uF  \|_{\Lip, s_2}  \leq \e \, .   
\end{equation}
\end{itemize}
We  implement the general strategy proposed in \cite{BBField},
used also in \cite{BBM-auto}, \cite{BM16},
where, instead of inverting $ d_i {\cal F}(\l; \ui ) $ 
(where all the $ (\theta, y, z ) $  components are coupled, see \eqref{operatorF}) 
we invert the linear operator $ {\mathbb D } (\ui) $ in  \eqref{operatore inverso approssimato},
which has a  triangular form. 
The operator $ {\mathbb D } (\ui) $ is found by a natural geometrical construction.
We define  the ``error function''\index{Error function}
\begin{equation} \label{def Zetone}
Z(\vphi) :=  (Z_1, Z_2, Z_3) (\vphi) := {\cal F}(\lambda; \ui) (\vphi) =
\om \cdot \pa_\vphi \ui(\vphi) - X_{K}(\ui(\vphi)) \, . 
\end{equation}
Notice that, if $ Z = 0 $ then the torus $ \ui $ is invariant for $ X_{K} $; 
in general, we say that $ \ui $ is ``{\it approximately invariant}''\index{Approximately invariant torus}, up to order $O(Z)$. 
Given  $ \ui (\vphi)  $ satisfying \eqref{ansatz 0}
we first construct an isotropic torus $i_\delta(\vphi) $ 
which is close to $ \ui $, see  \eqref{stima y - y delta} and \eqref{2015-2}. By  \eqref{stima toro modificato}, 
$ {\cal F}(i_\delta ) $ is also $ O(Z) $. 
Since the torus $i_\d$ is isotropic, the diffeomorphism $ (\phi, \zeta, w) \mapsto G_\delta(\phi, \zeta, w)$ 
defined in \eqref{trasformazione modificata simplettica} is symplectic.  
In these coordinates, the torus $i_\delta$ reads $(\phi, 0, 0)$,  
and the transformed Hamiltonian system becomes \eqref{sistema dopo trasformazione inverso approssimato}, where, by  \eqref{apprVF}  
the terms $ \pa_\phi {\mathtt K}_{00}, {\mathtt K}_{10} - \omega, {\mathtt K}_{01}$ are $O(Z)$. 
Neglecting such terms in the linearized operator \eqref{lin idelta} at $ (\phi,0,0)$, 
we obtain the linear operator ${\mathbb D}(\ui) $ in \eqref{operatore inverso approssimato}.

The main result of this section is the following Proposition. 

\begin{proposition}\label{Prop:sec6}
Let  $ \ui ( \vphi ) $ be a torus of the form \eqref{torus-linearized}, defined for all  $ \l \in \Lambda_\uF $, 
 satisfying \eqref{ansatz 0}. Then 
\begin{itemize}
\item {\bf (isotropic torus)}
there is  an isotropic torus $ i_\delta (\vphi) = (\uth (\vphi), y_\delta (\vphi), \uz(\vphi) ) $ satisfying, 
for some $ \bt := \bt(\es,\t_1) > 0 $, 
\begin{align} 
\label{stima y - y delta}
\| y_\delta - \uy \|_{\Lip, s} & \lesssim_s \| Z \|_{\Lip, s + \bt} + \| Z \|_{\Lip, s_0 + \bt} \|  \uF \|_{\Lip, s + \bt} 
\\
\label{stima toro modificato}
\| {\cal F}(i_\delta ) \|_{\Lip, s} & \lesssim_s  \| Z \|_{\Lip, s + \bt}  +  \| Z \|_{\Lip, s_0 + \bt} \| \uF \|_{\Lip, s + \bt} \, . 
\end{align}
Given another $ {\underline i}' $ satisfying \eqref{ansatz 0} we have  
\be\label{derivata i delta}
\| i_\d ( \ui ) -  i_\d ( \ui' ) \|_{s_1} \lesssim_{s_1} \| \uF - \uF' \|_{s_1+1} \, .  
\ee 
\item 
 {\bf (symplectic diffeomorphism)}
the change of variable $ G_\delta : (\phi, \ac , w) \to (\theta, y, z) $ 
of the phase space $\T^\es \times \R^\es \times H_{\St}^\bot$ defined by
\begin{equation}\label{trasformazione modificata simplettica}
\begin{pmatrix}
\theta \\
y \\
z
\end{pmatrix} := G_\delta \begin{pmatrix}
\phi \\
\ac \\
w
\end{pmatrix} := 
\begin{pmatrix}
\!\!\!\!\!\!\!\!\!\!\!\!\!\!\!\!\!\!\!\!\!\!\!\!\!\!\!\!\!\!\!\!
\!\!\!\!\!\!\!\!\!\!\!\!\!\!\!\!\!\!\!\!\!\!\!\!\!\!\!\!\!\!\!\!
\!\!\!\!\!\!\!\!\!\!\!\!\!\!\!\!\!\!\!\!\!\!\!\!\!\!\! \uth (\phi) \\
y_\delta (\phi) + [\pa_\phi \uth (\phi)]^{-\top} \ac - \big[ (\pa_\teta \tilde{\uz}) (\uth (\phi)) \big]^\top J w \\
\!\!\!\!\!\!\!\!\!\!\!\!\!\!\!\!\!\!\!\!\!\!\!\!\!\!\!\!\!
\!\!\!\!\!\!\!\!\!\!\!\!\!\!\!\!\!\!\!\!\!\!\!\!\!\!\!
\!\!\!\!\!\!\!\!\!\!\!\!\!\!\!\!\!\!\!\!\!\!\!\! \uz(\phi) + w
\end{pmatrix} 
\end{equation}
where $ \tilde{\uz} (\theta) := \uz (\uth^{-1} (\theta))$, is symplectic. 

In the new coordinates $ (\phi, \ac, w )$,  the isotropic torus $ i_\delta $ is the trivial embedded torus
$ (\vphi , 0, 0 ) $, i.e. 
\be\label{toro-new-coordinates}
i_\delta (\vphi) = G_\delta (\vphi , 0, 0 ) \, .
\ee
The linearized diffeomorphism $ DG_\delta(\vphi,0,0) $ satisfies, for all $ s \geq s_0 $, 
\be\label{DG delta}
\begin{aligned} 
  \|DG_\delta(\vphi,0,0) [\, \widehat \imath \, ] \|_{\Lip,s} + 
    \big \| \big( DG_\delta(\vphi,0,0) \big)^{-1} [ \, \widehat \imath \, ] \big\|_{\Lip, s} 
   &  \lesssim_s \| \, \widehat \imath \, \|_{\Lip, s} \\ 
& \ \,    +  \| \uF \|_{\Lip, s + 2}  \| \, \widehat \imath \, \|_{\Lip, s_0} 
\end{aligned}
\ee
and
\be
\begin{aligned}
\| D^2 G_\delta(\vphi,0,0)[\widehat \imath_1, \widehat \imath_2] \|_{\Lip,s}
& \lesssim_s  \| \, \widehat \imath_1\|_{\Lip,s}  \| \, \widehat \imath_2 \|_{\Lip, s_0} 
+ \| \, \widehat \imath_1\|_{\Lip, s_0}  \| \, \widehat \imath_2 \|_{\Lip, s} \\
& \quad +  \| \uF \|_{\Lip, s + 3} \| \, \widehat \imath_1 \|_{\Lip, s_0} 
 \| \, \widehat \imath_2\|_{\Lip, s_0}\, .  \label{DG2 delta} 
\end{aligned}
\ee
\item {\bf (Transformed Hamiltonian)}
Under the symplectic change of variables $ G_\d $, the Hamiltonian 
vector field $ X_{K} $ (the Hamiltonian $  K $ is defined in \eqref{primalinea}) transforms into 
\be\label{new-Hamilt-K}
X_{\mathtt K} = (D G_\d)^{-1} X_{K} \circ G_\d \qquad {\rm where} \qquad {\mathtt K} := K \circ G_\d  \, .
\ee
The Hamiltonian $ {\mathtt K} $ is reversible, i.e. $ {\mathtt K} \circ \tilde S = {\mathtt K} $.  
The $ 2 $-jets of the Taylor expansion of the Hamiltonian $ {\mathtt K} $ at the trivial torus $ (\phi , 0, 0 ) $, 
\begin{align} 
{\mathtt K} (\phi, \ac , w)
& =  {\mathtt K}_{00}(\phi ) + {\mathtt K}_{10}(\phi ) \cdot \ac + ({\mathtt K}_{0 1}(\phi ), w)_{L^2(\T_x)} + 
\frac12 {\mathtt K}_{2 0}(\phi) \ac \cdot \ac
\nonumber \\ & 
\quad +  \big( {\mathtt K}_{11}(\phi) \ac , w \big)_{L^2(\T_x)} 
+ \frac12 \big( {\mathtt K}_{02}(\phi) w , w \big)_{L^2(\T_x)} + {\mathtt K}_{\geq 3}(\phi, \ac, w)  
\label{KHG}
\end{align}
where $ {\mathtt K}_{\geq 3} $ collects the terms at least cubic in the variables $ (\ac , w )$, 
satisfy the following properties: 

i) The vector field 
\begin{align}\label{apprVF}
X_{\mathtt K} (\phi, 0, 0) & = 
\left(
\begin{array}{c}
\!\! {\mathtt K}_{10}(\phi) \!\!   \\
\!\!  -  \partial_{\phi} {\mathtt K}_{00}(\phi )  \!\! \\
\!\! J {\mathtt K}_{01} (\phi ) \!\!
 \end{array}
\right)  = 
\left(
\begin{array}{c}
\!\! \om  \!\! \\
\!\! 0 \!\! \\
\!\! 0 \!\!
 \end{array}
\right)
 - \big( DG_\d (\phi, 0, 0) \big)^{-1} Z_\d (\phi)  
\end{align}
where $ Z_\d (\phi)  := {\cal F}(i_\d)(\phi) $. 
The functions
$ {\mathtt K}_{00} : \T^\es \to \R $, $ {\mathtt K}_{0 1} : \T^\es \to \R^\es $ and 
$  {\mathtt K}_{0 1} : \T^\es \to H_{\mathbb S}^\bot $, that we regard as an element of  
$ {\mathcal H}^s (\T^\es \times \T^d, \R^2 ) $,   satisfy the estimate
\be\label{K 00 10 01}
\begin{aligned}
\|  \partial_\phi {\mathtt K}_{00} \|_{\Lip, s} 
+ \|  {\mathtt K}_{10} - \om  \|_{\Lip, s} +  \| {\mathtt K}_{0 1} \|_{\Lip,s}
& \lesssim_s  \| Z \|_{\Lip, s + \bt} \\ 
& \ \, +  \| Z \|_{\Lip, s_0 + \bt} \| \uF \|_{\Lip, s + \bt} \, .
\end{aligned}
\ee
ii) The average 
$ \langle {\mathtt K}_{20} \rangle    := (2 \pi)^{-\es} \int_{\T^\es}  {\mathtt K}_{20} (\phi) d \phi $ satisfies 
\begin{align}\label{stime coefficienti K 20 11 bassa}
& \| \langle {\mathtt K}_{20}  \rangle - \e^{2} \Ab \|_{\Lip}  \lesssim  \e^{4} 
\end{align}
 where $ \Ab  $ is the  twist  matrix  in  \eqref{def:AB}, and 
\begin{align}
& \| {\mathtt K}_{20} \zeta \|_{\Lip,s}  \lesssim_s  
\e^{2} \big( \| \ac \|_{\Lip,s} + \| \uF \|_{\Lip, s + \bt} \| \ac \|_{\Lip, s_0} \big) 
\label{tame:K20-appl} \\ 
\label{stime coefficienti K 11 alta}  
& \| {\mathtt K}_{11} \ac \|_{\Lip, s}
\lesssim_s \e^{2} \big( \| \ac \|_{\Lip,s} 
+ \| \uF \|_{\Lip, s + \bt} \| \ac \|_{\Lip, s_0} \big) 
\\
\label{stime coefficienti K 11 alta trasposto}  
& \| {\mathtt K}_{11}^\top w \|_{\Lip, s}
 \lesssim_s \e^{2} \big(  \| w \|_{\Lip, s} 
+ \| \uF\|_{\Lip, s + \bt} \| w \|_{\Lip, s_0} \big) \, . 
\end{align}
\item  {\bf (Linearized operator in the new coordinates)}
The linearized 
operator  
$$
 \Dom - d_{(\phi, \zeta, w)} X_{{\mathtt K}}( \phi,0,0 ) 
 $$ 
 is 
\begin{equation}\label{lin idelta}
\begin{pmatrix}
\widehat \phi  \\
\widehat \ac    \\ 
\widehat w 
\end{pmatrix} \mapsto
\begin{pmatrix}
\!\!\!\!\!\!\!\!\!\! \!\!\!\!\!\!\!\!\!\!\!\!\!\!\!\!\! \Dom \widehat \phi - \partial_\phi {\mathtt K}_{10}(\phi)[\widehat \phi \, ]  - 
{\mathtt K}_{2 0}(\phi)\widehat \ac - {\mathtt K}_{11}^\top (\phi) \widehat w \\
 \Dom  \widehat \ac + \partial_{\phi\phi} {\mathtt K}_{00}(\phi)[\widehat \phi] + 
[\partial_\phi {\mathtt K}_{10}(\phi)]^\top \widehat \ac + 
[\partial_\phi  {\mathtt K}_{01}(\phi)]^\top \widehat w   \\ 
\!\!\!\!\!\!\!\!\!\! \!\!\!\! \!\! \!\!\Dom  \widehat w - J 
\{ \partial_\phi {\mathtt K}_{01}(\phi)[\widehat \phi]  + {\mathtt K}_{11}(\phi) \widehat \ac + 
{\mathtt K}_{02}(\phi) \widehat w \}
\end{pmatrix} \! .  \hspace{-5pt}
\end{equation}
\end{itemize}
\end{proposition}
In order to find an approximate inverse of the linear operator in \eqref{lin idelta} 
 it is sufficient to invert  the operator
\begin{equation}\label{operatore inverso approssimato} 
{\mathbb D}   \begin{pmatrix}
\widehat \phi \\
\widehat \ac   \\
\widehat w
\end{pmatrix} := 
{\mathbb D}(\ui) 
  \begin{pmatrix}
\widehat \phi \\
\widehat \ac   \\
\widehat w
\end{pmatrix}
 := 
 \begin{pmatrix}
 \!\!\!\!\!\! \!\!\Dom \widehat \phi  - 
{\mathtt K}_{20}(\phi) \widehat \ac  - {\mathtt K}_{11}^\top(\phi) \widehat w\\
\!\! \!\!\!\!\!\!\!\!\!\!\!\!\!\!\!\!\!\!\!\!\!\!\!\!\!\!\!\!\!\!\!\!\!\!\!\!\! \!\!\!\!\!\!\!\!\!\!\!\!\!\!\!\!\!\!\!\!\!\!\!\! \Dom  \widehat \ac  \\
\Dom  \widehat w  - J {\mathtt K}_{02}(\phi) \widehat w    -J {\mathtt K}_{11}(\phi)\widehat \ac  
\end{pmatrix}
\end{equation}
which is obtained by neglecting in \eqref{lin idelta} 
the terms $ \partial_\phi {\mathtt K}_{10}  $, $ \partial_{\phi \phi} {\mathtt K}_{00} $, 
$ \partial_\phi {\mathtt K}_{00} $, $ \partial_\phi {\mathtt K}_{01} $,
which vanish if $ Z = 0 $ by \eqref{K 00 10 01}. 
The linear operator $ {\mathbb D } (\ui) $ 
can be inverted in a ``triangular" way. 
Indeed the second component in  \eqref{operatore inverso approssimato} for the action variable is decoupled
from the others. Then  one inverts the operator in the 
third component, 
 i.e. the operator 
\be\label{Lomega def}
{\cal L}_\omega := {\cal L}_\omega (\ui) := \Pi_{\St}^\bot \big(\Dom   - J {\mathtt K}_{02}(\phi) \big)_{|{H_{\St}^\bot}}  \, , 
\ee
and finally the first one.
The invertibility properties of $ {\cal L}_\om $ will be obtained in Proposition \ref{prop:inv-ap-vero} using the results 
of Chapters \ref{sec:6n}-\ref{sec:proof.Almost-inv}. 
We now provide the explicit  expression of ${\cal L}_\omega (\ui) $.

\begin{lemma} {\bf (Linearized operator in the normal directions)} \label{thm:Lin+FBR}
The linear operator  $ {\mathtt K}_{02}(\phi) := {\mathtt K}_{02}(\ui; \phi)  $  has the  form
\be\label{lin:normal-f-ge}
{\mathtt K}_{02}(\phi) =  D_V + \e^2 {\mathtt B}(\phi)  + {\mathtt r}_\e (\phi)  \, , 
\ee
where $ {\mathtt B} := {\mathtt B}(\e, \lambda) $ is 
the self-adjoint operator 
\be\label{form-of-B-ge}
{\mathtt B} 
\begin{pmatrix}
Q  \\
P 
\end{pmatrix}
:= \begin{pmatrix}
\Pi^\bot_{\mathbb S} D_V^{-1/2} \big( 3 a(x) ( v(\phi, 0 , \xi) )^2 + \e 4  a_4(x) ( v(\phi, 0, \xi) )^3 ) D_V^{-1/2} Q \big)   \\
0  
\end{pmatrix}   
\ee
with the function  $ v $  defined in \eqref{defv}, 
the functions $ a(x)$ and  $ a_4(x) $ are in \eqref{nonlinearity1},
the vector $ \xi = \xi (\lambda) \in [1,2]^\es $ in \eqref{xil}, and $ {\mathtt r}_\e  := {\mathtt r}_\e (\uF ) $
is a self-adjoint  remainder   satisfying
\begin{align}
& | {\mathtt r}_\e |_{\Lip,  +, s_1}  \lesssim_{s_1}   \e^4  \, , \quad \label{estimate:rep-ge}  \\
& | {\mathtt r}_\e |_{\Lip,  +,s}  \lesssim_s  \e^2 (\e^2 + \| \uF  \|_{\Lip, s+2}) \, .   \label{estimate:rep-ge0}
\end{align}
Moreover, given another torus $ \ui' = (\vphi, 0, 0) + \uF' $ satisfying \eqref{ansatz 0},  we have
\be\label{estimate:rep-ge0-diff}
| {\mathtt r}_\e - {\mathtt r}_\e' |_{+, s_1}  \lesssim_{s_1}   \e^2  \| \uF  - \uF' \|_{s_1+2} \, .   
\ee
\end{lemma}

The next Chapters \ref{sec:6n}-\ref{sec:proof.Almost-inv} will be devoted to obtain an approximate right inverse 
of $ {\cal L}_\om (\ui) $, as stated in Proposition \ref{prop:inv-ap-vero}.
In Chapter \ref{sec:6n} we shall conjugate $ {\cal L}_\om (\ui) $ to an operator (see \eqref{transf-op-1}-\eqref{newA+} and 
\eqref{coniugazione-yes})
which is in a suitable form to apply Proposition \ref{prop-cruciale}, and so proving Proposition \ref{prop:inv-ap-vero}.  
Proposition \ref{prop-cruciale}  is proved in Chapters \ref{sec:splitting} and \ref{sec:proof.Almost-inv}.

The rest of this chapter is devoted to the proof of Proposition \ref{Prop:sec6} and Lemma \ref{thm:Lin+FBR}. 

\section{Proof of Proposition \ref{Prop:sec6}}
\label{sec:decoupling}

By \eqref{dioep}, for all $ \lambda \in \Lambda $ the  frequency vector
$ \om =  (1+ \e^2 \l ) \bar \om_\e $ satisfies the  Diophantine condition 
\be\label{dioph}
 | \om \cdot \ell | \geq \frac{\g_2}{ \langle \ell \rangle^{ \t_1}} \, , \quad \forall \ell \in \Z^\es \setminus \{ 0 \} \, , 
 \quad {\rm where} \quad \gamma_2 = \gamma_1 / 2  = \gamma_0 / 4 \, . 
\ee
We recall that the constant $ \g_0 $ in \eqref{diop} depends only on the potential $ V(x) $, and it is considered 
as a fixed $ O(1) $ quantity, and thus we shall not track its dependence in the estimates.  

An invariant torus $ \ui  $ for the Hamiltonian vector field $X_{K}$,  supporting a  Diophantine flow,  
is isotropic\index{Isotropic torus} (see e.g. Lemma 1 in \cite{BBField} 
and Lemma \ref{lem:iso}), 
namely the pull-back $ 1$-form $ \ui^* \form $ is closed, 
where $ \form $\index{Liouville $ 1$-form} is the Liouville 1-form defined in \eqref{Lambda 1 form}. 
This is equivalent to say that the 2-form 
$$ 
\ui^* {\cal W} =  \ui^* d \form  = d (\ui^* \form )= 0 
$$ 
vanishes,  where $ {\cal W} = d \form $ is  defined in \eqref{2form-y}. 
Given an ``approximately invariant" torus $ \ui $, 
 the 1-form $ \ui^* \form $ is only  ``approximately closed".
In order to make this statement quantitative we consider
\begin{equation}\label{coefficienti pull back di Lambda}
\begin{aligned}
& \ui^* \form = \sum_{k = 1, \ldots, \es} a_k (\vphi) d \vphi_k \,, \\
& a_k(\vphi) :=  \big( [\pa_\ph \uth (\vphi)]^\top \uy (\vphi)  \big)_k 
+ \frac12 ( \partial_{\vphi_k} \uz (\ph), J \uz (\ph) )_{L^2(\T_x)}
\end{aligned}
\end{equation}
and we quantify how small is  the pull-back $ 2$-form
\begin{equation} \label{def Akj} 
\begin{aligned}
 \ui^* {\cal W} = d \, \ui^* \form  = \sum_{k, j =1, \ldots, \es, k <  j} 
& A_{k j}(\vphi) d \vphi_k \wedge d \vphi_j\,, \\
&  A_{k j} (\vphi) :=   \partial_{\vphi_k} a_j(\ph) - \partial_{\vphi_j} a_k(\ph)  \, , 
\end{aligned}
\end{equation} 
in terms of the error function\index{Error function} $Z(\vphi) $ defined in \eqref{def Zetone}.

\begin{lemma} 
The coefficients $ A_{kj}  $ in \eqref{def Akj} satisfy 
\begin{equation}\label{stima A ij}
\| A_{k j} \|_{\Lip, s} \lesssim_s \| Z \|_{\Lip, s+ \tau_1 + 1} + \| Z \|_{\Lip, s_0+1} 
\|  \uF \|_{\Lip, s+  \t_1 + 1} \,.
\end{equation}
\end{lemma}

\begin{pf}
The coefficients $ A_{kj}  $  satisfy the identity (see \cite{BBField}, Lemma 5, and \eqref{coefficientiAkj}) 
$$ 
 \Dom A_{k j} 
=  {\cal W}\big( \pa_\ph Z(\vphi) \underline{e}_k ,  \pa_\ph \ui(\vphi)  \underline{e}_j \big) 
+  {\cal W} \big(\pa_\ph \ui(\vphi) \underline{e}_k , \pa_\ph Z(\vphi) \underline{e}_j \big)  
$$
where  $ \underline{e}_k  $ denote the $ k $-th vector of the canonical basis of $ \R^\es $. 
Then by \eqref{ansatz 0} we get 
$$
\| \Dom A_{k j} \|_{\Lip, s}
\lesssim_s \| Z \|_{\Lip, s+1} + \| Z \|_{\Lip, s_0 + 1} \| \uF \|_{\Lip, s + 1}   \, . 
$$
Notice that  the functions $ A_{kj} (\vphi) $ defined in  \eqref{def Akj} have zero mean value in $ \vphi $, so that 
$$
A_{k j}(\vphi)  =   (\om \cdot \pa_\vphi )^{-1} \big({\cal W}\big( \pa_\ph Z(\vphi) \underline{e}_k ,  \pa_\ph \ui(\vphi)  \underline{e}_j \big) 
+  {\cal W} \big(\pa_\ph \ui(\vphi) \underline{e}_k , \pa_\ph Z(\vphi) \underline{e}_j \big)   \big) \, .
$$
Now, since $ \om $ is Diophantine according to  \eqref{dioph}, by \eqref{diof-est}  we have 
$$ 
\| (\om \cdot \pa_\vphi)^{- 1} g \|_{\Lip,s} \leq C  \| g\|_{\Lip, s + \tau_1} 
$$
where  the  constant $ \gamma_2 = \g_0 / 4 $ is included in $ C $ 
because it is considered a fixed constant $ O(1) $.
By the expression of $ {\cal W} $ in \eqref{2form-y}, the tame estimate \eqref{inter:pro-Lip}
and \eqref{ansatz 0} we deduce \eqref{stima A ij}. 
\end{pf}

We now modify the approximate torus $ \ui $ to obtain an isotropic torus $ i_\d  = i_\d  (\ui) $ nearby,  
which is 
still approximately invariant. We denote the Laplacian  
$$ 
\Delta_\vphi := \sum_{k = 1, \ldots, \es} \partial_{\vphi_k}^2 \, . 
$$ 

\begin{lemma}\label{toro isotropico modificato} {\bf (Isotropic torus)} 
The torus $ i_\delta(\vphi) :=  
(\uth (\vphi), y_\delta(\vphi), \uz (\vphi) ) $ defined by 
\begin{equation}\label{y 0 - y delta}
y_\d (\vphi)  := \uy (\vphi)  -  [\pa_\ph \uth (\vphi)]^{- T}  \rho(\vphi) \, ,  \ \rho = (\rho_j)_{j \in {\mathbb S}} \,,  \quad 
\rho_j(\vphi) := \Delta_\vphi^{-1} \sum_{k = 1, \ldots, \es} \partial_{\vphi_k} A_{k j}(\vphi) 
\end{equation}
 is {\it isotropic}. 
There is $ \bt := \bt(\es,\t_1) $ such that \eqref{stima y - y delta}-\eqref{stima toro modificato} hold. Moreover
\be\label{2015-2}
 \| y_\delta - \uy \|_{\Lip, s} \lesssim_s \| \uF \|_{\Lip, s+1} \\
\ee
and  \eqref{derivata i delta} holds. 
 \end{lemma}
Along the section we denote 
by $ \bt := \bt(\es, \tau_1 ) $ possibly different (larger) ``loss of derivatives"  constants. 
\\[1mm]
\begin{pf}
The proof of the isotropy of the torus $i_\d$ is in Lemma 6 of \cite{BBField}, see Lemma
\ref{modified}. 
Let us prove the bounds \eqref{stima y - y delta}-\eqref{derivata i delta}  and \eqref{2015-2}.
First notice that, since the map $A \mapsto A^{-1}$  is $C^\infty$ on the open set  of invertible matrices 
$ \{ A \in M_\es (\R) \, : \  \det A \neq 0 \}$, we derive from 
\eqref{ansatz 0} and Lemma \ref{Moser norme pesate} that, for $ \e $ small, the map 
$$ 
D \uth^{-1} : \vphi \mapsto [D \uth (\vphi)]^{-1}
$$ 
satisfies the tame estimates
\be \label{Duth-1}
\| D \uth^{-1} \|_{\Lip , s } \lesssim_s 1 + \| \uF \|_{s+1} \, , \quad \forall s \geq s_0 \, . 
\ee
Then \eqref{stima y - y delta}  and \eqref{2015-2} follow
by \eqref{y 0 - y delta}, \eqref{Duth-1},  \eqref{coefficienti pull back di Lambda}, \eqref{def Akj}, 
 \eqref{stima A ij} and \eqref{ansatz 0}. Moreover, we have that the difference 
$$
{\cal F}(i_\delta ) - {\cal F}(\ui ) =  
\begin{pmatrix} 0 \\ \Dom (y_\delta  - \uy )  \\ 0 \end{pmatrix}\, 
+ \e^2 \big( X_R(i_\delta ) - X_R(\ui) \big)
$$
and  \eqref{stima toro modificato} follows by  
\eqref{stima y - y delta}, Lemma\index{Moser estimates for composition operator} \ref{Moser norme pesate} and \eqref{ansatz 0}. 
Finally the bound  \eqref{derivata i delta} follows by 
\eqref{y 0 - y delta}, \eqref{def Akj}, \eqref{coefficienti pull back di Lambda}, \eqref{ansatz 0}.
\end{pf}

 It is proved in \cite{BBField} (see Lemma \ref{lemma:symplectic}) that the diffeomorphism 
$ G_\delta : (\phi, \ac , w) \to (\theta, y, z) $ 
defined in  \eqref{trasformazione modificata simplettica} is symplectic
because  the torus $ i_\d $ is isotropic (Lemma \ref{toro isotropico modificato}).
By construction, \eqref{toro-new-coordinates} holds.
Since $ G_\d $ is symplectic the Hamiltonian system generated by $ K $  transforms 
as in \eqref{new-Hamilt-K} into the
Hamiltonian system with Hamiltonian $ {\mathtt K} = K \circ G_\delta  $. 
By \eqref{parity solution} the transformation $ G_\d $ in  \eqref{trasformazione modificata simplettica} 
is also reversibility preserving and so the Hamiltonian 
$ {\mathtt K} $ is reversible, 
i.e. $ {\mathtt K} \circ \tilde S = {\mathtt K} $.  

\begin{lemma}
The tame estimates \eqref{DG delta} and \eqref{DG2 delta} hold. 
\end{lemma}

\begin{pf}
We write  \eqref{trasformazione modificata simplettica} as 
$$
 G_\delta \begin{pmatrix}
\phi \\
\ac \\
w
\end{pmatrix} := \left(
\begin{array}{l}
 \uth (\phi) \\
y_\delta (\phi) + M(\phi) \ac - \sum_{j=1}^\es (  m_j (\phi) , Jw )_{L^2_x} {\underline e}_j\\
 \uz(\phi) + w
\end{array} \right)
$$
where $( \underline{e}_j) $ denotes  the canonical basis of $\R^\es$ and we set 
$$ 
M(\phi) := [\pa_\phi \uth (\phi)]^{-\top} \, , 
\quad 
m_j(\phi) : =  (\partial_{\theta_j} {\tilde{\uz}}) (\uth (\phi)) = [(\partial_\vphi \uth)^{-\top} (\phi) \nabla \uz (\phi)]_j \, .
$$ 
The tame estimate \eqref{Duth-1} implies  
\be \label{Duth2}
\| M \|_{\Lip, s} \lesssim_s 1 + \|\uF\|_{\Lip , s+1} 
\ee
and using \eqref{inter:pro-Lip} and \eqref{ansatz 0} we have
\begin{align}   \label{g(1)est}
\| m_j \|_{\Lip ,s}   =\|(\partial_{\theta_j} \tilde{\uz}) (\uth (\phi))\|_{\Lip ,s} 
& \lesssim_s 
\|\uz\|_{\Lip, s+1} 
+ \|\uz \|_{\Lip, s_0+1} \| \underline \vartheta \|_{\Lip, s+1}  \nonumber  \\ 
& \lesssim_s \|\uF\|_{\Lip , s+1}  \, .
\end{align}
Now 
$$
DG_\d (\vphi , 0 , 0) [\widehat{\imath}(\vphi)]:= DG_\d (\vphi , 0 , 0)
\begin{pmatrix}  \widehat{\phi } (\vphi)\\  \widehat{\zeta } (\vphi)\\
\widehat{w} (\vphi) \end{pmatrix}=
\begin{pmatrix}  \widehat{a} (\vphi)\\  \widehat{b } (\vphi)\\
\widehat{c} (\vphi) \end{pmatrix} 
$$
where 
\be \label{exDG}
\begin{aligned}
& \widehat{a} := \partial_\vphi  \uth (\vphi)  [\widehat{\phi}]  
 \  , \\
& \widehat{b} :=\partial_\vphi y_\d (\vphi)  [\widehat{\phi}]  
+ M(\vphi) [\widehat{\zeta} ]- \sum_{j=1}^\es ( m_j (\vphi) , J\widehat{w} )_{L^2_x} {\underline e}_j
\  , \\
& \widehat{c} := \partial_\vphi \uz (\vphi) [\widehat{\phi} ]  +
\widehat{w} 
\end{aligned}
\ee
and 
$$
D^2G_\d (\vphi , 0 , 0) [\widehat{\imath}_1 (\vphi), \widehat{\imath}_2 (\vphi)]=
\begin{pmatrix}  \widehat{\a} (\vphi)\\  \widehat{\b } (\vphi)\\
\widehat{\gamma} (\vphi) 
\end{pmatrix} 
$$
where 
\be
\begin{aligned} \label{exD2G}
&\widehat{\alpha} := \partial_\vphi^2\uth  [\widehat{\phi}_1 , \widehat{\phi}_2]  \, ,  \\
&\widehat{\beta} := \partial_\vphi^2 y_\d [\widehat{\phi}_1 , \widehat{\phi}_2]
+ \partial_\vphi M[\widehat{\phi}_1] \widehat{\zeta}_2 + \partial_\vphi M[\widehat{\phi}_2] \widehat{\zeta}_1  \\
& \qquad -\sum_{j=1}^\es \Big( \big( \partial_\vphi m_j [\widehat{\phi}_1] , J\widehat{w}_2 \big)_{L^2_x} +
\big( \partial_\vphi m_j  [\widehat{\phi}_2] , J\widehat{w}_1 \big)_{L^2_x}    \Big) {\underline e}_j  \, ,   \\
&\widehat{\gamma} := \partial_\vphi^2 \underline{z} [\widehat{\phi}_1 , \widehat{\phi}_2]  \, .  
\end{aligned}
\ee
The tame estimates \eqref{DG delta} and \eqref{DG2 delta}
are a consequence of \eqref{exDG},  \eqref{exD2G}, 
 \eqref{inter:pro-Lip}, \eqref{ansatz 0}, \eqref{2015-2} and
 \eqref{Duth2}, \eqref{g(1)est}.  
\end{pf}

Then we consider the Taylor expansion \eqref{KHG} 
of the Hamiltonian $ {\mathtt K} $ at the trivial torus $ (\phi , 0, 0 ) $. 
Notice that the Taylor coefficient $  {\mathtt K}_{00}(\phi ) \in \R $,  
$ {\mathtt K}_{10}(\phi ) \in \R^\es $,  
$ {\mathtt K}_{01}(\phi ) \in H_{\St}^\bot$, 
$ {\mathtt K}_{20}(\phi) $ is a $\es \times \es$ real matrix, 
$ {\mathtt K}_{02}(\phi)$ is a linear self-adjoint operator of $ H_{\St}^\bot $ and 
$ {\mathtt K}_{11}(\phi) \in {\cal L}(\R^\es, H_{\St}^\bot )$. 

The Hamilton equations associated to \eqref{KHG}  are 
\begin{equation}\label{sistema dopo trasformazione inverso approssimato}
\begin{cases}
\dot \phi \hspace{-30pt} & = {\mathtt K}_{10}(\phi) +  {\mathtt K}_{20}(\phi) \ac + 
{\mathtt K}_{11}^\top (\phi) w + \partial_{\ac} {\mathtt K}_{\geq 3}(\phi, \ac, w)
\\
\dot \ac \hspace{-30pt} & = 
- \partial_\phi {\mathtt K}_{00}(\phi ) - [\partial_{\phi} {\mathtt K}_{10}(\phi )]^\top  \ac - 
[\partial_{\phi} {\mathtt K}_{01}(\phi )]^\top  w  
\\
& \quad -
\partial_\phi \big( \frac12 {\mathtt K}_{2 0}(\phi) \ac \cdot \ac + ( {\mathtt K}_{11}(\phi) \ac , w )_{L^2(\T_x)} + 
\frac12 ( {\mathtt K}_{02}(\phi) w , w )_{L^2(\T_x)} + {\mathtt K}_{\geq 3}(\phi, \ac, w) \big)
\\
\dot w \hspace{-30pt} & = J \big( {\mathtt K}_{01}(\phi) + 
{\mathtt K}_{11}(\phi) \ac +  {\mathtt K}_{0 2}(\phi) w + \nabla_w {\mathtt K}_{\geq 3}(\phi, \ac, w) \big) 
\end{cases} 
\end{equation}
where $ \partial_{\phi} {\mathtt K}_{10}^\top $ is the $ \es \times \es $ transposed matrix and 
$$ 
\partial_{\phi} {\mathtt K}_{01}^\top (\phi) \, ,  
\ {\mathtt K}_{11}^\top (\phi) : {H_{\St}^\bot \to \R^\es} \, , \quad 
\forall \phi \in \T^\es \, ,
$$
 are defined by the 
duality relation 
\be\label{duality-K01T}
( \partial_{\phi} {\mathtt K}_{01} [\widehat \phi ],  w)_{L^2_x}  = \widehat \phi \cdot [\partial_{\phi} {\mathtt K}_{01}]^\top w  \, ,
\quad \forall \widehat \phi \in \R^\es, \ w \in H_{\St}^\bot \, ,
\ee
where $ \cdot $ denotes the scalar product in $ \R^\es $. 
The transposed operator  $ {\mathtt K}_{11}^\top $ is  similarly defined and it turns out to be  
the following operator:  for all  $ w \in H_{\St}^\bot $, 
and denoting $\underline{e}_k$ the $k$-th vector of the canonical basis of $\R^\es$, 
\begin{equation} \label{K11 tras}
{\mathtt K}_{11}^\top(\phi) w =  \sum_{k = 1, \ldots, \es} 
\big( {\mathtt K}_{11}^\top(\phi) w \cdot \underline{e}_k\big) \underline{e}_k   =
\sum_{k = 1, \ldots, \es}  
\big( w, {\mathtt K}_{11}(\phi) \underline{e}_k  \big)_{L^2(\T_x)}  \underline{e}_k  \, \in \R^\es \, .  
\end{equation}
The terms $ \pa_\phi {\mathtt K}_{00} $, $ {\mathtt K}_{10} - \om $, 
$ {\mathtt K}_{01} $  in the Taylor expansion \eqref{KHG} vanish 
if $ Z  = 0 $.

\begin{lemma} \label{coefficienti nuovi} 
\eqref{apprVF} and  \eqref{K 00 10 01} hold.
\end{lemma}

\begin{pf}
Formula \eqref{apprVF} is proved in  Lemma 8 of \cite{BBField} (see Lemma \ref{lem:sc}).
Then \eqref{ansatz 0},  \eqref{stima y - y delta}, \eqref{stima toro modificato},  \eqref{DG delta} imply \eqref{K 00 10 01}.
\end{pf}

Notice that, if $ {\cal F} (\ui) = 0 $, namely $ \ui (\vphi )$ is an invariant torus for the 
Hamiltonian vector field $ X_K $, supporting  quasi-periodic solutions with frequency $ \om $, 
then, by \eqref{apprVF}-\eqref{K 00 10 01},  
the Hamiltonian in \eqref{KHG} simplifies to the KAM (variable coefficients) normal form 
\be\label{KAM:NF}
{\mathtt K}  = const + \om \cdot \ac + \frac12 {\mathtt K}_{2 0}(\phi) \ac \cdot \ac 
+ \big( {\mathtt K}_{11}(\phi) \ac , w \big)_{L^2(\T)} 
+ \frac12 \big( {\mathtt K}_{02}(\phi) w , w \big)_{L^2(\T)} + {\mathtt K}_{\geq 3} \, . 
\ee 
We now estimate  $ {\mathtt K}_{20}, {\mathtt K}_{11}$ in \eqref{KHG}. 

\begin{lemma} \label{lemma:Kapponi vari}
\eqref{stime coefficienti K 20 11 bassa}-\eqref{stime coefficienti K 11 alta trasposto}  hold. 
\end{lemma}

\begin{pf}
\\[1mm]
{\sc Proof of \eqref{stime coefficienti K 20 11 bassa}-\eqref{tame:K20-appl}. } 
By Lemma 9 of \cite{BBField} (see Lemma \ref{smalldepd}) 
and  the form of $ K $ in \eqref{primalinea},  we have
\begin{align}
{\mathtt K}_{2 0}(\phi)  & = [\pa_\phi \uth (\phi)]^{-1} \partial_{yy} K (i_\delta(\phi))  [\pa_\ph \uth (\phi)]^{-\top} \nonumber \\
&  = \e^2 [\pa_\phi \uth (\phi)]^{-1} \partial_{yy} R(i_\delta(\phi)) 
[\pa_\phi \uth (\phi)]^{-\top} \label{intemK20} \\
& =  \e^2 \partial_{yy} R(i_0(\phi))  +  r_{20} \label{K20-r20} 
\end{align}
where $ i_0(\phi) = (\phi, 0, 0 ) $ and  
$$
\begin{aligned}
r_{20} & := 
\e^2 \Big( [\pa_\phi \uth (\phi)]^{-1} \partial_{yy} R(i_\delta(\phi)) 
[\pa_\phi \uth (\phi)]^{-\top} - \partial_{yy} R(i_\delta(\phi)) \Big)  \\ 
& \quad +  \e^2 \Big(\partial_{yy} R(i_\delta(\phi)) - \partial_{yy} R(i_0(\phi)) 
\Big) \, . 
\end{aligned}
$$
By  Lemma \ref{Moser norme pesate}, \eqref{ansatz 0}, \eqref{2015-2} we have
$$ 
\| \partial_{yy} R (i_\d (\phi))- \partial_{yy} R (i_0(\phi)) \|_{\Lip , s_0} \leq C(s_1) \e^2 
$$ 
and, using also
$ \| (\partial_\phi \uth (\phi))^{-1} - {\rm Id}_{ \es }\|_{\Lip , s_0} \leq $ $ \|\uF \|_{\Lip, s_0+1} \lesssim_{s_1} \e^2 $, that
$$
\|  [\pa_\phi \uth (\phi)]^{-1} \partial_{yy} R(i_\delta(\phi)) 
[\pa_\phi \uth (\phi)]^{-\top} - \partial_{yy} R(i_\delta(\phi)) \|_{\Lip,s_0} \leq C(s_1) \e^2 \, . 
$$
Therefore
 \be\label{est:r20} 
 \| r_{20} \|_{\Lip, s_0} \leq C(s_1) \e^4 \, .
 \ee
Moreover, by \eqref{intemK20} and  Lemma\index{Moser estimates for composition operator} \ref{Moser norme pesate} the norm of $ {\mathtt K}_{20}$ (which is the sum 
of the norms of its $ \es \times \es $ 
matrix entries)  
satisfies 
$$
 \| {\mathtt K}_{20} \|_{\Lip,s}  \lesssim_s  \e^2 (1 +  \| \uF \|_{\Lip,s+ \bt})  
$$ 
and  \eqref{tame:K20-appl} follows by 
the tame estimates \eqref{inter:pro-Lip} for the product of functions. 
 
Next, recalling the expression of $ R $ in \eqref{restoNN}, with $ G$ as in \eqref{svilG1}, by 
computations similar to those in Section \ref{sec:shifted-tan}, it turns out  that the average with respect to 
$ \phi $ of $ \partial_{yy} R(i_0(\phi)) $ is 
\be\label{eq:R-da-twist}
\langle \partial_{yy} R(i_0(\phi)) \rangle = \Ab + r   \quad {\rm with} \quad \| r \|_{\Lip} \leq C \e^2  
\ee
where $ \Ab  $ is the twist matrix defined in \eqref{def:AB} 
(in particular  there is no contribution from $ a_4 (x) u^ 4 $). 

The estimate \eqref{stime coefficienti K 20 11 bassa} follows by \eqref{K20-r20}, 
\eqref{est:r20}   and \eqref{eq:R-da-twist}. 
\\[1mm]
{\sc Proof of \eqref{stime coefficienti K 11 alta}-\eqref{stime coefficienti K 11 alta trasposto}. } 
By Lemma 9 of \cite{BBField} (see Lemma \ref{smalldepd}) and the form of $ K $ in \eqref{primalinea} we have
\begin{align}
{\mathtt K}_{11}(\phi) & = 
\partial_{y} \nabla_z K (i_\delta(\phi)) [\pa_\ph \uth  (\phi)]^{-\top} 
+ J (\pa_\theta {\tilde \uz}) (\uth (\phi)) (\partial_{yy} K) (i_\delta(\phi)) [\pa_\ph \uth  (\phi)]^{-\top}   \nonumber\\
&  = 
\e^2 \partial_{y} \nabla_z R(i_\delta(\phi)) [\pa_\ph \uth  (\phi)]^{-\top} 
+ \e^2 J (\pa_\theta {\tilde \uz}) (\uth(\phi)) (\partial_{yy} R) (i_\delta(\phi)) [\pa_\ph \uth  (\phi)]^{-\top}\,, 
\end{align}
and using \eqref{ansatz 0}, 
\eqref{2015-2}, we deduce \eqref{stime coefficienti K 11 alta}. 
The bound \eqref{stime coefficienti K 11 alta trasposto} for $ {\mathtt K}_{11}^\top$ follows by \eqref{K11 tras}
and \eqref{stime coefficienti K 11 alta}.
\end{pf}

Finally formula \eqref{lin idelta} is obtained  linearizing \eqref{sistema dopo trasformazione inverso approssimato}.

\section{Proof of Lemma  \ref{thm:Lin+FBR}} \label{sec:Lomega}

We have to compute the quadratic term $ \frac12 ( {\mathtt K}_{02}(\phi) w, w )_{L^2(\T_x)} $  
in the Taylor expansion \eqref{KHG}
of  the Hamiltonian $ {\mathtt K} (\phi, 0, w)$. 
The operator $ {\mathtt K}_{02}(\phi)$ is
\begin{equation}
\label{trasformata forma normale}
\begin{aligned}
{\mathtt K}_{02}(\phi) & = \partial_w \nabla_w {\mathtt K} (\phi, 0, 0) \\
& = \partial_w \nabla_w (K \circ G_\delta)(\phi, 0, 0) =
D_V + \e^2 \partial_w \nabla_w (R \circ G_\delta)(\phi, 0, 0)
\end{aligned}
\ee
where  the Hamiltonian $ K  $ is defined in \eqref{primalinea} and $ G_\delta $ is the symplectic diffeomorphism
\eqref{trasformazione modificata simplettica}. 
 Differentiating with respect to $w$ the Hamiltonian 
$$ 
(R \circ G_\delta)(\phi, \ac, w) =  R(\uth (\phi), y_\delta(\phi) + L_1(\phi) \ac + L_2(\phi) w, \uz(\phi) + w) 
$$
where, for brevity, we set
\be\label{def:L1L2}
L_1(\phi) := [\partial_\phi \uth (\phi)]^{- T} \, , \quad  L_2(\phi) := - [\partial_\theta \tilde \uz(\uth (\phi))]^\top J \, ,
\quad \tilde{\uz} (\theta) = \uz (\uth^{-1} (\theta)) \, , 
\ee
(see \eqref{trasformazione modificata simplettica}),  we get 
$$
 \nabla_w (R \circ G_\delta)(\phi, \ac, w) = 
 L_2(\phi)^\top \partial_y R(G_\delta(\phi, \ac, w)) + \nabla_z R(G_\delta(\phi, \ac, w)) \, .
$$  
Differentiating such identity with respect to $ w $, and  recalling \eqref{toro-new-coordinates}, we get 
\be\label{P G delta}
 \partial_w \nabla_w (R \circ G_\delta)(\phi, 0, 0) = 
\partial_z \nabla_z R(i_\delta(\phi)) + r(\phi) 
\ee
with a self adjoint remainder $ r(\phi) := r_1(\phi) + r_2(\phi) + r_3 (\phi) $ given by  
\be\label{formr1r2r3}
\begin{aligned}
& r_1(\phi) := L_2(\phi)^\top \partial_{yy} R(i_\delta(\phi)) L_2(\phi)\,, \\ 
& r_2(\phi):=  L_2(\phi)^\top 
\nabla_z  \partial_{y} R(i_\delta(\phi)) \,, \\
& r_3(\phi) := \partial_y \nabla_z R(i_\delta(\phi))L_2(\phi) \, .  
\end{aligned}
\ee
Each operator $ r_1, r_2, r_3 $
is the composition of at least one operator with ``finite rank $ \R^\es $ in the space variable", 
and therefore 
it has the ``finite dimensional" form
\begin{equation}\label{forma buona resto}
r_l (\phi)[h] = \sum_{j = 1, \ldots, \es } \big(h  \,,\,g_j^{(l )} (\phi , \cdot) \big)_{L^2_x} 
\chi_j^{(l )} (\phi , \cdot) \,, \quad \forall h \in H_{{\mathbb S}}^\bot \, ,  \ l  = 1, 2, 3 \, , 
\end{equation}
for functions $ g_j^{(l )} (\phi, \cdot ), \chi_j^{(l )} (\phi , \cdot) \in H_{{\mathbb S}}^\bot  $. 
Indeed, writing the operator $L_2(\phi) : H_{{\mathbb S}}^\bot \to \R^\es $ as  
$$ 
L_2(\phi)[h] = \sum_{j  = 1, \ldots, \es } 
\big(h\,,\, L_2(\phi)^\top[ {\underline e}_j ] \big)_{L^2_x} {\underline e}_j \, , \quad 
 \forall h \in H_{{\mathbb S}}^\bot \, , 
$$
 we get, by \eqref{formr1r2r3}, 
\begin{align}\label{forma:r1}
& r_1(\phi)[h] = \sum_{j  = 1, \ldots, \es} \big(h\,,\, L_2(\phi)^\top[{\underline e}_j ] \big)_{L^2_x} 
A_1[{\underline e}_j ]\,,\quad A_1 := L_2(\phi)^\top \partial_{yy} R(i_\delta(\phi)) \, , \\
& 
r_2(\phi)[h]  =  \sum_{j  = 1, \ldots, \es} \big(h, A_2^\top[{\underline e}_j] \big)_{L^2_x} L_2(\phi)^\top[{\underline e}_j] \, ,
\quad A_2 := \pa_z \pa_{y} R(i_\delta(\phi)) \, , \label{forma:r2} \\
&  
 r_3(\phi)[h] = \sum_{j  = 1, \ldots, \es} \big(h, L_2(\phi)^\top [ {\underline e}_j] \big)_{L^2_x} A_3[{\underline e}_j] \, ,  \quad  
 A_3 := \partial_y \nabla_z R(i_\delta(\phi)) = A_2^\top \, . \label{forma:r3}
\end{align} 

\begin{lemma}\label{lem:gjhj}
For all $ l = 1, 2, 3 $, $ j = 1, \ldots, \es $, we have, for all $  s \geq s_0 $, 
\begin{equation}\label{stime gj chij}
\begin{aligned}
\| g_j^{(l )} \|_{\Lip, s} + \| \chi_j^{(l )} \|_{\Lip, s}  & \lesssim_s
1  
+ \| \uF \|_{\Lip , s+1}  \\ 
\min\{ \| g_j^{(l )} \|_{\Lip, s}, \| \chi_j^{(l )} \|_{\Lip, s}\} & \lesssim_s
\| \uF\|_{\Lip, s + 1} \, . \\
\end{aligned}
\end{equation}
\end{lemma}

\begin{pf}
Recalling the expression of $ L_2 $ in 
\eqref{def:L1L2} and \eqref{forma:r1}-\eqref{forma:r3} we have:
\\[1mm]
{\bf i)}  
$ g_j^{(1)}(\phi)=g_j^{(3)}(\phi)=\chi_j^{(2)}(\phi)=L_2 (\phi)^\top [\underline{e}_j]=J (\partial_{\theta_j} \tilde{\uz}) (\uth (\phi)) $.
Therefore by \eqref{g(1)est},
$$
\| g_j^{(1)} \|_{\Lip,s} = \| g_j^{(3)} \|_{\Lip,s} = \| \chi_j^{(2)} \|_{\Lip,s} 
\lesssim_s  \|\uF\|_{\Lip , s+1} \, . 
$$
{\bf ii)} $ \chi^{(1)}_j = L_2(\phi)^\top \partial_{yy} R (i_\d (\phi)) [ \underline{e}_j] $. Then, recalling \eqref{def:L1L2} and
using Lemma \ref{Moser norme pesate}, we have 
\begin{eqnarray*}
\| \chi^{(1)}_j \|_{\Lip ,s} 
&\lesssim_s &  \|\partial_{\theta} \tilde{\uz} (\uth (\phi))\|_{\Lip ,s} (1+ \| \ddI \|_{\Lip,s_0} ) + 
\|\partial_{\theta} \tilde{\uz} (\uth (\phi))\|_{\Lip ,s_0}  \| \ddI \|_{\Lip,s} \nonumber \\
&\stackrel{\eqref{ansatz 0}\eqref{g(1)est}\eqref{2015-2}}{\lesssim_s} &  \| \uF \|_{\Lip,s+1}    \, . 
\end{eqnarray*} 
{\bf iii)}
$ \chi_j^{(3)}=g_j^{(2)}=\partial_{y_j} \nabla_z R (i_\d ((\phi)) $. Then, using Lemma\index{Moser estimates for composition operator} \ref{Moser norme pesate}
and \eqref{2015-2}, we get 
$$
\|\chi_j^{(3)}\|_{\Lip, s}=\|g_j^{(2)}\|_{\Lip , s} {\lesssim_s} \, 1+ \|\ddI\|_{\Lip, s} \lesssim_s 1 + \| \uF \|_{\Lip, s+1}  \, . 
$$
Items $ {\bf i)}-{\bf iii)} $ imply the estimates \eqref{stime gj chij}.
\end{pf}

We now use Lemmata \ref{gchi-r} and \ref{lem:gjhj} 
to derive bounds on the decay norms of the remainders 
$ r_l $ defined in \eqref{forma buona resto}.  Recalling Definition \ref{def:decay-sub}, 
we have to estimate  the norm $ | r_l |_{\Lip,+,s} =  | D_m^{1/2} r_l \Pi^\bot_{\mathbb S} D_m^{1/2} |_{\Lip,s}  $ 
of the 
extended operators,  
acting on the whole $ {\cal H}^0 = L^2 (\T^d \times \T^\es; \C^2 ) $, defined by 
\begin{align}  
D_m^{1/2} r_l \Pi^\bot_{\mathbb S} D_m^{1/2} h & := 
\sum_{j  = 1, \ldots, \es}  \big(\Pi^\bot_{\mathbb S}  D_m^{1/2} h\,,\,  g_j^{(l)} \big)_{L^2_x} 
(D_m^{1/2} \chi_j^{(l )}) \nonumber \\ & = 
 \sum_{j  = 1, \ldots, \es}  \big(  h\,,\, D_m^{1/2} g_j^{(l)} \big)_{L^2_x} 
(D_m^{1/2} \chi_j^{(l )}) \,, \quad  h \in {\cal H}^0      \, ,  \label{con-plus}
\end{align}
where we used that $g_j^{(l)} \in H^{\bot}_{\mathbb S} $. Then the decay norm of $ r = r_1 + r_2 + r_3 $ satisfies 
\begin{align}
  | r  |_{\Lip,+, s} &  \lesssim 
 \max_{ l=1,2,3}   | D_m^{1/2} r_l \Pi^\bot_{\mathbb S} D_m^{1/2}  |_{\Lip,s} \nonumber \\ 
&  \stackrel{\eqref{con-plus}, \eqref{boundsimpleL}} {\lesssim_s} \!\!\!\!\!\! 
\!\!\!\!\max_{ l=1,2,3, 
 j = 1, \ldots, \es} 
 \| D_m^{1/2} g_j^{(l )} \|_{\Lip, s} \| D_m^{1/2} \chi_j^{(l )} \|_{\Lip, s_0} +  
 \| D_m^{1/2} g_j^{(l )} \|_{\Lip, s_0} \| D_m^{1/2} \chi_j^{(l )} \|_{\Lip, s} \nonumber \\
&  \stackrel{\eqref{stime gj chij}, \eqref{ansatz 0}} {\lesssim_s} 
\| \uF \|_{\Lip, s + 2} \, .  \label{restimr}
\end{align}
Finally, by \eqref{trasformata forma normale}, \eqref{P G delta} we have
\begin{align}
{\mathtt K}_{02}(\phi) & = D_V  + \e^2 \partial_z \nabla_z R(i_\delta(\phi)) + \e^2 r(\phi) \nonumber \\
& = D_V  + \e^2 \partial_z \nabla_z R(\phi, 0, 0 ) + \e^2 
\big( \partial_z \nabla_z R(i_\delta (\phi)) -   \partial_z \nabla_z R(\phi, 0, 0) \big)  + \e^2 r(\phi)   \nonumber \\
& = D_V  +  \e^2 {\mathtt B} + {\mathtt r}_\e \label{prre1} 
\end{align}
where $ {\mathtt B} $ is defined in \eqref{form-of-B-ge} and 
\be\label{repvero}
{\mathtt r}_\e :=  \e^2 \partial_z \nabla_z R(\phi, 0, 0) - \e^2 {\mathtt B} + 
\e^2 \big( \partial_z \nabla_z R(i_\delta (\phi)) -   \partial_z \nabla_z R(\phi, 0, 0) \big)  + \e^2 r(\phi) \, .
\ee
This is formula \eqref{lin:normal-f-ge}.  
We now prove that $ {\mathtt r}_\e $  satisfies \eqref{estimate:rep-ge0}. 
In \eqref{restimr} we have yet estimated $ |r|_{\Lip,  +, s} $.  
Recalling Definition \ref{def:decay-sub}, 
and the expression of $ R $ in \eqref{restoNN} (see also \eqref{gradRQ}) 
we have to estimate the 
decay norm of the extended operator, acting on the whole
 $ {\cal H}^0 = L^2 (\T^d \times \T^\es; \C^2 ) $, defined as 
\be \label{exprRzz}
\begin{aligned}
& \partial_z \nabla_z R ( i_\d (\phi) ) \Pi^\bot_{\mathbb S}
\begin{pmatrix}
h_1  \\
h_2 
\end{pmatrix}
:=  \\
& \begin{pmatrix}
\Pi^\bot_{\mathbb S} D_V^{-\frac12}  (\pa_u g)(\e, x,  v( \uth (\phi), y_\d (\phi), \xi) + D_V^{-\frac12} \uQ (\phi) )  D_V^{-1/2} \Pi^\bot_{\mathbb S} h_1    \\
0  
\end{pmatrix}   
\end{aligned}
\ee
where $ g(\e, x, u) = \pa_u G (\e, x, u) $ is the nonlinearity in \eqref{nonlinearity:gep}.
Hence, by Proposition \ref{prop:decay-prod}
and Lemma \ref{pisig}, we have
\begin{align} \label{idi0}
& |\partial_z \nabla_z R(i_\delta (\phi)) -  \partial_z \nabla_z R(\phi, 0, 0)|_{\Lip,+, s} \lesssim_s \nonumber \\
& \big\| (\pa_u g)(\e, x,  v(\uth (\phi), y_\d (\phi), \xi) + D_V^{-1/2} \uQ (\phi) ) - 
 (\pa_u g)(\e, x,  v(\phi, 0, \xi) )  \big\|_{\Lip , s} \nonumber \\
 & \stackrel{Lemma  \, \ref{Moser norme pesate}, \eqref{ansatz 0}} {\lesssim_s}  \| \ddI \|_{\Lip ,s}  \stackrel{ \eqref{2015-2}}{\lesssim_s}  \| \uF \|_{\Lip, s + 1} \, .
\end{align}
Moreover, recalling  \eqref{nonlinearity:gep1}
and  the definition of $ {\mathtt B} $ in \eqref{form-of-B-ge},  we have 
\be\label{prre2}
|\partial_z \nabla_z R ( \phi, 0, 0 )-  {\mathtt B} |_{\Lip,+,s} \lesssim_s \e^2 \, . 
\ee
In conclusion, the operator 
$ {\mathtt r}_\e  $ defined in \eqref{repvero} satisfies, by 
\eqref{prre2}, \eqref{idi0},  \eqref{restimr}, the estimate
  \eqref{estimate:rep-ge0}. 
In particular, 
\eqref{estimate:rep-ge} holds by \eqref{ansatz 0}.

There remains to prove \eqref{estimate:rep-ge0-diff}. 
Let $ \ui' =(\vphi,0,0)+{\uF}'$ be another torus embedding satisfying \eqref{ansatz 0} and let 
$ {\mathtt r}'_\e  $ be the associated remainder 
in 
\eqref{repvero}.  We have
$$
{\mathtt r}'_\e (\phi) - {\mathtt r}_\e (\phi) = \e^2 \big(\pa_z \nabla_z R (i'_\d (\phi))- \pa_z \nabla_z R (i_\d (\phi)) \big) +
\e^2 (r'(\phi) - r(\phi)) \, .
$$
By Lemma \ref{Moser norme pesate} and the expression \eqref{exprRzz} of $\pa_z \nabla_z R$, 
\be \label{626-1}
|\pa_z \nabla_z R (i'_\d (\phi))- \pa_z \nabla_z R (i_\d (\phi))|_{+,s_1} 
\lesssim_{s_1} \|\ddI'-\ddI\|_{s_1} \stackrel{\eqref{derivata i delta}}{\lesssim_{s_1}} 
\|\uF'-\uF\|_{s_1 +1} \, . 
\ee 
Moreover let  ${g'}_j^{(l)}$ and ${\chi'}_j^{(l)}$
be the functions obtained in \eqref{forma buona resto}-\eqref{forma:r3} from ${i}'$. 
Using again
Lemma \ref{Moser norme pesate} as well as \eqref{derivata i delta}  and \eqref{2015-2}, we obtain, for any 
$ j = 1,2,3 $, $ j  = 1, \ldots, \es  $, 
$$
 \| {g'}_j^{(l)}- g_j^{(l )} \|_s +\| {\chi'}_j^{(l )}-  \chi_j^{(l )}
\|_s  
\lesssim_s 
\| {\uF'} - \uF \|_{s + 1}+ \| \uF \|_{s + 1} \| {\uF'} - \uF\|_{s_0 + 1}  \, . 
$$
Hence, by Lemma \ref{gchi-r}, $|r'-r|_{+,s_1} \lesssim_{s_1} \| \uF' - \uF \|_{s_1+2}$. With
\eqref{626-1}, this gives \eqref{estimate:rep-ge0-diff}.

\chapter{Splitting of low-high normal subspaces up to $ O(\e^4)$}\label{sec:6n}

The main result of this chapter is Proposition \ref{prop:op-averaged}.
Its goal is to transform 
the linear operator $ {\cal L}_\om $  in \eqref{Lomega def},
into 
a  form (see \eqref{newA+}) 
suitable to 
apply Proposition \ref{prop-cruciale} in the next chapter, which will enable 
to prove the existence of   
 an approximate right inverse of $ {\cal L}_\om $ for most values of the parameter  $ \l  $. 
 
In the next section we fix the set ${\mathbb M} $ in the splitting 
$ H_{\mathbb S}^\bot = H_{\mathbb M} \oplus H_{\mathbb M}^\bot $.

\section{Choice of $ {\mathbb M}$ }

We first remind that, by Lemma \ref{thm:Lin+FBR}, 
the linear operator $ {\cal L}_\om $ defined in \eqref{Lomega def},
acting in the normal subspace $ {H_{\St}^\bot} $, has the  form
\be\label{lin:normal-f}
{\cal L}_\om  =
\om \cdot \partial_\vphi - J (D_V + \e^2 {\mathtt B} + {\mathtt r}_\e ) \, , \quad  \om = (1+\e^2 \lambda) \bar \om_\e \, ,
\ee
where $ {\mathtt B}  = {\mathtt B}(\e, \lambda) $ is 
the self-adjoint operator  in \eqref{form-of-B-ge}
and the self-adjoint 
remainder $ {\mathtt r}_\e  $ satisfies \eqref{estimate:rep-ge}. 
Recalling \eqref{form-of-B-ge}, \eqref{defv}, \eqref{xil}, we derive, for $ \| a_4 \|_{L^\infty} \leq 1 $, 
the following bounds
for  the $ L^2 $-operatorial norm of $ {\mathtt B} $: 
\be\label{estimate:B-dB} 
\| {\mathtt B} \|_{0}  \leq C ( \| a \|_{L^\infty}  + \e) \, , \quad
\| {\mathtt B} \|_{\lip, 0}  \leq C ( \| a \|_{L^\infty} + \e) \| \Ab^{-1}  \| 
\ee
where $ \| \Ab^{-1}  \|  $ is some  norm of the inverse  twist matrix $ \Ab^{-1} $. 

Dividing \eqref{lin:normal-f} by $ 1 + \e^2 \l  $ we now consider the operator 
\be\label{def:primo-op}
\begin{aligned}
 \frac{{\cal L}_\om}{1 + \e^2 \l } = & 
\ \bar \om_\e \cdot \partial_\vphi - J \Big( {\mathtt A} +  \frac{{\mathtt r}_\e }{1+\e^2 \lambda} \Big)  \, , \\
& {\mathtt A} := {\mathtt A}(\lambda) = \frac{D_V}{1+ \e^2 \l} +  
\frac{\e^2 {\mathtt B}}{1+ \e^2 \l}   \, .
\end{aligned} 
\ee
We denote by $ \varrho $ the self-adjoint operator 
\be\label{def:varrho0}
\varrho := \frac{\e^2 {\mathtt B}}{1+ \e^2 \l} 
\ee
and, according to the splitting  $ H_{\mathbb S}^\bot = H_{\mathbb M} \oplus H_{\mathbb M}^\bot $, and 
taking in $ H_{\mathbb M} $ the basis 
$$ 
\{ ( \Psi_j(x) ,0) \, ,  (0, \Psi_j(x) ) \}_{j \in {\mathbb M}} \, ,
$$
we represent ${\mathtt A} $ as (recall that $ D_V \Psi_j = \mu_j \Psi_j $ and the notation \eqref{decoFGsotto})
\be\label{def:primo-op-bis}
{\mathtt A} = 
 \begin{pmatrix}
 {\rm Diag}_{j \in {\mathbb M} } \frac{{\mu}_j}{1+ \e^2 \l }  {\rm Id}_2 &  0   \\
  0 &  \frac{D_V}{1+  \e^2 \l}  \\
\end{pmatrix} + 
 \begin{pmatrix}
 \varrho_{{\mathbb M}}^{{\mathbb M}} &  \varrho_{{\mathbb M}}^{{\mathbb M}^c}    \\
 \varrho_{{\mathbb M}^c}^{{\mathbb M}}   &  \varrho_{{\mathbb M}^c}^{{\mathbb M}^c}   \\
\end{pmatrix} \, .
\ee
Recalling  \eqref{form-of-B-ge}, 
since the functions $ a, a_4  \in C^\infty (\T^d ) $, using \eqref{defv}, \eqref{xil}, and Proposition 
\ref{prop:decay-prod}, 
we derive that the operator $ \varrho $ defined in \eqref{def:varrho0} satisfies
\be\label{est:varrho}
| \varrho |_{\Lip, +,s} \leq C(s) \e^2  \, . 
\ee
In the next lemma we fix the subset of indices $ {\mathbb M} $. 
We recall the notation  $ {\mathtt A}^{{\mathbb M}^c}_{{\mathbb M}^c} =  
\Pi_{H_{{\mathbb M}^c}} {\mathtt A}_{| H_{{\mathbb M}^c} } = 
\Pi_{H_{\mathbb M}^\bot} {\mathtt A}_{| H_{\mathbb M}^\bot } $.

\begin{lemma}\label{choice:M} {\bf (Choice of $ \mathbb M $)}
Let 
\be\label{def:Omega-n}
{\mathtt \Omega}_j (\e, \l) := [\Bb \Ab^{-1} \bar \om_\e]_j + 
\frac{ \mu_j - [ \Bb \Ab^{-1} \bar \mu]_j}{1+ \e^2 \l} \, , \quad j \in {\mathbb M } \, , 
\ee
where $ \Ab $, $ \Bb $ are the ``Birkhoff" matrices\index{Birkhoff matrices} defined in \eqref{def:AB}, \eqref{def G}, 
the  $ \mu_j $ are the unperturbed frequencies defined in \eqref{auto-funzioni},  and the vectors 
$ \bar \mu $, $ \bar \om_\e $ are defined  respectively in \eqref{unp-tangential}, \eqref{def omep}. 
There is a constant $ C  >  0 $ such that, if $ {\mathbb M} \subset \N  $ contains the subset  
$ {\mathbb F} $  defined in \eqref{taglio:pos-neg}, i.e. $ {\mathbb F} \subset {\mathbb M} $, and 
\be\label{cond:suM}
\min_{j \in {\mathbb M}^c} \mu_j \geq  C \, \big( 1 + \e^2 \|a\|_{L^\infty} 
+ \| a \|_{L^\infty} \| \Ab^{-1} \| + \e \| \Ab^{-1} \| \big)  \, ,
\ee
then 
\be\label{negative-on-up}
\pa_\lambda {\mathtt A}^{{\mathbb M}^c}_{{\mathbb M}^c} 
\leq -  \Big( \max_{j \in {\mathbb F}} |\pa_\l {\mathtt \Omega}_j | + \e^2  \Big){\rm Id} \, . 
\ee
\end{lemma}

\begin{pf}
Differentiating the expression of ${\mathtt A} $ in 
\eqref{def:primo-op} we get
$$
\pa_\lambda {\mathtt A} := - \frac{\e^2 D_V }{(1+ \e^2 \l)^2} + 
\frac{\e^2 \pa_\lambda {\mathtt B}}{1+ \e^2 \l} -   \frac{\e^4 {\mathtt B}}{(1+ \e^2 \l)^2} \, . 
$$
Then by \eqref{estimate:B-dB} and the fact that 
$  [D_V]_{{\mathbb M}^c}^{{\mathbb M}^c} \geq   \min_{j \in {\mathbb M}^c} \mu_j  $, we get 
\begin{align}\label{lower-bound-der}
\pa_\lambda {\mathtt A}^{{\mathbb M}^c}_{{\mathbb M}^c} 
& 
\leq   -   \frac{\e^2  }{(1+ \e^2 \l)^2} \min_{j \in {\mathbb M}^c} \mu_j  + 
 C \e^4 (\|a\|_{L^\infty} + \e) + C \e^2  (\| a \|_{L^\infty}+ \e)  \| \Ab^{-1} \| \, . 
\end{align}
Now, recalling \eqref{def:Omega-n}, we have 
\be\label{upp-bound-dero}
\pa_\l {\mathtt \Omega}_j ( \e, \l) = \frac{- \e^2}{(1+ \e^2 \l )^2}
\big( \mu_j - [ \Bb \Ab^{-1} \bar \mu]_j \big) 
\ee
and so \eqref{lower-bound-der} and \eqref{upp-bound-dero} imply 
$$
\begin{aligned}
\pa_\lambda {\mathtt A}^{{\mathbb M}^c}_{{\mathbb M}^c}  + 
\max_{j \in {\mathbb F}}  |\pa_\l {\mathtt \Omega}_j | 
& \leq 
 - \frac{\e^2  }{(1+ \e^2 \l)^2} \Big(  \min_{j \in {\mathbb M}^c} \mu_j  - 
\max_{j \in {\mathbb F}} 
\big| \mu_j - [ \Bb \Ab^{-1} \bar \mu]_j \big| \Big) \\
& \quad + C \e^4 (\|a\|_{L^\infty} + \e) + 
C \e^2 (\|a\|_{L^\infty} + \e)\| \Ab^{-1} \|  \, .
\end{aligned}
$$
Thus  \eqref{negative-on-up} holds taking $ {\mathbb M} $ large enough  such that (recall that $ \mu_j \to + \infty $)
\be\label{penu1}
\min_{j \in {\mathbb M}^c} \mu_j \geq
\max_{j \in {\mathbb F}} 
\big| \mu_j - [ \Bb \Ab^{-1} \bar \mu]_j \big| +  C( 1 + \e^2 \|a\|_{L^\infty} + \| a \|_{L^\infty} \| \Ab^{-1} \|
+ \e \| \Ab^{-1} \| )  \, .
\ee
By \eqref{taglio:pos-neg}, \eqref{def:gs}, the fact that $ \mu_j \geq \sqrt{\beta } $ (see
\eqref{auto-funzioni}), we get 
\be\label{penu2}
\begin{aligned}
\max_{j \in {\mathbb F}} 
\big| \mu_j - [ \Bb \Ab^{-1} \bar \mu]_j \big| \leq
\gap & =  
\max_{j \in {\mathbb S}^c}\{   ( \Bb \Ab^{-1} \bar \mu)_j  -  \mu_j  \} \\ 
& \leq C(  \| a \|_{L^\infty} \| \Ab^{-1}\| + 1 ) 
\end{aligned}
\ee
for some $ C := C(\beta, {\mathbb S}) $, having used that 
\eqref{def G} we have 
$ \sup_{k \in {\mathbb S}} |G^j_k| \leq C \| a \|_{L^\infty} \beta $. 
By \eqref{penu1}-\eqref{penu2}  and  \eqref{cond:suM}  we deduce
\eqref{negative-on-up}. 
\end{pf}

In the sequel of the monograph the subset $ {\mathbb M} $ is kept fixed.
Note that the condition \eqref{cond:suM}  can be fulfilled taking $ {\mathbb M}  $ large enough because  $ \mu_j \to + \infty $ as $ j \to + \infty $.

In the next part of the chapter we perform one step of averaging 
to eliminate, as much as possible, the terms of order $ O(\e^2) $  of 
$ \varrho_{\mathbb M}^{\mathbb M} $, $ \varrho_{\mathbb M}^{{\mathbb M}^c} $, 
$ \varrho_{{\mathbb M}^c}^{{\mathbb M}} $ in \eqref{def:primo-op-bis}. 

\section{Homological equations}

According to the splitting 
$ H = H_{\mathbb M} \oplus H_{\mathbb M}^\bot $ we  consider   the linear map
$$
{\mathtt S} \mapsto J \ppamuvphi {\mathtt S}  + [J{\mathtt  S}, J D_V] 
$$
where $ \bar \mu \in \R^\es $ is the unperturbed tangential frequency vector defined in \eqref{unp-tangential}, and 
\be\label{forma cal-S0}
\begin{aligned}
&  {\mathtt S}(\vphi)= 
\begin{pmatrix}
{\mathtt d} (\vphi) & {\mathtt a}(\vphi)^* \\ 
{\mathtt a}(\vphi) & 0
\end{pmatrix} \in {\cal L} (H_{\mathbb S}^\bot) \, , \ \forall \vphi \in \T^\es \, ,  \\
& {\mathtt d}(\vphi) = {\mathtt d}^*(\vphi) \in {\cal L} (H_{\mathbb M}) \, ,  \   
{\mathtt a}(\vphi) \in {\cal L} (H_{\mathbb M}, H_{\mathbb M}^\bot) 
\end{aligned}
\ee
is  self-adjoint. 

Since $ D_V $ and $ J $ commute, we have 
\be\label{def:homo-matrix-1}
\begin{aligned}
& J \ppamuvphi {\mathtt S}  + [J{\mathtt S}, J D_V] \\
& = 
\begin{pmatrix}
J \ppamuvphi {\mathtt d}  + D_V {\mathtt d}  + J {\mathtt d} J D_V  & 
J \ppamuvphi  {\mathtt a}^*   + D_V  {\mathtt a}^*  + J {\mathtt a}^* J D_V    \\
J \ppamuvphi  {\mathtt a}  + D_V  {\mathtt a} + J {\mathtt a} J D_V & 0   
\end{pmatrix}  \, .
\end{aligned}
\ee
Recalling the definition of $ \Pi_{\mathtt D} $ in \eqref{def:pro+} (with $ {\mathbb F} \rightsquigarrow {\mathbb M} $) and 
of $ \Pi_{\mathtt O}  $ in \eqref{lem:off-dia} 
we decompose the self-adjoint operator 
$ \varrho = \e^2  {\mathtt B}  (1+ \e^2 \lambda)^{-1} $ defined in  \eqref{def:varrho0} as 
\be\label{def:varrho}
\varrho = \Pi_{\mathtt D} \varrho + \Pi_{\mathtt O} \varrho \, .
\ee
The term $ \Pi_{\mathtt O} \varrho $ has the form 
\be\label{decomposition-rho1-1}
\begin{aligned}
& \Pi_{\mathtt O}  {\varrho}(\vphi)=\begin{pmatrix} \varrho_1 (\vphi)& \varrho_2 (\vphi)^* \\
\varrho_2 (\vphi) & 0    \end{pmatrix} \in {\cal L} (H_{\mathbb S}^\bot ) \, , \\  
& \varrho_1 (\vphi) \in {\cal L} (H_{\mathbb M}) \, ,  \   \varrho_2 (\vphi) \in {\cal L} (H_{\mathbb M}, H_{\mathbb M}^\bot) \, , 
\end{aligned}
\ee
where $ \varrho_1 (\vphi) = \varrho_1^* (\vphi) $, 
and, recalling \eqref{new normal form D+}, \eqref{deco:M+M-},  the $ \vphi $-average 
\be\label{varrho1-NF}
[\widehat{ \varrho_1} ]^j_j(0) =
\frac{1}{(2\pi)^\es} \int_{\T^\es} ( \widehat{ \varrho_1} )^j_j( \vphi ) \, d \vphi
 \in M_- \, , \quad \ \forall \, j \in {\mathbb M} \, . 
\ee
The aim is to  solve the ``homological" equation\index{Homological equations}
\be \label{linhom-1}
J \ppamuvphi {\mathtt S}   + [J{\mathtt S} , J D_V ]= J \Pi_{\mathtt O}  \varrho  
\ee
which, recalling \eqref{def:homo-matrix-1}, \eqref{decomposition-rho1-1}, amounts to solve the decoupled pair of equations
\begin{align} \label{lineqd-0}
& J \ppamuvphi {\mathtt d}  + D_V {\mathtt d}  + J {\mathtt d} J D_V  = J \varrho_1  \\
&  \label{eqhoma-0}
J \ppamuvphi {\mathtt a} + D_V {\mathtt a} + J {\mathtt a} J D_V =J \varrho_2 \, . 
\end{align}
Note that, taking the adjoint equation of \eqref{eqhoma-0}, multiplying
by  $ J $ on the left and the right, and since $J $ and $ D_V $ commute, we obtain 
also  the equation in top right  in \eqref{linhom-1}, see \eqref{def:homo-matrix-1} and \eqref{decomposition-rho1-1}. 

The arguments of this section are similar  to those  developed in Section \ref{sec:homo-split},
actually simpler because the equation  \eqref{eqhoma-0} has constant coefficients in $ \vphi $, unlike the
corresponding equation  \eqref{eqhoma} where the operator  $ V_0 $ depends on $ \vphi \in \T^\es $. 
Thus in the sequel we shall often  refer to Section \ref{sec:homo-split}.

We first find a solution $ {\mathtt d} $ 
of the equation \eqref{lineqd-0}.
We recall that a linear operator
$ {\mathtt d}(\vphi ) \in {\cal L}( H_{\mathbb M}) $ is represented 
by a finite dimensional square matrix  $ ( {\mathtt d}_{i}^j (\vphi)  )_{i, j \in {\mathbb M}}  $ with  entries 
$ {\mathtt d}_{i}^j (\vphi) \in {\cal L}( H_j, H_i) \simeq {\rm Mat}_2 (\R) $.
To solve \eqref{lineqd-0} we use
the second order Melnikov non-resonance conditions  \eqref{2Mel+}-\eqref{2Mel rafforzate}, that
depend just on the unperturbed linear frequencies defined by \eqref{auto-funzioni}. 

\begin{lemma} \label{homdiag-1} {\bf (Homological equation \eqref{lineqd-0})}
Assume 
the second order Melnikov non resonance conditions \eqref{2Mel+}-\eqref{2Mel rafforzate}.
Then the equation  \eqref{lineqd-0} has a solution $ {\mathtt d} (\vphi ) = ( {\mathtt d}_i^j (\vphi) )_{i, j \in {\mathbb M}} $,
$ {\mathtt d}(\vphi ) = {\mathtt d}^*(\vphi ) $, satisfying
\be\label{loss-su-d}
\| {\mathtt d}_i^j \|_{\Lip, H^s(\T^\es)} \leq C  \| (\varrho_1)_i^j \|_{\Lip, H^{s+ 2 \tau_0}(\T^\es)} \, , \quad \forall i,j \in {\mathbb M} \, . 
\ee
\end{lemma}

\begin{pf}
Since the symplectic operator $ J $ leaves invariant each subspace $ H_{j} $ and recalling 
\eqref{def:DV},  the equation \eqref{lineqd-0} is equivalent to 
$$
J \ppamuvphi {\mathtt d}_i^j (\vphi) + 
\mu_i  {\mathtt d}_i^j (\vphi)  + \mu_j J {\mathtt d}_i^j (\vphi) J = J (\varrho_1 (\vphi) )_i^j \, , \  \forall i,j \in {\mathbb M} \, ,
\   \  
J = 
\begin{pmatrix}
 0 & 1   \\
 -1 & 0   \\
\end{pmatrix} \, ,
$$
and, by a Fourier series expansion with respect to the variable $ \vphi \in \T^\es $ writing
$$
d_i^j (\vphi) = \sum_{\ell \in \Z^\es} \widehat{d}_i^j (\ell) e^{\ii \ell \cdot \vphi} \, , \quad
\widehat{d}_i^j (\ell) \in {\rm Mat}_2 (\C ) \,, \quad \ov{\widehat{d}_i^j (\ell)} = \widehat{d}_i^j (-\ell) \, , 
$$
to
\be \label{lineqd3-1}
\ii (\bar \mu \cdot \ell ) J {  \widehat{\mathtt d}}_i^j (\ell) + \mu_i \widehat{  {\mathtt d}}_i^j (\ell) + 
\mu_j  J \widehat{  {\mathtt d}}_i^j (\ell) J 
= J [\widehat{ \varrho_1}]_i^j (\ell) \, ,  \
\forall i,j \in {\mathbb M} \, ,  \  \ell \in \Z^\es \, . 
\ee
Using the second order Melnikov non resonance conditions \eqref{2Mel+}-\eqref{2Mel rafforzate},
and since, by \eqref{varrho1-NF}, the Fourier coefficients 
$ [\widehat{ \varrho_1}]_j^j (0)  \in M_- $, 
the equations \eqref{lineqd3-1} can be solved  arguing as in Lemma \ref{homdiag} and
the estimate \eqref{loss-su-d} follows by standard arguments.
\end{pf}

We now solve  the equation \eqref{eqhoma-0} in the unknown  $ {\mathtt a} \in {\cal L}(H_{\mathbb M} , H_{\mathbb M}^\bot) $.
We use again the second order Melnikov non-resonance conditions \eqref{2Mel+}-\eqref{2Mel rafforzate}.
 
\begin{lemma} \label{homdiag-2-1} {\bf (Homological equation \eqref{eqhoma-0})}
Assume 
the second order Melnikov non resonance conditions \eqref{2Mel+}-\eqref{2Mel rafforzate}.
Then the homological  
equation  \eqref{eqhoma-0} has a solution $ {\mathtt a} \in 
L^2 (\T^\es, {\cal L}(H_{\mathbb M} , H_{\mathbb M}^\bot)) $ 
satisfying
\be\label{stima:a-av1}
|{\mathtt a} |_{\Lip, s} \leq C(s)  \| \varrho_2 \|_{\Lip, s+ 2 \tau_0} \, . 
\ee
\end{lemma}

\begin{pf}
Writing $ {\mathtt a}^j (\vphi) := {\mathtt a (\vphi)}_{|H_{j}} $,  
$ \varrho_2^j (\vphi) := (\varrho_2 (\vphi))_{|H_j} \in {\cal L}(H_j, H_{\mathbb M}^\bot) $ 
 and recalling  \eqref{def:DV}, the equation \eqref{eqhoma-0} amounts to 
\be \label{eqhomai-zero}
J \ppamuvphi {\mathtt a}^j  (\vphi) + D_V  {\mathtt a}^j (\vphi) + 
\mu_j  J {\mathtt a}^j  (\vphi) J  = J \varrho_2^j (\vphi)  \, , \quad \forall j \in {\mathbb M} \, . 
\ee
Writing, by  a Fourier series expansion with respect to $ \vphi \in \T^\es $, 
$$
\begin{aligned}
& {\mathtt a}^j  (\vphi)  = \sum_{\ell \in \Z^\es} {\widehat {\mathtt a}^j}  (\ell) e^{\ii \ell \cdot \vphi } \, , 
\ {\widehat {\mathtt a}^j}  (\ell) = \frac{1}{(2\pi)^\es} \int_{\T^\es} {\mathtt a}^j  (\vphi ) e^{-\ii \ell \cdot \vphi} \, d \vphi \, , \\
& 
\varrho_2^j (\vphi)   = \sum_{\ell \in \Z^\es} J {\widehat \varrho_2^j}  (\ell) e^{\ii \ell \cdot \vphi } \, , 
\end{aligned}
$$
the equation \eqref{eqhomai-zero} amounts  to 
\be \label{eqhomai-0}
\ii  {\bar \mu} \cdot \ell  J {\widehat {\mathtt a}^j}  (\ell) + 
D_V  {\widehat {\mathtt a}^j} (\ell) + \mu_j  J {\widehat {\mathtt a}^j} (\ell) J  = 
J {\widehat \varrho_2^j} (\ell)  \, , \quad \forall j \in {\mathbb M} \, ,  \  \ell \in \Z^\es  \, . 
\ee
According to the $ L^2 $-orthogonal splitting $H_{\mathbb M}^\bot = \oplus_{j \in {\mathbb M}^c} H_j $
the linear operator 
$ {\widehat {\mathtt a}^j}  (\ell) $, which maps $ H_j $ into the complexification 
of 
$ H_{\mathbb M}^\bot $ and satisfies 
$ \ov{{\widehat {\mathtt a}^j}  (\ell)}=  {\widehat {\mathtt a}^j}  (-\ell) $, 
is identified (as in \eqref{a1-a4}-\eqref{per-lemma7.2} with 
index $ k \in {\mathbb M}^c $) 
with a sequence of $ 2 \times 2 $ matrices
\be\label{seque-ajk-reality}
( {\widehat {\mathtt a}^j_k} (\ell) )_{ k \in {\mathbb M}^c } \, , 
\quad 
{\widehat {\mathtt a}^j_k} (\ell) \in 
{\rm Mat}_2 ( \C) \, , 
\quad
\ov{{\widehat {\mathtt a}^j_k} (\ell)} = {\widehat {\mathtt a}^j_k} (-\ell) \, .
\ee
Similarly  $ {\widehat \varrho_2^j} (\ell) \equiv ({[\widehat \varrho_2]_k^j} (\ell))_{k \in {\mathbb M}^c } $. Thus 
\eqref{eqhomai-0}
amounts to the following sequence of equations
\be \label{eqhomai-0-ink}
\begin{aligned}
& \ii  {\bar \mu} \cdot \ell  J {\widehat {\mathtt a}^j_k}  (\ell) + 
\mu_k  {\widehat {\mathtt a}^j_k} (\ell) + \mu_j  J {\widehat {\mathtt a}^j_k} (\ell) J  = 
J {[\widehat \varrho_2]_k^j} (\ell)  \, , \\
&  \quad j \in {\mathbb M} \, ,  \ k \in {\mathbb M}^c \, , \  \ell \in \Z^\es \, , 
\    J = 
\begin{pmatrix}
 0 & 1   \\
 -1 & 0   \\
\end{pmatrix}
 \, . 
 \end{aligned}
\ee
Note that the equation \eqref{eqhomai-0-ink} is like \eqref{lineqd3-1}. Since 
$ j \neq k $ (indeed $ j \in {\mathbb M} $, $ k \in {\mathbb M}^c $),
by the second order Melnikov non resonance conditions \eqref{2Mel+}-\eqref{2Mel rafforzate},  
each equation \eqref{eqhomai-0-ink} has a unique solution for any $ \varrho_2 $, 
the reality condition \eqref{seque-ajk-reality} holds,   and 
\be\label{small-div-av2}
\|{\widehat {\mathtt a}^j_k}  (\ell)\| \leq C  \langle \ell \rangle^{\t_0} \| {[\widehat \varrho_2]_k^j} (\ell) \| \, . 
\ee
We now estimate $ |{\mathtt a} |_s $. By Lemma \ref{Lemma:rho} we have 
$ |{\mathtt a} |_s  \simeq_s \| {\mathtt a} \|_s $
(we  identify each $  {\mathtt a}^j  \in  {\cal L}(H_j, H_{\mathbb M}^\bot) $
with a function  of 
$ H_{{\mathbb M}}^\bot \times H_{{\mathbb M}}^\bot $ as in 
\eqref{ident-H-4-me}, \eqref{ident-H-4}).
Given a function 
$$
 u (\vphi, x ) =  \sum_{\ell \in \Z^\es, k \in {\mathbb M}^c}  u_{\ell, k} e^{\ii \ell \cdot \vphi }\Psi_k (x)  \in H_{\mathbb M}^\bot 
 $$  
the Sobolev norm $ \| u \|_s $ defined in \eqref{def:Hs} is equivalent, using  \eqref{equiv-norms-s}, \eqref{spaces-cal-Hs},  
to 
\be\label{equiv-in-HM}
\begin{aligned}
\| u \|_s^2 & \simeq_s \| u \|_{L^2_\vphi (H^s_x \cap  H_{\mathbb M}^\bot)}^2 + 
\| u \|_{H^s_\vphi (L^2_x \cap  H_{\mathbb M}^\bot)}^2 \\ 
& \simeq_s
\sum_{\ell \in \Z^\es, k \in {\mathbb M}^c} 
\big( \mu_k^{2s} + \langle \ell \rangle^{2s}  \big)  | u_{\ell, k}|^2 \, . 
\end{aligned}
\ee
In conclusion \eqref{equiv-in-HM}  and \eqref{small-div-av2} and the fact that by Young inequality
$$
\begin{aligned}
& \langle \ell \rangle^{2\tau_0}  \mu_k^{2s} 
 \leq \frac{\langle \ell \rangle^{2\tau_0 p }}{p}   + \frac{\mu_k^{2s q}}{q} 
 \lesssim_{s, \tau_0} 
  \langle \ell \rangle^{2(\tau_0+s)}   + \mu_k^{2(\tau_0+s)} \, , \\
& \qquad \qquad p := \frac{2(\tau_0+s)}{2 \tau_0} \, ,
 \ q := \frac{2(\tau_0+s)}{2 s} \, , 
 \end{aligned}
 $$
 imply 
$ |{\mathtt a} |_{s} \leq C(s)  \| \varrho_2 \|_{s+ \tau_0} $. 
The estimate  \eqref{stima:a-av1} for the Lipschitz norm follows as usual. 
\end{pf}

\section{Averaging step} 

We consider the family of invertible symplectic transformations
\be\label{def:P0}
{\mathtt P} (\vphi) := e^{J {\mathtt S}(\vphi) } \, , \quad  {\mathtt P}^{-1} (\vphi) := e^{- J {\mathtt S}(\vphi) } \, , 
\quad \vphi \in \T^\es \, , 
\ee
where $ {\mathtt S} := {\mathtt S} (\vphi) $ is the self-adjoint operator in 
$ {\cal L} (H_{\mathbb S}^\bot) $  of the form \eqref{forma cal-S0} with $ {\mathtt d} (\vphi) $, $ {\mathtt a}(\vphi)  $  defined in Lemmata  \ref{homdiag-1}, \ref{homdiag-2-1}.
By Lemma \ref{Lemma:rho}, 
the estimates \eqref{loss-su-d}, \eqref{stima:a-av1}
and \eqref{est:varrho} 
 imply
\be\label{stima:S}
| {\mathtt S} |_{\Lip, +, s} \lesssim_s | {\mathtt S} |_{\Lip, s + \frac12} 
\lesssim_s  \| \varrho\|_{\Lip, s+ 2 \tau_0 + 1}  \leq 
C(s) \e^2  
\ee
and 
the transformation $ {\mathtt P} (\vphi) $ in  \eqref{def:P0} satisfies, for $ \e $ small,  the estimates 
\be\label{tame-P0}
| {\mathtt P} |_{\Lip, s_1} \leq 2 \, , \quad  | {\mathtt P} |_{\Lip, s} \, , \ | {\mathtt P}^{-1} |_{\Lip,s}  \leq C(s) 
\, , \quad \forall s \geq s_1 \, . 
\ee 
In the next proposition, which is the main result of this chapter,  
we conjugate the whole operator 
$ \bar \om_\e \cdot \partial_\vphi - J ( {\mathtt A} +  \frac{{\mathtt r}_\e }{1+\e^2 \lambda}) $
 defined in \eqref{def:primo-op}  by $ {\mathtt P}(\vphi ) $. 

\begin{proposition}\label{prop:op-averaged} {\bf (Averaging) }
Assume 
the second order Melnikov non resonance conditions \eqref{2Mel+}-\eqref{2Mel rafforzate}
where the set $ {\mathbb M} $ is fixed in Lemma \ref{choice:M}. 
Let $ {\mathtt P} (\vphi) $ be the  symplectic transformation  \eqref{def:P0} of $ {H_{\St}^\bot} $, 
where $ {\mathtt S} (\vphi) $ is the self-adjoint operator of the form 
\eqref{forma cal-S0} with $ {\mathtt d} (\vphi) $, $ {\mathtt a}(\vphi)  $  defined in Lemmata  \ref{homdiag-1}, \ref{homdiag-2-1}.
Then the conjugated operator 
\be\label{transf-op-1}
 {\mathtt P}^{-1} (\vphi) \Big[ \ppavphi - J \Big( {\mathtt A} + 
\frac{{\mathtt r}_\e}{ 1+\e^2 \lambda }  \Big) \Big]
 {\mathtt P} (\vphi)  =  \ppavphi - J {\mathtt A}^+ 
\ee
acting in the normal subspace $ {H_{\St}^\bot} $, 
has the following form, with respect to the splitting $ H_{\mathbb S}^\bot = H_{\mathbb M} \oplus H_{\mathbb M}^\bot $, 
\be\label{newA+}
{\mathtt A}^+ = {\mathtt A}_0 + \varrho^+ \, , \quad 
{\mathtt A}_0 := 
\frac{D_V}{1+ \e^2 \l} + \Pi_{\mathtt D} \varrho 
= 
\begin{pmatrix}
D(\e,\l)  & 0 \\ 
0 & {\mathtt A}_{{\mathbb M}^c}^{{\mathbb M}^c} 
\end{pmatrix}   \, ,  
\ee
where, in the basis  $ \{ (\Psi_j,0),(0,\Psi_j) \}_{j \in \mathbb M}$ of $H_{\mathbb M}$, 
the operator $ D(\e,\l) $ is represented by the diagonal matrix 
\begin{align} \label{shifted-normal-freq}
D(\e,\l) = 
{\rm Diag}_{j \in {\mathbb M}} {\mathtt \Omega}_j ( \e, \l)  {\rm Id}_2 \, , \quad 
{\rm Id}_2 = 
\begin{pmatrix}
1  & 0   \\
 0 & 1  
\end{pmatrix} 
\end{align}
with $ {\mathtt \Omega}_j (\e, \l ) $ defined in \eqref{def:Omega-n}, 
$ \mathtt A $ is defined in \eqref{def:primo-op},  
 and
\be\label{est:rho+}
 | \varrho^+ |_{\Lip, +,s}  \lesssim_{s}   \e^4 + | {\mathtt r}_\e |_{\Lip, +,s}   \, .
\ee
Moreover, given  another self-adjoint operator $ {\mathtt r}_\e' $  satisfying \eqref{estimate:rep-ge}, we have that 
\be\label{est:rho+-variation}
 | \varrho^+ - (\varrho^+)' |_{+,s_1}  \lesssim_{s_1}   | {\mathtt r}_\e - {\mathtt r}_\e' |_{+,s_1}   \, .
\ee
\end{proposition}

The rest of this chapter is dedicated to prove Proposition \ref{prop:op-averaged}.
We first study the conjugated operator 
$$ 
{\mathtt P}^{-1} (\vphi) ( \ppavphi - J{\mathtt A} ) {\mathtt P} (\vphi)  
$$ 
where $ {\mathtt A} $ is defined in \eqref{def:primo-op}.
We have the Lie series expansion\index{Lie expansion}
\begin{align} \label{Lieseries-1}
{\mathtt P}^{-1} (\vphi) ( \ppavphi - J{\mathtt A}  ) {\mathtt P} (\vphi) & =
(\ppavphi - J {\mathtt A}) 
+ {\rm Ad}_{(-J{\mathtt S})} ({\mathtt X}_0) + \sum_{k \geq 2} \frac{1}{k!} 
{\rm Ad}^k_{(-J{\mathtt S})} ({\mathtt X}_0) 
\end{align}
where  $ {\mathtt X}_0 := \ppavphi - J {\mathtt A}  $.
Recalling that $ {\mathtt A} = (1+ \e^2 \l )^{-1} D_V + \varrho $
by \eqref{def:primo-op}-\eqref{def:varrho0},  
we expand the commutator as
\begin{align}
{\rm Ad}_{(-J {\mathtt S})} (X_0) 
&=  \ppavphi (J {\mathtt S} ) + [J{\mathtt S},J {\mathtt A}]    \label{sviluppo-Ad-1}  \\
&  = 
J (\ppavphi {\mathtt S}) + (1+ \e^2 \l )^{-1} [J{\mathtt S},J  D_V] + [J{\mathtt S}, J\varrho] 
\nonumber \\
& =  J (\ppamuvphi {\mathtt S}) +   [J{\mathtt S},J  D_V] 
- \e^2 \l (1+ \e^2 \l )^{-1} [J{\mathtt S},J  D_V] + [J{\mathtt S}, J\varrho] \nonumber \\
& \quad + J (\bar \om_\e - \bar \mu) \! \cdot \! \pa_\vphi {\mathtt S} 
\nonumber \\
& \stackrel{\eqref{linhom-1}, \eqref{def omep}} =  J \Pi_{\mathtt O} \varrho 
- \e^2 \l (1+ \e^2 \l )^{-1} [J{\mathtt S},J  D_V]  + [J{\mathtt S}, J\varrho] + 
\e^2 (J \zeta  \! \cdot \! \pa_\vphi) {\mathtt S}  \, . \nonumber 
\end{align}
As a consequence of  \eqref{Lieseries-1}, \eqref{sviluppo-Ad-1}, 
\eqref{def:varrho}, \eqref{def:primo-op}-\eqref{def:varrho0} we deduce \eqref{transf-op-1} with 
\be\label{newrho-1}
 {\mathtt A}^+    = 
{\mathtt A_0}  + \varrho^+ \, , \quad {\mathtt A_0} :=  (1+ \e^2 \l)^{-1} D_V + \Pi_{\mathtt D} \varrho  
\ee
and 
\begin{align}
- J \varrho^+  & := - \e^2 \l (1+ \e^2 \l )^{-1} [J{\mathtt S},J  D_V]  + 
 [J{\mathtt S}, J\varrho] + \e^2 (J \zeta \! \cdot \! \pa_\vphi) {\mathtt S} \\
& \quad \ - (1 + \e^2 \l )^{-1} {\mathtt P}^{-1} J {\mathtt r}_\e {\mathtt P} + \sum_{k \geq 2} \frac{1}{k!} 
{\rm Ad}^k_{(-J{\mathtt S})} ( {\mathtt X}_0)  \, . \nonumber 
\end{align}
The estimate \eqref{est:rho+} for $ \varrho^+ $ follows by \eqref{stima:S}, \eqref{tame-P0},  
\eqref{est:varrho} and the same arguments used in Lemmata \ref{Ads1}-\ref{AdS}.  
Similarly we deduce \eqref{est:rho+-variation}.

Now, recalling  \eqref{def:pro+}-\eqref{new normal form D+}  (with $ {\mathbb F} \rightsquigarrow {\mathbb M} $,
$ {\mathbb G} \rightsquigarrow {\mathbb M}^c $), 
and \eqref{def:primo-op}-\eqref{def:varrho0}, 
we have that the operator $ {\mathtt A}_0 $ in \eqref{newrho-1} can be decomposed,
with respect to the splitting  $ H_{\mathbb M} \oplus H_{{\mathbb M}^c} $,  as
$$ 
{\mathtt A}_0  = \frac{D_V}{1+ \e^2 \l} + \Pi_{\mathtt D} \varrho  = 
 \begin{pmatrix}
D(\e,\l)  & 0 \\ 
0 &  {\mathtt A}_{{\mathbb M}^c}^{{\mathbb M}^c}  
\end{pmatrix} \, , 
$$
where, 
taking in $ H_{\mathbb M} $ the basis 
$ \{ ( \Psi_j(x) ,0) \, ,  (0, \Psi_j(x) ) \}_{j \in {\mathbb M}} $, we have 
\begin{align} \label{D-computed}
D(\e,\l) & 
= \frac{ [D_V]^{\mathbb M}_{\mathbb M}}{1 + \e^2 \l } 
+ 
{\rm Diag}_{j \in {\mathbb M}} \big(  \pi_+ [ {\widehat {\varrho}_j^j}(0)] \big) \nonumber \\
& 
 \stackrel{\eqref{def:varrho0}} = 
\frac{1}{1+  \e^2 \l } {\rm Diag}_{j \in {\mathbb M}} \begin{pmatrix}
 { \mu_j}  & 0   \\
 0 & { \mu_j}    
\end{pmatrix}  + \e^2 
{\rm Diag}_{j \in {\mathbb M}} \Big( \frac{\pi_+ [ {\widehat{\mathtt B}_j^j}(0) ]}{1+  \e^2 \l} \Big) 
\, . 
\end{align}
We now prove that $ D(\e,\l) $ has the form  \eqref{shifted-normal-freq}. 

\begin{lemma} {\bf (Shifted normal frequencies)}\index{Shifted normal frequencies} \label{lem:ns}
The operator $ D(\e, \l) $  in \eqref{D-computed} has the form  \eqref{shifted-normal-freq}
with $ {\mathtt \Omega}_j ( \e, \l) $  defined in \eqref{def:Omega-n}. 
 \end{lemma}

\begin{pf}
By \eqref{D-computed} and recalling the definition of $ \pi_+ $ in  \eqref{proiettore:sim}, the  operator $ D(\e, \l) $ is 
\be\label{corrected-normal-fre}
D(\e, \l) = 
{\rm Diag}_{j \in {\mathbb M}} \frac{\mu_j + \e^2 {\mathtt b}_j }{1+  \e^2 \l }  {\rm Id}_2 \, , 
\quad 
{\mathtt b}_j := \frac12 {\rm Tr} ({\widehat{\mathtt B}_j^j} (0) )  \, ,  
\ee
and, by the definition of $ \mathtt B $  in \eqref{form-of-B-ge}, we have 
\begin{align}\label{new-rj+}
{\rm Tr} ( {\widehat {\mathtt B}_j^j }(0)) 
& =  \frac{3}{(2 \pi)^{|{\mathbb S}|}}\int_{\T^\es}
\big( \Psi_j, D_V^{-1/2} \big( a(x) ( v(\vphi, 0, \xi) )^2   D_V^{-1/2} \Psi_j \big) \big)_{L^2 (\T^d)} \, d \vphi \, , \\
& +  \frac{4 \e}{(2 \pi)^{|{\mathbb S}|}}\int_{\T^\es}
\big( \Psi_j, D_V^{-1/2} \big( a_4(x) ( v(\vphi, 0, \xi) )^3   D_V^{-1/2} \Psi_j \big) \big)_{L^2 (\T^d)} \, d \vphi \, . 
\label{new-rj+-bis}
\end{align}
We first compute the term \eqref{new-rj+}.
Expanding the expression of  $ v $ in \eqref{defv} we get 
\begin{align}
\eqref{new-rj+} & = \frac{3}{(2 \pi)^{|{\mathbb S}|}}  \int_{\T^{\es+d}} ( D_V^{-1/2} \Psi_j (x) )   a(x) ( v(\vphi, 0, \xi) )^2   
D_V^{-1/2} \Psi_j (x) \, d \vphi dx  \nonumber \\
  & =   \frac{3}{(2 \pi)^{|{\mathbb S}|}}  \int_{\T^{\es+d}} \frac{\Psi_j (x) }{\sqrt{\mu_j}  } 
 a(x) \Big( \sum_{k \in \mathbb S}  \mu_k^{-1/2} \sqrt{2 \xi_k } \cos \vphi_k \Psi_k (x)   \Big)^2  
 \frac{\Psi_j (x) }{\sqrt{\mu_j}}   \, d \vphi dx \nonumber \\
  & =   \frac{3}{(2 \pi)^{|{\mathbb S}|}} \sum_{k_1, k_2 \in \mathbb S}   \int_{\T^{\es+d}} 
 a(x) \ \mu_{k_1}^{-1/2} \mu_{k_2}^{-1/2} \sqrt{2 \xi_{k_1} } \sqrt{2 \xi_{k_2} } 
 \cos \vphi_{k_1}\cos \vphi_{k_2} \nonumber \\ 
& \qquad\qquad\qquad\qquad\qquad\qquad\qquad
\qquad\qquad \times \Psi_{k_1} (x) \Psi_{k_2}(x)   \frac{\Psi_j^2 (x)}{\mu_j}   \, d \vphi dx \nonumber \\
 & = \frac{3}{(2 \pi)^{|{\mathbb S}|}} \sum_{k \in \mathbb S}   \int_{\T^d} 
 a(x) \ \mu_{k}^{-1}  2 \xi_{k}   \Psi_{k}^2 (x)  \Psi_j^2 (x) \mu_j^{-1}    dx \int_{\T^{|\mathbb S|} } \cos^2 \vphi_{k}\, d \vphi 
 \nonumber \\
 & = \frac{6}{(2 \pi)^{|{\mathbb S}|}}    \mu_j^{-1}  
 \sum_{k \in \mathbb S}  \mu_{k}^{-1}  \xi_{k}  \int_{\T^d}  a(x) \  
 \Psi_{k}^2 (x)  \Psi_j^2 (x) dx \, (2 \pi )^{|{\mathbb S}|-1} \, \int_{\T } \cos^2 \teta \, d \teta \nonumber
 \\
 & =  3  \mu_j^{-1} 
 \sum_{k \in \mathbb S}    \mu_{k}^{-1}    
 ( \Psi_j^2, a(x)  \Psi_{k}^2 )_{L^2} 
   \xi_{k} = 2 (\Bb \xi)_j  \label{formula-aj}
 \end{align}
using \eqref{def:AB}-\eqref{def G}.  
On the other hand, the term in \eqref{new-rj+-bis} is equal to zero, because
\begin{align}
 \eqref{new-rj+-bis}
  & =   \frac{4 \e}{(2 \pi)^{|{\mathbb S}|}}  \int_{\T^{d}} \int_{\T^{\es}} \frac{\Psi_j (x) }{\sqrt{\mu_j}  } 
 a_4 (x) \big( v(\vphi, 0, \xi)  \big)^3  
 \frac{\Psi_j (x) }{\sqrt{\mu_j}}   \, d \vphi dx  \label{termine-secondo-zero}
 \end{align}
and,  by \eqref{simmetry:v}, the integral 
$$
\int_{\T^\es} \big( v(\vphi, 0, \xi)  \big)^3   d \vphi = 
\int_{\T^\es} \big( v(\vphi + \vec \pi, 0, \xi)  \big)^3 d \vphi = - \int_{\T^\es} \big( v(\vphi, 0, \xi)  \big)^3   d \vphi \, ,
$$
is equal to zero. 

In conclusion, we deduce by \eqref{corrected-normal-fre}, \eqref{new-rj+},  
\eqref{new-rj+-bis}, \eqref{formula-aj}, \eqref{termine-secondo-zero} that 
$ {\mathtt b}_j   
 = (\Bb \xi)_j  $, and, inserting  the value 
$ \xi := \xi (\l) = \e^{-2} \Ab^{-1} ( ( 1+ \e^2 \lambda) \bar \omega_\e - \bar \mu ) $ 
defined in \eqref{xil}, we get that the eigenvalues of $ D(\e, \lambda ) $ in \eqref{corrected-normal-fre}, are equal to 
$$
 \frac{\mu_j + \e^2 {\mathtt b}_j }{1+  \e^2 \l } = {\mathtt \Omega}_j (\e, \l)  
$$
with  $ {\mathtt \Omega}_j (\e, \l)  $ defined in \eqref{def:Omega-n}. 
 \end{pf}
 
The new operator $  \ppavphi - J {\mathtt A}^+ $
obtained in Proposition \ref{prop:op-averaged} 
is in a suitable form to admit, for most values of the parameter  $ \l  $, 
an approximate right inverse according to Proposition \ref{prop-cruciale} in the next chapter.

\chapter{Approximate right inverse in  normal directions}\label{sec:7}

The goal of this chapter is to state 
the crucial Proposition \ref{prop-cruciale} for the existence of an approximate right inverse  
for a class of quasi-periodic linear Hamiltonian operators, 
acting in the normal subspace $ H_{\mathbb S}^\bot $,  of  the form
$\bar \om_\e \cdot \partial_\vphi - J  (A_0  + \rho) $
where $ A_0 $ and $\rho $ are self-adjoint operators, 
 $ A_0 $ is a split admissible  operator 
 according to 
 Definition \ref{def:calC}, and 
$ \rho $ is  ``small", see \eqref{rho-R0:small}.
We shall use Proposition \ref{prop-cruciale} to construct the sequence of
approximate solutions along  
the iterative nonlinear Nash-Moser scheme of Chapter \ref{sec:NM}, more precisely to prove Proposition \ref{prop:inv-ap-vero}.

\section{Split admissible operators}

We first define the following class of  split admissible operators $ A_0 $. 

\begin{definition} \label{def:calC}
{\bf (Split admissible operators)}\index{Split admissible operator}
Let $ C_1, c_1, c_2 >  0 $ be constants. We denote by $ {\cal C}(C_1,c_1,c_2)$ the class of 
self-adjoint operators  
\be\label{defA0}
A_0(\e, \lambda , \varphi)= \frac{D_V}{1+  \e^2 \l} + R_0(\e, \l, \varphi) \, , \qquad 
D_V = (-\Delta +V(x) )^{1/2} \, , 
\ee
acting on $ H_{\mathbb S}^\bot $,  defined for all $ \l \in \wtilde \Lambda \subset \Lambda $, 
that satisfy  
\begin{enumerate}
\item \label{item1R0C1e2}
$| R_0  |_{\Lip, +, s_1} \leq C_1 \e^2 $, 
\item  $A_0$ is block diagonal 
with respect to the splitting $ H_{\mathbb S}^\bot  = H_{\mathbb F} \oplus H_{{\mathbb G}}  $, i.e. $ A_0 $ has the form  (see \eqref{decoFGsotto}) 
\be\label{form:A0}
A_0 := A_0( \e , \l , \vphi )= 
\begin{pmatrix}
D_0(\e , \lambda) & 0  \\
0  &  V_0(\e, \l, \varphi)  
\end{pmatrix}  \, , 
\ee
and, moreover,  
in the basis of the eigenfunctions $ \{ (\Psi_j,0) $, $(0,\Psi_j)\}_{ j \in {\mathbb F}} $ (see \eqref{Fj-eigenfunctions}),  
the operator $ D_0 $ is represented by the diagonal matrix
\be\label{form:D0}
D_0 := D_0 (\e, \lambda) = {\rm Diag}_{j \in {\mathbb F}} \, \mu_j (\e,\l) {\rm Id}_2 \, , \quad \mu_j (\e,\l) \in \R \, .
\ee
The eigenvalues $ \mu_j (\e, \l) $ satisfy 
\be\label{muje2}
|\mu_j ( \e, \lambda ) - \mu_j | \leq C_1 \e^2  \, , 
\ee 
where $ \mu_j  $ are defined in  \eqref{auto-funzioni},  and 
\begin{align}
& {\mathfrak d}_\l (\mu_i - \mu_j) (\e, \l) \geq c_2 \e^2 \quad {or} \quad {\mathfrak d}_\l (\mu_i - \mu_j) (\e, \l) \leq - c_2 \e^2
\, , \quad  i \neq  j   \, , \label{Hyp2} \\
&
{\mathfrak d}_\l (\mu_i + \mu_j) (\e, \l) \geq c_2 \e^2 \quad {or} \quad {\mathfrak d}_\l (\mu_i + \mu_j) (\e, \l) \leq - c_2 \e^2 \, , 
   \label{Hyp3}  \\
 & \forall  j \in {\mathbb F}\, , \ 
 \Big( c_2 \e^2 \leq {\mathfrak d}_\l \mu_j (\e, \l)  \leq c_2^{-1} \e^2 \  {\rm or} \ 
  - c_2^{-1} \e^2 \leq {\mathfrak d}_\l \mu_j (\e, \l)  \leq - c_2 \e^2 \Big)    \,   , \label{Hyp4} 
\end{align} 
and, in addition, for all $ j \in {\mathbb F} $, 
\be\label{Hyp1}
\begin{cases}
{\mathfrak d}_\l  \big( V_0 (\e , \l) + \mu_j (\e, \l) {\rm Id}  \big) \leq - c_1 \e^2 \  \cr
{\mathfrak d}_\l \big( V_0 (\e , \l) - \mu_j (\e, \l) {\rm Id} \big) \leq - c_1 \e^2 \, . 
 \end{cases}
\ee
\end{enumerate}
\end{definition}

We now verify that the operator $ {\mathtt A}_0 $ defined in \eqref{newA+}
is a split admissible operator
according to Definition \ref{def:calC}.

\begin{lemma}\label{A-0-splitted} {\bf (${\mathtt A}_0 $ is split admissible)}
There exist  positive constants $ C_1, c_1, c_2 $ such that 
the operator  $ {\mathtt A}_0 $ defined in  \eqref{newA+}-\eqref{shifted-normal-freq}  belongs to the class $ {\cal C}(C_1, c_1, c_2) $  of split admissible operators
introduced  in Definition  \ref{def:calC}.  Notice that $ {\mathtt A}_0 $ is defined for all $ \lambda \in \Lambda $. 
\end{lemma}

\begin{pf}
Remind that the decomposition in  \eqref{newA+} refers to the splitting 
$ H_{\mathbb S}^\bot  = H_{\mathbb M} \oplus H_{{\mathbb M}^c} $, where the finite set 
$ {\mathbb M} \supset {\mathbb F} $ has been fixed in 
Lemma \ref{choice:M}, whereas  the decomposition in \eqref{form:A0}  concerns 
the splitting 
$ H_{\mathbb S}^\bot  = H_{\mathbb F} \oplus H_{{\mathbb G}} $.

By \eqref{newA+} the operator $ {\mathtt A}_0 $ has the form \eqref{defA0}
with $ R_0 = \Pi_{\mathtt D} \varrho $ and \eqref{estimate-deco-resto},  \eqref{est:varrho} imply 
$$ 
| R_0  |_{\Lip, +, s_1}  = |  \Pi_{\mathtt D} \varrho  |_{\Lip, +, s_1} \leq C_1 \e^2 \, . 
$$ 
In addition, with respect to the splitting $  H_{\mathbb S}^\bot   = H_{\mathbb F} \oplus H_{{\mathbb G}} $,
and taking in  $ H_{\mathbb F} $ the basis of the eigenfunctions $ \{ (\Psi_j,0) $, $(0,\Psi_j)\}_{ j \in {\mathbb F}} $,
the operator $ {\mathtt A}_0 $ in \eqref{newA+}-\eqref{shifted-normal-freq}
has the form \eqref{form:A0}-\eqref{form:D0} (recall that $ {\mathbb F} \subset {\mathbb M} $)  
with 
$$
D_0 (\e, \lambda) = {\rm Diag}_{j \in {\mathbb F}} {\mathtt \Omega}_j (\e, \l)   {\rm Id}_2 \, , 
\quad \mu_j (\e,\l) = {\mathtt \Omega}_j (\e, \l)  \, ,
$$
and the operator $ V_0 $, which acts in $ H_{\mathbb G} $,  admits the decomposition, with respect to 
the splitting $ H_{\mathbb G} = H_{{\mathbb M} \setminus {\mathbb F}} \oplus H_{\mathbb M^c} $,
and taking in  $ H_{{\mathbb M} \setminus {\mathbb F}} $ 
the basis of the eigenfunctions $ \{ (\Psi_j,0) $, $(0,\Psi_j)\}_{ j \in {{\mathbb M} \setminus {\mathbb F}}} $,
\be\label{def:V0a}
V_0 = 
\begin{pmatrix}
{\rm Diag}_{j \in {\mathbb M} \setminus {\mathbb F}} {\mathtt \Omega}_j (\e, \l)  {\rm Id}_2 & 0  \\
 0  & {\mathtt A}_{{\mathbb M}^c}^{{\mathbb M}^c} \\
\end{pmatrix} \,  .
\ee
By 
\eqref{def:Omega-n} and \eqref{def omep} we have, for all $ j \in {\mathbb M} $,  
$ {\mathtt \Omega}_j (\e, \l)  = {\mathtt \Omega}_j (0, \lambda) + O(\e^2) $,  where
$ {\mathtt \Omega}_j (0, \l ) = \mu_j  $ are the unperturbed linear frequencies defined in \eqref{auto-funzioni}, 
and the derivative of the functions ${\mathtt \Omega}_j (\e, \l)  $ (which are defined for all $ \lambda \in \Lambda$)
is 
\begin{align}
\pa_\l {\mathtt \Omega}_j (\e, \l) =  - \frac{ \e^2 }{(1+ \e^2 \l)^2} 
( \mu_j - [ \Bb \Ab^{-1} \bar \mu]_j ) , \quad \forall j \in {\mathbb M} \label{der-pert-1passo} \, . 
\end{align}
Then, by \eqref{der-pert-1passo} and the 
definition of $ {\mathbb F}$ and $ {\mathbb G}$ 
in \eqref{taglio:pos-neg}, for any $ j \in {\mathbb M} \setminus {\mathbb F} \subset {\mathbb G} $, we have
\begin{align}
\pa_\l {\mathtt \Omega}_j  + \max_{j \in {\mathbb F}} | \pa_\lambda {\mathtt \Omega}_j |
& = - \frac{ \e^2 }{(1+ \e^2 \l)^2}  \Big( ( \mu_j - [ \Bb \Ab^{-1} \bar \mu]_j ) - 
\max_{j \in {\mathbb F}} | \mu_j -  ( \Bb \Ab^{-1} \bar \mu)_j| \Big) \nonumber \\
& <  - \frac{ \e^2 }{(1+ \e^2 \l)^2} \mathfrak c \label{eq:Omega-nor} 
\end{align}
for some $ \mathfrak c  >  0 $. 
By \eqref{eq:Omega-nor} (recall that $ \mu_j (\e,\l) = {\mathtt \Omega}_j (\e, \l)$) and \eqref{negative-on-up} 
we conclude that, for $ \e $ small,  
 property \eqref{Hyp1} holds for some $ c_1 >  0 $. 
The other properties \eqref{Hyp2}-\eqref{Hyp4} follow, for $ \e $ small,  
by   \eqref{der-pert-1passo} and the  assumptions \eqref{non-reso}-\eqref{non-reso1}. 
\end{pf}

\section{Approximate right inverse}
 
The main result of this chapter is the following proposition which provides an approximate 
solution of the linear equation 
$$
{\mathfrak L}   h = g  
\qquad {\rm where} \qquad 
{\mathfrak L}  := \bar {\om }_\e \cdot \partial_\varphi - J (A_0 + \rho) \, . 
$$

\begin{proposition}\label{prop-cruciale} 
{\bf (Approximate right inverse in normal directions)}\index{Approximate right inverse in normal directions}
Let  $ \bar \om_\e \in \R^\es $ be  $(\gamma_1 , \tau_1) $-Diophantine and satisfy 
property $ {\bf (NR)}_{\gamma_1, \tau_1} $ in Definition \ref{NRgamtau} with $ \g_1, \t_1 $ fixed in 
\eqref{def:tau1}. Fix $ s_1 > s_0 $ according to Proposition \ref{propmultiscale} (more precisely Proposition \ref{propinv}), 
and $ s_1 < s_2 < s_3 $ such that 
\be\label{choice s2 s3} 
(i) \ s_2-s_1 >  300 (\tau' +3 s_1 +3) \, , \qquad  (ii) \ s_3-s_1\geq  3 (s_2 -s_1) 
\ee
where the constant $\tau' $ appears in the multiscale Proposition \ref{propmultiscale} and 
we assume that $ \t' > 2 \t + 3 $. 

Then there is  $ \e_0 > 0  $  such that, $\forall \e \in (0, \e_0)$, for each self-adjoint operator 
$$ 
A_0 =  \frac{D_V}{1+ \e^2 \l} + R_0 
$$ 
acting in $ H_{\mathbb S}^\bot $, 
belonging to the class $ {\cal C} (C_1,c_1,c_2) $ of  split admissible operators (see Definition \ref{def:calC}), 
for any self-adjoint operator $ \rho  \in L^2( \T^\es, {\cal L}( H_{\mathbb S}^\bot  )) $, 
defined in $ \wtilde \Lambda \subset \Lambda $, 
satisfying
\be\label{rho-R0:small}
|  \rho  |_{\Lip, +, s_1}  \leq  \e^3 \, , \quad 
 |  R_0  |_{\Lip, +, s_2} + |  \rho |_{\Lip, +, s_2}   \leq  \e^{-1} \, , 
\ee
there  are closed subsets $ {\bf \Lambda} (\e; \cc ,A_0, \rho) \subset \wtilde \Lambda $,  $1/2 \leq \cc \leq 5/6$,  satisfying 
\begin{enumerate}
\item $  {\bf \Lambda}  (\e ; \cc, A_0, \rho) \subseteq 
{\bf \Lambda}  (\e ; \cc', A_0, \rho) $, for all $1/2 \leq \cc \leq \cc' \leq 5/6$, 
\item \label{AApro2}
$ \big| [ {\bf \Lambda} (\e; 1/2, A_0, \rho)]^c \cap \wtilde \Lambda \big| \leq b(\e) $ where $ \lim_{\e \to 0} b(\e)=0 $, 
\item \label{AApro3}
if $ A_0' = (1+ \e^2 \l)^{-1} D_V  + R_0' \in {\cal C} (C_1,c_1,c_2) $ and $\rho'$ satisfy  
\be\label{smallR0R0'}
|R_0' -R_0|_{+, s_1} + |\rho' - \rho|_{+, s_1} \leq \delta \leq \e^3 \, , 
\ee
for all $\l \in \wtilde{\Lambda} \cap \wtilde{\Lambda}'  \subset \Lambda$, 
then, for all $ (1/2) + \d^{2/5} \leq \cc \leq 5/6$, 
\be\label{finadav}
\big| \wtilde{\Lambda}' \cap  [{\bf \Lambda} (\e; \cc, A'_0, \rho')]^c \cap  {\bf \Lambda} (\e ; \cc - \d^{2/5} , A_0, \rho) \big| \leq \delta^{\alpha/3} ; 
\ee 
\end{enumerate}
and, for any $ \nu \in (0, \e) $, there exists a linear operator 
$$ 
{\mathfrak L}^{-1}_{approx} := {\mathfrak L}^{-1}_{approx, \nu} \in {\cal L} ({\mathcal H}^{s_3} \cap H_{\mathbb S}^\bot) 
$$ 
such that, 
for any function $ g : \wtilde \Lambda  \to  {\mathcal H}^{s_3} \cap H_{\mathbb S}^\bot $ satisfying 
\be\label{g:small-large}
\| g \|_{\Lip, s_1} \leq \e^2 \nu  \,  , \quad 
|  R_0  |_{\Lip, +, s_3}  + |\rho|_{\Lip, +, s_3}+ \| g \|_{\Lip, s_3}  \leq \e^2 \nu^{-1} \ ,
\ee
the function  $ h := {\mathfrak L}^{-1}_{approx}  g $,  $ h  :  {\bf \Lambda} (\e; 5/6, A_0, \rho ) 
 \to  {\mathcal H}^{s_3} \cap H_{\mathbb S}^\bot $  
satisfies  
\be\label{h:s1s2s3}
 \|h\|_{\Lip, s_1} \leq \e^2 \nu^{ \frac45} \, , \quad
  \|h\|_{\Lip, s_3} \leq  \e^2 \nu^{ - \frac{11}{10}} \, , 
\ee
and 
\be\label{sol:almost-inv}
(\bar {\om }_\e \cdot \partial_\varphi - J (A_0 + \rho)) h = g +  r 
\ee
with 
\be\label{estimate:r-final}
 \|r\|_{\Lip, s_1} \leq \e^2 \nu^{3/2} \, . 
\ee
Furthermore, setting $ Q' := 2(\tau'+ \varsigma s_1 ) +3 $ (where $ \varsigma = 1/ 10 $  and 
$ \tau' $ is given by Proposition \ref{propmultiscale}), for all $ g \in {\mathcal H}^{s_0 + Q'} \cap  H_{\mathbb S}^\bot $, 
\be\label{con-loss}
 \| {\mathfrak L}^{-1}_{approx} \, g \|_{\Lip, s_0} \lesssim_{s_1} \| g \|_{\Lip, s_0 + Q'} \, . 
\ee
\end{proposition}

The required bounds  in \eqref{rho-R0:small}-\eqref{g:small-large}
will be verified along the nonlinear Nash-Moser scheme of Chapter \ref{sec:NM}. 
Proposition \ref{prop-cruciale} will be  applied to the operator
$  \bar \om_\e \cdot \pa_{\vphi}  - J {\mathtt A}^+ $ where $ {\mathtt A}^+ = {\mathtt A}^0 + \varrho^+ $ is defined in Proposition \ref{prop:op-averaged}, to prove Proposition \ref{prop:inv-ap-vero}. Notice that 
$ {\mathtt A}^0 $ is split admissible 
by Lemma \ref{A-0-splitted} and  the coupling operator 
$ \varrho^+ $ satisfies  $ |  \varrho^+  |_{\Lip, +, s_1}  \lesssim_{s_1} \e^4 \leq \e^3 
$ by \eqref{est:rho+} and \eqref{estimate:rep-ge}.  

Proposition \ref{prop-cruciale} is proved in 
Chapter \ref{sec:proof.Almost-inv} using the results of  Chapters
\ref{sec:multiscale} and \ref{sec:splitting}.   

\smallskip

We remark that the value of $ \e_0 $ given in  Proposition \ref{prop-cruciale} may depend on $s_3$,
because, in the estimates of the proof, 
there are quantities of the form $ C(s_3) \e^a $, $ a > 0 $,
and we choose $ \e $ small so that $ C(s_3)  \e^a < 1 $. However,  such terms
appear by quantities  $ C(s_3)  \nu^{a_1}  $, $ a_1 > 0 $, and so we may require 
  $ C(s_3)  \nu^{a_1} < 1  $
assuming only $\nu$ small enough, 
i.e. $\nu \leq \wtilde{\nu} (s_3)$ (see Remark \ref{rem:nuep}). As a result, the following
more specific statement holds.

\begin{proposition}  \label{crucial-s}
The conclusion of  Proposition \ref{prop-cruciale}  can be modified from  \eqref{finadav}
in the following way.

For any $s \geq s_3$, there is $\wtilde{\nu} (s) >  0 $  such that:  for any $0<\nu < \min (\e, \wtilde{\nu} (s) )$, 
there exists a linear operator 
$$ 
{\mathfrak L}^{-1}_{approx} := {\mathfrak L}^{-1}_{approx, \nu, s} \in {\cal L} ({\mathcal H}^{s} \cap H_{\mathbb S}^\bot) 
$$ 
such that, 
for any function $ g : \wtilde \Lambda  \to  {\mathcal H}^{s} \cap H_{\mathbb S}^\bot $ satisfying 
\be\label{g:small-large-inf}
\| g \|_{\Lip, s_1} \leq \e^2 \nu  \,  , \quad 
|  R_0  |_{\Lip, +, s}  + |\rho|_{\Lip, +, s}+ \| g \|_{\Lip, s}  \leq \e^2 \nu^{-1} \ ,
\ee
the function  $ h := {\mathfrak L}^{-1}_{approx}  g $,  $ h  :  {\bf \Lambda} (\e; 5/6, A_0, \rho ) 
 \to  {\mathcal H}^{s} \cap H_{\mathbb S}^\bot $  
satisfies  
\be\label{h:s1s2s3-inf}
 \|h\|_{\Lip, s_1} \leq \e^2 \nu^{ \frac45} \, , \quad
  \|h\|_{\Lip, s} \leq  \e^2 \nu^{ - \frac{11}{10}} \, , 
\ee
and 
\be\label{sol:almost-inv-i}
(\bar {\om }_\e \cdot \partial_\varphi - J (A_0 + \rho)) h = g +  r 
\ee
with 
\be\label{estimate:r-final-inf}
 \|r\|_{\Lip, s_1} \leq \e^2 \nu^{3/2} \, . 
\ee
Furthermore, setting $ Q' := 2(\tau'+ \varsigma s_1 ) +3 $ (where $ \varsigma = 1/ 10 $  and 
$ \tau' $ is given by Proposition \ref{propmultiscale}), for all $ g \in {\mathcal H}^{s_0 + Q'} \cap  H_{\mathbb S}^\bot $, 
\be\label{con-loss-inf}
 \| {\mathfrak L}^{-1}_{approx} \, g \|_{\Lip, s_0} \lesssim_{s_1} \| g \|_{\Lip, s_0 + Q'} \, . 
\ee
\end{proposition}

\chapter{Splitting between low-high normal subspaces}\label{sec:splitting}

The main result of this chapter is Corollary \ref{cor:split} below. Its goal is to block-diagonalize
a quasi-periodic  Hamiltonian operator of the form
$ \bar \om_\e \cdot \partial_\vphi - J  (A_0 + \rho ) $ 
according to the splitting $ H_{\mathbb S}^\bot = H_{\mathbb F} \oplus H_{\mathbb G} $,  
up to a very small coupling term, 
see the conjugation  \eqref{goal-trasf-iterata} where $ \rho_\ind $ is small according to  \eqref{form-A'-m}.

\section{Splitting step and corollary}

 The proof of Corollary  \ref{cor:split} is based on an iterative application of the following 
Proposition. 


\begin{proposition} \label{prophomeq} {\bf (Splitting step)}
Let  $ \bar \om_\e \in \R^\es $ be  $(\gamma_1 , \tau_1)$-Diophantine and satisfy 
property $ {\bf (NR)}_{\gamma_1, \tau_1} $ in Definition \ref{NRgamtau} 
with $ \g_1, \t_1 $ fixed in  \eqref{def:tau1}.
Assume 
$ s_2 - s_1 >  120 (\tau' +3 s_1 +3) $, 
i.e. condition \eqref{choice s2 s3}-($i$). 
Then, given $C_1>0$, $c_1>0$, $c_2>0$, there is  $ \e_0 > 0  $  such that, $\forall \e \in (0, \e_0)$, for each self-adjoint operator 
$$ 
A_0 
=  \frac{D_V}{1+ \e^2 \l} + R_0 =  
\begin{pmatrix}
{D}_0(\e , \lambda) & 0  \\
0  &  V_0(\e, \l, \varphi)  
\end{pmatrix} 
$$ 
belonging to the class $ {\cal C} (C_1,c_1,c_2) $ of split admissible operators of Definition \ref{def:calC} (see \eqref{defA0}-\eqref{form:A0}), 
defined in $ \wtilde \Lambda $,  
with  
$$
D_0 (\e, \lambda) = {\rm Diag}_{j \in {\mathbb F}} \, \mu_j (\e,\l) {\rm Id}_2 \, , \quad
\mu_j (\e,\l) \in \R \, , 
$$ 
as in \eqref{form:D0}, there 
are 
\begin{itemize}
\item
closed subsets $ \Lambda (\e; \cc ,A_0) \subset \wtilde \Lambda $,  $1/2 \leq \cc \leq 1$,  satisfying the  properties
\begin{enumerate}
\item $  \Lambda  (\e ; \cc, A_0) \subseteq  \Lambda  (\e ; \cc', A_0) $ 
for all  $1/2 \leq \cc \leq \cc' \leq 1$, 
\item \label{sp:item2}
$ \big| [ \Lambda (\e; 1/2, A_0)]^c \cap \wtilde \Lambda \big| \leq b(\e) $ 
where $ \lim_{\e \to 0} b(\e)=0 $, 
\item  \label{sp:item3}
if $ A_0' = (1+ \e^2 \l)^{-1} D_V  + R_0' \in {\cal C} (C_1,c_1,c_2) $ with $ |R'_0 - R_0|_{+, s_1} \leq \delta \leq \e^{5/2} $ for all
$\l \in \wtilde{\Lambda} \cap \wtilde \Lambda' $, 
then, for all $ (1/2) + \sqrt{\d} \leq \cc \leq 1$, 
\be\label{Cantor:intersect0}
\big| \wtilde{\Lambda}' \cap  [\Lambda (\e; \cc, A'_0)]^c \cap  \Lambda (\e ; \cc -\sqrt{\d} , A_0) \big| \leq \delta^\alpha \, , \quad \a > 0 \, ,  
\ee 
\end{enumerate}
\end{itemize}
and, for each self-adjoint operator 
$ \rho $ acting in $ H_{\mathbb S}^\bot $, defined for  $ \l \in \wtilde \Lambda $,  satisfying 
\begin{align}\label{rho-R0:small-0}
& |  \rho  |_{\Lip, +, s_1}   \leq \delta_1 \leq  \e^3 \, , \qquad   \delta_1 (|  R_0  |_{\Lip, +, s_2} 
 + |  \rho  |_{\Lip, +, s_2} )  \leq  \e \, , 
\end{align}
there exist
\begin{itemize}
\item 
 a symplectic linear  invertible transformation
\be\label{def:P}
e^{J {\cal S}  (\vphi)} \in {\cal L}( H_{\mathbb S}^\bot  ) \, , \quad \vphi \in \T^\es \, , 
\ee
defined for all $ \l \in \Lambda(\e;1,A_0)   $, where $ {\cal S}  (\vphi) := {\cal S}(\e, \l )  (\vphi)$ 
is a self-adjoint operator in $ {\cal L}( H_{\mathbb S}^\bot  ) $, satisfying 
the estimates \eqref{stima-S}-\eqref{estanys} below; 

\item 
a self-adjoint  operator $ A^+ $ of the form 
\be\label{new-A'}
A^+  =  
\frac{D_V}{1+ \e^2 \l } + R_0^+ + \rho^+  =
\begin{pmatrix}
 {D}^+_0(\e , \lambda) & 0  \\
0  &  V^+_0(\e, \l, \varphi)  
\end{pmatrix} + \rho^+  \, , 
\ee
defined for all $ \l \in \Lambda(\e;1,A_0)   $, with 
\be
\label{form-A'}
R^+_0, \rho^+ \  L^2{\it-self-adjoint} \, , \ \quad
D^+_0(\e, \l)={\rm Diag}_{j \in {\mathbb F}} \, \mu^+_j (\e,\l) {\rm Id}_2 \, , 
\ee
and $ \mu_j^+ (\e,\l) \in \R $;
\end{itemize} 
such that,  for all $ \l \in \Lambda(\e;1,A_0)  $, 
\be\label{goal-trasf}
\big( \bar \om_\e \cdot \partial_\vphi - J  A_0 - J \rho  \big) e^{J {\cal S} (\vphi)}  = 
e^{J {\cal S}  (\vphi)}   \big( \bar \om_\e \cdot \partial_\vphi - J A^+ \big) \, , 
\ee
and the following estimates hold: 
\begin{itemize}
\item The self-adjoint operator ${\cal S} $ satisfies 
\be \label{stima-S}
| {\cal S}|_{\Lip, s_1+1} \leq \d_1^{\frac78} \, , \
| {\cal S} |_{\Lip, s_2+1}
\leq \d_1^{- \frac14}\big( |R_0|_{\Lip, +,s_2}  + 
|  \rho  |_{\Lip, +, s_2} \big) + \d_1^{- \frac34} \, ,
\ee
and more generally, for all $s \geq s_2$,
\be \label{estanys}
  | {\cal S} |_{\Lip, s+1} 
\leq C(s) \Big[\d_1^{- \frac14}\big( |R_0|_{\Lip, +,s}  +  
|  \rho  |_{\Lip, +, s} \big) + \d_1^{- \frac34} \d_1^{- 3 \varsigma 	\frac{s-s_2}{ s_2 - s_1 }} \Big] \, ;
\ee
\item The operators in \eqref{new-A'}-\eqref{form-A'} satisfy 
\begin{align} 
&  \label{new-aut-V+}
|\mu_j^+ (\e, \l) -  \mu_j (\e, \l) |_{\Lip} \leq \d_1^{3/4} = o(\e^2) \, ,   \\
&  \|  V_0^+ -  V_0\|_{\Lip, 0} \leq \d_1^{3/4} = o (\e^2) \, , \label{new-aut-V+1}
\end{align}
and 
\begin{align}
|  R^+_0  - R_0  |_{\Lip, +, s_1} & \leq  \d_1^{3/4} \, , \quad |\rho^+  |_{\Lip, +, s_1} \leq  \d_1^{3/2} / 2  \, ,
\label{propo1} \\
| R^+_0 |_{\Lip, +, s_2}+ | \rho^+  |_{\Lip, +, s_2} 
& \leq \d_1^{-1/4} \big(| R_0  |_{\Lip, +, s_2}+ | \rho |_{\Lip, +, s_2} \big)
+ \d_1^{-3/4} \label{propo2-0}  \\
&  \stackrel{\eqref{rho-R0:small-0}} \leq \e \delta_1^{-3/2} \, , \label{propo2} 
  \end{align}
 and, more generally, for all  $ s \geq s_2 $,
\be\label{R0+s}
 | R^+_0  |_{\Lip, +, s}+ | \rho^+  |_{\Lip, +, s} \lesssim_s  
 \d_1^{-\frac14} \big(| R_0 |_{\Lip, +, s}+ | \rho |_{\Lip, +, s} \big)
+ \d_1^{- \frac34} \d_1^{- 3 \varsigma 	\frac{s-s_2}{ s_2 - s_1 }} 
\ee 
where $\varsigma = 1 /10 $ is fixed in \eqref{def:varsigma}. 
\end{itemize}
At last, given $ A'_0 = (1+ \e ^2 \l)^{-1} D_V + R'_0 \in {\cal C}(C_1,c_1,c_2)$ and $\rho'$ satisfying 
\eqref{rho-R0:small-0}, 
then, for all $\l \in \Lambda (\e; 1 , A_0) \cap \Lambda (\e; 1 , A'_0)  $, we have the estimates
\begin{align} \label{pertR0+}
& |{R'_0}^+ - R_0^+ |_{+, s_1} \leq |R'_0-R_0|_{+, s_1} + C |\rho' - \rho|_{+, s_1} \\
& \label{pertrho+}
|{\rho'}^+ - \rho^+|_{+, s_1} \leq   \d_1^{1/2} |R'_0-R_0|_{+, s_1} + \d_1^{-1/20}  |\rho' - \rho|_{+, s_1}  \, . 
\end{align}
\end{proposition}

Note that, according to \eqref{goal-trasf}, \eqref{new-A'}, \eqref{propo1},   
the new coupling term $  \rho^+  $ is much smaller than  $ \rho $, 
in low norm  $ | \  |_{\Lip, +, s_1} $. 
Moreover, by \eqref{propo2},  the new $ \rho^+ $ satisfies 
also the second assumption \eqref{rho-R0:small-0} with $ \d_1^{3/2} $ instead of $ \d_1 $. Notice also that the Cantor sets 
$  \Lambda (\e; \cc ,A_0)  $ depend only on $ A_0 $ and not on $ \rho $. 
Finally we point out that  \eqref{pertR0+}-\eqref{pertrho+} 
will be used for the measure estimate of  Cantor sets of ``good'' 
parameters $\lambda$, in
relation with property \ref{sp:item3} (see \eqref{Cantor:intersect0}): for this application,
an estimates of low norms $ | \ |_{+,s_1} $, without the control of the Lipschitz dependence,
is enough.

Applying iteratively Proposition \ref{prophomeq}  we deduce the following corollary.


\begin{corollary} {\bf (Splitting)} \label{cor:split}
Let  $ \bar \om_\e \in \R^\es $ be  $(\gamma_1, \tau_1)$-Diophantine and satisfy 
property $ {\bf (NR)}_{\gamma_1, \tau_1} $ in Definition \ref{NRgamtau}
with $ \g_1, \t_1 $ fixed in \eqref{def:tau1}. 
Let 
$$ 
A_0 = \frac{D_V}{1+ \e^2 \l} + R_0 
$$ 
be a self-adjoint operator in  the class of split admissible  operators $ {\cal C} (C_1,c_1,c_2) $
(see Definition \ref{def:calC}) and $ \rho  $ be a 
self-adjoint operator in $ {\cal L}( H_{\mathbb S}^\bot ) $, defined for 
$ \l \in \wtilde \Lambda $, satisfying \eqref{rho-R0:small}. 
Then there exist
\begin{itemize}
\item 
closed sets  $ \Lambda_\infty (\e; \eta,  A_0, \rho) \subset \wtilde \Lambda $, $ 1/2 \leq \eta \leq 5/6$,  satisfying the properties 
\begin{enumerate}
\item  \label{item2-91} $  \Lambda_\infty  (\e ; \cc, A_0, \rho) \subseteq  \Lambda_\infty  (\e ; \cc', A_0, \rho) $
 for all $1/2 \leq \cc \leq \cc' \leq 5/6 $, 
\item \label{item2-92}
$ \big| \big[ \Lambda_\infty (\e; 1/2, A_0, \rho) \big]^c \cap \wtilde \Lambda  \big| \leq b_1 (\e) $ 
where $ \lim_{\e \to 0} b_1(\e)=0 $, 
\item \label{item2-93}
if $ A_0' = (1+ \e^2 \l)^{-1} D_V  + R_0' \in {\cal C} (C_1,c_1,c_2) $ and $\rho'$ satisfy  
\be\label{R0R0'vicini}
|R_0-R'_0|_{+, s_1} + |\rho - \rho'|_{+, s_1} \leq \delta \leq \e^2 
\ee
for all $ \l \in \wtilde{\Lambda} \cap  \wtilde{\Lambda}' $, 
then, for all $ (1/2) + \d^{2/5} \leq \cc \leq 5/6$, 
\be\label{Cantor:intersect}
\big| \wtilde{\Lambda}' \cap  [\Lambda_\infty (\e; \cc, A'_0, \rho')]^c \cap  \Lambda_\infty 
(\e ; \cc - \d^{2/5} , A_0, \rho) \big| \leq \delta^{\alpha/2} ; 
\ee 
\end{enumerate}
\item 
a sequence of symplectic linear invertible transformations 
\be\label{def:Pind}
{\cal P}_{0} := {\rm Id} \, , \quad {\cal P}_{\ind}  := {\cal P}_\ind (\e, \l) (\vphi)= 
e^{J {\cal S}_1 (\e, \l) (\vphi)} \ldots e^{J {\cal S}_\ind (\e, \l) (\vphi)} \, , \ \  \ind \geq 1 \, , 
\ee
defined for all $ \l \in \Lambda_\infty (\e; 5/6, A_0,\rho)  $, acting in $ H_{\mathbb S}^\bot $, satisfying, for 
\be\label{def:d1}
\d_1 = \e^3 \, , 
\ee
the estimates 
\begin{align}\label{Pins1}
& {\rm for} \ \ind \geq 1 , \, 
|  {\cal P}_\ind^{\pm 1}  -  {\cal P}_{\ind-1}^{\pm 1} |_{\Lip, +,s_1}  
\leq \dnew^{ (\frac32)^{\ind-1} \frac34}, \, 
 |  {\cal P}_\ind^{\pm 1}  - {\rm Id} |_{\Lip, +,s_1} \leq 2 \dnew^{\frac34}\, , 
\\ 
& |  {\cal P}_\ind^{\pm 1} |_{\Lip, +, s_2}  
 \leq (C(s_2))^\ind   \dnew^{-(3/2)^{\ind -1} \frac{3}{4}} \big[ \e
\dnew^{-\frac{1}{2} } +1 \big]  \, ,\label{Pins2}
\end{align}
and, more generally, for all $ s \geq s_2 $, 
\be \label{Pinsany} 
 |  {\cal P}_\ind^{\pm 1} |_{\Lip, +,s}    \leq 
 (C(s))^\ind   \d_1^{-(\frac32)^{\ind -1} (\frac{3}{4}+ \alpha (s) )} \big[ 
 (|R_0|_{\Lip, +,s} + |\rho |_{\Lip, +,s}) \d_1^{\frac{1}{2} + \frac{2\alpha (s) }{3}} +1 \big]  
\ee
where
\be\label{def:as}
 \a (s) := 3 \varsigma \frac{s-s_2}{s_2-s_1} \, ;
\ee
\item 
a sequence of self-adjoint block diagonal operators of the form 
\be\label{defAm}
A_\ind  =  \frac{D_V}{1+ \e^2 \l}  + R_\ind   
:=  \begin{pmatrix}
 D_\ind (\e , \lambda) & 0  \\
0  &  V_\ind (\e, \l, \varphi)  
\end{pmatrix}  \, , \ \  \ind \geq 1 \, , 
\ee
defined for $ \l \in \Lambda_\infty (\e; 5/6, A_0,\rho) $, 
belonging to the class $ {\cal C} (2C_1,c_1/2, c_2/2) $ of  split admissible operators, with 
\begin{align}\label{form-A'-m0}
D_\ind (\e, \l)={\rm Diag}_{j \in {\mathbb F}} \, \mu_j^{(\ind)} (\e,\l) {\rm Id}_2	\, , 
\end{align}
satisfying 
\begin{align}
&  \label{der-lambda-m}
 | \mu_j^{(\ind)} (\e, \l) -  \mu_j (\e, \l)  |_{\Lip} =O(\d_1^{3/4}) = o(\e^2) \, , \\
 & 
 \|  V_\ind -  V_0\|_{\Lip, 0} = O(\d_1^{3/4})= o (\e^2) \, ,  \label{der-lambda-mV} \\
& |  R_\ind   |_{\Lip, +, s_1}  \leq C_1 \e^2 +  2\d_1^{3/4} \leq 2C_1 \e^2 \, , \label{stime-rm-low} \\ 
& | R_\ind  |_{\Lip, +, s_2}  \ll (C(s_2))^\ind   \d_1^{-(\frac32)^{\ind -1} \frac{3}{4} } 
\big[ \e \d_1^{- \frac{1}{2} } +1 \big] \, , \label{stime-rm-high}
\end{align} 
and, more generally, for all $ s \geq s_2 $, 
\be\label{Rrn-tame0}
 | R_\ind  |_{\Lip, +, s}  \leq  
 (C(s))^\ind   \d_1^{-(\frac32)^{\ind -1} (\frac{3}{4}+ \alpha (s) )} 
 \big[ \big( |R_0|_{\Lip, +,s} + |\rho|_{\Lip, +,s} \big)
\d_1^{\frac{1}{2} + \frac{2\alpha (s) }{3}} +1 \big] 
\ee
where $ \a (s) $ is defined in \eqref{def:as}; 
\item a sequence of $ L^2 $ self-adjoint operators  
$ \rho_\ind \in {\cal L}( H_{\mathbb S}^\bot ) $, $  \ind \geq 1 $, 
defined for $ \l \in \Lambda_\infty (\e; 5/6, A_0,\rho)  $, satisfying  
\begin{align}\label{form-A'-m}
| \rho_\ind  |_{\Lip, +, s_1}   \leq \d_1^{(\frac32)^\ind} \, , \quad
 | \rho_\ind |_{\Lip, +, s_2}   \ll  (C(s_2))^\ind   \d_1^{-(\frac32)^{\ind -1} \frac{3}{4} } [\e \d_1^{- \frac12} +1]\, , 
\end{align} 
and, more generally, for all $ s \geq s_2 $, 
\be\label{Rrn-tame}
  | \rho_\ind  |_{\Lip, +, s}  \leq  
 (C(s))^\ind   \d_1^{-(\frac32)^{\ind -1} (\frac{3}{4}+ \alpha (s) )} \big[ (|R_0|_{\Lip, +,s} + |\rho|_{\Lip, +,s})
\d_1^{\frac{1}{2} + \frac{2\alpha (s) }{3}} +1 \big]  
\ee
where $ \a (s) $ is defined in \eqref{def:as}; 
\end{itemize}
such that, for all $ \l \in \Lambda_\infty (\e; 5/6, A_0,\rho) $,  
\be\label{goal-trasf-iterata}
\big( \bar \om_\e \cdot \partial_\vphi - J  A_0  -J \rho  \big) {\cal P}_\ind (\vphi) = 
{\cal P}_\ind (\vphi) \big( \bar \om_\e \cdot \partial_\vphi - J A_\ind -J \rho_\ind \big) \, . 
\ee
At last, given 
$$
A'_0= \frac{D_V}{1+\e^2 \l} + R'_0 \in {\cal C}(C_1, c_1, c_2)
$$ 
and a
self-adjoint operator $ \rho' \in  {\cal L}( H_{\mathbb S}^\bot ) $ satisfying 
\eqref{rho-R0:small}, if 
\be \label{pertA0}
|A'_0 - A_0|_{+, s_1} + |\rho'-\rho|_{+, s_1} \leq \d \leq \e^3 \, , \quad \forall \l \in \wtilde \Lambda \cap \wtilde \Lambda' \, ,  
\ee
then,  for all $  \l \in \Lambda_\infty(\e; 5/6, A_0, \rho) \cap \Lambda_\infty(\e; 5/6, A'_0, \rho')  $, for all $ \ind \in \N $,     
\be \label{pertAn}
|A'_\ind - A_\ind|_{+, s_1} \leq \d^{4/5} \, . 
\ee
\end{corollary}

Note that each operator 
$ A_\ind $ in \eqref{defAm} is block-diagonal according to the splitting 
$ H_{\mathbb S}^\bot = H_{\mathbb F} \oplus H_{\mathbb G} $ in \eqref{orthogonal-deco-FS}, i.e. it
 has the same form as $ A_0 $ in \eqref{form:A0} but  the coupling term 
$ \rho_\ind $ in \eqref{goal-trasf-iterata} is much smaller than $ \rho $ (in the low norm $ | \ |_{\Lip, +, s_1} $), compare 
the first inequality in \eqref{form-A'-m}  (where $ \d_1 = \e^3 $ by \eqref{def:d1})
and the first inequality in \eqref{rho-R0:small}.  
The tame estimates \eqref{Pinsany}  and \eqref{Rrn-tame} for all $ s \geq s_2 $, 
will be used in the proof of Proposition \ref{prop-cruciale},
which provides an approximate right  inverse required  in the Nash-Moser nonlinear iteration in Chapter \ref{sec:NM}.   

\smallskip

\begin{pf}  
Let us define the sequences of  real numbers 
$(\d_\ind)_{\ind \geq 1}$ and $(\eta_\ind)_{\ind \geq 0}$ by 
\be\label{def:indices}
\d_\ind := \d_1^{(\frac32)^{\ind-1}}, \ \d_1 = \e^3 \, ,   \quad {\rm and} \quad 
\eta_0 := 0 \, , \   \eta_{\ind+1} := \eta_\ind + \d_{\ind+1}^{ \frac38} =	\eta_\ind + \d_1^{ \frac38 (\frac32)^{\ind}} \, . 
\ee
We shall prove by induction the following statements: for any $ \ind \in \N $ 
\begin{itemize}
\item $ \bf (P)_\ind $ there exist: 
\begin{itemize}
\item[(i)] 
symplectic linear invertible transformations 
$ {\cal P}_0, \ldots, {\cal P}_\ind $ of the form \eqref{def:Pind},
defined respectively for $ \l $ in  decreasing subsets 
$$
 \Lambda_\ind  (\e; 5/6, A_0, \rho) \subset \ldots 
 \subset  \Lambda_1  (\e; 5/6, A_0, \rho) \subset \Lambda_0 := \wtilde  \Lambda \, , 
$$
satisfying  \eqref{Pins1}-\eqref{Pinsany} at any order $ k = 0, \ldots, \ind $.   
The sets $\Lambda_\ind (\e; \cc, A_0, \rho) $, $1/2 \leq \eta \leq 5/6 $, 
are defined
inductively by 
\be\label{set-cor-split}
\Lambda_0 := \wtilde \Lambda  \quad {\rm and} \quad  
 \Lambda_\ind  (\e; \eta, A_0, \rho) := \bigcap_{k=0}^{\ind -1}  \Lambda (\e ; \eta + \eta_{k}, A_{k}) \, , 
 \ \ind \geq 1 \, ,
\ee
where
the sets $ \Lambda (\e ; \eta + \eta_{k}, A_{k}) $
 are those defined by  Proposition \ref{prophomeq}. Notice that, for $ \e $ small enough, 
 $ \eta_k \leq 1/ 6 $ for any $ k \geq 0 $ (and so $ \eta + \eta_{k} \leq 1 $). 
\item[(ii)] 
Self-adjoint split admissible  operators $ A_0,  \ldots, $ $ A_\ind  $ 
as in \eqref{defAm}-\eqref{form-A'-m0},  in the class $ {\cal C}(2C_1,c_1/2, c_2/2) $ (see Definition \ref{def:calC}), 
satisfying  \eqref{der-lambda-m}-\eqref{Rrn-tame0} at any order $ k = 0, \ldots, \ind $, and 
such that the conjugation identity  \eqref{goal-trasf-iterata} holds for all $ \lambda $ in  $ 
\Lambda_\ind (\e; 5/6, A_0, \rho) $, with $\rho_k$ satisfying  \eqref{form-A'-m}-\eqref{Rrn-tame}
for $ k = 0, \ldots, \ind $. 
\item[(iii)] 
Moreover  we have, $ \forall \ind \geq 1 $,   
\be\label{An-An-1} 
 |A_{\ind} - A_{\ind -1}|_{\Lip, +, s_1}  \leq  \d_\ind^{\frac34} = \d_1^{\frac34 (\frac32)^{\ind-1}}
 \  {\rm on} \  \Lambda_\ind (\e; 5/6, A_0, \rho) \, . 
\ee
\end{itemize}
\end{itemize}
{\bf Initialization.} 
The statement $ {\bf (P)}_0$-(i) holds
with  $ {\cal P}_0 := {\rm Id} $, and \eqref{Pins1}-\eqref{Pinsany} trivially hold.
The conjugation identity  \eqref{goal-trasf-iterata} at $ \ind = 0 $ trivially holds with 
$ \rho_0 := \rho $. 
Then, in order to prove $ {\bf (P)}_0$-(ii), it is sufficient to notice that 
the self-adjoint operator $ A_0 \in  {\cal C} (C_1,c_1,c_2)$ 
 has the form \eqref{defAm}, \eqref{form-A'-m0} with 
 $ \mu_j^{(0)} (\e, \lambda) =  \mu_j (\e, \lambda) $, 
 and \eqref{der-lambda-m}-
 \eqref{stime-rm-low} hold.
The estimate \eqref{stime-rm-high} (which for $\ind=0$ is
$ |R_0|_{\Lip,+,s_2} \ll \d_1^{-1/2} [\e \d_1^{-1/2} +1]$) 
and \eqref{form-A'-m} are consequences of Assumption \eqref{rho-R0:small}, since $\d_1=\e^3$.  
Finally notice that \eqref{Rrn-tame0} and \eqref{Rrn-tame} are  trivially satisfied.
\\[1mm]
{\bf Induction.} 
Next assume that $ {\bf (P)}_\ind $ holds.   
In order to define 
$ {\cal P}_{\ind +1} $ and $ A_{\ind + 1 } $
we  apply the ``splitting step" Proposition \ref{prophomeq}  
with $ A_0, \rho $ replaced by  $A_\ind , \rho_\ind $.
In fact, by the inductive assumption  $ \bf (P)_\ind $, 
the operator $A_\ind $  in \eqref{defAm} 
belongs to  the class
 $ {\cal C}(2C_1,c_1/2, c_2/2)$ of split admissible  operators, according to Definition \ref{def:calC}.
Moreover, by \eqref{form-A'-m}, \eqref{stime-rm-high},  
we have 
\be\label{smallness-n}
| \rho_\ind |_{\Lip, +, s_1}   \leq 
 \d_1^{(\frac32)^\ind}   \stackrel{\eqref{def:indices}} = \d_{\ind+1} \, , \qquad 
\d_{\ind+1}  \big( |R_\ind|_{\Lip,+,s_2} + |\rho_\ind   |_{\Lip,+,s_2} \big) \leq \e \, ,
\ee
which is \eqref{rho-R0:small-0} with  $ \rho_\ind $, $ R_\ind $, $  \d_{\ind+1} $ instead of 
$  \rho $, $ R_0 $, $  \d_1 $.

We define, for $1/2 \leq \eta \leq 5/6$,  the set 
\be\label{def:Lambdan+1}
 \Lambda_{\ind +1}(\e; \eta, A_0, \rho)  :=  \Lambda_{\ind } (\e; \eta, A_0, \rho)  \cap  \Lambda \big( \e ; \cc + \eta_{\ind}, A_{\ind} \big) \, . 
\ee
in agreement  with \eqref{set-cor-split} at $ \ind + 1 $.

By \eqref{smallness-n}, Proposition \ref{prophomeq} implies the existence of a self-adjoint operator 
$ {\cal S}_{\ind+1} \in {\cal L}( H_{\mathbb S}^\bot ) $, defined for 
$\l \in \Lambda_{\ind +1}(\e; 5/6, A_0, \rho) \subset  \Lambda \big( \e ; 1, A_{\ind} \big)$, satisfying (see \eqref{stima-S}-\eqref{estanys}) 
\be
\begin{aligned} \label{estimateS-n+1}
& | {\cal S}_{\ind +1}|_{\Lip,s_1 +1} \leq \d_{\ind +1}^{\frac78}=\d_1^{\frac{7}{8} (\frac{3}{2})^\ind} \, , \\
& | {\cal S}_{\ind +1}|_{\Lip,s_2 +1} \leq  \d_{\ind +1}^{-\frac{1}{4} } 
 \big( |R_\ind |_{\Lip,+,s_2} + |\rho_\ind  |_{\Lip,+,s_2} \big) + \d_{\ind +1}^{-\frac{3}{4}}  \, , \\
 & | {\cal S}_{\ind +1} |_{\Lip, s+1} 
\leq C(s) \Big[\d_{\ind +1}^{- \frac14}\big( |R_\ind|_{\Lip, +,s}  +  |  \rho_\ind  |_{\Lip, +, s} \big) + \d_{\ind +1}^{- \frac34} \d_{\ind +1}^{- 3 \varsigma 	\frac{s-s_2}{ s_2 - s_1 }} \Big] \, ,
 \end{aligned}
 \ee
such that, for any $ \l \in \Lambda_{\ind +1}(\e; 5/6, A_0, \rho)$,  
we have 
\be\label{ind-conj}
\big( \bar \om_\e \cdot \partial_\vphi - J  A_\ind  -J \rho_\ind  \big) e^{J {\cal S}_{\ind +1}} = 
e^{J {\cal S}_{\ind +1}} \big( \bar \om_\e \cdot \partial_\vphi - J A_{\ind+1}  -J \rho_{\ind+1} \big) 
\ee
where the operator 
\be\label{def:An+1}
A_{\ind+1} =\dps \frac{D_V}{1+\e^2 \l} + R_{\ind +1} 
\ee
is block-diagonal as in \eqref{defAm}-\eqref{form-A'-m0}.
By  \eqref{propo1}-\eqref{propo2}, \eqref{R0+s}, we have
\be  \label{Rrhon+1}
\begin{aligned}   
|  R_{\ind+1}  - R_\ind  |_{\Lip, +, s_1} & \leq  \d_{\ind+1}^{3/4}=\d_1^{\frac34 (\frac32)^\ind} \, , \\
 |\rho_{\ind +1}  |_{\Lip, +, s_1} &\leq  \frac12 \, \d_1^{(\frac{3}{2})^{\ind +1}}  \, ,  \\
 | R_{\ind +1} |_{\Lip, +, s_2}+ | \rho_{\ind +1}  |_{\Lip, +, s_2} 
& \leq \d_{\ind +1}^{-\frac14} \big(| R_\ind  |_{\Lip, +, s_2}+ | \rho_\ind |_{\Lip, +, s_2} \big)
+ \d_{\ind +1}^{-\frac34}   \, ,  \\
 | R_{\ind +1}  |_{\Lip, +, s}+ | \rho_{\ind +1}  |_{\Lip, +, s} &\lesssim_s  
 \d_{\ind +1}^{- \frac14} \big(| R_\ind |_{\Lip, +, s}+ | \rho_\ind |_{\Lip, +, s} \big)
+ \d_{\ind +1}^{-\frac34} \d_{\ind +1}^{- 3 \varsigma 	\frac{s-s_2}{ s_2 - s_1 }}  \, . 
\end{aligned} 
\ee
In particular, we have the first bound in \eqref{form-A'-m} at the step $  \ind + 1 $.
\\[1mm]
{\sc The symplectic transformation  $ {\cal P}_{\ind +1} $.} 
We define, for $\l \in \Lambda_{\ind +1}(\e; 5/6, A_0, \rho) $,
the  symplectic linear invertible transformation  
\be\label{def:Pind+1}
{\cal P}_{\ind +1} := {\cal P}_{\ind }e^{J {\cal S}_{\ind +1}} \stackrel{\eqref{def:Pind}} = 
e^{J {\cal S}_{1}} \ldots  e^{J {\cal S}_{\ind }}  e^{J {\cal S}_{\ind +1}}  
\ee
which has the  form \eqref{def:Pind} at order $ \ind + 1 $.
By the inductive assumption $ \bf (P)_\ind $, the conjugation identity 
\eqref{goal-trasf-iterata} holds  for all $ \l \in \Lambda_\ind (\e; \eta, A_0, \rho)$, 
and we deduce, by \eqref{ind-conj}, that, for all $ \lambda $ in the set 
$\Lambda_{\ind+1} (\e; \eta, A_0, \rho) $ defined in \eqref{def:Lambdan+1}, we have 
\begin{align*}
\big( \bar \om_\e \cdot \partial_\vphi - J  A_0  -J \rho  \big) {\cal P}_{\ind +1} & \stackrel{\eqref{def:Pind+1} } =
\big( \bar \om_\e \cdot \partial_\vphi - J  A_0  -J \rho  \big) {\cal P}_{\ind} e^{J {\cal S}_{\ind +1}} \\
& \stackrel{\eqref{goal-trasf-iterata}} = 
{\cal P}_{\ind} \big( \bar \om_\e \cdot \partial_\vphi - J  A_\ind  -J \rho_\ind  \big)  e^{J {\cal S}_{\ind +1}} \\
& \stackrel{\eqref{ind-conj}} = 
{\cal P}_{\ind} e^{J {\cal S}_{\ind +1}} \big( \bar \om_\e \cdot \partial_\vphi - J  A_{\ind +1} -J \rho_{\ind +1}  \big) \\
&  \stackrel{\eqref{def:Pind+1} } = {\cal P}_{\ind +1 }  \big( \bar \om_\e \cdot \partial_\vphi - J  A_{\ind +1} -J \rho_{\ind +1}  \big) 
\end{align*}
which is  \eqref{goal-trasf-iterata} at the step $ \ind + 1 $.

We have $ {\cal P}_{\ind +1} -  {\cal P}_{\ind } = {\cal P}_{\ind } ( e^{J {\cal S}_{\ind +1}}- {\rm Id})   $. By 
\eqref{estimateS-n+1}, 
using \eqref{special-form}-\eqref{special-form-Lip}  and the definition of $ \d_{\ind +1} $ in \eqref{def:indices},
\be \label{ejsn+1}
|e^{J {\cal S}_{\ind +1}}-{\rm Id}|_{\Lip,+,s_1 } \lesssim |J {\cal S}_{\ind +1}|_{\Lip,+,s_1}  
\lesssim |J {\cal S}_{\ind +1}|_{\Lip,s_1 +1}  
\lesssim \d_{\ind +1}^{7/8}   \, .
\ee
Moreover, by the second inequality of \eqref{Pins1}, $ |{\cal P}_{\ind}|_{\Lip,+, s_1 } \leq 2 $.
Hence  
$$ 
| {\cal P}_{\ind +1} -  {\cal P}_{\ind }|_{\Lip,+,s_1} \leq \d_{\ind+1}^{3/4} 
$$ 
which is the first
inequality in \eqref{Pins1} at the step $ \ind + 1 $ (we obtain in the same way the estimate for 
$ {\cal P}_{\ind +1}^{-1} $). As a consequence,
$$
|{\cal P}_{\ind +1}-{\rm Id}|_{\Lip,+, s_1 } \leq \sum_{k = 1}^{\ind+1} \d_k^{3/4}  \leq 2  \d_1^{3/4}  
$$
which is the second inequality in \eqref{Pins1} at the step $ \ind +1 $ (we obtain in the same way the estimate for 
$ {\cal P}_{\ind +1}^{-1} $). 
Estimates  \eqref{Pins2}-\eqref{Pinsany} at  the step $ \ind +1 $
are proved below. 
\\[1mm]
{\sc $A_{\ind+1} $ in \eqref{def:An+1}  
 is a split admissible operator  in $ {\cal C}(2C_1,c_1/2, c_2/2)$}, 
 see Definition \ref{def:calC}.
By  \eqref{Rrhon+1} we have  
$$
|A_{\ind+1} - A_{\ind}|_{\Lip,+,s_1}  = |R_{\ind+1} - R_{\ind}|_{\Lip,+,s_1}    \leq \d_{\ind+1}^{\frac34} = \d_1^{\frac34 (\frac32)^\ind}   
\quad {\rm on} \  \ \Lambda_{\ind +1} (\e; 5/6, A_0, \rho) \, , 
$$
which is \eqref{An-An-1} at the step  $ \ind + 1 $. As a consequence
\be\label{stime-rm-low-n+1}
| R_{\ind+1} - R_{0}|_{\Lip,+,s_1}  = 
|A_{\ind+1}-A_{0}|_{\Lip,+,s_1} \leq \sum_{k=1}^{\ind +1} \d_1^{\frac{3}{4} (\frac{3}{2})^{k-1} }  \leq 2 \d_1^{3/4} 
\ee
and  \eqref{stime-rm-low} at order $\ind +1$ follows:
in fact, $ |R_0|_{\Lip, +,s_1} \leq C_1 \e^2 $  by item \ref{item1R0C1e2} of Definition \ref{def:calC} and
$$
|R_{\ind +1}|_{\Lip, +,s_1} \leq |R_0|_{\Lip, +,s_1} + 2 \d_1^{3/4} 
\leq C_1 \e^2 + 2 \e^{9/4} \leq 2C_1 \e^2 
$$
for $ \e $ small enough, since $\d_1 = \e^3 $.  Recalling 
\eqref{defAm}, the estimate 
\eqref{der-lambda-mV} at order $\ind +1$ is also 
a direct consequence of \eqref{stime-rm-low-n+1}, as well as 
$$ 
\|D_{\ind +1} -D_0 \|_{\Lip, 0} = O(\d_1^{3/4} )=o(\e^2) \, , 
$$ 
which implies 
\eqref{der-lambda-m} at order $\ind +1$. 

Now, since $A_0 \in {\cal C}(C_1,c_1,c_2)$ (Definition \ref{def:calC}), 
either ${\mathfrak d}_\lambda (\mu_i-\mu_j) (\e,\l) \geq c_2 \e^2$ 
or ${\mathfrak d}_\lambda (\mu_i-\mu_j) (\e,\l) \leq -c_2 \e^2$, see \eqref{Hyp2}. 
Thus, by \eqref{der-lambda-m} at order $\ind +1$, we deduce that
${\mathfrak d}_\lambda (\mu^{(\ind +1)}_i-\mu^{(\ind +1)}_j) (\e,\l) \geq c_2 \e^2/2$ in the first case and
${\mathfrak d}_\lambda (\mu^{(\ind +1)}_i-\mu^{(\ind +1)}_j) (\e,\l) \leq -c_2 \e^2/2$  in the second case,
for $\e$ small enough.

In a similar way  \eqref{der-lambda-m} provides properties \eqref{Hyp3}-\eqref{Hyp4} with constant $c_2/2$ instead of $c_2$ at order
$\ind +1$, and  
 \eqref{der-lambda-m}-\eqref{der-lambda-mV} provide \eqref{Hyp1} with constant $c_1/2$ instead of $c_1$ at order
$\ind +1$.
\\[1mm]
{\sc Estimates of $|R_{\ind +1}|_{\Lip,+,s}$ and $|\rho_{\ind +1}|_{\Lip,+,s}$.} 
We first consider the case  $ s = s_2 $. By \eqref{Rrhon+1} and  \eqref{stime-rm-high},
\eqref{form-A'-m}, and recalling the definition of $ \d_{\ind+1} $ in \eqref{def:indices}, 
\begin{align*}
|R_{\ind +1}|_{\Lip,+,s_2} +|\rho_{\ind +1}|_{\Lip,+,s_2} &\leq 
 \d_{\ind +1}^{- \frac14} \big( |R_{\ind }|_{\Lip,+,s_2} +|\rho_{\ind }|_{\Lip,+,s_2} \big) + \d_{\ind +1}^{- \frac34} 
 \\
&\ll    (C(s_2))^{\ind +1} \d_{\ind +1}^{- \frac34}  (\e \d_1^{- \frac12}+1) \, , 
\end{align*}
which gives \eqref{stime-rm-high} and the second bound in \eqref{form-A'-m}  at the step $ \ind + 1 $. 

To estimate the $s$-norms for any $s \geq s_2$, let us introduce the notation 
\be\label{uinds}
u_\ind (s)  := | R_\ind  |_{\Lip, +, s} + | \rho_\ind |_{\Lip, +, s} \quad {\rm and} \quad  
 \alpha (s)  := 3 \varsigma \frac{s-s_2}{ s_2 - s_1 } \quad {\rm as \ in \ }   \eqref{def:as} \, . 
\ee 
By \eqref{Rrhon+1}, we have the inductive bound
$$
u_{\ind+1} (s) = 
| R_{\ind +1} |_{\Lip, +, s} + | {\rho}_{\ind +1} |_{\Lip, +, s}   \leq C'(s) \big[ \d_{\ind+1}^{- \frac14} u_\ind (s) + \d_{\ind+1}^{- \frac34} \, \d_{\ind+1}^{- \alpha (s)} \big] \, . 
$$
Then, since 
\eqref{Rrn-tame0} and \eqref{Rrn-tame} 
hold at order $\ind$, we obtain
\begin{align*}
u_{\ind +1} (s) &\leq C'(s) \big[ \d_{\ind+1}^{- \frac14} (C(s))^\ind \d_\ind^{-\frac34 -\alpha(s)} \big( u_0(s) \d_1^{\frac12 + \frac{2 \alpha (s)}{3}} +1 \big)
+ \d_{\ind+1}^{- \frac34} \d_{\ind+1}^{- \alpha (s)} \big] \\
&\leq  (C(s))^{\ind + 1} \d_{\ind+1}^{-\frac34 -\alpha(s)}
\big( u_0(s) \d_1^{\frac12 + \frac{2 \alpha (s)}{3}} +1 \big)
\end{align*}
for $C(s)$ large enough, and using that $ \d_{\ind +1} = \d_\ind^{3/2}$. Hence the estimates 
\eqref{Rrn-tame0} and \eqref{Rrn-tame}
are proved also at order $\ind +1$. 
\\[1mm]
{\sc Estimates of $ |{\cal P}_{\ind +1}^{\pm 1}|_{\Lip,+,s}$.} 
We first consider the case $ s = s_2 $. By  \eqref{AkLip} and the first 
estimate in \eqref{estimateS-n+1} 
we obtain
 \begin{align}
 |e^{J {\cal S}_{\ind +1}}|_{\Lip, s_2+1} 
 & \leq 
1+   C(s_2)  \sum_{k \geq 1}  \frac{(C(s_2))^{k-1}}{k!} |{\cal S}_{\ind +1}|^{k-1}_{\Lip, s_1} 
 |{\cal S}_{\ind +1}|_{\Lip, s_2+1} \nonumber \\
 & 
  \leq  1+ C'(s_2) |{\cal S}_{\ind +1}|_{\Lip, s_2+1}  \, . \label{C2s2}
 \end{align}
Hence by \eqref{def:Pind+1} and \eqref{inter-norma+s}, we get  
\begin{align}
|{\cal P}_{\ind +1}|_{\Lip,+, s_2 } & \leq  
 C (s_2) \big( |{\cal P}_{\ind }|_{\Lip,+,s_2 } |e^{J {\cal S}_{\ind +1}}|_{\Lip, s_1+1} + 
|{\cal P}_{\ind }|_{\Lip,+,s_1 } |e^{J {\cal S}_{\ind +1}}|_{\Lip, s_2+1} \big) \nonumber \\
& \stackrel{\eqref{ejsn+1}, \eqref{C2s2}} \leq C (s_2) \big( |{\cal P}_{\ind }|_{\Lip,+,s_2 }  + 
|{\cal P}_{\ind }|_{\Lip,+,s_1 } | {\cal S}_{\ind +1}|_{\Lip,s_2 +1} \big) \label{qqs2}
\end{align}
for some new constant $ C(s_2) $. 
Therefore, by \eqref{qqs2},  properties \eqref{Pins2}, \eqref{Pins1} at order $\ind$, 
the second inequality in \eqref{estimateS-n+1}, 
and recalling the definition of $ \delta_{\ind}$ in \eqref{def:indices}, we have 
\begin{align*}
|{\cal P}_{\ind +1}|_{\Lip,+, s_2 } & \leq   
 C(s_2) \Big( (C(s_2))^\ind \d_\ind^{- \frac34}(\e \d_1^{- \frac12}+1)+ 2 \d_{\ind +1}^{- \frac14} 
\big( |R_\ind |_{\Lip,+,s_2} + |\rho_\ind |_{\Lip,+,s_2} \big) +\d_{\ind +1}^{- \frac34} \Big) \\
& \stackrel{\eqref{stime-rm-high},\eqref{form-A'-m}}{\leq}  
(C(s_2))^{\ind +1} \d_{\ind +1}^{- \frac34}  (\e \d_1^{- \frac12}+1) 
\end{align*}
provided the constant $C(s_2)$ is large enough, proving \eqref{Pins2} at the step $ \ind + 1 $. 
We obtain in the same way the estimate for  $ {\cal P}_{\ind +1}^{-1} $. 

Now, for any $ s \geq s_2 $, we derive by 
the last inequality in \eqref{estimateS-n+1}, \eqref{uinds}, 
\eqref{Rrn-tame0}, \eqref{Rrn-tame}, 
\begin{align}
|{\cal S}_{\ind +1}|_{\Lip,s+1} & \leq C(s) \big[ \d_{\ind+1}^{-\frac14} u_\ind (s) 
+ \d_{\ind+1}^{- \frac34 -\alpha (s)} \big] \nonumber \\
& \leq (C(s))^{\ind + 1} \d_{\ind+1}^{-\frac34 -\alpha(s)}
\big( u_0(s) \d_1^{\frac12 + \frac{2 \alpha (s)}{3}} +1 \big) \, . \label{eq:prop-Sn+1}
\end{align}
As in the case $s= s_2$ (see \eqref{C2s2})  we obtain
$$
|e^{J {\cal S}_{\ind +1}}|_{\Lip,s+1} \leq 1 + C(s) 
|{\cal S}_{\ind +1}|_{\Lip, s+1} \, . 
$$
Since $ {\cal P}_{\ind +1} = {\cal P}_\ind e^{J {\cal S}_{\ind +1}} $ we derive
the bound  \eqref{Pinsany} on $|{\cal P}^{\pm 1}_{\ind+1}|_{\Lip, +,s}$, exactly as in the case $s=s_2$,
using 
the interpolation inequality \eqref{inter-norma+s}, the inductive assumptions \eqref{Pinsany} and \eqref{Pins1}, 
\eqref{eq:prop-Sn+1}, and taking the constant $C(s)$ of \eqref{Pinsany} large enough.
\\[1mm]
This completes the iterative proof of $({\bf P}_\ind)_{\ind \geq 0}$. 
\\[2mm]
{\sc The Cantor-like sets $\Lambda_\infty (\e;\eta,  A_0, \rho)  $. }
We define, for $1/2 \leq \eta \leq 5/6 $, the  set
\be\label{def:Lambda-infty}
\Lambda_\infty (\e;\eta,  A_0, \rho) := \bigcap_{\ind=0}^\infty \Lambda_\ind (\e;\eta,  A_0, \rho)
= \bigcap_{k=0}^\infty \Lambda(\e; \eta + \eta_k,A_k) 
\ee
where 
$ \Lambda_\ind (\e;\eta,  A_0, \rho) $ are defined in  \eqref{set-cor-split}
and $ \Lambda (\e ; \eta + \eta_{k}, A_{k}) $ are  
 defined by  Proposition \ref{prophomeq}. We recall that the sequence $(\eta_k)$ is defined in \eqref{def:indices}.
 
 The sets $ \Lambda_\infty (\e;\eta,  A_0, \rho) $ satisfy  Property \ref{item2-91} of Corollary \ref{cor:split} as 
the sets $ \Lambda (\e ; \eta + \eta_{k}, A_{k}) $  satisfy Property 1 of   Proposition \ref{prophomeq}.
 \\[1mm]
{\sc Proof of Property \ref{item2-92} 
for the sets $\Lambda_\infty (\e;\eta,  A_0, \rho) $. }  
The complementary set of $ \Lambda_\infty (\e; \eta , A_0 , \rho)  $ may be decomposed
as (the sets $ \Lambda_\ind (\e;\eta,  A_0, \rho) $ in \eqref{set-cor-split} are decreasing in $ \ind $)
\begin{align}\label{compl-set-Linfty}
\Lambda_\infty (\e; \eta , A_0 , \rho)^c 
&  \stackrel{\eqref{def:Lambda-infty}} = \bigcup_{\ind=0}^\infty \Lambda_\ind (\e;\eta,  A_0, \rho)^c  \nonumber \\
&  \stackrel{ \Lambda_0  = \wtilde \Lambda  } = \wtilde \Lambda^c \cup
\bigcup_{k=0}^\infty \Big( \Lambda_k (\e;\eta,  A_0, \rho) \cap \Lambda_{k+1} (\e;\eta,  A_0, \rho)^c \Big) \nonumber \\
& \stackrel{\eqref{set-cor-split}} = \wtilde \Lambda^c \cup 
\bigcup_{k=0}^\infty \Big( \Lambda_k (\e;\eta,  A_0, \rho) \cap 
\big( \Lambda_k (\e;\eta,  A_0, \rho) \cap \Lambda (\e; \eta + \eta_k, A_k)  \big) ^c \Big)  \nonumber \\
& =  \wtilde \Lambda^c \cup  \bigcup_{k = 0}^{\infty} 
\Big(  \Lambda_k (\e; \eta , A_0 , \rho) \cap \Lambda (\e; \eta + \eta_k, A_k)^c  \Big) \nonumber \\
&  \stackrel{ \Lambda_0  = \wtilde \Lambda  } =  \wtilde \Lambda^c \cup  \Lambda (\e; \eta, A_0)^c 
 \bigcup_{k = 1}^{\infty} 
\Big(  \Lambda_k (\e; \eta , A_0 , \rho) \cap \Lambda (\e; \eta + \eta_k, A_k)^c  \Big)
\, . 
\end{align} 
By  Property \ref{sp:item2} of Proposition \ref{prophomeq} we have  
\be\label{Lam012}
 |\Lambda(\e;1/2,A_0)^c \cap \wtilde \Lambda | \leq b(\e) \quad {\rm with} \ \lim_{\e \to 0} b(\e) = 0 \, . 
\ee 
Moreover,  for $ k \geq 1 $, we have,  
by the definition of $ \Lambda_{k}  (\e; \eta , A_0 , \rho) $ in   \eqref{set-cor-split},
\begin{align}
& \Lambda_{k} \big(\e; \cc , A_0, \rho \big)   \bigcap  \Lambda \Big(\e; \cc + \eta_{k},A_{k} \Big)^c = \nonumber \\ 
& 
\Lambda_{k} \big(\e; \cc , A_0, \rho \big)  
 \bigcap  \Lambda \Big(\e; \cc + \eta_{k},A_{k} \Big)^c 
 \bigcap \Lambda \big(\e; \cc +\eta_{k-1},A_{k-1} \big) \, . 
 \label{trivial:incl}
\end{align}
Now, since, for all $ k \geq 1 $ we have $ | A_{k} - A_{k-1} |_{+,s_1} \leq \d_{k}^{3/4}$  on $\Lambda_{k} $ 
(see \eqref{An-An-1}) and 
$ \eta_{k} = \eta_{k-1} + \d_{k}^{3/8} $ (see \eqref{def:indices}), 
we deduce by  \eqref{Cantor:intersect0} and \eqref{trivial:incl}  the Lebesgue measure estimate,
\be \label{mesetilde}
 \Big| \Lambda_{k} \big(\e; \frac12, A_0, \rho \big)   \bigcap  \Lambda \Big(\e; \frac12 + \eta_{k},A_{k} \Big)^c \Big| 
  \leq \d_{k}^{3\alpha /4} \, , \quad \forall k \geq 1 \, .
\ee
In conclusion, by \eqref{compl-set-Linfty}, \eqref{Lam012},
\eqref{mesetilde} we obtain,  recalling 
\eqref{def:indices},
$$
\begin{aligned}
|\Lambda_\infty(\e; 1/2, A_0, \rho)^c \cap \wtilde \Lambda | & \leq  
b(\e)+ \sum_{k=1}^\infty \d_1^{\frac{3\alpha }{4}(\frac{3}{2})^{k-1}} \\
& \leq  b( \e  ) + 2 \d_1^{\frac{3\alpha }{4}} = b( \e  ) + 2 \e^{\frac{9 \alpha }{4}} =: b_1(\e) 
\end{aligned}
$$
since $\d_1 = \e^3 $.
Property \ref{item2-92} of Corollary \ref{cor:split} is proved. 
\\[1mm]
{\sc Proof of \eqref{pertAn}}. Actually we prove that, if $A'_0 \in {\cal C}(C_1,c_1,c_2)$, $(R'_0,\rho')$ satisfy 
\eqref{rho-R0:small} as $(R_0, \rho)$  and $(A'_0, \rho')$ satisfy  \eqref{pertA0}, then  
\be \label{pertAn-finita}
|A'_\ind - A_\ind|_{+, s_1} \leq \d^{4/5} \, , \quad 
\forall \l \in \Lambda_\ind (\e; 5/6, A_0, \rho) \cap  \Lambda_\ind (\e; 5/6, A_0', \rho') 	\, , 
\ee
where the sets $\Lambda_\ind $ are defined in \eqref{set-cor-split}. 
Let 
\be\label{def-Deltan}
\Delta_\ind^{(1)}:= |R_\ind' -R_\ind|_{+, s_1} \, ,   \quad 
\Delta_\ind^{(2)}:= |\rho_\ind' -\rho_\ind|_{+, s_1}
\, , \quad \Delta_\ind := \Delta_\ind^{(1)} + \Delta_\ind^{(2)} \, . 
\ee
The assumption \eqref{pertA0} means that $ \Delta_0 \leq \delta $ on $ \wtilde \Lambda \cap  \wtilde \Lambda' $.
We recall that $(A_{\ind +1}, \rho_{\ind +1})$ (resp. $(A'_{\ind +1}, \rho'_{\ind +1})$) is built 
applying Proposition \ref{prophomeq} to $(A_\ind , \rho_\ind)$ (resp. $(A'_\ind , \rho'_\ind)$) instead of 
$(A_0 , \rho_0)$ (resp. $(A'_0 , \rho'_0)$). Hence, by \eqref{pertR0+}-\eqref{pertrho+}
(with $ \d_1 $ replaced by  $ \d_{\ind + 1 } $),   for all $\ind \in \N$,  for all 
$ \l \in \Lambda_{\ind+1} (\e; 5/6 , A_0, \rho) \cap  \Lambda_{\ind+1} (\e; 5/6 , A_0', \rho') $, 
we have the iterative inequalities 
\be\label{induz-n+1}
\Delta_{\ind+1}^{(1)} \leq \Delta_\ind^{(1)} + C \Delta_\ind^{(2)}  \qquad
{\rm and} \qquad
\Delta_{\ind +1}^{(2)} \leq \d_{\ind +1}^{1/2} \Delta_\ind^{(1)} + \d_{\ind +1}^{-1/20} \Delta_\ind^{(2)} \, .
\ee
As a consequence  $ \Delta_{\ind +1} \leq \d_{\ind +1}^{- \frac{1}{19} } \Delta_\ind$ and, 
recalling the definition of $ \d_\ind $ in  \eqref{def:indices}, we deduce that, for any $\ind \geq 1  $, 
for all $ \l \in \Lambda_\ind (\e; \eta, A_0, \rho) \cap  \Lambda_\ind (\e; \eta, A_0', \rho') $, we have 
\be\label{Delta-ind}
\Delta_\ind \leq (\d_\ind \ldots \d_1)^{- \frac{1}{19} } \Delta_0 \leq \d_1^{-\frac1{19} (1+ \ldots + (\frac32)^{\ind-1})} \d
\leq \d_1^{-\frac{2}{19} (\frac32)^{\ind}} \d \, . 
\ee
Let $ \ind_0 \geq 1  $  be the integer such that 
\be\label{n0picco}
\d_1^{ (\frac32)^{\ind_0}} < \d < \d_1^{ (\frac32)^{\ind_0-1}}  \, . 
\ee
Thus, by \eqref{Delta-ind} and the second inequality in \eqref{n0picco}, we get
\be\label{nsmallern0}
\forall \ind \leq \ind_0 \, , \quad \Delta _\ind \leq \d_1^{-\frac2{19} (\frac32)^{\ind_0}} \d
\leq \d^{- \frac{3}{19}} \d = \d^{ \frac{16}{19} } \, . 
\ee
On the other hand,  recalling \eqref{def-Deltan} we bound,  
using  the first estimate in \eqref{form-A'-m}, 
\be\label{stima>n0+1}
\forall \ind \geq \ind_0 +1 \, , \  \ \Delta_\ind^{(2)}\leq |\rho_\ind|_{+, s_1} + |\rho'_\ind|_{+, s_1} 
 \leq 2 \d_1^{(\frac32)^\ind} \, .
\ee
Applying iteratively the first inequality in \eqref{induz-n+1} we get, $ \forall \ind \geq \ind_0 +1 $, 
\begin{eqnarray}
|A'_\ind - A_\ind|_{+, s_1} = \Delta_\ind^{(1)} & 
\leq &  \Delta_{\ind_0}^{(1)} + C \big( \Delta_{\ind_0}^{(2)}
+ \ldots + \Delta_{\ind-1}^{(2)}\big)  \nonumber \\
& \stackrel{\eqref{nsmallern0}, \eqref{stima>n0+1}}\leq  & \d^{\frac{16}{19} } +2 C ( \d_1^{(\frac32)^{\ind_0}} + \ldots + \d_1^{(\frac32)^{\ind-1}} ) \nonumber \\
& \leq & \d^{ \frac{16}{19} } + 3 C \d_1^{(\frac32)^{\ind_0}} \nonumber \\
& \stackrel{\eqref{n0picco}} \leq & \d^{ \frac{16}{19} } + 3 C \d \leq \d^{ \frac45 }  \label{ngreatern0}
\end{eqnarray}
for $ \d $ small enough. In conclusion \eqref{nsmallern0} and \eqref{ngreatern0} 
imply 
\eqref{pertAn-finita}. 
\\[1mm] 
{\sc Proof of Property \ref{item2-93} 
for the sets $\Lambda_\infty (\e;\eta,  A_0, \rho) $. } 
By \eqref{compl-set-Linfty} (with $ A_0', \rho' $ instead of $ A_0, \rho $) and \eqref{def:Lambda-infty}
we deduce, for all $ (1/2) + \delta^{2/5} \leq \eta \leq 5/ 6 $, the inclusion 
\be\label{M-included-Mn}
{\cal M} :=  \wtilde \Lambda' \cap 
\Lambda_\infty (\e; \eta , A'_0 , \rho')^c \cap  \Lambda_\infty (\e;  \eta-\d^{2/5} , A_0 , \rho)  
\subset  {\cal M}_0   \bigcup \Big( \bigcup_{\ind \geq 1} {\cal M}_\ind \Big) 
\ee
where
$$
 {\cal M}_0  :=  \wtilde \Lambda' \cap  \Lambda (\e; \eta, A'_0)^c \cap \Lambda (\e; \eta - \d^{2/5}, A_0) 
$$
and
\begin{align*} 
& {\cal M}_\ind  := \\ 
& \Lambda_\ind (\e; \eta , A'_0 , \rho') \bigcap \Lambda (\e;  \eta + \eta_\ind, A'_\ind)^c 
 \bigcap  \Lambda_\ind (\e;  \eta -\d^{2/5} , A_0 , \rho) 
\bigcap \Lambda (\e; \eta + \eta_{\ind} -\d^{2/5}, A_{\ind})  \, .
\end{align*}
By  the assumption 
\eqref{R0R0'vicini}  we have $  | R_0' - R_0|_{+, s_1}  \leq \d $ on $ \wtilde \Lambda \cap \wtilde \Lambda' $, 
and therefore property \ref{sp:item3} of Proposition \ref{prophomeq}
and the fact that $ \d^{1/2} \leq \d^{2/5} $ imply that, for
all $ (1/2) + \delta^{2/5} \leq \eta \leq 5/ 6 $, 
\begin{align}
|{\cal M}_0|  
& \leq  | \wtilde \Lambda' \cap  \Lambda (\e; \eta, A'_0)^c \cap \Lambda (\e; \eta - \d^{1/2}, A_0)  | 
\leq \d^\alpha \, . \label{estimM0}
\end{align}
For $\ind \geq 1$,  we have,  by 
\eqref{set-cor-split}, the inclusion  
$$
{\cal M}_{\ind} \subset \Lambda_\ind (\e; \eta , A'_0 , \rho') 
\cap \Lambda (\e; \eta + \eta_{\ind-1}, A'_{\ind-1}) \cap \Lambda (\e; \eta + \eta_\ind, A'_\ind)^c \, .
$$
and therefore, since $ | R'_\ind  - R'_{\ind-1}|_{+, s_1} =  | A'_\ind  - A'_{\ind-1}|_{+, s_1}  \leq \d_\ind^{3/4} $ 
on $ \Lambda_\ind (\e; \eta , A'_0 , \rho') $ by \eqref{An-An-1}, 
the estimate \eqref{Cantor:intersect0} and 
$ \eta_{\ind} = \eta_{\ind-1} + \d_{\ind}^{3/8} $, imply
\be\label{estimMn}
|{\cal M}_{\ind}| \leq \d_\ind^{3\alpha /4} \, .
\ee
On the other hand, by \eqref{pertAn-finita},   $|A_\ind - A'_\ind|_{+, s_1} \leq \d^{4/5}$  for any 
$$ \l \in {\cal M}_n \subset 
\Lambda_\ind (\e; \eta , A'_0 , \rho') \bigcap \Lambda_\ind (\e;  \eta , A_0 , \rho) 
$$ 
and we deduce, by \eqref{Cantor:intersect0}, 
the measure estimate
\be\label{estim:Nlow}
|{\cal M}_\ind| \leq \big| 
\Lambda (\e;  \eta + \eta_\ind, A'_\ind)^c \bigcap
 \Lambda (\e;  \eta + \eta_{\ind} -\d^{2/5}, A_{\ind}) \big| \leq \d^{4\alpha /5} \, .
\ee
Finally  \eqref{M-included-Mn}, \eqref{estimM0}, \eqref{estimMn}, \eqref{estim:Nlow} imply the measure estimate 
$$
|{\cal M}| \leq  \d^{\alpha} + \sum_{\ind \geq 1} \min (  \d_{\ind}^{3\alpha /4} , \d^{4\alpha /5}) \leq \d^{\alpha /2}   
$$
for $ \d $ small, proving \eqref{Cantor:intersect}. 
The proof of Corollary \ref{cor:split} is complete. 
\end{pf}

The rest of the chapter is devoted to the Proof of Proposition \ref{prophomeq}.

\section{The linearized homological equation}

We  consider   the linear map
$$
{\cal S} \mapsto J \ppavphi {\cal S}  + [J{\cal S}, J {A_0}] 
$$
where $ {\cal S} := {\cal S} (\vphi) $, $\vphi \in \T^\es $, has the form 
\be\label{forma cal-S}
\begin{aligned}
& \qquad {\cal S}(\vphi)= \begin{pmatrix}
d(\vphi) & a(\vphi)^* \\ a(\vphi) & 0
\end{pmatrix} \in {\cal L} (H_{\mathbb S}^\bot) \, , \quad \\
& d(\vphi) = d^* (\vphi ) \in {\cal L} (H_{\mathbb F}) \, ,  \quad   a(\vphi) \in {\cal L} (H_{\mathbb F}, H_{\mathbb G}) \, , 
\end{aligned}
\ee
and it is self-adjoint. 

Recalling \eqref{form:A0} and using that $ D_0 $ and $ J $ commute, $ J^2 = - {\rm Id} $, we have 
\be\label{def:homo-matrix}
\begin{aligned}
& J \ppavphi {\cal S}  + [J{\cal S}, J {A_0}]   = \\ 
& \begin{pmatrix}
J \ppavphi d  + D_0 d  + J d J D_0  & 
J \ppavphi  a^*   + J a^* J V_0  + D_0  a^*   \\
J \ppavphi  a  - J V_0 J a+ J a J D_0 & 0   
\end{pmatrix}  \, .
\end{aligned}
\ee
The key step in the proof 
of the splitting Proposition   \ref{prophomeq} (see Section \ref{sec:proof-splitting})
is, given  $  \rho(\vphi) $ of the form 
\be\label{decomposition-rho1}
\begin{aligned}
& \qquad \rho(\vphi) = \begin{pmatrix} \rho_1 (\vphi)& \rho_2 (\vphi)^* \\
\rho_2 (\vphi) & 0    \end{pmatrix} \in {\cal L} (H_{\mathbb S}^\bot ) \, , \\
& \rho_1 (\vphi) = \rho_1^* (\vphi) \in {\cal L} (H_{\mathbb F}) \, ,  \   \rho_2 (\vphi) \in {\cal L} (H_{\mathbb F}, H_{\mathbb G}) \, ,
\end{aligned}
\ee
to solve  (approximately) the ``homological" equation
\be \label{linhom}
 J \ppavphi {\cal S}   + [J{\cal S} , J {A_0} ]= J \rho  \, .  
\ee
The equation \eqref{linhom}  amounts to solve the pair of decoupled equations\index{Homological equations}
\begin{align} 
& J \ppavphi d + D_0 d + J d J D_0= J \rho_1 \, , \label{lineqd} \\
& J \ppavphi a -J V_0 J a + J a J D_0 =J \rho_2 \, . \label{eqhoma} 
\end{align}
Note that, taking the adjoint equation of \eqref{eqhoma}, multiplying
by  $ J $ on the left and the right, since $ V_0 $ and $ D_0 $ are self-adjoint
and  $J D_0 J=- D_0 $, $ J^* = - J $,   we obtain
$$
J \ppavphi a^* +  J a^* J V_0 +  D_0 a^* =J \rho_2^* \, ,
$$
which is the equation in the top right  in \eqref{linhom}, \eqref{def:homo-matrix}, \eqref{decomposition-rho1}. 

\begin{remark}
We  shall solve only approximately the homological equation \eqref{linhom} up to terms which are Fourier supported on high frequencies.
The main reason is that 
the multiscale Proposition
\ref{propmultiscale} provides tame estimates of the inverses of {\it finite} dimensional restrictions of infinite dimensional operators. 
This is sufficient for proving Proposition \ref{prophomeq}.
\end{remark}

We shall decompose  
an operator $ \rho $ of the form  \eqref{decomposition-rho1}, as well as $ {\cal S} $ in \eqref{forma cal-S}, in the following way.  
The operator $\rho_1 \in {\cal L} (H_{\mathbb F}) $ 
can be  represented as a finite dimensional self-adjoint square matrix 
$((\rho_1)_i^j )_{i,j \in {\mathbb F}}$ 
with entries $(\rho_1)_i^j \in {\cal L}(H_j , H_i) $. Using in each subspace $ H_{j} $, $ j \in {\mathbb F} $,   
the basis $((\Psi_j,0), (0,\Psi_j))$, see \eqref{Fj-eigenfunctions},
we identify each operator 
$ (\rho_1)_{i}^j (\vphi) \in {\cal L}( H_{j}, H_i) $ with a $ 2 \times 2 $-real matrix that we still denote
by $ (\rho_1)_{i}^j (\vphi) \in {\rm Mat}_2 (\R) $. We shall also Fourier expand 
\be\label{rho1:Fourier}
 (\rho_1  )_i^j (\vphi)  =  \sum_{\ell \in \Z^\es} [\widehat{\rho_1}]_i^j(\ell) e^{\ii \ell \cdot \vphi} \, , \ \
[\widehat{\rho_1}]_i^j(\ell) \in {\rm Mat}_2 (\C ) \,, \quad
\ov{[\widehat{\rho_1}]_i^j(\ell)} = [\widehat{\rho_1}]_i^j(-\ell)\,  , 
\ee
where $ [\widehat{\rho_1}]_i^j(0)  $ is the average 
\be\label{average-vphi}
[\widehat{\rho_1}]_i^j(0)  = \frac{1}{(2 \pi)^{\es}} \int_{\T^\es} (\rho_1  )_i^j (\vphi) \, d \vphi \in   {\rm Mat}_2 (\R ) \, . 
\ee
The operator $ \rho_2 \in {\cal L} (H_{\mathbb F}, H_{\mathbb G}) $ is identified,  as in 
\eqref{a-deco-sez}, with  $ (\rho_2^j )_{j \in {\mathbb F} } $ 
where  $ \rho_2^j \in {\cal L} (H_j, H_{\mathbb G}) $, which, using in $ H_j $    
the basis $((\Psi_j,0), (0,\Psi_j))$, can be identified with a vector of 
$ H_{\mathbb G} \times H_{\mathbb G} $.  

We also recall that, for a $ \vphi $-dependent family of operators $ \rho (\vphi) $, $ \vphi \in \T^\es $, of the form \eqref{decomposition-rho1},  we have 
the estimates \eqref{special-form}-\eqref{special-form-Lip}.

The next lemma provides an approximate solution of the homological equation \eqref{linhom}.


\begin{lemma} \label{reshomlin} {\bf (Homological equations)}
Given $C_1>0$, $c_1>0$, $c_2>0$, there is $ \e_1 > 0 $ such that 
$ \forall \e \in (0, \e_1)$, for each split admissible operator
$A_0 \in {\cal C} (C_1,c_1,c_2) $ (see Definition \ref{def:calC}), defined for $ \lambda \in \wtilde \Lambda $,
there are closed subsets  $ \Lambda (\e; \cc,A_0) \subset \wtilde \Lambda $,
$ 1/2 \leq \eta \leq 1 $,  satisfying 
the properties 1-3 of Proposition \ref{prophomeq}, such that,  if $\rho  \in {\cal L}(H_{\mathbb S}^\bot ) $ 
has the form \eqref{decomposition-rho1}, $ \rho $  is  
$ L^2 $-{\it self-adjoint}, defined for  $ \l \in \wtilde \Lambda $,  and 
satisfies 
\begin{align}\label{prop:rho-plus}   
& |  \rho  |_{\Lip, +, s_1} \leq \delta_1^{\frac{9}{10}} \, , \ \d_1 \leq  \e^3 \, , \qquad 
 \delta_1^{\frac{11}{10}} (|  R_0  |_{\Lip, +, s_2} +  |  \rho |_{\Lip, +, s_2}  )  \leq  \e \, , \\
& \label{propertyM-}
 {[\widehat{\rho_1}]^j_j} (0) \in M_-  \  ({\it recall} \ \eqref{average-vphi}, \eqref{deco:M+M-}),   \  \ \forall j \in {\mathbb F} \, ,   \  \forall \l \in 
 \wtilde \Lambda  \, , 
\end{align}
then there is  a linear self-adjoint operator 
$ {\cal S} := {\cal S}(\e, \l)( \vphi) \in {\cal L} (H_{\mathbb S}^\bot ) $  
of the form \eqref{forma cal-S}, defined for all $ \l \in \Lambda(\e; 1, A_0) $, 
satisfying \eqref{stima-S}-\eqref{estanys},  such that 
\begin{align}\label{approxlin}
& |J \ppavphi {\cal S}  + [J{\cal S}, J {A_0}] - J \rho|_{\Lip, +,s_1}  \leq \d_1^{\frac74} \, , 
 \\
& | J\ppavphi {\cal S} + [J{\cal S}, J {A_0}] |_{\Lip, +,s_2}
 \leq \d_1^{-\frac14} \big( |R_0|_{\Lip, +,s_2} + |\rho|_{\Lip, +,s_2} \big) +
\d_1^{-\frac34} \, , \label{31S}
\end{align}
and, more generally, for all $s \geq s_2 $, 
\be\label{tameS-any-s}
 | J\ppavphi {\cal S} + [J{\cal S}, J {A_0}] |_{\Lip, +,s} \leq 
C(s)  \Big[\d_1^{- \frac14} \big(| R_0 |_{\Lip, +, s}+ | \rho |_{\Lip, +, s} \big)
+ \d_1^{- \frac34} \d_1^{- 3 \varsigma \frac{s-s_2}{s_2-s_1}} \Big] \, . 
\ee
At last, denoting by $ {\cal S}_{A_0, \rho} $ and $ {\cal S}_{A_0', \rho'} $ the operators 
defined as above  
associated,  respectively, to $ (A_0, \rho) $ and  $( A'_0, \rho') $,  
we have, for all  $\l \in \Lambda(\e; 1, A_0) \cap \Lambda(\e;1, A'_0)$, 
\be  \label{pertS}
|{\cal S}_{A_0, \rho} - {\cal S}_{A'_0, \rho'}|_{+, s_1}  \leq 
\d_1^{\frac34} |A_0 -A'_0|_{+, s_1} + \d_1^{-\frac{1}{30}} |\rho -\rho'|_{+, s_1} \, , 
\ee
and
\begin{align} \label{perteq}
& \big|(J \ppavphi {\cal S}_{A_0,\rho}  + [J{\cal S}_{A_0,\rho}, J {A_0}] - J \rho)-
( J \ppavphi {\cal S}_{A'_0,\rho'}  + [J{\cal S}_{A'_0,\rho'}, J {A_0}] - J \rho' ) \big|_{ +,s_1} \nonumber \\
& \leq 
\d_1^{\frac34} |A_0 -A'_0|_{+, s_1} + \d_1^{- \frac{1}{30}} |\rho -\rho'|_{+, s_1} \, . 
\end{align}
\end{lemma} 

The proof of Lemma \ref{reshomlin}  is given  in the next section.

\section{Solution of homological equations: proof of Lemma \ref{reshomlin} } \label{sec:homo-split} 


\noindent
{\bf Step 1: approximate solution of the homological equation \eqref{lineqd}.}  
We represent a  linear operator
$ d(\vphi ) \in {\cal L}( H_{\mathbb F}) $ 
by a finite dimensional square matrix  $ ( d_{i}^j (\vphi)  )_{i, j \in {\mathbb F}}  $ with  entries 
$ d_{i}^j (\vphi) \in {\cal L}( H_j, H_i ) \simeq {\rm Mat}_2 (\R) $.
Since the symplectic operator $ J $ leaves invariant each subspace $ H_{j} $ and 
$$ 
D_0 = {\rm Diag}_{ j \in \mathbb F} \, \mu_j(\e,\l) {\rm Id}_2 
$$ 
(see  \eqref{form:D0}),  the equation \eqref{lineqd} is equivalent to 
$$
\begin{aligned}
& J \ppavphi d_i^j (\vphi) + 
\mu_i(\e,\l) d_i^j (\vphi)  + \mu_j(\e,\l) J d_i^j (\vphi) J = J (\rho_1 )_i^j (\vphi) \, , \\ 
& \qquad  \qquad   \forall i,j \in {\mathbb F} \qquad {\rm where} \qquad
J = 
\begin{pmatrix}
 0 & 1   \\
 -1 & 0   \\
\end{pmatrix} \, ,
\end{aligned}
$$
and, by a Fourier series expansion with respect to the variable $ \vphi \in \T^{\es}$, writing
\be\label{svil-dij-vphi}
d_i^j (\vphi) = \sum_{\ell \in \Z^\es} \widehat{d}_i^j (\ell) e^{\ii \ell \cdot \vphi} \, , \quad
\widehat{d}_i^j (\ell) \in {\rm Mat}_2 (\C ) \,, \quad \ov{\widehat{d}_i^j (\ell)} = \widehat{d}_i^j (-\ell) \, , 
\ee
to  
\be \label{lineqd3}
\begin{aligned}
 \ii (\bar \om_\e \cdot \ell ) J \widehat{d}_i^j (\ell) + \mu_i(\e,\l) \widehat{d}_i^j (\ell) + \mu_j(\e,\l) J \widehat{d}_i^j  (\ell) J 
 = J [\widehat{\rho_1}]_i^j(\ell) \, , \\
 \forall i,j \in {\mathbb F} \, ,  \ell \in \Z^\es \, .
\end{aligned} 
\ee
In order to solve \eqref{lineqd3} we have to study the linear operator 
\be\label{op-F-homo}
 T_{ij\ell} :  {\rm Mat}_2 (\C ) \to  {\rm Mat}_2 (\C )\, , \quad 
{\mathtt d} \mapsto \ii \bar \om_\e \cdot \ell \, J {\mathtt d}  + \mu_i (\e, \l) {\mathtt d}   + \mu_j (\e, \l)
J {\mathtt d}  J  \, .
\ee
In the basis $ (M_1,M_2,M_3,M_4)$ of $ {\rm Mat}_2 (\C ) $ defined in \eqref{def:E1E2}-\eqref{def:E3E4}
 the linear operator $ T_{ij\ell} $  is represented  by the following 
 self-adjoint matrix (recall that $JM_l =M_l J$ for $ l =1,2$ and
$JM_l = -M_l J$ for $ l =3,4$, and $ J M_1 = - M_2 $, $ J M_3 = - M_4 $)
\be\label{representation-matrix-M1M4}
\begin{pmatrix}
  (\mu_i - \mu_j)(\e, \l) &  \ii \bar \om_\e \cdot \ell  & 0 & 0  \\
 -  \ii \bar \om_\e \cdot \ell  &    (\mu_i - \mu_j)(\e, \l)   &  0 & 0  \\
 0 & 0  &     (\mu_i + \mu_j)(\e, \l)  &   \ii \bar \om_\e \cdot \ell \\
 0 & 0  &  -  \ii \bar \om_\e \cdot \ell &   (\mu_i + \mu_j)(\e, \l)
\end{pmatrix} \, .
\ee
As a consequence  the eigenvalues of  $ T_{ij \ell} $ are
\be\label{eige:Tijl}
\pm \bar \om_\e \cdot \ell + \mu_i (\e, \l) - \mu_j(\e, \l) \, , \quad
  \pm \bar \om_\e \cdot \ell + \mu_i(\e, \l) + \mu_j (\e, \l) \, . 
\ee
To impose non-resonance conditions we define
the sets, for $ 1/2 \leq \eta \leq 1 $,  
\begin{align}\label{def:Lambda1}
\Lambda^1 (\e;\cc,A_0) & := \Big\{ \l \in \wtilde \Lambda \, : \, 
 | \bar \om_\e \cdot \ell \pm \mu_j(\e,\l)\pm \mu_i(\e,\l)| \geq \frac{\gamma_1}{2 \cc \langle \ell \rangle ^{\tau}}, \\
&  \qquad \qquad  \qquad \forall (\ell , i, j)   \in \Z^\es \times {\mathbb F} \times {\mathbb F} \, , \  (\ell, i, j) \neq (0,j,j) 
  \Big\} \nonumber 
\end{align}
where the constant $ \g_1 = \g_0 / 2 $  (recall that $ \g_0 $ is fixed in \eqref{diop}) and
\be\label{def>tau} 
\tau \geq (3/2) \tau_1 +3 + \es  
\ee 
where $ \tau_1 $ is defined in \eqref{def:tau1}.
The inequalities in \eqref{def:Lambda1} are 
second-order Melnikov\index{Second Melnikov non-resonance conditions} non-resonance 
concerning only a finite number of normal frequencies.

\begin{remark}
Since $ \mu_j (\e, \l) > c_0 > 0 $, for any $ j \in \N $,  if $  \gamma_1 \leq  2 c_0 $, then 
 the inequality 
$ | \bar \om_\e \cdot \ell + \mu_j(\e,\l) +  \mu_i(\e,\l)| 
\geq \frac{\gamma_1}{2 \cc \langle \ell \rangle ^{\tau}} $ in \eqref{def:Lambda1} holds  for $ \ell = 0 $ and $ j = i  $, for all
$ \l \in \Lambda $, $ \eta \in [1/2, 1 ]$.   
\end{remark}

In the next lemma we find a solution $ d_N $
of the projected homological equation 
\be\label{eq:approx-eq-1}
J \ppavphi d_N + D_0 d_N + J d_N J D_0 = {\it \Pi}_N J \rho_1 
\ee
where here the projector  $ {\it \Pi}_N $ applies to functions depending only on the variable $ \vphi $, namely 
\be\label{def:PiN-time}
{\it \Pi}_N  : \quad h(\vphi) = \sum_{\ell \in \Z^\es} h_{\ell} e^{\ii \ell \cdot \vphi } \quad \mapsto \quad 
({\it \Pi}_N h)(\vphi) := \sum_{|\ell |  \leq N} h_{\ell} e^{\ii \ell \cdot \vphi } \, . 
\ee

\begin{lemma} \label{homdiag} {\bf (Homological equation \eqref{lineqd})}
Let $ \rho  $ be a  self-adjoint operator 
of the form \eqref{decomposition-rho1}, defined for $ \l \in \wtilde \Lambda $, satisfying  
\eqref{prop:rho-plus}, and  assume that the average $ [\widehat{ \rho_1}]^j_j (0)  \in M_- $, 
$ \forall j \in {\mathbb F}$, $ \forall \l $, i.e. \eqref{propertyM-} holds.
Let  
\be \label{defNsp:lemma1} 
N \in \Big[ \d_1^{-\frac{3}{s_2-s_1}} -1, \d_1^{-\frac{3}{s_2-s_1}} +1 \Big]\, .
\ee
Then, for all $\l \in  \Lambda^1 (\e;1,A_0) $ (set defined in \eqref{def:Lambda1}) 
there exists $ d_N := d_N (\lambda; \vphi )  \in {\cal L} ( H_{\mathbb F} ) $, $ d_N  = d_N^* $, 
satisfying 
\begin{align}\label{stima-di-d1}
 | d_N |_{\Lip, s_1+1} \ll \d_1^{\frac78} \, , \quad 
 | d_N |_{\Lip , s+1} \ll \d_1^{- \frac{1}{40}} |\rho|_{\Lip ,s} \, , \ \forall s \geq s_1 \, , 
\end{align}
which solves the projected homological equation \eqref{eq:approx-eq-1}, and which is 
an approximate solution of the homological equation \eqref{lineqd} in the sense that 
\begin{align}\label{soluzione:dN}
|J \ppavphi d_N + D_0 d_N  + J d_N J D_0- J \rho_1|_{\Lip, s_1}   \ll \d_1^{7/4} \, . 
\end{align}
At last, denoting by $ d'_N $ the operator associated to $ D'_0, \rho' $ as above, we have, 
for all $ \l \in \Lambda^1 (\e;1, A_0') \cap \Lambda^1 (\e;1,A_0)$,  
\be \label{pertdN}
| d_N' -d_N|_{s_1} \leq \d_1^{ \frac78} |D_0' - D_0|_{s_1} + \d_1^{- \frac{1}{40}} |\rho_1' - \rho_1|_{s_1} 
\ee
and
\begin{align}
& \big| \big(J \bar \om_\e \cdot \pa_\vphi d_N' + D_0' d_N'  + J d_N' J D_0'- J \rho_1' \big)-
\big(J \bar \om_\e \cdot \pa_\vphi  d_N + D_0 d_N  + J d_N J D_0- J \rho_1 \big) 
\big|_{s_1} \nonumber \\
& \leq |\rho_1'- \rho_1 |_{s_1} \, . \label{perteq1}
\end{align}
\end{lemma}

\begin{pf} We split the proof in some steps. 
\\[1mm]
{\sc Step $ 1$.}
{\it For any $ (\ell, i,j ) \in \Z^\es \times {\mathbb F} \times {\mathbb F} $, 
for all $ \l \in \Lambda^1 (\e;1,A_0) $,  
there exists a solution 
$ {\widehat d}_i^j (\ell) $ of the equation \eqref{lineqd3},
satisfying the reality condition \eqref{svil-dij-vphi} and  
\be\label{trivialbo}
\| {\widehat d}_i^j (\ell) \|  
\leq C \langle \ell \rangle^\tau   \| [ {\widehat \rho_1}]_i^j (\ell) \| 
\ee
where $ \| \ \| $ denotes some norm in $ {\rm Mat}_2 (\C ) $.}

Recalling \eqref{eige:Tijl} and 
the definition of $  \Lambda^1 (\e;1,A_0) $ in \eqref{def:Lambda1}, we have that,
for all $ \l \in \Lambda^1 (\e;1,A_0) $, 
the operators $ T_{ij \ell} $ defined in \eqref{op-F-homo} are isomorphisms of $ {\rm Mat}_2 (\C ) $
for $ \ell \neq 0  $ or $ i \neq  j $, and 
\be\label{Tijellinv}
\| T_{ij\ell}^{-1} \| \leq C \langle \ell \rangle^{\tau}  
\ee
where $ \| \ \| $ denotes  some norm in $ {\cal L}({\rm Mat}_2 (\C )) $ and 
$ C := C(\g_1 ) $.  Hence, in this case, there exists a unique solution 
$ {\widehat d}_i^j (\ell) =
T_{ij\ell}^{-1} \big(J[ {\widehat \rho_1}]_i^j (\ell)$ of the equation \eqref{lineqd3}.

Let us consider the case  $\ell=0$ and $i=j$.  
By \eqref{representation-matrix-M1M4}, the linear operator  
$T_{jj0}$ is represented, in the basis
$(M_1,M_2,M_3,M_4)$ of ${\rm Mat}_2 (\C)$, by the matrix 
$$
2 \mu_j (\e,\l) \begin{pmatrix} 0_2 & 0 \\ 0 & {\rm Id}_2 \end{pmatrix} \, .
$$ 
Moreover
there is $c_0>0$ such that $\mu_j(\e,\l) > c_0$  for any $j \in \N$, $\e \in (0,\e_0)$,
$\l \in \Lambda$.
Hence the range of $T_{jj 0} $, is the subspace 
$ M_- = {\rm Span}(M_3,M_4)$ of ${\rm Mat}_2 (\C)$ defined in 
\eqref{deco:M+M-}, \eqref{def:E3E4},  and
$$
\forall \wtilde{\rho} \in M_-  \, , \quad \exists \, ! \, \wtilde{d} \in M_- \quad {\rm solving} \quad
T_{jj0}\wtilde{d}=\wtilde{\rho}  \quad {\rm with } \quad 
\| \wtilde{d} \| \lesssim \| \wtilde{\rho} \| \, ; 
$$
with a little abuse of notation, we denote 
$\wtilde{d}$ by $T_{jj0}^{-1} (\wtilde{\rho})$.
By the assumption $ [\widehat{ \rho_1}]^j_j (0) \in M_-$ and  the fact that 
 $M_-$ is stable under the action of $J$, we have that 
 $ J [\widehat{ \rho_1}]^j_j ( 0 )$ is in $ M_- $ and there is a 
  unique solution $ {\widehat d}_j^j(0) \in M_-$ of  the equation \eqref{lineqd3}.
 
To conclude the proof of Step 1, it remains to show that  
the matrices 
$ {\widehat d}_i^j (\ell) \in {\rm Mat}_2 (\C ) $ satisfy the reality condition \eqref{svil-dij-vphi}. 
It is consequence of $T_{ij(-\ell)} =\overline{T_{ij\ell}}$, the fact that 
$ [ {\widehat \rho_1}]_i^j (\ell) $ satisfies the reality condition \eqref{rho1:Fourier},  and 
 the uniqueness of the solutions of $ T_{ij\ell} ( {\widehat d}_i^j (\ell)  ) =  J [ {\widehat \rho_1}]_i^j (\ell)$ 
 (with the condition  $d_i^j (\ell)  \in M_-$ if $ \ell = 0 $ and $i=j$). 
\\[1mm]
{\sc Step 2. Definition of $ d_N $ and proof that $ d_N = d_N^* $. }
We define $ d_N := {\it \Pi}_N d (\vphi) $ where 
$ {\it \Pi}_N $ is the projector defined in \eqref{def:PiN-time} and
\be\label{defdN}
d (\vphi) = (d_i^j (\vphi) )_{i,j \in {\mathbb F}} \, , \quad 
d_i^j (\vphi)  = \sum_{\ell \in \Z^\es} \widehat{d}_i^j (\ell) e^{\ii \ell \cdot \vphi} \, , 
\quad \widehat{d}_i^j (\ell) = 
T_{ij\ell}^{-1} \big(J[ {\widehat \rho_1}]_i^j (\ell)  \, , 
\ee
is the unique solution of \eqref{lineqd}. 
Note that, taking the adjoint equation of \eqref{lineqd}, multiplying
by  $ J $ on the left and the right, since $ V_0 $ and $ D_0 $ are self-adjoint,
the operators 
 $ D_0 $ and $ J $ commute, $ J^* = - J $, and using that $ \rho_1 $ is self-adjoint, 
we obtain
$$
 J \ppavphi d^* + D_0 d^* + J d^* J D_0= J \rho_1 \, . 
$$
By uniqueness $ d^* = d $ is the unique solution of \eqref{lineqd}. 
Then also $ d_N^* = d_N $. 
\\[1mm]
{\sc Step $ 3 $. Proof of \eqref{stima-di-d1}.}
By \eqref{trivialbo}\index{Smoothing operators} (and \eqref{special-form}, \eqref{smoothingS1S2-Lip})  we have 
\be\label{es:int0}
  | d_N |_{s +1} \lesssim N^{\tau+1}   |\rho_1 |_{s}  \, . 
\ee
We claim also the estimate 
\be\label{es:int1}
| d_N |_{\lip,s+1} \lesssim N^{2 \t+1}  | \rho_1 |_{\Lip,s} 
\ee
which implies, together with \eqref{es:int0},  
\be\label{es:int-total}
| d_N |_{\Lip, s+1} \lesssim N^{2 \t+1}  | \rho_1 |_{\Lip,s} \, . 
\ee
To prove \eqref{es:int1} notice that, by \eqref{Tijellinv}, 
\begin{align}
\| T_{ij \ell}^{-1} (\l_1) - T_{ij \ell}^{-1} (\l_2) \| & = \| T_{ij \ell}^{-1} (\l_1) \big(  T_{ij \ell} (\l_2) - T_{i,j, \ell} (\l_1) \big) 
T_{ij\ell}^{-1} (\l_2) \big) \| \nonumber		\\
& \lesssim  \langle \ell \rangle^{2 \tau}  \| T_{ij \ell} (\l_2) - T_{ij \ell} (\l_1) \| 
\lesssim \langle \ell \rangle^{2 \tau}  \e^2 | \l_2 - \l_1 | \label{diffTijell}
\end{align}
recalling the definition of  $ T_{ij \ell} $  in \eqref{op-F-homo}, and \eqref{Hyp4}. 
By \eqref{Tijellinv} and \eqref{diffTijell} we  estimate the $ \lip $ seminorm of 
$ d_i^j (\ell) = T_{ij\ell}^{-1} \big( [\widehat \rho_1]_j^j (\ell) \big) $ as 
\begin{align*}
\| d_i^j (\ell) \|_{\lip}  & 
\lesssim \langle \ell \rangle^{2 \tau}  \e^2 \| [\widehat \rho_1]_j^j (\ell) \|  +   
\langle \ell \rangle^{ \tau}  \|  [\widehat \rho_1]_j^j (\ell) \|_{\lip}  
\end{align*}
and \eqref{es:int1} follows. 
The estimates \eqref{es:int-total}, 
\eqref{defNsp:lemma1} and \eqref{choice s2 s3}-(i), finally imply 
$$ 
| d_N |_{\Lip, s+1} \leq C 
N^{2 \tau+1}   | \rho_1 |_{\Lip, s} \ll \d_1^{- \frac{1}{40}}  | \rho_1 |_{\Lip, s} 
$$ 
which gives the second inequality in \eqref{stima-di-d1}. 
By \eqref{prop:rho-plus} we deduce 
$$
 | d_N |_{\Lip, s_1+1} \ll \d_1^{- \frac{1}{40}}  \d_1^{\frac{9}{10}} = \d_1^{\frac78} 
 $$
which is the first inequality in \eqref{stima-di-d1}.
\\[1mm]
{\sc Step $ 4 $. Proof of  \eqref{soluzione:dN}.}
By the definition of $ d_N $ we have 
\be  \label{eq1approx}
J \ppavphi d_N  + D_0 d_N + J d_N J D_0- J \rho_1   = - {\it \Pi}_N^\bot  J \rho_1 \, .
\ee
Then, recalling \eqref{special-form}, 
\begin{align*}
| J \ppavphi d_N  + D_0 d_N + J d_N J D_0- J \rho_1 |_{\Lip, s_1} 
& 
= | {\it \Pi}_N^\bot  J \rho_1 |_{\Lip, s_1}  \\
& \stackrel{\eqref{smoothingS1S2-Lip},\eqref{prop:rho-plus}} 
{\lesssim_{s_2}} N^{-(s_2-s_1)} \d_1^{- \frac{11}{10}} \\
& \stackrel{\eqref{defNsp:lemma1}} 
{\lesssim_{s_2}} \d_1^3 \d_1^{- \frac{11}{10}} \ll \d_1^{7/4} \, . 
\end{align*}
{\sc Step 5. Proof of \eqref{pertdN}-\eqref{perteq1}}. 
We denote by $ T_{ij\ell}' $ the linear operator as in \eqref{op-F-homo}
associated to $ A_0' $, i.e. with $ \mu_j' (\e, \l) $ instead of $ \mu_j (\e, \l) $.
Arguing as for \eqref{diffTijell} we get that, 
if $ \l \in \Lambda^1 (\e;1, A_0) \cap \Lambda^1 (\e;1,A'_0)$,  then 
\be\label{diff-eigenval}
\big\|  T_{ij\ell}^{-1} - (T_{ij\ell}')^{-1} \big\|
\leq C  \la \ell \ra^{2\tau}    \max_{j\in \mathbb F} \big| \mu_j (\e, \l) - \mu'_j (\e, \l) \big|  \, .
\ee
Let us denote $d_N := d_{N,D_0,\rho} $ and $ d_N' := d_{N,D'_0, \rho'} $ where we highlight the dependence of $ d_N $ 
with respect to $ D_0 $ and  $ \rho $. 
Using that $d_N$ depends linearly on $\rho$, actually only on $ \rho_1 $, see \eqref{defdN}, we get 
\begin{align*}
|d_{N,D_0,\rho} - d_{N,D'_0, \rho'}|_{ s_1} & \leq  |d_{N,D_0,\rho} - d_{N,D'_0, \rho}|_{ s_1} +
|d_{N,D'_0, \rho-\rho'}|_{s_1} \\
& \stackrel{\eqref{diff-eigenval}, \eqref{es:int0}} \lesssim  
 N^{2\tau}  |\rho_1 |_{s_1}  \max_{j \in \mathbb F} |\mu_j (\e, \l) - \mu'_j (\e, \l) | + 
 N^{\tau+1}  |\rho_1 -\rho_1'|_{s_1} 
\end{align*}
and therefore, by \eqref{prop:rho-plus}, 
\eqref{defNsp:lemma1}, 
\eqref{choice s2 s3}-(i), $ | \rho_1 |_{+,s} \simeq  | \rho_1 |_{s} $  we conclude that
\begin{align*}
|d_{N,D_0,\rho} - d_{N,D'_0, \rho'}|_{s_1}  &  \leq  \d_1^{- \frac{1}{40}} 
 \d_1^{ \frac{9}{10}} |D_0-D'_0|_{ s_1} + \d_1^{- \frac{1}{40}} |\rho_1 - \rho_1'|_{s_1} \\
&\leq  \d_1^{\frac45} |D_0-D'_0|_{s_1} + \d_1^{- \frac{1}{40}} |\rho_1 - \rho_1'|_{s_1}   
\end{align*}
proving \eqref{pertdN}. Finally \eqref{perteq1} is an immediate consequence of \eqref{eq1approx}.
\end{pf}

We now prove {\it measure estimates} for the sets $\Lambda^1 (\e;\cc,A_0)  $
defined in \eqref{def:Lambda1}. 

\begin{lemma}\label{lemma:measure1}
$ |[\Lambda^1 (\e;\cc, A_0)]^c \cap \wtilde \Lambda | \leq  \e $ for all $ 1/2 \leq \eta \leq 1 $.  
\end{lemma}

\begin{pf}
For $ \ell \in \Z^\es $, $ \cc \in [1/2,1]$, we define the set 
\be\label{Def:Lambda-l1}
\begin{aligned}
\Lambda^1_\ell (\e;\cc, A_0) 
& := \Big\{ \l \in \wtilde \Lambda \, : \, 
 | \bar \o_\e \cdot \ell \pm \mu_j(\e,\l)\pm \mu_i(\e,\l)| \geq \frac{\gamma_1}{ 2 \cc \langle \ell \rangle^{\tau}} \, , \\
& \qquad \ 
\forall (i,j) \in  {\mathbb F}  \times {\mathbb F}  \ {\rm with }  \   | \ell | + | i-j| \neq 0 \Big\} 
\end{aligned} 
\ee
where $ \gamma_1 = \g_0 / 2 $ and $ \tau > \tau_1 $, see  \eqref{def>tau}. 
By the unperturbed second order Melnikov non-resonance conditions \eqref{2Mel+}-\eqref{2Mel}, the Diophantine property \eqref{dioep}, 
since $ \bar \o_\e = \bar \mu + \e^2 \zeta  $ (see \eqref{def omep}) and
$ \mu_j(\e,\l) = \mu_j  + O(\e^2) $ (see \eqref{muje2}), we have
$$ 
{\rm if} \ \e^2 \langle \ell \rangle^{\tau_0 + 1 } \g_0^{-1} \leq c  \ {\rm small \ enough} 
\quad \Rightarrow
\quad  \Lambda^1_\ell (\e;1/2, A_0) = \wtilde \Lambda \, . 
$$ 
Hence, for $ \e $ small enough,
\be\label{prima-inclusione}
\Lambda^1 (\e;\cc, A_0) = \bigcap_{\ell \in \Z^\es} \Lambda^1_\ell (\e;\cc, A_0)  
= \bigcap_{|\ell | > c_1 \e^{- 2/ \tau_0}} \Lambda^1_\ell (\e;\cc, A_0) 
\ee
where $ c_1 := (c \g_0)^{1/\tau_0} $. 
Now, using  \eqref{Hyp2}-\eqref{Hyp3}, we deduce that 
the complementary set $ [\Lambda^1_\ell (\e;\cc, A_0)]^c $ is included in the union of $ 4{\mathtt f}^2 $ intervals of
length $  \dps \frac{ \g_1}{ \eta \langle \ell 
\rangle^\tau c_2 \e^2} $ where $ {\mathtt f} $ denotes the cardinality  of $ {\mathbb F} $. Hence
its measure satisfies
$$
|[\Lambda^1_\ell (\e;\cc, A_0)]^c \cap \wtilde \Lambda | 
\leq  \frac{4 \g_1 {\mathtt f}^2}{ \eta \langle \ell \rangle^\tau c_2 \e^2} \, , 
$$ 
and, for any $ L > 0  $, $ 1/ 2 \leq \eta \leq 1 $, $ \tau > \es $, we have 
\be \label{estsum}
 \sum_{| \ell | \geq L} \big| [\Lambda^1_\ell (\e;\cc, A_0)]^c \big| \leq
\sum_{| \ell | \geq L} \frac{4 \g_1 {\mathtt f}^2}{ \eta \langle \ell \rangle^\tau c_2 \e^2} \leq   \frac{ C}{ L^{\tau-\es}  \e^2} \, . 
\ee
The lemma follows by \eqref{prima-inclusione}, \eqref{estsum} with 
$ L = c_1 \e^{-2/\tau_0} $ and \eqref{def>tau}. 
\end{pf}

\begin{remark}\label{rem:meas1s}
The measure of the set 
$ [\Lambda^1 (\e;\cc, A_0)]^c   $ can be made
smaller than $ \e^p $ as $\e$ tends to $0$, for any $ p $, if we take the exponent 
 $ \tau $  large enough. This is analogous to the situation described in Remark \ref{rem:meas1}. 
\end{remark}

\begin{lemma}\label{lemma:hom1-var} 
Assume that $|A'_0-A_0|_{+,s_1} \leq \d \leq \e^3 $ on 
$ \wtilde{\Lambda} \cap \wtilde{\Lambda}' $. 
Then, 
for $ \cc \in [(1/2) + \sqrt{\d}, 1] $,  we have 
\begin{align}
|\wtilde{\Lambda}' \cap [\Lambda^1  (\e;\cc, A'_0)]^c \cap \Lambda^1 (\e;\cc-\sqrt{\d}, A_0)| 
& \leq \d^{\frac{1}{12}} \, . \label{Cantor-sovra1}
\end{align}
\end{lemma}

\begin{pf} 
We first prove  the estimate
\be\label{midif0}
\big| [\Lambda^1 (\e;\cc, A'_0)]^c  \bigcap \Lambda^1 (\e;\cc-\sqrt{\d}, A_0) \big| 
\lesssim  \min ( \e^{-2} \d^{\frac12-  \frac{\es}{2\tau}}, \e) \, .
\ee
Denoting $  \Lambda^1_\ell (\e;\cc, A'_0) $  the set \eqref{Def:Lambda-l1} associated to
$ A_0' $, we claim that there exists $ c(s_0) > 0 $ such that  
\be \label{prdis1} 
[ \Lambda^1_\ell (\e;\cc, A'_0)]^c \cap \Lambda^1 (\e;\cc-\sqrt{\d}, A_0) = \emptyset \, , \quad {\rm if} \ 
\langle \ell \rangle^\tau \leq c(s_0) \g_1 \slash (\sqrt{\d}) \, .
\ee
Indeed, denoting by  $ \mu'_j (\e,\l) \in \R $ 
the eigenvalues of $ D'_0=\Pi_{\mathbb F} A'_0 |_{H_{\mathbb F}} $, 
since $ \| A'_0-A_0 \|_0 \leq C(s_0) | A'_0-A_0 |_{+, s_0}  \leq C(s_0) \d $ (see \eqref{incl-plus}), we have  
\be \label{pert1}
 |\mu'_j(\e,\l)-\mu_j(\e,\l)| \leq C(s_0) \d \, , \quad \forall j \in {\mathbb F} \,  .
\ee
If  $ \l \notin \Lambda^1_\ell (\e;\cc, A'_0) $ then, recalling \eqref{Def:Lambda-l1}, there exist
$ i , j \in {\mathbb F} $ (with $i\neq j$ if $\ell=0$), 
signs $\ep_i=\pm 1$ and $\ep_j=\pm 1$ such that 
$$
| \bar \o_\e \cdot \ell + \ep_i  \mu'_i(\e,\l) + \ep_j  \mu'_j (\e, \l)| <  \frac{\g_1}{2 \cc \langle \ell \rangle^\tau} 
$$ 
and so, by \eqref{pert1}, 
$$
| \bar \o_\e \cdot \ell + \ep_i   \mu_i(\e,\l) +  \ep_j  \mu_j (\e, \l)| <  \frac{\g_1}{2 \cc \langle \ell \rangle^\tau} 
+ 2 C(s_0)  \d < \frac{\g_1}{2( \cc -\sqrt{\d}) \langle \ell \rangle^\tau} \, , 
$$
 for all $ (1/ 2) + \sqrt{\d} \leq \cc \leq 1 $,  
$ \langle \ell \rangle^\tau \leq c(s_0) \g_1 \slash \sqrt{\d} $, for some 
 $ c(s_0) $ small enough,  proving \eqref{prdis1}.
Hence
\begin{align*}
[\Lambda^1 (\e;\cc, A'_0)]^c  \bigcap \Lambda^1 (\e;\cc-\sqrt{\d}, A_0) & \subset 
 \bigcup_{\ell \in \Z^\es} [\Lambda^1_\ell (\e;\cc, A_0')]^c  
 \cap \Lambda^1 (\e;\cc-\sqrt{\d}, A_0)
\\
& \stackrel{\eqref{prdis1}} \subset
\bigcup_{| \ell | \geq (\frac{c(s_0) \g_1}{\sqrt{\d}})^{1/\tau}} [\Lambda^1_\ell  (\e;\cc, A'_0)]^c  
\end{align*}
which implies \eqref{midif0}
by 
\eqref{estsum}  with $ L = (\g_1 / \sqrt{\d} )^{1/\tau} $
(applied to $ A'_0 $) and Lemma \ref{lemma:measure1}. 

Finally, by \eqref{def>tau} and since $\tau_1 \geq \es $, we have  $\es /2\tau < 1/4$  
and so
\be\label{dis12} 
\min(\e^{-2} \d^{\frac12 -  \frac{\es}{2\tau} } , \e) \leq \min(\e^{-2} \d^{\frac14} , \e) 
\leq \d^{\frac{1}{12}}  
\ee
where the last inequality follows distinguishing the cases $ \e \leq \d^{\frac{1}{12}}$ and $ \e > \d^{\frac{1}{12}} $. 
By \eqref{midif0} and \eqref{dis12}, the estimate \eqref{Cantor-sovra1} is proved. 
\end{pf}


\smallskip

\noindent
{\bf Step 2: approximate solution of the homological equation \eqref{eqhoma}. } 
We decompose $ a \in {\cal L}(H_{\mathbb F} , H_{\mathbb G})$ as  $ a =(a^j)_{j \in {\mathbb F}} $   
and   $ \rho_2 \in {\cal L}(H_{\mathbb F} , H_{\mathbb G})$ as  $ \rho_2  =(\rho_2 ^j)_{j \in {\mathbb F}} $   
where $ a^j := a_{|H_{j}} \in {\cal L}(H_j, H_{\mathbb G}) $ and $ \rho_2^ j := (\rho_2)_{|H_j} \in {\cal L}(H_j, H_{\mathbb G}) $.  
 Recalling the form of $ D_0 $ in \eqref{form:D0}, the equation \eqref{eqhoma}
is equivalent to 
\be \label{eqhomai}
T_j (a^j) = J \rho_2^j \, , \quad  \forall j \in {\mathbb F} \, , 
\ee
where we define $T_j$ as the linear operator which maps  
$ a^j  :  \T^{\es} \to  {\cal L}(H_j, H_{\mathbb G}) $  to
\be\label{eqhomai1}
T_j (a^j) := J \ppavphi a^j -J V_0 J a^j + \mu_j (\e, \lambda) J a^j J  \, .
\ee
We shall consider $T_j$ as an unbounded linear operator
from the space $ L^2( \T^{\es}, {\cal L}( H_j, H_{\mathbb G} ) )$ to itself.

In the sequel $ i_{\mathbb G} $ denotes the injection 
\be\label{injectionG}
i_{\mathbb G} : H_{\mathbb G} \hookrightarrow  H \, . 
\ee
Before  applying the multiscale Proposition \ref{propmultiscale} we need to 
extend the linear operator $ T_j $   defined in \eqref{eqhomai1}, 
which acts on  $ L^2(\T^{\es}, {\cal L}(H_j,  H_{\mathbb G})) $, 
to a linear operator $ T^\sharp_j $  acting on the {\it whole}
space $ L^2(\T^{\es}, {\cal L}( H_j,  H)) 
\simeq L^2(\T^{\es}, H \times H)$ (by \eqref{ident-H-4-me})
and satisfying the properties of Definition  \ref{definition:Xr}.

We define, for any $ a^j \in L^2(\T^{\es}, {\cal L}( H_j, H)) $, the operator
\be\label{Ext10}
\begin{aligned}
& T^\sharp_j (a^j )   := \\
&  J \ppavphi  a^j   +  
\frac{D_V}{1+ \e^2 \l} \Pi_{\mathbb S \cup  \mathbb F} 
a^j + \frac{\co}{1+ \e^2 \l} \Pi_{ \mathbb S} a^j 
- Ji_{\mathbb G}  V_0   \Pi_{\mathbb G} J a^j + \mu_j (\e, \lambda) J \Pi_{\mathbb G} a^j J  
\end{aligned}
\ee
where $ \Pi_{\mathbb S \cup  \mathbb F}, \Pi_{\mathbb G}, \Pi_{\mathbb S}  $ 
are the $ L^2$-orthogonal projectors on the subspaces 
$ H_{\mathbb S \cup  \mathbb F}, H_{\mathbb G}, H_{\mathbb S} $ and $ \co > 0 $ 
is a positive constant that we fix according to \eqref{diof:co}. 
Clearly $ T^\sharp_j $  is an extension of $ T_j $, i.e.
$$ 
T^\sharp_j (a^j ) = T_j (a^j ) \, ,  \quad \forall  a^j \in 
L^2(\T^{\es}, {\cal L}(H_j,  H_{\mathbb G})) \subset 
L^2(\T^{\es}, {\cal L}(H_j,  H))\, .
$$ 
Recalling  the decomposition \eqref{ortho-operator}-\eqref{ort-deco-op},  i.e.
\be\label{deco:HjFSG}
\begin{aligned}
{\cal L}( H_j,   H ) & = {\cal L}( H_j,  H_{\mathbb S \cup  \mathbb F} ) \oplus 
{\cal L}( H_j,  H_{\mathbb G} ) \, , \\
a^j (\vphi) & = \Pi_{\mathbb S \cup  \mathbb F} a^j (\vphi) + \Pi_{{\mathbb G}} a^j  (\vphi) \, ,
\end{aligned}
\ee
we may write the operator $ T^\sharp_j $ in \eqref{Ext10} as
\be\label{Ext1}
T^\sharp_j (a^j )  =  J \ppavphi  a^j   +  
\frac{D_V}{1+ \e^2 \l  } a^j +  \mu_j (\e, \lambda) J \Pi_{\mathbb G} a^j J  +
 \frac{\co}{1+ \e^2 \l} \Pi_{ \mathbb S} a^j + {\cal R}_j a^j
\ee
where, using that $ J $ commutes with $ D_V $ and $ \Pi_{\mathbb G} $,  $ J^2 = - {\rm Id} $, 
\be\label{ext-wtildeT}
{\cal R}_j a^j := J \Big(  \frac{D_V}{ 1+ \e^2 \l } - i_{\mathbb G} V_0  \Big) \Pi_{\mathbb G} J a^j  =
- JR_0 \Pi_{\mathbb G} Ja^j 
\ee 
and we used \eqref{defA0}-\eqref{form:A0} to obtain the last equality.
Notice that,
according to the decomposition \eqref{deco:HjFSG}, 
the operator $ T^\sharp_j $ is represented by the block-diagonal matrix of operators
\be\label{TB1B2}
T^\sharp_j = 
\begin{pmatrix}
 J \ppavphi + \frac{D_V}{1+ \e^2 \l} + \frac{\co}{1+ \e^2 \l} \Pi_{\mathbb S}
 & 0    \\
 0 & T_j       
\end{pmatrix} 
\, ,
\ee
and, according to the further splitting 
$$
\begin{aligned}
 {\cal L}( H_j, H ) & = {\cal L}( H_j, H_{\mathbb S} ) \oplus {\cal L}( H_j, H_{\mathbb F} )  \oplus 
{\cal L}( H_j, H_{\mathbb G} ) \, , \\ 
 a^j (\vphi) & = \Pi_{\mathbb S} a^j (\vphi) +  \Pi_{\mathbb F} a^j  (\vphi) +  \Pi_{{\mathbb G}} a^j  (\vphi) \, , 
\end{aligned}
$$
recalling \eqref{Ext10}, by the matrix of operators
\begin{align} \label{three-split}
T_j^\sharp & 
 = J \ppavphi +
\left(
\begin{array}{ccc}
\frac{D_V}{1+ \e^2 \l} + \frac{\co \, {\rm Id}}{1+ \e^2 \l} &  0 & 0  \\
0  &   \frac{D_V}{1+ \e^2 \l} &  0 \\
0  &  0 &   -J V_0 J + \mu_j (\e, \lambda) {\cal J}
\end{array}
\right)
\end{align}
where ${\cal J} $ is defined in \eqref{def:opJ}.
Note also that 
\be\label{commTjTtj} 
T_j \Pi_{\mathbb G} = \Pi_{\mathbb G} T^\sharp_j \, .
\ee
Now, given 
$ g^j := J \rho_2^j : \T^\es \to  {\cal L} (H_j , H_{\mathbb G} ) $ 
we define
\be\label{extended-datum}
g_j^\sharp : \T^\es \to  {\cal L} (H_j , H ) \simeq H \times H \, , \quad 
g_j^\sharp (\vphi) := i_{\mathbb G} g^j (\vphi) \, , 
\ee
 (recall \eqref{ident-H-4-me}, \eqref{injectionG})
and we look for an (approximate) solution $ a_j^\sharp  : 
\T^\es \to {\cal L} (H_j , H) \simeq H \times H$ 
of  the equation 
\be \label{invTj}
 T^\sharp_j ( a_j^\sharp )= g_j^\sharp \, .
\ee
As already said, $ T^\sharp_j $ is an (unbounded) operator on 
$ L^2(\T^\es, {\cal L} (H_j,  , H))$ and an extension of $T_j$.
It is important to notice  that, since the subspaces   
$ L^2(\T^\es, {\cal L} (H_j,  H_{\mathbb S \cup \mathbb F}))$ and 
$ L^2(\T^\es, {\cal L} (H_j,  H_{{\mathbb G}})) $  are 
invariant under the operator $ T^\sharp_j$ (see e.g. \eqref{TB1B2}), a solution of  \eqref{invTj} 
satisfies $T_j ( \Pi_{\mathbb G} a_j^\sharp  ) = g_j^\sharp $, by \eqref{commTjTtj}, and therefore 
 $ a^j :=  \Pi_{\mathbb G}  a_j^\sharp $ solves  
equation \eqref{eqhomai} for the non extended operator $ T_j $. Actually in 
Lemma \ref{lem:homoeq2} we  find an approximate solution of the equation 
\eqref{eqhomai}. 

In the sequel 
$ T_j^\sharp $ is regarded as an operator acting  on $ L^2(\T^\es ; H \times H) $: in fact,
 using the basis $((\Psi_j,0), (0,\Psi_j))$ of $H_j$, see \eqref{Fj-eigenfunctions}, 
 we have the identification
\eqref{ident-H-4-me}
and  $ L^2(\T^\es, {\cal L}( H_j, H )) \simeq L^2(\T^\es ; H \times H )$. 
We  will apply the multiscale Proposition\index{Multiscale proposition} \ref{propmultiscale} 
 (in case-($ii$) in \eqref{phase-space-multi})
to the operator 
$$
{\cal L}_{r, \mu} = (1+ \e^2 \l) T^\sharp_j  \, .
$$ 
Actually, recalling \eqref{Ext1}, 
$ \Pi_{\mathbb G} =  \Pi_{{\mathbb F} \cup {\mathbb S}}^\bot $, 
and the definition of ${\cal J} $ in \eqref{def:opJ}, we have that  
\be\label{applicazione-ii}
\begin{aligned}
& {\cal L}_{r, \mu} = (1+ \e^2 \l) T^\sharp_j = J \om \cdot \partial_\vphi + X_{r, \mu} \, , \\
& 
X_{r, \mu} = D_V + \co \Pi_{\mathbb S} + \mu  {\cal J}  \Pi_{\mathbb S \cup \mathbb F}^\bot + r \, , 
\end{aligned}
\ee
has the form \eqref{def:Lr}  with $ \omega = (1+ \e^2 \l) \bar \om_\e $ and  
$X_{r, \mu} $ as in \eqref{def:Ar} with 
\be\label{rquRj}
\mu = (1+ \e^2 \l) \mu_j (\e, \lambda) \, ,
\qquad
r = (1+ \e^2 \l ) {\cal R}_j 	\stackrel{\eqref{ext-wtildeT}} 
= -(1+\e^2 \l) JR_0 \Pi_{\mathbb G}J  \, .
\ee
We now prove that  $ X_{r, \mu} $ in \eqref{applicazione-ii} satisfies the properties 
 stated in Definition \ref{definition:Xr},
beginning with its self-adjointness. 

\begin{lemma}\label{Ttilde-self-ad}
 {\bf (Self-adjointness)}
The operator  $ X_{r,\mu}  $ defined in  \eqref{applicazione-ii}-\eqref{rquRj} 
is  self-adjoint  with respect to the scalar product 
$ \langle \ \, , \ \rangle_0  $ in \eqref{sc-pr-operators-vphi}.
\end{lemma}

\begin{pf}
The self-adjointness of $D_V$ and $\Pi_{\mathbb S}$ with respect to 
$ \langle \ \, , \ \rangle_0  $ directly 
follows by the fact that  $ D_V $ and $\Pi_{\mathbb S}$ are $ L^2 $-
self-adjoint.
Now let $ a := a (\vphi ) $,  $ b  := b (\vphi ) $ belong to 
$ L^2(\T^\es, {\cal L}( H_j,  , H))$.  We obtain,
using that $\Pi^\bot_{{\mathbb S}\cup {\mathbb F}}$
is $ L^2$-self-adjoint, that 
$$
\begin{aligned}
\langle a, {\cal J}  \Pi^\bot_{{\mathbb S}\cup {\mathbb F}} b\rangle_0 
& \stackrel{\eqref{def:opJ}}  = \int_{\T^\es} {\rm Tr} (  (J \Pi^\bot_{{\mathbb S}\cup {\mathbb F}} b J  )^* a ) d \vphi \\ 
& = \int_{\T^\es}
{\rm Tr} (   J b^* \Pi^\bot_{{\mathbb S}\cup {\mathbb F}} J a  ) d \vphi  = \int_{\T^\es} {\rm Tr} (  J   b^*   J  \Pi^\bot_{{\mathbb S}\cup {\mathbb F}} a ) d \vphi
\end{aligned}
$$
because  $J $ and $\Pi^\bot_{{\mathbb S}\cup {\mathbb F}} $ commute. Thus, using that
$ {\rm Tr}(AB) = {\rm Tr}(BA) $, we deduce that 
$$
\begin{aligned}
\langle a, {\cal J}  \Pi^\bot_{{\mathbb S}\cup {\mathbb F}} b \rangle_0 
& = \int_{\T^\es} {\rm Tr} (    b^*   J  \Pi^\bot_{{\mathbb S}\cup {\mathbb F}} a J ) d \vphi 
 \stackrel{\eqref{def:opJ}, \eqref{sc-pr-operators-vphi}} = 
\langle   {\cal J}  \Pi^\bot_{{\mathbb S}\cup {\mathbb F}}  a,  b \rangle_0 \, .  
\end{aligned}
$$
Hence the operator ${\cal J}  \Pi^\bot_{{\mathbb S}\cup {\mathbb F}} $ is self-adjoint
 with respect to $ \langle \ \, , \ \rangle_0  $. 

 It remains to prove that the operator $ {\cal R}_j $ in \eqref{ext-wtildeT} is 
 self-adjoint  with respect to $ \langle \ \, , \ \rangle_0  $. 
Since $R_0$ is self-adjoint and $R_0 (H_{\mathbb G}) \subset H_{\mathbb G}$, we deduce that $ B := R_0 \Pi_{\mathbb G} $ is $ L^2 $-self-adjoint, i.e. it satisfies  
$B^*=B$, and, recalling \eqref{ext-wtildeT}, we obtain
$$
\langle a, {\cal R}_j  b \rangle_0  
= -\int_{\T^\es} {\rm Tr} (  (J B  J b)^* a) d \vphi =- \int_{\T^\es} {\rm Tr} (  b^* J B  J a ) d \vphi
 = \langle {\cal R}_j  a, b \rangle_0 \, .  
 $$
This completes the proof of the  lemma. 
\end{pf}

\begin{lemma}\label{lem:off-diag}
  {\bf (Off-diagonal decay)}
The operator ${\cal R}_j $ defined in \eqref{ext-wtildeT} satisfies 
\be\label{lem:9.10}
|{\cal R}_j|_{\Lip,s} \leq C(s) |R_0|_{\Lip,s}  \, ,  \qquad 
|{\cal R}_j|_{\Lip, +,s} \leq C(s) |R_0|_{\Lip, +,s} \, ,
\ee
in particular $ |{\cal R}_j|_{\Lip, +,s_1} \leq C_1' \e^2 $ for some constant  
$ C_1' $ depending only on $C_1$.
\end{lemma} 

\begin{pf}
We identify $ a(\vphi)  \in {\cal L} ( H_j, H ) $ with the vector $ (a^{(1)}(\vphi), a^{(2)}(\vphi), a^{(3)}(\vphi), a^{(4)}(\vphi)) $ in $  H \times H $
as in  \eqref{ident-H-4}. Then,  using \eqref{actionJ}, we have 
$$
{\cal R}_j a = \big( -JR_0 \Pi_{\mathbb G} (a_2,-a_1) , -JR_0 \Pi_{\mathbb G} (a_3, -a_4) \big) \, . 
$$
Hence  
$$ 
|{\cal R}_j|_{\Lip,s} \sim | R_0 \Pi_{\mathbb G}|_{\Lip,s} \quad {\rm and } \quad 
|{\cal R}_j|_{\Lip, +,s} \sim | R_0 \Pi_{\mathbb G}|_{\Lip, +,s}  \, . 
$$
Lemmata   \ref{normprod} and \ref{pisig}  imply that  
$$ 
\begin{aligned}
& |R_0 \Pi_{\mathbb G}|_{\Lip,s} 
\stackrel{\eqref{tame-s-decayLip}} {\lesssim_s} |\Pi_{{\mathbb G}}|_{\Lip, s} |R_0|_{\Lip,s} 
\stackrel{\eqref{sepr2}} {\lesssim_s} |R_0|_{\Lip,s}  \\ 
& |R_0 \Pi_{\mathbb G}|_{\Lip, +,s} 
\stackrel{\eqref{inter-norma+s}}
{\lesssim_s} |\Pi_{{\mathbb G}} |_{\Lip, s+ \frac12} |R_0|_{\Lip, +,s} 
\stackrel{\eqref{sepr2}} {\lesssim_s} |R_0|_{\Lip, +,s} 
\end{aligned}
$$
proving \eqref{lem:9.10}. Since $ |R_0|_{\Lip, +,s_1} \leq C_1 \e^2 $
by Definition \ref{def:calC}, the second estimate in \eqref{lem:9.10} implies 
$ |{\cal R}_j|_{\Lip, +,s_1} \leq C_1' \e^2 $.
\end{pf}

\begin{lemma}\label{pos:def-var}
 {\bf (Sign condition)}
The operator $ X_{r, \mu}$ defined in \eqref{applicazione-ii} satisfies, for some $ c >  0 $ depending on the constant $c_1$ in \eqref{Hyp1},  
$$ 
{\mathfrak d}_\l \Big(\frac{X_{r, \mu} }{1+\e^2 \l}\Big) \leq - c \, \e^2 {\rm Id} \, .
$$
\end{lemma}

\begin{pf}
According to \eqref{three-split}
the operator  $ \frac{X_{r, \mu} }{1+\e^2 \l} $ is represented by the matrix of operators
$$
\left(
\begin{array}{ccc}
\frac{D_V + \co \, {\rm Id}}{1+ \e^2 \l} &  0 & 0  \\
0  &    \frac{D_V}{1+ \e^2 \l} &  0 \\
0  &  0 &   -J V_0  J + \mu_j (\e, \lambda) \,  {\cal J}
\end{array}
\right) 
$$
where  $ {\cal J} a = JaJ $.
We clearly have 
\be\label{c1c2}
\partial_\l  \frac{D_V + \co \, {\rm Id}}{1+ \e^2 \l}  
 = -\frac{\e^2 (D_V + \co \, {\rm Id})}{(1+\e^2 \l)^2}  \leq -\frac{\e^2 D_V}{(1+\e^2 \l)^2}   
\leq - c\e^2 {\rm Id}
\ee
for some $c>0$. Then it is sufficient to prove that, for all 
$ a \in L^2(\T^\es, {\cal L}(H_j,  H_{{\mathbb G}}))$, 
and for all $\l_1, \l_2 \in \wtilde{\Lambda}$ with $\l_1 \neq \l_2$,  using the notation \eqref{part-quo}, 
$$
\Big\la -J \frac{\Delta V_0}{\Delta \l}   J a + \frac{\Delta  \mu_j (\e, \lambda)}{\Delta \l}  \, J aJ,a \Big\ra_0 \leq 
- c \e^2 \| a \|_0^2 \, ,
$$
for the scalar product $\langle \  ,  \  \rangle_0$ defined  in \eqref{sc-pr-operators-vphi}.
We have 
\begin{align} 
\Big\la -J \frac{\Delta V_0}{\Delta \l}   J a + \frac{\Delta  \mu_j (\e, \lambda)}{\Delta \l}  \, J aJ,a \Big\ra_0
\nonumber
&=  \Big\la  \frac{\Delta V_0}{\Delta \l}    J a, Ja  \Big\ra_0 - 
\frac{\Delta  \mu_j (\e, \lambda)}{\Delta \l}  \, \big\la aJ, J a \big\ra_0 \\ \nonumber
&\leq  \Big\la  \frac{\Delta V_0}{\Delta \l}   J a, Ja  \Big\ra_0 +
\Big| \frac{\Delta  \mu_j (\e, \lambda)}{\Delta \l} \Big| \| aJ \|_0  \| J a \|_0 \\
& =   \Big\la  \frac{\Delta V_0}{\Delta \l}   J a, Ja  \Big\ra_0 +
\Big| \frac{\Delta  \mu_j (\e, \lambda)}{\Delta \l} \Big| \| J a \|_0^2   
\label{pacru}
\end{align}
noting that $ \| aJ \|_0 = \| Ja \|_0 = \| a \|_0$
from the definition \eqref{sc-pr-operators-vphi} of the scalar product 
$\langle \  ,  \  \rangle_0 $. 
Let 
$$
(Ja)^{(1)} (\vphi)=J a(\vphi) (\Psi_j,0) \, , \quad 
(Ja)^{(2)} (\vphi)=J a(\vphi) (0, \Psi_j) \, , \quad 
(Ja)^{(i)} \in L^2(\T^\es , H_{\mathbb G}) \, . 
$$
Using the notations $\la \ , \ \ra_0$ and $\|\ \|_0$  for the scalar product and its associated norm both in 
$ L^2(\T^\es, {\cal L}(H_j, H_{{\mathbb G}}))$ and in  
$L^2(\T^\es ; H_{{\mathbb G}})$, we obtain by \eqref{pacru}
\be  \label{positivityf}
\begin{aligned}
& \Big\la -J \frac{\Delta V_0}{\Delta \l}   J a + \frac{\Delta  \mu_j (\e, \lambda)}{\Delta \l}  \, J aJ,a \Big\ra_0 \\
& \leq 
\sum_{i=1}^2  \Big\la  \frac{\Delta V_0}{\Delta \l}   (J a)^{(i)}, (Ja)^{(i)} \Big\ra_0 +
\Big| \frac{\Delta  \mu_j (\e, \lambda)}{\Delta \l} \Big| \| (J a)^{(i)} \|_0^2 \, . 
\end{aligned}
\ee
Now the assumption \eqref{Hyp1} implies that for all $h \in L^2(\T^\es , H_{\mathbb G})$
\be \label{positivity3} 
\Big\la  \frac{\Delta V_0}{\Delta \l}  h, h  \Big\ra_0 + 
\Big| \frac{\Delta  \mu_j (\e, \lambda)}{\Delta \l}\Big| \| h \|_0^2 \leq 
- c_1 \e^2 \| h \|_0^2 \, .
\ee
The estimates \eqref{c1c2},  \eqref{positivityf} and \eqref{positivity3}
imply the lemma. 
\end{pf}

\smallskip

By Lemmata \ref{Ttilde-self-ad}, \ref{lem:off-diag}, \ref{pos:def-var}
we  apply the multiscale Proposition \ref{propmultiscale} 
to the operator $ {\cal L}_{r, \mu} = (1+ \e^2 \l) T^\sharp_j $ 
in \eqref{applicazione-ii} where $ r $ is given in \eqref{rquRj}. 
As a consequence
there exist, for any $ j \in {\mathbb F} $,  closed subsets   
\be\label{def:setsLj}
\Lambda^2_j (\e; \cc, A_0) \, , \   1/2 \leq \eta \leq 1 \, , \ {\rm  satisfying \ Properties} \, 
1{\mathrm -}3 \ {\rm of \ Proposition \ \ref{prophomeq}} \, , 
\ee
and $ \bar N \in \N $, such that, for all $ \forall N \geq \bar N $, 
for all $\l \in \Lambda^2_j (\e;1, A_0) $, there are 
operators $ (T_j^\sharp)_N^{-1}  $ defined in Proposition \ref{propmultiscale} as
\be\label{TN_tr}
\begin{aligned} 
& i)  \ {\rm the \  right \ inverse \ of } \ {\it \Pi}_N (T_j^\sharp)_{|{\mathcal H}_{2N}} \  \hbox{} \ {\rm if}  \ \bar N \leq N < N(\e) \, , \\
& ii)  \ {\rm the \ inverse \ of } \ {\it \Pi}_N (T_j^\sharp)_{|{\mathcal H}_N}  \ {\rm if}  \ N \geq N(\e) \, ,
\end{aligned}
\ee
where  $ {\mathcal H}_N $ are the 
finite dimensional subspaces  defined in \eqref{def:EN-tr}.  
By \eqref{est:Right+1}, \eqref{Lip-sDM}, \eqref{rquRj}-\eqref{lem:9.10} and
$ |R_0|_{\Lip, +,s_1} \leq C_1 \e^2 $
by Definition \ref{def:calC}, 
 the operators $ (T_j^\sharp)_N^{-1}  $ satisfy 
the following tame estimates: 
$ \forall s \geq s_0 $, 
\begin{align} 
&  |  (T_j^\sharp)_N^{-1} |_{\Lip, s_1}  \leq C(s_1)   N^{2(\tau' + \loss s_1 ) + 3} \label{Tn-1Lips1} \\
& | (T_j^\sharp)_N^{-1} |_{\Lip, s} 
\leq C(s)   N^{2(\tau' + \loss s_1) +3}  \big(   N^{\loss (s - s_1)} + | R_0 |_{\Lip, +,s} \big)   \label{Tn-1Lipss} \, .
\end{align}
To justify that the sets $\Lambda^2_j (\e; \cc, A_0) $ satisfy Property $3$ of Proposition
\ref{prophomeq}, we remark that if the operators $A_0$ and $A'_0$ satisfy
$|R'_0-R_0|_{+,s_1} \leq \d \leq \e^{5/2}$, then 
$\mu , \mu', r,r'$ defined in \eqref{rquRj} satisfy 
$$
|r'-r|_{+,s_1} + |\mu'-\mu| \leq C\d \leq \e^2 \, ,
$$
and refer to \eqref{prop:spos} in Proposition \ref{propmultiscale} 
(note that the exponent $\alpha$ of Property $3$ of Proposition
\ref{prophomeq} may be taken
slightly smaller than in \eqref{prop:spos}).

We now obtain an approximate solution of the homological equation 
\eqref{eqhomai}, finding an approximate solution of 
\eqref{invTj}. 

\begin{lemma} {\bf (Homological equation \eqref{eqhomai})}\label{lem:homoeq2}
Let
\be \label{defNsp} 
N \in \Big[ \d_1^{-\frac{3}{s_2-s_1}} -1, \d_1^{-\frac{3}{s_2-s_1}} +1 \Big]\, . 
\ee 
Then, for all $ j \in {\mathbb F} $, 
for all $\l \in \Lambda^2_j (\e;1, A_0) $ $($set introduced in \eqref{def:setsLj}$)$, the function $ a^j $ defined as  
\be \label{homdefa} 
a^j := \Pi_{\mathbb G}  b_j^\sharp   \, ,   \ \ 
b_j^\sharp := (T_j^\sharp)_N^{-1} g_{j,N}^\sharp \, , \ \   g_{j,N}^\sharp := {\it \Pi}_N g_j^\sharp \, , \ \ 
g_j^\sharp := i_{\mathbb G}  g^j \, , \  g^j = J (\rho_2)_{|H_j}  
\ee
(recall the notation  \eqref{extended-datum}, \eqref{decomposition-rho1}), satisfies 
\begin{align}\label{stima-di-a}
& | a^j |_{\Lip, s_1+1}  \ll  \d_1^{\frac78} \, , \\
& | a^j |_{\Lip, s+1}  \lesssim_s   \d_1^{- \frac{1}{40}} \big( 
|\rho|_{\Lip, s} + | R_0 |_{\Lip, +,s} \big) + 
\d_1^{\frac12} \d_1^{- 3 \varsigma \frac{s-s_2}{ s_2 - s_1 }} \, , \ \  \forall s \geq s_2 \, , \label{stima-di-a1}
\end{align}
and it  is an approximate solution of the homological equation \eqref{eqhomai} in the sense that  
\be\label{approx-sol-2}
 |T_j a^j -g^j |_{\Lip, s_1+1}  \ll \d_1^{7/4}  \, . 
 \ee 
At last, denoting by 
$ a^j := a^j_{ A_0, \rho} $ and $ (a^j)' := a^j_{A'_0, \rho'} $ the operators 
defined as above 
associated to $ (A_0, \rho) $ and  $ (A'_0, \rho') $ respectively, 
for all $\l \in \Lambda^2 (\e; 1, A_0) \cap \Lambda^2 (\e; 1, A'_0)$, we have 
\be \label{perta}
\big|  a^j_{ A_0, \rho} - a^j_{A'_0, \rho'} \big|_{s_1 +1} \ll 
\d_1^{\frac45} |A_0 - A'_0|_{+, s_1}  + \d_1^{- \frac{1}{40}} |\rho -\rho'|_{+, s_1} \,   
\ee
and
\be \label{perteq2}
\big| \big( T_{j, A_0}  a^j_{A_0,\rho} -g^j_\rho \big) - 
\big( T_{j,A'_0} a^j_{A'_0,\rho'} -g^j_{\rho'} \big) \big|_{s_1+ \frac12} 
\ll \d_1^{\frac45} |A_0 - A'_0|_{+, s_1}  + \d_1^{- \frac{1}{40}} |\rho -\rho'|_{+, s_1}  \, .
\ee
\end{lemma}

\begin{pf}
The function $ b_j^\sharp $ defined in  \eqref{homdefa} satisfies, by \eqref{TN_tr}, 
the property 
\be\label{b-cases}
\begin{aligned} 
& i)  \ {\it \Pi}_{2N} b_j^\sharp = b_j^\sharp \quad  {\rm if} \quad  \bar N \leq N < N(\e) \\
& ii)\  {\it \Pi}_{N} b_j^\sharp = b_j^\sharp   \quad {\rm if}  \quad N \geq N(\e) \, ,
\end{aligned}
\ee
and it is a solution of $ {\it \Pi}_N T_j^\sharp b_j^\sharp = {\it \Pi}_N g_{j}^\sharp $. Then,
up to terms  which are Fourier supported on high harmonics, it is an approximate solution of the equation 
\eqref{invTj},  more precisely, 
\be\label{svizN-1}
T_j^\sharp b_j^\sharp = {\it \Pi}_N T_j^\sharp b_j^\sharp + {\it \Pi}_N^\bot T_j^\sharp b_j^\sharp = 
{\it \Pi}_N g_{j}^\sharp + {\it \Pi}_N^\bot
T_j^\sharp b^\sharp_j = g_j^\sharp + z_N 
\ee
where 
\be\label{svizN1}
 z_N := - {\it \Pi}_N^\bot g_j^\sharp + {\it \Pi}_N^\bot T_j^\sharp b_j^\sharp \, .
\ee
Applying $ \Pi_{\mathbb G}  $ to both sizes in \eqref{svizN-1}, we deduce,
by \eqref{commTjTtj} and \eqref{homdefa}, 
 that $ a^j := \Pi_{\mathbb G}  b_j^\sharp $ is an approximate solution of the homological equation 
\eqref{eqhomai},  in the sense that  
\be\label{soluzione-approssimata-vera2}
T_j a^j - J (\rho_2)_{|H_j}  =  
\Pi_{\mathbb G} z_N \, , \qquad J (\rho_2)_{|H_j} = J \rho_2^j = g^j  \, . 
\ee
We now prove \eqref{stima-di-a}-\eqref{approx-sol-2}.
\\[1mm]
{\sc Estimates of $ b_j^\sharp $.} By \eqref{ident-H-4-me} we identify 
the $ \vphi $-dependent family of operators $ g_j^\sharp = i_{\mathbb G} J (\rho_2)_{|H_j} $ in $ L^2(\T^\es , {\cal L} (H_j, H) )$ defined in \eqref{homdefa}
with a function of $ L^2(\T^\es , H \times H )$. 
We have the equivalence of the norms 
\be\label{normeqs} 
| g_j^\sharp |_{s} \sim_s \| g_j^\sharp \|_{s}  \, . 
\ee
In fact, let us define, for all $ \vphi \in \T^\es $,  the functions 
$$
g_j^{\sharp , 1} (\vphi) := g_j^\sharp  (\vphi)\big[ (\Psi_j,0)\big] \in H \quad
{\rm and}  \quad 
g_j^{\sharp , 2} (\vphi) := g_j^\sharp  (\vphi) \big[ (0,\Psi_j)\big] \in H \, . 
$$
We have to prove that  $| g_j^\sharp |_{s} \sim_s \| g_j^{\sharp , 1}  \|_{s} + \| g_j^{\sharp , 2}  \|_{s}$.
By Definition \ref{def:decay-sub}, 
$ | g_j^\sharp |_{s}=| \wtilde{g_j^\sharp} |_{s}$, where the operator
$ \wtilde{g_j^\sharp}$ is defined on the whole $L^2(\T^\es, H )$ by 
$$
 \wtilde{g_j^\sharp} \big[ (h^{(1)}, h^{(2)}) \big] := \sum_{i=1}^2
 g_j^{\sharp , i} \big( h^{(i)},\Psi_j \big)_{L^2_x} \, .
$$
Now, by Lemma  \ref{gchi-r}, 
$$
| \wtilde{g_j^\sharp} |_{s} \leq C(s) 	\big( 
\| g_j^{\sharp , 1}  \|_{s} + \| g_j^{\sharp , 2} \|_{s} \big) \, . 
$$
Conversely
$$
\| g_j^{\sharp , 1}  \|_{s}= \| \wtilde{g_j^\sharp} \big[ (\Psi_j , 0) \big] \|_s
\lesssim_s | \wtilde{g_j^\sharp} |_{s} \| \Psi_j \|_s \lesssim_s | \wtilde{g_j^\sharp} |_{s}
$$
and similarly $\| g_j^{\sharp , 2} \|_s \lesssim_s | \wtilde{g_j^\sharp} |_{s}$. This
proves the norm equivalence \eqref{normeqs}.

\smallskip
We have  $ g_{j,N}^\sharp  = {\it \Pi}_N g_{j,N}^\sharp $,  
and\index{Smoothing operators} \eqref{smoothingS1S2-Lip}, \eqref{homdefa},   \eqref{Tn-1Lips1}, \eqref{prop:rho-plus},
imply
\begin{align}
| b_j^\sharp |_{\Lip, s_1+1} & 
\leq N | b_j^\sharp |_{\Lip, s_1} \nonumber \\
& \lesssim_{s_1} N  |(T_j^\sharp)_N^{-1}|_{\Lip, s_1} 
| g_{j,N}^\sharp |_{\Lip, s_1} \nonumber \\
&  \lesssim_{s_1}  N^{2(\tau'+ \loss  s_1 + 2)}  | \rho |_{\Lip, s_1}  
 \lesssim_{s_1}  N^{Q} \d_1^{\frac{9}{10}}   \label{estimate:bs1}
\end{align}
having set 
\be\label{def:Q}
Q:= 2 (\tau'+ \loss s_1 + 2) \quad {\rm with} \quad \loss = 1/10 \quad {\rm  (as \ in \ \eqref{def:varsigma})} \, . 
\ee
Similarly, for $ s \geq s_1 $, using \eqref{smoothingS1S2-Lip}, \eqref{homdefa},  
the tame estimate \eqref{tame-s-decayLip}, 
  \eqref{Tn-1Lips1}-\eqref{Tn-1Lipss}, \eqref{prop:rho-plus}, we get 
\begin{align}
| b_j^\sharp |_{\Lip, s+1} & \leq N | b_j^\sharp |_{\Lip, s} \nonumber \\
& \lesssim_s N  |(T_j^\sharp)_N^{-1}|_{\Lip, s} | g_{j,N}^\sharp |_{\Lip, s_1} +
N  |(T_j^\sharp)_N^{-1}|_{\Lip, s_1} | g_{j,N}^\sharp |_{\Lip, s} \nonumber \\
& \lesssim_s N^{Q} \big( N^{\loss (s- s_1)} + |R_0|_{\Lip,+,s} \big) \d_1^{\frac{9}{10}} +
N^{Q}  | g_{j,N}^\sharp |_{\Lip, s} \nonumber \\
&\lesssim_s  N^{Q} \big(
N^{\loss(s-s_1)} \d_1^{\frac{9}{10}}  +  |R_0|_{\Lip, +,s} \d_1^{\frac{9}{10}}  
+  | \rho |_{\Lip,s} \big)  \, .
\label{estimate:bs}
\end{align}
{\sc Estimates of $ z_N $ defined in \eqref{svizN1}.}
We first claim that 
\begin{align}\label{cancellazione-bot} 
& |{\it \Pi}_N^\bot T_j^\sharp {\it \Pi}_{N} |_{\Lip, s} \leq 
|{\it \Pi}_N^\bot T_j^\sharp  {\it \Pi}_{2N} |_{\Lip, s} \lesssim_s N+|R_0|_{\Lip, s} \, . 
\end{align}
Indeed, recalling the expression \eqref{Ext1} of  $ T_j^\sharp $,  
setting  $ {\cal J} (a^j ) := J  a^j J  $, and writing $ D_V = D_m + (D_V - D_m) $, 
we get 
\be
\begin{aligned}
{\it \Pi}_N^\bot T_j^\sharp {\it \Pi}_{2N} & =  {\it \Pi}_N^\bot  J \ppavphi  {\it \Pi}_{2N} 
+ \frac{1}{1+ \e^2 \l} \big( {\it \Pi}_N^\bot D_m {\it \Pi}_{2N} + {\it \Pi}_N^\bot ( D_V - D_m) {\it \Pi}_{2N} \big) \\
& \quad +  \mu_j (\e, \lambda)  {\it \Pi}_N^\bot {\cal J} \Pi_{\mathbb G} {\it \Pi}_{2N}  + {\it \Pi}_N^\bot {\cal R}_j {\it \Pi}_{2N} +
{\it \Pi}_N^\bot  \frac{\co}{1+ \e^2 \l} \Pi_{ \mathbb S} {\it \Pi}_{2N} 
 \label{deco-tre-pezzi} \, .
\end{aligned}
\ee
In view of \eqref{deco-tre-pezzi}, 
Lemmas \ref{pisig}, \ref{lem:off-diag} and Proposition \ref{lemma-DVvsDm} imply \eqref{cancellazione-bot}. 

Then, by \eqref{svizN1}, \eqref{b-cases},  \eqref{homdefa}, \eqref{cancellazione-bot},
 \eqref{tame-s-decay}, $ | R_0|_{\Lip,s_1} \leq 1 $, 
we obtain
\begin{align}
|z_N |_{\Lip, s} &\leq   | g_j^\sharp |_{\Lip, s} + | {\it \Pi}_N^\bot T_j^\sharp b_j^\sharp |_{\Lip, s} \nonumber \\
&  \lesssim_s  | \rho |_{\Lip, s} + 
N | b_j^\sharp |_{\Lip, s} + (N+ |R_0|_{\Lip, s}) | b_j^\sharp |_{\Lip, s_1} 
\nonumber \\
& 
\stackrel{\eqref{estimate:bs}, \eqref{estimate:bs1}} 
{\lesssim_s} N^{Q} \big( | \rho |_{\Lip, s} +  |R_0|_{\Lip, +,s} \d_1^{\frac{9}{10}}  + 
N^{\loss(s-s_1)} \d_1^{\frac{9}{10}} \big) \, . \label{estimate:zN}
\end{align}
By \eqref{svizN1} we have  $ {\it \Pi}_N z_N  = 0 $, and, 
using \eqref{estimate:zN} and the assumption \eqref{prop:rho-plus}, we derive that
\begin{align}\label{estimatezN:bs1}
|z_N|_{\Lip, s_1+1} & \leq N^{-(s_2-s_1-1)} |z_N |_{\Lip, s_2} \nonumber  \\
&  \lesssim_{s_2} 
 N^{Q+1-(s_2-s_1)} \d_1^{- \frac{11}{10}} + N^{Q+1 -(1-\loss) (s_2 -s_1)} \d_1^{\frac{9}{10}}  \, .
\end{align}
{\sc Proof of \eqref{stima-di-a}-\eqref{approx-sol-2}.} 
By the choice of $ N $ in \eqref{defNsp}, the condition \eqref{choice s2 s3}-($i$), and \eqref{def:Q} we get 
\begin{align}
& N^{Q+1}=o(\d_1^{- \frac{1}{40}}) \ , \label{stN1}\\
& N^{Q+1-(s_2-s_1)} \d_1^{- \frac{11}{10}} = o(\d_1^{ \frac74}) \, , \  
N^{Q+1 - (1- \loss) (s_2-s_1)} \d_1^{\frac{9}{10}} = 
o(\d_1^{\frac74}) \, , \label{stN2} \\
& \forall s \geq s_2 \, , \quad N^{Q+1 + \loss (s-s_1)} \d_1^{\frac{9}{10}} = N^{Q+1+ \loss (s_2-s_1)} \d_1^{\frac{9}{10}} N^{ \loss (s-s_2)}
\leq \d_1^{\frac12} \d_1^{-3 \loss \frac{s-s_2}{s_2-s_1}} \, . \label{stN3} 
\end{align}
Then, by \eqref{estimate:bs1},  \eqref{stN1} and  \eqref{estimatezN:bs1}, \eqref{stN2},  
for $\d_1$ small enough, we deduce the bounds 
\be\label{stima-di-b}
| b_j^\sharp |_{\Lip, s_1+1}  \ll  \d_1^{7/8},  \qquad  
|z_N|_{\Lip, s_1+1} \ll  \d_1^{7/4} \, . 
\ee
In addition, \eqref{estimate:bs}, \eqref{stN1}, \eqref{stN3} imply  the estimate in high norm 
\be\label{estimate:bs-high}
 | b_j^\sharp |_{\Lip, s+1}  \leq   \d_1^{- \frac{1}{40}} \big( |\rho|_{\Lip, s} +
|R_0|_{\Lip, +,s} \big) + \d_1^{\frac12} \d_1^{- 3 \varsigma 	\frac{s-s_2}{ s_2 - s_1 }} \, . 
\ee
By  Lemma \ref{pisig}, the function
$ a^j := \Pi_{\mathbb G} b_j^\sharp $ satisfies the same estimates as $ b_j^\sharp $ in  \eqref{stima-di-b}, 
\eqref{estimate:bs-high}.  
In particular \eqref{stima-di-a}-\eqref{stima-di-a1} hold and, by \eqref{soluzione-approssimata-vera2}, 
Lemma \ref{pisig},  \eqref{stima-di-b}, we get  
$$
 |T_j a^j -g^j|_{\Lip, s_1+1} =   |\Pi_{\mathbb G} z_N |_{\Lip, s_1+1}
\ll \d_1^{7/4} 
 $$ 
 proving \eqref{approx-sol-2}.
\\[1mm]  
{\bf Proof of \eqref{perta}-\eqref{perteq2}.}
We denote by $ ({a^j})'$, ${b^\sharp_j}'$ the functions 
obtained in \eqref{homdefa} from $(A'_0, \rho')$ instead of $ A_0, \rho $, and, similarly,  
by ${T^\sharp_j}' $ the  linear operator defined in \eqref{Ext1} from $(A_0', \rho') $. 
Recall that to define $T^\sharp_j$ we applied Proposition \ref{propmultiscale} with
$\mu$ and $ r $ defined in \eqref{rquRj}.  In particular, we have
$$
| \mu - \mu'| \lesssim   |A_0-A'_0|_{+,s_1 }  \, ,  \quad
|r-r'|_{+,s_1}  \lesssim   |R_0-R'_0|_{+,s_1 } =  |A_0-A'_0|_{+,s_1 } \, .
$$
Hence, by \eqref{MSLip1}, \eqref{MSLip2} in Proposition  \ref{propmultiscale},
for any $\l \in \Lambda^2_j (\e; 1, A_0) \cap   \Lambda^2_j (\e; 1, A'_0)$, we have
\begin{equation}  \label{ultima-differenza}
\big| (T_j^\sharp)_N^{-1} - ({T_j^\sharp}')_N^{-1} \big|_{s_1}
 \leq   N^{2(\tau' + \varsigma s_1) +3  } |A_0-A'_0|_{+, s_1}  \, . 
\end{equation}
By \eqref{homdefa} and Lemma \ref{pisig}  we have 
 \begin{align}
 |({a^j})' - a^j |_{s_1+1} & \lesssim_{s_1}  |{b_j^\sharp}' - b_j^\sharp|_{s_1+1} \nonumber \\
 & \stackrel{\eqref{b-cases}, \eqref{tame-s-decay}} {\lesssim_{s_1}} 
   N \big| ({T_j^\sharp}')_N^{-1} - 
 (T_j^\sharp)_N^{-1} \big|_{s_1} |{g_{j,N}^{\sharp'}}|_{s_1} +
 N \big| (T_j^\sharp)_N^{-1} \big|_{s_1} 
 | {g_{j,N}^{\sharp'}} -{g_{j,N}^\sharp}  |_{s_1} \nonumber  \\
 & \stackrel{\eqref{ultima-differenza}, \eqref{Tn-1Lips1},  \eqref{prop:rho-plus}} 
 {\lesssim_{s_1}}   N^{2(\tau' + \varsigma s_1 )+4} |A'_0 - A_0|_{+, s_1} \d_1^{\frac{9}{10}} +
 N^{2(\tau' + \varsigma s_1 )+4} |\rho-\rho'|_{s_1}  \nonumber \\
 & \stackrel{\eqref{def:Q}}  {\lesssim_{s_1}}   N^{Q} \d_1^{\frac{9}{10}} |A'_0 - A_0|_{+, s_1}  +
 N^{Q} |\rho-\rho'|_{+, s_1}   \label{disuvoluta}  \\
 &\stackrel{\eqref{stN1}}{ \ll}  \d_1^{\frac78} |A'_0-A_0|_{+, s_1}  +
 \d_1^{- \frac{1}{40}} |\rho-\rho'|_{+, s_1}  \nonumber
 \end{align}
 proving \eqref{perta}. 
In addition 
$$ 
 T_j^\sharp b_j^\sharp = g_j^\sharp + z_N \, ,  \quad 
 {T_j^\sharp}' {b_j^\sharp}' = {g_j^\sharp}' + z'_N \, , 
 $$    
where, by \eqref{svizN1},  
\begin{align}
|z_N-z'_N|_{s_1+ \frac12}  &  =  
\big| \big(- {\it \Pi}_N^\bot g_j^\sharp + {\it \Pi}_N^\bot T_j^\sharp b_j^\sharp \big) -
\big( - {\it \Pi}_N^\bot {g_j^\sharp}' + 
{\it \Pi}_N^\bot {T_j^\sharp}' {b_j^\sharp}' \big) \big|_{s_1+ \frac12} \nonumber \\
 & \leq    |{g_j^\sharp}-{g_j^\sharp}'|_{s_1+ \frac12} + 
\big| {\it \Pi}_N^\bot (T_j^\sharp -{T_j^\sharp}') {b_j^\sharp}   \big|_{s_1+ \frac12} + 
\big|  {\it \Pi}_N^\bot {T_j^\sharp}' {\it \Pi}_{2N} (b_j^\sharp - {b_j^\sharp}')  \big|_{s_1+ \frac12} \nonumber \\
&\stackrel{\eqref{Ext1}, \eqref{ext-wtildeT}, 
\eqref{cancellazione-bot}}{\lesssim_{s_1}}   |{g_j^\sharp}-{g_j^\sharp}'|_{s_1+ \frac12} + 
\big| A_0-A'_0 \big|_{s_1 + \frac12}  \big|{b_j^\sharp}   \big|_{s_1+ \frac12} \nonumber  \\ 
& \qquad \qquad \qquad \quad + 
\big( N+ \big|  R_0 \big|_{\Lip , s_1+\frac12 } \big) \big| b_j^\sharp - {b_j^\sharp}'\big|_{s_1+\frac12}\nonumber \\
&  \lesssim_{s_1} |{g_j^\sharp}-{g_j^\sharp}'|_{s_1+ \frac12} + 
\big| A_0-A'_0 \big|_{+,s_1 }  \big|{b_j^\sharp}   \big|_{s_1+ \frac12} + 
(N+ \big|  R_0 \big|_{\Lip , +,s_1 } ) 
\big| b_j^\sharp - {b_j^\sharp}'\big|_{s_1+\frac12}\nonumber \\
 & \stackrel{\eqref{homdefa}, \eqref{stima-di-b}, \eqref{disuvoluta}} {\lesssim_{s_1}} 
|\rho -\rho'|_{+, s_1} + \d_1^{\frac78} |A_0-A'_0|_{+, s_1} \nonumber \\ 
& \qquad \qquad \qquad \quad 
+ N^{Q+1} \d_1^{\frac{9}{10}} |A'_0-A_0|_{+,s_1} + N^{Q+1} |\rho - \rho'| _{+, s_1}
\nonumber \\
& \stackrel{\eqref{stN1}}{\ll}   \d_1^{\frac45} |A'_0-A_0|_{+, s_1}  +
 \d_1^{- \frac{1}{40}} |\rho-\rho'|_{+, s_1}  \, . \label{ultimissima-2}
\end{align}  
Finally, recalling \eqref{soluzione-approssimata-vera2} and Lemma \ref{pisig}, 
 we obtain  
$$
\big| (T_j a^j -g^j) -(T'_j (a^j)' - (g^j)') \big|_{s_1+\frac12} = |\Pi_{\mathbb G} (z_N - z'_N)|_{s_1+\frac12} 
\lesssim_{s_1}  |z_N - z'_N   |_{s_1+\frac12} 
$$  
and \eqref{ultimissima-2} implies the estimate \eqref{perteq2}.
\end{pf}

\smallskip

\noindent
{\bf Step 3:  Conclusion of the proof of Lemma \ref{reshomlin}.} 
We consider the sets 
\be\label{def:Cantor-homol} 
\Lambda (\e;\cc,A_0) := \bigcap_{j \in \mathbb F} \Lambda^2_j(\e;\cc,A_0) \bigcap \Lambda^1(\e;\cc,A_0)   \, , 
\quad  \frac12 \leq \eta \leq 1 \, , 
\ee
where  the sets $ \Lambda^1(\e;\cc,A_0) $ are defined in \eqref{def:Lambda1} and the sets $ \Lambda^2_j(\e;\cc,A_0) $ in \eqref{def:setsLj}. 
By Lemmata \ref{lemma:measure1} and \ref{lemma:hom1-var}
and 
 \eqref{def:setsLj}, the sets $\Lambda (\e;\cc,A_0) $, $ 1/ 2 \leq \eta \leq 1 $, 
satisfy the properties 1-3 listed in Proposition
\ref{prophomeq} for some $ 0 < \a < 1 / 12 $. 

For all $ \l \in  \Lambda (\e;1,A_0) $ we define  the self-adjoint operator 
\be\label{definition di SN}
{\cal S}(\vphi) := \begin{pmatrix}
d_N (\vphi) & a^*(\vphi) \\
a(\vphi) & 0
\end{pmatrix}  
\ee
where $ d_N (\vphi) = d_N^* (\vphi)$ is defined in Lemma \ref{homdiag} and   
$ a := ( a^j)_{j \in {\mathbb F}} $ by Lemma \ref{lem:homoeq2}. 
The operator  $ {\cal S} $ defined in \eqref{definition di SN} satisfies the estimates \eqref{stima-S}-\eqref{estanys}
by  \eqref{stima-di-d1} 
and \eqref{stima-di-a}-\eqref{stima-di-a1}.
\\[1mm]
{\sc Proof of \eqref{approxlin}.}
The estimate  \eqref{approxlin} follows by \eqref{soluzione:dN}, 
\eqref{special-form}-\eqref{special-form-Lip},  and \eqref{approx-sol-2}, recalling 
the definition of $ T_j $ in \eqref{eqhomai1} and of $ g^j  $ in \eqref{homdefa}.  
\\[1mm]
{\sc Proof of \eqref{31S}-\eqref{tameS-any-s}.}
By \eqref{def:homo-matrix}, \eqref{eq:approx-eq-1}, Lemma \ref{lem:homoeq2},
\eqref{special-form}-\eqref{special-form-Lip}  and the fact that
$$
T_j(a_j)=T_j(\Pi_{\mathbb G} b_j^\sharp)=\Pi_{\mathbb G } T_j^\sharp (b_j^\sharp) \, ,
$$
(see \eqref{homdefa}, \eqref{commTjTtj}), 
we derive 
\begin{align} 
& | J\ppavphi {\cal S} + [J{\cal S}, J {A_0}] |_{\Lip, +,s} \nonumber \\
& \leq  | {\it \Pi}_N \rho_1|_{\Lip,+,s} + \sum_{j \in \mathbb F} 
\big|T_j (a_j)\big|_{\Lip, s+\frac12}  \nonumber
 \\ 
 & \lesssim_s |  \rho_1|_{\Lip,+,s} + \sum_{j \in \mathbb F} \big| {\it \Pi}_N T_j^\sharp (b_j^\sharp ) \big|_{\Lip,s+\frac12}   + 
 \sum_{j \in \mathbb F} \big| {\it \Pi}_N^\bot  T_j^\sharp ( b_j^\sharp) \big|_{\Lip,s+\frac12}  \nonumber \\
 & \lesssim_s   |  \rho_1|_{\Lip, +,s} +  \big| {\it \Pi}_N \rho_2 \big|_{\Lip,s + \frac12} + 
 \sum_{j \in \mathbb F} \big| {\it \Pi}_N^\bot T_j^\sharp (b_j^\sharp) \big|_{\Lip, s+ \frac12} \,. \label{remain1} 
 \end{align}
 Now 
 \begin{align}
 \big| {\it \Pi}_N^\bot T_j^\sharp (b_j^\sharp) \big|_{\Lip, s+ \frac12}  & \stackrel{\eqref{cancellazione-bot}} {\lesssim_s} 
 \big( N+ |  R_0 |_{\Lip , s+\frac12} \big) \|  b_j^\sharp\|_{\Lip , s_1} +
  \big( N+ |  R_0 |_{\Lip , s_1} \big) \|  b_j^\sharp\|_{\Lip , s + \frac12} \nonumber  \\
  &\lesssim_s   \big( N+ \big|  R_0  \big|_{\Lip , +, s} \big) \|  b_j^\sharp\|_{\Lip , s_1} +
  (N+ \e^2) \|  b_j^\sharp\|_{\Lip , s + \frac12}  \nonumber \\
  & \stackrel{\eqref{estimate:bs1}, \eqref{estimate:bs}}{\lesssim_s} 
  \big( N+ |  R_0 |_{\Lip , +, s} \big)  N^Q \d_1^{9/10} \nonumber \\
  & \qquad \qquad +
  N^{Q+1}\big((N^{\varsigma (s-s_1)} +  |  R_0  |_{\Lip , +, s} ) \d_1^{9/10}+ |\rho|_{\Lip,+,s}  \big) \nonumber \\
  &\lesssim_s  N^{Q+1} \big(|  R_0 |_{\Lip , +, s}  +   |\rho|_{\Lip,+,s} \big) + N^{Q+1}\d_1^{9/10}N^{\varsigma (s-s_1) }\nonumber \\
 & \lesssim_s  \d_1^{-1/40} \big( |R_0|_{\Lip,+,s} + |\rho|_{\Lip,+,s} \big) + \d_1^{7/8} \d_1^{- 3 \varsigma \frac{s-s_2}{ s_2 - s_1 }} \, .
  \label{remain2}
 \end{align}
 Estimates \eqref{31S} (for $\e$ small enough) and \eqref{tameS-any-s} are immediate consequences of \eqref{remain1} and \eqref{remain2}.
 \\[1mm]
 {\sc Proof of \eqref{pertS}.}
By \eqref{definition di SN}, Lemma \ref{homdiag}, Lemma \ref{Lemma:rho},
and Lemma \ref{lem:homoeq2},  for any 
$ \l \in \Lambda (\e ; 1, A_0) \cap \Lambda (\e ; 1, A'_0)$, we have 
\begin{align*}
|{\cal S}_{A_0, \rho} -  {\cal S}_{A_0', \rho'} |_{+, s_1} & \lesssim  |d_N -d'_N|_{s_1+\frac12} +
\sum_{j \in \mathbb F} |a_j - a_j'|_{s_1+\frac12} \\
&  \leq \d_1^{\frac45} |A_0-A_0'|_{+, s_1}  +
 \d_1^{-\frac{1}{30}} |\rho-\rho'|_{+, s_1} 
\end{align*}
proving \eqref{pertS}. 
\\[1mm]
{\sc Proof of \eqref{perteq}.}
Also \eqref{perteq}  is a consequence of Lemma \ref{homdiag} and Lemma \ref{lem:homoeq2}, more precisely of 
 \eqref{perteq1}, \eqref{special-form}-\eqref{special-form-Lip} and \eqref{perteq2}.
This concludes the proof of Lemma \ref{reshomlin}.

\section{Splitting step: Proof of Proposition \ref{prophomeq}}\label{sec:proof-splitting}

Recalling   \eqref{def:pro+}, \eqref{new normal form D+}, \eqref{lem:off-dia} we decompose
the coupling term $ \rho \in {\cal L}(H_{\mathbb S}^\bot) $ as
\be\label{deco-rho}
\rho = \Pi_{\mathtt D} \rho + \Pi_{\mathtt O} \rho \, . 
\ee
The operator $ \Pi_{\mathtt O} \rho $ has the form \eqref{decomposition-rho1}, satisfies 
\eqref{propertyM-} and, by 
\eqref{estimate-deco-resto} and \eqref{rho-R0:small-0}, it satisfies also \eqref{prop:rho-plus}.
We then apply Lemma \ref{reshomlin} with $ A_0 $ and $ \rho \rightsquigarrow  \Pi_{\mathtt O} \rho $.  
It provides the  existence of closed subsets $ \Lambda (\e;\eta, A_0) \subset \wtilde{\Lambda} $, $ 1/ 2 \leq \eta \leq 1 $, 
satisfying  the properties 1-3 of Proposition \ref{prophomeq}, 
and a self-adjoint operator $ {\cal S}(\vphi) := {\cal S}(\e,\l)(\vphi) \in {\cal L}
(H_{\mathbb S}^\bot) $ 
of the form \eqref{forma cal-S}, defined 
for all $ \l \in \Lambda(\e; 1, A_0) $, such that 
\be\label{approxlin-applicato}
|J \ppavphi {\cal S}  + [J{\cal S}, J {A_0}] - J \Pi_{\mathtt O} \rho |_{\Lip, +,s_1}  
 \leq \d_1^{\frac74} 
\ee
(see \eqref{approxlin}) and  \eqref{stima-S}-\eqref{estanys}, \eqref{31S}, \eqref{tameS-any-s} hold. 

We now conjugate the Hamiltonian operator 
\be\label{A=A0+rho} 
\ppavphi - J A \, , \quad A =  A_0 + \rho 
\stackrel{\eqref{defA0}} = \frac{D_V}{1+ \e^2 \l} + R_0 +  \rho  \, ,  
\ee
by 
the symplectic linear invertible transformations $ e^{J{\cal S} ( \vphi)} $.
We first notice that the  conjugated operator  
\be\label{transf-op}
 e^{- J {\cal S} (\vphi)} \big( \ppavphi - JA  (\vphi) \big) e^{ J {\cal S} (\vphi)}  = \ppavphi - JA^+ (\vphi) 
\ee
is Hamiltonian, because  
$ e^{J{\cal S} ( \vphi) } $ 
is symplectic and Lemma \ref{lem:sym} implies that 
the (unbounded) operator
$ A^+(\vphi) $ 
is self-adjoint.  
We are going to prove that (see \eqref{new:A+}, \eqref{new-A'})
 $$ 
A^+ (\vphi) =  \frac{D_V}{1+ \e^2 \l} + R_0^+ (\vphi) + \rho^+ (\vphi)  \, ,
$$ 
 where 
\be\label{defR0'}
 R_0^+ :=  R_0 +  \Pi_{\mathtt D} \rho  \, ,
 \ee 
 and 
the estimates \eqref{propo1}-\eqref{propo2}, \eqref{R0+s} of Proposition \ref{prophomeq} hold. First notice that the bounded operator $\rho^+ (\vphi) $ is self-adjoint because 
the operators $A^+ (\vphi) $, $ D_V$ and $ R_0^+ (\vphi) $ are self-adjoint.
\\[1mm]
{\bf Lie series expansion.} 
Consider the $ 1$-parameter family of operators\index{Lie expansion} 
\be\label{pathXt}
X_t := e^{-t J {\cal S} (\vphi)} \big( \ppavphi - JA \big) e^{t J {\cal S} (\vphi)}\, , 
\quad t\in [0,1] \, ,
\ee
connecting $ X_0 :=  \ppavphi - JA $ to 
\be\label{transf-op-new}
X_1 :=  e^{- J {\cal S} (\vphi)} (\ppavphi - JA)e^{J {\cal S} (\vphi)} 
\stackrel{\eqref{transf-op}} = \ppavphi - JA^+  \, . 
\ee 
Since the path  $ t \mapsto X_t $  solves the problem
\be\label{diff:path}
\begin{cases}
\frac{d X_t}{dt} = -J{\cal S} X_t + X_t J {\cal S} = [-J {\cal S}, X_t] = {\rm Ad}_{(-J {\cal S})} (X_t) 
 \cr
 X_0 = \ppavphi - JA \, , 
 \end{cases}
\ee
we deduce  from \eqref{transf-op-new}, \eqref{pathXt}, \eqref{diff:path} the Lie series expansion
\begin{align} \label{Lieseries}
\ppavphi - JA^+  & 
= X_{1} =  \sum_{k \geq 0} \frac{1}{k!} {\rm Ad}^k_{(-J {\cal S})} (X_0) \nonumber \\
& = X_0 + {\rm Ad}_{(-J{\cal S})} (X_0) + \sum_{k \geq 2} \frac{1}{k!} {\rm Ad}^k_{(-J{\cal S})} (X_0) \, . 
\end{align}
Now 
\begin{align}
{\rm Ad}_{(-J{\cal S})} (X_0)& = X_0 J{\cal S} -J {\cal S} X_0 \nonumber \\
&= J {\cal S} \, \ppavphi + \ppavphi (J{\cal S})-JAJ{\cal S} - (J{\cal S} \, 
\ppavphi - J {\cal S} JA) \nonumber \\
&=  \ppavphi (J{\cal S}) + [J{\cal S},JA] \nonumber \\
&  
= J \, \ppavphi {\cal S}  + [J{\cal S},JA_0] + [J{\cal S}, J\rho]  \label{sviluppo-Ad} 
\end{align}
recalling that $ A = A_0 + \rho $, see \eqref{A=A0+rho}.
Comparing \eqref{Lieseries} and \eqref{sviluppo-Ad} we obtain
\begin{align}
-JA^+ & = - JA+ J \ppavphi {\cal S} + [J{\cal S},JA_0] + [J{\cal S}, J\rho] + \sum_{k \geq 2} \frac{1}{k!} {\rm Ad}^k_{(-J{\cal S})} (X_0) \nonumber \\
&  = - JA_0 - J \Pi_{\mathtt D} \rho  - J \rho^+ \label{newA'}
\end{align}
where,  recalling the decomposition \eqref{deco-rho}, 
\be \label{newrho}
- J \rho^+ := \big( J \ppavphi {\cal S} + [J{\cal S},JA_0] - J\Pi_{\mathtt O} \rho \big) +
 [J{\cal S}, J\rho] + \sum_{k \geq 2} \frac{1}{k!} {\rm Ad}^k_{(-J{\cal S})} (X_0) \, . 
\ee
Note that, by \eqref{approxlin-applicato}, the first addendum in  the expression of 
$ - J \rho^+ $ in \eqref{newrho} is very small (in low norm) and the others terms in
\eqref{newrho} are ``quadratic" in $ \rho $ (note that the term $ \rho^+ $ satisfies the 
estimate \eqref{picco-rho+}). By \eqref{newA'} we have that 
\begin{align}\label{new:A+}
A^+ = A_0 + \Pi_{\mathtt D} \rho + \rho^+ 
& \stackrel{ \eqref{defA0}}   = \frac{D_V}{1+ \e^2 \l} +  R_0 + \Pi_{\mathtt D} \rho + \rho^+ \nonumber \\
& \stackrel{\eqref{defR0'}} = \frac{D_V}{1+ \e^2 \l}  + R_0^+   + \rho^+ \, .
\end{align}
\begin{lemma}
$ A^+ $  has the form \eqref{new-A'}-\eqref{form-A'} and \eqref{new-aut-V+}-\eqref{new-aut-V+1} holds. 
\end{lemma}

\begin{pf}
By \eqref{new:A+}, \eqref{form:A0}, \eqref{def:pro+}, \eqref{new normal form D+}, 
\eqref{form:D0},  the operator $ A^+ $ has the form \eqref{new-A'} 
with 
$$ 
\begin{aligned}
& D_0^+ = D_0 + D_+ ( \rho_{\mathbb F}^{\mathbb F} ) = 
{\rm Diag}_{j \in {\mathbb F}} \big(  \mu_j (\e, \l) {\rm Id}_2 
+ \pi_+ [{\widehat \rho}^j_j (0)] \big) 
 \, , \\ 
& V_0^+ = V_0 + \rho_{\mathbb G}^{\mathbb G} \, . 
\end{aligned} 
$$
Since  $ \rho^j_j (\vphi )  $ is a $ 2 \times 2 $ symmetric matrix,  \eqref{proiettore:sim} implies 
\eqref{form-A'}  with 
$$ 
\mu_j^+ (\e, \l) = \mu_j (\e, \l) + \frac12 {\rm Tr}[ {\widehat \rho}^j_j (0)]  \, . 
$$
Hence
$$
 | \mu_j^+ (\e, \l) -  \mu_j (\e, \l) |_{\Lip} +  \| V_0^+ - V_0 \|_{\Lip ,0 } \lesssim \|  \rho \|_{\Lip,0}  
$$
and by \eqref{rho-R0:small-0} the bounds \eqref{new-aut-V+}-\eqref{new-aut-V+1} follow. 
\end{pf}
 
The first  estimate in \eqref{propo1} follows by 
\be \label{estpro1}
|R^+_0 - R_0 |_{\Lip, +,s_1} \stackrel{ \eqref{defR0'}} = 
|\Pi_{\mathtt D} \rho|_{\Lip, +, s_1} \stackrel{\eqref{estimate-deco-resto}, \eqref{rho-R0:small-0}} 
{\lesssim_{s_1}} \d_1 \leq \d_1^{\frac34} 
\ee
for $ \d_1 \leq \e^3 $ small. 
For the estimate of $|\rho^+ |_{\Lip, +, s_1} $ we use the following lemma.

\begin{lemma}  \label{Ads1}
We have
\begin{align}
& |{\rm Ad}_{(-J{\cal S})} (X_0)|_{\Lip, +,s_1} \lesssim_{s_1} \d_1 \label{len1} \\   
& |{\rm Ad}^k_{(-J{\cal S})} (X_0)|_{\Lip, +,s_1} \leq \d_1^{1+\frac34 (k-1)} \, ,  \label{len2} 
\  \forall k \geq 2 \, .
\end{align}
\end{lemma}

\begin{pf}  
By \eqref{sviluppo-Ad} and
\eqref{inter-norma+s}, we have
\begin{align*}
|{\rm Ad}_{(-J{\cal S})} (X_0)|_{\Lip, +,s_1}&= 
\big| J \ppavphi {\cal S} + [J{\cal S},JA_0] + [J{\cal S}, J\rho] \big|_{\Lip, +,s_1} \\
&\leq  \big| J \ppavphi {\cal S} + [J{\cal S},JA_0]-\Pi_{\mathtt O} \rho \big|_{\Lip, +,s_1} \\ 
& \quad +  \big| \Pi_{\mathtt O} \rho \big|_{\Lip, +,s_1} +
C(s_1) |J{\cal S}|_{\Lip, s_1+ \frac12} |\rho |_{\Lip, +,s_1} \\
& 
\stackrel{
\eqref{approxlin-applicato}, \eqref{estimate-deco-resto},  
\eqref{rho-R0:small-0}, \eqref{stima-S}} \leq 
 \d_1^{\frac74} + C(s_1) \d_1 + C (s_1) \d_1^{\frac78+1} \\
 & \lesssim_{s_1} \d_1  
\end{align*}
for  $ \d_1 $ small, proving \eqref{len1}.  
In order to prove \eqref{len2}, 
let 
$ m_k := |{\rm Ad}^k_{(-J{\cal S})} (X_0)|_{\Lip, +,s_1} $.
For $k\geq 1$, 
we get, for $ \d_1$ small, 
$$
\begin{aligned}
m_{k+1}  = \big| [-J{\cal S}, {\rm Ad}^k_{(-J{\cal S})} (X_0)] \big|_{\Lip, +,s_1} 	
& \stackrel{\eqref{inter-norma+s}} \leq C(s_1) |{\cal {\cal S}}|_{\Lip, s_1+ \frac12}\,  m_k \\
&  \stackrel{\eqref{stima-S}} \leq  C(s_1) \d_1^{\frac78} m_k \\
&  \ll \d_1^{\frac34} m_k \, . 
\end{aligned}
$$
This  iterative inequality and \eqref{len1} imply \eqref{len2}. 
\end{pf}

\smallskip

We derive by \eqref{newrho} the bound
\begin{align}
|\rho^+ |_{\Lip, +,s_1} & 
\leq | J \ppavphi {\cal S} + [J{\cal S},JA_0] - J\Pi_{\mathtt O} \rho|_{\Lip, +,s_1}
 + |[J{\cal S}, J\rho]|_{\Lip, +,s_1} \nonumber \\
 & \quad  +  \sum_{k \geq 2} \frac{1}{k!} \big| {\rm Ad}^k_{(-J{\cal S})} (X_0)\big|_{\Lip, +,s_1} \nonumber \\
 & \stackrel{\eqref{approxlin-applicato}, \eqref{inter-norma+s}, \eqref{rho-R0:small-0}, 
 \eqref{stima-S}, \eqref{len2}} 
 \leq \d_1^{\frac74} + C(s_1) \d_1^{\frac78 +1} +   \d_1^{\frac74} 
 \ll  \d_1^{\frac32}  
 \label{picco-rho+}
\end{align}
for $ \d_1 \leq \e^3 $ small.  This proves the second estimate in \eqref{propo1}. 

There remains to estimate the high norms  $ |R^+_0|_{\Lip, +,s} + |\rho^+ |_{\Lip, +,s}$ for $ s \geq s_2 $. 
First, by \eqref{defR0'} and \eqref{estimate-deco-resto}, 
\be\label{estimate-R0S}
|R^+_0|_{\Lip, +,s}=|R_0+ \Pi_{\mathtt D} \rho|_{\Lip, +,s} \leq |R_0|_{\Lip, +,s} + C(s) |\rho|_{\Lip, +,s} 
\ee
which implies the bound for $ |R^+_0|_{\Lip, +,s} $ in \eqref{R0+s} (actually \eqref{estimate-R0S}
 is much better than the estimate \eqref{R0+s} for $ |R^+_0|_{\Lip, +,s} $).  
 For $|\rho^+|_{\Lip, +,s}$ we use the following lemma. 

\begin{lemma} \label{AdS}
For $k\geq 1$,
\be
\begin{aligned}\label{induz-norma-S}
|{\rm Ad}^k_{(-J{\cal S})} (X_0)|_{\Lip, +,s} \leq (C(s))^k 
{\cal M}_s  \d_1^{3 \frac{k-1}{4}}   
\end{aligned}
\ee
where
\be\label{la-prima-S}
{\cal M}_s  :=   \d_1^{- \frac14} \big( |R_0|_{\Lip, +,s}  + |\rho|_{\Lip, +,s} \big) 
 + \d_1^{- \frac34} \d_1^{- 3 \varsigma 	\frac{s-s_2}{ s_2 - s_1 }}   
\, . 
\ee
\end{lemma}
\begin{pf}
By \eqref{sviluppo-Ad} and \eqref{inter-norma+s} we have 
\begin{align}
|{\rm Ad}_{(-J{\cal S})} (X_0)|_{\Lip, +,s}&\leq 
\big| J \ppavphi {\cal S} + [J{\cal S},JA_0] \big|_{\Lip, +,s} + |[J{\cal S}, J\rho] \big|_{\Lip, +,s} \nonumber \\
& \leq  \big| J \ppavphi {\cal S} + [J{\cal S},JA_0] \big|_{\Lip, +,s} + 
C(s) \big| {\cal S} \big|_{\Lip, s_1+ \frac12} |\rho|_{\Lip, +,s}  \nonumber \\
& \quad +C(s)  |{\cal S}|_{\Lip, s+ \frac12} |\rho|_{\Lip, +,s_1}  \nonumber \\
& \stackrel{ \eqref{tameS-any-s}, \eqref{stima-S}, \eqref{estanys}} 
\leq C(s) \Big[ \d_1^{- \frac14} \big( |R_0|_{\Lip, +,s} + |\rho|_{\Lip, +,s} \big) + \d_1^{-\frac34} \d_1^{- 3 \varsigma 	\frac{s-s_2}{ s_2 - s_1 }} \Big] \,. \nonumber 
\end{align} 
Hence, the estimate 
\eqref{induz-norma-S} is proved for $ k = 1 $, recall \eqref{la-prima-S}. 
In order to prove \eqref{induz-norma-S} for $ k \geq 2 $, let 
$ M_k :=|{\rm Ad}^k_{(-J{\cal S})} (X_0)|_{\Lip, +,s} $. 
For $ k \geq 1 $,  we have, 
by \eqref{inter-norma+s}, 
\begin{align}
M_{k+1} 
& =  \big| [-J {\cal S}, {\rm Ad}^k_{(-J{\cal S})} (X_0)] \big|_{\Lip, +,s} \nonumber \\
& \leq C(s)  |{\cal S}|_{\Lip, s_1+ \frac12} M_k + C(s) |{\cal S}|_{\Lip, s+ \frac12} 
|{\rm Ad}^k_{(-J{\cal S})} (X_0)|_{\Lip, +,s_1}  \nonumber \\
& \stackrel{\eqref{stima-S},  \eqref{estanys}, \eqref{la-prima-S}, \eqref{len2}} 
\leq C(s) \big( \d_1^{\frac78} M_k + {\cal M}_s  \d_1^{1+\frac34 (k-1)} \big) \, . 
\label{iteMk}
\end{align}
 Then \eqref{induz-norma-S} follows 
by iteration from  \eqref{iteMk}, 
provided $\varepsilon $ is small enough (independently of $s$).
\end{pf}

Finally, by \eqref{newrho} and \eqref{sviluppo-Ad}  we get 
\be\label{formula-per-resto}
-J \rho^+ = - J\Pi_{\mathtt O} \rho  +  \sum_{k \geq 1} \frac{1}{k!} {\rm Ad}^k_{(-J{\cal S})} (X_0) 
\ee
so that
\begin{align}\label{estimate-rho'S}
|\rho^+|_{\Lip, +,s} & \leq |\Pi_{\mathtt O} \rho|_{\Lip, +,s} + 
\sum_{k \geq 1} \frac{1}{k!} \big| {\rm Ad}^k_{(-J{\cal S})} (X_0) \big|_{\Lip, +,s} \nonumber \\
& \stackrel{\eqref{estimate-deco-resto}, \eqref{induz-norma-S}} 
\leq C(s) |\rho|_{\Lip, +,s}  + {\cal M}_s \d_1^{-3/4} \sum_{k\geq 1} \frac{\big( C(s) \d_1^{3/4}  \big)^k}{k!}  \nonumber \\
& \leq C(s) |\rho|_{\Lip, +,s}  + {\cal M}_s C(s)  \frac{e^{C(s) \d_1^{3/4}}-1}{C(s) \d_1^{3/4}} \nonumber \\
& \stackrel{\eqref{la-prima-S}} {{\lesssim}_s}  \d_1^{-1/4} \big( |R_0|_{\Lip, +,s} + |\rho|_{\Lip, +,s} \big) + \d_1^{-3/4} \d_1^{- 3 \varsigma 	\frac{s-s_2}{ s_2 - s_1 }}  \, . 
\end{align}
The estimate \eqref{R0+s} is a consequence of \eqref{estimate-R0S} and \eqref{estimate-rho'S}. We can prove 
in the same way  \eqref{propo2} taking $\d_1$ small 
enough (depending on $s_2$).  

At last, by \eqref{defR0'} and \eqref{estimate-deco-resto} we get 
$$
|{R'_0}^+ - R_0^+   |_{+, s_1} \leq | R'_0 - R_0  |_{+, s_1} +  C | \rho' - \rho  |_{+, s_1}  
$$ 
for some positive $ C := C(s_1)  $, proving  \eqref{pertR0+}.  

\begin{lemma}
\eqref{pertrho+}  holds. 
\end{lemma}

\begin{pf}
By \eqref{formula-per-resto} we have
\begin{align}
- J (\rho^+ - (\rho^+)') 
& =  - J\Pi_{\mathtt O} (\rho - \rho')   +   {\rm Ad}_{(-J{\cal S})} (X_0) - {\rm Ad}_{(-J{\cal S}')} (X'_0) \nonumber \\
&  \quad + \sum_{k \geq 2} \frac{1}{k!} \big( {\rm Ad}^k_{(-J{\cal S})} (X_0)  - 
{\rm Ad}^k_{(-J{\cal S'})} (X'_0) \big) \, .
\label{prima-dis1} 
\end{align}
Set  
\be\label{def:Vk} 
V_k :={\rm Ad}^k_{(-J{\cal S})} (X_0) - {\rm Ad}^k_{(-J{\cal S}')} (X'_0) \, ,
\qquad 
v_k := |V_k|_{+, s_1} \, , \quad  k \geq 1 \, . 
\ee
By \eqref{perteq} in Lemma \ref{reshomlin} (that we applied with $ \Pi_{\mathtt O} \rho $ instead of $ \rho $), and recalling \eqref{sviluppo-Ad},  we have  
\begin{align}
| V_1 - J \Pi_{\mathtt O}(\rho-\rho')|_{+, s_1}   
& \leq \d_1^{\frac34} |A_0-A'_0|_{+,s_1} + 
\d_1^{- \frac{1}{30}} |\Pi_{\mathtt O} (\rho -\rho') |_{+, s_1}  \nonumber \\
& \quad  + | [J {\cal S}, J \Pi_{\mathtt O} \rho] - [J {\cal S}', J \Pi_{\mathtt O} \rho'] |_{+,s_1} \label{pezdim1} \\
& \lesssim_{s_1} \d_1^{\frac34} |A_0-A'_0|_{+,s_1} +  \d_1^{- \frac{1}{30}} |\rho -\rho'|_{+, s_1}   \label{pezdim2}
\end{align}
because the  term  in \eqref{pezdim1} is bounded as 
\begin{align}
 | [J {\cal S}, J \Pi_{\mathtt O} \rho] - [J {\cal S}', J \Pi_{\mathtt O} \rho'] |_{+,s_1} 
& \stackrel{\eqref{normprod1}} {\lesssim_{s_1}} | {\cal S} - {\cal S}' |_{ s_1+ \frac12 } | \Pi_{\mathtt  O}\rho |_{+,s_1} 
+    | {\cal S}' |_{s_1 + \frac12}  | \Pi_{\mathtt O} (\rho - \rho' )|_{+, s_1  } \nonumber \\
&\stackrel{\eqref{estimate-deco-resto}}{ \lesssim_{s_1} }
| {\cal S} - {\cal S}' |_{+, s_1} | \rho |_{+,s_1} 
+    | {\cal S}' |_{s_1+ \frac12}  | \rho - \rho' |_{+, s_1 } \nonumber \\
& \stackrel{ 
\eqref{pertS}, \eqref{rho-R0:small-0}, \eqref{stima-S}} 
{\lesssim_{s_1}}  	\big(\d_1^{\frac34} |A_0-A'_0|_{+,s_1} + 
\d_1^{- \frac{1}{30}} |\rho -\rho'|_{+, s_1} \big) \d_1 \nonumber \\
& \qquad \qquad  \qquad  +  \d_1^{ \frac78} | \rho - \rho' |_{+, s_1}   \nonumber 
\end{align}
(the estimate \eqref{stima-S} is applied to $ {\cal S}' $). 
By \eqref{pezdim2} and \eqref{estimate-deco-resto} we deduce that 
\be\label{stimav1}
v_1 := | V_1 |_{+, s_1} \leq  C(s_1)  \big(\d_1^{\frac34} |A_0-A'_0|_{+,s_1} + \d_1^{- \frac{1}{30}} |\rho -\rho'|_{+, s_1} \big) \, .
\ee
Now, recalling \eqref{def:Vk}, 
$$
V_{k+1} =  {\rm Ad}_{(-J{\cal S})} (V_k) +
\big[ J{\cal S}' - J {\cal S} , {\rm Ad}^k_{(-J{\cal S}')}(X'_0) \big] \, .
$$
Hence, for $k \geq 1$, we get, setting $ m'_k := |  {\rm Ad}^k_{(-J{\cal S}')}(X'_0) |_{+, s_1} $, 
\begin{align*}
v_{k+1} := | V_{k+1}|_{+, s_1} & 
\lesssim_{s_1}  |{\cal S}|_{+, s_1} v_k +   |{\cal S}-{\cal S}'|_{+, s_1} m'_k \\
& \leq  C(s_1)  \d_1^{\frac78} v_k + C(s_1)  \d_1^{1+ \frac34 (k-1)} \big( \d_1^{\frac34} |A_0-A'_0|_{+,s_1} +
\d_1^{- \frac{1}{30}} |\rho -\rho'|_{+, s_1} \big) 
\end{align*}
by \eqref{stima-S}, \eqref{pertS} and Lemma
\ref{Ads1}. From the previous iterative estimate and \eqref{stimav1}, we may derive
\be\label{terza-disu3}
\sum_{k \geq 2} v_k \leq \d_1^{\frac34} \big( |A_0-A'_0|_{+,s_1} +
|\rho -\rho'|_{+, s_1} \big) \, . 
\ee
Finally, recalling \eqref{prima-dis1} and \eqref{def:Vk}, we get 
$$
\begin{aligned}
|\rho^+ -{\rho'}^+|_{+, s_1} &  \leq \big| V_1- J \Pi_{\mathtt O}(\rho-\rho') \big|_{+, s_1} +  \sum_{k \geq 2} v_k \\
& \stackrel{\eqref{pezdim2}, \eqref{terza-disu3}} \leq 
\d_1^{\frac12} |A_0-A'_0|_{+,s_1} +
\d_1^{- \frac{1}{20}} |\rho -\rho'|_{+, s_1} 
\end{aligned}
$$
which is \eqref{pertrho+}. 
\end{pf}

\chapter{
Construction of  approximate right inverse 
}
\label{sec:proof.Almost-inv}

The goal of this chapter is to prove Proposition \ref{prop-cruciale}. 

\section{Splitting of low-high normal subspaces} \label{s101}

The first step in order to find an approximate solution 
of the equation  
\be\label{sistema-da-riso}
\big( \bar {\om }_\e \cdot \partial_\varphi - J (A_0 + \rho ) \big) h = g \, , 
\ee
see \eqref{sol:almost-inv}, is to apply the splitting Corollary \ref{cor:split}, 
whose  assumptions 
are verified by the hypothesis of  Proposition \ref{prop-cruciale}.
As a consequence 
there are sets $\Lambda_\infty (\e; \eta, A_0, \rho) \subset \tilde \Lambda $, 
$ 1/2 \leq \eta \leq 5/ 6 $, 
and a sequence of symplectic linear transformations  
$ {\cal P}_{\ind} (\vphi) \in {\cal L}( H_{\mathbb S}^\bot ) $, $ \ind \geq 1 $, 
satisfying \eqref{Pins1}-\eqref{Pinsany}, such that, 
for all $ \l \in \Lambda_\infty (\e; 5/6, A_0, \rho) $, the conjugation
 \eqref{goal-trasf-iterata} holds for any $ \ind \geq 1 $, namely 
\be\label{goal-trasf-iterata-ancora}
\big( \bar \om_\e \cdot \partial_\vphi - J  A_0  -J \rho  \big) {\cal P}_\ind (\vphi) = 
{\cal P}_\ind (\vphi) \big( \bar \om_\e \cdot \partial_\vphi - J A_\ind  -J \rho_\ind \big)  \, , 
\ee
where  the operators $ A_\ind $ have the form described in \eqref{defAm}-\eqref{form-A'-m0} and  satisfy
\eqref{der-lambda-m}-\eqref{Rrn-tame0}, in particular $A_\ind \in {\cal C}(2C_1, c_1/2, c_2/2)$ is 
a split admissible operator according to  Definition \ref{def:calC}.  
Moreover  the coupling operator $ \rho_\ind \in {\cal L}( H_{\mathbb S}^\bot )$ in 
\eqref{goal-trasf-iterata-ancora} satisfies the estimates \eqref{form-A'-m}-\eqref{Rrn-tame}, in 
particular it is small in the low norm $ | \ |_{\Lip, +, s_1} $.

Thus, given a function $ g $ satisfying \eqref{g:small-large}, 
in order to 
find a function $ h $  such that  \eqref{sol:almost-inv} holds with a remainder $ r $ satisfying 
\eqref{estimate:r-final}, 
we make the change of variables
\be\label{def:change-var}
g' (\vphi):= {\cal P}_\ind^{-1} (\vphi ) g (\vphi) \in H_{\mathbb S}^\bot \, , \qquad h' (\vphi) := {\cal P}_\ind^{-1} (\vphi ) h(\vphi) \in H_{\mathbb S}^\bot \, , 
\ee
and we look for an approximate solution of the equation
\be\label{eq:vera-n}
\big( \bar \om_\e \cdot \partial_\vphi - J A_\ind  -J \rho_\ind  \big) h' = g'  \, . 
\ee
In Sections \ref{sec:9.1}  and  \ref{sec:9.2}  we shall  solve approximately the more general equation
\be\label{eq:astratta}
\big( \bar \om_\e \cdot \partial_\vphi - J \FA (\e, \l, \vphi) -J \varrho  \big) h' = g'  
\ee
with a split admissible operator $ \FA \in {\cal C}(2C_1, c_1/2, c_2/2)$ as in Definition \ref{def:calC}, i.e. of the form 
\be \label{estFR}
\FA  (\e,\l, \varphi)= \frac{D_V}{1+ \e^2 \l} + \FR (\e,\l, \varphi) \, ,  \quad 
|\FR|_{\Lip, +, s_1} \leq 2C_1 \e^2 \, , 
\ee
and  a self-adjoint operator $ \varrho \in {\cal L} (H_{\mathbb S}^{\bot}) $
satisfying suitable smallness conditions (see  \eqref{smallness:rho}).

Multiplying by $ J $ both sizes of the equation \eqref{eq:astratta}, we look for an approximate solution $ h' $ of
\be\label{def:calLD}
{\cal L} h' = J g'  \qquad  {\rm with} \qquad 
{\cal L} :=  {\cal L}_D + \varrho \, , \quad {\cal L}_D  := J \bar \om_\e \cdot \partial_\vphi  + \FA     \, .
\ee
Notice that the operator $ {\cal L}_D $ 
is block-diagonal according to the splitting  $ H_{\mathbb S}^\bot  = H_{\mathbb F} \oplus H_{{\mathbb G}}  $, 
as the split admissible operator $ \FA \in  {\cal C}(2C_1, c_1/2, c_2/2) $. 

In Section \ref{sec:9.1} we shall first find an approximate right  inverse of $ {\cal L}_D $.
Then in Section   \ref{sec:9.2} we shall obtain, by a Neumann series argument, an approximate right  inverse of the whole 
$ {\cal L} =  {\cal L}_D + \varrho  $ under a suitable smallness condition on $ \varrho $, see \eqref{smallness:rho}. 

Finally in Section \ref{sub:choice:N-ind} we shall apply
the results of Sections \ref{sec:9.1} and \ref{sec:9.2} to obtain an approximate solution of 
the equation \eqref{eq:astratta}
with $\FA = A_\ind$ and $\varrho=\rho_\ind $, namely of equation \eqref{eq:vera-n}, and, ultimately, of \eqref{sistema-da-riso}.  
The number  $ \ind $ of splitting steps will be chosen  in Section \ref{sub:choice:N-ind} large enough (see
\eqref{defNsp-nu}-\eqref{choice:ind}),  so that   
the coupling term $ \rho_\ind $  is 
a small perturbation  of  $ {\cal L}_D $, satisfying the  smallness condition \eqref{smallness:rho}. 

\section{Approximate right inverse of $ {\cal L}_D $}\label{sec:9.1} 

We first obtain an approximate right  inverse of the block diagonal operator $ {\cal L}_D $ introduced in \eqref{def:calLD}.

\begin{proposition}\label{lem:app-inv-LD}
{\bf (Approximate right inverse of $ {\cal L}_D $)}.   
Let ${\FA}  $ be a split admissible operator\index{Split admissible operator} in $ {\cal C} (2C_1,c_1/2 , c_2/2)$ according to Definition \ref{def:calC}.
Then there are closed subsets $ {\mathtt \Lambda} (\e;\cc,\FA ) $, $ 1/2 \leq \eta \leq 1 $, 
satisfying  Properties
1-3 of  Proposition  \ref{propmultiscale},  
and $ \bar N \in \N $ such that, for all $ N \geq \bar N $, 
there exists a linear operator $ {\cal I}_D := {\cal I}_{D,N} $, 
defined for 
$\l \in   {\mathtt \Lambda} (\e;\cc,\FA ) $, with the following properties:
\begin{itemize}
\item setting
\be \label{defQprime} 
Q' := 2(\tau'+ \varsigma s_1 ) +3 \, ,  
\ee
where $ \tau' $ is given by 
Proposition \ref{propmultiscale}, we have 
\begin{align}
   \| {\cal I}_D g' \|_{\Lip, s_0} & \lesssim_{s_0} N^{Q'} \| g' \|_{\Lip, s_0}  \, , \quad 
\| {\cal I}_D g' \|_{\Lip, s_1} \lesssim_{s_1} N^{Q'} \| g' \|_{\Lip, s_1} \label{est:Is1} \, , 
\end{align}
and, $ \forall s \geq s_1 $, 
\begin{align}
 \| {\cal I}_D g' \|_{\Lip, s} & \leq  C N^{Q'}  \|  g' \|_{\Lip, s}  + 
 C(s) N^{Q'}  ( N^{\loss (s-s_1)}+ | \FR |_{\Lip, +,s}) \|  g' \|_{\Lip, s_1}  \label{est:Iss}
\end{align} 
where $ C$ is a constant independent of $ s $ (it depends on $ s_1 $)
and $\FR$ is introduced in \eqref{estFR}.
Furthermore 
\be\label{est:s1-with-loss} 
\| {\cal I}_D g' \|_{\Lip, s_0} \lesssim_{s_1} \| g' \|_{\Lip, s_0+ Q' }   \, . 
\ee
\item 
$ {\cal I}_D$  is an approximate right inverse of $ {\cal L}_D $, in the sense  that
\begin{align}\label{stima:error}
& \| ({\cal L}_D {\cal I}_D - {\rm Id}) g' \|_{\Lip, s_1} \nonumber \\ 
& \qquad  \lesssim_{s_3}   N^{Q'+1-(s_3-s_1)} 
\big( \|g' \|_{\Lip, s_3} +  \big( N^{\varsigma (s_3-s_1)} 
+    |   \FR |_{\Lip,+, s_3} \big)   \| g' \|_{\Lip, s_1} \big) \, . 
\end{align}
\end{itemize}
\end{proposition}

The rest of this section is devoted to the proof of Proposition \ref{lem:app-inv-LD}. 

Recalling Definition \ref{def:calC}, the split admissible\index{Split admissible operator} 
operator   $ \FA   \in {\cal C}(2C_1, c_1/2, c_2/2) $ 
in \eqref{estFR}
is block-diagonal according to the splitting  $ H_{\mathbb S}^\bot  = H_{\mathbb F} \oplus H_{{\mathbb G}}  $, i.e. 
\be \label{proprFA}
\FA ( \e , \l , \vphi )= 
\begin{pmatrix}
{\mathfrak D}(\e , \lambda) & 0  \\
0  & {\mathfrak  W}(\e, \l, \varphi)  
\end{pmatrix}   \, ,  
\ee
and, moreover,  
in the basis of the eigenfunctions $ \{ (\Psi_j,0) $, $(0,\Psi_j)\}_{ j \in {\mathbb F}} $ (see \eqref{Fj-eigenfunctions}),  
the operator $ \FD $ is represented by the diagonal matrix
\be\label{def:Dmj}
\FD :={\rm Diag}_{j \in \mathbb F} \, {\Fm}_j (\e,\l) {\rm Id}_2 \, ,
\ee
where  each $ {\Fm}_j (\e,\l) \in \R $ satisfies (see \eqref{muje2})
$$ 
|{\Fm}_j (\e,\l) - \mu_j | \leq C_1 \e^2 \, , 
$$
and the estimates (see \eqref{Hyp4})
\be \label{proprFA2}
\forall  j \in {\mathbb F} 
 \quad \Big( \frac{c_2}{2} \e^2 \leq {\mathfrak d}_\l {\Fm}_j (\e, \l)  \leq 2 c_2^{-1} \e^2 \quad {\rm or} \quad  
  - 2 c_2^{-1} \e^2 \leq {\mathfrak d}_\l {\Fm}_j (\e, \l)  \leq - \frac{c_2}{2} \e^2 \Big) 
\ee
and, by \eqref{Hyp1}, 
\be \label{proprFA21}
{\mathfrak d}_\l \FW (\e,\l) \leq -\frac{c_1}{2}  \e^2 {\rm Id} \, . 
\ee
By \eqref{proprFA}, 
in order  to find an approximate solution of 
$ ( J \bar \om_\e \cdot \partial_\vphi + \FA   ) h' = g' $, where  $ g' \in L^2 (\T^\es, H_{\mathbb S}^\bot)$, 
we have to solve approximately the  pair of {\it decoupled} equations 
\begin{align}
& \big( J \bar \om_\e \cdot \partial_\vphi + \FD (\e, \l ) \big) h'_{\mathbb F} = g'_{\mathbb F}  \label{le2eq-I} \\
& \big( J \bar \om_\e \cdot \partial_\vphi + \FW (\e, \l , \vphi) \big) h'_{\mathbb G} = g'_{\mathbb G}  \label{le2eq-II}
\end{align}
having used the notation 
\be\label{def:PiG-pr}
\begin{aligned}
h' = 
\begin{pmatrix}
h'_{\mathbb F}     \\
h'_{\mathbb G}   
\end{pmatrix} \,  , &  \ 
g' = \begin{pmatrix}
g'_{\mathbb F}     \\
g'_{\mathbb G}   
\end{pmatrix} \, ,  \\
  h'_{\mathbb F} := \Pi_{\mathbb F} h' \, , \ h'_{\mathbb G} := \Pi_{\mathbb G} h' \, , & 
\ g'_{\mathbb F} := \Pi_{\mathbb F} g' \, , \ g'_{\mathbb G} := \Pi_{\mathbb G} g' \, .
  \end{aligned}
\ee
{\sc Solution of \eqref{le2eq-I}.} Recalling \eqref{proprFA} and \eqref{def:Dmj}, and  decomposing 
\be\label{hpegp}
\begin{aligned} 
& h'_{\mathbb F} = \sum_{j \in {\mathbb F}} h'_j (\vphi) \Psi_j (x) \, , \quad  h'_j (\vphi)  \in \R^2 \, , \\
& g'_{\mathbb F} =  \sum_{j \in {\mathbb F}} g'_j (\vphi) \Psi_j (x) \, , \quad   g'_j (\vphi) \in \R^2 \, , 
\end{aligned}
\ee
the  equation  \eqref{le2eq-I} reduces to the decoupled system of scalar equations
\be\label{inv-1-de}
\big( J \bar \om_\e \cdot \partial_\vphi + \Fm_j (\e, \l ) \big) h'_j (\vphi) = g'_j (\vphi) \, , \quad j \in {\mathbb F} \, .  
\ee
Setting, for $ 1/2 \leq \eta \leq 1 $, 
\begin{align}\label{def:Lambda1-L}
{\mathtt \Lambda}^1 (\e;\cc,\FA) & := \Big\{ \l \in \wtilde{\Lambda} \, : \, 
 | \bar \om_\e \cdot \ell + \Fm_j(\e,\l)| \geq \frac{\gamma_1}{2 \cc \langle \ell \rangle ^{\tau}},  \ 
 \forall (\ell ,  j) \in \Z^\es \times {\mathbb F}   \Big\}  
\end{align}
with  $ \g_1 = \g_0/2 $ and $\tau$ satisfying \eqref{def>tau}, we have the following usual lemma. 

\begin{lemma} \label{homdiag-2} {\bf (Solution of \eqref{le2eq-I}) }
For all $\l  \in {\mathtt \Lambda}^1 (\e;1,\FA) $
  the  equation  \eqref{le2eq-I}  
 has a solution $ h_{\mathbb F}'  = \sum_{j \in {\mathbb F}} h'_j (\vphi) \Psi_j (x) $, written as in \eqref{hpegp}, satisfying
\be \label{easy}
\| h_j' \|_{\Lip, H^s(\T^\es)} \leq C \| g_j' \|_{\Lip, H^{s+ 2 \tau}(\T^\es)} \, , \quad \forall j \in {\mathbb F} \, . 
\ee
Moreover, for all $N \in \N\backslash \{0\}$,
\be \label{easy-N}
\begin{aligned}
\| {\it \Pi}_N h_j' \|_{\Lip, H^s(\T^\es)} & \leq C \| {\it \Pi}_N g_j' \|_{\Lip, H^{s+ 2 \tau}(\T^\es)}  \\
& \leq CN^{2\tau }  \| {\it \Pi}_N g_j' \|_{\Lip, H^{s}(\T^\es)}  \, , \ \  \forall j \in {\mathbb F} \, . 
\end{aligned}
\ee

\end{lemma}

\begin{pf}
By the Fourier expansion 
\begin{align}\label{svil-hj-primog}
& g'_j  (\vphi) = \sum_{\ell \in \Z^\es} \widehat{g'_j} (\ell) e^{\ii \ell \cdot \vphi} \, , \quad
\widehat{g'_j} (\ell) \in \C^2 \,, \quad \ov{ \widehat{g'_j} (\ell)} = \widehat{g'_j} (-\ell) \, , \\
& \label{svil-hj-primoh}
h'_j  (\vphi) = \sum_{\ell \in \Z^\es} \widehat{h'_j}(\ell) e^{\ii \ell \cdot \vphi} \, , \quad
\widehat{h'_j}(\ell) \in \C^2 \,, \quad \ov{ \widehat{h'_j} (\ell)} = \widehat{h'_j} (-\ell) \, , 
\end{align}
each  equation \eqref{inv-1-de} amounts to 
\be\label{ciascuna-amount}
M_{j,\ell} \widehat{h'_j}(\ell) = \widehat{g'_j}(\ell) \, , \ \ell \in \Z^{\es} \, , \quad {\rm where} \quad 
M_{j,\ell} := \ii \bar{\om}_\e \cdot \ell J + \Fm_j {\rm Id}_2 \, , \quad j \in {\mathbb F } \, , 
\ee
are  $2\times 2$ self-adjoint matrices with eigenvalues  $\Fm_j \pm \bar{\om}_\e \cdot \ell $.
As a consequence, for any $\l  $ in the set $  {\mathtt \Lambda}^1 (\e;1,\FA) $ defined in 
\eqref{def:Lambda1-L}, the matrices $ M_{j,\ell} $ are invertible and   
 $ \| M_{j, \ell}^{-1} \| \leq C \la \ell \ra^\tau $, for some positive constant 
 $ C := C(\gamma_0) $,  
for any
$ \ell \in \Z^{|\mathbb S|}$. Hence the equation  \eqref{inv-1-de} has the unique solution
\be\label{h1:stima1}
h'_j(\varphi) =\sum_{\ell \in \Z^\es} M_{j,\ell}^{-1} \widehat{g'_j} (\ell ) e^{\ii \ell \cdot \varphi}
\qquad {\rm satisfying} \qquad 
\| h_j' \|_{ H^s(\T^\es)} \leq C \| g_j' \|_{ H^{s+  \tau}(\T^\es)} \, . 
\ee 
In addition  $ \widehat{h'_j} (\ell ) = M_{j,\ell}^{-1} \widehat{g'_j} (\ell ) $ 
satisfies the reality condition \eqref{svil-hj-primoh} since,  
 taking the complex conjugated equation in \eqref{ciascuna-amount},  we obtain,
 by \eqref{svil-hj-primog}, 
$$ 
M_{j, - \ell} \ov{ \widehat{h'_j} (\ell)} = \widehat{g'_j} (-\ell) 
$$ 
and, therefore,
by uniqueness $\ov{ \widehat{h'_j} (\ell)} = \widehat{h'_j} (-\ell)$. 

Moreover, since
$ 
\| M_{j, \ell}^{-1} \| \leq C \la \ell \ra^\tau$ and  
$ \| M_{j, \ell} \|_{\lip} \simeq | \Fm_j |_{\lip} \leq 2c_2^{-1} \e^2 $ (see \eqref{proprFA2}), 
we get 
\begin{align}
\| h_j' \|_{\lip ,  H^s(\T^\es)} &\lesssim \Big( \sum_{\ell \in \Z^\es} 
\langle \ell \rangle^{2s} \big( \| M_{j,\ell}^{-1} \|_\lip |\widehat{g'_j}(\ell) | + 
 \| M_{j,\ell}^{-1} \| |\widehat{g'_j}(\ell) |_\lip \big)^2 \Big)^{1/2} \nonumber \\
 &\lesssim \Big( \sum_{\ell \in \Z^\es} \langle \ell \rangle^{2s} \big( \| M_{j,\ell}^{-1} \|^2 \| M_{j,\ell} \|_\lip |\widehat{g'_j}(\ell) | + 
\| M_{j,\ell}^{-1} \| |\widehat{g'_j}(\ell) |_\lip \big)^2 \Big)^{1/2} \nonumber \\
&\lesssim \Big( \sum_{\ell \in \Z^\es} \langle \ell \rangle^{2s} \big( \la \ell \ra^{2\t} \e^2 |\widehat{g'_j}(\ell) | + 
\la \ell \ra^{\t} |\widehat{g'_j}(\ell) |_\lip \big)^2 \Big)^{1/2} \nonumber \\
&\lesssim \e^2  \| g_j' \|_{ H^{s+ 2 \tau}(\T^\es)} +  \| g_j' \|_{\lip, H^{s+  \tau}(\T^\es)}  \, . \label{h2:stima2}  
\end{align}
The bounds \eqref{h1:stima1} and \eqref{h2:stima2} imply \eqref{easy}.
The estimate \eqref{easy-N} is obtained in the same way, just considering sums over $| \ell | \leq N$.
\end{pf}

Moreover, 
with arguments similar to those used in Lemmas \ref{lemma:measure1} and  \ref{lemma:hom1-var}, 
using \eqref{proprFA2}, the fact that $ {\Fm}_j (\e,\l) = \mu_j + O(\e^2 )  $, and the unperturbed first order Melnikov
non-resonance conditions \eqref{1Mel},  we deduce the following measure estimate.

\begin{lemma}\label{lem:meas1} {\bf (Measure estimate)} 
Let $ \tau \geq (3/2) \tau_0 +3 + \es $ ($ \t_0 $ is the Diophantine exponent in 
\eqref{1Mel}).
Then the sets $ {\mathtt \Lambda}^1 (\e;\cc, \FA) $ defined in \eqref{def:Lambda1-L} satisfy,
for $\e$ small enough (depending on $C_1, c_1$),
\be\label{meas:1-nu}
|[{\mathtt \Lambda}^1 (\e;\cc, \FA)]^c \cap \wtilde{\Lambda} | \leq   \e \, , \quad \forall 1/2 \leq \eta \leq 1 \, .  
\ee
Moreover, if  $| \FA' - \FA|_{+,s_1} \leq \d \leq \e^3 $ on $\wtilde{\Lambda} \cap \wtilde{\Lambda}' $, 
then, for $ \cc \in [(1/2) + \sqrt{\d}, 1] $,  
\be\label{Cantor-sovra1-sec9}
|\wtilde{\Lambda}' \cap [ {\mathtt \Lambda}^1 (\e;\cc, \FA')]^c \cap {\mathtt \Lambda}^1 (\e;\cc-\sqrt{\d}, \FA)| 
\leq \d^{\frac{1}{12}} \, . 
\ee
\end{lemma}

\begin{remark}
The measure of the set 
$ [\Lambda^1 (\e;\cc, \FA)]^c   $ is 
smaller than $ \e^p $, for any $ p $, taking the exponent 
 $ \tau $ in \eqref{def:Lambda1-L} large enough. This is analogous to the situation described in Remark \ref{rem:meas1s}. 
\end{remark}

\noindent
{\sc Approximate solution of \eqref{le2eq-II}.} 
We now solve approximately the  equation  \eqref{le2eq-II}, that we write as  
\be \label{eqhomai-dir}
T ( h'_{\mathbb G} ) =  g'_{\mathbb G} \, ,   \qquad T  := J \bar \om_\e \cdot \partial_\vphi + \FW (\e, \l , \vphi) \, , 
\ee
where we  regard $ T $ as an unbounded operator of 
 $ L^2(\T^\es, H_{\mathbb G}) $.  
We first extend the operator $ T $ to an unbounded  linear operator $ T^\sharp $  acting on the whole 
space 
$$ 
L^2(\T^\es, H) \, , \quad 
H  = H_{{\mathbb S} \cup {\mathbb F}} \oplus  H_{ {\mathbb G}} = 
 H_{\mathbb S} \oplus H_{\mathbb F} \oplus  H_{ {\mathbb G}} \, , 
 $$ 
 (see \eqref{H:SFG}), by defining
\begin{align}
T^\sharp  & :=  J \ppavphi     +  \frac{D_V}{1+ \e^2 \l} \Pi_{\mathbb S \cup  \mathbb F}  + 
 \frac{\co}{1+ \e^2 \l} \Pi_{\mathbb S}  
+  i_{\mathbb G}  \FW (\e, \l, \vphi )\Pi_{\mathbb G} \nonumber \\
&  \stackrel{\eqref{proprFA}, \eqref{estFR}} =  J \ppavphi   
+  \frac{D_V}{1+ \e^2 \l}   +  \frac{\co}{1+ \e^2 \l} \Pi_{\mathbb S}   +   \FR (\e, \l, \vphi) \Pi_{\mathbb G}   \label{def:T-tilde}
\end{align}
where $ \co > 0 $ is a positive constant that we fix according to \eqref{diof:co},
$ i_{\mathbb G} $ is the canonical injection defined in \eqref{injectionG} and 
$ \FR (\e, \l, \vphi) $ is given  in \eqref{estFR}.

According to the decomposition  $ H = H_{\mathbb S \cup \mathbb F} \oplus H_{{\mathbb G}} $  the operator
$ T^\sharp $ is represented by the matrix of operators
\be\label{wT}
T^\sharp = 
\begin{pmatrix}
 J \ppavphi +   \frac{D_V}{1+ \e^2 \l} +
\frac{\co}{1+ \e^2 \l} \Pi_{\mathbb S}
 & 0    \\
 0 & T       
\end{pmatrix}  \, , 
 \qquad T \Pi_{\mathbb G} = \Pi_{\mathbb G} T^\sharp   \, , 
\ee
and, according to the decomposition  $ H = H_{\mathbb S} \oplus
H_{ \mathbb F}  \oplus H_{{\mathbb G}} $, and recalling \eqref{eqhomai-dir},   
by the matrix of operators
\be\label{three-split0}
T^\sharp = J \ppavphi  + 
\left(
\begin{array}{ccc}
\frac{D_V + \co \, {\rm Id}}{1+ \e^2 \l} &  0 & 0  \\
0  &   \frac{D_V}{1+ \e^2 \l} &  0 \\
0  &  0 &  \FW (\e, \l , \vphi)
\end{array}
\right) \,.
\ee
We look for an approximate solution of
\be\label{aps10}
T^\sharp {h^\sharp}' =  { g^\sharp_{\mathbb G}}'  \qquad {\rm where} \qquad 
 { g^\sharp_{\mathbb G}}' := i_{\mathbb G}  g'_{\mathbb G}   \, . 
\ee
With this aim we apply the multiscale Proposition\index{Multiscale proposition}   \ref{propmultiscale} 
(in case \eqref{phase-space-multi}-($i$) and \eqref{op:i}) to the operator 
\be\label{applicazione-i}
{\cal L}_{r} = (1+ \e^2 \l) T^\sharp \stackrel{\eqref{def:T-tilde}} =
J \om \cdot \partial_\vphi + D_V + \co \Pi_{\mathbb S} + r 
= J \om \cdot \partial_\vphi + X_r 
\ee
where 
$$ 
\omega = (1+ \e^2 \l) \bar \om_\e \, , \quad  r = (1+ \e^2 \l ) \FR \Pi_{\mathbb G} \, ,
\quad  X_r = D_V + \co \Pi_{\mathbb S} + r  \, .  
$$
The operator $ X_r $ belongs to the class
$ {\mathfrak C} (C_1, c_1) $ (see Definition \ref{definition:Xr}-($i$)) for some 
positive constant $ C_1, c_1 > 0 $. 
Indeed 
$ X_r $ is  self-adjoint  with respect to the 
scalar product  $ L^2 (\T^\es, H ) $ (argue as in Lemma \ref{Ttilde-self-ad}) 
and $ r = (1+ \e^2 \l ) \FR \Pi_{\mathbb G} $ has off-diagonal decay by \eqref{estFR}.
More precisely, arguing as in Lemma \ref{lem:off-diag}, we obtain
\be\label{rep+s}
| r |_{\Lip,+, s} \stackrel{\eqref{inter-norma+s}} {\lesssim_{s}}
(1+ \e^2 \l ) |\FR |_{\Lip,+, s} |\Pi_{\mathbb G}|_{\Lip, s + \frac12} 
\stackrel{\eqref{sepr2}}  {\lesssim_{s}}  |\FR |_{\Lip,+, s} \, , 
\ee
and therefore, by \eqref{estFR}, $ | r |_{\Lip,+, s_1} \leq C(s_1) \e^2  $. 
Finally, recalling  \eqref{applicazione-i} and \eqref{three-split0}, 
according to the decomposition  $ H = H_{\mathbb S} \oplus
H_{ \mathbb F}  \oplus H_{{\mathbb G}} $, we represent  
$$
 \frac{X_r}{1+ \e^2 \lambda} = 
\left(
\begin{array}{ccc}
\frac{D_V + \co \, {\rm Id}}{1+ \e^2 \l} &  0 & 0  \\
0  &   \frac{D_V}{1+ \e^2 \l} &  0 \\
0  &  0 &  \FW (\e, \l , \vphi)
\end{array}
\right) 
$$
and 
we deduce by the assumption   \eqref{proprFA21}, 
and arguing  as in Lemma \ref{pos:def-var}, that 
$$ 
 {\mathfrak d}_\lambda  \Big( \frac{X_r}{1+ \e^2 \lambda}  \Big) \leq  - c \, \e^2 \, {\rm Id} 
$$ 
for some constant $ c >  0 $. 
We have proved that   
$X_r $ is in ${\mathfrak C} (C_1, c_1)$ according to Definition \ref{definition:Xr}-($i$), 
for suitable positive constants $ C_1, c_1 $.

Thus  Proposition \ref{propmultiscale} implies the existence of $\e_0>0$ (depending only on fixed positive constants as 
$C_1,c_1, \gamma_1, \tau_1$), of closed subsets  
\be\label{def:setsLj1}
{\mathtt \Lambda}^2 (\e; \cc, \FA) \, , \   1/2 \leq \eta \leq 1 \, , \ {\rm  satisfying \ Properties} \, 
1-3 \ {\rm of \ Proposition \ \ref{propmultiscale}} \, , 
\ee
and $ \bar N \in \N $, 
such that,  for all 
$ \l \in {\mathtt \Lambda}^2 (\e;1, \FA)$, $ \forall N \geq \bar N  $,  
\begin{align}
 &  \forall s \geq s_0 \, ,   \  | (T_N^\sharp)^{-1} |_{\Lip,  s} 
\leq C(s)  N^{Q'} 
 \big(   N^{\loss (s - s_1)} + | \FR |_{\Lip, +,s} \big)   \label{Tn-1Lipss-bis} \\
&  |  (T_N^\sharp)^{-1} |_{\Lip,  s_0}  \leq 
C(s_0) N^{Q'} \, ,  \qquad  |  (T_N^\sharp)^{-1} |_{\Lip,  s_1}  \leq 
C(s_1) N^{Q'} \label{Tn-1Lips1-bis} 
\end{align}
(see \eqref{est:Right+1}, \eqref{Lip-sDM}, \eqref{rep+s}) where  $Q'$ is defined in \eqref{defQprime}  and  $ (T_N^\sharp)^{-1}  $ denotes 
\be\label{da-dove}
\begin{aligned} 
& i) \  {\rm the \  right \ inverse \ of } \ {\it \Pi}_N (T^\sharp)_{|{\mathcal H}_{2N}} \  \hbox{defined in Proposition \ref{propmultiscale}} \ {\rm if}  \ \bar N \leq N < N(\e) \, , \\
& ii) \ {\rm the \ inverse \ of } \ {\it \Pi}_N (T^\sharp)_{|{\mathcal H}_N}  \ {\rm if}  \ N \geq N(\e) \, ,
\end{aligned}
\ee
and $ {\mathcal H}_N $ are the 
finite dimensional subspaces  defined in \eqref{def:EN-tr}.
We lay the stress on the fact that $\bar N$ can be regarded as a fixed constant, being independent of $\e$.
Note also that \eqref{Tn-1Lips1-bis} is a straightforward consequence of \eqref{Tn-1Lipss-bis}, by \eqref{estFR}.

Now, given $ g'_{\mathbb G}  \in L^2 (\T^\es, H_{\mathbb G}) $, 
we  define the following approximate solution of the equation  \eqref{aps10},  
\be \label{homdefa-nu} 
h'_{\mathbb G} := \Pi_{\mathbb G} {h^\sharp}'   \, ,   \qquad
{h^\sharp}' := (T_N^\sharp)^{-1} {{\mathtt g}_N^\sharp}  \, , \qquad  {\mathtt g}_N^\sharp 
:= {\it \Pi}_N  { g^\sharp_{\mathbb G}}' \, , \quad  { g^\sharp_{\mathbb G}}' := i_{\mathbb G}  g'_{\mathbb G} \, , 
\ee
and notice that, by \eqref{da-dove}, 
\be\label{hsN2N}
{h^\sharp}' \in {\cal H}_{2N}  \ {\rm if}   \  N < N(\e) \, ,  \qquad  {\rm and} \qquad  {h^\sharp}' \in {\cal H}_{N}
 \ {\rm   if}  \  N \geq N(\e) \, .
\ee
Finally, 
for all  $ \l $ in 
\be\label{sets:L-inverse}
 {\mathtt \Lambda} (\e;\cc,\FA) :=  {\mathtt \Lambda}^{1} (\e;\cc,\FA) \cap {\mathtt \Lambda}^{2} (\e;\cc,\FA) \, , \quad 
  1/2 \leq \eta \leq 1 \, , 
\ee
where $ {\mathtt \Lambda}^{(1)} $, ${\mathtt \Lambda}^{(2)} $ are the sets introduced in \eqref{def:Lambda1-L}, \eqref{def:setsLj1}, 
 we define the following  approximate right  inverse $ {\cal I}_D $ of 
 the operator $ {\cal L}_D $:  given $ g' \in L^2 (\T^\es, H_{\mathbb S}^\bot)  $, let  
\be\label{la-soluz-inv}
h' := {\cal I}_D g' := {\mathtt \Pi}_N h'_{{\mathbb F}}   + h'_{\mathbb G} \, , 
\qquad  {\mathtt \Pi}_N h'_{{\mathbb F}} :=  \sum_{j \in {\mathbb F} }  ({\it \Pi}_N h'_j (\vphi) ) \Psi_j (x)   \, , 
\ee
where $ h_{\mathbb F}'  $ is the solution of equation  \eqref{le2eq-I} given 
in Lemma \ref{homdiag-2},  $ h'_{\mathbb G} $ is defined in \eqref{homdefa-nu} and  
the projector  $ {\it \Pi}_N $ applies to functions depending only on the variable $ \vphi $ as in 
\eqref{def:PiN-time}.

\begin{lemma}
The operator $ {\cal I}_D $ defined in \eqref{la-soluz-inv}  satisfies \eqref{est:Is1}-\eqref{est:Iss}. 
\end{lemma}

\begin{pf}
We first estimate the function 
$  h'_{\mathbb G} = \Pi_{\mathbb G}{h^\sharp}' 
= \Pi_{\mathbb G} (T_N^\sharp)^{-1} {{\mathtt g}_N^\sharp} $ in \eqref{homdefa-nu}. 
Since $ |\Pi_{\mathbb G}|_{\Lip,s} \leq C(s) $ by \eqref{sepr2},
the estimate 
\eqref{opernorm-Lip+} implies that, 
$ \forall s \geq s_0 $,
\begin{align}
   \| h'_{\mathbb G} \|_{\Lip, s}  
&  = \| \Pi_{\mathbb G}{h^\sharp}' \| _{\Lip, s} 
= \| \Pi_{\mathbb G} (T_N^\sharp)^{-1} {{\mathtt g}_N^\sharp} \|_{\Lip, s}  \label{hiG-0}  \\
& 
 \leq C \| (T_N^\sharp)^{-1} {{\mathtt g}_N^\sharp} \| _{\Lip, s} + 
C(s) \| (T_N^\sharp)^{-1} {{\mathtt g}_N^\sharp} \| _{\Lip, s_1} \nonumber \\
& \stackrel{\eqref{opernorm-Lip+}}  \leq C |(T_N^\sharp)^{-1}|_{\Lip , s_1}     \| {\mathtt g}_N^\sharp\| _{\Lip, s}
+ C(s) |(T_N^\sharp)^{-1}|_{\Lip , s}  \| {\mathtt g}_N^\sharp\| _{\Lip, s_1} \nonumber \\
& \stackrel{\eqref{Tn-1Lipss-bis},\eqref{Tn-1Lips1-bis}}  \leq C N^{Q'}  \|  {{\mathtt g}_N^\sharp} \|_{\Lip, s} +   C(s) N^{Q'}  
(N^{\loss (s - s_1)}+ | \FR |_{\Lip, +,s}) \|  {{\mathtt g}_N^\sharp} \|_{\Lip, s_1} \nonumber \\
& \stackrel{\eqref{homdefa-nu},\eqref{smoothingS1S2-Lip}}  
\leq C N^{Q'}  \|  g'_{\mathbb G} \|_{\Lip, s} +   C(s) N^{Q'}  
( N^{\loss (s - s_1)}+ | \FR |_{\Lip, +,s}) \|  g'_{\mathbb G} \|_{\Lip, s_1}  \label{hiG}
\end{align}
where\index{Smoothing operators} $ C $ is a positive constant which depends on $ s_1 $.  
Moreover \eqref{hiG-0},  \eqref{sepr2}, \eqref{opernorm-Lip}, 
 \eqref{Tn-1Lips1-bis}, \eqref{homdefa-nu},\eqref{smoothingS1S2-Lip} imply 
\begin{align}
& \label{other}
\| h'_{\mathbb G} \|_{\Lip, s_1 } \lesssim_{s_1} \| {h^\sharp}' \| _{\Lip, s_1} 
\lesssim_{s_1} N^{Q'}  \|  g'_{\mathbb G} \|_{\Lip, s_1} \, , \\
& \| h'_{\mathbb G} \|_{\Lip, s_0 } \lesssim_{s_0} \| {h^\sharp}' \| _{\Lip, s_0} 
 \lesssim_{s_0} N^{Q'}  \|  g'_{\mathbb G} \|_{\Lip, s_0} \, .  \label{other-0}
\end{align} 
In addition, Lemma \ref{homdiag-2}   and the fact that $\Psi_j (x) \in C^\infty (\T^d)$ for all $j \in {\mathbb F}$
imply that 
$ {\mathtt \Pi}_N h'_{\mathbb F}  $ defined in \eqref{la-soluz-inv} satisfies, for $s \geq s_0$, 
\begin{align} 
\| {\mathtt \Pi}_N h'_{\mathbb F} \|_{\Lip, s} & \stackrel{\eqref{intbasic+Lip}} \leq 
\max_{j \in {\mathbb F}}  \big( C  \| \Psi_j \|_{s_1}   \| {\it \Pi}_N h'_j\|_{\Lip , H^s(\T^\es)}
+ C(s)  \| \Psi_j \|_s  \| {\it \Pi}_N h'_j\|_{\Lip , H^{s_1} (\T^\es)} \big)  \nonumber \\
& \stackrel{\eqref{easy-N}}{\leq }   N^{2\tau} \big( C  \max_{j \in {\mathbb F}} \|  g'_j\|_{\Lip , H^s(\T^\es)} 
+ C(s)   \max_{j \in {\mathbb F}} \|  g'_j\|_{\Lip , H^{s_1}(\T^\es)}\big)\nonumber \\
& \leq    N^{2\tau} \big( C \| g'_{\mathbb F} \|_{\Lip, s}  + C(s) \| g'_{\mathbb F} \|_{\Lip, s_1} \big) \label{PiNhF}
\end{align}
using that $ g'_j(\vphi) =\la g'_{\mathbb F} (\vphi , \cdot) ,  \Psi_j \ra_{L^2_x }$, see \eqref{hpegp}.  

In conclusion, by \eqref{hiG}, 
\eqref{PiNhF}, the fact that $2\tau < Q'$, \eqref{other}, \eqref{other-0}, \eqref{estFR}, 
the function $ h' := {\cal I}_D g' $ defined in \eqref{la-soluz-inv} satisfies \eqref{est:Iss} and \eqref{est:Is1}. 
\end{pf}

\begin{lemma}\label{loss-regularity-low}
The operator $ {\cal I}_D $ defined in \eqref{la-soluz-inv}  satisfies \eqref{est:s1-with-loss}. 
\end{lemma}

\begin{pf}
Recalling \eqref{la-soluz-inv},  
that $\Psi_j (x) \in C^\infty (\T^d)$ for all $j \in {\mathbb F} $, 
and arguing as in \eqref{PiNhF}, we obtain 
$$ 
\begin{aligned}
\| {\mathtt \Pi}_N h'_{\mathbb F} \|_{\Lip, s_0} 
\lesssim_{s_0}  \max_{j \in {\mathbb F}}   \| h'_{j} \|_{\Lip, H^{s_0}(\T^\es)} 
& \stackrel{ \eqref{easy}} 
{\lesssim_{s_0}}  \max_{j \in {\mathbb F}}   \| g'_{j} \|_{\Lip, H^{s_0+ 2 \tau}(\T^\es)} \\
& \lesssim_{s_1} \| g'_{\mathbb F} \|_{\Lip, s_0 + Q' } 
\end{aligned}
$$
where $ Q' = 2(\tau'+ \varsigma s_1 ) +3  $ is defined in \eqref{defQprime}.
Hence, in order to prove \eqref{est:s1-with-loss},  it is sufficient to show that 
\be\label{wwwa}
\| h'_{\mathbb G} \|_{\Lip, s_0} \lesssim_{s_1} \| g'_{\mathbb G} \|_{\Lip, s_0 + Q' } \, .
\ee
We use a dyadic decomposition argument. First, given an integer $ N > \bar N $, we define a sequence $(M_p)_{0 \leq p \leq q}$
of positive integers, $ q \geq 1 $,  by
$$
M_0 := \bar N   \, , \quad    M_p := 2M_{p-1} \, , \forall p \in [ \! [ 1, q-1 ] \! ] \quad
{\rm and} \quad 2 M_{q-1} \leq M_q:=N < 4 M_{q-1} \, , 
$$
so that 
$$ 
 [ \! [ 0, N ] \! ] = [ \! [ 0, M_0 ] \! ] \cup \ldots \cup [ \! [ M_{q-1}, M_q ] \! ] \, . 
 $$
For $ 0 \leq p \leq q $, we set 
$ {\it \Pi}_p  := {\it \Pi}_{M_p}  $
and we define the dyadic projectors  
$$
\Delta_0 := {\it \Pi}_{0} \, , \quad 
\Delta_p := {\it \Pi}_{p} - {\it \Pi}_{{p-1}}  = {\it \Pi}_{p}  {\it \Pi}_{{p-1}}^\bot \, , \ 
 \forall p \in [ \! [ 1, q ] \! ] \, , \quad 
{\it \Pi}_{{p-1}}^\bot := {\rm Id} - {\it \Pi}_{{p-1}} \, .
$$
For any function $ h \in {\mathcal H}_N $ (recall that $ N = M_q $) 
we consider its dyadic decomposition 
\be\label{dyadic:deco}
h = \sum_{p=0}^q h_p  \qquad 
{\rm where} \qquad h_p :=   \Delta_p h \, , 
\ee
and, for $ 0 \leq p \leq q $, we denote
\be\label{sum-prime-p}
H_p :=  \sum_{\ell = 0}^p h_\ell =  {\it \Pi}_p h \, . 
\ee
In order to estimate $ h'_{\mathbb G} := \Pi_{\mathbb G} {h^\sharp}' $ 
we recall that, by \eqref{homdefa-nu} and \eqref{da-dove},  we have
\be\label{eq:hsN}
{\it \Pi}_N T^\sharp  {h^\sharp}'  =  {{\mathtt g}_N^\sharp}  \, , 
\quad \forall N \geq \bar N \, , 
\ee
and the function $ {h^\sharp}' $ satisfies \eqref{hsN2N}. 
The key estimate is the following:
\begin{itemize}
\item
Let $ g_p^\sharp := \Delta_p {{\mathtt g}_N^\sharp} $. 
Then, each function $ h_p^{\sharp'} = \Delta_p h^{\sharp'}  $,  $  p \in [\![ 0 ,q]\!] $, 
 satisfies 
\be
\begin{aligned} \label{dyadinterm}
\| h_p^{\sharp'} \|_{\Lip,s_0} & \lesssim_{s_1} 
\bar N^{Q'} \big\| g_{p-1}^\sharp + g_p^\sharp + g_{p+1}^\sharp + g_{p+2}^\sharp \big\|_{\Lip, s_0+ Q'} \\ 
& \qquad +  M_p^{-1}   \big( \| {h^\sharp}' \|_{\Lip,s_0} + \| g_{\mathbb G}' \|_{\Lip,s_0} \big) 
\end{aligned}
\ee
with the convention that  $ g_l^\sharp :=0$ for $ l<0$ or $ l>q$. 
\end{itemize}

We split the proof of \eqref{dyadinterm} in different cases. 
\\[2mm]
{\bf Case I:} $   0 \leq p \leq q-3  $ and  $ M_{p+1} \leq N(\e)$. Applying in \eqref{eq:hsN} the projectors  
$ {\it \Pi}_{p+1} $ and $ {\it \Pi}_{p+2} $ for $  0 \leq p \leq q-3  $,  
 and using the splitting  $ {\rm Id} = {\it \Pi}_{p+1}  +  {\it \Pi}_{p+1}^\bot  $, 
we get
\be\label{def:gNsharp}
{\it \Pi}_{p+2}  T^\sharp {\it \Pi}_{p+1} {h^\sharp}' 
+ {\it \Pi}_{p+2}  T^\sharp  {\it \Pi}_{p+1}^\bot   {h^\sharp}' 
 = G_{p+2}^\sharp  \qquad {\rm where} \qquad 
 G_{p+2}^\sharp := {\it \Pi}_{p+2}  {{\mathtt g}_N^\sharp} \, ,
\ee 
and therefore 
\be \label{def:gNsharp1}
{\cal T}_{p+1}^\sharp  H_{p+1}^{\sharp'}  
 = G_{p+2}^\sharp -  {\it \Pi}_{p+2}  T^\sharp  {\it \Pi}_{p+1}^\bot   {h^\sharp}' 
\ee
where 
${\cal T}_{p+1}^\sharp := {\it \Pi}_{p+2} (T^\sharp)_{|{\cal H}_{M_{p+1}}}  $
and 
$ H_{p+1}^{\sharp'} := {\it \Pi}_{p+1} {h^\sharp}'  $. 

We claim that $ {\cal T}^\sharp_{p+1} $ has a left inverse. 
Indeed, by \eqref{def:setsLj1}-\eqref{da-dove} (applied with $ M_{p+1} $ instead of $ N $)
each operator $ {\it \Pi}_{p+1} (T^\sharp)_{|{\cal H}_{M_{p+2}}}$ has a right inverse 
$$ 
R_{p+1} := (T_{M_{p+1}}^\sharp)^{-1} : {\cal H}_{M_{p+1}} \mapsto {\cal H}_{M_{p+2}}  \,.
$$ 
Taking the adjoints in the identity $ \Pi_{p+1} T^\sharp R_{p+1} = {\rm Id}_{{\cal H}_{M_{p+1}}} $, 
we obtain 
$$
R_{p+1}^* {\cal T}^\sharp_{p+1} = {\rm Id}_{{\cal H}_{M_{p+1}}} \qquad {\rm where} \qquad
{\cal T}_{p+1}^\sharp = {\it \Pi}_{p+2} (T^\sharp)_{|{\cal H}_{M_{p+1}}} \, , 
$$  
i.e. $ R_{p+1}^* : {\cal H}_{M_{p+2}}  \to {\cal H}_{M_{p+1}}  $ is a left inverse of  $ {\cal T}^\sharp_{p+1} $, that we denote by 
 $({\cal T}_{p+1}^\sharp)^{-1}  := R_{p+1}^* $. By 
\eqref{A-s-adj} and since 
 $ R_{p+1} = (T_{M_{p+1}}^\sharp)^{-1} $ satisfies 
\eqref{Tn-1Lips1-bis} we deduce  that
\be\label{new38}
\big| ({\cal T}_{p+1}^\sharp)^{-1} \big|_{\Lip, s_1} \leq C(s_1) M_{p+1}^{Q'} \, . 
\ee
Applying the left inverse $ ({\cal T}_{p+1}^\sharp)^{-1} $ 
in \eqref{def:gNsharp1}  
we deduce that $ H_{p+1}^{\sharp'} = {\it \Pi}_{p+1} {h^\sharp}'   $ may be expressed as
$$
 H_{p+1}^{\sharp'}  = 
({\cal T}_{p+1}^{\sharp})^{-1} G_{p+2}^\sharp - 
({\cal T}_{p+1}^\sharp)^{-1} {\it \Pi}_{p+2}  T^\sharp  {\it \Pi}_{p+1}^\bot   {h^\sharp}'  \, . 
$$
Finally, applying the projector $ \Delta_p $ we get,
setting $ g_p^\sharp := \Delta_p {{\mathtt g}_N^\sharp} $, 
\begin{align}
h_p^{\sharp'} := \Delta_p H_{p+1}^{\sharp'} 
&  =  \Delta_p ({\cal T}_{p+1}^\sharp)^{-1} (g_{p-1}^\sharp + g_p^\sharp + g_{p+1}^\sharp  + g_{p+2}^\sharp) \label{dyadic-1}\\
& \quad
  +  \Delta_p ({\cal T}_{p+1}^\sharp)^{-1} G_{p-2}^\sharp \label{dyadic-2} \\
& \quad -   \Delta_p ({\cal T}_{p+1}^\sharp)^{-1} {\it \Pi}_{p+2}  T^\sharp  {\it \Pi}_{p+1}^\bot   {h^\sharp}' \, .  \label{dyadic-3}
\end{align}
For $ p = 0, 1 $ the previous formula holds with  $ g_{-1}^\sharp := 0 $ and 
$ G_{-2}^\sharp := G_{-1}^\sharp := 0 $. In particular for $ p = 0, 1 $ the term  \eqref{dyadic-2}
is not present. 
 
We now estimate separately the terms in \eqref{dyadic-1}-\eqref{dyadic-3}. 
\\[1mm]
{\sc Estimate of \eqref{dyadic-1}. }
By \eqref{opernorm-Lip}, \eqref{new38} and since $ M_{p+1} = 2 M_p $ we have, for 
$ p = 0, \ldots, q -3$,  
\begin{align}
&  \| \Delta_p ({\cal T}_{p+1}^\sharp)^{-1} (g_{p-1}^\sharp + g_p^\sharp + g_{p+1}^\sharp + g_{p+2}^\sharp )\|_{\Lip,s_0} 
\nonumber \\
& \qquad \lesssim_{s_1}
M_p^{Q'} \| g_{p-1}^\sharp + g_p^\sharp + g_{p+1}^\sharp + g_{p+2}^\sharp\|_{\Lip, s_0} \nonumber \\
& \qquad \lesssim_{s_1} \bar N^{Q'} \| g_{p-1}^\sharp + g_p^\sharp + g_{p+1}^\sharp + g_{p+2}^\sharp\|_{\Lip, s_0+ Q'} \, . \label{pezzo-diadico1}
\end{align}
Note that we make appear the multiplicative constant $\bar N^{Q'}$ to deal with the cases $p=0,1$, where $g_0^\sharp$ is in the last sum  (we use $M_0 , M_1 \leq 4 \bar N$). 
Notice that the constant $C(s_1)$ in  \eqref{pezzo-diadico1} does not depend on $\bar N$.
\\[1mm]
{\sc Estimate of \eqref{dyadic-2}. } 
Denoting  
\be\label{def:Bp}
B_p := \big\{ i \in \Z^{\es+d} \, : \, |i| \leq M_p   \big\} \, , 
\ee
we have, if $ 2 \leq p \leq q - 3 $ (for $ p = 0, 1 $ this term is not present)
\begin{align}
\| \Delta_p ({\cal T}_{p+1}^\sharp)^{-1} G_{p-2}^\sharp \|_{\Lip,s_0} &
\stackrel{ \eqref{opernorm-Lip}} {\lesssim_{s_0}}
\big| \Delta_p ({\cal T}_{p+1}^\sharp)^{-1}_{| {\cal H}_{M_{p-2}}}  
\big|_{\Lip,s_0} \| G_{p-2}^\sharp \|_{\Lip,s_0} \nonumber \\
& \stackrel{\eqref{Sm1}, \eqref{def:gNsharp}} {\lesssim_{s_0}} {\mathtt d}(B_{p-2}, B_{p-1}^c)^{-(s_1-s_0)} 
\big| \Delta_p ({\cal T}_{p+1}^\sharp)^{-1}_{| {\cal H}_{M_{p-2}}}  \big|_{\Lip,s_1} \| {{\mathtt g}_N^\sharp} \|_{\Lip,s_0} \nonumber \\
& \stackrel{\eqref{def:Bp}, \eqref{new38}, \eqref{homdefa-nu}} 
{\lesssim_{s_1}} M_p^{-(s_1-s_0)} M_p^{Q' }  \| g_{\mathbb G}' \|_{\Lip,s_0}  \, , 
\label{pezzo-diadico2}
\end{align}
keeping in mind that $M_{l+1}=2 M_l$ for $0 \leq l \leq p$.  
\\[2mm]
{\sc Estimate of \eqref{dyadic-3}. } By \eqref{opernorm-Lip}  we have
\begin{align}
& \| \Delta_p ({\cal T}_{p+1}^\sharp)^{-1} {\it \Pi}_{p+2}  T^\sharp  {\it \Pi}_{p+1}^\bot   {h^\sharp}' \|_{\Lip,s_0} \nonumber \\
& \quad  \lesssim_{s_0}
\big| \Delta_p  ({\cal T}_{p+1}^\sharp)^{-1} {\it \Pi}_{p+2}  T^\sharp  {\it \Pi}_{p+1}^\bot  \big|_{\Lip,s_0} \| {h^\sharp}' \|_{\Lip,s_0} 
\nonumber \\
&\quad  \stackrel{\eqref{Sm1} } {\lesssim_{s_0}}  {\mathtt d}(B_{p}, B_{p+1}^c)^{-(s_1-s_0)} 
\big| ({\cal T}_{p+1}^\sharp)^{-1} {\it \Pi}_{p+2}  T^\sharp  {\it \Pi}_{p+1}^\bot  \big|_{\Lip,s_1} \| {h^\sharp}'  \|_{\Lip,s_0} \nonumber \\
& \quad  \stackrel{\eqref{def:Bp}, \eqref{new38}} {\lesssim_{s_1}}  M_p^{-(s_1-s_0)} M_p^{Q'+1 }  \| {h^\sharp}'  \|_{\Lip,s_0}  \, ,
\label{pezzo-diadico3}
\end{align}
using that 
$ | {\it \Pi}_{p+2}  T^\sharp |_{\Lip,s_1} \lesssim_{s_1} M_{p+2} \lesssim_{s_1} M_p$. 
Now, since $ s_1 $ is large according to \eqref{s1},  
and  $ Q' = 2(\tau'+ \varsigma s_1 ) + 3 $,
 we have that   $  Q' - (s_1 - s_0 )  < -1 $ 
 and by \eqref{dyadic-1}-\eqref{dyadic-3} and
\eqref{pezzo-diadico1}, \eqref{pezzo-diadico2}, \eqref{pezzo-diadico3} we obtain the estimate
\eqref{dyadinterm}, for any 
$0 \leq p \leq q-3 $ and $M_{p+1} \leq N(\e)$. 
\\[1mm]
{\bf Case II:  $ 0\leq p \leq q-3$} and  $ M_{p+1} > N(\e)$.  We have just to replace \eqref{def:gNsharp} with
$$
{\it \Pi}_{p+1}  T^\sharp {\it \Pi}_{p+1} {h^\sharp}' 
+ {\it \Pi}_{p+1}  T^\sharp  {\it \Pi}_{p+1}^\bot   {h^\sharp}' 
 = G_{p+1}^\sharp 
$$
and apply the inverse $({\cal T}_{p+1}^\sharp)^{-1}$ of ${\cal T}_{p+1}^\sharp :=T_{M_{p+1}}^\sharp$ (recall that by \eqref{da-dove} this operator admits an inverse), which 
satisfies \eqref{Tn-1Lips1-bis}. 
Then $ H_{p+1}^{\sharp'} := {\it \Pi}_{p+1} {h^\sharp}'  $ satisfies 
$$
H_{p+1}^{\sharp'}  = 
({\cal T}_{p+1}^{\sharp})^{-1} G_{p+1}^\sharp - 
({\cal T}_{p+1}^\sharp)^{-1} {\it \Pi}_{p+1}  T^\sharp  {\it \Pi}_{p+1}^\bot   {h^\sharp}'   
$$
and, applying $\Delta_p$, we derive the estimate \eqref{dyadinterm} in the same way (notice that 
since the functions $ g_{p-1}^\sharp, \ldots, g_{p+2}^{\sharp} $ are orthogonal
for the scalar product associated to the $s_0$-norm, 
we have $ \| g_{p-1}^\sharp + \ldots + g_{p+1}^{\sharp}\|_{\Lip, s_0} \leq 
\| g_{p-1}^\sharp + \ldots + g_{p+1}^{\sharp} + g_{p+2}^{\sharp}\|_{\Lip, s_0}  $). 
\\[2mm]
{\bf Case III: $ q-2 \leq p \leq q$}. By \eqref{homdefa-nu} we write 
$$
\begin{aligned} 
h_p^{\sharp'} & = \Delta_p  h^{\sharp'}  
= \Delta_p (T_N^\sharp )^{-1} {\mathtt g}_N^\sharp \\ 
& =\Delta_p (T_N^\sharp )^{-1} (g_{p-1}^\sharp + \ldots +
g_{q}^\sharp) + 
\Delta_p (T_N^\sharp )^{-1} {\it \Pi}_{p-2}{\mathtt g}_N^\sharp \, .
\end{aligned}
$$
Now recalling \eqref{def:Bp}, and  since $p-2 \geq q-4$,  $ M_q = N $, we have that 
$$
{\mathtt d}(B_{p-1}^c , B_{p-2}) \geq M_{p-2} \geq N/32 \, . 
$$ 
Hence, arguing as above,   we get by \eqref{Tn-1Lips1-bis}, \eqref{Sm1}, 
\eqref{opernorm-Lip},  
\begin{align*}
\| h_p^{\sharp'} \|_{\Lip , s_0}& \lesssim_{s_1} 
 N^{Q'} ( \| g_{p-1}^\sharp \|_{\Lip, s_0} + \ldots + \| g_{q}^\sharp \|_{\Lip, s_0})
+ N^{-(s_1-s_0)} | (T_N^\sharp )^{-1}|_{\Lip, s_1} \| {\mathtt g}_N^\sharp \|_{\Lip, s_0} \\
&  \lesssim_{s_1}   \big( \| g_{p-1}^\sharp \|_{\Lip, s_0+Q'} + \ldots + \| g_{q}^\sharp \|_{\Lip, s_0+Q'} \big)
+ N^{-1} \| g_{\mathbb G}' \|_{\Lip, s_0} \, .
\end{align*}
In conclusion, since $N^{-1} \leq M_p^{-1} $, the estimate 
\eqref{dyadinterm} is proved for any $p \in [\![ 0 ,q]\!]$. 

\medskip

By \eqref{dyadinterm} we have 
\begin{align}
\| h_p^{\sharp'} \|_{\Lip,s_0}^2 & \lesssim_{s_1} 
\bar N^{2Q'} \big( \| g_{p-1}^\sharp \|_{\Lip, s_0+ Q'}^2 + \| g_p^\sharp \|_{\Lip, s_0+ Q'}^2 + 
 \| g_{p+1}^\sharp \|_{\Lip, s_0+ Q'}^2 + \| g_{p+2}^\sharp \|_{\Lip, s_0+ Q'}^2  \big) 
 \nonumber \\ 
 & \quad \ +  M_p^{-2}    \| {h^\sharp}' \|_{\Lip,s_0}^2  +  M_p^{-2}    \|  g_{\mathbb G}' \|_{\Lip,s_0}^2 \, . \label{ultimo:sq}
\end{align}
Taking the sum for $ 0 \leq p \leq q $ of  \eqref{ultimo:sq}, 
recalling that $ g_p^\sharp := \Delta_p {{\mathtt g}_N^\sharp} $, 
the definition of $ {{\mathtt g}_N^\sharp}  $ in 	\eqref{homdefa-nu}, 
and that $ M_0 = \bar N $, we obtain
$$
\| {h^\sharp}' \|_{\Lip,s_0}^2 \lesssim_{s_1}  
\bar N^{2Q'} \| g_{\mathbb G}' \|_{\Lip, s_0+ Q'}^2  
+   \bar N^{-2}    \| g_{\mathbb G}' \|_{\Lip,s_0}^2  
+  \bar N^{-2}    \| {h^\sharp}' \|_{\Lip,s_0}^2 \,.
$$
For $ \bar N^{-1}    \leq \delta (s_1) $  small enough, this implies 
$$
 \| {h^\sharp}' \|_{\Lip,s_0} \lesssim_{s_1} \bar N^{Q'} \| g_{\mathbb G}' \|_{\Lip, s_0+ Q'}  \, . 
 $$
Then the function $ h'_{\mathbb G} := \Pi_{\mathbb G} {h^\sharp}' $
satisfies, by Lemma \ref{pisig},  the estimate  
$$ 
\| h'_{\mathbb G}  \|_{\Lip,s_0} 
\lesssim_{s_1} \bar N^{Q'} \| g_{\mathbb G}' \|_{\Lip, s_0+ Q'}  \, .
$$
Since $\bar N $ is  a fixed constant, depending on $ s_1 $,  
this inequality implies \eqref{wwwa}. The proof of the lemma is complete. 
\end{pf}

We now prove that $ {\cal I}_D $ is an approximate  right  inverse of $ {\cal L}_D $ satisfying \eqref{stima:error}. 

\begin{lemma}
\eqref{stima:error} holds. 
\end{lemma}

\begin{pf}
We have to estimate
\be\label{forma:error}
({\cal L}_D {\cal I}_D  - {\rm Id})g' \stackrel{\eqref{la-soluz-inv}} = {\cal L}_D h'  - g' 
\stackrel{\eqref{def:calLD}} = (J \bar \om_\e \cdot \partial_\vphi + \FA) h' - g' \, .
\ee
Recalling the definition of $ h' := {\cal I}_D g'  $ in  
\eqref{la-soluz-inv},  \eqref{le2eq-I}-\eqref{le2eq-II}, 
Lemma \ref{homdiag-2}, \eqref{eqhomai-dir}
  \eqref{wT},  \eqref{da-dove}, \eqref{homdefa-nu},   we have 
\be\label{ancora-forma-error}
 (J \bar \om_\e \cdot \partial_\vphi + \FA ) h' - g' = 
 \begin{pmatrix}
  - {\mathtt \Pi}_N^\bot g_{\mathbb F}'   \\
 \Pi_{\mathbb G}  {\it \Pi}_N^\bot \big({T^\sharp} \, {h^\sharp}'  - {g^\sharp_{\mathbb G}}' \big)   
\end{pmatrix} \, . 
\ee
Hence by \eqref{forma:error}  and \eqref{ancora-forma-error} we have 
\begin{align}
\big\| \big( J \bar \om_\e \cdot \partial_\vphi + \FA \big) h'  -  g'  
\big\|_{\Lip, s_1} & \lesssim_{s_1} 
\|   {\mathtt \Pi}_N^\bot g_{\mathbb F}'  \|_{\Lip, s_1} + \|{\it \Pi}_N^\bot ( {T^\sharp} \, { h^\sharp}' -  {g^\sharp_{\mathbb G}}' ) \|_{\Lip, s_1} \label{est:err1}	\\
& \lesssim_{s_1}
\|   {\mathtt \Pi}_N^\bot g_{\mathbb F}'  \|_{\Lip, s_1} + \| {\it \Pi}_N^\bot {T^\sharp} \, {h^\sharp}' \|_{\Lip, s_1} + \|  
{\it \Pi}_N^\bot {g^\sharp_{\mathbb G}}'  \|_{\Lip, s_1}  
\nonumber \\
& \lesssim_{s_1}  N^{- (s_3- s_1)} \|  g'  \|_{\Lip, s_3} + N^{- (s_3- s_1)} \|  {\it \Pi}_N^\bot {T^\sharp} \, {h^\sharp}' \|_{\Lip, s_3} 
\nonumber \end{align}
by the smoothing property\index{Smoothing operators} \eqref{smoothingS1S2-Lip}. 

Recalling the form of $ T^\sharp $ in \eqref{def:T-tilde}, and writing 
$ D_V =  D_m +  (D_V - D_m ) $, we obtain  
\begin{align}
{\it \Pi}_N^\bot T^\sharp \, { h^\sharp}' & = {\it \Pi}_N^\bot T^\sharp \, {\it \Pi}_{2N} {h^\sharp}' \nonumber \\
&  = 
 {\it \Pi}_N^\bot J \ppavphi  {\it \Pi}_{2N} {h^\sharp}' +
\frac{1}{1+ \e^2\l } 
\big( {\it \Pi}_N^\bot  D_m  {\it \Pi}_{2N}  {h^\sharp}' 
+ {\it \Pi}_N^\bot ( D_V- D_m )  {h^\sharp}' \big) 
\nonumber \\
& \quad 
+  \frac{1}{1+ \e^2\l } {\it \Pi}_N^\bot \co \Pi_{\mathbb S}  {\it \Pi}_{2N} {h^\sharp}'  
+ {\it \Pi}_N^\bot  \FR \Pi_{\mathbb G}  {h^\sharp}'   \, . \label{altri-res}
\end{align}
Notice that, by \eqref{hsN2N}, if $ N \geq N(\e) $ then
$ {h^\sharp}'  \in {\cal H}_N $ and the first and second terms in
\eqref{altri-res}  are zero.

In conclusion, by \eqref{est:err1}, \eqref{altri-res}, using \eqref{opernorm-Lip},
 \eqref{sepr2}, we get  
\begin{align}
&  \|  (J  \bar \om_\e \cdot \partial_\vphi + \FA  ) h'  -   g'   \|_{\Lip, s_1} \label{ais3last}\\
&
 \lesssim_{s_3}  N^{- (s_3 - s_1) } \big( \| g' \|_{\Lip, s_3} +  
\|   ( D_V- D_m )  {h^\sharp}' \|_{\Lip, s_3} +
\|   \FR \Pi_{\mathbb G}  {h^\sharp}' \|_{\Lip, s_3}  + N \|{h^\sharp} \|_{\Lip,s_3} \big)  
\nonumber \\
& \stackrel{\eqref{DeltaV2}, \eqref{estFR}} {\lesssim_{s_3}} 
N^{- (s_3 - s_1) } \big( \| g' \|_{\Lip, s_3} +  
N \|   { h^\sharp}' \|_{\Lip, s_3} +
|   \FR |_{\Lip, s_3} \|   {h^\sharp}' \|_{\Lip, s_1}  \big)  \nonumber \\
& \stackrel{\eqref{hiG-0}, \eqref{hiG}, \eqref{other}} {\lesssim_{s_3}}
N^{- (s_3 - s_1) }  N^{Q'+1} \big(\| g' \|_{\Lip, s_3} +   
 ( N^{\varsigma (s_3-s_1)} +      | \FR |_{\Lip, +, s_3} ) \|   g' \|_{\Lip, s_1})  \big) 
 \nonumber 
\end{align}
(and recall that $ g_{\mathbb G}' = \Pi_{\mathbb G} g' $ by \eqref{def:PiG-pr}).
Recalling \eqref{forma:error}, the estimate \eqref{ais3last}  proves \eqref{stima:error}. 
\end{pf}

Now by Lemma \ref{lem:meas1} and \eqref{def:setsLj1},  ${\mathtt \Lambda}^1 (\e;\cc,\FA) $ and 
$ {\mathtt \Lambda}^2 (\e;\cc,\FA)$ satisfy  properties 1-3 of 
Proposition \ref{propmultiscale}. Since these properties are preserved
under finite intersection, this implies the following lemma. 

\begin{lemma} \label{measureL}
The sets ${\mathtt \Lambda} (\e;\cc,\FA)$ defined in \eqref{sets:L-inverse} satisfy properties 1-3 of 
Proposition \ref{propmultiscale}. 
\end{lemma}

\section{Approximate right  inverse of $ {\cal L} = {\cal L}_D + \varrho $}\label{sec:9.2}

In this section we construct  an approximate right inverse 
of $ {\cal L} = {\cal L}_D + \varrho $ by a perturbative 
Neumann series argument, for $ \varrho \in {\cal L} ( H_{\mathbb S}^\bot) $ small,  using the approximate right inverse $ {\cal I}_D $ of $ {\cal L}_D $ found in  Proposition \ref{lem:app-inv-LD}.
 
We denote 
\be\label{Rrn-tame-nu-Un}
 U (s) :=   U_{\FA, \varrho } (s) := | \FR  |_{\Lip, +, s} + | \varrho  |_{\Lip, +, s}  \, . 
\ee

\begin{proposition} \label{lemma:ap-inv_L} {\bf (Approximate right inverse of $ {\cal L} $)}.   
There is $c_0>0$ (depending on $s_1 $)  such that, for  $ N \geq \bar{N}$, if  
$ \varrho $ satisfies 
\be\label{smallness:rho}
 | \varrho |_{\Lip, s_1} N^{Q'} = | \varrho |_{\Lip, s_1} N^{2(\t' + \varsigma s_1) + 3} \leq c_0 \, , 
\ee
then 
$ \forall \l \in {\mathtt \Lambda} (\e;\cc,\FA) $, $  1/2 \leq \eta \leq 1 $ (where  the set 
$ {\mathtt \Lambda} (\e;\cc,\FA)$ is defined in Proposition \ref {lem:app-inv-LD}), 
the operator $ {\rm Id}  + {\cal I}_D \varrho$ is invertible and 
\be\label{def:inv-ap-I}
{\cal I} := {\cal I}_N :=  ({\rm Id}  + {\cal I}_D \varrho )^{-1} {\cal I}_D  
\ee
is an approximate  right inverse of $ {\cal L} $, in the sense that
\be \label{error:ap-I}
\begin{aligned}
& \| ({\cal L} \, {\cal I} - {\rm Id}) g'  \|_{\Lip, s_1}  \lesssim_{s_3} \\
& N^{ Q'+1 - (s_3-s_1)} \big(  \| g' \|_{\Lip, s_3}  + 
 (N^{ \varsigma(s_3- s_1)}+ N^{Q'} U (s_3) ) \| g' \|_{\Lip, s_1} \big) \, . 
\end{aligned}
\ee
Moreover the operator $ {\cal I}$ satisfies
\begin{align}
& \| {\cal I} g' \|_{\Lip, s_1} \lesssim_{s_1} N^{Q'} \| g' \|_{\Lip, s_1} \, , \label{est:Is1-tot} \\
& \| {\cal I} g' \|_{\Lip, s} \leq C N^{Q'}   \|  g' \|_{\Lip, s}  +  C(s)  N^{Q'} 
\big( N^{\varsigma (s-s_1) } + 
N^{Q'} U (s) \big)  \| g' \|_{\Lip, s_1} 
\, , \quad \forall s \geq s_1 \, , \label{est:Iss-tot}
\end{align} 
where the constant $ C $ is independent of $ s $ (it depends on $ s_1 $). 
Furthermore
\be\label{new:cal-loss}
 \| {\cal I} g' \|_{\Lip, s_0} \lesssim_{s_1} \| g' \|_{\Lip, s_0+ Q'} \, .
\ee
\end{proposition}

The proof of Proposition \ref{lemma:ap-inv_L} is given in the rest of this section.  

We first justify that the operator
${\rm Id}  + {\cal I}_D \varrho$ is invertible and provide appropriate estimates for its inverse, as an application of Lemma
\ref{pert+}. By  \eqref{est:Is1}  and \eqref{opernorm-Lip} we have 
\be \label{IDrhos1}
\begin{aligned}
& \| {\cal I}_D \varrho h'  \|_{\Lip , s_0 } \lesssim_{s_0} N^{Q'} \| \varrho h' \|_{\Lip , s_0} \lesssim_{s_0} 
N^{Q'} |\varrho |_{\Lip , s_0 } \| h' \|_{\Lip , s_0 } \, , \\
& \| {\cal I}_D \varrho h'  \|_{\Lip , s_1 } \lesssim_{s_1} N^{Q'} \| \varrho h' \|_{\Lip , s_1 } \lesssim_{s_1} 
N^{Q'} |\varrho |_{\Lip , s_1 } \| h' \|_{\Lip , s_1 } \, , 
\end{aligned}
\ee
and, by  \eqref{est:Iss}, \eqref{Rrn-tame-nu-Un}, 
\begin{align}  
\| {\cal I}_D \varrho h'  \|_{\Lip , s} & \leq  CN^{Q'} \| \varrho h' \|_{\Lip , s} + C(s) N^{Q'} 
 \big( N^{\varsigma (s-s_1)} + U (s) \big) \| \varrho h' \|_{\Lip , s_1 } \nonumber \\
 & \stackrel{\eqref{opernorm-Lip+}} 
 \leq  CN^{Q'}  |\varrho |_{\Lip , s_1 }  \| h' \|_{\Lip , s} + C(s) N^{Q'}  |\varrho |_{\Lip , s }  \| h' \|_{\Lip , s_1} \nonumber \\
 &  \qquad + C(s) N^{Q'} |\varrho |_{\Lip , s_1 }  \big( N^{\varsigma (s-s_1)} + 
 U (s) \big)   \| h' \|_{\Lip , s_1}  \label{IDerhos}
\end{align}
where $ C = C(s_1) $.
Hence, there is a positive constant $c_0$ (depending on $ s_1 $) such that, if
\be \label{cond:rhos1}
 | \varrho |_{\Lip, s_1} N^{Q'} \leq c_0 \, , 
\ee
then the operator ${\cal I}_D \varrho $ satisfies, by \eqref{IDrhos1}, \eqref{IDerhos}  and
recalling also \eqref{Rrn-tame-nu-Un}, 
\begin{align}
& \| {\cal I}_D \varrho h'  \|_{\Lip , s_0 } \leq \frac12 \|h'\|_{\Lip , s_0 } \label{IDrhos0+} \\
& \| {\cal I}_D \varrho h'  \|_{\Lip , s_1 } \leq \frac12 \|h'\|_{\Lip , s_1 } \label{IDrhos1+} \\
& \| {\cal I}_D \varrho h'  \|_{\Lip , s } \leq \frac12 \|h'\|_{\Lip , s} + C(s)  \big( N^{\varsigma (s-s_1)} + N^{Q'} U(s) \big) \|h'\|_{\Lip , s_1 }  \, . 
\label{IDrhos+}
\end{align}
By \eqref{IDrhos0+}-\eqref{IDrhos+}  and Lemma \ref{pert+} (applied
with $ R  = {\cal I}_D \varrho $ and $ E = L^2 (\T^\es, H_{\mathbb S}^\bot)$), 
the operator $ {\rm Id}  + {\cal I}_D \varrho $ is invertible and 
its inverse satisfies  the tame estimates 
\begin{align}
&  \| ({\rm Id}  + {\cal I}_D \varrho)^{-1} g' \|_{\Lip, s_0} \leq 2 \|  g' \|_{\Lip, s_0}  \label{inv-t-s1-s0}\\ 
& \| ({\rm Id}  + {\cal I}_D \varrho)^{-1} g' \|_{\Lip, s_1} \leq 2 \|  g' \|_{\Lip, s_1}  \label{inv-t-s1}\\ 
& \| ({\rm Id}  + {\cal I}_D \varrho )^{-1} g' \|_{\Lip, s} \leq
2 \| g' \|_{\Lip, s} + C(s) \big( N^{\varsigma (s-s_1)} + N^{ Q'} U (s) \big) \| g' \|_{\Lip, s_1}  
\label{inv-t-ss}
\end{align}
for all $ s \geq s_1 $. 
We now estimate the operator $ {\cal I} $ defined in \eqref{def:inv-ap-I}. 

\begin{lemma}
The operator $ {\cal I} =  ({\rm Id}  + {\cal I}_D \varrho )^{-1} {\cal I}_D   $ 
defined in \eqref{def:inv-ap-I} satisfies
\eqref{est:Is1-tot}-\eqref{est:Iss-tot} and \eqref{new:cal-loss}. 
\end{lemma}

\begin{pf}
 By   \eqref{inv-t-s1} and \eqref{est:Is1} we have    
$$
\| {\cal I} g' \|_{\Lip, s_1} \leq 2  \| {\cal I}_D g' \|_{\Lip, s_1} \lesssim_{s_1} N^{Q'} \| g' \|_{\Lip, s_1} 
$$
which is \eqref{est:Is1-tot}. In addition \eqref{inv-t-ss}, \eqref{est:Is1}, \eqref{est:Iss},
\eqref{Rrn-tame-nu-Un}, imply 
\begin{align*}
\| {\cal I} g' \|_{\Lip, s} & 
\leq 2 \| {\cal I}_D g' \|_{\Lip, s} +  C(s)
\big( N^{\varsigma (s-s_1)} + N^{ Q' } U (s) \big) \| {\cal I}_D g' \|_{\Lip, s_1}  \\
&  \leq C N^{Q'}  \|  g' \|_{\Lip, s}  + 
C(s)  N^{Q'}  \big( N^{\varsigma (s-s_1) } + 
N^{ Q'} U (s) \big)  \| g' \|_{\Lip, s_1} 
\end{align*}
proving \eqref{est:Iss-tot}. Finally, by \eqref{inv-t-s1-s0}, 
$$
\begin{aligned}
\| {\cal I} g' \|_{\Lip, s_0} = \| ({\rm Id}  + {\cal I}_D \varrho)^{-1}  {\cal I}_{D} g' \|_{\Lip, s_0}
&  \leq  2 \|  {\cal I}_{D} g' \|_{\Lip, s_0} \\
& \stackrel{\eqref{est:s1-with-loss}} {\lesssim_{s_1}}  \| g' \|_{\Lip, s_0+ Q'} 
\end{aligned}
$$
which is \eqref{new:cal-loss}.  
\end{pf}

We now prove that  ${\cal I}  $ is an  approximate right inverse of $ {\cal L } $ satisfying  \eqref{error:ap-I}.

\begin{lemma} \label{approxinvs3}
\eqref{error:ap-I} holds.
\end{lemma}

\begin{pf} 
Recalling that $ {\cal L} = {\cal L}_D + \varrho $ by \eqref{def:calLD}, and 
setting  
\be\label{defh'new}  
h' := {\cal I} g' \stackrel{\eqref{def:inv-ap-I}} = ({\rm Id}  + {\cal I}_D \varrho)^{-1} ({\cal I}_D g') \, , 
\ee
we have 
\begin{align}
({\cal L} \, {\cal I}  - {\rm Id}) g' & = 
{\cal L}_D h' + \varrho h' - g'  \nonumber \\
&= {\cal L}_D ({\rm Id} 
+ {\cal I}_D \varrho)h'- {\cal L}_D {\cal I}_D \varrho h' + \varrho h' -g' \nonumber \\
& = {\cal L}_D {\cal I}_D  g' - {\cal L}_D {\cal I}_D \varrho h' + \varrho h' -g' \nonumber \\
& =  ( {\cal L}_D {\cal I}_D - {\rm Id} ) (g'- \varrho h' ) \, . \label{LIDI}
\end{align}
Then we estimate \eqref{LIDI} as 
\begin{align}
& \| ({\cal L} {\cal I}  - {\rm Id}) g' \|_{\Lip, s_1} = 
 \| ( {\cal L}_D {\cal I}_D - {\rm Id} ) (g'- \varrho h' ) \|_{\Lip, s_1}  \nonumber \\
& \stackrel{\eqref{stima:error}, \eqref{Rrn-tame-nu-Un}} 
 {\lesssim_{s_3}} N^{Q'+1-(s_3-s_1)} 
\Big(   \| g' \|_{\Lip, s_3} + \| \varrho h' \|_{\Lip, s_3}  \nonumber \\
& \qquad \qquad \qquad \qquad \quad \qquad + 
 ( N^{\varsigma (s_3-s_1)} +    U (s_3) ) ( \| g' \|_{\Lip, s_1} + \| \varrho h' \|_{\Lip, s_1} ) \Big) \nonumber \\
& 
\stackrel{\eqref{opernorm-Lip}, 
\eqref{defh'new}, \eqref{est:Is1-tot}} {\lesssim_{s_3}} N^{Q'+1-(s_3-s_1)} 
 \Big( \| g' \|_{\Lip, s_3} + |\varrho|_{\Lip , s_1} \|h'\|_{\Lip, s_3} +  |\varrho|_{\Lip , s_3} N^{Q' }\|g'\|_{\Lip, s_1} \nonumber \\
 &\qquad \qquad \qquad \qquad  \qquad \qquad \quad + 
  ( N^{\varsigma (s_3-s_1)} +    U (s_3))  
  \big( \|g'\|_{\Lip, s_1} + |\varrho|_{\Lip , s_1} N^{Q' }\|g'\|_{\Lip, s_1} \big) \Big) \nonumber \\
 &
 \stackrel{\eqref{defh'new}, \eqref{est:Iss-tot}, \eqref{smallness:rho}, \eqref{Rrn-tame-nu-Un}} 
 {\lesssim_{s_3}} N^{Q'+1-(s_3-s_1)} 
\Big( \| g' \|_{\Lip, s_3} + 	\big( N^{\varsigma (s_3-s_1) }+ N^{ Q'} U (s_3) \big) 
\| g' \|_{\Lip, s_1} \Big) 
\nonumber 
\end{align}
which proves \eqref{error:ap-I}. 
\end{pf}

\section{Approximate right inverse of $ \bar \om_\e \cdot \pa_\vphi - J (A_0  + \rho) $}\label{sub:choice:N-ind}

In this section we   complete the proof of  Proposition \ref{prop-cruciale}. 
As said in Section \ref{s101}, by the hypotheses of  Proposition \ref{prop-cruciale},
the assumptions of Corollary \ref{cor:split} are satisfied by $(A_0, \rho)$.
Corollary \ref{cor:split} then provides, 
 for all $ \l $ in $ \Lambda_\infty (\e; 5/6, A_0, \rho) $, the conjugation
 (see \eqref{goal-trasf-iterata-ancora})
$$
\big( \bar \om_\e \cdot \partial_\vphi - J  A_0  -J \rho  \big) {\cal P}_\ind (\vphi) = 
{\cal P}_\ind (\vphi) \big( \bar \om_\e \cdot \partial_\vphi - J A_\ind  -J \rho_\ind \big)  \, ,  
$$
where $ A_\ind $ 
is a self-adjoint block diagonal operator of the form 
$$
A_\ind  =  \frac{D_V}{1+ \e^2 \l}  + R_\ind  \, , 
$$
see \eqref{defAm}-\eqref{form-A'-m0}, satisfying
\eqref{der-lambda-m}-\eqref{Rrn-tame0}, in particular $A_\ind \in {\cal C}(2C_1, c_1/2, c_2/2)$ is 
a split admissible operator according to  Definition \ref{def:calC}.  
The sequence $({\cal P}_\ind (\vphi))$ of symplectic transformations satisfies \eqref{Pins1}-\eqref{Pinsany}
and $R_\ind$, $\rho_\ind$ satisfy \eqref{stime-rm-low}-\eqref{Rrn-tame}; in particular
\be \label{raprhon}
|\rho_\ind|_{\Lip,+, s_1} \leq \d_1^{(\frac32)^\ind} \, .
\ee
We  define, for $1/2 \leq \eta \leq 5/6 $,  the sets 
\be\label{set-mathtt-Lambda}
{\bf \Lambda} (\e; \eta, A_0, \rho) := \Lambda_\infty (\e; \eta , A_0, \rho) \bigcap 
\Big( \bigcap_{\ind \geq 0} {\mathtt \Lambda} (\e;\cc + \eta_\ind ,A_\ind) \Big) 
\ee
where  $  \Lambda_\infty $ is defined in Corollary \ref{cor:split},  ${\mathtt \Lambda} $ 
in Proposition \ref{lem:app-inv-LD} (i.e. \eqref{sets:L-inverse})
and the sequence $(\eta_\ind)$  is defined in  \eqref{def:indices}.
Note that we can apply Proposition \ref{lem:app-inv-LD} with $\FA=A_\ind$ for any $\ind \geq 0$,  because $A_\ind \in {\cal C}(2C_1, c_1/2, c_2/2)$.
In Lemma \ref{lem:finame} below we shall prove that the
 sets $ {\bf \Lambda} (\e; \eta, A_0, \rho)  $ satisfy properties 1-3 of 
Proposition \ref{prop-cruciale}. 
\\[2mm]
{\sc Choice of the cut-off $ N $ and of the number $ \ind $ of splitting steps.}  
\\[1mm]
For any $ \nu \in (0, \e ) $, we choose $N \in \N$ such that
\be \label{defNsp-nu} 
N \in \Big[ \nu^{-\frac{3}{s_3-s_1}} - \frac12, \nu^{-\frac{3}{s_3-s_1}} + \frac12 \Big) 
\ee 
and  the number  $ \ind \in \N $ of splitting steps 
in Corollary \ref{cor:split}  as  
$$
\ind := \min \Big\{ k \in \N \, : \,   \d_1^{(\frac32)^{k }} N^{Q'}=\d_1^{(\frac32)^{k }} N^{2(\t' + \varsigma s_1) + 3} \leq c_0 \Big\} 
$$
where $ \d_1 = \e^3 $ (see \eqref{def:d1} and \eqref{rho-R0:small-0}) and $c_0$ is the 
strictly positive constant of Proposition  \ref{lemma:ap-inv_L}. 
Hence
\be\label{choice:ind}
\d_1^{(\frac32)^{\ind }} N^{Q'} \leq c_0 \, , \qquad \hbox{and, if} \ \ \ind \geq 1 \, , \quad  
\d_1^{(\frac32)^{\ind -1}} N^{Q'} > c_0 \, .
\ee
\begin{remark} \label{rem:nuep}
In the sequel we suppose that $ \nu $ is small enough, possibly depending on $ s_3 $
(the aim is to obtain estimates without any multiplicative constant depending on $s_3$). 
Since $\nu \in (0,\e)$, any smallness condition for $\nu$ is satisfied if $\e$ is small enough,
possibly depending on $s_3$, which is a  large, but fixed, constant. 
For this reason, often  we will  not explicitly indicate such
dependence.  However it is useful to point out that, for $\e$ fixed, 
our estimates still hold for $\nu$
small enough. In fact, in Proposition \ref{crucial-s}, we  extend some estimates on the 
$s_3$-norm to estimates on the $s$-norm for any $s\geq s_3$, with a smallness condition on $\e$
which does not depend on $s$ but only on $s_3$, while the smallness condition on $\nu$  depends on $s$.
\end{remark}

In the sequel of this section the notation $a \ll b$ means that $a/b \to 0$ as 
$\nu \to 0$. 

By \eqref{defNsp-nu} that,   
for $ \nu $ small enough,   $ N \geq \bar{N}$.

By \eqref{defNsp-nu},  \eqref{choice s2 s3} and recalling that $ Q' = 2(\t' + \varsigma s_1) + 3 $ and 
$ \varsigma = 1 /10 $, we  have 
\begin{align}\label{N-versus-vu}
& N^{Q'+1} \ll  \nu^{- \frac{1}{20}} \, , \quad N^{\varsigma (s_3- s_1)} \ll \nu^{- \frac25} \, , 
\quad N^{-(1- \varsigma)(s_3-s_1)} \ll \nu^{\frac52} \, \\
& \quad N^{Q'+1-(s_3-s_1)}  \leq \nu^{- \frac{1}{20}} N^{-(s_3 -s_1)} \lesssim_{s_3} \nu^{- \frac{1}{20}} \nu^3 \ll \nu^{\frac{11}{4}} \, .
\label{N-versus-vu-bis}
\end{align} 
Moreover, by \eqref{defNsp-nu}, \eqref{choice:ind},  \eqref{choice s2 s3},
we have 
\begin{align}
& \d_1^{(\frac32)^{\ind }}  \lesssim   \nu^{\frac{3 Q'}{s_3-s_1}}  \ll \nu^{-\frac1{20}}\, , 
\quad
 \ind \leq \frac{1}{\ln (3/2)} \ln  \ln \nu^{-1} \, , \nonumber \\
&  [C(s)]^\ind = e^{\ind \ln C(s) }  \leq e^{ \frac{\ln C(s)}{\ln (3/2)} \ln \ln \nu^{-1}  } 
\leq \big( \ln \nu^{-1}  \big)^{ \frac{\ln C(s)}{\ln (3/2)} } \, , \label{first-disu}
\end{align}
and,  for $\ind \geq 1$,
\be\label{second-disu}
\d_1^{- ( \frac32)^{\ind -1}}  \lesssim \nu^{\frac{ - 3 Q'}{s_3-s_1}} \ll  \nu^{-\frac1{20}}  \, .
\ee
Recalling \eqref{Rrn-tame-nu-Un},  we introduce the notation
\be\label{Rrn-tame-nu-Un-ind}
U_\ind (s) = U_{A_\ind, \rho_\ind} (s) = |R_\ind|_{\Lip,+, s} + |\rho_\ind|_{\Lip,+, s} \, .
\ee
The bounds \eqref{Rrn-tame0}, \eqref{Rrn-tame} 
where  $ \a (s) = 3 \varsigma \frac{s-s_2}{s_2-s_1} $  
and \eqref{first-disu}, \eqref{second-disu} provide 
the estimate 
 \begin{eqnarray}\label{Rrn-tame-nu}
U_\ind (s_3) 
& \lesssim_{s_3} &  2 ( \ln (\nu^{-1}) )^{ \frac{\ln C(s_3)}{\ln (3/2)} }   \nu^{ - \frac{3Q'}{s_3-s_1} ( \frac34 + 
3\varsigma \frac{s_3-s_2}{s_2-s_1}  )} \big[ \big( |R_0|_{\Lip, +,s_3} + |\rho|_{\Lip, +,s_3} \big)
\d_1^{\frac{1}{2} + \frac{2\alpha (s_3) }{3}} +1 \big] \nonumber \\
& \stackrel{\eqref{g:small-large}}  {\lesssim_{s_3}} & 2 ( \ln (\nu^{-1}) )^{ \frac{\ln C(s_3)}{\ln (3/2)} }  
 \nu^{ - \frac{3Q'}{s_3-s_1} ( \frac34 + 3\varsigma \frac{s_3-s_2}{s_2-s_1}  )} \big[ \nu^{-1} \e^2 +1 \big] \nonumber  \\ 
 &  \lesssim_{s_3}   &  2 ( \ln (\nu^{-1}) )^{ \frac{\ln C(s_3)}{\ln (3/2)} }  
  \nu^{ - \frac{9Q'}{4(s_3-s_1)} } \nu^{ - \frac{9Q' \varsigma }{s_2-s_1} }  \big[ \nu^{-1} \e^2 +1 \big] \nonumber \\
 & \stackrel{ \eqref{choice s2 s3}}\leq & C(s_3) \nu^{- \frac{11}{10}}  \, .   \label{Uns3}
\end{eqnarray} 
The estimate \eqref{Pinsany} provides similarly
\be \label{Pns3}
|{\cal P}_\ind^{\pm 1} |_{\Lip, +, s_3} \leq C(s_3) \nu^{-\frac{11}{10}} \, . 
\ee
\noindent
{\sc Solution of \eqref{sol:almost-inv}.} 
\\[1mm]
By  \eqref{raprhon} and 
\eqref{choice:ind},  
the smallness condition \eqref{smallness:rho} is satisfied by $ \varrho = \rho_\ind $. Then Proposition \ref{lemma:ap-inv_L} applies
to the operator ${\cal L}:=J \bar {\om }_\e \cdot \partial_\varphi + A_\ind  + \rho_\ind $,
implying the existence of an approximate right inverse ${\cal I} := {\cal I}_{N,\ind}$, defined for $\l \in {\bf \Lambda} (\e; \eta, A_0, \rho) $ 
(see \eqref{set-mathtt-Lambda}),  satisfying  \eqref{error:ap-I}-\eqref{new:cal-loss}.
The operator $ {\cal I} $ satisfies 
 \be \label{s1Ig}
 \| {\cal I} \wtilde{g} \|_{\Lip, s_1} 
\stackrel{ \eqref{est:Is1-tot}} {\lesssim_{s_1}} N^{Q'} \|\wtilde{g} \|_{\Lip, s_1}
 \stackrel{\eqref{N-versus-vu}} \leq \nu^{-\frac1{20}} \|\wtilde{g} \|_{\Lip, s_1} 
 \ee
and  
\begin{align}
\| {\cal I} \wtilde{g} \|_{\Lip, s_3}  & 
\stackrel{ \eqref{est:Iss-tot}} \leq 
C N^{Q'}  \|\wtilde{g} \|_{\Lip, s_3} + C(s_3) \big(  N^{\varsigma (s_3-s_1)} + N^{Q'} U_\ind (s_3) \big)  \|\wtilde{g} \|_{\Lip, s_1}
\nonumber \\
&  \stackrel{\eqref{N-versus-vu}, \eqref{Uns3}} \leq  \nu^{- \frac{1}{20}} \| \wtilde{g} \|_{\Lip, s_3} +
\big( \nu^{- \frac25} + \nu^{- \frac{1}{20}}  \nu^{- \frac{11}{10}} \big) \| \wtilde{g} \|_{\Lip, s_1} \nonumber \\
&\leq   \nu^{- \frac{1}{20}} \| \wtilde{g} \|_{\Lip, s_3} + \nu^{- \frac{6}{5}} \| \wtilde{g}\|_{\Lip, s_1} \, .  \label{estInu}
\end{align}
In addition estimates \eqref{error:ap-I}, \eqref{N-versus-vu}-\eqref{N-versus-vu-bis},  \eqref{Rrn-tame-nu} give
\begin{align}  
\| ({\cal L I} - {\rm Id})  \wtilde{g} \|_{\Lip, s_1} & \leq    
\nu^{\frac{11}{4}}  \Big(  \| \wtilde{g} \|_{\Lip, s_3}   +
\big( \nu^{-\frac25} + \nu^{-\frac1{20} } \nu^{-\frac{11}{10}} \big) \| \wtilde{g} \|_{\Lip, s_1} \Big)
\nonumber \\
&\leq 
\nu^{\frac{11}{4}} \|\wtilde{g} \|_{\Lip , s_3} + \nu^{\frac32}  \|\wtilde{g}  \|_{\Lip, s_1}   \, .   \label{estapproxnu}
\end{align}
Now  let $ g \in H_{\mathbb S}^\bot $ satisfy the assumption \eqref{g:small-large} and 
consider the function 
 $ g' = {\cal P}_\ind^{-1} (\vphi ) g $ introduced in  \eqref{def:change-var}. 
We define, for any $ \l $ in $ {\bf \Lambda} (\e;\eta , A_0, \rho) $,  
the approximate solution $ h' = {\cal I} ( J g' ) $ of the equation \eqref{def:calLD}
where $ {\cal I} $ is  the  approximate right inverse, 
obtained in Proposition \ref{lemma:ap-inv_L}, of the operator 
$ {\cal L } = J \bar {\om }_\e \cdot \partial_\varphi + A_\ind  + \rho_\ind  $, and
we consider the  ``approximate solution" $ h $ of the equation \eqref{sistema-da-riso} as
\be\label{def:sol-finale}
h := {\cal P}_\ind (\vphi) h'  \qquad {\rm where} 
\qquad h' := {\cal I} (J g' ) \, , \quad  g' = {\cal P}_\ind^{-1} (\vphi) g \, .
\ee
It means that we define the approximate right inverse 
of $ \bar {\om }_\e \cdot \partial_\varphi - J (A_0 + \rho) $ as
\be\label{prop 8.1:approximate-inv}
 {\mathfrak L}^{-1}_{approx}  :=  {\cal P}_\ind (\vphi)  {\cal I} J {\cal P}_\ind (\vphi)^{-1}  \, . 
\ee
We claim that the function $ h $ defined in \eqref{def:sol-finale}
is the required solution of \eqref{sol:almost-inv} with a remainder $ r $ satisfying 
\eqref{estimate:r-final}. 

\begin{lemma} \label{esthnu}
The approximate solution $ h $ defined in \eqref{def:sol-finale} 
satisfies \eqref{h:s1s2s3}  and  \eqref{con-loss}. 
\end{lemma}

\begin{pf}
We have 
\begin{align}\label{bound:hs1-ok}
\| h \|_{\Lip, s_1} & 
\stackrel{ \eqref{def:sol-finale} } =  \| {\cal P}_\ind \, {\cal I} (J  {\cal P}_\ind^{-1}  g ) \|_{\Lip, s_1} 
\stackrel{ \eqref{Pins1}, \eqref{opernorm-Lip}} {\lesssim_{s_1}}  \|{\cal I} (J  {\cal P}_\ind^{-1} g)  \|_{\Lip, s_1} \nonumber \\
& \stackrel{ \eqref{s1Ig}  }  {\lesssim_{s_1}} \nu^{ - \frac{1}{20}} \|  {\cal P}_\ind^{-1}  g \|_{\Lip, s_1} 
\stackrel{ \eqref{Pins1} } {\lesssim_{s_1}} \nu^{ - \frac{1}{20}} \| g \|_{\Lip, s_1}  
\lesssim_{s_1}  \e^2 \nu^{\frac{19}{20}} 
\end{align}
by  the assumption $ \| g \|_{\Lip, s_1} \leq \e^2 \nu $, see \eqref{g:small-large}. Therefore
$ \| h \|_{\Lip, s_1}  \leq \e^2 \nu^{4/5} $, for $\nu$ small enough (depending on $s_1$), 
 proving the first inequality in  \eqref{h:s1s2s3}. 

Let us now estimate $\| h \|_{\Lip, s_3}$. We have 
 \begin{align}
\|  g' \|_{\Lip, s_3} & \stackrel{\eqref{def:sol-finale}} = 
\|  {\cal P}_\ind^{-1} g \|_{\Lip, s_3} \nonumber \\
& \stackrel{\eqref{opernorm-Lip}} {\lesssim_{s_3}}  | {\cal P}_\ind^{-1}|_{\Lip,  s_3} \| g \|_{\Lip, s_1} +
| {\cal P}_\ind^{-1}|_{\Lip,  s_1} \| g \|_{\Lip, s_3} \nonumber \\ 
& \stackrel{\eqref{Pns3}, \eqref{Pins1}} {\lesssim_{s_3}} 
 \nu^{- \frac{11}{10}} \| g \|_{\Lip, s_1} + \| g \|_{\Lip, s_3} \nonumber \\
& \stackrel{\eqref{g:small-large}}  {\lesssim_{s_3}}   \e^2 (\nu^{- \frac{1}{10}} + \nu^{-1}) 
 \lesssim_{s_3} \e^2 \nu^{-1} \, . \label{stima:ulti}
\end{align}
Then  \eqref{estInu}, \eqref{stima:ulti} and $ \| g' \|_{\Lip, s_1} \lesssim_{s_1} \| g \|_{\Lip, s_1} \lesssim_{s_1} \e^2 \nu $,  imply 
\begin{align} 
\|  {\cal I} (J g' ) \|_{\Lip, s_3} &
\lesssim_{s_3} \e^2 (\nu^{-\frac1{20}} \nu^{- 1}+ \nu^{-  \frac{6}{5}} \nu) 
\lesssim_{s_3}  \e^2  \nu^{-\frac{21}{20}}  \label{eq:inter-fi}
\end{align}
and finally, using \eqref{opernorm-Lip},  \eqref{Pins1},  \eqref{Pns3},\eqref{eq:inter-fi} and
$ \|{\cal I} (J g' ) \|_{\Lip, s_1} \lesssim_{s_1}  \e^2 \nu^{\frac{19}{20}}  $  (see \eqref{bound:hs1-ok}), we conclude that 
\begin{align}\label{bound:hs3}
\| h \|_{\Lip, s_3} = \| {\cal P}_\ind  \, {\cal I} (J g' ) \|_{\Lip, s_3} \lesssim_{s_3} 
\e^2 \big(  \nu^{- \frac3{20}} + \nu^{- \frac{21}{20}} \big) \leq  \e^2 \nu^{-\frac{11}{10}} 
\end{align}
for $\nu$ small enough (depending on $s_3$). This proves the second inequality in  \eqref{h:s1s2s3}.   

At last, by \eqref{Pins1} and \eqref{new:cal-loss}, for any $g \in {\mathcal H}^{s_0 +Q'} \cap H^{\bot}_{\mathbb S}$, 
$$
\| {\mathfrak L}^{-1}_{approx} g \|_{\Lip, s_0} =\| {\cal P}_\ind {\cal I}  J {\cal P}^{-1}_\ind g\|_{\Lip, s_0}
\lesssim_{s_1} \| J {\cal P}^{-1}_\ind g \|_{\Lip, s_0+Q'}  \lesssim_{s_1} \| g \|_{\Lip, s_0+Q'}  
$$
using that $s_0 + Q' \leq s_1$.  This proves  \eqref{con-loss}.
\end{pf}

We finally estimate the error $ r $ in the equation  \eqref{sol:almost-inv}. 

\begin{lemma} 
$ \|r \|_{\Lip, s_1} \leq \e^2 \nu^{3/2} $ , where 
$r = \big( \bar \om_\e \cdot \partial_\vphi - J (A_0  + \rho) \big)h  -   g  $.
\end{lemma}

\begin{pf}
By  \eqref{goal-trasf-iterata-ancora}, \eqref{def:sol-finale} and recalling that 
$ {\cal L} = J \bar {\om }_\e \cdot \partial_\varphi + A_\ind  + \rho_\ind$ we have 
\begin{align} 
r& = \big( \bar \om_\e \cdot \partial_\vphi - J (A_0  + \rho) \big)h  -   g \nonumber  \\
& = {\cal P}_\ind \big[ \big( \bar \om_\e \cdot \partial_\vphi - J (A_\ind  + \rho_\ind ) \big) h' - g' \big] \nonumber  \\
&    = - {\cal P}_\ind J \big[ {\cal L} h' - J g' \big]  \nonumber \\ 
&  = - {\cal P}_\ind J \big[ ({\cal L} {\cal I}  - {\rm Id}) J g'
\big] \label{forma-r}
\, .
\end{align}
Now, by the estimates in the proof of Lemma \ref{esthnu}, see \eqref{stima:ulti}, 
we know that $g'={\cal P}_\ind^{-1} g$ satisfies
\be \label{estg'nu}
\|g' \|_{\Lip , s_3} \lesssim_{s_3}  \e^2 ( \nu^{-1} + \nu^{- \frac{1}{10}}) 
\lesssim_{s_3} \e^2 \nu^{-1} \, , 
\qquad
\| g' \|_{\Lip , s_1} \lesssim \e^2 \nu \, . 
\ee
Hence  by  \eqref{forma-r}, \eqref{opernorm-Lip},   \eqref{Pins1},  we get 
\begin{align}
\| r \|_{\Lip, s_1} &  \lesssim   \|  ({\cal L} {\cal I}  - {\rm Id}) J g' \|_{\Lip , s_1} \nonumber \\
&  \stackrel{ \eqref{estapproxnu}} 
\lesssim \nu^{\frac{11}{4}} \| g'  \|_{\Lip, s_3} +  \nu^{ \frac32} \| g' \|_{\Lip, s_1}   \nonumber \\
& \stackrel{ \eqref{estg'nu}} {\lesssim_{s_3}} \e^2 ( \nu^{\frac74} + \nu^{\frac52})  \leq \e^2 \nu^{\frac32} \nonumber 
\end{align} 
proving the lemma. 
\end{pf}


The next lemma completes  the proof of  Proposition \ref{prop-cruciale}. 

\begin{lemma} {\bf (Measure estimates)}\label{lem:finame}
The sets ${\bf \Lambda} (\e;\eta,A_0,\rho)$, $ 1/2 \leq \eta \leq 5/ 6 $, 
defined in \eqref{set-mathtt-Lambda} satisfy properties 1-3 of 
Proposition \ref{prop-cruciale}. 
\end{lemma}

\begin{pf}    
Property 1 follows immediately because the sets 
$  \Lambda_\infty $ defined in Corollary \ref{cor:split}, 
 and ${\mathtt \Lambda} $  in Proposition \ref{lem:app-inv-LD} are increasing in $ \eta $.
\\[1mm]
For the proof of properties \ref{AApro2} and \ref{AApro3} , we observe that whereas $A_0$ is defined for 
any $\l \in \wtilde{\Lambda}$, by Corollary \ref{cor:split} the operators $A_\ind$, $\ind \geq 1$, are defined for 
$ \l \in \Lambda_\infty(\e ; 5/6, A_0, \rho) \subset \wtilde{\Lambda}$. Hence for $\ind \geq 1$,
the set ${\mathtt \Lambda} (\e , \eta , A_\ind)$ ($1/2 \leq \eta \leq 1$) is considered as 
a subset of $\Lambda_\infty(\e ; 5/6, A_0, \rho)$, and we shall apply Proposition \ref{lem:app-inv-LD}
(more precisely Lemma \ref{measureL}) in this setting.
\\[1mm]
{\sc Proof of property \ref{AApro2} of Proposition \ref{prop-cruciale}.}
Definining
\be\label{def:Laminf}
{\mathtt \Lambda}_\ind  (\e;\cc , A_0, \rho) := \Lambda_\infty (\e;\cc,A_0,\rho) 
\bigcap \Big( \bigcap_{k=0}^{\ind-1} {\mathtt \Lambda} (\e; \cc + \eta_k , A_k)  \Big)  \, ,  \quad \forall  \ind \geq 1 \, ,
\ee
we write the set in \eqref{set-mathtt-Lambda} (recall that $ \eta_0 = 0 $) as 
\be\label{BigL}
{\bf \Lambda} (\e; \cc, A_0, \rho) = 
\Lambda_\infty (\e; \cc , A_0, \rho) 
\bigcap  {\mathtt \Lambda} (\e;  \cc ,A_0) 
\bigcap \Big( \bigcap_{\ind \geq 1} {\mathtt \Lambda}_\ind (\e;\cc , A_0, \rho) \Big) \, .
\ee
Then its complementary set may be decomposed as (arguing as in \eqref{compl-set-Linfty}) 
\begin{align}
{\bf \Lambda}(\e; \cc, A_0, \rho)^c  & =\Lambda_\infty (\e; \cc, A_0, \rho)^c   \bigcup
{\mathtt \Lambda} (\e; \cc, A_0)^c  \nonumber \\
& \quad \bigcup 
\Big(  \bigcup_{\ind \geq 1} \big( {\mathtt \Lambda}_\ind (\e; \cc, A_0, \rho) \bigcap
{\mathtt \Lambda} (\e; \cc +\eta_\ind , A_\ind)^c \big)  \Big) \, . \label{inclusL}
\end{align}
By property \ref{item2-92} of Corollary \ref{cor:split} we have 
\be \label{est1L}
\big|  \Lambda_\infty (\e; 1/2, A_0, \rho)^c  \cap  \wtilde{\Lambda} \big| \leq b_1(\e)  \quad {\rm with} \quad  
\lim_{\e \to 0} b_1(\e)=0 \, .
\ee
By Lemma \ref{measureL},
\be \label{est2L}
\big|   {\mathtt \Lambda} (\e ; 1/2, A_0)^c  \cap \wtilde{\Lambda} \big| \lesssim  \e \, . 
\ee
Moreover, for $\ind \geq 1$,   by \eqref{def:Laminf}, we have 
$$ 
\begin{aligned} 
& {\mathtt \Lambda}_\ind (\e; \cc , A_0, \rho) \cap
{\mathtt \Lambda} (\e; \cc +\eta_\ind , A_\ind)^c \\
& \subset 
\Lambda_\infty (\e ; 5/6 , A_0, \rho) \cap
{\mathtt \Lambda} (\e; \cc +\eta_{\ind -1} , A_{\ind -1}) \cap
{\mathtt \Lambda} (\e; \cc +\eta_\ind , A_\ind)^c \, . 
\end{aligned}
$$
Then,  since, for all $ \ind \geq 1 $, we have 
$ |A_\ind - A_{\ind -1}|_{+,s_1} \leq \d_\ind^{3/4} $  on $ \Lambda_\infty (\e ; 5/6 , A_0, \rho)  $ (see \eqref{An-An-1}) 
and $  \cc + \eta_{\ind -1} =  \cc + \eta_\ind - \d_\ind^{3/8} $ (see \eqref{def:indices}),  
we deduce by Lemma \ref{measureL} that, for any $1/2 \leq \cc \leq 5/6$,
\be \label{est3L}
\big| {\mathtt \Lambda}_\ind (\e; \cc, A_0, \rho) \bigcap
{\mathtt \Lambda} (\e ; \cc +\eta_\ind , A_\ind)^c   \big|
\leq 
 \d_\ind^{3\alpha /4} \, .
\ee
In conclusion, by \eqref{inclusL}  and \eqref{est1L},  \eqref{est2L},  \eqref{est3L} at $ \eta = 1 / 2 $, we deduce 
$$
\begin{aligned}
\big|   {\bf \Lambda}(\e; 1/2, A_0, \rho)^c  \cap \wtilde{\Lambda}  \big| 
& \leq
b_1(\e) + C \e + \sum_{\ind \geq 1} \d_\ind^{3\alpha/4} \\
& \leq 
b_1(\e) + C \e + c \d_1^{3\alpha/4} \leq b(\e) 
\end{aligned}
$$
with $\lim_{\e \to 0} b(\e)=0 $, since  $\d_1=\e^3 $. This proves property \ref{AApro2} of 
Proposition \ref{prop-cruciale} for the sets ${\bf \Lambda} (\e; \cc, A_0, \rho)$.
\\[1mm] 
{\sc Proof of property \ref{AApro3} of Proposition \ref{prop-cruciale}.}
By  \eqref{inclusL} 
(with $ A_0', \rho' $ instead of $ A_0, \rho $) and  \eqref{BigL} we deduce,  
for all $ (1/2) + \d^{2/5} \leq \cc \leq 5/6 $,  the inclusion 
\begin{align} \label{inclMfin1}
{\cal M} & := \wtilde{\Lambda}' \bigcap {\bf \Lambda} (\e ; \cc, A'_0 , \rho')^c
\bigcap {\bf \Lambda} (\e; \cc- \d^{2/5} , A_0, \rho) \nonumber \\
& \subset
{\cal M}_\infty \bigcup {\cal M}_0 \bigcup  \Big( \bigcup_{\ind \geq 1} {\cal M}_\ind   \Big) \, ,
\end{align}
where
\begin{align*}
& {\cal M}_\infty := \wtilde{\Lambda'} \cap [\Lambda_\infty (\e ; \cc , A'_0, \rho')]^c \cap \Lambda_\infty (\e ; \cc - \d^{2/5} , A_0, \rho) \\
& {\cal M}_0 := \wtilde{\Lambda}' \cap {\mathtt \Lambda} (\e ; \cc, A'_0 )^c
\cap {\mathtt \Lambda} (\e; \cc- \d^{2/5} , A_0) 
\\
& {\cal M}_\ind := {\mathtt \Lambda}_\ind (\e; \cc, A'_0, \rho') \cap {\mathtt \Lambda} (\e; \cc + \eta_\ind , A'_\ind)^c
\cap {\mathtt \Lambda}(\e; \cc + \eta_\ind-\d^{2/5} , A_\ind) \, , \quad {\rm for} \ \ind \geq 1 \, . 
\end{align*} 
By  \eqref{smallR0R0'} and noting that  $ \d^{1/2} \leq \d^{2/5} $, we deduce by
\eqref{Cantor:intersect} that, for all $ (1/2) + \d^{2/5} \leq \cc \leq 5/6 $, 
\be\label{Minfty1}
| {\cal M}_\infty |  \leq \d^{\alpha /2} \, .
\ee
Moreover Lemma \ref{measureL}  implies that
\be\label{Minfty2}
| {\cal M}_0  | \leq \big| \wtilde{\Lambda}' \cap {\mathtt \Lambda} (\e ; \cc, A'_0 )^c
\cap {\mathtt \Lambda} (\e; \cc- \d^{1/2} , A_0)   \big| \leq \d^\alpha \, .
\ee
For $ \ind \geq 1 $, we have,  by \eqref{est3L} (applied to $(A'_0, \rho')$)
\be\label{Minfty3}
| {\cal M}_\ind | \leq \big| {\mathtt \Lambda}_\ind (\e; \cc, A'_0, \rho') \cap {\mathtt \Lambda} (\e; \cc + \eta_\ind , A'_\ind)^c   \big|
\leq  \d_\ind^{3 \alpha /4}  \, .
\ee
On the other hand,  assumption  \eqref{smallR0R0'} implies \eqref{pertA0}, and so 
\eqref{pertAn} holds, i.e.  
{ $|A_\ind - A'_\ind|_{+, s_1} \leq \d^{4/5}$  for any $ \l $ in 
$ {\cal M}_\ind \subset 
\Lambda_\infty (\e; 5/6 , A'_0 , \rho') \bigcap \Lambda_\infty (\e;  5/6 , A_0 , \rho) $, 
and we deduce, by Lemma \ref{measureL},  the measure estimate}
\be\label{Minfty2.5}
\begin{aligned}
| {\cal M}_\ind | & 
\leq  \big|  \Lambda_\infty (\e ; 5/6 , A'_0, \rho')  \cap
{\mathtt \Lambda} (\e; \cc + \eta_\ind , A'_\ind)^c \cap {\mathtt \Lambda} (\e; \cc + \eta_\ind - \d^{2/5}, A_\ind) \big| \\
& \lesssim \d^{4\alpha /5} \, .
\end{aligned}
\ee
Finally \eqref{inclMfin1}, \eqref{Minfty1}, \eqref{Minfty2}, \eqref{Minfty2.5}, \eqref{Minfty3} imply the measure estimate
\be \label{me1}
| {\cal M} | \leq 
\d^{\alpha/2} + \d^\alpha + C \sum_{\ind \geq 1}  \min (\d^{4\alpha /5} , \d_\ind^{3\alpha/4} ) \, . 
\ee
Now, using that  $\dps \d_\ind =\d_1^{(  \frac32 )^{\ind-1}} $ with $ \d_1 = \e^3 $,  
there is a constant $C>0$ such that 
$$
\d^{\frac{4\a}{5}} \leq \d_\ind^{ \frac{3\a}{4}} \quad  \Longrightarrow  \quad  \ind \leq C \ln ( \ln (\d^{-1})) \, .
$$
Hence
\begin{align}  \label{me2-10}
 \sum_{\ind \geq 1}  \min (\d^{ \frac{4\alpha}{5}} , \d_\ind^{ \frac{3\alpha}{4}} )  &\leq C  \ln ( \ln (\d^{-1})) \d^{ \frac{4\a}{5}}  +
 \sum_{\ind \geq 1 ,   \d_\ind^{ \frac{3\a}{4} } < \d^{ \frac{4\a}{5}}}  \d_\ind^{ \frac{3\alpha}{4}}   \nonumber \\
 &\leq  
 C  \ln ( \ln (\d^{-1})) \d^{ \frac{4\a}{5}}  + C   \d^{ \frac{4\a}{5}} \, .  
\end{align}
The estimates \eqref{me1} and \eqref{me2-10} imply, for $ \d $ small, that  
$$
| {\cal M} | \leq \d^{\alpha /3} \, ,  
$$ 
which is \eqref{finadav}. 
\end{pf}
\\[3mm]
{\sc Proof of Proposition \ref{crucial-s}.}  
The result of Proposition \ref{prop-cruciale} holds for any fixed $s_3$ 
satisfying \eqref{choice s2 s3}.  
The sets $ {\bf \Lambda} (\e; \cc ,A_0, \rho) $ do not depend on the choice of $s_3$.
On the other hand, the smallness condition for $\e$ may depend on $s_3$ but, as explained in Remark \ref{rem:nuep},  what is really used in the proof of Proposition \ref{prop-cruciale} is that
$\nu$ is small enough. This naturally leads to Proposition \ref{crucial-s}, the proof of which is exactly the same as for Proposition \ref{prop-cruciale}, 
replacing $s_3$ by $ s $, and 
the smallness condition on $\e$ by a smallness condition on $\nu$ (depending on $s$). 
Notice that the approximate inverse ${\mathfrak L}^{-1}_{approx} $ 
depends not only on $ \nu $, but also on $ s$, in particular through the choice of 
$ N $, that, replacing $ s_3 $ with $ s $ in   \eqref{defNsp-nu}, results in
$$
N \in \Big[ \nu^{-\frac{3}{s-s_1}} - \frac12, \nu^{-\frac{3}{s-s_1}} + \frac12 \Big) \, . 
$$

\chapter{Proof of the Nash-Moser Theorem} 
\label{sec:NM}
 
In this chapter we finally prove the Nash-Moser 
Theorem \ref{thm:NM}, finding,  by an  iterative scheme 
(see Section \ref{sec:NM-ite}),  a solution of the nonlinear equation
$ {\cal F} (\lambda;  i ) = 0 $ where $ \cal F $ is the operator  
defined in \eqref{operatorF}. 

By the procedure described in 
Chapter \ref{sezione almost approximate inverse} (see Proposition \ref{Prop:sec6})
 in order to  find an approximate inverse for the linearized operator 
 at an approximate solution $ \ui $,
 and thus implement a convergent Nash-Moser scheme, 
 it is sufficient to prove the existence of  
an approximate right inverse of the operator $ {\mathbb D } (\ui) $ defined in  \eqref{operatore inverso approssimato}. This is achieved in 
Proposition \ref{Prop:inversione}.

\section{Approximate right inverse of $  {\cal L}_\omega $}

We first give the key result about the existence of an approximate right inverse\index{Approximate right inverse in normal directions} 
of the operator $ {\cal L}_\om :=  {\cal L}_\omega(\ui ) $ 
defined in \eqref{Lomega def}, acting on the normal subspace $ H_{\mathbb S}^\bot $. 
We recall that $ \ui (\vphi ) = (\vphi,0,0) + \uF (\vphi) $ is defined for all $ \l \in \Lambda_{\uF}  $. 

\begin{proposition} \label{prop:inv-ap-vero} {\bf (Approximate right inverse of $ {\cal L}_\omega (\ui)  $)}
Let  $ \bar \om_\e \in \R^\es $ be  $(\gamma_1 , \tau_1) $-Diophantine and satisfy 
property $ {\bf (NR)}_{\gamma_1, \tau_1} $ in Definition \ref{NRgamtau} with $ \g_1, \t_1 $ fixed in 
\eqref{def:tau1}. Assume  \eqref{choice s2 s3}. 

Then there is  $ \e_0 > 0  $  such that, $ \forall \e \in (0, \e_0)$, for all $ \uF $,  defined for all $ \l \in  \Lambda_{\uF} $,
satisfying $  \| \uF \|_{\Lip, s_2+2} \leq \e$, 
there  are closed subsets $ {\bf \Lambda} (\e; \cc , \uF  ) \subset \Lambda_{\uF} $,  $ 1/2 \leq \cc \leq 5/6 $,  satisfying 
\begin{enumerate}
\item $  {\bf \Lambda}  (\e ; \cc, \uF ) \subseteq  {\bf \Lambda}  (\e ; \cc', \uF) $ 
for all  $1/2 \leq \cc \leq \cc' \leq 5/6 $, 
\item \label{mmqf2}
$ \big| [ {\bf \Lambda} (\e; 1/2, \uF)]^c \cap \Lambda_{\uF} \big| \leq b(\e) $ where $ \lim_{\e \to 0} b(\e) = 0 $, 
\item 
if $ \| \uF' - \uF \|_{s_1+2}  \leq \delta \leq \e^{3/2} $ for all
$ \l \in \Lambda_{\uF} \cap \Lambda_{\uF'} $, 
then, for all $ (1/2) + \d^{2/5} \leq \cc \leq 5/6 $, 
\be\label{exp-misu}
\big| \Lambda_{\uF'} \cap  [{\bf \Lambda} (\e; \cc, \uF')]^c \cap  {\bf \Lambda} (\e ; \cc -\d^{2/5} , \uF) \big| \leq 
\delta^{\alpha/3} ; 
\ee 
\end{enumerate}
and, for any $ \nu \in (0, \e^{\frac32}) $ such that 
 $\| \uF \|_{\Lip, s_3+4} \leq \e  \nu^{-\frac9{10}} $, 
there exists a linear operator 
$$ 
{\cal L}^{-1}_{appr} :=  {\cal L}^{-1}_{appr, \nu} 
$$ 
such that,   for 
any function $ g : \Lambda_{\uF} \to {\mathcal H}^{s_3+2} \cap H_{\mathbb S}^\bot $ satisfying 
\be\label{g:small-large-IA}
\| g \|_{\Lip, s_1} \leq  \e^2 \nu  \,  , \quad    \| g \|_{\Lip, s_3+2}  \leq  \e^2 \nu^{- \frac{9}{10}} \ ,
\ee
the function  $ h := {\cal L}^{-1}_{appr} g $, $ h : {\bf \Lambda} (\e; 5/6, \uF ) \to  
{\mathcal H}^{s_3+2} \cap H_{\mathbb S}^\bot $  
satisfies
\be\label{h:s1s2s3-IA-per}
 \| h \|_{\Lip, s_1} \leq C(s_1) \e^2 \nu^{\frac45} \, , \quad 
 \| h \|_{\Lip, s_3+2} \leq C(s_3)  \e^2 \nu^{ - \frac{11}{10} } \, , 
\ee
and, setting  $ {\cal L}_\omega := {\cal L}_\omega(\ui ) $  defined in \eqref{Lomega def},  we have 
\be\label{estimate:r-final-IA-per}
 \| {\cal L}_\omega h - g \|_{\Lip, s_1} \leq C(s_1) \e^2 \nu^{\frac32} \,  . 
\ee
Furthermore, 
setting $ Q' := 2(\tau'+ \varsigma s_1 ) +3 $ (where $ \varsigma =1/ 10 $, see \eqref{def:varsigma}) and 
$ \tau' $ are given by Proposition \ref{propmultiscale}), for all $ g \in {\mathcal H}^{s_0 + Q'} \cap  H_{\mathbb S}^\bot $, 
\be\label{new:estimate-cal-L-1}
\| {\cal L}^{-1}_{appr} g \|_{\Lip, s_0} \lesssim_{s_1} \| g \|_{\Lip, s_0 + Q'} \, .
\ee
We underline that the estimate  \eqref{new:estimate-cal-L-1} is independent of $ \nu $ which defines $ {\cal L}_{appr}  = {\cal L}_{appr, \nu} $. 
\end{proposition}

\begin{remark} \label{rem:s3s}
We lay the stress on the fact that along the proof, any smallness condition
depending on $s_3$ will concern $\nu$ rather than $\e $ (while
$ \e $ is small depending on $ s_1 $ and $ s_2 $). Since $\nu \in (0, \e^{\frac32})$, 
these conditions obviously will be satisfied for $\e$ small enough (depending on $s_3$). 
However to obtain $C^\infty$ solutions in Section \ref{Cinfty}, we shall need to replace $s_3$ 
with some $s$ that is slowly increasing along the Nash Moser scheme, and we shall use that, 
for $\e$ fixed, the estimates of Proposition \ref{prop:inv-ap-vero} hold for any $\nu$ small
enough (depending on $s$).  
\end{remark}

\begin{pf}
The proposition is a consequence of Lemma \ref{thm:Lin+FBR}, Proposition \ref{prop:op-averaged}, Lemma
\ref{A-0-splitted} and Proposition \ref{prop-cruciale}. 
By \eqref{Lomega def}, \eqref{lin:normal-f-ge} and recalling that $ \om = (1 + \e^2 \l ) \bar \om_\e $ we write the operator ${\cal L}_\omega $ acting on $ H_{\mathbb S}^\bot $, 
$$
{\cal L}_\omega = (1 + \e^2 \l ) \Big[  \bar \om_\e \cdot \partial_\vphi   - J  \Big( 
{\mathtt A } + \frac{{\mathtt r}_\e}{1 + \e^2 \l}  \Big)  \Big]  \quad {\rm where} \quad 
{\mathtt A}  = \frac{D_V}{1+ \e^2 \l} +  
\frac{\e^2 {\mathtt B}}{1+ \e^2 \l} \, ,   
$$
as in \eqref{def:primo-op}. Notice that 
  the term $ {\mathtt r}_\e = {\mathtt r}_\e (\uF)$ 
 depends on the torus 
 $ \ui $ at which we linearize, see \eqref{torus-linearized}, 
 unlike $ {\mathtt B} $ in \eqref{form-of-B-ge}, and thus $ {\mathtt A} $.  
 Moreover, by 
 \eqref{estimate:rep-ge0} and  the assumptions of the proposition
 $$ 
 \| \uF \|_{\Lip, s_2+2} \leq \e \, , \quad 
 \| \uF \|_{\Lip, s_3+4} \leq \e  \nu^{-\frac9{10}} \, , 
 $$ 
 we get
 \begin{align}
& | {\mathtt r}_\e |_{\Lip,+,s_2} \leq C(s_2)  \e^2 (\e^2 + \e)  \ll \e^{\frac52}   \label{estrep} \\
&  | {\mathtt r}_\e |_{\Lip,+,s_3+2} \leq   C(s_3) \e^2 (\e^2 + \e \nu^{-\frac9{10}})  \ll \e^{\frac52} \nu^{-1} \, . \label{estrep1}
  \end{align} 
Applying Proposition \ref{prop:op-averaged} we get 
\be\label{coniugazione-yes}
{\cal L}_\omega =  (1 + \e^2 \l )
 {\mathtt P} (\vphi)  \big[ \ppavphi - J ({\mathtt A}_0 + \varrho^+ )  \big] {\mathtt P}^{-1} (\vphi) \,  ,
\ee
 where, as in \eqref{newA+}, \eqref{def:varrho0},  
 \be\label{def:matA0}
   {\mathtt A}_0 :=  \frac{D_V}{1+ \e^2 \l}  + R_0 \, , \quad R_0 = \Pi_{\mathtt D} \varrho \, ,
 \quad  \varrho = \frac{\e^2}{1 + \e^2 \l} \mathtt B  \, , 
\ee
and  the coupling term  $ \varrho^+ $ satisfies \eqref{est:rho+}, i.e. 
\be\label{utile-ripet}
| \varrho^+ |_{\Lip, +,s}  \leq C(s) ( \e^4 + | {\mathtt r}_\e |_{\Lip, +,s} ) \, . 
\ee
By \eqref{estimate-deco-resto} and \eqref{est:varrho}    we have
\be\label{stimaR0n} 
|R_0|_{\Lip , + , s}  = |\Pi_{\mathtt D} \varrho |_{\Lip , + , s}  \leq C(s) \e^2  \, .
\ee
We want to apply Proposition \ref{prop-cruciale} to the operator
\be\label{applic-P}
 \bar \om_\e \cdot \pa_{\vphi} - J ({\mathtt A}_0 + \varrho^+ )  \qquad 
 {\rm with} \qquad   \bg =  \frac{{\mathtt P}^{-1} (\vphi)  g}{1 + \e^2 \lambda } \, , 
\ee 
namely with the substitutions $ A_0 \leadsto {\mathtt A}_0 $, $ R_0 = \Pi_{\mathtt D} \varrho $, $ \rho \leadsto \varrho^+ $ and 
$ g \leadsto \bg $. 
First of all Lemma \ref{A-0-splitted} proves that 
$ {\mathtt A}_0  $ is a split admissible operator in the class $ {\cal C} (C_1,c_1,c_2) $. 
Then, by \eqref{utile-ripet}, \eqref{estimate:rep-ge}, \eqref{stimaR0n},   \eqref{estrep}-\eqref{estrep1},  
we verify that
\be\label{rho-R0:small-bis}
\begin{aligned}
& |  \varrho^+  |_{\Lip, +, s_1} 
\leq  \e^3 \, ,  \\
& |  R_0  |_{\Lip, +, s_2} + |  \varrho^+ |_{\Lip, +, s_2}  
  \leq  \e^{-1} \, ,   \\
& \nu \big( |  R_0  |_{\Lip, +, s_3+2}  + | \varrho^+ |_{\Lip, +, s_3+2}  \big) \ll \e^2  \, , 
\end{aligned}
\ee
for $ \nu$ small depending on $ s_3 $ (for the last estimate we 
use the first inequality in \eqref{estrep1}),  
which implies the conditions required in \eqref{rho-R0:small} and \eqref{g:small-large} with $ s_3 \leadsto s_3 + 2 $. 
Notice also that the function  $ \bg  $ defined in \eqref{applic-P} satisfies, 
by \eqref{tame-P0} and \eqref{opernorm-Lip}, the estimate
$$ 
\| \bg \|_{\Lip,s} \leq C(s) \| g \|_{\Lip,s} \, , \quad \forall s \geq s_1 \, , 
$$ 
and therefore, by \eqref{g:small-large-IA}, the function $ \bg $ satisfies the assumption
$$ 
\| \bg \|_{\Lip,s_1} \leq \e^2 C(s_1)  \nu  \, , \quad  \| \bg \|_{\Lip,s_3+2} \leq  C(s_3) \e^2 \nu^{-\frac9{10}} \leq  \e^2 ( C(s_1)  \nu)^{-1} 
$$ 
required in \eqref{g:small-large} with $ \nu \leadsto C(s_1) \nu $ and
$ s_3 \leadsto s_3 + 2 $, for $ \nu $ small depending on $ s_3 $. 
Note also that, since $ \nu \in (0, \e^{3/2})$, we have $ C(s_1) \nu \in (0,\e)$, 
as required in Proposition \ref{prop-cruciale}.
Therefore Proposition \ref{prop-cruciale} applies and there  are closed subsets (independent of $g$)
\be \label{def:Canti}
{\bf \Lambda} (\e; \cc , \uF  ) := {\bf \Lambda} (\e; \cc , {\mathtt A}_0, \varrho^+) \subset \Lambda_\uF \, , 
\quad 1/2 \leq \cc \leq 5/6 \, , 
\ee
satisfying the properties 1-3 listed in Proposition \ref{prop-cruciale} with $ \wtilde \Lambda = \Lambda_\uF $, 
and an operator $ {\mathfrak L}_{approx}^{-1} = {\mathfrak L}_{approx, \nu}^{-1}  $ such that the  function  
$$ 
\bh  := {\mathfrak L}_{approx}^{-1} \breve g \, , \  \breve h  :  {\bf \Lambda} (\e; 5/6 , \uF  ) 
 \to  {\mathcal H}^{s_3+2} \cap H_{\mathbb S}^\bot \, , 
 $$   
  satisfies \eqref{h:s1s2s3} with $ s_3 \leadsto s_3 + 2 $, i.e. 
\be\label{h:s1s2s3-here}
 \| \bh \|_{\Lip, s_1} \lesssim_{s_1} \e^2 \nu^{ \frac45} \, , 
 \quad   \| \bh \|_{\Lip, s_3+2} \lesssim_{s_1} \e^2 \nu^{ -  \frac{11}{10} } \, , 
\ee
and, 
see \eqref{sol:almost-inv}, \eqref{estimate:r-final}, 
\be\label{sol:almost-inv-here}
(\bar {\om }_\e \cdot \partial_\varphi - J ({\mathtt A}_0 +  \varrho^+)) \bh =  
\bg  +  \br \quad {\rm with} \quad  \| \br\|_{\Lip, s_1} \lesssim_{s_1} \e^2 \nu^{\frac32} \, .  
\ee
Set  
\be\label{defined:L-1}
{\cal L}_{appr}^{-1} :=   \frac{1}{1 + \e^2 \lambda} \, {\mathtt P} (\vphi)  {\mathfrak L}_{approx}^{-1}   {\mathtt P}^{-1} (\vphi)  
 \,, \quad 
h := {\mathtt P} (\vphi) \breve h =  {\cal L}_{appr}^{-1}  g   \, . 
\ee 
By  \eqref{sol:almost-inv-here}-\eqref{defined:L-1} and 
\eqref{coniugazione-yes}, \eqref{applic-P}  we get 
\be\label{ciochere} 
{\cal L}_\omega h - g  = (1 + \e^2 \l)  {\mathtt P} (\vphi) \br \, . 
\ee
By \eqref{defined:L-1}, \eqref{tame-P0}, \eqref{opernorm-Lip}, \eqref{ciochere} we get 
\begin{align*}
& \| h \|_{\Lip,s}  =  \|   {\mathtt P} (\vphi) \breve h \|_{\Lip,s}  \leq C(s) \| \bh \|_{\Lip,s}  \\
&  \| {\cal L}_\omega h - g \|_{\Lip,s} 
= \| (1 + \e^2 \l)  {\mathtt P} (\vphi) \br \|_{\Lip,s}  \leq C(s) \|  \br  \|_{\Lip,s} \, , 
\end{align*}
and, from the estimates of $ \breve h $ and $ \breve r $
in \eqref{h:s1s2s3-here}, \eqref{sol:almost-inv-here}, we deduce \eqref{h:s1s2s3-IA-per}, \eqref{estimate:r-final-IA-per}. 
Furthermore the estimate \eqref{new:estimate-cal-L-1} follows by 
\eqref{defined:L-1}, \eqref{tame-P0}, \eqref{opernorm-Lip} and \eqref{con-loss}. 

Let us finally prove  that  the sets 
$ {\bf \Lambda} (\e; \cc , \uF  ) $ defined in \eqref{def:Canti}
satisfy the properties  1-3 listed in the statement of the Proposition  \ref{prop:inv-ap-vero}.
Items 1-2 are immediate consequences the corresponding ones in Proposition   \ref{prop-cruciale}.  For proving 
item 3, first notice  that, $ \varrho $ (defined in \eqref{def:varrho0}), 
$ R_0 = \Pi_{\mathtt D} \varrho $ and hence $ {\mathtt A}_0 $ (defined in \eqref{def:matA0}) do not depend on the
torus $ \uF $ at which we linearize.  
Then,  if $ \| \uF - \uF' \|_{s_1+2} \leq \delta \leq \e^{3/2} $,  we have  
$$ 
\begin{aligned} 
| {\mathtt A}'_0 - {\mathtt A}_0|_{+,s_1} + | \varrho^+ - (\varrho^+)' |_{+,s_1}
& = | \varrho^+ - (\varrho^+)' |_{+,s_1} \\ 
& \stackrel{\eqref{est:rho+-variation}}  {\lesssim_{s_1}}  | {\mathtt r}_\e -{\mathtt r}'_\e|_{+,s_1} \\
& \stackrel{\eqref{estimate:rep-ge0-diff}} {\lesssim_{s_1}}    \e^2 \| \uF - \uF' \|_{s_1+2}  \leq  \e^{\frac32} \delta \leq \e^3 \, ,
\end{aligned}
$$
so that condition \eqref{smallR0R0'} of Proposition  \ref{prop-cruciale} is satisfied.
Hence, by the property 3 (see \eqref{finadav}) in Proposition  \ref{prop-cruciale},
and the inclusion 
$
{\bf \Lambda} (\e; \cc - \d^{\frac25} , \uF  ) \subset 
{\bf \Lambda} (\e; \cc - (\e^{\frac32} \d)^{\frac25} , \uF  ) \, 
$, we have 
$$
\big|  \Lambda_{\uF'} \cap [{\bf \Lambda} (\e; \cc , \uF ' )  ]^c \cap 
{\bf \Lambda} (\e; \cc - \d^{\frac25} , \uF  )       \big| \leq (\e^{\frac32} \d)^{\frac{\a}{3}} \, .
$$
This proves property 3 for the sets ${\bf \Lambda} (\e; \cc , \uF )$ in Proposition 
\ref{prop:inv-ap-vero}.
\end{pf}

\section{Nash-Moser iteration}\label{sec:NM-ite}

In Chapter \ref{sec:thmNM} we observed that the nonlinear operator $ {\cal F} $ defined in \eqref{operatorF}
evaluated at the trivial torus $ i_0 (\vphi) := (\vphi,0,0)  $ satisfies (see \eqref{0sol})
\be\label{Fi0s1}
 \| {\cal F}(\l;  i_0 ) \|_{\Lip,s} \leq C(s) \e^2 \quad {\rm on } \ \Lambda \, , \quad \forall s \geq s_0 \, , 
\ee
 and, in Lemma \ref{NM:step1} we have defined,  for all 
 $ \l \in \Lambda $, a torus embedding $ i_1 (\vphi) $ such that  (see \eqref{approx1})
\be\label{eq:F1e4}
\| {\cal F}(\l; i_1 ) \|_{\Lip,s} \leq C(s) \e^4 \quad {\rm on } \ \Lambda \, ,  \quad \forall s \geq s_0 \, . 
\ee
In this section we define a sequence of torus embeddings $ (i_n ( \vphi) )_{n\geq 2}$, $ n \geq 2 $, 
defined for $ \l $ belonging to a decreasing sequence of 
subsets $ {\bf \Lambda}_n \subset \Lambda $ which converges, for all $ \l $ belonging to the 
intersection $ \cap_{n \geq 2} {\bf \Lambda}_n $ to a solution $ i_\infty (\vphi) $  of $ {\cal F} (\l; i_\infty ) = 0 $.

Fixing $ \nu_1 = C(s_1) \e^2 $ such that $\| {\cal F} (i_1)\|_{\Lip,s_1} \leq \e^2 \nu_1 $ (see \eqref{eq:F1e4}), 
we define the decreasing sequence
$ (\nu_n)_{n \geq 1} $  by 
\be\label{def:nun-NM}
\nu_n :=  \nu_1^{ {\mathfrak q}^{n-1} }\, , \qquad \nu_{n+1} = \nu_n^{\mathfrak q} \, , \qquad 
\mathfrak q := \frac{3}{2} - \s_\star  \, ,
\ee
for some $ \s_\star  \in ( 0, 1/4) $. 
Theorem \ref{thm:NM}\index{Nash-Moser theorem} will be a consequence of the following result. 

\begin{theorem}\label{thm:NM-ite} {\bf (Nash-Moser)}
Let  $ \bar \om_\e \in \R^\es $ be  $(\gamma_1 , \tau_1) $-Diophantine and satisfy 
property $ {\bf (NR)}_{\gamma_1, \tau_1} $ in Definition \ref{NRgamtau} with $ \g_1, \t_1 $ fixed in 
\eqref{def:tau1}. Assume  \eqref{choice s2 s3} and 
$ s_2 - s_1 \geq \underline{\tau} + 2 $, $ s_1 \geq s_0+2+\underline{\tau} +Q'  $, 
where $ \underline{\tau} $ is the loss of derivatives  defined in Proposition
\ref{Prop:sec6} and $ Q' := 2(\tau'+ \varsigma s_1 ) +3 $ is defined in  
Proposition \ref{lem:app-inv-LD}.

Then, there exists $ s_3 $ large enough, a constant
$ \sigma_\star := \sigma_\star (s_3) >  0 $, satisfying $ \sigma_\star (s_3) \to 0 $ as $ s_3 \to + \infty $, and 
 $ \e_0 > 0 $, such that, defining the sequence $(\nu_n)$ by \eqref{def:nun-NM},
 for all $ 0 < \e \leq \e_0 $,  for all  $ n \geq 1 $, there exist 
\begin{enumerate}
\item a subset $ {\bf \Lambda}_n \subseteq {\bf \Lambda}_{n-1}  $, $ {\bf \Lambda}_1 := {\bf \Lambda}_0 := \Lambda $,
satisfying 
\be\label{meas:Ln-NM}
\begin{aligned}
& 
| {\bf \Lambda}_1 \setminus {\bf \Lambda}_{2} | \leq b(\e) \quad {\rm with} \quad \lim_{\e \to 0} b(\e) = 0 \, , \\ 
&  | {\bf \Lambda}_{n-1} \setminus {\bf \Lambda}_{n}| \leq \nu_{n-2}^{ \alpha_*} \, , \ \forall n \geq 3 \, ,
\end{aligned} 
\ee
where $ \a_* =  \alpha / 4 $ and $ \a > 0  $ is the exponent in \eqref{exp-misu}, 
\item  a torus $ i_n (\vphi) = (\vphi,0,0) + \fracchi_n (\vphi) $,
 defined for all $ \l \in {\bf \Lambda}_n  $, 
satisfying
\begin{align}
& 
\| \fracchi_n  \|_{\Lip,s_1+2} \leq 
 C(s_1) \e^{2}  \, , \quad  
\| \fracchi_n  \|_{\Lip,s_2+2} \leq \e  \, , \quad 
\| \fracchi_n  \|_{\Lip,s_3+2} \leq 
\e^2 \nu_n^{- \frac45} \, , \label{NM:indu-line2}  \\
&  
\| i_n - i_{n-1}  \|_{\Lip,s_1+2} \leq C(s_1) \nu_{n-1} + \e^2 \nu_{n-1}^{\frac{4}{5}-\sigma_\star} \, , \ n \geq 2 \, ,  \label{vicin-in+1-in} \, \\
&  
\| i_n - i_{n-1}  \|_{\Lip,s_2+2} \leq  \e \nu_{n-1}^{\frac{1}{5}} \, , \quad n \geq 2  \label{vicin-in+1-in-s2} \, , 
\end{align}
such that 
\be\label{NM:indu} 
 \| {\cal F}(i_n )\|_{\Lip,s_1} \leq \e^2 \nu_n  \,  .   
\ee
\end{enumerate}
\end{theorem}

The rest of this section is devoted to the proof of Theorem \ref{thm:NM-ite}. 

\medskip

\noindent
{\bf First step.} For $ n = 1 $ 
the torus $ i_1 (\vphi ) $ defined in Lemma \ref{NM:step1}, for all $\lambda \in \Lambda = {\bf \Lambda}_1 $, 
satisfies  \eqref{estimate-i1} and \eqref{approx1}, which imply
\eqref{NM:indu-line2} and  \eqref{NM:indu}  at $ n = 1 $,  for $ \e $ small enough (recall 
that $ \nu_1 = C(s_1) \e^2 $).  
\\[1mm]
\noindent
{\bf Iteration.} We now proceed by induction. 
Assume that we have already defined subsets $ {\bf \Lambda}_n \subseteq 
 {\bf \Lambda}_{n-1} \subseteq  \ldots \subseteq \Lambda $ satisfying \eqref{meas:Ln-NM},
and $ n $ tori  $ i_1, \ldots, i_{n} $, of the form  
$ i_n (\vphi) = (\vphi,0,0) + \fracchi_n (\vphi) $, satisfying 
\eqref{NM:indu-line2}-\eqref{vicin-in+1-in-s2}  and  such that \eqref{NM:indu} holds. 
We are going to define a subset  $  {\bf \Lambda}_{n+1} \subseteq {\bf \Lambda}_n  $ and,
for all  $ \l \in {\bf \Lambda}_{n+1} $, 
the subsequent better approximate torus embedding $ i_{n+1} $ which 
 satisfies  \eqref{NM:indu-line2}-\eqref{vicin-in+1-in-s2} and  \eqref{NM:indu} 
 at order $ n + 1 $. 

We shall define $ i_{n+1} $ by a Nash-Moser type iterative scheme, using 
Propositions \ref{Prop:sec6} and \ref{prop:inv-ap-vero}.  
\begin{itemize}
\item
{\bf Notation.} 
From now, $\sigma$ denotes an arbitrarily small, strictly positive, constant. 
In the sequel, when it appears in some inequality  (at the exponent), it means
that for any $\s>0$, if $s_3$ has been chosen large enough and $\nu_n$ is small enough  (depending on $s_3$ and $\s$), the inequality holds.  In particular, since $\nu_n \leq \nu_1 \leq C (s_1) \e^2$, the inequality holds for
$\e$  small enough (depending on $s_3$). 
\end{itemize}

\noindent
{\bf Step 1. Regularization.}
We first consider the regularized approximate torus
\be\label{trunc-torus}
\breve \imath_n (\vphi) := (\vphi,0,0) + \breve \fracchi_n  (\vphi)\, , \quad \breve \fracchi_n := {\it \Pi}_{N_n} \fracchi_n \, , 
\ee
where $ {\it \Pi}_N $ is the Fourier projector defined in \eqref{def:projectorN} and 
\be\label{def:Nn}
N_n \in \Big[ \nu_n^{-\frac{3}{s_3-s_1}} -1, \nu_n^{-\frac{3}{s_3-s_1}} +1 \Big] \, . 
\ee
\begin{lemma} \label{1122}
The torus $ \breve \imath_n (\vphi) := (\vphi,0,0) + \breve \fracchi_n (\vphi) $ satisfies 
\begin{align}
& \| \breve \fracchi_n  \|_{\Lip,s_1+2} \leq C(s_1) \e^2 \, , 
\quad  \| \breve \fracchi_n  \|_{\Lip,s_2+2} \leq \e\, ,  \quad
\| \breve \fracchi_n  \|_{\Lip,s_3+ \bt + 6} \leq \e^2 \nu_n^{- \frac{4}{5} -  \sigma} \, , \label{NM:indu-line2-trunc} \\
& \| i_n - \breve \imath_n \|_{\Lip, s_1+2} \ll \e^2 \nu_n^2 \, , \quad 
 \| i_n - \breve \imath_n \|_{\Lip, s_3+2} \leq \e^2 \nu_n^{- \frac45} \, , \label{NM:int2} \\
& \| {\cal F}( \breve \imath_n )\|_{\Lip,s_1} \leq 2 \e^2 \nu_n  \,  , 
\quad \| {\cal F}( \breve \imath_n )\|_{\Lip,s_3+2+\bt } \leq \e^2 \nu_n^{-  \frac{4}{5} - 2 \sigma} \, , 
\label{NM:indu-trunc} 
\end{align}
where $ \sigma > 0 $ can be taken arbitrarily small taking $ s_3 - s_1 $ large enough.
\end{lemma}

\begin{pf}
The first two estimates in \eqref{NM:indu-line2-trunc} follow by \eqref{NM:indu-line2} and\index{Smoothing operators}  \eqref{smoothingS1S2-Lip}. 
The third estimate in \eqref{NM:indu-line2-trunc} follows  by
$$
\begin{aligned}
\| \breve \fracchi_n  \|_{\Lip,s_3+ \bt + 6}  
& \stackrel{\eqref{smoothingS1S2-Lip}} \leq N_n^{\bt + 4} \|  \fracchi_n  \|_{\Lip,s_3+2} \\
& \stackrel{\eqref{def:Nn}, \eqref{NM:indu-line2}} \leq 
C \nu_n^{- \frac{3(\bt + 4)}{s_3-s_1}} \e^2 \nu_n^{- \frac{4}{5}}  \leq \e^2 \nu_n^{- \frac{4}{5} -  \sigma}  
\end{aligned}
$$
for $ 3 ( \underline{\tau} + 4)/(s_3 - s_1) < \sigma $.
\\[1mm]
{\sc Proof of \eqref{NM:int2}.} 
We have
$$
\begin{aligned}
\| i_n - \breve \imath_n \|_{\Lip,s_1 +2} \stackrel{ \eqref{trunc-torus}} 
= \| {\it \Pi}_{N_n}^\bot \fracchi_n  \|_{\Lip,s_1+2} & \stackrel{\eqref{smoothingS1S2-Lip}} 
\leq N_n^{-(s_3-s_1)} \| \fracchi_n  \|_{\Lip,s_3+2}  \\
& \stackrel{\eqref{def:Nn}, \eqref{NM:indu-line2}} 
{ \lesssim_{s_3} } \nu_n^{3} \e^2 \nu_n^{- \frac45}
\end{aligned}
$$
which implies the first bound in \eqref{NM:int2}.  Similarly 
\eqref{trunc-torus}, \eqref{smoothingS1S2-Lip} and  \eqref{NM:indu-line2} imply the second estimate in \eqref{NM:int2}. 
\\[1mm]
{\sc Proof of \eqref{NM:indu-trunc}.} Recalling  the definition of the operator $ {\cal F } $
in \eqref{operatorF}, 
and using Lemma\index{Moser estimates for composition operator} \ref{Moser norme pesate}, 
\eqref{NM:indu-line2}, \eqref{NM:indu-line2-trunc}, we have
\begin{align*}
\| {\cal F}( \breve \imath_n ) \|_{\Lip,s_1} 
& \leq \| {\cal F}( i_n ) \|_{\Lip,s_1} + C (s_1) \| i_n - \breve \imath_n \|_{\Lip, s_1+2} \\
& \stackrel{\eqref{NM:indu}, \eqref{NM:int2}}  \leq    2\e^2 \nu_n 
\end{align*}
proving the first inequality in \eqref{NM:indu-trunc}. 
Finally Lemma \ref{Moser norme pesate} 
and \eqref{NM:indu-line2-trunc} imply 
\begin{align*}
\| {\cal F}( \breve \imath_n ) \|_{\Lip,s_3+2+ \bt} 
& \leq \| {\cal F}( i_0 ) \|_{\Lip,s_3+2+ \bt} +
C(s_3) \| \breve \fracchi_n \|_{\Lip, s_3+4+\bt} \\
& \stackrel{\eqref{Fi0s1}, \eqref{NM:indu-line2-trunc}} {\lesssim_{s_3}} \e^2  + 
 \e^2 \nu_n^{- \frac{4}{5} -  \sigma}  
\leq  \e^2 \nu_n^{- \frac{4}{5} - 2 \sigma} 
\end{align*}
for $ \e $ small, proving the second inequality in \eqref{NM:indu-trunc}. 
\end{pf}

\medskip
\noindent
{\bf Step 2. Isotropic torus.} 
We associate to the regularized torus $ \breve \imath_n  $  defined in \eqref{trunc-torus}
 the isotropic torus\index{Isotropic torus}  $ i_{n, \delta}   $ defined by Proposition \ref{Prop:sec6} with $ \ui = \breve \imath_n  $. Notice that,  by \eqref{NM:indu-line2-trunc}, the torus  
 $ \breve \imath_n $ satisfies the assumption 
 \eqref{ansatz 0} required in Proposition \ref{Prop:sec6}.
 The torus  $ i_{n, \delta}  $ is an approximate solution ``essentially as good" as $ i_n $
(compare \eqref{NM:indu} and \eqref{NM:indu-iso}).  

\begin{lemma}\label{lem:ison}
The isotropic torus $ i_{n, \delta } (\vphi) = (\vphi, 0,0) +  \fracchi_{n,\d} (\vphi) $ defined by 
Proposition \ref{Prop:sec6} with $ \ui = \breve \imath_n  $,  satisfies
\begin{align}\label{eq:diff-con-isot}
& \| i_{n,\d} - \breve \imath_n \|_{\Lip,s_0+2}  \lesssim_{s_0}  \e^2 \nu_n \, , \qquad 
\| i_{n,\d} - \breve \imath_n \|_{\Lip,s_1+4}  \ll  \e^2 \nu_n^{1-  \sigma} \, , 
 \\
& \label{bound-int-iso}
\| i_{n,\d} - \breve \imath_n   \|_{\Lip, s_3+4}  \ll \e^2 \nu_n^{ - \frac{4}{5} -  2 \sigma} \, ,
\quad \| \fracchi_{n,\d}  \|_{\Lip, s_3+4}  \ll \e^2 \nu_n^{ - \frac{4}{5} - 2 \sigma} \, , 
\end{align}
and 
\be
\begin{aligned}\label{NM:indu-iso}
& \| {\cal F}(i_{n, \d} )\|_{\Lip,s_0}    \lesssim_{s_0} \e^2 \nu_n \, , \\
& \| {\cal F}(i_{n, \d} )\|_{\Lip,s_1+2}   \leq  \e^2 \nu_n^{1-  \sigma}   \\
& \| {\cal F}(i_{n, \d} )\|_{\Lip,s_3+2}   \leq \e^2 \nu_n^{-  \frac{4}{5} - 2 \sigma} \, .
\end{aligned}
\ee
\end{lemma}

\begin{pf}

\noindent
{\sc Proof of \eqref{eq:diff-con-isot}.}
Since $ s_0 + 2 +  \bt \leq s_1 $, we have by \eqref{NM:indu-line2-trunc},
$$
\| \breve \fracchi_n \|_{\Lip , s_0+2+\bt} \leq C(s_1) \e^2 \, .
$$
Hence by
\eqref{stima y - y delta} (with $ s = s_0 +2$), we have  
 \be \label{eq:diff-iso0}
 \| i_{n,\d} - \breve \imath_n \|_{\Lip,s_0+2} \lesssim_{s_0} \|  {\cal F}(\breve \imath_n ) \|_{\Lip,s_0 + 2 + \bt}
 \lesssim_{s_0} \|  {\cal F}(\breve \imath_n ) \|_{\Lip,s_1}  
 \stackrel{\eqref{NM:indu-trunc}}{\lesssim_{s_0}} \e^2 \nu_n 
 \ee
 which is the first estimate in \eqref{eq:diff-con-isot}. 
Similarly, since $s_1+2+\bt \leq s_2$, we have by the second estimate in 
\eqref{NM:indu-line2-trunc},
$$
\| \breve \fracchi_n  \|_{\Lip , s_1+4+\bt} \leq \| \breve \fracchi_n \|_{\Lip , s_2+2} \leq \e \, . 
$$
Hence  by \eqref{stima y - y delta} (with $s=s_1+4$),  
\be\label{eq:diff-iso}
\| i_{n,\d} - \breve \imath_n \|_{\Lip,s_1+4} \lesssim_{s_1} \|  {\cal F}(\breve \imath_n ) \|_{\Lip,s_1 + 4 + \bt} \, .
\ee
Now, by the interpolation inequality \eqref{inter-Sobo-Lip} 
we have 
\be\label{15-6-2017}
\|  {\cal F}(\breve \imath_n ) \|_{\Lip,s_1 + 4 + \bt} 
 \lesssim_{s_3} \|  {\cal F}(\breve \imath_n ) \|_{\Lip,s_1}^{\theta} \|  {\cal F}(\breve \imath_n ) \|_{\Lip,s_3}^{1-\theta} 
\ee
where 
\be\label{def:theta-int}
\theta := 1  - \frac{ 4 + \bt  }{s_3 - s_1 } \, , \quad 1 - \theta := \frac{4 + \bt }{s_3 - s_1  } \, . 
\ee
Therefore \eqref{15-6-2017}, \eqref{NM:indu-trunc}, \eqref{def:theta-int} imply 
\be\label{15-6-2017-bis}
\|  {\cal F}(\breve \imath_n ) \|_{\Lip,s_1 + 4 + \bt} \lesssim_{s_3} 
(\e^2 \nu_n)^{\theta} (\e^2 \nu_n^{-  \frac{4}{5} -  2 \sigma} )^{1-\theta} 
\ll  \e^2 \nu_n^{1-  \sigma}  
\ee
for $ s_3 - s_1 $ large enough and $ \nu_n $ small.
Inequalities  \eqref{eq:diff-iso} and \eqref{15-6-2017-bis} prove the second estimate in \eqref{eq:diff-con-isot}.
\\[1mm]
{\sc Proof of \eqref{bound-int-iso}.}
By
\eqref{2015-2} we have 
$$
\| \fracchi_{n,\d} - \breve \fracchi_n \|_{\Lip, s_3+4}  \lesssim_{s_3}
\|  \breve \fracchi_n \|_{\Lip, s_3+5}   
\stackrel{\eqref{NM:indu-line2-trunc} }\ll \e^2 \nu_n^{ - \frac{4}{5} -  2 \sigma} 
$$
for $ \nu_n $ small enough,  
proving the first inequality in \eqref{bound-int-iso}. The second  inequality in \eqref{bound-int-iso}
follows as a consequence of the first one and the third estimate in \eqref{NM:indu-line2-trunc}. 
\\[1mm]
{\sc Proof of \eqref{NM:indu-iso}.}
Recalling the definition of the operator $ {\cal F } $ in \eqref{operatorF}
 and using Lemma\index{Moser estimates for composition operator} \ref{Moser norme pesate}, 
 \eqref{NM:indu-line2-trunc}, we have 
$$
\| {\cal F} ( i_{n,\d}) \|_{\Lip, s} \leq  \| {\cal F} ( \breve \imath_{n}) \|_{\Lip, s} + 
C(s) \|  i_{n,\d} -  \breve \imath_{n} \|_{\Lip, s+2}
$$
both for $ s = s_0 $ and $ s = s_1 + 2 $, and therefore
$$
\begin{aligned}
& \| {\cal F} ( i_{n,\d}) \|_{\Lip, s_0}  
 \stackrel{\eqref{NM:indu-trunc} \eqref{eq:diff-con-isot}} {\lesssim_{s_0}}
\e^2 \nu_n \, , \\
& \| {\cal F} ( i_{n,\d}) \|_{\Lip, s_1 +2}   \stackrel{\eqref{15-6-2017-bis} \eqref{eq:diff-con-isot}} {\leq }
\e^2 \nu_n^{1-\s}  \, , 
\end{aligned}
$$
proving the first two estimates in \eqref{NM:indu-iso}.
Finally
\begin{align}
 \| {\cal F}(i_{n, \d} )\|_{\Lip,s_3+2} 
& \leq \| {\cal F}(i_{0} )\|_{\Lip,s_3+2} + \| {\cal F}(i_{n, \d} ) - {\cal F}(i_{0} ) \|_{\Lip,s_3+2} \nonumber  \\
& \stackrel{\eqref{Fi0s1}, Lemma \, \ref{Moser norme pesate}} {\lesssim_{s_3}} \e^2 + \| \fracchi_{n,\d} \|_{\Lip, s_3+4}  \nonumber	\\
& \stackrel{\eqref{bound-int-iso}} \leq  \e^2 \nu_n^{-  \frac{4}{5} - 2 \sigma}
\end{align}
for $ \nu_n $ small,  which is the third inequality in \eqref{NM:indu-iso}. 
\end{pf}

\medskip

\noindent
{\bf Step 3. Symplectic diffeomorphism.} We apply the symplectic change of variables 
$ G_{n,\delta} $ defined in \eqref{trasformazione modificata simplettica} with $ (\uth (\phi), y_\d (\phi), \uz (\phi) ) 
=  i_{n,\delta} (\phi) $, 
which transforms the isotropic torus $ i_{n,\delta} $ into (see \eqref{toro-new-coordinates})
\be\label{torus:origin}
G_{n,\delta}^{-1}(i_{n,\delta} (\vphi) )= (\ph, 0 , 0 ) \, . 
\ee 
It conjugates the Hamiltonian vector field $ X_{K} $ associated to $  K $  defined in \eqref{primalinea},  with the Hamiltonian vector field
(see \eqref{new-Hamilt-K})
\be\label{new-Hamilt-Kn}
X_{{\mathtt K}_n} = (D G_{n, \d})^{-1} X_{K} \circ G_{n, \d} 
\qquad {\rm where} \qquad {\mathtt K}_n := K \circ G_{n, \d}  \, .
\ee
Denote by $  {\mathtt u} := (\phi, \ac, w) $  
the symplectic coordinates induced by the diffeomorphism $ G_{n, \d} $ in \eqref{trasformazione modificata simplettica}.  
Under the symplectic map $ G_{n,\delta} $, the nonlinear operator ${\cal F}$ in \eqref{operatorF} is transformed into 
\be \label{trasfo imp}
{\mathtt F}_n ({\mathtt u} (\vphi) ) := 
 \Dom {\mathtt u} (\vphi) - X_{{\mathtt K}_n} ( {\mathtt u} (\vphi))   
 =  (D G_{n, \delta}( {\mathtt u}  (\vphi) ))^{-1} {\cal F}(G_{n, \delta}(  {\mathtt u} (\vphi) ) ) \, .
\ee
By \eqref{trasfo imp} and \eqref{torus:origin} (see also \eqref{apprVF}) we have that 
\be\label{mathtt-Fn-new}
{\mathtt F}_n ( \ph, 0 , 0  ) = 
(\om,0,0) - X_{{\mathtt K}_n} ( \vphi,0,0) =   
 (D G_{n, \delta}( \vphi, 0, 0 ))^{-1} {\cal F} \big( i_{n, \delta} (\vphi) \big) \, .
\ee
Hence by \eqref{DG delta}, for $s\geq s_0$,
\begin{align} 
\| {\mathtt F}_n ( \ph, 0 , 0 )\|_{\Lip,s}  & \lesssim_s   \| {\cal F} ( i_{n, \delta} ) \|_{\Lip , s}
+ \|  \breve \fracchi_n \|_{\Lip, s +2}   \| {\cal F} ( i_{n, \delta} ) \|_{\Lip , s_0} \nonumber \\
 &  \stackrel{\eqref{NM:indu-iso}}{\lesssim_s}    \| {\cal F} ( i_{n, \delta} ) \|_{\Lip , s}
+ \|  \breve \fracchi_n \|_{\Lip, s +2}  \, \e^2 \nu_n   \, . \label{estFn1}
\end{align}
By \eqref{estFn1}, \eqref{NM:indu-iso} and  \eqref{NM:indu-line2-trunc}, we have 
\be\label{NM:indu-new}
\begin{aligned}
& \| {\mathtt F}_n ( \ph, 0 , 0 )\|_{\Lip,s_0} \lesssim_{s_0} \e^2 \nu_n   \\
& \| {\mathtt F}_n ( \ph, 0 , 0 )\|_{\Lip,s_1+2} \leq  \e^2 \nu_n^{1-  2 \sigma}  \\
& \| {\mathtt F}_n ( \ph, 0 , 0 ) \|_{\Lip,s_3+2} \leq  \e^2 \nu_n^{-  \frac{4}{5} - 3 \sigma} \, .
\end{aligned}
\ee
{\bf Step 4. Approximate solution of the linear equation associated to a Nash-Moser step in new coordinates.}
In order  to  look for a better approximate zero  
\be\label{def:un+1}
 {\mathtt u}_{n+1} (\vphi) = (\vphi,0,0) + {\mathtt h}_{n+1}(\vphi) 
\ee
of the operator $ {\mathtt F}_n  (  {\mathtt u} )  $  defined in \eqref{trasfo imp}, we expand
\be \label{Ham-new-expa}
 {\mathtt F}_n (  {\mathtt u}_{n+1}  ) =  {\mathtt F}_n ( \vphi,0,0  ) + 
 \big( \Dom - d_{\mathtt u} X_{{\mathtt K}_n} (\vphi,0,0 ) \big)  {\mathtt h}_{n+1}  +
 Q_{n} ( {\mathtt h}_{n+1})
\ee
where $Q_{n} ( {\mathtt h}_{n+1})$ is a quadratic  remainder, 
and we want to solve 
(approximately) the linear equation
\be \label{Newtonnewvar}
{\mathtt F}_n ( \vphi,0,0  ) + 
 \big( \Dom - d_{\mathtt u} X_{{\mathtt K}_n} (\vphi,0,0 ) \big)  {\mathtt h}_{n+1} =0 \, .
\ee
This is the goal of the present step. 
Notice that the operator  $ \Dom  - d_{\mathtt u} X_{{\mathtt K}_n}( \vphi,0,0 ) $ is provided by \eqref{lin idelta},
and we decompose it as
\begin{equation}\label{splitting linearizzato nuove coordinate}
\om \! \cdot \! \pa_\vphi 
- d_{\mathtt u} X_{{\mathtt K}_n}( \vphi,0,0 ) 
= \mathbb{D}_n + {\mathtt R}_{Z_n}  
\end{equation}
where $ \mathbb{D}_n $ is  the operator in \eqref{operatore inverso approssimato} with  $ \ui $ replaced by 
$ \breve \imath_n $, i.e. 
\begin{equation}\label{def:operatore inverso approssimato proiettato}
{\mathbb D}_n 
 \begin{pmatrix}
\widehat \phi \\
 \widehat \ac \\
\widehat w
\end{pmatrix}  := \begin{pmatrix}
 \omega \cdot \partial_\vphi \widehat \phi  - 
{\mathtt K}_{20}^{(n)}(\vphi) \widehat \ac  - [{\mathtt K}_{11}^{(n)}]^\top(\vphi) \widehat w\\
\!\!\!\!\!\! \!\!\!\!\!\!\!\!\!\!\!\!\!\!\!\!\!\!\!\!\!\!\!\!\!\!\!\!\!\!\!\!\!\!\!\! \!\!\!\!\!\!\!\!\!\!\!\!\!\!\!\!\!\!
\!\!\!\!\!\! \! 
\omega \cdot \partial_\vphi \widehat \ac  \\
\!\!\!\! \!\!\!\!\!\! \!\!\!\!\!\!\!\!\!\!\!\! \!\!\!\!\!\!\!\!\!\!\!\! \!\!\!
{\cal L}_{\omega}^{(n)} \widehat w -J {\mathtt K}_{11}^{(n)}(\vphi)\widehat \ac  
\end{pmatrix} 
\end{equation}
with $ {\cal L}_{\omega}^{(n)} := {\cal L}_{\omega} ( \breve \imath_n )  $ (see \eqref{Lomega def}), the functions 
$ {\mathtt K}_{20}^{(n)} (\vphi) $, ${\mathtt K}_{11}^{(n)} (\vphi) $ are the Taylor coefficients 
in \eqref{KHG}
of the Hamiltonian $ {\mathtt K}_n $,  and 
\be\label{term:RZ}
 {\mathtt R}_{Z_n} 
  \begin{pmatrix}
\widehat \phi \\
 \widehat \ac \\
\widehat w
\end{pmatrix} := \begin{pmatrix}
 - \partial_\phi {\mathtt K}_{10}^{(n)}(\vphi) [\widehat \phi ] \\

 \partial_{\phi \phi} {\mathtt K}^{(n)}_{00} (\vphi) [ \widehat \phi ] + 
 [\partial_\phi {\mathtt K}^{(n)}_{10}(\vphi)]^\top \widehat \ac + 
 [\partial_\phi {\mathtt K}^{(n)}_{01}(\vphi)]^\top \widehat w  \\
 - J \{ \partial_{\phi} {\mathtt K}^{(n)}_{01}(\vphi)[ \widehat \phi ] \}
 \end{pmatrix}\, . 
\ee
In the next proposition we define $  {\mathtt h}_{n+1} $ as an approximate solution of the linear equation 
$ {\mathtt F}_n ( \vphi,0,0  ) +  \mathbb{D}_n  {\mathtt h}_{n+1} = 0 $.

\begin{proposition}\label{Prop:inversione}{\bf (Approximate right inverse of $ {\mathbb D}_n $)} \index{Approximate right inverse of $ {\mathbb D}_n $}
For all $  \l $ in the set 
\be\label{def:Lambda-n+1} 
{\bf \Lambda}_{n+1} := {\bf \Lambda}_{n} \cap 
{\bf \Lambda} (\e; \eta_n, \breve \fracchi_n ) \, , \quad \ \eta_1 := 1/2 \, ,  \quad 
\eta_n  := \eta_{n-1} + \nu_{n-1}^{\frac14} \, ,  \ n \geq 2 \, ,  
\ee
where $ {\bf \Lambda} (\e; \eta, \uF ) $ is defined in Proposition \ref{prop:inv-ap-vero}, 
there exists $  {\mathtt h}_{n+1} $ 
satisfying 
\begin{equation} \label{stima T 0 b} 
 \|   {\mathtt h}_{n+1} \|_{\Lip, s_1+2} \lesssim_{s_1}  \nu_n + 
 \e^2 \nu_n^{\frac{4}{5} - 5 \sigma} \, , \quad  \|   {\mathtt h}_{n+1} \|_{\Lip, s_3+2} \leq 
 \e^2 \nu_n^{- \frac{11}{10}-2\s} \, ,
\end{equation}
such that 
\be\label{def:rn+1-si}
{\mathtt r}_{n+1} := {\mathtt F}_n ( \vphi,0,0  ) + {\mathbb D}_n \,   {\mathtt h}_{n+1} 
\ee
satisfies 
\be \label{ap-in-no0}
 \| {\mathtt r}_{n+1} \|_{\Lip, s_1} \leq  \e^2 \nu_n^{\frac{3}{2}-6 \sigma} \, . 
\ee
\end{proposition}

\begin{pf}  
Recalling \eqref{def:operatore inverso approssimato proiettato}
we look for an approximate solution 
$  {\mathtt h}_{n+1} = (\widehat \phi, \widehat \ac, \widehat w ) $ of 
\begin{equation}\label{operatore inverso approssimato proiettato}
{\mathbb D}_n 
 \begin{pmatrix}
\widehat \phi \\
 \widehat \ac \\
\widehat w
\end{pmatrix}  = \begin{pmatrix}
\omega \cdot \partial_\vphi \widehat \phi  - 
{\mathtt K}_{20}^{(n)}(\vphi) \widehat \ac  - [{\mathtt K}_{11}^{(n)}]^\top(\vphi) \widehat w\\
\!\! \!\!\!\!\!\!\!\!\!\! \!\! \!\!\!\!\!\!\!\!\!\! \!\! \!\!\!\!\!\!\!\!\!\! \!\!\!\!\!\!\!\!\!\!\!\!\!\!\!\!
\!\!\!\!\!\!\!\!\!\!\!\!\!\! \omega \cdot \partial_\vphi \widehat \ac  \\
\!\!\!\!\!\!\!\!\!\!\!\! \!\!\!\!\!\!\!\!\!\! \!\! \!\!\!\!\!\!\!\!\!\! \!\! \!  {\cal L}_{\omega}^{(n)} \widehat w -J {\mathtt K}_{11}^{(n)}(\vphi)\widehat \ac  
\end{pmatrix} =
 \begin{pmatrix}
{\mathtt f}_1^{(n)} \\
{\mathtt f}_2^{(n)} \\
{\mathtt f}_3^{(n)}
\end{pmatrix} 
\end{equation}
where
$ ( {\mathtt f}_1^{(n)}, {\mathtt f}_2^{(n)}, {\mathtt f}_3^{(n)}) :=  - {\mathtt F}_n ( \vphi,0,0  ) $. Since
$ X_{{\mathtt K}_n} $ is a reversible vector field, its 
components $ {\mathtt f}_1^{(n)}, $ $ {\mathtt f}_2^{(n)},  $ $ {\mathtt f}_3^{(n)}$ satisfy  the reversibility property 
\begin{equation}\label{parita g1 g2 g3}
{\mathtt f}_1^{(n)}(\vphi) = {\mathtt f}_1^{(n)}(- \vphi)\,,\quad 
{\mathtt f}_2^{(n)}(\vphi) = - {\mathtt f}_2^{(n)}(- \vphi)\,,\quad 
{\mathtt f}_3^{(n)}(\vphi) = - (S {\mathtt f}_3^{(n)})(- \vphi)\,.
\end{equation}
We solve approximately \eqref{operatore inverso approssimato proiettato} in a triangular way.
\\[1mm]
{\sc Step 1. Approximate solution of the second equation  in \eqref{operatore inverso approssimato proiettato}}, i.e. 
$ \omega \cdot \partial_\vphi  \widehat \ac  = {\mathtt f}_2^{(n)}   $.
We solve the equation 
\be\label{appro-eq-2}
\omega \cdot \partial_\vphi  \widehat \ac  = {\it \Pi}_{N_n} {\mathtt f}_2^{(n)}(\vphi)
\ee
where $ N_n $ is defined in \eqref{def:Nn} and the projector $ {\it \Pi}_{N_n} $ 
applies to functions depending only on the variable $ \vphi $ as in \eqref{def:PiN-time}. 
Notice that  
 $ \omega $ is a Diophantine vector by \eqref{dioph} for all $ \lambda \in \Lambda $, and that,
by \eqref{parita g1 g2 g3}, the $\vphi$-average of $ {\it \Pi}_{N_n} {\mathtt f}_2^{(n)}  $ is zero. 
The solutions of \eqref{appro-eq-2} are
\begin{equation}\label{soleta}
\widehat \ac = \widehat \ac_0 + [\widehat \ac] \, , \quad 
[\widehat \ac] := (\omega \cdot \partial_\vphi )^{-1}  {\it \Pi}_{N_n} {\mathtt f}_2^{(n)}
 \, , \quad \widehat \ac_0 \in \R^\es \, , 
\end{equation}
where  $(\omega \cdot \partial_\vphi )^{-1} $ is defined in \eqref{op-inv-KAM}  and   $ \widehat \ac_0 $ is a free parameter that we fix below in \eqref{sol alpha}.
By \eqref{soleta}, \eqref{dioph}, \eqref{diof-est}, 
we have
\be\label{eq:acs}
\| [ \widehat \ac ] \|_{\Lip, s} \lesssim \| {\it \Pi}_{N_n} {\mathtt f}_2^{(n)} \|_{\Lip,s+ \tau_1} 
\stackrel{\eqref{smoothingS1S2-Lip}} \lesssim N_n^{\tau_1}  \|  {\mathtt f}_2^{(n)} \|_{\Lip,s}
\ee
($ \g_2 = \g_0 / 4 $ is considered  as a fixed constant).
Therefore, using \eqref{eq:acs}, \eqref{NM:indu-new}, \eqref{def:Nn}, and taking $ s_3 - s_1 $ large enough, we have 
for $\nu_n$ small enough (depending on $s_3$)
\be\label{est:zeta}
\| [ \widehat \ac ] \|_{\Lip, s_1+2} 
 \leq \e^2 \nu_n^{1- 3 \sigma} \, , \quad 
\| [ \widehat \ac ] \|_{\Lip, s_3+2} 
\leq \e^2 \nu_n^{- \frac{4}{5} - 4 \sigma} \, . 
\ee
{\sc Step 2. Approximate solution of the third equation  in \eqref{operatore inverso approssimato proiettato}}, i.e. 
\be\label{ultima-eq}
{\cal L}_\om^{(n)} \widehat w = {\mathtt f}_3^{(n)}   + J {\mathtt K}_{11}^{(n)} \widehat \ac_0 
 + J {\mathtt K}_{11}^{(n)} [ \widehat \ac  ] \, . 
\ee
For some $\nu \in (0,\e^{\frac32})$ such that 
\be \label{conditionnu}
\| \breve \fracchi_n \|_{\Lip, s_3+4} \leq \e \nu^{- \frac{9}{10}} 
\ee
and which will be chosen later in \eqref{choice:nul}, we consider the operator
$ {\cal L}^{-1}_{appr, \nu} $ defined in Proposition \ref{prop:inv-ap-vero} (with $ \uF  \leadsto \breve \fracchi_n $). 
Notice that, by \eqref{NM:indu-line2-trunc}, also the assumption 
$ \| \breve \fracchi_n \|_{\Lip, s_2+2} \leq \e $ required in 
Proposition \ref{prop:inv-ap-vero} with $ \uF  \leadsto \breve \fracchi_n $ holds. 
Moreover ${\cal L}^{-1}_{appr, \nu} $  satisfies \eqref{new:estimate-cal-L-1} (independently of $ \nu $).
We define the approximate solution $\widehat w$ of \eqref{ultima-eq}, 
\be\label{def:Lappro}
\widehat w := {\cal L}^{-1}_{appr,\nu} \big( {\mathtt f}_3^{(n)}    + J {\mathtt K}_{11}^{(n)} [ \widehat \ac  ] \big)  +
 {\cal L}^{-1}_{appr, \nu} J {\mathtt K}_{11}^{(n)} \widehat \ac_0   \, . 
\ee 
{\sc Step 3. Approximate solution of the first equation  in \eqref{operatore inverso approssimato proiettato}.}
With \eqref{def:Lappro}, the first equation in \eqref{operatore inverso approssimato proiettato} can be written as
\be\label{equazione psi hat}
\begin{aligned}
\omega \cdot \partial_\vphi \widehat \phi &  = T_n \widehat \ac_0  
 + g_n   
\end{aligned}
\ee
where
\be\label{def:Tn-gn}
\begin{aligned}
& T_n  := 
 {\mathtt K}_{20}^{(n)} + 
 [{\mathtt K}_{11}^{(n)}]^\top  {\cal L}^{-1}_{appr} J {\mathtt K}_{11}^{(n)}   \, , \\
& g_n := {\mathtt f}_1^{(n)}  + {\mathtt K}_{20}^{(n)} [\widehat \zeta] + 
[{\mathtt K}_{11}^{(n)}]^\top  {\cal L}^{-1}_{appr} \big( {\mathtt f}_3^{(n)}  
 + J {\mathtt K}_{11}^{(n)} [ \widehat \ac  ] \big)    \, . 
 \end{aligned}
\ee
We  have 
to choose the constant 
$ \widehat \ac_0 \in \R^\es $ such that the right hand side in \eqref{equazione psi hat}  has zero average.  
By \eqref{stime coefficienti K 20 11 bassa},  the matrix 
$$ 
\langle {\mathtt K}_{20}^{(n)} \rangle := 
(2 \pi)^{-\es} \int_{\T^\es}  {\mathtt K}_{20}^{(n)} (\vphi) d \vphi 
$$ 
is close to the invertible twist matrix 
$ \e^2 \Ab $ (see \eqref{A twist}) and therefore by  \eqref{A twist}  there is a constant $ C $ such that
\be\label{averageK20}
\| \langle {\mathtt K}_{20}^{(n)}  \rangle^{-1} \| \leq C \e^{-2}  \, . 
\ee
By \eqref{def:Tn-gn}  the operator $ \langle T_n \rangle : \zeta_0 \in \R^\es \mapsto 
 \langle  T_n \zeta_0 \rangle  \in \R^\es $ is decomposed as
\be\label{Tnsvi}
\langle T_n \rangle := 
\langle  {\mathtt K}_{20}^{(n)} \rangle+ 
\langle [{\mathtt K}_{11}^{(n)}]^\top  {\cal L}^{-1}_{appr} J {\mathtt K}_{11}^{(n)}  \rangle \, . 
\ee
Now, for $ \zeta_0  \in \R^\es $, 
\begin{align}
| \langle  [{\mathtt K}_{11}^{(n)}]^\top  {\cal L}^{-1}_{appr} J {\mathtt K}_{11}^{(n)}  \zeta_0 \rangle  |_\Lip
& \leq \|   [{\mathtt K}_{11}^{(n)}]^\top  {\cal L}^{-1}_{appr} J {\mathtt K}_{11}^{(n)}   \zeta_0  \|_{\Lip , s_0} \nonumber \\
& \stackrel{\eqref{stime coefficienti K 11 alta trasposto},\eqref{NM:indu-line2-trunc}}{\leq} 
\e^2    \|    {\cal L}^{-1}_{appr} J {\mathtt K}_{11}^{(n)}  \zeta_0   \|_{\Lip , s_0} \nonumber \\
& \stackrel{\eqref{new:estimate-cal-L-1}}{\leq } C(s_1) \e^2   \|  {\mathtt K}_{11}^{(n)}  \zeta_0   \|_{\Lip , s_0+Q'} \nonumber \\
&\stackrel{ \eqref{stime coefficienti K 11 alta}}{\leq} C(s_1) \e^4
\big(   | \zeta_0|_\Lip +   \|  \breve \fracchi_n  \|_{\Lip, s_0 + \bt + Q'}   | \zeta_0|_\Lip  \big)   \nonumber \\
& \stackrel{\eqref{NM:indu-line2-trunc}}{\leq}  C(s_1) \e^4  | \zeta_0|_\Lip   \label{pert4} 
\end{align}
using the fact that $s_0 + \bt + Q' \leq s_1$.  By \eqref{Tnsvi}, \eqref{averageK20}, 
\eqref{pert4}, we deduce that, for $ \e $ small enough, $ \langle T_n \rangle$
is invertible and 
\be\label{Tn-1:estimate}
\| \langle T_n \rangle^{-1} \| \leq 2 C \e^{-2} \, . 
\ee
Thus, in view of \eqref{equazione psi hat},  we define 
\begin{equation}\label{sol alpha}
\widehat \ac_0  := - \langle T_n \rangle^{-1} \langle g_n \rangle  \in \R^\es \, .
\end{equation}
We now estimate the function $ g_n $ defined in \eqref{def:Tn-gn}. 
By  
\eqref{tame:K20-appl}, \eqref{stime coefficienti K 11 alta trasposto}, \eqref{NM:indu-line2-trunc}, \eqref{new:estimate-cal-L-1},  and using that $s_1 > s_0 + Q' +\bt $,  
\begin{align}
\| g_n \|_{\Lip , s_0}  &\leq  \| {\mathtt f}_1^{(n)}  \|_{\Lip ,s_0} + \| {\mathtt K}_{20}^{(n)} [\widehat \zeta] \|_{\Lip ,s_0} + 
\|[{\mathtt K}_{11}^{(n)}]^\top  {\cal L}^{-1}_{appr} \big( {\mathtt f}_3^{(n)}  
 + J {\mathtt K}_{11}^{(n)} [ \widehat \ac  ] \big) \|_{\Lip ,s_0} \nonumber \\
 & \lesssim_{s_0}  \| {\mathtt f}_1^{(n)}  \|_{\Lip ,s_0} + \e^2 \| [\widehat \zeta] \|_{\Lip ,s_0} + 
\e^2 \big(  \|{\mathtt f}_3^{(n)} \|_{\Lip ,s_0+Q'} 
 + \| {\mathtt K}_{11}^{(n)} [ \widehat \ac  ] \|_{\Lip ,s_0+Q'} \big) \, . \nonumber 
\end{align}
Hence by \eqref{est:zeta}, 
\eqref{stime coefficienti K 11 alta}, \eqref{NM:indu-line2-trunc},
\begin{align} \label{gns0}
\| g_n \|_{\Lip , s_0}  & \lesssim_{s_0} 
\| {\mathtt f}_1^{(n)}  \|_{\Lip ,s_0} + \e^4 \nu_n^{1-3\s} + \e^2
\big(   \|{\mathtt f}_3^{(n)} \|_{\Lip ,s_0+Q'}  +  \e^4 \nu_n^{1-3\s} \big) \nonumber \\
& \stackrel{\eqref{NM:indu-new}}{\lesssim_{s_1}} \e^2 \nu_n + \e^4 \nu_n^{1-3\s} +
\e^2 \big(  \e^2 \nu_n^{1- 2 \sigma} + \e^4  \nu_n^{1-3\s} \big) \nonumber \\
& \lesssim_{s_1} \e^2 \nu_n+ \e^4 \nu_n^{1- 3 \sigma}   \, .
\end{align}
Then the constant  $ \widehat \ac_0 $ defined in \eqref{sol alpha} satisfies, by \eqref{Tn-1:estimate}, \eqref{gns0}, 
\be \label{estzeta0}
| \widehat \ac_0 |_{\Lip}   
 \lesssim_{s_1} \nu_n+ \e^2 \nu_n^{1- 3 \sigma}  \, . 
\ee
Using   \eqref{NM:indu-new}, \eqref{stime coefficienti K 11 alta}, \eqref{estzeta0}, \eqref{est:zeta} 
and \eqref{NM:indu-line2-trunc}, 
 we have 
\be\label{est:g3-new}
\begin{aligned}
& \| {\mathtt f}_3^{(n)}  \|_{\Lip,s_1} + \| J {\mathtt K}_{11}^{(n)} \widehat \ac_0 \|_{\Lip,s_1} + 
\| J {\mathtt K}_{11}^{(n)} [\widehat \ac ] \|_{\Lip,s_1} \leq \e^2 \nu_n^{1 - 3 \sigma} \, , \\
& \| {\mathtt f}_3^{(n)} \|_{\Lip,s_3+2} + \| J {\mathtt K}_{11}^{(n)} \widehat \ac_0 \|_{\Lip,s_3+2} +
\| J {\mathtt K}_{11}^{(n)} [\widehat \ac] \|_{\Lip,s_3+2} \leq \e^2 \nu_n^{- \frac{4}{5} - 5 \sigma}
\end{aligned}
\ee
for $ \nu_n $ small enough.

We now choose the constant $ \nu $ in \eqref{conditionnu} as 
\be\label{choice:nul}
\nu :=  \nu_n^{1-3\s} 
\ee
so that $ \nu \in (0, \e^{3/2}) $ and,  by  the third inequality in \eqref{NM:indu-line2-trunc},
condition \eqref{conditionnu} is satisfied, provided $\s$ is chosen small enough (and hence $s_3$ is large enough). 
We apply Proposition \ref{prop:inv-ap-vero} with
$ g = {\mathtt f}_3^{(n)} + J {\mathtt K}_{11}^{(n)} \widehat \ac_0 +J
 {\mathtt K}_{11}^{(n)} [\widehat \ac ] $  (and  $ \uF  \leadsto \breve \fracchi_n $) 
 noting that  
\eqref{est:g3-new} implies \eqref{g:small-large-IA}  with $ \nu $ defined in \eqref{choice:nul}, again provided  $ \s $  is chosen small enough.
As a consequence, by \eqref{h:s1s2s3-IA-per} the function $\widehat w $ defined in \eqref{def:Lappro} satisfies (for $\nu_n$ small enough)
\be\label{est:hatw}
\begin{aligned}  
 &\|  \widehat w \|_{\Lip, s_1} \leq C(s_1) \e^2 \nu_n^{\frac{4}{5}(1 -3\s)} \leq  \e^2 \nu_n^{\frac{4}{5}-3\s}  \, , \\ 
 &\| \widehat w \|_{\Lip, s_3+2} \leq  C(s_3)   \e^2 \nu_n^{ - \frac{11}{10} (1-3\s)}  \leq  \e^2 \nu_n^{ - \frac{11}{10} +3\s} \, . 
\end{aligned} 
\ee
By \eqref{est:hatw} and interpolation inequality \eqref{inter-Sobo-Lip}, arguing as above \eqref{15-6-2017-bis}, we obtain
\be \label{ws1+2}
\|  \widehat w \|_{\Lip, s_1+2} \leq  \e^2 \nu_n^{\frac{4}{5}-4\s} \, ,
\ee
for $s_3$ large enough and $\nu_n$ small enough (depending on $s_3$).
Moreover, by \eqref{estimate:r-final-IA-per}, 
\begin{equation}\label{normalw}
{\cal L}_\om^{(n)}  \widehat w - \big( {\mathtt f}_3^{(n)} 
+ J {\mathtt K}_{11}^{(n)} \widehat \ac_0 
+ J {\mathtt K}_{11}^{(n)} [\widehat \ac] \big) = \rf_{n+1} \, , 
\end{equation}
where $ \rf_{n+1}  $ satisfies 
\be\label{est:frak-r}
 \| \rf_{n+1} \|_{\Lip, s_1} \leq  C(s_1)  \e^2 \nu_n^{\frac{3}{2} (1-3\s)} \leq \e^2 \nu_n^{\frac{3}{2}-5\s} \, .
\ee
Now, since $ \widehat \ac_0 $ has been chosen in \eqref{sol alpha} so that $ T_n \widehat \ac_0 + g_n $
has zero mean value,  the equation \eqref{equazione psi hat} has the approximate solution 
\begin{equation}\label{sol psi}
\widehat \phi := (\omega \cdot \partial_\vphi)^{-1} {\it \Pi}_{N_n} ( T_n \widehat \ac_0 + g_n)  \, . 
\end{equation} 
Recalling the definition of $ g_n $ and $ T_n $  in \eqref{def:Tn-gn} we have
\be\label{esprgn}
g_n + T_n  \widehat \zeta_0 =  {\mathtt f}_1^{(n)} + {\mathtt K}_{20}^{(n)} \widehat{\ac} + [{\mathtt K}_{11}^{(n)}]^\top \widehat{w}  
\ee
where $ \widehat w $  is defined in \eqref{def:Lappro}.
By \eqref{esprgn} we estimate 
\begin{align}
\|   g_n + T_n  \widehat \zeta_0 \|_{\Lip , s_1} & 
\stackrel{\eqref{tame:K20-appl},  \eqref{stime coefficienti K 11 alta trasposto},\eqref{NM:indu-line2-trunc}} \leq \| {\mathtt f}_1^{(n)} \|_{\Lip , s_1} + C(s_1) \e^2 
\big( \|  \widehat \zeta \|_{\Lip , s_1} + \|  \widehat w \|_{\Lip , s_1}  \big)  \nonumber \\
& \stackrel{\eqref{NM:indu-new},  \eqref{est:zeta}, \eqref{estzeta0}, \eqref{est:hatw}} 
\leq \e^2 \nu_n^{1-2\s} + C(s_1) \e^2 \big( \nu_n + \e^2 \nu_n^{1-3\s} + \e^2 \nu_n^{\frac45 - 3 \s}   \big) \nonumber \\
\label{gns1}
& \leq  \e^2 \nu_n^{\frac45 - 3 \s} 
\end{align}
and, similarly, by \eqref{tame:K20-appl},  \eqref{stime coefficienti K 11 alta trasposto},\eqref{NM:indu-line2-trunc},  
\eqref{NM:indu-new},  \eqref{est:zeta}, \eqref{estzeta0}, \eqref{est:hatw}, we get  
\begin{align}
\|   g_n + T_n  \widehat \zeta_0 \|_{\Lip , s_3} & 
\leq \| {\mathtt f}_1^{(n)} \|_{\Lip , s_3} + C(s_3) \e^2 
\big( \|  \widehat \zeta \|_{\Lip , s_3} + \|  \widehat w \|_{\Lip , s_3}  \big) \nonumber \\
& \qquad \qquad \quad \quad + C(s_3) \e^2 \|\breve \fracchi_n \|_{\Lip , s_3+ \bt } 
\big(   \|  \widehat \zeta \|_{\Lip , s_0} + \|  \widehat w \|_{\Lip , s_0}  \big) \nonumber \\
& 
\leq \e^2 \nu_n^{-\frac45-3\s} + C(s_3) \e^2 \nu_n^{-\frac{11}{10}+3\s} + C(s_3)   \e^2 \nu_n^{-\frac45-\s}    \nu_n^{\frac45-3\s} \nonumber \\
\label{gns3}
& \leq  \e^2 \nu_n^{-\frac{11}{10}} \, . 
\end{align}
The function $\widehat \phi$ defined in \eqref{sol psi} satisfies
by \eqref{dioph}, \eqref{diof-est}, \eqref{smoothingS1S2-Lip}, \eqref{def:Nn}, 
for $ s_3 - s_1 $ large enough,  and  \eqref{gns1}, \eqref{gns3}, 
\be\label{est:hatphi}
\| \widehat \phi  \|_{\Lip, s_1} \leq  \e^2 \nu_n^{\frac{4}{5}-  4 \s} \, , \quad 
\| \widehat \phi  \|_{\Lip, s_3+2} \leq  \e^2 \nu_n^{- \frac{11}{10}-\s} \, . 
\ee
Using again the interpolation inequality \eqref{inter-Sobo-Lip}, we deduce by 
\eqref{est:hatphi},  the estimate
\be \label{phis1+2}
\| \widehat \phi  \|_{\Lip, s_1+2} \leq  \e^2 \nu_n^{\frac{4}{5}-  5 \s} \, ,
\ee
for $s_3$ large enough and $\nu_n$ small enough.
Moreover  the remainder 
\be\label{def:r3n+1} 
{\mathfrak r}_3^{(n+1)} := {\it \Pi}_{N_n}^\bot (  T_n \widehat \ac_0 + g_n ) 
\ee
 satisfies
\be\label{est:r3-nel}
\|  {\mathfrak r}_3^{(n+1)}  \|_{\Lip, s_1} \stackrel{\eqref{smoothingS1S2-Lip}} 
\leq N_n^{-(s_3-s_1)}  \|  T_n \widehat \ac_0 + g_n \|_{\Lip, s_3}  \stackrel{\eqref{def:Nn}, \eqref{gns3}} 
 \leq  \e^2 \nu_n^{ \frac{3}{2}}   \, . 
\ee
{\sc Step 4. Conclusion.} We set 
$ {\mathtt h}_{n+1} := (\widehat \phi,  \widehat \ac, \widehat w ) $
defined in  \eqref{sol psi},  \eqref{soleta}, \eqref{sol alpha} and \eqref{def:Lappro}. 
The estimates in  \eqref{stima T 0 b} follow by
\eqref{est:hatphi}-\eqref{phis1+2}, \eqref{est:zeta}, \eqref{estzeta0}, \eqref{est:hatw}-\eqref{ws1+2} and recalling that 
$ \nu_n \leq \nu_1 = C(s_1) \e^2 $.
Finally, by \eqref{operatore inverso approssimato proiettato}, 
\eqref{soleta}, \eqref{sol alpha}, \eqref{normalw},  we have that 
\be\label{forma-resto}
{\mathbb D}_n 
\begin{pmatrix}
\widehat \phi \\
 \widehat \ac \\
\widehat w
\end{pmatrix}
 -   \begin{pmatrix}
{\mathtt f}_1^{(n)} \\
{\mathtt f}_2^{(n)} \\
{\mathtt f}_3^{(n)}
\end{pmatrix}  = 
\begin{pmatrix}
-{\mathfrak r}_3^{(n+1)} \\
{\it \Pi}_{N_n}^\bot {\mathtt f}_2^{(n)} \\
\rf_{n+1}   
\end{pmatrix}
\ee
where $ {\mathfrak r}_3^{(n+1)} $ is defined in \eqref{def:r3n+1} and 
$ \rf_{n+1} $  in \eqref{normalw}. 
The estimate  \eqref{ap-in-no0}  follows by \eqref{forma-resto}, \eqref{est:r3-nel}, the bound 
$$
\| {\it \Pi}_{N_n}^\bot {\mathtt f}_2^{(n)} \|_{\Lip , s_1} \stackrel{\eqref{smoothingS1S2-Lip}} 
 \leq
N_n^{-(s_3-s_1)}  \| {\mathtt f}_2^{(n)} \|_{\Lip , s_3}  \stackrel{\eqref{def:Nn}, \eqref{NM:indu-new}}{\ll} \e^2 \nu_n^{\frac32}
$$
and \eqref{est:frak-r}. 
\end{pf}

The function $  {\mathtt h}_{n+1} $ defined in Proposition \ref{Prop:inversione} is an approximate solution of equation 
\eqref{Newtonnewvar}, according to the following corollary.  

\begin{corollary}\label{cor:NM}
The function
\be\label{def:rn+1'}
{\mathtt r}_{n+1}' := {\mathtt F}_n (\vphi,0,0) +
\big( \om \! \cdot \! \pa_\vphi 
- d_{\mathtt u} X_{{\mathtt K}_n}( \vphi,0,0 )  \big)  {\mathtt h}_{n+1}
\ee
satisfies 
\be\label{fina-n+1}
\| {\mathtt r}_{n+1}'  \|_{\Lip,s_1} \leq \e^2 \nu_n^{\frac32 - 8 \sigma} \, . 
\ee
\end{corollary}

\begin{pf}
The term $ {\mathtt r}_{n+1}'  $ in \eqref{def:rn+1'} is,
by  \eqref{splitting linearizzato nuove coordinate} and \eqref{def:rn+1-si}, 
\be\label{rn+1++}
{\mathtt r}_{n+1}'  = {\mathtt r}_{n+1} +  {\mathtt R}_{Z_n}   {\mathtt h}_{n+1} \, .
\ee
By the expression of $  {\mathtt R}_{Z_n}  $ in \eqref{term:RZ}, using the tame 
estimate \eqref{inter:pro-Lip} for the product of functions and \eqref{duality-K01T}, we get 
\begin{align}
\|  {\mathtt R}_{Z_n}   {\mathtt h}_{n+1}\|_{\Lip,s_1} &\lesssim_{s_1}
\big( \| {\mathtt K}_{10}^{(n)} - \omega \|_{\Lip, s_1+1} + \| {\mathtt K}_{01}^{(n)} \|_{\Lip , s_1+1} + \| \partial_\phi {\mathtt K}_{00}^{(n)} \|_{\Lip , s_1+1} \big)
\| {\mathtt h}_{n+1}\|_{\Lip , s_1}  \nonumber \\
& \stackrel{\eqref{K 00 10 01},  \eqref{NM:indu-line2-trunc}}{\lesssim_{s_1}}
  \| {\mathtt F}_n (\vphi,0,0) \|_{\Lip , s_1 +1+ \underline{\tau}} 
\| {\mathtt h}_{n+1}\|_{\Lip , s_1} \label{estRZn}  
\end{align} 
having used 
that  $ s_1+ 1+ \bt \leq s_2$ and the estimate $\| \breve \fracchi_n \|_{\Lip , s_2} \leq \e $  in 
\eqref{NM:indu-line2-trunc}.  Now, by \eqref{NM:indu-new} and the 
interpolation inequality \eqref{inter-Sobo-Lip}, 
we get, for $ s_3 - s_1 $ large enough,  
\be\label{Fn-inte}
\| {\mathtt F} _n (\vphi,0,0) \|_{\Lip , s_1 +1+ \underline{\tau}} \leq \e^2 \nu_n^{1- 3 \s} \, .
\ee
So we derive, by \eqref{estRZn}, \eqref{Fn-inte} and  \eqref{stima T 0 b},  the estimate
\be\label{RZ:new}
\|  {\mathtt R}_{Z_n}   {\mathtt h}_{n+1}\|_{\Lip,s_1} \leq \e^2 \nu_n^{ \frac{9}{5} - 8 \sigma} \, . 
\ee
In conclusion \eqref{rn+1++}, \eqref{ap-in-no0},  \eqref{RZ:new} 
imply \eqref{fina-n+1}. 
\end{pf}
\\[1mm]
{\bf Step 5. Approximate solution $ i_{n+1}$. } 
Finally, for all $  \l $ in the set  $ {\bf \Lambda}_{n+1} $ introduced in \eqref{def:Lambda-n+1}, we define
the new approximate solution of the Nash-Moser iteration in the original coordinates as
\be\label{def:in+1}
i_{n+1} := i_{n, \d} + h_{n+1} \, , \quad h_{n+1}  :=  DG_{n,\delta} (\vphi,0,0)  {\mathtt h}_{n+1}  \, ,
\ee
where $ i_{n, \d} $ is the isotropic torus defined in Lemma \ref{lem:ison} and $  {\mathtt h}_{n+1} $ is 
the function defined in Proposition \ref{Prop:inversione}. 
By \eqref{DG delta},  \eqref{NM:indu-line2-trunc} coupled with the inequality $s_1+4 \leq s_2+2$ and \eqref{stima T 0 b},
we have 
\begin{align}\label{Hn+1:s1}
& \| h_{n+1} \|_{\Lip, s_1+2}  
\lesssim_{s_1}  \|  {\mathtt h}_{n+1} \|_{\Lip, s_1+2} \lesssim_{s_1} \nu_n +  \e^2 \nu_n^{\frac45 - 5 \sigma} \\
& \| h_{n+1} \|_{\Lip, s_3+2}  \lesssim_{s_3}  \|  {\mathtt h}_{n+1} \|_{\Lip, s_3+2} + \| \breve \fracchi_n \|_{\Lip , s_3 +4}   \|  {\mathtt h}_{n+1} \|_{\Lip, s_0}  \leq  \e^2 \nu_n^{- \frac{11}{10} - 3 \sigma}  \label{Hn+1:s3}
\end{align}
for $\nu_n$ small enough.
\begin{lemma}\label{lem:rhon+1}
The term
\be\label{def:varn+1}
\varrho_{n+1}  := {\cal F} (i_{n,\d}) + d_i  {\cal F} (i_{n,\d}) h_{n+1} 
\ee
satisfies 
\be\label{new-restn}
\| \varrho_{n+1} \|_{\Lip,s_1} \leq \e^2 \nu_n^{\frac32 - 9 \sigma} \, . 
\ee
\end{lemma}

\begin{pf}
Differentiating the identity (see \eqref{trasfo imp})
\be\label{laid1}
D G_{n, \delta}( {\mathtt u}  (\vphi) ) {\mathtt F}_n ({\mathtt u} (\vphi) ) =  {\cal F}(G_{n, \delta}(  {\mathtt u} (\vphi) ) )
\ee
 we obtain
$$
\begin{aligned}
& D^2 G_{n,\delta} ({\mathtt u} (\vphi))  [  {\mathtt h},  
{\mathtt F}_n ({\mathtt u} (\vphi) )] +   D G_{n,\delta} ( {\mathtt u} (\vphi) ) 
d_{\mathtt u}  {\mathtt F}_n ({\mathtt u} (\vphi) ) [  {\mathtt h}  ] \\ 
& = 
d_i {\cal F} (G_{n, \delta}({\mathtt u} (\vphi)) ) DG_{n,\delta} ( {\mathtt u} (\vphi)  ) [   {\mathtt h}  ] \, . 
\end{aligned}
$$
For $ {\mathtt u}  (\vphi)  = (\vphi,0,0)$ and $  {\mathtt h} =  {\mathtt h}_{n+1} $ this 
gives, recalling
\eqref{torus:origin},  \eqref{def:in+1},  
\be\label{dd2}
\begin{aligned}
& D^2 G_{n,\delta} (\vphi,0,0)  [  {\mathtt h}_{n+1},  
{\mathtt F}_n (\vphi,0,0 )] +  DG_{n,\delta} ( \vphi,0,0) d_{\mathtt u}  {\mathtt F}_n (\vphi,0,0) [  {\mathtt h}_{n+1}  ] \\
& = d_i {\cal F} (i_{n, \delta}(\vphi) ) h_{n+1} \, . 
\end{aligned}
\ee
By \eqref{laid1} we have $ {\cal F} (i_{n,\d}) =  D G_{n, \delta}( \vphi,0,0 ) {\mathtt F}_n (\vphi,0,0) $ and
therefore, by \eqref{dd2}, 
\begin{align*}
&  {\cal F} (i_{n,\d})  + d_i {\cal F} (i_{n,\delta}) h_{n+1} \\
&  = 
  D G_{n, \delta}( \vphi,0,0 ) \big( {\mathtt F}_n (\vphi,0,0) + d_{\mathtt u}  {\mathtt F}_n (\vphi,0,0) [  {\mathtt h}_{n+1}  ] \big) 
 + D^2 G_{n,\delta} (\vphi,0,0)  [  {\mathtt h}_{n+1},  {\mathtt F}_n (\vphi,0,0 )] \\
& \stackrel{\eqref{def:rn+1'}, \eqref{Ham-new-expa}}=   D G_{n, \delta}( \vphi,0,0 ) {\mathtt r}_{n+1}'   
  + D^2 G_{n,\delta} (\vphi,0,0)  [  {\mathtt h}_{n+1},  {\mathtt F}_n (\vphi,0,0 )] \, . 
\end{align*}
In conclusion the term $\varrho_{n+1} $ in \eqref{def:varn+1} satisfies 
\begin{align}
\| \varrho_{n+1} \|_{\Lip,s_1} &\leq 
\| D G_{n, \delta}( \vphi,0,0 ) {\mathtt r}_{n+1}'  \|_{\Lip,s_1} 
 + \|D^2 G_{n,\delta} (\vphi,0,0)  [  {\mathtt h}_{n+1},  {\mathtt F}_n (\vphi,0,0 )]\|_{\Lip,s_1} 
  \nonumber \\
 &\stackrel{ \eqref{DG delta} , \eqref{DG2 delta} } {\lesssim_{s_1}} 
 \big(   1+\|  \breve \fracchi_n  \|_{\Lip , s_1+2}  \big)  \| {\mathtt r}_{n+1}'\|_{\Lip , s_1} 
 \nonumber \\
 & \qquad \qquad \ + 
\big(   1+\|  \breve \fracchi_n  \|_{\Lip , s_1+3}  \big)  \| {\mathtt h}_{n+1}\|_{\Lip , s_1}  \|  {\mathtt F}_n (\vphi,0,0 )\|_{\Lip , s_1}  \nonumber  \\
& \stackrel{\eqref{NM:indu-line2-trunc}}{\lesssim_{s_1} } \| {\mathtt r}_{n+1}'\|_{\Lip , s_1} 
+ \| {\mathtt h}_{n+1}\|_{\Lip , s_1}  \|  {\mathtt F}_n (\vphi,0,0 )\|_{\Lip , s_1}  
\label{lem:fin}
\end{align}
where we use that $s_1+3 \leq s_2 +2$.  Hence by  \eqref{lem:fin}, \eqref{fina-n+1}, 
 \eqref{stima T 0 b}, \eqref{NM:indu-new}, we get
$$
\| \varrho_{n+1} \|_{\Lip,s_1}  \lesssim_{s_1} \e^2 \nu_n^{\frac32 - 8\s} + (\nu_n + \e^2 \nu_n^{\frac45 - 5\s}) \e^2 \nu_n^{1 - 2\s} 
$$ 
which gives \eqref{new-restn}. 
\end{pf}

\begin{lemma}\label{lem:Fin+1}
$ \| {\cal F} ( i_{n+1}) \|_{\Lip,s_1} \leq \e^2 \nu_n^{\frac32 - 10 \sigma} $. 
\end{lemma}

\begin{pf}
By \eqref{def:in+1} we have 
\be\label{Tay1}
{\cal F}( i_{n+1}) = {\cal F}( i_{n, \d} + h_{n+1}) = 
{\cal F}( i_{n, \d})  + d_i {\cal F}( i_{n, \d}) h_{n+1} + Q( i_{n,\d}, h_{n+1}) \, 
\ee
and,  by the form of $ {\cal F} $ in \eqref{operatorF}, 
\begin{eqnarray*}
Q( i_{n,\d}, h_{n+1})& =&\wtilde{\cal F}( i_{n, \d} + h_{n+1}) - 
\wtilde{\cal F}( i_{n, \d})  - d_i \wtilde{\cal F}( i_{n, \d}) h_{n+1} \\
&=&\int_0^1 (1- \tau ) (D^2_i \wtilde{\cal F})( i_{n, \d} + \tau h_{n+1}) [h_{n+1},h_{n+1}] \, d \tau 
\end{eqnarray*}
where $\wtilde{\cal F}$ is the ``nonlinear part'' of $\cal F$ defined as 
$$
\wtilde{\cal F} (i) := \left(
\begin{array}{c}
  - \e^2 (\partial_y R) ( i(\vphi), \xi   )   \\
 \e^2 (\partial_\teta R) ( i(\vphi), \xi  )  \\
- \e^2 \big( 0, (\nabla_Q R) (i(\vphi), \xi) \big)
\end{array}
\right) 
$$
and $R$ is the Hamiltonian defined in \eqref{restoNN} and \eqref{defv}.  We have
\begin{align}
\| Q( i_{n,\d}, h_{n+1})\|_{\Lip , s_1} & \lesssim_{s_1}   \e^2 \big(1+ \|\fracchi_{n,\d}\|_{\Lip , s_1} + \| h_{n+1} \|_{\Lip , s_1}\big)
\| h_{n+1} \|_{\Lip , s_1}^2 \nonumber \\
&\stackrel{\eqref{NM:indu-line2-trunc}  , \eqref{eq:diff-con-isot}}{\lesssim_{s_1} }  \e^2  \| h_{n+1} \|_{\Lip , s_1}^2 \nonumber \\
\label{Qest1}
& \stackrel{\eqref{Hn+1:s1}}{\lesssim_{s_1}}  \e^2 \nu_n^{\frac85 - 10 \s} \, . 
\end{align}
The lemma follows by \eqref{Tay1},  Lemma \ref{lem:rhon+1} and \eqref{Qest1}. 
\end{pf}
\\[1mm]
{\bf Step 6. End of the Induction.} 
We are now ready to prove  that the approximate solution  $i_{n+1} $ defined in \eqref{def:in+1} for all
$ \l $ in  the set $ {\bf \Lambda}_{n+1} \subset {\bf \Lambda}_{n} $ introduced in \eqref{def:Lambda-n+1} satisfies 
 the estimates  \eqref{NM:indu-line2}-\eqref{vicin-in+1-in-s2} and  \eqref{NM:indu}  at step $ n + 1 $. 
 By Lemma \ref{lem:Fin+1} and \eqref{def:nun-NM} we obtain
 $$
 \| {\cal F} ( i_{n+1}) \|_{\Lip,s_1} \leq \e^2 \nu_n^{\frac32 - 10 \sigma} = \e^2 \nu_{n+1} $$  
 with
 \be \label{defsigstar} \s_\star := 10 \sigma \, , 
\ee
proving \eqref{NM:indu}$_{n+1}$. Moreover we have 
\begin{align} \label{diff-i-n+1-in}
\| i_{n+1} - i_n \|_{\Lip,s_1+2} 
& \quad \leq \|  i_{n+1} - i_{n,\d} 
\|_{\Lip,s_1+2} +  \| i_{n,\d} - \breve \imath_n \|_{\Lip,s_1+2} 
 +  \| \breve \imath_n - i_{n}  \|_{\Lip,s_1+2} \nonumber  \\
& \! \! \! \! \! \! \stackrel{\eqref{def:in+1}, \eqref{Hn+1:s1},  \eqref{eq:diff-con-isot}, \eqref{NM:int2}} 
 {\lesssim_{s_1}}   \nu_n  +   \e^2 \nu_n^{\frac{4}{5} - 5 \sigma}  + 
\e^2 \nu_n^{1 - \sigma} + \e^2 \nu_n^{2} \nonumber \\
&\quad   \leq C(s_1) \nu_n  +   \e^2 \nu_n^{\frac{4}{5} - 6 \sigma}  
\end{align}
proving \eqref{vicin-in+1-in}$_{n+1}$ since $ \s_\star := 10 \s $. Similarly  we get 
 \begin{align} 
\| i_{n+1} - \breve \imath_n \|_{\Lip,s_3+2} 
& \leq \| i_{n+1} - i_{n,\d}  \|_{\Lip,s_3+2} +  \| i_{n,\d} - \breve \imath_n \|_{\Lip,s_3+2} 
  \nonumber  \\
& 
\stackrel{\eqref{def:in+1}, \eqref{Hn+1:s3}, 
 \eqref{bound-int-iso}} \leq  \e^2 \nu_n^{- \frac{11}{10} - 4 \s} \label{nuovas} \\
 & \stackrel{\eqref{def:nun-NM}, \eqref{defsigstar}} \ll \e^2 \nu_{n+1}^{- \frac{4}{5}}  \label{diffns3new}
\end{align}
and, by \eqref{diffns3new} and \eqref{NM:int2},
\be \label{diffns3}
\| i_{n+1} - i_n \|_{\Lip,s_3+2} 
 \leq \| i_{n+1} - \breve \imath_{n}  \|_{\Lip,s_3+2} +  \|  \breve \imath_n - i_{n}  \|_{\Lip,s_3+2} 
 \ll \e^2 \nu_{n+1}^{- \frac{4}{5}}  
\ee
provided $\s$ is chosen small enough.

The bound  \eqref{vicin-in+1-in-s2}$_{ n + 1} $ for $ \|  i_{n+1} - i_n  \|_{\Lip, s_2} $  follows by interpolation:
setting
$$
s_2 + 2 = \theta s_1 + (1-  \theta) s_3 \, , \quad 
\theta = \frac{s_3-s_2-2}{s_3 - s_1}  \, , \quad 
1- \theta = \frac{s_2+2-s_1}{s_3-s_1}  \, ,  
$$
we have, using $\nu_n=O(\e^2)$, 
  \begin{align}
\| i_{n+1} - i_{n} \|_{\Lip, s_2+2} & 
\stackrel{\eqref{inter-Sobo-Lip}} {\lesssim_{s_1, s_3}} \| i_{n+1} - i_{n} \|_{\Lip, s_1}^{\theta} 
\| i_{n+1} - i_{n} \|_{\Lip, s_3}^{1- \teta} 
\nonumber \\ 
& \stackrel{\eqref{diff-i-n+1-in}, \eqref{diffns3}} {\lesssim_{s_1, s_3}} 
\big( \max \{\nu_n,   \e^2 \nu_n^{\frac{4}{5} - 6\sigma}\}  \big)^{\theta}   
\big( \e^2 \nu_n^{- \frac{11}{10} - 4\s} \big)^{1- \theta} 
  \nonumber \\
 & {\lesssim_{s_1, s_3}} \, 
 \big(  \e^{\frac32}  \nu_n^{\frac{1}{4}}  \big)^{\theta}   \big( \e^\frac32 \nu_n^{- \frac{11}{10} - 3 \s}
 \big)^{1- \theta} 
 \leq \e^{\frac32} \nu_n^{\frac14 -\sigma} \label{in:s2+2} 
 \end{align}
for $ s_3  $ large enough. 
The bounds  for $ {\fracchi}_{n+1}  $ follow by a telescoping argument.  
Using  the first estimate in \eqref{NM:indu-line2}$_1$ and  
\eqref{vicin-in+1-in}$_{k}$ for all $ 1 \leq k \leq n+ 1 $ we get 
$$
\begin{aligned}
\| {\fracchi}_{n+1} \|_{\Lip, s_1+2}  \leq  \| {\fracchi}_{1} \|_{\Lip, s_1+2} + 
\sum_{k= 1}^{n+1} \| {\fracchi}_{k+1} - {\fracchi}_{k} \|_{\Lip, s_1+2} 
& \lesssim_{s_1}   \e^2  + \e^2 \nu_1^{\frac45 - \sigma_\star} \\
& \leq C(s_1) \e^2 
\end{aligned}
$$
since $ \nu_1 = C(s_1) \e^2 $. This proves the first inequality in \eqref{NM:indu-line2}$_{n+1}$.
Moreover 
\begin{align} \label{n+1-s3+1} 
\| {\fracchi}_{n+1} \|_{\Lip, s_3+2} 
& \   \leq  \| \breve {\fracchi}_{n} \|_{\Lip, s_3+2} +  \| i_{n+1} - \breve \imath_n \|_{\Lip, s_3+2} \nonumber \\ 
&\! \! \! \! \! \!  \stackrel{\eqref{NM:indu-line2-trunc} \eqref{nuovas}}{\leq} \e^2 \nu_n^{- \frac45-\s} +  \e^2 \nu_n^{- \frac{11}{10} - 4 \sigma} \nonumber  \\
&  \leq  2 \e^2 \nu_n^{- \frac{11}{10} - 4 \sigma} \stackrel{\eqref{def:nun-NM}} 
\leq \e^2 \nu_{n+1}^{- \frac45} 
\end{align}
which is the third estimate in \eqref{NM:indu-line2}$_{n+1} $.  The
second estimate in \eqref{NM:indu-line2}$_{n+1} $ follows 
similarly by  \eqref{estimate-i1} 
and  \eqref{vicin-in+1-in-s2}$_k $, for all 
$ k \leq n+ 1 $.

In order to complete the proof of Theorem \ref{thm:NM-ite} it remains to prove the measure estimates 
\eqref{meas:Ln-NM}.

\begin{lemma}
The sets $ {\bf \Lambda}_n $ defined  iteratively 
in \eqref{def:Lambda-n+1} with $ {\bf \Lambda}_1 := \Lambda $, satisfy \eqref{meas:Ln-NM}. 
\end{lemma}

\begin{pf}
The estimate 
$ | {\bf \Lambda}_1 \setminus {\bf \Lambda}_{2} | \leq b(\e)$ with $ \lim_{\e \to 0} b(\e) = 0 $
follows by item \ref{mmqf2}  of Proposition \ref{prop:inv-ap-vero} with 
 $ \Lambda_{\uF} =   \Lambda_{ \breve \fracchi_1} = \Lambda $.
In order to prove \eqref{meas:Ln-NM}  for $ n \geq 3 $  
notice that, 
by the definition of $ {\bf \Lambda}_{n+1} $ (and  $ {\bf \Lambda}_{n} $) in \eqref{def:Lambda-n+1}, we have
\be\label{set-in-set}
{\bf \Lambda}_n \setminus {\bf \Lambda}_{n+1} = 
{\bf \Lambda}_{n-1} \bigcap 
 {\bf \Lambda} (\e; \eta_{n-1}, \breve \fracchi_{n-1} ) \bigcap 
 [{\bf \Lambda} (\e; \eta_{n}, \breve \fracchi_{n} ) ]^c \, , \quad n \geq 2 \, . 
\ee
In addition,  for all $ \l \in   {\bf \Lambda}_{n} $, we have, for $ \b := {\frac{3}{4}}$,  
\begin{align}  
\| \breve \fracchi_{n} - \breve \fracchi_{n-1} \|_{s_1+2} &\stackrel{\eqref{trunc-torus}}=
\| {\it \Pi}_{N_n}  \fracchi_{n} - {\it \Pi}_{N_n}  \breve \fracchi_{n-1} \|_{s_1+2} \nonumber \\
&\leq  \| {\it \Pi}_{N_n}  (\fracchi_{n} - \fracchi_{n-1}) \|_{s_1+2} +  \| {\it \Pi}_{N_n}  (\fracchi_{n-1} -\breve \fracchi_{n-1}) \|_{s_1+2} 
\nonumber \\
&\stackrel{ \eqref{smoothingS1S2}}{\leq } N_n^2 \| i_n -i_{n-1} \|_{s_1}   +  \| i_{n-1}- \breve \imath_{n-1} \|_{s_1+2}  \nonumber \\
&\stackrel{\eqref{vicin-in+1-in} , \eqref{NM:int2}}{\leq} \nu_n^{-\frac6{s_3-s_1}} \big( C(s_1) \nu_{n-1} + \e^2 \nu_{n-1}^{\frac45 - \s_\star}     \big) 
+ \e^2 \nu_{n-1}^2 \nonumber \\
& \ll    \nu_{n-1}^{\b}  \leq \e^{\frac32} \, ,  \quad \forall n \geq 2 \, ,
\label{frann-1}
\end{align}
since $ \nu_{n-1} \leq C(s_1) \e^2 $. 
Now by \eqref{def:Lambda-n+1}, we have $ \eta_n - \eta_{n-1}  = \nu_{n-1}^{\frac14} \geq \nu_{n-1}^{\frac25 \beta} $. 
Hence the estimates 
\eqref{exp-misu} and \eqref{frann-1} 
imply that the set of \eqref{set-in-set} satisfies the measure estimate
$$ 
|{\bf \Lambda}_n \setminus {\bf \Lambda}_{n+1}| \leq \nu_{n-1}^{\b \a / 3 }  \leq \nu_{n-1}^{\a /4 }  \, ,
$$ 
which proves \eqref{meas:Ln-NM}. 
\end{pf}
\\[1mm]
{\bf Proof of Theorem \ref{thm:NM}.}\index{Nash-Moser theorem}  The torus embedding  
$$
i_\infty := i_0 + ( i_1 - i_0) + \sum_{n \geq 2} (i_{n} - i_{n-1})   
$$
is defined for all $ \lambda $ in $ {\cal C}_\infty := \cap_{n \geq 1} {\bf \Lambda}_n  $, 
and, by   \eqref{estimate-i1}   and \eqref{vicin-in+1-in}-\eqref{vicin-in+1-in-s2}, 
 it is convergent in $ \| \ \|_{\Lip,s_1} $ and $ \| \ \|_{\Lip,s_2} $-norms with  
\begin{align*}
\| i_\infty - i_0 \|_{\Lip,s_1} & 
\leq C(s_1) \e^2  \, , \quad 
 \| i_\infty - i_0 \|_{\Lip,s_2}  \leq \e \, , 
\end{align*}
proving \eqref{bound-i-infty}.
By \eqref{NM:indu} we deduce that
$$
\forall \l \in  {\cal C}_\infty = \cap_{n \geq 1} {\bf \Lambda}_n \, ,  \quad 
 {\cal F} (\l; i_\infty (\l) ) = 0  \, . 
$$ 
Finally, by \eqref{meas:Ln-NM} we deduce \eqref{C-infty-AM}, i.e. 
that $ {\cal C}_\infty $ is a set of asymptotically full measure. 

\begin{remark}\label{rem:Cq}
The previous result  holds if  
the nonlinearity $ g (x, u) $   and the potential $ V(x) $ in \eqref{NLW1} are of class 
$ C^q $ for some $ q $ large enough, depending on $ s_3 $.  
\end{remark}

\section{$C^\infty$ solutions} \label{Cinfty} 

In this section we prove  the last statement of Theorem \ref{thm:NM} 
about $ C^\infty $ solutions.  

\smallskip

By Proposition \ref{crucial-s} and a simple modification of the proof of Proposition 
\ref{prop:inv-ap-vero} which substitutes any $s\geq s_3$ to $s_3$ (see Remark \ref{rem:s3s}), we obtain the following result.

\begin{proposition} \label{prop:inv-ap-vero-s}
Same assumptions as in Proposition \ref{prop:inv-ap-vero}. From
\eqref{exp-misu}, the conclusion can be modified in the following way:  
 for any  $s \geq s_3$, there  is $\wtilde{\nu} (s)$ such that the following holds.
 
 For any $ \nu \in (0, \e^\frac{3}{2}) \cap (0, \wtilde{\nu} (s))$ such that 
 $\| \uF \|_{\Lip, s+4} \leq \e  \nu^{-\frac9{10}} $, 
there exists a linear operator $ {\cal L}^{-1}_{appr} :=  {\cal L}^{-1}_{appr, \nu,s} $ such that,   for 
any function $ g : \Lambda_{\uF} \to {\mathcal H}^s \cap H_{\mathbb S}^\bot $ satisfying 
\be\label{g:small-large-IA-inf}
\| g \|_{\Lip, s_1} \leq  \e^2 \nu  \,  , \quad    \| g \|_{\Lip, s}  \leq  \e^2 \nu^{- \frac{9}{10}} \ ,
\ee
the function  $ h := {\cal L}^{-1}_{appr} g $, $ h : {\bf \Lambda} (\e; 5/6, \uF ) \to  
{\mathcal H}^{s} \cap H_{\mathbb S}^\bot $  
satisfies
\be\label{h:s1s2s3-IA-per-inf}
 \| h \|_{\Lip, s_1} \leq C(s_1) \e^2 \nu^{\frac45} \, , \quad 
 \| h \|_{\Lip, s+2} \leq   C(s) \e^2 \nu^{ - \frac{11}{10} } \, , 
\ee
and, setting  $ {\cal L}_\omega := {\cal L}_\omega(\ui ) $,  we have 
\be\label{estimate:r-final-IA-per-inf}
 \| {\cal L}_\omega h - g \|_{\Lip, s_1} \leq C(s_1) \e^2 \nu^{\frac32} \,  . 
\ee
Furthermore, 
setting $ Q' := 2(\tau'+ \varsigma s_1 ) +3 $ (where $ \varsigma = 1/ 10 $ and 
$ \tau' $ are given by Proposition \ref{propmultiscale}), for all $ g \in {\mathcal H}^{s_0 + Q'} \cap  H_{\mathbb S}^\bot $, 
\be\label{new:estimate-cal-L-1-inf}
\| {\cal L}^{-1}_{appr} g \|_{\Lip, s_0} \lesssim_{s_1} \| g \|_{\Lip, s_0 + Q'} \, .
\ee
\end{proposition}

Note that in Proposition \ref{prop:inv-ap-vero-s}  the sets ${\bf \Lambda} (\e; 5/6, \uF )$ do not depend on $s$ and are 
the same as in Proposition \ref{prop:inv-ap-vero}. 

Thanks to Proposition \ref{prop:inv-ap-vero-s} we can modify the Nash-Moser scheme
 in order to keep
along the iteration the control of higher and higher Sobolev 
norms. This new scheme relies on the
following result, which is used in the  iteration. We recall that the sequence $(\nu_n)$ is defined 
in  \eqref{def:nun-NM}, and the sequence $(\eta_n)$ in \eqref{def:Lambda-n+1}.

\begin{proposition} \label{NM:s3s}
For any $s \geq s_3 $, there is $\wtilde{\nu}'(s) >0$ with the following property. 
Assume that for some $n$ such that $\nu_n \leq \wtilde{\nu}'(s)$, there is a map
$ \fracchi_n  :  \l \mapsto   \fracchi_n  (\l)$, defined for $\l$ in some set
${\bf \Lambda}_n$, such that,
\be  \label{it-hyp-s}
\| \fracchi_n \|_{\Lip, s_1+2} \leq  C(s_1) \e^2  \  ,  \quad
\| \fracchi_n \|_{\Lip, s_2+2}  \leq  \e  \   , \quad 
\| \fracchi_n \|_{\Lip, s+2}  \leq \e^2 \nu_n^{-\frac45} \, , 
\ee
and 
$$
 \| {\cal F}(i_n )\|_{\Lip,s_1} \leq \e^2 \nu_n    \quad {\rm where} \quad
 i_n(\vphi)=(\vphi,0,0) + \fracchi_n (\vphi) \, .$$
 Let 
\be\label{def:Nns}
 N_{n,s} \in  \Big[ \nu_n^{-\frac3{s-s_1}} -1,   \nu_n^{-\frac3{s-s_1}} +1 \Big]  \quad 
 {\rm and}  \quad  \breve \fracchi_{n,s} := {\it \Pi}_{N_{n,s}} \fracchi_n \, .
 \ee
 Then there exists a map $\l \mapsto \fracchi_{n+1} (\l)$, defined for 
 $\l \in {\bf \Lambda }_n \cap {\bf \Lambda}(\e; \eta_n ,  \breve \fracchi_{n,s})$, which 
 satisfies
 \begin{align}  \label{it-concl}
&  \| i_{n+1} - i_{n} \|_{\Lip, s_1+2} \leq  C(s_1) \nu_n + \e^2 \nu_n^{\frac45-\sigma_\star}, 
 \quad
  \| i_{n+1} - i_{n} \|_{\Lip, s_2} \leq   \e^2 \nu_n^{\frac15} \, , \\
&  \label{it-concl-s}
  \| \fracchi_{n+1} \|_{\Lip, s+3} \leq  \e^2 \nu_{n+1}^{-\frac45}  
\end{align}
 and 
 $$
  \| {\cal F}(i_{n+1} )\|_{\Lip,s_1} \leq \e^2 \nu_{n+1}  \,  ,
 $$
 where $\s_\star$ is chosen as in Theorem \ref{thm:NM-ite}.
\end{proposition}

\begin{pf}
The proof follows exactly the steps of the inductive part in the proof of Theorem \ref{thm:NM-ite},
with $N_{n,s} \sim \nu_n^{-\frac3{s-s_1}}$ and with the same small exponent $\s$. We have just to observe that
in Step 1 (regularization)
\begin{eqnarray*}
\| \breve \fracchi_n \|_{\Lip , s+ \bt +7} & \stackrel{\eqref{smoothingS1S2-Lip}}  
\leq & N_{n,s}^{\bt +5} \|  \fracchi_n \|_{\Lip , s+ 2}  \\
&\stackrel{\eqref{def:Nns}, \eqref{it-hyp-s} } \leq& C(s) \nu_n^{- \frac{3(\bt +5)}{s-s_1}} \e^2 \nu_n^{-\frac45 }  \leq \e^2 \nu_n^{-\frac45 -\s} \,  ,
\end{eqnarray*}
for $3(\bt +5)/(s_3-s_1) < \s $  and $\nu_n$ small enough (depending on $s$). 
As a result  
$$
\| \breve \fracchi_n  \|_{\Lip,s +1 + \bt + 6} \leq \e^2 \nu_n^{- \frac{4}{5} -  \sigma} \, , 
$$
i.e. in the third estimate of \eqref{NM:indu-line2-trunc} in Lemma \ref{1122},
$s+1$ can be substituted to $s_3$. 
 From this point, we use the substitution
$s_3 \leadsto s+1$  in all the estimates of the induction where $s_3$ appears, 
except the second estimate of \eqref{NM:int2} and \eqref{diffns3}, where we keep $s_3$. 
Note that these last two estimates are useful only to obtain the bound
\eqref{in:s2+2}  on the norm 
$\| i_{n+1}-i_n \|_{\Lip ,s_2+2}$ by an interpolation argument.  All the estimates where we apply
the substitution $s_3 \leadsto s+1$ hold provided that $\nu_n$ is smaller than some possibly very small, but positive constant,  depending on $s$. In particular, the analogous of Proposition 
\ref{Prop:inversione} uses Proposition \ref{prop:inv-ap-vero-s}, and 
the second estimate of  \eqref{stima T 0 b}  is replaced by 
$$
\|   {\mathtt h}_{n+1} \|_{\Lip, s+3} \leq 
 \e^2 \nu_n^{- \frac{11}{10}-2\s} \, .
$$
To prove \eqref{it-concl-s}, we use  \eqref{n+1-s3+1}  with 
the substitution $s_3 \leadsto s+1$. Note that in this estimate we need only a bound
for  $\| \breve \fracchi_n \|_{\Lip , s+3}$, not for  $\|  \fracchi_n \|_{\Lip , s+3}$.
\end{pf} 

\medskip

We now consider a non-decreasing sequence $(p_n)$ of integers with the following properties:
\begin{itemize}
\item[(i)]  $p_1 = 0$ and $\lim_{n \to \infty} p_n=\infty$ ;
\item[(ii)]  $ \forall n \geq 1 \, , \  \nu_n \leq \wtilde{\nu}'(s_3+p_n)$ ;
\item[(iii)]  $\forall n \geq 1 \, , \  p_{n+1}=p_n \   {\rm or } \   p_{n+1}=p_n+1$ .
\end{itemize}
The sequence $(p_n)$ can be defined iteratively in the following way:  $p_1 := 0 $, so that
property (ii) is satisfied for $n=1$, provided $\e$ is small enough (the smallness condition depending on $s_3$). 
Once $p_n$ is defined (satisfying (ii)), we choose
$$
p_{n+1} := \begin{cases}
p_n +1 \quad \ {\rm if} \   \nu_{n+1} \leq \wtilde{\nu}'(s_3 + p_n +1)\\
p_n \quad \qquad  {\rm otherwise} \, .  
\end{cases}  
$$
Then property (ii) is satisfied at step $n+1$ in both cases, because $(\nu_n)$ is decreasing.  The sequence $(p_n)$ 
satisfies (i) because $\lim_{n \to \infty} \nu_n=0$ (so that the sequence cannot be stationary) and (iii) by definition.

Starting, as in Section \ref{sec:NM-ite}, with the torus $i_1(\vphi)$ defined in Lemma \ref{NM:step1} for all 
$\l \in \Lambda={\bf \Lambda}_1 $ and using repeatedly Proposition \ref{NM:s3s}, we obtain the following theorem. 


\begin{theorem}\label{thm:NM-ite-inf} {\bf (Nash-Moser $C^\infty$)}
Let  $ \bar \om_\e \in \R^\es $ be  $(\gamma_1 , \tau_1) $-Diophantine and satisfy 
property $ {\bf (NR)}_{\gamma_1, \tau_1} $ in Definition \ref{NRgamtau} with $ \g_1, \t_1 $ fixed in 
\eqref{def:tau1}. Assume  \eqref{choice s2 s3} and 
$ s_2 - s_1 \geq \underline{\tau} + 2 $, $ s_1 \geq s_0+2+\underline{\tau} +Q'  $, 
where $ \underline{\tau} $ is the loss of derivatives  defined in Proposition
\ref{Prop:sec6} and $ Q' := 2(\tau'+ \varsigma s_1 ) +3 $ is defined in   Proposition \ref{lem:app-inv-LD}.

Then, for all $ 0 < \e \leq \e_0 $ small enough, 
 for all  $ n \geq 1 $, there exist 
\begin{enumerate}
\item a subset $ {\bf \Lambda}_n \subseteq {\bf \Lambda}_{n-1}  $, $ {\bf \Lambda}_1 := {\bf \Lambda}_0 := \Lambda $,
satisfying 
\be\label{meas:Ln-NM-inf}
\begin{aligned}
& | {\bf \Lambda}_1 \setminus {\bf \Lambda}_{2} | \leq b(\e) \quad {\rm with} \quad \lim_{\e \to 0} b(\e) = 0 \, , \\
& 
| {\bf \Lambda}_{n-1} \setminus {\bf \Lambda}_{n}| \leq \nu_{n-2}^{ \alpha_*} \, , \ \forall n \geq 3 \, , 
\end{aligned}
\ee
where $ \a_* =  \alpha / 4 $ and $ \a > 0  $ is the exponent in \eqref{exp-misu}, 
\item  a torus $ i_n (\vphi) = (\vphi,0,0) + \fracchi_n (\vphi) $,
 defined for all $ \l \in {\bf \Lambda}_n  $, 
satisfying
\begin{align}
& 
\| \fracchi_n  \|_{\Lip,s_1} \leq 
 C(s_1) \e^{2}  \, , \quad  
\| \fracchi_n  \|_{\Lip,s_2+2} \leq \e  \, , \quad 
\| \fracchi_n  \|_{\Lip,s_3 + p_n  +2} \leq 
\e^2 \nu_n^{- \frac45} \, , \label{NM:indu-line2-inf}  \\
&  
\| i_n - i_{n-1}  \|_{\Lip,s_1+2} \leq C(s_1) \nu_{n-1} + \e^2 \nu_{n-1}^{\frac{4}{5}-\sigma_\star} \, , \ n \geq 2 \, ,  \label{vicin-in+1-in-inf} \, \\
&  
\| i_n - i_{n-1}  \|_{\Lip,s_2} \leq  \e \nu_{n-1}^{\frac{1}{5}} \, , \quad n \geq 2 \, ,  \label{vicin-in+1-in-s2-inf} \,
\end{align}
and
\be\label{NM:indu-inf} 
 \| {\cal F}(i_n )\|_{\Lip,s_1} \leq \e^2 \nu_n  \,  .  
\ee
\end{enumerate}
\end{theorem}

As said at the end of Section \ref{sec:NM-ite}, \eqref{vicin-in+1-in-inf} and \eqref{vicin-in+1-in-s2-inf}  imply
that the sequence $(i_n)$ converges in $\|  \  \|_{\Lip , s_1}$  and $\|  \  \|_{\Lip , s_2}$ norms 
to a map $i_\infty : \l \mapsto i_\infty (\l)$, defined on 
$$
{\cal C}_\infty = \cap_{n \geq 1} {\bf \Lambda}_n \, ,
$$
which satisfies 
$$
\forall \l \in  {\cal C}_\infty \, ,  \quad 
 {\cal F} (\l; i_\infty (\l) ) = 0  \, . 
$$
There remains to justify that 
$$
\|  \fracchi_\infty \|_{\Lip , s} =\| i_\infty - i_0 \|_{\Lip , s} < \infty \, , \quad \forall s > s_2 \, . 
$$
For a given $s > s_2$, let us fix $\ov{n}$ large enough, so that
\be  \label{sbars}
s \leq \frac34 s_1 + \frac14 \ov{s}   \qquad {\rm with }  \ \   \ov{s}:= s_3 + p_{\ov{n}} +2 \, .
\ee
Then, by  \eqref{NM:indu-line2-inf}, for all $k \geq \ov{n}+1$, 
\be  \label{diff-s}
\|\fracchi_k - \fracchi_{k-1}\|_{\Lip , \ov{s}} \leq \|\fracchi_k \|_{\Lip , \ov{s}} + \|\fracchi_{k-1}\|_{\Lip , \ov{s}} \leq 2 \e^2 \nu_k^{-\frac45}  \, .
\ee
Hence, by  interpolation,
\begin{align}
\|\fracchi_k - \fracchi_{k-1}\|_{\Lip , {s}} & 
\stackrel{\eqref{sbars} } \leq \|\fracchi_k - \fracchi_{k-1}\|_{\Lip , \frac{3s_1+\ov{s}}{4}} 
\nonumber \\
&\stackrel{\eqref{inter-Sobo-Lip}} 
\leq C(\ov{s}) \|\fracchi_k - \fracchi_{k-1}\|_{\Lip , s_1}^{\frac34}  \|\fracchi_k - \fracchi_{k-1}\|_{\Lip , \ov{s}}^{\frac14} \nonumber  \\
& \stackrel{\eqref{vicin-in+1-in-inf}, \eqref{diff-s}}{\leq } C(\ov{s}) 
 \nu_{k-1}^{\frac3{5}- \s_\star}  \nu_k^{-\frac1{5}}  \leq 
 \nu_{k-1}^{\frac3{10} - \s_\star}  \label{ffina}
\end{align} 
using that $\nu_k=\nu_{k-1}^{\mathfrak q} $, with $ {\mathfrak q} <3/2$, see 
\eqref{def:nun-NM}.  Moreover,  recalling that
$ \ov{s}:= s_3 + p_{\ov{n}} +2 $ (see \eqref{sbars}), we deduce by \eqref{NM:indu-line2-inf} that 
\be\label{bounfIn}
\| \fracchi_{\ov{n}}\|_{\Lip , \ov{s}} < \infty \,	 .
\ee
In conclusion,  \eqref{bounfIn} and \eqref{ffina}, imply, since $ (3/10) -\s_\star>0$, that
$$
\| \fracchi_{\ov{n}}\|_{\Lip , {s}}  + \sum_{k \geq \ov{n} +1}  \|\fracchi_k - \fracchi_{k-1}\|_{\Lip , {s}} < \infty \, ,
$$
and the sequence $(\fracchi_n)_{n \geq \ov{n}}$ converges in $ H^s_\vphi \times  H^s_\vphi \times 
( {\mathcal H}^s  \cap H_{\mathbb S}^\bot ) $ to  $ \fracchi_\infty=i_\infty - i_0$. 
Since $s$ is arbitrary, we conclude that $i_\infty (\l)$ is $C^\infty$ for any $\lambda \in {\cal C}_\infty$.

\chapter{Genericity of the assumptions}
\label{section:gener1} 

The aim of this chapter is to prove the genericity result\index{Genericity result} stated in Theorem \ref{Prop:genericity}. 

\section{Genericity of non-resonance and non-degeneracy  conditions}

We fix $s>d/2$, so that we have the compact embedding $H^s (\T^d) \hookrightarrow C^0 (\T^d)$. 
We denote by $B(w,r)$ the open ball of center $ w $ and radius $r$ in $H^s (\T^d)$.

Recalling Definition \ref{Def:dense-open} of $C^\infty$-dense open sets,  
it is straightforward to check the following lemma. 
\begin{lemma} The following properties hold:
\begin{enumerate}
\item a finite intersection of $C^\infty$-dense open subsets of ${\cal U}$ is $C^\infty$-dense open
in ${\cal U}$;
\item
a countable intersection of $C^\infty$-dense open subsets of ${\cal U}$ is $C^\infty$-dense in ${\cal U}$;
\item 
given three open subsets ${\cal W} \subset {\cal V} \subset {\cal U}$ of $H^s(\T^d)$
(resp.  $H^s(\T^d) \times H^s(\T^d)$), if ${\cal W}$ is $C^\infty$-dense in ${\cal V}$ and ${\cal V}$
is $C^\infty$-dense in ${\cal U}$, then ${\cal W}$ is $C^\infty$-dense in ${\cal U}$. 
\end{enumerate}
\end{lemma}

Moreover we have the following useful result.

\begin{lemma} \label{analyt}
Let $U$ be a connected open subset of $H^s (\T^d)$ and let $ f : U \to \R$ be
a real analytic function. If $f \not \equiv 0$, then 
$$
Z(f)^c:= \big\{ w \in U \ : \  f(w)\neq 0 \big\} 
$$ 
is a $C^\infty$-dense open subset of $U$.
\end{lemma}  

\begin{pf}   
Since $f : U \to \R$ is continuous, $Z(f)^c$ is an open subset of $U$. 
Arguing by contradiction, we assume that
$Z(f)^c$ is not $C^\infty$-dense in $ U $. 
Then there are $w_0 \in U$, $s' \geq s$ and $\ep>0$ such that: for all 
$h \in C^\infty (\T^d)$ satisfying $ \| h \|_{H^{s'}} < \ep$, 
we have that $f(w_0+h)=0 $. 
Let $\rho>0$ be such that the ball $B(w_0,\rho) \subset U$. 
We claim that $f$ vanishes on $B(w_0, \rho) \cap C^\infty ( \T^d ) $.
Indeed, if $h \in C^\infty(\T^d)$ satisfies 
$\|h\|_{H^s} < \rho$ 
we have the segment  $[w_0, w_0+h] \subset U$, and the map $ \varphi : t\mapsto f(w_0+th)$ is real
 analytic on an open interval of $ \R $
which contains $[0,1]$. Moreover $\varphi (t) $ vanishes on the whole interval
$ | t | <  \ep \| h \|_{H^{s'}}^{-1} $. 
Hence $\varphi $ vanishes everywhere, and $f(w_0+h)=\varphi(1)=0$. 

Now, by the fact that $C^\infty (\T^d)$ is a dense subset of $H^s (\T^d)$, and
$  f $ is continuous,
we  conclude the $f$ vanishes on the whole ball $B(w_0, \rho)$. 

Let 
$$
V := \big\{ w \in U \ : \ f \ \hbox{vanishes on some open neighborhood of} \ w \big\} \, .
$$ 
From the previous argument, 
$V$ is not empty, and it is by definition an open subset of $U$. Let us prove that it is closed in $U$ too.
Assume that the sequence $(w_n)$ $H^s$-converges to some $w$ and satisfies: $w_n \in V$ for all $n$. 
Then there is $r>0$ such that, for $n$ large enough, the ball $B(w_n,r)$ contains $w$ and is 
included in $U$. Using the same argument as before (with $w_n$ instead of $w_0$), we can conclude
that $f$ vanishes on the whole ball $B(w_n,r)$, hence $w \in V$. 

Since $U$ is connected, we finally obtain $V=U$, which contradicts the hypothesis $f \not \equiv 0 $.
\end{pf}

The proof  of Theorem \ref{Prop:genericity} uses results and arguments provided by Kappeler-Kuksin \cite{KaKu} that we recall below.  
We first introduce some preliminary information.  For a real valued 
potential $V \in H^s(\T^d)$, we denote by $(\lambda_j(V))_{j \in \N}$ the sequence of the eigenvalues\index{Eigenvalues of Schr\"odinger operator}  of 
the Schr\"odinger operator $-\Delta + V(x) $,\index{Schr\"odinger operator} written
in increasing order and counted with multiplicity
$$
\lambda_0(V) < \lambda_1 (V) \leq \lambda_2 (V) \leq \ldots \, . 
$$ 
These eigenvalues are Lipschitz-continuous functions of the potential, namely
\be \label{Lipdep}
|\lambda_j(V_1) - \lambda_j (V_2) | \leq \| V_1-V_2 \|_{L^\infty(\T^d)} \lesssim \| V_1-V_2 \|_{H^s(\T^d)} \, . 
\ee  
{\bf Step 1.} {\bf The construction of  Kappeler-Kuksin \cite{KaKu}.} For $J \subset  \N$, define
the set of potentials $ V := V(x) $ in $ H^s (\T^d) $ such that the eigenvalues $ \lambda_j (V) $, $ j \in J $, 
of $ - \Delta + V(x) $ are simple, i.e. 
$$
E_{J} := \Big\{ V  \in H^s(\T^d)  \, : \, \lambda_j (V) \ \text{is a simple eigenvalue of }
- \Delta + V(x) \, ,  \ \forall j \in 
J \Big\} \, . 
$$
The set $E_{[ \! [ 0,N ] \! ]}$ will be simply denoted by $ E_N$. 
Since the eigenvalues $ \lambda_j(V) $ of $ - \Delta + V(x) $ are simple on $ E_J $, it turns out that 
each function
$$ 
\lambda_j : E_J  \to \R \, , \quad  j \in {J} \, , 
$$ 
is {\it real analytic}. 
Moreover the corresponding eigenfunctions $ \Psi_j := \Psi_j (V) $, normalized with   $ \| \Psi_j \|_{L^2} = 1 $,  
can  locally be expressed as real analytic functions of the potential $ V \in E_J $.

By  \cite{KaKu}, Lemma 2.2, we have that,  for any $J \subset \N$ finite,
\begin{enumerate}
\item   $ E_{J}  $ is an  {\it open}, {\it dense} and {\it connected} subset of $H^s(\T^d)$,
\item 
$H^s(\T^d) \backslash E_J $ is a real analytic variety. 
This implies that for all $ V \in H^s(\T^d) \backslash E_J$, there are $r>0$
and real analytic functions $f_1, \ldots , f_s$ on the open ball $B(V,r)$ such that
\be\label{real-anal-var}
(i) \quad B(V,r) \backslash E_J \subset \bigcup_{i=1}^s f_i^{-1}(0)  \qquad  \quad
(ii)  \quad \forall i \in [ \![   1,s ] \! ] \, , \  f_{i|B(V,r)} \not \equiv 0 \, . 
\ee
\end{enumerate}
Recalling Definition \ref{Def:dense-open} we prove this further lemma. 
\begin{lemma} \label{ej}
Let $ {\cal P} $ be the set defined in \eqref{pos:poten}.
The subset $E_J \cap {\cal P} \subset H^s(\T^d) $ is connected and $C^\infty$-dense open in ${\cal P}$.
\end{lemma}

\begin{pf}
We first prove that  $E_J \cap {\cal P}$ is $C^\infty$-dense open in ${\cal P}$.
By 	\eqref{real-anal-var} 
\be\label{triv:inclu}
\bigcap_{i=1}^s B(V,r) \backslash f_i^{-1}(0) \subset B(V,r) \cap E_J \, .
\ee
By Lemma \ref{analyt} 
each set $ B(V,r) \backslash f_i^{-1}(0)$, $ i =1, \ldots, s $,  is $C^\infty$-dense open in 
$B(V,r)$ (since $ f_{i|B(V,r)} \not \equiv 0$), as well as their intersection, and therefore \eqref{triv:inclu} implies that 
$ B(V,r)  \cap E_J $
is $C^\infty$-dense in $B(V,r)$. 
Thus  $E_J$ is $C^\infty$-dense open in $H^s(\T^d)$. 
Finally, since ${\cal P} $ is an open
subset of $H^s(\T^d)$, the set $E_J \cap {\cal P}$ is $C^\infty$-dense open in ${\cal P}$.

There remains to 
justify that $E_J \cap {\cal P}$ is connected. Let $V_0, V_1 \in E_J \cap {\cal P}$. Since $E_J$ 
is an open connected subset of $H^s(\T^d)$, it is arcwise connected: there is a continuous path
$ \gamma : [0,1]  \to E_J $ such that $\gamma (0)=V_0$, $\gamma (1)=V_1$. 
Notice that, if $V \in E_J$, then, for all $m\in \R $, the potential $V+m \in E_J$. 
Let $\lambda_0 (t)$ be the smallest eigenvalue of $-\Delta + \gamma (t)$. The map 
$t \mapsto \lambda_0(t)$ is continuous. Since $V_0, V_1 \in
{\cal P}$, we have $\l_0(0)>0$, $\l_0(1)>0$. Choose any continuous map $ \mu : [0,1] \to ]0, + \infty[$ 
such that $\mu (0)= \lambda_0(0)$ and $\mu(1)=\l_0(1)$ and define $m : [0,1] \to \R$ by $m(t) := \mu (t) - \l_0(t)$. 
Then $\gamma + m$ is a continuous path in $E_J \cap {\cal P}$ connecting $V_0$ and $V_1$. In 
conclusion, $E_J \cap {\cal P}$ is arcwise connected.
\end{pf}

We have the following result, which is 
 Lemma 2.3 in \cite{KaKu}, with some new estimates 
on the eigenfunctions, i.e. items ($iii$)-($iv$).

\begin{lemma} \label{convL4}
Fix $J \subset \N$ , $J$ finite. There is a sequence $(q_n)_{n \in \N} $  of $C^\infty$ positive potentials with the following properties: 
\begin{itemize}
\item[(i)] $ \forall n \in \N $,  the potential $ q_n $ is in $ E_J$. More precisely, for each $j \in J$, the sequence of eigenvalues 
$(\l_{n,j})_n := (\l_j (q_n))_n $ converges to some $\l_j>0$, with
$\l_j \neq \l_k$ if $j,k \in J$, $j\neq k$. 
\item[(ii)] Let $\pm \Psi_{n,j}$ be the eigenfunctions
of $ -\Delta + q_n (x) $ with   $ \| \Psi_{n,j} \|_{L^2} = 1$. 
For each $j \in J$, the sequence $ (\Psi_{n,j})_n \to \Psi_j$ weakly in $H^1 (\T^d)$, hence strongly in $L^2(\T^d)$, 
where the functions $\Psi_j \in L^\infty (\T^d) $ 
have disjoint essential supports. 
\item[(iii)] For each $j \in J$, the sequence $  (\Psi_{n,j})_n \to \Psi_j$  strongly in $L^q(\T^d)$ for any $q \geq 2$.
\item[(iv)]  For any $\rho \in L^\infty(\T^d)$,  
$ \forall j,k \in J $, 
$$ 
\lim_{n \to \infty} \dps \int_{\T^d} \rho (x) \Psi^2_{n,j} (x) \Psi^2_{n,k} (x) \, dx
= \d^j_k \int_{\T^d} \rho(x) \Psi^4_j (x) \, dx
$$
 where $ \d^j_k $ is the Kronecker delta with values $ \d^j_k := 0 $, if $ k \neq j $, and $  \d^j_j  := 1 $. 
\end{itemize}
\end{lemma}
\begin{pf}
Let $M\in \N$ be such that $J \subset [ \! [ 0,M]\! ]$. It is enough to prove the lemma for 
$J = [ \! [ 0,M]\! ]$. 
We recall the construction  in Lemma 2.3 of \cite{KaKu}.
Choose disjoints open balls 
$B_j := B(x_j, r_j)$, $ 0 \leq j \leq M$, of decreasing radii $r_0 > \ldots > r_M $ 
in such a way that, denoting  by $\l_j $ 
the smallest Dirichlet eigenvalue of  $ -\Delta $ on $ B_j $, it results
$$ 
\l_0 < \ldots  < \l_M < \l_0^{(2)} 
$$ 
where $ \l_0^{(2)} $ is the second Dirichlet 
eigenvalue of $ -\Delta $ on $ B_1 $.
Define a sequence of $ C^\infty $ positive potentials such that 
$$
q_n(x) :=
\begin{cases}
n \, ,  \quad
\forall x \in \T^d \backslash \bigcup_{i=0}^M B_i \, ,  \cr
0 \, , \quad  \forall 
x \in B(x_i,r_i-\ep_n) \, ,  
\end{cases}
\quad {\rm with} \quad \lim_{n \to \infty} \ep_n =0  \, , \quad 0 \leq q_n\leq n \, .
$$
It is proved in  \cite{KaKu}, Lemma 2.3, that  properties $(i)-(ii)$ hold 
with functions $ \Psi_j $, $ j \in [\![ 0,M]\!] $, in $ H^1 (\T^d) $, 
satisfying  ${\rm suppess} (\Psi_j) = B_j$ and such that $\Psi_{j|B_j} \in H^1_0(B_j)$ is
an eigenfunction of $ - \Delta $ associated to the eigenvalue $\lambda_j $.
Such functions $ \Psi_j $ are in $ L^\infty (\T^d) $ because they are smooth in each ball $ B_j $ and vanish outside.

To prove ($iii$), it is sufficient to prove the $ L^p $ bounds
\be\label{bound-eige}
\forall j \in [\! [ 0,M ] \! ] \, , \ \forall p\geq 2 \, ,  \quad 
\sup_n \| \Psi_{n,j} \|_{L^p(\T^d)} \leq C_{j,p} < + \infty \, . 
\ee
Indeed, since $(\Psi_{n,j})_n $ converges to $\Psi_j$ 
in $L^2(\T^d)$ by item ($ii$), the bound \eqref{bound-eige}, the fact that $ \Psi_j \in L^\infty (\T^d) $, and 
 H\"older inequality, imply that the sequence $(\Psi_{n,j})_n$ converges to $\Psi_j$ 
in $L^q(\T^d)$ for any $q$. 

We fix $ 0 \leq j \leq M$ and, for simplicity, 
we write $ \Psi_{n,j} = \Psi_n $  in what follows. 
 To  prove \eqref{bound-eige} we perform a bootstrap argument 
for the  $ L^p $ norms of the solutions of the elliptic eigenvalue equation 
\be\label{ei-eq}
-\Delta \Psi_n (x) + q_n (x)  \Psi_n (x) = \l_n \Psi_n (x)  \, .
\ee
Remark that the eigenfunctions $ \Psi_n  $ are in $ C^\infty (\T^d) $, but we shall not perform Schauder estimates
because we want  bounds independent of the potentials 
$ q_n (x) $ which are unbounded in Sobolev spaces.  
We multiply  \eqref{ei-eq} 
by  $|\Psi_n|^{r-2} \Psi_n$, $r\geq 2 $,  and integrate by parts on $\T^d$, obtaining
\be\label{int-par-int}
\int_{\T^d} \nabla \Psi_n (x) \cdot \nabla (|\Psi_n (x) |^{r-2} \Psi_n (x)) + 
q_n (x) |\Psi_n (x) |^r \, dx= \l_n \int_{\T^d}  |\Psi_n (x)|^r \, dx  \, .
\ee
Now 
$$ 
\begin{aligned}
 \nabla \Psi_n \cdot \nabla (|\Psi_n|^{r-2} \Psi_n)=(r-1) |\Psi_n|^{r-2} |\nabla \Psi_n|^2 
& =(r-1) | |\Psi_n|^{\frac{r}{2}-1} \nabla \Psi_n|^2 \\
& =K_r |\nabla z_n|^2   
\end{aligned}
$$
where 
$$
K_r := 4r^{-2} (r-1) \, , \qquad 
z_n := z_n (x) := |\Psi_n (x)|^{\frac{r}{2}-1}  \Psi_n (x) \, .
$$
 Hence, by 
\eqref{int-par-int} and since $q_n \geq 0$,
$$
K_r  \int_{\T^d} |\nabla z_n(x)|^2 \, dx   \leq      \int_{\T^d} K_r |\nabla z_n(x)|^2 + q_n (x) |z_n (x)|^2 \leq \l_n  
\int_{\T^d}  |z_n(x)|^2 \, dx \,  ,
$$
which gives 
$$ 
\| z_n \|_{H^1} \leq (K_r^{-1} \l_n +1)^{1/2} \| z_n \|_{L^2} \leq C_r \| z_n \|_{L^2}
$$ 
because the sequence $(\l_n)_n $ 
is bounded, see item ($i$). The continuous Sobolev embedding $H^1(\T^d) \hookrightarrow L^{d_*} (\T^d)$ with 
$ d_*=\dps \frac{2d}{d-2}$, implies that $ \| z_n \|_{L^{d_*}} \leq C_r' \| z_n \|_{L^2}$, i.e.
\be \label{iterLr}
\| \Psi_n \|_{L^{\frac{rd}{d-2}}} \leq C_r \| \Psi_n \|_{L^{r}} \, . 
\ee
Iterating \eqref{iterLr} (and starting from $r=2$), we obtain that the sequence $( \| \Psi_n \|_{L^p})_n $
is bounded for any $p$. Note that if $d\leq 2$, we obtain \eqref{bound-eige} in one step only since
in this case $H^1$ is continuously embedded in $L^p$ for any $p$. 

$(iv)$ is a straightforward consequence of the convergence of the sequence 
$(\Psi_{n,j})_n $ to $\Psi_j $ in $L^4(\T^d)$ 
for all $ j \geq M$ and of the fact that the functions $\Psi_j $ ($0\leq j \leq M$) have disjoint essential 
supports.
\end{pf}

As a corollary we deduce the following lemma. 

\begin{lemma}\label{KaKu}
{\bf (Lemma 2.3 in \cite{KaKu})}
There is a $ C^\infty $ potential $ q(x) $  such that all the eigenvalues $ \l_j (q) $, $ j \in {\mathbb S} $, are simple
(therefore $ q $ is in $ E_{\mathbb S}  $), and the 
corresponding $ L^2$-normalized eigenfunctions $ \Psi_j (q) $, $ j \in {\mathbb S} $,
have the property that $ (\Psi_j^2(q))_{j \in {\mathbb S}} $ are linearly independent.
\end{lemma} 

Consider the real analytic map
$$
\Lambda: E_{\mathbb S} \to \R^\es \, , \quad \Lambda (V) := ( \l_j(V))_{j \in {\mathbb S}} \, .  
$$

\begin{lemma} \label{cor:N1S}
There is a $|\mathbb S|$-dimensional
linear subspace $E$ of $C^\infty(\T^d)$ such that 
\be\label{def:cal-N1}
{\cal N}^{(1)}_{\mathbb S} := \big\{ V\in E_{\mathbb S} \  :   \   d\Lambda (V)_{|E} \ \hbox{is an isomorphism} \big\}
\ee
is a  $C^\infty$-dense open subset of $E_{\mathbb S}$, 
thus of $H^s(\T^d)$. 
\end{lemma}

\begin{pf}
We follow \cite{KaKu}.  For any $ \widehat V \in H^s (\T^d) $ we have that the differential
$$
d\Lambda(q) [ \widehat{V}] = \big( ( \widehat{V} , \Psi_j^2 )_{L^2} \big)_{j \in {\mathbb S}} \, . 
$$
Since the $\Psi_j^2(q) $ are linearly independent by Lemma \ref{KaKu},  $ d\Lambda(q) $ is onto, and 
there is a $\mathbb |{\mathbb S}|$-dimensional linear subspace $E$ of $C^\infty(\T^d)$
such that $d\Lambda (q)_{|E}$ is an isomorphism. Let $(g_1, \ldots , g_\es)$ be a basis of $E$. 
For any $V\in E_{\mathbb S}$, denote by $A_V$  the $|{\mathbb S}| \times |{\mathbb S}|$-matrix
whose columns are given by $d\Lambda (V) [g_j] $, so that 
$$
{\cal N}^{(1)}_{\mathbb S} = \big\{ V\in E_{\mathbb S} \  :   \   d\Lambda (V)_{|E} \ \hbox{is an isomorphism} \big\}
= \big\{ V\in E_{\mathbb S} \ : \ \det A_V \neq 0 \big\} \, . 
$$ 
 $ E_{\mathbb S}$ is a connected open subset of $H^s(\T^d)$ and the map $V \mapsto \det A_V$ is real analytic
on $E_{\mathbb S}$ and does not vanish at $q$. Hence  Lemma \ref{analyt} implies that ${\cal N}^{(1)}_{\mathbb S}$  
is a \dcopen subset of $E_{\mathbb S}$.
\end{pf}

For any $V \in {\cal P}$, defined in \eqref{pos:poten},  all the eigenvalues $ \lambda_j  $ of $-\Delta +V (x) $ are strictly positive, and therefore we deduce the following lemma:  

\begin{lemma}\label{lem:isoN1}
The map 
\be\label{def:Upsilon}
\bar \mu : E_{\mathbb S} \cap {\cal P} \to \R^\es \, , \quad 
\bar \mu (V ) := ( \mu_j (V))_{j \in {\mathbb S}} = \big( {\l_j^{\frac12} (V)} \big)_{j \in {\mathbb S}}\, , 
\ee
is  real analytic and, for any $V \in {\cal N}^{(1)}_{\mathbb S}$, the differential  $d \bar \mu (V)_{|E}$ is an isomorphism
onto $\R^{|\mathbb S|}$. 
\end{lemma}

\begin{Remark} \label{extNJ}
More generally, for any finite $J \subset \N$, there is a \dcopen subset ${\cal N}_J^{(1)}$ of $H^s(\T^d)$ such that,
for all $V \in {\cal N}_J^{(1)} \cap {\cal P}$, the linear map $ (d\mu_j (V))_{j \in J} : H^s(\T^d) \to \R^{|J|}$ is onto. 
\end{Remark}

The $ |\mathbb S|$-dimensional
linear subspace $E$ of $C^\infty(\T^d)$ defined in Lemma \ref{cor:N1S}
is the same subspace that appears in the statement of Theorem \ref{Prop:genericity}.

\medskip

\noindent
{\bf Step 2.} {\bf Genericity  of the twist condition \eqref{A twist}} 
\\[1mm]
Our aim is to prove that the twist matrix $\Ab$ defined  in \eqref{def:AB} is 
invertible for $(V,a)$ belonging to some $C^\infty$-dense open subset of $E_{\mathbb S} \times H^s (\T^d)$.
Note that $\Ab$ is invertible if and only if  
 $ \det G \neq 0  $ where the matrix $ G   := (G^j_k(V, a))_{j,k\in {\mathbb S}} $ is defined in \eqref{def G}. 
The matrix $ G $ depends 
linearly on 
the coeffient $ a(x) $,  
and nonlinearly on 
the potential $ V(x) $, through 
the eigenfunctions $ \Psi_j := \Psi_j (V) $ defined in \eqref{auto-funzioni}.  
By previous considerations, the functions $\Psi_j^2 (x) $, $ j \in {\mathbb S}$, 
where the eigenfunctions $\Psi_j (x)$  are normalized by the condition  $ \| \Psi_j \|_{L^2} = 1 $, depend analytically 
on 
the potential $ V \in E_{\mathbb S} $, and so each 
map 
\be\label{anal:entries}
\begin{aligned}
& \qquad \qquad G_k^j  : H^s (\T^d) \times H^s (\T^d) \to \R \, , \\
& (V, a) \mapsto G_k^j (V,a)=\dps \frac{3}{4} (2- \d_k^j) ( \Psi_j^2, a \Psi_k^2)_{L^2} \, , \ \  \forall j, k \in {\mathbb S} \, , 
\end{aligned}
\ee
is real analytic on $E_{\mathbb S} \times H^s (\T^d)$, as well as 
the map $(V,a) \mapsto \det G (V,a) $. 

\begin{lemma}\label{lem:gene-a}
The set 
$$ 
{\cal N}^{(2)}_{\mathbb S} := \{ V \in E_{\mathbb S} \ : \  \det G (V,1) \neq  0 \} 
$$
is a $C^\infty$-dense open subset of $E_{\mathbb S} $.
\end{lemma}

\begin{pf}
Consider the sequence $(q_n)_n $ of $C^\infty$ potentials provided by
Lemma \ref{convL4} with $J={\mathbb S} $. By property $(iv)$ of that lemma, taking the limit
for $ n \to \infty $ in
 \eqref{anal:entries}, we get 
$$
\lim_{n \to \infty} G_k^j (q_n , 1)=\frac{3}{4} \d^{j}_k (\Psi_j^2, \Psi_k^2)_{L^2}
$$
and therefore
\be\label{det:G}
\lim_{n \to \infty} \det { G} (q_n,1) = (3/4)^{\es} \prod_{j \in {\mathbb S}} \int_{\T^d} \Psi_{j}^4 (x) \, dx =: \rho \neq 0 \, . 
\ee
In particular,  there is a potential $ q_n (x)  $ in $ E_{\mathbb  S} $ such that
$ \det G(q_n,1) \not\equiv 0$.  
Since the map $V \mapsto \det G(V,1)$ is real analytic on the open and connected  subset 
$E_{\mathbb S} $ of $ H^s (\T^d) $, 
Lemma \ref{analyt} implies that
$ {\cal N}_{\mathbb S}^{(2)}  $ 
is $C^\infty$-dense open in $ E_{\mathbb S}$. 
\end{pf}

We deduce the following corollary. 

\begin{corollary} \label{twist2}
The set  
\be\label{def:G2}
{\cal G}^{(2)} := \big\{  (V,a) \in  E_{\mathbb S} \times H^s (\T^d)  \  : \ 
\det G (V,a) \neq  0 \big\} 
\ee
is a $C^\infty$-dense open subset of $ E_{\mathbb S} \times H^{s} (\T^d)$, thus of  $ H^s (\T^d) \times H^{s} (\T^d) $. 
\end{corollary}

\begin{pf}
By Lemma \ref{lem:gene-a}, for each potential $ V \in  {\cal N}^{(2)}_{\mathbb S} $, 
we have that $ \det G(V,1) \neq 0 $ and, 
since the function  $ (V, a) \mapsto \det G(V,a)$ is real analytic 
on the open and connected  subset  
$ E_{\mathbb S} \times H^s (\T^d) $,
Lemma \ref{analyt} implies that 
 $ {\cal G}^{(2)}   $ is a  $C^\infty$-dense open   subset   of $ E_{\mathbb S} \times H^{s} (\T^d) $. 
\end{pf}

\begin{remark}\label{cor:twist}
With similar arguments we deduce that, 
for each potential $ V \in  {\cal N}_{\mathbb S}^{(2)} $, the set
$ 
\big\{ a \in     H^{s} (\T^d) \ : \  \det G (V,a) \neq 0  \big\} $
is a  $C^\infty$-dense open   subset   of $ H^{s} (\T^d) $.   
\end{remark}

\noindent
{\bf Step 3.}  {\bf Genericity  of the non-degeneracy conditions \eqref{non-reso}-\eqref{non-reso1}} 
\\[1mm]
\indent
Let $ M \in \N $ such that $ {\mathbb S} \cup {\mathbb F}   \subset [\![ 0, M ]\! ] $. We define  
\be\label{def:GM}
{\cal G}_M := \Big\{ (V,a) \in  (E_{\mathbb S}\cap {\cal P} ) \times  H^s(\T^d) \  :  \  \text{the following conditions hold} 
\ee
\begin{enumerate}
\item $\big( \det \Ab \, \mu_j-[ \Bb \Ab^\sharp \bar{\mu}]_j \big) - 
 \big( \det \Ab \, \mu_k -[ \Bb \Ab^\sharp \bar{\mu}]_k\big) \neq 0 ,  \forall j,k \in {\mathbb S}^c, j,k \leq M, 
j \neq k $, 
\item $\big( \det \Ab \, \mu_j-[ \Bb \Ab^\sharp \bar{\mu}]_j \big) + 
 \big( \det \Ab \, \mu_k -[ \Bb \Ab^\sharp \bar{\mu}]_k\big) \neq 0  \,   \forall j,k \in {\mathbb S}^c \, , \,  j,k \leq M $ \Big\}
\end{enumerate}
where 
$ \Ab := \Ab (V,a) $, $ \Bb := \Bb (V, a) $ are the Birkhoff matrices   introduced
in \eqref{def:AB} 
and where $ \Ab^\sharp $ denotes the comatrix of $ \Ab $. 
Notice that 
\be\label{matrix-adjugate}
{\rm for \ any} \    (V,a) \in {\cal G}^{(2)} \, , \ \text{defined in \eqref{def:G2}} \, , \   \Ab^{-1} = \Ab^\sharp / \det \Ab 
\ee 
so that conditions 1, 2, above imply the non-degeneracy conditions 
\eqref{non-reso}-\eqref{non-reso1} for ${\mathbb S} \cup {\mathbb F} \subset [\![ 0,M ]\!]$.

\begin{proposition} \label{127}
The set ${\cal G}_M$ defined by \eqref{def:GM} is a \dcopen subset of ${\cal P} \times H^s(\T^d)$.  
As a result, $ {\cal G}_M \cap {\cal G}^{(2)}$ is a \dcopen subset of ${\cal P} \times H^s(\T^d)$ and, 
 for any $(V,a) \in {\cal G}_M\cap {\cal G}^{(2)}$, the conditions \eqref{non-reso}-\eqref{non-reso1} hold, provided that 
${\mathbb S} \cup {\mathbb F} \subset [\![ 0,M ]\!]$. 
\end{proposition}

\begin{pfn}{\sc of Proposition \ref{127}.}
Define
$$
 {\cal I}_M := \big\{  (j,k,\s) \in  
 \big([\![ 0,M ]\!] \cap {\mathbb S}^c \big)^2 \times \{ \pm 1 \}  \ : \  j \neq k \  {\rm or } \ \s =1  \big\} \, , 
$$
and, for any $(j,k,\s) \in {\cal I}_M$, 
the real  function $F_{j,k,\s}$ on ${\cal P}\times H^s(\T^d)$  by 
\be  \label{non-dege} 
F_{j,k,\s}(V,a) := \big( \det \Ab \, \mu_j-[ \Bb \Ab^\sharp \bar{\mu}]_j \big) + 
\s \big( \det \Ab \, \mu_k -[ \Bb \Ab^\sharp \bar{\mu}]_k\big) \, .
\ee
\begin{lemma}  \label{nonresder}
The set 
\be\label{NSM:pot}
{\cal N}_{{\mathbb S}, M} := \big\{ V \in E_M \cap {\cal P}  \   :  \  
\forall (j,k,\s) \in  {\cal I}_M \, , \   F_{j,k,\s}(V,1)  \neq 0 \big\}
\ee
is an open and $C^\infty$-dense subset    of ${\cal P}$.
\end{lemma}

\begin{pf}
It is enough to prove that 
\be \label{nonzero}
\hbox{for each} \ (j,k,\s) \in {\cal I}_M \, , \   \hbox{there exists}  \  V\in E_M \cap {\cal P}
\   \hbox{such that} \
F_{j,k,\s}(V,1) \neq 0 \, .
\ee
Indeed, since, for all $ ( j, k, \s ) \in  {\cal I}_M$, the function $F_{j,k,\s}(\cdot ,1)$ is real analytic on
the open connected subset $ E_M \cap {\cal P}$ of $H^s(\T^d)$, by \eqref{nonzero},  
Lemma \ref{analyt} and the finiteness of ${\cal I}_M$, we  conclude 
that ${\cal N}_{{\mathbb S}, M} $ is \dcopen in $E_M \cap {\cal P}$. 
Hence, since $E_M \cap {\cal P}$ is \dcopen in ${\cal P}$, 
the set ${\cal N}_{{\mathbb S}, M} $ is \dcopen in $ {\cal P}$.

To prove \eqref{nonzero}, we consider  the sequence $(q_n)$ of potentials 
provided by Lemma \ref{convL4}, with $J=[ \! [ 0,M ] \! ]$. In particular, $q_n \in E_M \cap {\cal P}$.
By Lemma \ref{convL4}-(i)  we get  
\be\label{convg-1}
\forall j \in [\![ 0,M ]\!] \, , \  \mu_{n,j} \stackrel{n \to + \infty}\longrightarrow \mu_j=\sqrt{\l_j}>0 \  \ {\rm with} \  \  \mu_j \neq \mu_k \ \  {\rm for} 
\ \ j \neq k \, .
\ee
Moreover, by Lemma \ref{convL4}-($i$)-($ii$), for all $ k \in {\mathbb S} $, $ j \in  [\![ 0,M ]\!] \cap {\mathbb S}^c $, 
the  matrix elements  (recall   \eqref{def:AB}, \eqref{def G} and that $ j \neq k $) 
\be\label{convg-2}
\big[\Bb (q_n, 1) \big]^k_j  = \frac32 
\mu_{n,j}^{-1}   ( \Psi_{n,j}^2, \Psi_{n,k}^2 )_{L^2}  \mu_{n,k}^{-1} \to 0  \quad 
{\rm as} \ n \to + \infty \,   .
\ee
In addition, by  \eqref{def:AB}, \eqref{det:G} and \eqref{convg-1}, we have 
\be\label{convg-3}
\det \Ab (q_n,1)  = \Big( \prod_{j \in {\mathbb S}} \mu_{n,j}^{-1} \Big)^2
\det G (q_n,1) \stackrel{n \to + \infty}\longrightarrow \Big( \prod_{j \in {\mathbb S}} \mu_{j}^{-1} \Big)^2 \rho =: \rho_1  >  0 \, . 
\ee
Therefore $(F_{j,k,\s} (q_n,1))$  defined in \eqref{non-dege} 
converges to $ \rho_1 ( \mu_j  + \s  \mu_k  ) $, and 
since 
$ (\mu_j)_{j \in [\![ 0,M ]\!]} $ are distinct and strictly positive, 
$ \rho_1 ( \mu_j  + \s  \mu_k  ) \neq 0 $ for $(j,k,\s) \in {\cal I}_M $. This  implies \eqref{nonzero}. 
\end{pf}

We have the following corollary. 

\begin{corollary}\label{cor:non-deg2}
For each potential $  V(x)  \in {\cal N}_{{\mathbb S},M} $ (defined in \eqref{NSM:pot}), the set 
$$
{\cal G}_{V,{\mathbb S},M}:= \bigcap_{(j,k,\s) \in {\cal I}_M} \big\{ a\in H^s(\T^d) \ : \ 
F_{j,k,\s} (V,a) \neq 0   \big\}
$$
is  $C^\infty$-dense in $ H^s(\T^d) $.
\end{corollary}

\begin{pf}
By Lemma \ref{nonresder},
for each potential $  V  \in {\cal N}_{{\mathbb S},M} $, for each $(j,k,\s) \in {\cal I}_M$,
we have that $F_{j,k,\s} (V, 1)  \neq 0$. Since the function 
$ a \mapsto F_{j,k,\s} (V,a)$ is real analytic on $ H^s (\T^d) $,  
 Lemma \ref{analyt}, and the fact that  ${\cal I}_M$ is finite, imply that 
$ {\cal G}_{V,{\mathbb S},M} $
is a \dcopen subset of $ H^s (\T^d) $. 
\end{pf}

\noindent
{\sc Proof of Proposition \ref{127} concluded.}
The set  ${\cal G}_{ M}$ defined in \eqref{def:GM}   is  clearly open. 
Moreover ${\cal G}_{ M}$
 is  $C^\infty$-dense in $(E_{\mathbb S} \cap {\cal P}) \times H^s(\T^d)$, 
 by Corollary \ref{cor:non-deg2}, and because 
${\cal N}_{{\mathbb S}, M}$ is $C^\infty$-dense in $E_{\mathbb S} \cap {\cal P}$.
Hence ${\cal G}_{ M}$ is a $C^\infty$-dense open subset of $(E_{\mathbb S} \cap {\cal P}) \times H^s(\T^d)$.
Since $E_{\mathbb S}$ is \dcopen in 
${\cal P}$,  the set ${\cal G}_{ M}$ is a $C^\infty$-dense open subset of  $ {\cal P} \times H^s(\T^d)$.
By Corollary \ref{twist2}, so is ${\cal G}_{ M} \cap {\cal G}^{(2)}$, and by
\eqref{matrix-adjugate},  we deduce 
the last claim of Proposition \ref{127}.
\end{pfn}

\medskip
\noindent
{\bf  Step 4. Genericity  of finitely many
first and second Melnikov conditions \eqref{1Mel},
\eqref{2Mel+}\text{-}\eqref{2Mel rafforzate}}  
\\[1mm]
\indent
Let $ L, M  \in \N$ with  ${\mathbb S} \subset [\![ 0, M ]\! ]$. Consider 
the following Conditions:  
\begin{itemize}
\item[($C1$)]  $\bar{\mu}\cdot \ell + \mu_j \neq 0 \, , \   \forall \ell \in \Z^{\mathbb S} ,  | \ell | \leq L , j \in {\mathbb S}^c$, 
\item[($C2$)]  $\bar{\mu}\cdot \ell + \mu_j -\mu_k \neq 0 \, , \   \forall (\ell,j,k) \neq (0,j,j)  \in 
\Z^{\mathbb S} \times ([\![ 0,M] \!] \cap {\mathbb S}^c) \times 
{\mathbb S}^c  \ ,  |\ell| \leq L$ , 
\item[($C3$)]  $\bar{\mu}\cdot \ell + \mu_j +\mu_k \neq 0 \, , \   \forall (\ell,j,k)   \in \Z^{\mathbb S} \times 
([\![ 0,M] \!] \cap {\mathbb S}^c) \times 
{\mathbb S}^c  \ ,  |\ell | \leq L $ .  
\end{itemize} 
Note that ($C1$), ($C2$), ($C3$) correspond to the 
Melnikov conditions \eqref{1Mel},  \eqref{2Mel+}-\eqref{2Mel rafforzate}.

We denote by ${\cal E}_L^{(1)}$, respectively  ${\cal E}_{L,M}^{(2)}$, ${\cal E}_{L,M}^{(3)}$,
the set of potentials $V \in E_{\mathbb S} \cap {\cal P}$
satisfying conditions ($C1$), respectively ($C2$), ($C3$). 

\begin{lemma} \label{finiten}
Let $ L, M \in \N$. The set of potentials 
\be\label{def:ELM} 
{\cal E}_{L,M} := {\cal E}_L^{(1)} \cap   {\cal E}_{L,M}^{(2)} \cap  {\cal E}_{L,M}^{(3)} =
\big\{ V \in E_{\mathbb S} \cap  {\cal P} \ : \,  \hbox{ ($C1$),($C2$),($C3$) hold} \, \big\} 
\ee
is a $C^\infty$-dense open subset of $E_{\mathbb S}\cap {\cal P}$, thus of ${\cal P}$.
\end{lemma}

\begin{pf}
First note that for any potential $ q \in {\cal P}$,  Weyl asymptotic formula
about the distribution of the eigenvalues of $ - \Delta + q(x) $ implies that 
\be\label{weil:mezzo} 
 C_1 j^{1/d} \leq \mu_j \leq C_2 j^{1/d} \, , \quad  \forall j \in \N \, ,  
\ee
for some positive constants $C_1, C_2 $, which is uniform on some open neighborhood $B_q$  of $q$.  
Hence Condition ($C1$) above may be violated by some $V \in B_q$ for $j \leq C (|\bar{\mu}| L)^d$ 
only. Similarly, Conditions ($C2$) and ($C3$) may be violated by some $V \in B_q$ for $ k \leq C (M + (|\bar{\mu}| L)^d) $ only. 
Hence the inequalities in Conditions ($C1$)-($C3$) above are locally finitely many so that 
it enough to
check that for any  $ \ell \in \Z^{|{\mathbb S}|}$, $j \in {\mathbb S}^c$, $k \in {\mathbb S}^c$,
the sets 
\be\label{sets:reso-g}
\begin{aligned}
 {\cal E}^{(1)}_{\ell,j} & := \big\{ V \in E_{\mathbb S} \cap {\cal P} \ : \ \bar{\mu}  \cdot \ell + \mu_j \neq 0 \big\}  \  , \\
 {\cal E}^{(2)}_{\ell,j,k} & := \big\{ V \in E_{\mathbb S} \cap {\cal P} \ : \ \bar{\mu}\cdot \ell + \mu_j -\mu_k \neq 0 \big\}  \  {\rm with }\  \ell \neq 0 \  
{\rm or} \  j\neq k , \\
 {\cal E}^{(3)}_{\ell,j,k} & := \big\{ V \in E_{\mathbb S} \cap {\cal P} \ : \ \bar{\mu}\cdot \ell + \mu_j +\mu_k \neq 0 \big\}
\end{aligned}
\ee
are \dcopen  in $E_{\mathbb S} \cap {\cal P}$. 
The sets in \eqref{sets:reso-g} are open, 
since each $\mu_j$ depends continuously on $V \in E_{\mathbb S} \cap {\cal P}$. 
We now prove  the $C^\infty$-density of ${\cal E}^{(2)}_{\ell,j,k}$, with 
$ (\ell,j,k) \neq (0,j,j) $.
Notice that
the map  
$$
{\Upsilon}_{j,k} : E_{{\mathbb S} \cup \{ j \}\cup \{ k \}} \cap {\cal P} \to \R^{|{\mathbb S}|+2} \, , \quad 
V \mapsto \Upsilon_{j,k} (V)   := (\bar{\mu}, \mu_j, \mu_k ) \, , 
$$ 
is real analytic  and, by Remark \ref{extNJ},  the differential 
$ d {\Upsilon}_{j,k} (V)$ is onto for 
any $V$ belonging to  some \dcopen subset ${\cal N}^{(1)}_{{\mathbb S}\cup \{j\} \cup \{ k \} }$ of $ H^s (\T^d) $.
Hence the map
$ V \mapsto  \bar{\mu}  \cdot \ell + \mu_j - \mu_k $ is real analytic on 
the connected  open set $E_{{\mathbb S} \cup \{ j \} \cup \{ k \}} \cap {\cal P}$ 
and does not vanish everywhere for $ (\ell,j,k) \neq (0,j,j) $. Therefore Lemma \ref{analyt} implies that
 ${\cal E}^{(2)}_{\ell,j,k}$ is
\dcopen in $E_{{\mathbb S} \cup \{ j \} \cup \{ k \}} \cap {\cal P}$  , hence also in $E_{\mathbb S} \cap {\cal P} $. 
The same arguments can be applied to  ${\cal E}^{(1)}_{\ell,j}$ and ${\cal E}^{(3)}_{\ell,j,k}$. 
\end{pf}

\noindent
{\bf Step 5.} {\bf Genericity  of the Diophantine conditions \eqref{diop}, \eqref{NRgamma0}, 
of the first Melnikov conditions \eqref{1Mel} and the second Melnikov conditions
\eqref{2Mel+}\text{-}\eqref{2Mel rafforzate}} 
\\[1mm]
Consider the set of $ (V, a) $ defined by 
\be\label{def:W1}
{\cal G}^{(3)} := \big( ({\cal N}^{(1)}_{\mathbb S} \cap {\cal P}) \times H^s(\T^d) \big) \bigcap {\cal G}^{(2)}  
\ee
where $  {\cal N}^{(1)}_{\mathbb S}  $ is the set of potentials
 defined in \eqref{def:cal-N1} and $ {\cal G}^{(2)} $ is 
 the set of $(V,a)$ defined in  \eqref{def:G2} (for which the twist condition \eqref{A twist} holds). 
 \begin{lemma}\label{lem:G3-dense}
 ${\cal G}^{(3)}$ is a $C^\infty$-dense open subset of ${\cal P} \times H^s (\T^d)$.
 \end{lemma}
 
 \begin{pf}
By Lemma \ref{cor:N1S} and  Corollary \ref{twist2}. 
 \end{pf}

We also remind that  the map 
$$
\bar{\mu} : E_{\mathbb S} \cap {\cal P} \to \R^\es \, , \quad
V \mapsto \bar{\mu}(V) :=(\mu_j(V))_{j \in {\mathbb S}}
$$ 
defined in \eqref{def:Upsilon} is  real analytic, 
and, by Lemma \ref{lem:isoN1}, for any $V \in {\cal N}^{(1)}_{\mathbb S}$, the differential  $d \bar \mu (V)_{|E}$ is an isomorphism,
where $E$ is the  $ |{\mathbb S}|$-dimensional subspace of $ C^\infty (\T^d) $ defined  in Lemma \ref{cor:N1S}. 

We  fix some $(\bar{V} , \bar{a}) \in {\cal G}^{(3)}$ and,
according to the decomposition 
$$ 
H^s(\T^d)=E \oplus F \, , \quad {\rm where } \quad  
F := E^{\bot_{L^2}} \cap H^s (\T^d) \, , 
$$ 
we write  uniquely 
$$
\bar{V}=\bar{v}_1 + \bar{v}_2 \, , \quad \bar{v}_1 \in E \, , \ \bar{v}_2 \in F \, ,  
$$ 
i.e. $ \bar{v}_1 $ is the projection of  $\bar{V} $ on $ E $ and  $ \bar{v}_2 $ is the projection of  $\bar{V} $ on $ F $. 
\begin{lemma}\label{prelim-A} 
(i) There are open balls $B_1 \subset E $ (the subspace $ E \simeq \R^\es $), $B_2 \subset F \subset H^s (\T^d) $, 
centered respectively  at $\bar{v}_1 \in E $, $\bar{v}_2 \in F $,    
such that,  for all $ v_2  \in B_2 $, 
the map 
$$
u_{v_2} : B_1 \subset E \to \R^\es \, , \quad 
v_1 \mapsto u_{v_2} (v_1 ) := \bar \mu ( v_1 + v_2)
$$ 
is a $C^1$ diffeomorphism from $B_1$ onto its image
$ \bar \mu ( B_1 + v_2) =:  {\cal O}_{v_2} $ which is an open bounded subset 
of $\R^{|{\mathbb S}|}$, the closure of which is included in $(0, +\infty)^{|{\mathbb S}|}$. 

(ii) There is a constant $ \barK > 0 $ such that 
the inverse functions  
$$ 
u_{v_2}^{-1} ( \cdot ) :  {\cal O}_{v_2} \subset \R^\es \to E 
$$ 
are  $ \barK $-Lipschitz continuous. 

(iii) There is an open  ball $B \subset H^s(\T^d)$ centered at $\bar{a}$ such that
(possibly after reducing $B_1$ and $B_2$) the  
neighborhood  
$ {\cal U}_{\bar{V}, \bar{a}}:= (B_1+B_2) \times B $ of $ (\bar{V}, \bar{a})  $ 
has closure  contained in  $ {\cal G}^{(3)}$, 
and there is a constant $ C_{\bar{V}, \bar{a}}  > 0 $ such that 
\be\label{def:Vb}
\|a\|_{L^\infty} (1+ \| \Ab^{-1}\|)  \leq C_{\bar{V}, \bar{a}} \, , \quad \forall (V,a) \in {\cal U}_{\bar{V}, \bar{a}} \, , 
\ee
where $ \Ab $ is the Birkhoff matrix in \eqref{def:AB}-\eqref{def G}. 
\end{lemma}

\begin{pf}
Since $\bar{V} = \bar v_1 + \bar v_2 $ is in $  {\cal N}^{(1)}_{\mathbb S}$, the differential 
$d \bar \mu ( \bar v_1 + v_2 )_{|E }$ is an isomorphism and  the local inversion 
theorem implies item ($i$)  of the lemma. 
Items ($ii$)-($iii$) are straightforward taking $ B_1 $ and $ B_2 $ small enough. 
\end{pf}

In the sequel we shall always restrict to the neighborhood   
$ {\cal U}_{ \bar{V}, \bar{a}} $ of $ ( \bar V, \bar a)  \in {\cal G}^{(3)} $ provided by  Lemma \ref{prelim-A}-($iii$).

For any $j \in  \N $, $ v_2 \in B_2  $, we consider the $ C^1 $ function 
\be\label{def-Theta}
\Xi_j : {\cal O}_{v_2} \subset \R^\es \to \R \, , \quad 
\om \mapsto  \Xi_j (\om) := \mu_j (u_{v_2}^{-1} (\om )  + v_2 ) \, .
\ee
Note that, for any $j \in {\mathbb S}$, it results $\Xi_j(\om)=\om_j$.

\begin{lemma}\label{Lip-ov-K}
There is a constant $ \bar{K} := \bar{K}_{\bar V} >  0 $ 
such that, for any potential $ v_2 \in B_2 \subset F \subset H^s (\T^d) $,   
each function $ \Xi_j $ defined in \eqref{def-Theta} is $\bar {K}$-Lipschitz continuous. 
Moreover there are constants $C, C' > 0 $, such that, for all $ v_2 \in B_2 $,  
\be \label{weil} 
C j^{1/d} \leq \Xi_j (\om) \leq C' j^{1/d} \, , \quad \forall j \in \N \, , \ \om \in {\cal O}_{v_2} \, .  
\ee
\end{lemma}
 
\begin{pf}
By \eqref{Lipdep}, $ \lambda_j = \mu_j^2 $,   and the fact that
$$ 
\a:= \inf_{V \in B_1 + B_2} \inf_{j \in \N} \mu_j (V) >  0 
$$ 
we get that 
\be\label{mu-Lip}
\forall V,V'  \in B_1+B_2 \, , \  |\mu_j (V) - \mu_j(V')|  \leq  \frac{1}{2\alpha} \| V-V' \|_{H^s} \, .
\ee 
Since, by Lemma \ref{prelim-A} the functions $ u_{v_2}^{-1} $ are $ \barK $-Lipschitz continuous,  
\eqref{mu-Lip} implies that the function $ \Xi_j $ defined in \eqref{def-Theta} 
 is $\dps \frac{\barK}{2\alpha }$-Lipschitz continuous. 
 The bound \eqref{weil} follows by   \eqref{weil:mezzo}. 
\end{pf}

Notice that, in view of \eqref{cond:suM}, 
the set ${\mathbb M} := {\mathbb M}_{V,a} \subset  \N \backslash {\mathbb S}$ associated to $(V,a)$ 
in Lemma \ref{choice:M} 
has  an upper bound $C(V,a)$ that depends only on $\| \Ab^{-1}\|$ and $\|a\|_{L^\infty}$. 
Hence 
$ \mathbb M $ can be taken constant for all 
$ (V,a) $ in a neighborhood of $ ( \bar{V}, \bar{a}) $, namely,  by \eqref{def:Vb},
\be\label{choice-of-bar-M}
\exists \bar{M} \in \N \ \hbox{such that},  \quad
\forall (V,a) \in {\cal U}_{\bar{V}, \bar{a}} \, , \quad  {\mathbb M}_{V,a}  \cup {\mathbb S} \subset [ \! [ 0, \bar{M} ] \! ] \, . 
\ee

\begin{lemma} \label{genmeasure}
Fix an integer $\bar{L} \geq 4  \bar{K}$ where the constant 
$ \bar{K} $ is defined in Lemma \ref{Lip-ov-K} and let $ \bar{M} \in \N $ 
be as in \eqref{choice-of-bar-M}.
Given a potential  $ v_2 \in B_2  \subset F \subset H^s (\T^d) $, 
we define, for any $\gamma>0$,  the subset of potentials 
${\cal G} ( \gamma , v_2) \subset B_1 \subset E \subset C^\infty $ of all 
the $v_1 \in B_1$ such that, for $V=v_1+v_2$, the following Diophantine conditions hold:
\begin{enumerate} 
\item  
$|\bar{\mu} \cdot \ell | \geq \frac{\gamma}{\la  \ell \ra^{\tau_0}}  \, , \ \forall  \ell \in \Z^{|{\mathbb S}|} \backslash 
\{ 0 \} \, ;$   
\item  
$\Big| n + \sum_{i, j \in {\mathbb S}, i \leq j}  p_{ij}  \mu_i  \mu_j \Big| \geq  \frac{\g}{ \langle p \rangle^{ \t_0}} \, , \quad 
\forall (n, p) \in \Z \times \Z^{\frac{\es(\es+1)}{2}} \setminus \{(0,0)\} \, ;$ 
\item 
$| \bar \mu \cdot \ell + \mu_j | \geq \frac{\g}{ \langle \ell \rangle^{\tau_0}} \, , \quad \forall \ell \in \Z^\es \, ,\   
|  \ell | \geq \bar{L}  \, , \ j \in  {\mathbb S}^c := \N \setminus {\mathbb S} \, ;$
\item 
\be
\begin{aligned}
 | \bar \mu \cdot \ell +  \mu_j -  \mu_k | \geq \frac{\g}{ \langle \ell \rangle^{\tau_0}} \, , \ 
& \forall (\ell, j, k) \in  \Z^\es \times ([\![ 0, \bar{M} ] \! ] \cap  {\mathbb S}^c)   \times {\mathbb S}^c  \, , \\
& (\ell, j, k) \neq (0,j,j) \, ,  \  | \ell | \geq \bar{L} \, , 
\end{aligned}
\ee
\be
\begin{aligned}
 | \bar \mu \cdot \ell +  \mu_j +  \mu_k | \geq \frac{\g}{ \langle \ell \rangle^{\tau_0}} \, , \ 
& \forall (\ell, j, k) \in  \Z^\es \times ([\![ 0, \bar{M} ] \! ] \cap  {\mathbb S}^c) \times {\mathbb S}^c  \, , \\
&   |  \ell | \geq \bar{L} \, . 
\end{aligned}
\ee
\end{enumerate}
Then the measure (on the finite dimensional subspace $ E \simeq \R^\es $)
\be\label{B1minus}
\Big| B_1 \backslash \Big(  \bigcup_{\gamma >0} {\cal G} ( \gamma , v_2) \Big) 
 \Big| = 0 \, . 
\ee
\end{lemma}
\begin{pf}
In this lemma we denote by $ m_{\rm Leb}  $ the Lebesgue measure in $ \R^\es $.
\\[1mm]
i) Let ${\cal F}_{1,\g}$ be the set of the Diophantine frequency vectors $\om \in {\cal O}_{v_2}$ such that 
\be \label{diophomega}
|\om \cdot \ell  | \geq \frac{\gamma}{\la  \ell \ra^{\tau_0}}  \, , \ \forall  \ell \in \Z^{|{\mathbb S}|} \backslash 
\{ 0 \} \, .  
\ee
It is well known that, for $\tau_0> |{\mathbb S}|-1$, $m_{\rm Leb} ({\cal O}_{v_2} \backslash {\cal F}_{1,\g}) =O(\gamma) $.
\\[1mm]
ii) Let ${\cal F}_{2,\g}$ be the set of the frequency vectors $\om \in {\cal O}_{v_2}$ such that 
\be \label{diophomega2}
\Big| n + \sum_{i, j \in {\mathbb S}, i \leq j}  p_{ij} {\om}_i \om_j \Big| \geq  \frac{\g}{ \langle p \rangle^{ \t_0}} \, , \quad 
\forall (n, p) \in \Z \times \Z^{\frac{\es(\es+1)}{2}} \setminus \{(0,0)\} \, . 
\ee
Arguing as in Lemma \ref{lemma:rho-dioph}  we deduce that 
$m_{\rm Leb}({\cal O}_{v_2} \backslash {\cal F}_{2,\g})=O(\g)$, see also Lemma 6.3 in \cite{BB12}. 
\\[2mm]
iii) Let ${\cal F}_{3,\g}$ be the set of the frequency vectors $\om \in {\cal O}_{v_2}$ such that 
\be \label{kov1}
| \om \cdot \ell + \Xi_j (\om) | \geq \frac{\g}{ \langle \ell \rangle^{\tau_0}} \, , \quad \forall \ell \in \Z^\es \, ,\   |  \ell | \geq \bar{L} \, , 
\ j \in {\mathbb S}^c \, . 
\ee
Let us prove that $m_{\rm Leb} ({\cal O}_{v_2} \backslash {\cal F}_{3,\g})=O(\g)$. Define 
on ${\cal O}_{v_2}$ the map 
$ f_{\ell,j}(\om) :=  \om \cdot \ell + \Xi_j (\om) $. By \eqref{weil}, 
there is a constant $ C > 0 $ such
that if $j \geq C |\ell |^d$, then $f_{\ell,j} \geq 1$ on ${\cal O}_{v_2}$. Assume $j \leq C | \ell |^d$.  
 Since $\Xi_j$ is $\bar{K}$-Lipschitz continuous by Lemma \ref{Lip-ov-K}, we have, 
for $| \ell | \geq \bar{L} \geq 4 \bar{K}$, 
$$
\frac{\ell }{| \ell |} \cdot \partial_\om (\om \cdot \ell + 
\Xi_j (\om)) =|\ell|+ \frac{\ell }{| \ell |} \cdot \partial_\om \Xi_j (\om) \geq \frac{3|\ell|}{4} \, . 
$$
Hence, for $| \ell | \geq \bar{L}$, 
$$
m_{\rm Leb} \Big( \Big\{ \om \in {\cal O}_{v_2} \ : \  \exists j \in {\mathbb S}^c \, , \ 
|f_{\ell,j}(\om)| \leq \frac{\g}{ \la \ell \ra^{\tau_0}} \Big\}  \Big) 
\leq C |\ell|^d \frac{\g}{\la \ell \ra^{\tau_0+1}} \leq C \frac{\g}{\la \ell \ra^{\tau_0+1-d}}  \, . 
$$
Hence 
$$
m_{\rm Leb} ({\cal O}_{v_2} \backslash {\cal F}_{3,\g}) \leq C \sum_{|\ell | \geq \bar{L}} \frac{\g}{| \ell |^{\tau_0+1-d}} \leq C \g \, , 
$$ 
provided that $\tau_0>d+ |{\mathbb S}|-1$. \\[1mm]
iv) Let ${\cal F}_{4,\g}$, resp.  ${\cal F}_{5,\g}$,  be the set of the frequency vectors $\om \in {\cal O}_{v_2}$ 
such that 
\be \label{kov2}
\begin{aligned}
| \om \cdot \ell +  \Xi_j (\om) -  \Xi_k (\om) | \geq \frac{\g}{ \langle \ell \rangle^{\tau_0}} \, , 
& \  \forall (\ell, j, k) \in  \Z^\es \times ([\![ 0, \bar{M} ] \! ] \cap  {\mathbb S}^c)   \times {\mathbb S}^c \, , \\
& (\ell, j, k) \neq (0,j,j) \, ,  \  | \ell | \geq \bar{L}  \, ,
\end{aligned}
\ee 
respectively 
\be \label{kov3}
\begin{aligned}
& | \bar \mu \cdot \ell +  \Xi_j (\om) +  \Xi_k (\om) | \geq \frac{\g}{ \langle \ell \rangle^{\tau_0}} \, , \\ 
& \ \forall (\ell, j, k) \in  \Z^\es \times ([\![ 0, \bar{M} ] \! ] \cap  {\mathbb S}^c) \times {\mathbb S}^c  \, , \  |  \ell | \geq \bar{L} \, .
\end{aligned}
\ee
Define on ${\cal O}_{v_2}$ the map 
$f_{\ell,j,k}(\om) :=  \om \cdot \ell + \Xi_j (\om)-\Xi_k(\om) $. By \eqref{weil}, there is a constant $C$ such that, 
if $j \leq \bar{M}$ and $k \geq C (|\ell|^d + \bar{M})$, then
$f_{\ell,j,k} \geq 1$ on ${\cal O}_{v_2}$.  Moreover,  
since $\Xi_j$, $\Xi_k $  are $\bar{K}$-Lipschitz continuous by Lemma \ref{Lip-ov-K},
we deduce that, for $|\ell| \geq \bar{L} \geq 4 \bar{K}$, 
$$
\frac{\ell}{|\ell|} \cdot \partial_\om (\om \cdot \ell + \Xi_j(\om)-\Xi_k(\om)) =
|\ell|+ \frac{\ell }{| \ell |} \cdot \partial_\om \Xi_j (\om) - \frac{\ell }{| \ell |} \cdot \partial_\om \Xi_k (\om) \geq \frac{|\ell|}{2} \, . 
$$
Therefore, for $| \ell | \geq \bar{L}$,
\begin{align*}
& m_{\rm Leb} \Big( \Big\{ \om \in {\cal O}_{v_2} \ : \  \exists (j,k)  \in ([\![ 0, \bar{M} ] \!] \cap {\mathbb S}^c) \times  {\mathbb S}^c \, , \  |f_{\ell,j,k}(\om)| \leq \frac{\g}{\la \ell \ra^{\tau_0}} 
\Big\}  \Big) \\
& \leq C \bar{M}(|\ell|^d+ \bar{M}) \frac{\g}{\la \ell \ra^{\tau_0+1}} \\
& \leq C'(\bar{M}) \frac{\g}{\la \ell \ra^{\tau_0+1-d}}  \, . 
\end{align*}
Hence, as in iii), $m_{\rm Leb} ({\cal O}_{v_2} \backslash {\cal F}_{4,\g})  \leq C \g $, provided that $\tau_0>d+|{\mathbb S}|-1$. We have the similar estimate  
 $m_{\rm Leb} ({\cal O}_{v_2} \backslash {\cal F}_{5,\g})  \leq C \g $. \\[1mm]
 We have proved that $m_{\rm Leb} ({\cal O}_{v_2} \backslash {\cal F}_{i,\g}) =O( \g )$ for $i=1, \ldots , 5$. 
We conclude that   ${\cal F}(\g) :=\cap_{i=1}^5 {\cal F}_{i,\g}$ satisfies
$$
m_{\rm Leb} ({\cal O}_{v_2} \backslash {\cal F} (\g)) \leq \sum_{i=1}^5 m_{\rm Leb} ({\cal O}_{v_2} \backslash {\cal F}_{i,\g}) =O(\g) \, \quad 
\Longrightarrow \quad 
m_{\rm Leb} \big({\cal O}_{v_2} \backslash \bigcup_{\g >0} {\cal F} (\g) \big) = 0 \, . 
$$ 
Finally, since  by Lemma \ref{prelim-A}-($i$) 
the map $u_{v_2}$ is a diffeomorphism between $ B_1 $ and $ {\cal O}_{v_2} $, 
the measure
$$
 m_{\rm Leb} (B_1 \backslash {\cal G} (\g, v_2,a)) = m_{\rm Leb} ( u_{v_2}^{-1} ({\cal O}_{v_2} \backslash {\cal F} (\g)))=0 \, .
 $$
This completes the proof of the lemma. 
\end{pf}
\\[2mm]
{\bf Step 6.} {\bf Conclusion: Proof of Theorem \ref{Prop:genericity}}
\\[1mm]
 For any $ ( \bar{V}, \bar{a}) \in {\cal G}^{(3)} $ (see \eqref{def:W1}), set  
 introduced in \eqref{def:W1}, we define
\be\label{def:GbVbA}
{\cal G}_{\bar{V}, \bar{a}} := \big( {\cal E}_{\bar{L} , \bar{M}} \times H^s(\T^d) \big)
\bigcap {\cal G}_{\bar{M}} \bigcap {\cal U}_{\bar{V}, \bar{a}}  
\ee
where 
$ {\cal U}_{\bar{V}, \bar{a}} \subset {\cal G}^{(3)}  $ is the neighborhood of  $ (\bar{V}, \bar{a}) $ fixed in Lemma \ref{prelim-A}, and 
the sets
${\cal E}_{\bar{L} , \bar{M}} $, $ {\cal G}_{\bar{M}}$ are 
defined respectively in  
 \eqref{def:ELM}, \eqref{def:GM} with 
 integers $\bar{L}$, $\bar{M}$, associated to $(\bar{V}, \bar{a})$, that are fixed in 
 Lemma \ref{genmeasure} and  \eqref{choice-of-bar-M}.  

\begin{lemma}\label{GbVba-dense}
The set ${\cal G}_{\bar{V}, \bar{a}}$ is  $C^\infty$-dense open in ${\cal U}_{\bar{V}, \bar{a}}$. 
\end{lemma}

\begin{pf}
 By Lemma \ref{finiten} and Proposition \ref{127}. 
\end{pf}

Finally, we define the set $ \cal G $ of Theorem \ref{Prop:genericity} as
\be\label{def:setG}
{\cal G} := \bigcup_{(\bar{V}, \bar{a}) \in {\cal G}^{(3)}} {\cal G}_{\bar{V}, \bar{a}}  \, ,
\ee
where $ {\cal G}^{(3)} $ is defined in \eqref{def:W1}. 

\begin{lemma}
 ${\cal G}$ is a \dcopen subset of ${\cal P} \times H^s(\T^d)$. 
\end{lemma} 
 
 \begin{pf}
Since
 ${\cal G}_{\bar{V}, \bar{a}} $ is  open and $C^\infty$-dense  in ${\cal U}_{\bar{V}, \bar{a}}$
 by Lemma \ref{GbVba-dense},  the set $ {\cal G} $ defined in \eqref{def:setG} 
is open and $ C^\infty $-dense in 
$$
\bigcup_{(\bar{V} , \bar{a}) \in {\cal G}^{(3)}} {\cal U}_{\bar{V}, \bar{a}}={\cal G}^{(3)} \, . 
$$ 
(recall that $ {\cal U}_{\bar{V}, \bar{a}} \subset {\cal G}^{(3)} $ by Lemma \ref{prelim-A}). 
Now, by Lemma
\ref{lem:G3-dense}, 
 $ {\cal G}^{(3)} $ is a \dcopen subset of ${\cal P} \times H^s(\T^d)$  
and therefore
${\cal G}$ is a \dcopen subset of ${\cal P} \times H^s(\T^d)$.
\end{pf}

The next lemma completes the proof  of Theorem \ref{Prop:genericity}. 

\begin{lemma}
Let $ \wtilde{a} \in H^s(\T^d)$, $\wtilde{v}_2 \in E^{\bot_{L^2}} \cap H^s(\T^d)$, where $E$ is the finite dimensional linear subspace of
$C^{\infty}(\T^d)$ defined in Lemma \ref{cor:N1S}. Then 
\be  \label{mes0final}
\Big| \big\{ v_1 \in E \, : \,   (v_1+ \wtilde{v}_2, \wtilde{a}) \in {\cal G} \backslash 
\wtilde{\cal G} \big\} \Big| = 0 
\ee
where $\wtilde{\cal G} $ is defined in \eqref{def:tildeG}.
\end{lemma}

\begin{pf}
We may suppose  that 
$$
{\cal W}_{\wtilde{v}_2, \wtilde{a}} := \big\{ v_1 \in E \, : \,   
(v_1+ \wtilde{v}_2, \wtilde{a}) \in {\cal G} \big\} \neq \emptyset \, ,
$$
otherwise \eqref{mes0final} is trivial. Since  $ {\cal G} $ is open in $ H^s (\T^d ) $,
the set 
$ {\cal W}_{\wtilde{v}_2, \wtilde{a}} $ is an open subset of $ E $ and, in order to
deduce  \eqref{mes0final}, it is enough
to prove that, for any $\wtilde{v}_1 \in {\cal W}_{\wtilde{v}_2, \wtilde{a}}$, there is
an open neighborhood 
$ {\cal W}'_{\wtilde{v}_1}  \subset {\cal W}_{\wtilde{v}_2, \wtilde{a}} $ of 
$ \wtilde{v}_1$ such that
\be \label{mes0first}
\Big| \big\{ v_1 \in {\cal W}'_{\wtilde{v}_1} \, : \,   (v_1+ \wtilde{v}_2, \wtilde{a}) \notin 
\wtilde{\cal G} \big\} \Big|=0 \, . 
\ee 
Since $(\wtilde{v}_1 + \wtilde{v}_2,\wtilde{a}) \in {\cal G} $, 
by the definition of $\cal G$ in \eqref{def:setG}, there is 
$(\bar{V}, \bar{a}) \in {\cal G}^{(3)} $ such that 
$$
(\wtilde{v}_1 + \wtilde{v}_2,\wtilde{a}) \in {\cal G}_{\bar{V} , \bar{a}} 
\stackrel{ \eqref{def:GbVbA}} = 
\big({\cal E}_{\bar L , \bar M} \times H^s (\T^d) \big) \bigcap
 {\cal G}_{\bar M} \bigcap  {\cal U}_{\bar{V} , \bar{a}}
$$
where 
$$
{\cal U}_{\bar{V} , \bar{a}} = (B_1 + B_2) \times B
$$ 
is defined in Lemma \ref{prelim-A}.
Define
$$
{\cal W}'_{\wtilde{v}_1} := \big\{ v_1 \in E \, : \, 
(v_1+ \wtilde{v}_2, \wtilde{a}) \in {\cal G}_{\bar{V} , \bar{a}} \big\} \subset B_1 \, .
$$
Since ${\cal G}_{\bar{V} , \bar{a}}$ is open, ${\cal W}'_{\wtilde{v}_1} 
\subset  {\cal W}_{\wtilde{v}_2, \wtilde{a}} 
$ is an open neighborhood of $\wtilde{v}_1   \in {\cal W}_{\wtilde{v}_2, \wtilde{a}} $. 

Now for all $(V,a) \in {\cal G}_{\bar{V} , \bar{a}}$, 
the twist condition \eqref{A twist} and the non-degeneracy properties
 \eqref{non-reso}-\eqref{non-reso1} hold and, by the definition of 
${\cal E}_{\bar L , \bar M}$ in \eqref{def:ELM} and by \eqref{choice-of-bar-M},
there is $\g_0=\gamma_0(V,a)>0$ such that \eqref{1Mel} and \eqref{2Mel+}-\eqref{2Mel rafforzate} 
for all $| \ell | \leq \bar L$. 
Thus, recalling the definition of
 the sets ${\cal G}(\g, \wtilde{v}_2)$ in Lemma \ref{genmeasure} and 
$\wtilde{\cal G}$ in \eqref{def:tildeG}, and \eqref{choice-of-bar-M}, 
we have
$$
{\cal G}_{\bar{V} , \bar{a}} \bigcap \Big( \bigcup_{\g>0} {\cal G}(\g, \wtilde{v}_2)  \Big) 
\subset \wtilde{\cal G} \, , 
$$
so that
$$
\big\{ v_1 \in {\cal W}'_{\wtilde{v}_1} \, : \,   (v_1+ \wtilde{v}_2, \wtilde{a}) \notin 
\wtilde{\cal G} \big\} \subset B_1 \backslash 
\Big( \bigcup_{\g>0} {\cal G}(\g, \wtilde{v}_2)  \Big) \, . 
$$
Hence, by \eqref{B1minus}, the measure 
estimate \eqref{mes0first} holds. 
This completes the proof of \eqref{mes0final}.
\end{pf}


\appendix

\chapter{Hamiltonian and reversible PDEs}\label{App0}

In this appendix 
we first introduce the concept of 
Hamiltonian and/or reversible vector field.
Then we shortly  review the Hamiltonian and/or reversible 
structure of some  classical PDE.

\section{Hamiltonian and reversible vector fields}\index{Hamiltonian PDE}\label{App0-1}

Let $ E $ be a real Hilbert space  with scalar product $ \langle \  , \ \rangle $.  
Endow $ E $ with a constant exact symplectic $ 2 $-form
$$
\Omega (z,w) =  \langle \bar J z , w \rangle  \, , \quad \forall z, w \in E \, , 
$$
where $ \bar J : E \to E  $ is a non-degenerate, antisymmetric linear operator.
Then, given a Hamiltonian  function $ H :   {\cal D}(H) \subset E \to \R $, 
we associate the Hamiltonian system 
\be\label{HS-H}
u_t = X_H (u)  \qquad {\rm where} \qquad d H (u) [ \cdot ] = - \Omega ( X_H(u) , \, \cdot ) 
\ee
formally defines the Hamiltonian  vector field  $  X_H $.  
The vector field  $ X_H : E_1 \subset E \to E  $ is, in general, 
well defined  and smooth  only on a dense subspace  $  E_1 \subset E $.  
A continuous curve 
 $ [t_0, t_1] \ni t \mapsto u (t)  	\in E $ is a solution of  \eqref{HS-H}
 if it is $ C^1 $ as a map from $ [t_0, t_1] \mapsto  E_1 $ and  
$$
 u_t (t) = X_H ( u(t) ) \  {\rm in } \ E \, , \quad  \forall t \in [t_0, t_1 ]  \, . 
$$ 
If, for all $ u \in E_1 $, there is a  vector $ \nabla_u H ( u ) \in  E $ such that 
the differential writes 
\be\label{gradient-0}
d H(u ) [h] = \langle \nabla_u H (u), h \rangle \, , \quad \forall  h  \in E_1  \, ,
\ee
(since $ E_1 $ is in general not a Hilbert space  
with the scalar product $ \langle \,  , \,  \rangle $ then \eqref{gradient-0} does not follow by the Riesz theorem),
then the Hamiltonian vector field  
$ X_H  :  E_1  \mapsto  E  $
writes 
\be\label{HamS} 
X_H  =  J \nabla_u H  \, , \qquad J : = - {\bar J}^{-1} \, . 
\ee
In PDE applications that we shall review below, usually 
$ E = L^2 $ and the dense subspace $ E_1 $ belongs to 
the Hilbert scale
formed by the Sobolev spaces of periodic functions
$$
H^s :=  \Big\{  u(x) = 
\sum_{j \in \Z^d} u_j e^{\ii j \cdot x }  \, : \, \| u \|_s^2 :=  \sum_{j \in \Z^d }   |u_j|^2  (1 + |j|^{2s })  < + \infty  \Big\} 
$$
for $ s \geq 0 $,  
or by the spaces of analytic functions
$$
H^{\s,s} :=  \Big\{  u(x) = \sum_{j \in \Z^d} u_j \, 
e^{\ii j \cdot x }  \, : \, 
\| u \|_{\s,s}^2 := \sum_{j \in \Z^d } |u_j|^2 \, e^{2|j| \s} (1 + |j|^{2 s})  < + \infty   \Big\} 
$$
for $ \s > 0 $. 
For $ s > d / 2 $ the spaces $ H^s $ and $ H^{\s,s} $ are an algebra with respect to the product of functions. 
We refer to  \cite{Ku1} for a general functional setting of Hamiltonian PDEs  on scales of Hilbert spaces. 
\\[1mm]
{\bf Reversible vector field.} 
A vector field $ X $ is 
{\it  reversible}\index{Reversible PDE} \index{Reversible vector field} if there exists an involution $ S $\index{Involution} of the phase space, 
i.e a  linear operator of $ E $ satisfying $ S^2 =  {\rm Id} $,   such that   
\be\label{reversible-VF}
 X  \circ S = - S \circ X \, . 
 \ee
Such condition 
is equivalent to the relation 
$$
\Phi^t  \circ S = S \circ \Phi^{-t}
$$ 
for the flow $ \Phi^t $ associated to the vector field   $ X $. 
For a reversible equation 
it is natural to look for ``reversible" solutions $ u(t) $ of 
$ u_t = X(u) $, namely such that\index{Reversible solution}
$$
u(-t) = S u(t) \, . 
$$
If $ S $ is antisymplectic, i.e. $S^*\Omega = - \Omega$, then
a Hamiltonian vector field $ X = X_H $ is reversible if and only if the Hamiltonian 
$ H $ satisfies 
$$
H \circ S = H  \, .
$$
\begin{remark}
The possibility of developing KAM theory for reversible systems was first observed 
for finite dimensional systems  by Moser in \cite{M67}, see
\cite{Sev} for a complete presentation. In infinite dimension, 
 the first KAM results for reversible
PDEs have been obtained in \cite{ZGY}.
\end{remark}

We now present  some examples of Hamiltonian and/or reversible PDE. 

\section{Nonlinear wave and Klein-Gordon equations}

We consider the nonlinear wave 
equation (NLW) 
\begin{equation}\label{NLWeq}
y_{tt} - \Delta y + V(x) y = f(x, y) \, , \quad x \in \T^d := (\R/ 2 \pi \Z)^d\, ,  \quad y \in \R \, , 
\end{equation}
with a real valued  multiplicative 
potential $ V(x) \in \R $. If  $ V(x) = m $ is  constant, \eqref{NLWeq}  is 
also called a nonlinear Klein-Gordon equation.

The NLW equation \eqref{NLWeq} can be written as the first order Hamiltonian system
$$
\frac{d}{dt}
\left(
\begin{array}{ccc}
 y  \\
 p      
\end{array}
\right)
 = 
\left(
\begin{array}{ccc}
 p  \\
 \Delta y  - V(x) y +  f(x,y)     
\end{array}
\right) = 
\left(
\begin{array}{cc}
0 & {\rm Id} \\ 
- {\rm Id} & 0
\end{array}
\right) 
\left(
\begin{array}{ccc}
 \nabla_y H (y,p)  \\
  \nabla_p H (y,p)       
\end{array}
\right)
$$
where $ \nabla_y H $, $ \nabla_p H $ denote the $ L^2 (\T^d_x) $-gradient
of the Hamiltonian 
\be\label{H-wave}
H(y,p) :=
\int_{\T^d} \frac{p^2}{2} + \frac12 \big( (\nabla_x y)^2 + V(x) y^2 \big) + F(x,y)  \, d x 
\ee
with potential density 
$$
  F (x, y) := - \int_0^y f( x, z ) \, dz  \qquad {\rm and} \qquad \nabla_x y := (\pa_{x_1} y, \ldots, \pa_{x_d} y) \, . 
  $$ 
Thus  for the NLW equation $ E = L^2  \times L^2 $ with $ L^2 := L^2 (\T^d, \R) $, the symplectic matrix 
$$
J =  {\bar J} = 
\left(
\begin{array}{cc}
 0 &    {\rm Id}  \\
 - {\rm Id}  &  0 
\end{array}
\right)
$$
and  $ \langle \ , \ \rangle  $ is the 
$ L^2 $ real scalar product. 
The variables $( y, p)$ are ``Darboux coordinates". 
\begin{remark}\label{rem:1}
The  transposed operator $ {\bar J}^\top = - {\bar J} $  (with respect to $ \langle \, , \, \rangle $)
and 
$ {\bar J}^{-1} = {\bar J}^\top  $. 
\end{remark}

The Hamiltonian $ H $ is properly defined on the subspaces 
$  E_1 := H^s  \times H^s $, $ s \geq 2 $,  
so that the Hamiltonian vector field is a map 
$$ 
J \nabla_u H :  H^{s}  \times H^{s} \to H^{s-2}  \times H^{s-2} \subset  L^2  \times L^2 \, . 
$$
Note that the loss of two derivatives is only due
to the 
Laplace operator $ \Delta $ and that the Hamiltonian vector field 
generated by the nonlinearity is bounded, because the composition operator 
\be\label{composition-operator}
y(x) \mapsto  f(x, y(x)) 
\ee
is a map of $ H^{s}  $ into itself. For such a reason the PDE \eqref{NLWeq} is  {\it semi-linear}.  

If the 
nonlinearity $ f(x,u) $ is analytic 
in the variable $ u $,  and $ H^{s_0} $, 
$ s_0 > d / 2 $, in the space variable $ x $, 
then  the composition operator \eqref{composition-operator} is analytic 
from $ H^{s_0} $ in itself. 
It is 
only finitely many times
differentiable on $ H^{s_0} $ 
if the nonlinearity $ f (x, u) $ is  only (sufficiently)  many times  differentiable with respect to  $ u $. 

\begin{remark}
The  regularity of the vector field is relevant  for KAM theory:
for finite dimensional systems, it has been rigorously proved that,  
if the vector field is not sufficiently smooth, then all the invariant tori 
could be destroyed and only  discontinuous  Aubry-Mather invariant sets survive, see e.g. \cite{He1}. 
\end{remark}

If  the potential density $  F(x,y, \nabla_x y) $ in \eqref{H-wave}
 depends also on the first order
derivative $  \nabla_x y $,  
 we obtain a {\it quasi-linear} wave equation.  
For simplicity we write explicitly only 
the  equation 
in dimension $ d  = 1 $. 
 Given the    Hamiltonian 
$$
H(y,p) :=
\int_{\T} \frac{p^2}{2} + \frac12 \big( y_x^2 + V(x) y^2 \big) + F(x,y,  y_x)  \, d x 
$$
we derive the Hamiltonian wave equation 
$$
y_{tt} - y_{xx} + V(x) y = f(x,y,y_x, y_{xx}) 
$$
with nonlinearity 
(denoting by  $ (x, y, \zeta) $ the independent variables of the potential density
$ (x, y, \zeta) \mapsto F( x, y, \zeta ) $) given by 
\begin{align*}
& f(x,y,y_x, y_{xx})  =  - (\pa_y F)(x, y, y_x ) + \frac{d}{dx} \big\{   ({\pa_{ \zeta } F})(x, y, y_x ) \big\} \\
& = - {(\pa_y F)}(x, y, y_x ) + ({\pa_{ \zeta x } F})(x, y, y_x ) + 
  ({\pa_{\zeta y} F})(x, y, y_x ) y_x +  ({\pa_{\zeta \zeta} F})(x, y, y_x ) y_{xx} \, .
 \end{align*}
Note that $ f $ depends on all the derivatives $ y , y_x, y_{xx} $   
but  it is linear in the second order derivatives $ y_{xx} $, i.e. $ f $ is quasi-linear. 
The nonlinear 
 composition operator 
 $$ 
 y(x) \mapsto   f(x,y(x),y_x (x), y_{xx} (x) ) 
 $$ 
 maps  $ H^{s} \to  H^{s-2} $, i.e. loses two derivatives.
\\[1mm] 
\noindent
{\sc Derivative Wave equations}.
If the nonlinearity $ f(x, y,y_x ) $ depends only on first order derivatives, then  
the Hamiltonian structure of the wave equation is lost (at least the usual one). 
However such equation 
can admit a reversible structure.
Consider the derivative wave, or, better called, derivative Klein-Gordon equation  
$$
 y_{tt} - y_{xx} + m y = f(x, y,y_x, y_t)  \, , \quad   x \in \T \, , 
$$
where the nonlinearity  depends also 
on the first order space and time derivatives $ (y_x, y_t) $, and 
  write it  
 as the first order system
$$
\frac{d}{dt}
\left(
\begin{array}{ccc}
 y  \\
 p     
\end{array}
\right)
 = 
\left(
\begin{array}{ccc}
 p  \\
  y_{xx}  - m y +  f(x,y, y_x, p)     
\end{array}
\right) \, . 
$$
Its vector field $ X $ is {\it reversible} (see \eqref{reversible-VF}) with respect to the involution
$$
S : (y,p) \to (y,-p) \, , \qquad {\rm resp.} \  \ S : (y(x),p(x)) \to (y(-x), -p(-x)) \, , 
$$
assuming the reversibility 
condition
$$
 f(x, y, y_x, - p ) = f(x, y, y_x, p ) \, , \qquad {\rm resp.}  \ \   f(x, y, - y_x,  p ) = f(- x, y, y_x, - p ) \, . 
$$
KAM results have been obtained for reversible derivative wave equations in 
\cite{Berti-Biasco-Procesi-rev-DNLW}.

\section{Nonlinear Schr\"odinger equation }

Consider the Hamiltonian Schr\"odinger equation\index{Nonlinear Schr\"odinger equation} 
\begin{equation}\label{NLSeq}
\ii u_t - \Delta u + V(x) u = f(x, u ) \, , \quad x \in \T^d \, , \ u \in \C \, , 
\end{equation}
where $ f(x, u ) = \partial_{\bar u} F (x, u) $ and the potential $ F(x, u) \in \R  $, $ \forall u \in \C $,  is real valued. 
The NLS equation \eqref{NLSeq} can be written as the infinite dimensional complex 
Hamiltonian equation
$$
u_t =  \ii \nabla_{\bar u} H (u) \, , \quad 
H (u) := \int_{\T^d}    |\nabla u|^2 + V(x) |u|^2  - F(x, u ) \, dx \, .
$$
Actually \eqref{NLSeq}  is a real Hamiltonian PDE. 
In the variables $(a, b) \in \R^2 $, real and imaginary part of 
$$ 
u = a + \ii b \, , 
$$ 
denoting  the real valued potential $ W(a, b) := F(x, a + \ii b )  $ 
so that  
$$ 
\pa_{\bar u} F (x, a+ \ii b) :=  \frac12 (\pa_a + \ii \pa_b) W(a, b)  \, , 
$$ 
the NLS equation \eqref{NLSeq} reads
$$
\frac{d}{dt}
\left(
\begin{array}{ccc}
 a  \\
 b      
\end{array}
\right)
=  
\left(
\begin{array}{ccc} 
\ \Delta b - V(x) b + \frac12 \pa_b W(a, b) 
   \\
\! - \Delta a + V(x) a - \frac12 \pa_a W(a, b)
\end{array}
\right) = \frac12
\left(
\begin{array}{cc}
0 & - {\rm Id} \\ 
{\rm Id} & 0
\end{array}
\right) 
\left(
\begin{array}{ccc}
 \nabla_a H (a,b)  \\
  \nabla_b H (a,b)       
\end{array}
\right)  
$$
with real valued Hamiltonian $ H(a, b) := H(a + \ii b )$ and $ \nabla_a, \nabla_b $ denote the $L^2$-real gradients.
The variables $( a, b )$ are ``Darboux coordinates". 

\smallskip 
A simpler Hamiltonian pseudo-differential model  
 equation  often studied is 
\begin{equation}\label{NLSeq-conv}
\ii u_t - \Delta u + V * u = \pa_{\bar u} F(x, u ) \, , \quad x \in \T^d \, , \ u \in \C \, , 
\end{equation}
where the convolution potential 
$ V * u $ is the Fourier multipliers operator
$$
u(x)  = \sum_{j \in \Z^d} u_j e^{\ii j \cdot x} \quad \mapsto \quad 
V * u :=  \sum_{j \in \Z^d} V_j u_j e^{\ii j \cdot x}  
$$
with real valued Fourier multipliers $ V_j \in \R $. The Hamiltonian of \eqref{NLSeq-conv} is 
$$
H (u) := \int_{\T^d}    |\nabla u|^2 +   (V * u) \, \bar u  - F(x, u ) \, dx \, .
$$
Also for the NLS equation  \eqref{NLSeq}, 
the Hamiltonian vector field 
loses two derivatives because of 
the Laplace  operator $ \Delta $. On the other hand the nonlinear 
 Hamiltonian vector field $ \ii \pa_{\bar u} F(x,u) $
 is bounded, 
so that the PDE \eqref{NLSeq}, as well as \eqref{NLSeq-conv}, is a   semi-linear equation.  

If the nonlinearity 
depends also on first and second order derivatives,  we have,  
respectively, derivative NLS (DNLS) and  fully-non-linear (or quasi-linear) Schr\"odinger equations. For simplicity we 
present the model
equations only in dimension $ d = 1 $. 
For the DNLS equation
\be \label{DNLS}
\ii u_t + u_{xx} = f( x, u, u_x ) 
\ee
the Hamiltonian structure is lost (at least the usual one).  However \eqref{DNLS}
is reversible  with respect to the involution
$$
S : u \mapsto \bar u  
$$
(see \eqref{reversible-VF}) if 
the nonlinearity $ f $ satisfies the condition 
$$
f ( x, \bar u, {\bar u}_x ) = \overline{f (  x, u,  u_x )} \, .  
$$
KAM results for DNLS have been proved  in  \cite{ZGY}.  
On the other hand, fully nonlinear, or quasi-linear,  perturbed  NLS equations like 
$$
\ii u_t = u_{xx} + \e f(\om t, x, u, u_x, u_{xx})
$$
may be  both reversible or Hamiltonian. 
They have been considered in  \cite{FP}.

\section{Perturbed KdV equations}\label{KdV:int}

An important class of equations which arises in fluid mechanics concerns nonlinear perturbations
\begin{equation}\label{KdV-QLA}
u_t + u_{xxx} + \pa_x u^2 + {\cal N} (x, u, u_x, u_{xx}, u_{xxx}) = 0 \, , \quad x \in \T \, ,    
\end{equation}
of 
the KdV equation 
\begin{equation}\label{KdV}
u_t + u_{xxx} + \pa_x u^2  = 0 \, .     
\end{equation}
Here the unknown  $ u(x) \in \R $ is {\it real} valued. 
If the nonlinearity $ {\cal N} (x, u, u_x, u_{xx}, u_{xxx}) $ does not depend on $ u_{xxx} $ such equation is  semilinear.
The most general Hamiltonian (local) nonlinearity is 
\be\label{qlpertA}
{\cal N} (x, u, u_x, u_{xx}, u_{xxx}) := 
- \partial_x \big[ (\partial_u f)(x, u,u_x)  - \partial_{x} ((\partial_{u_x} f)(x, u,u_x)) \big]  
\ee
which is quasi-linear.  In this case  \eqref{KdV-QLA}  is the Hamiltonian PDE
\be\label{KdV Hamiltonian0}
u_t = \partial_x  \nabla_u H(u)  \, , \qquad H (u) = \int_\T  \frac{u_x^2}{2} - \frac{u^3}{3} + f(x, u,u_x) \, dx  \, , 
\ee
where $ \nabla_u H $ denotes the $L^2(\T_x)$  gradient.

The ``mass" $ \int_{\T} u(x) \, dx  $ is a prime integral of \eqref{KdV-QLA}-\eqref{qlpertA} and a 
natural phase space is 
$$
H^1_0 (\T_x) := \Big\{ u(x ) \in H^1 (\T, \R) \ : \ \int_{\T} u(x) dx = 0  \Big\} \, . 
$$
Thus for the Hamiltonian PDE \eqref{KdV Hamiltonian0} we have 
$$
E = L^2_0 (\T, \R) := \Big\{  u \in L^2(\T, \R) \, : \,  \int_\T u(x) \, dx = 0 \Big\} \, , \quad  \bar J = - \partial_x^{-1}  \, , 
$$
 and $ \langle \ , \ \rangle  $ is the 
$ L^2 $ real scalar product. 
Note that $ J = - {\bar J}^{-1} = \pa_x $ see \eqref{HamS}. 
The nonlinear Hamiltonian vector field  is defined and smooth on the subspaces   $ E_1 = H^s_0 ( \T ) $, $ s \geq 3 $, because 
$$ 
{\cal N} (x, u, u_x, u_{xx}, u_{xxx})  :  H^s_0 (\T)  \mapsto H^{s-3}_0 (\T)  \subset L^2_0 (\T) \, . 
$$   
Note that, if the Hamiltonian density $ f = f(x,u) $ does not depend on the first order derivative $ u_x $, the 
nonlinearity in \eqref{qlpertA} reduces to 
$ {\cal N} = - \partial_x  [(\partial_u f)(x, u)]  $ and so the 
PDE \eqref{KdV-QLA}-\eqref{qlpertA} is semilinear.  
\\[1mm]
{\sc Birkhoff coordinates.} 
The KdV equation \eqref{KdV}  has very special algebraic properties:
it possesses infinitely many analytic prime integrals in involution, i.e.  pairwise commuting, 
 and it is integrable in the strongest possible sense, namely it possesses global analytic action-angle variables, 
 called the  Birkhoff coordinates. 
The whole infinite dimensional phase space is foliated by quasi-periodic and almost-periodic solutions.  
The quasi-periodic solutions are called the ``finite gap" solutions. 
Kappeler and collaborators (see for example \cite{KaP}  and references therein) proved 
that there exists an analytic  symplectic diffeomorphism 
\begin{align*}
& \Psi^{-1} : H^N_0 (\T)  \to   \ell^2_{N + \frac 12} (\R) \times \ell^2_{N+ \frac12} (\R) \, , 
\end{align*}
from the Sobolev spaces $ H^N_0 (\T) := H^N (\T) \cap L^2_0 (\T) $, $ N \geq 0 $,   
to the spaces of  coordinates $(x, y) \in   \ell^2_{N + \frac 12} (\R) \times \ell^2_{N+ \frac12} (\R) $ equipped 
with the canonical symplectic form
$ {\mathop \sum}_{n \geq 1} d x_n \wedge d y_n $,  where 
$$
 \ell^2_s (\R )  := \Big\{  x := (x_n)_{n \in \N} \, , \ x_n \in \R 
 \, : \, {\mathop \sum}_{n \geq 1} x_n^2 n^{2s} < + \infty  \Big\} \, , 
$$
such that,  for $ N = 1 $, the KdV Hamiltonian 
$$
 H_{KdV} (u) = \int_\T  \frac{u_x^2}{2} - \frac{u^3}{3}  \, dx  \ , 
$$
expressed in the new coordinates, i.e. $ K := H_{KdV} \circ \Psi $,
 depends only on $ x_n^2 + y_n^2 $, $ n \geq 1 $ (actions).
The KdV equation appears therefore, in these new coordinates,  as an infinite chain of anharmonic oscillators, whose
frequencies 
$$ 
\omega_n (I) := \pa_{I_n} K (I) \, , \quad  I = (I_1, I_2, \ldots ) \, ,  \quad
 I_n = (x_n^2 + y_n^2)/ 2 \,
$$ 
depend on their amplitudes in a nonlinear and real analytic fashion.  
This  is the basis for studying  perturbations of the finite gap solutions of KdV. 

\begin{remark}
The existence of Birkhoff coordinates is much more than what is needed 
for local KAM perturbation theory. As noted by Kuksin in \cite{K2-KdV}-\cite{Ku1} 
it is sufficient that the unperturbed
torus is {\it reducible}, 
namely that the linearized equation at the quasi-periodic solution have,
in a suitable set of coordinates, constant coefficients. 
\end{remark}

Other integrable Hamiltonian PDEs which possess Birkhoff coordinates are the mKdV equation,
see \cite{KaT},
\be\label{mKdV}
u_t + u_{xxx} + \pa_x u^3 = 0 \, , \quad x \in \T \, , 
\ee
and the cubic $ 1 $-d NLS\index{Nonlinear Schr\"odinger equation} equation, see \cite{GrK}, 
\be\label{1d-NLS}
\ii u_t = u_{xx} + | u|^2 u \, , \quad x \in \T \, . 
\ee

\chapter{Multiscale Step}\label{App:mult}

In this appendix we provide the 
proof of the multiscale step proposition  \ref{propinvA} proved in \cite{BBo10}, which is used to prove Proposition
\ref{propinv}. 

\begin{itemize}
\item 
{\bf Notation.} In this appendix we use the notation of \cite{BBo10}, in particular be aware that 
the Definitions \ref{goodmatrixA}, \ref{regularsA}, \ref{ANregA} below   
of $N$-good/bad matrix, regular/singular sites, and $(A,N)$-good/bad site are different with respect to those 
introduced in Section \ref{sec:MS} which are used in the monograph. 
\end{itemize}
In order to give a self-contained presentation we first prove
the properties of the decay norms introduced in \cite{BBo10}, recalled in Section \ref{sec:DN}.

\section{Matrices with off-diagonal decay}\setcounter{equation}{0}\label{sec:off}

Let $ e_i = e^{\ii (\ell \cdot \vphi + j \cdot x)} $ for 
$ i := (\ell,j) \in \Z^b := \Z^\nu \times \Z^d $.  
In the vector-space $ {\mathcal H}^s = 
  {\mathcal H}^s ( \T^\nu \times \T^d;  \C^r)  $ defined in \eqref{def:Hs} 
with $ \nu =  \es $, 
  we consider 
  the basis 
\be\label{basisi}
e_k  
= e_i  e_{\mathfrak a}  \, , \quad k := (i, \mathfrak a) \in \Z^b \times {\frak I}\, , 
\ee
where 
$ e_{\mathfrak a} := (0, \ldots, \underbrace{1}_{{\mathfrak a}-th}, \ldots,0 )  \in \C^r $,
$  {\mathfrak a} = 1, \ldots, r $,  
denote the canonical basis of $ \C^r $, and 
$$ 
{\frak I} := \{1, \ldots, r \} \, . 
$$
Then we write any $ u \in {\mathcal H}^s  $ as
$$
u =  \sum_{k \in \Z^b \times {\frak I}}  u_{k} e_k  
\, ,  \ \  u_{k} \in \C \, . 
$$
For $ B \subset \Z^b \times {\frak I} $, we introduce the subspace
$$
{\mathcal H}^s_B := \Big\{ u \in {\mathcal H}^s \, : \, u_{k} = 0 \ {\rm if} 
\ k \notin B \Big\} \, .
$$
When $ B $ is finite, the space $ {\mathcal H}^s_B $ does not 
depend on $ s $ and will be denoted $ {\mathcal H}_B $. We define 
$$ 
\Pi_B : {\mathcal H}^s \to {\mathcal H}_B
$$ 
the $L^2$-orthogonal projector onto $ {\mathcal H}_B $. 

In what follows $ B, C, D, E $ are finite  subsets of $ \Z^b\times {\frak I} $.

We identify the space
$ \lin^B_C $ of the linear maps $L : {\mathcal H}_B \to {\mathcal H}_C $  with the space of matrices
$$
\matr^B_C := \Big\{ M = (M^{k'}_k)_{k' \in B, k \in C} \, , \ M^{k'}_{k} \in \C  \Big\} 
$$ 
according to the following usual definition. 

\begin{definition}\label{matrixrepre}
The matrix $ M \in \matr^B_C  $ represents the linear operator $ L \in  \lin^B_C $, if 
$$ 
\forall k'= (i', \mathfrak a') \in B, \, k=(i, \mathfrak a) \in C \, ,  \quad
\Pi_k L e_{k'}=M_k^{k'} e_k  \, , 
$$
where $ e_k $ are defined in \eqref{basisi} and $ M_k^{k'} \in \C $.
\end{definition}

\noindent
{\sc Example.} The multiplication operator  
for 
$$
\left(
\begin{array}{cc}
p(\vphi, x) & q(\vphi, x) \\
\ov{q}(\vphi, x) & p(\vphi, x)   
\end{array} \right) \, , 
$$
acting in 
$  {\mathcal H}^s ( \T^\nu \times \T^d;  \C^2)  $,  is represented by the matrix
\be\label{T-mult}
 T := (T_i^{i'})_{i, i' \in \Z^b} \qquad {\rm where} \qquad
T_i^{i'} = 
\left(
\begin{array}{ccc}
 p_{i-i'} &   q_{i-i'}    \\
 (\overline{q})_{i-i'} &      p_{i-i'}    
\end{array}\right) 
\ee
and $ p_i, q_i $ denote the Fourier coefficients of $ p(\vphi, x), q(\vphi, x) $.  
With the above notation, the set $ {\mathfrak I } = \{1, 2 \} $ and 
$$
T_{(i,1)}^{(i',1)} = p_{i-i'} \, , \ T_{(i,1)}^{(i',2)} = q_{i-i'} \, , \ T_{(i,2)}^{(i',1)} 
= (\overline{q})_{i-i'} \, , \ T_{(i,2)}^{(i',1)}  = p_{i-i'} \, .
$$
{\sc Notation.} 
For any subset $B$ of $\Z^b \times {\frak I}$, we denote by 
\be\label{proj}
\ov{B} := {\rm proj}_{\Z^b} B 
\ee
the projection of $B$ in $\Z^b$.

Given $ B \subset B' $, $ C \subset C' $ $ \subset \Z^b \times {\frak I}$ and $ M \in {\cal M}_{C'}^{B'} $ we can introduce
the restricted matrices  
\be\label{defMBC}
M^B_C := \Pi_C M_{|{\mathcal H}_B} \, ,  \quad
M_C := \Pi_C M \, , \quad M^B := M_{|{\mathcal H}_B} \, .
\ee
If $ D \subset {\rm proj}_{\Z^b} B'$, $E \subset {\rm proj}_{\Z^b} C'  $, then we define 
\be\label{BCtilde}
M_E^{D} \ \ {\rm  as } \ \ M_{C}^{B}
\quad {\rm where} \quad   {B} := (D \times {\frak I}) \cap B', \
 {C} := (E\times {\frak I}) \cap C' \, .
\ee
In the particular case  $ D = \{ i' \} $,  $ E := \{ i \} $, $i,i' \in \Z^b$,  we use the simpler notation  
\begin{align}\label{line1}
M_i & := M_{\{ i \}} \quad {\rm (it \ is \ either \ a \ line \ or \ a \ group \ of } \  2, \ldots, r \ { \rm lines \ of} \   M {\rm )} \, , \\
\label{line2}
M^{i'} & := M^{\{ i' \}} \ \   {\rm (it \ is \ either \ a \ column \ or \ a \ group \ of } \ 2, \ldots, r \ 
{\rm columns \ of} \ M {\rm )} \, , 
\end{align}
and 
\be\label{line3} 
M_i^{i'} := M_{\{ i \}}^{\{ i' \}} \, , 
\ee
it is a $ m \times m' $-complex matrix, where $ m \in \{1, \ldots, r \} $ 
(resp. $ m' \in \{1, \ldots, r \} $) is
the cardinality of $ C $ (resp. of $  B $) defined in \eqref{BCtilde} with  
 $ E := \{ i \} $ (resp. $ D = \{ i' \} $). 

We endow the vector-space of the $ m \times m' $, $ m, m' \in \{ 1, \ldots,  r \} $, 
complex matrices with a norm $ | \ | $ such that 
$$
|U W | \leq |U\| W| \, , 
$$
whenever the dimensions of the matrices make their multiplication possible, and $|U| \leq |V|$ if $U$ is a submatrix of $V$. 

\begin{remark}
The notation in \eqref{BCtilde}, \eqref{line1}, \eqref{line2}, \eqref{line3}, may be not very specific, but it is deliberate:  
it is convenient not to distinguish the index $ \mathfrak a \in {\frak I}$, which is irrelevant
in the definition of the $ s $-norms, in Definition \ref{defnormatrA}.
\end{remark}

We also set the $ L^2 $-operatorial norm 
\be\label{L2norm}
\| M^B_C \|_0 := \sup_{h \in {\mathcal H}_B, h \neq 0} \frac{\| M^B_C h \|_0}{\| h \|_0}  
\ee
where  $ \| \ \|_0 := \| \ \|_{L^2 } $.

\begin{definition} \label{defnormatrA} {\bf ($s$-norm)}
The $ s $-norm\index{Decay norm}  of a matrix $ M \in \matr^B_C $ is defined by
\be\label{defM}
\norma M \norma_s^2 :=  \sum_{n \in \Z^b} [M(n)]^2 \langle n \rangle^{2s}  
\ee
where $ \langle n \rangle := \max (|n|,1)$, 
\be\label{Mpri}
[M(n)] := \begin{cases}
\max_{i-i'=n, i \in {\ov C}, i' \in {\ov B}} |M^{i'}_i|  \ \ \ \ \quad   {\rm if}  \ \   n \in \ov{C}-\ov{B}\\
0 \qquad \qquad \qquad \qquad \ \ \quad  {\rm if}  \ \  n \notin \ov{C}-\ov{B} \, , 
\end{cases}
\ee
with $ \ov{B} := {\rm proj}_{\Z^b} B $, $ \ov{C} := {\rm proj}_{\Z^b} C $
(see \eqref{proj}).
\end{definition}

It is easy to check that $\nors{\ }$ is a norm on $ \matr_C^B $. It verifies $ \nors{\ } \leq \norma \ \norma_{s'} $, 
$ \forall s \leq s' $, and
$$ 
\forall M \in \matr^B_C \, , 
\quad \forall B' \subseteq B \ , \ C' \subseteq C \ , \quad \nors{M^{B'}_{C'}} \leq \nors{M} \, . 
$$
The $ s $-norm is designed to estimate the off-diagonal decay of matrices 
similar to the T\"oplitz matrix  which represents the multiplication operator
for a Sobolev function.

\begin{lemma}\label{lem:multi}  
The matrix $ T $ in \eqref{T-mult}
with $ (p, q) \in {\mathcal H}^s ( \T^\nu \times \T^d;  \C^2)  $ satisfies
\be\label{multiplication}
\nors{T} \lesssim \| (q,p) \|_s \, . 
\ee
\end{lemma}
\begin{pf}
By \eqref{Mpri}, \eqref{T-mult} we get
$$
[T (n)] := \max_{i-i'=n} \Big| \left(
\begin{array}{cc}
p_{i-i'} & q_{i-i'} \\
\ov{q_{i-i'}} & p_{i-i'}
\end{array}
\right)   \Big|  
\lesssim |p_n|+|q_n| \, .
$$
Hence, the definition in \eqref{defM} implies
$$
\norma T \norma_s^2 =  \sum_{n\in \Z^b} [T(n)]^2 \la n \ra^{2s}
\lesssim \sum_{n \in \Z^b} (|p_n|+|q_n|)^2 \la n \ra^{2s}
\lesssim \|(p,q)\|_s^2 
$$
and \eqref{multiplication} follows. 
\end{pf}

In order to prove that the matrices with finite 
$ s $-norm satisfy the interpolation inequalities \eqref{interpm}, and then the algebra property
\eqref{algebra}, the guiding principle
is the analogy between these matrices and the T\"opliz matrices which represent  
the multiplication operator  for functions.
We introduce the set  $ \mapp $ of 
the trigonometric polynomials with positive Fourier coefficients
$$
\begin{aligned}
\mapp 
:= \Big\{ & h = \sum h_{\ell,j} e^{\ii (\ell \cdot \varphi + j \cdot x)}  \ {\rm with} \ h_{\ell,j} \neq 0  \\
&  {\rm for \ a \ finite \ number \ of} \ (\ell,j) \ {\rm only \  and \ }  
h_{\ell,j} \in \R_+ \Big\}\, . 
\end{aligned}
$$
Note that the  sum and the product of two functions in $ \mapp $ remain in $ \mapp $. 

\begin{definition}\label{domina} 
Given $ M \in \matr^B_C $, $ h \in \mapp $, we say that $M$ is dominated by $ h $, and we write $ M \prec h $,
if 
\be\label{Msec}
[M(n)] \leq h_n \, , \quad \forall n \in \Z^b \, ,
\ee
in other words  if  $\  |M_i^{i'}| \leq h_{i-i'} $ , $ \forall i' \in {\rm proj}_{\Z^b} B $, $ i \in {\rm proj}_{\Z^b} C $. 
\end{definition}

It is easy to check  ($B$ and $C$ being finite) that
\be \label{normdom}
\begin{aligned}
 \norma M \norma_s   = & \min \Big\{ \|h\|_s \  : \  h\in \mapp \ , \  M \prec h \Big\} \quad {\rm and} \\ 
& \exists h \in \mapp \ ,  \  \forall s\geq 0 \, ,\  \nors{M} =  \|h\|_s \, .
\end{aligned} 
\ee

\begin{lemma} \label{domprod}
For $ M_1 \in \matr^C_D $, $M_2 \in \matr^B_C$, $M_3 \in \matr^C_D$, we have 
$$
\begin{aligned}
& \quad M_1\prec h_1 \, , \
M_2\prec h_2 \, , \ M_3\prec h_3 \quad \Longrightarrow \\
&  M_1+M_3 \prec h_1 + h_3  \quad {\rm and} \quad M_1M_2 \prec h_1 h_2 \, .
\end{aligned}
$$
\end{lemma}

\begin{pf}
Property $ M_1+M_3 \prec h_1 + h_3 $ is straightforward. 
For $i\in {\rm proj}_{\Z^b}  D$, $i'\in {\rm proj}_{\Z^b} {B}$, we have
\begin{eqnarray*}
| (M_1M_2)^{i'}_i |  = \Big| \sum_{q \in \ov{C} := {\rm proj}_{\Z^b} C} (M_1)^q_i (M_2)_q^{i'} \Big| &\leq &
\sum_{q \in \ov{C}} |(M_1)^q_i|  |(M_2)_q^{i'}| \\ 
& \leq &   \sum_{q \in \ov{C}} (h_1)_{i-q} (h_2)_{q-i'} \\
&\leq & \sum_{q \in \Z^b} (h_1)_{i-q} (h_2)_{q-i'}=(h_1h_2)_{i-i'}
\end{eqnarray*}
implying $ M_1M_2 \prec h_1 h_2  $ by Definition \ref{domina}. 
\end{pf}

We deduce from  \eqref{normdom}
and 
 Lemma 
\ref{lem:int-impro}, the following interpolation estimates. 

\begin{lemma} \label{prodest}
{\bf (Interpolation)} 
For all $ s \geq s_0 > (d+\nu)/2 $, there is $ C(s) \geq 1 $, such that, 
for any finite subset $ B, C, D \subset \Z^b \times {\frak I} $,
for all matrices $ M_1 \in \matr^C_D$, $ M_2 \in \matr^B_C $, 
\be \label{interpm}
\nors{M_1 M_2} \leq  C(s_0) \norso{M_1} \nors{M_2} + C(s) \nors{M_1} \norso{M_2} \, ,   
\ee
in particular, 
\be \label{algebra}
\norma M_1 M_2 \norma_s \leq  C(s) \nors{M_1} \nors{M_2} \, .
\ee
\end{lemma}

Note that the constant $ C(s) $ in Lemma \ref{prodest} is independent of $ B $, $ C $, $ D $. 

\begin{lemma} \label{prodest1}
For all $ s \geq s_0 > (d+\nu)/2 $, there is $ C(s) \geq 1 $, such that, 
for any finite subset $ B, C, D \subset \Z^b \times {\frak I}  $, 
for all $  M_1 \in \matr^C_D $, $ M_2 \in \matr^B_C $, we have
\be \label{norsoest}
\norso{M_1M_2} \leq C(s_0 ) \norso{M_1}\norso{M_2} \, , 
\ee
and, $ \forall M \in \matr^B_B $, $ \forall n \geq 1 $, 
\label{powerest} 
\be\label{Mnab}
\begin{aligned}
& \norma M^n \norma_{s_0} \leq 
( C(s_0))^{n-1} \norma M \norma_{s_0}^n \, , \\
& \norma M^n \norma_{s} \leq  C(s) (C(s_0))^{n-1} 
\norma M \norma_{s_0}^{n-1}  \norma M \norma_{s} \, , \ \forall s \geq s_0 \, .
\end{aligned}
\ee
\end{lemma}

\begin{pf}
The first estimate in \eqref{Mnab} follows by 
\eqref{algebra} with $ s = s_0 $ and 
the second estimate in (\ref{Mnab}) is obtained from (\ref{interpm}),  using $C(s) \geq 1$. 
\end{pf}

The $ s $-norm of a matrix  $ M \in \matr^B_C $ controls also the Sobolev $ {\mathcal H}^s $-norm. Indeed, 
we  identify $ {\mathcal H}_B $ with the space $ \matr_B^{\{0\}} $ of column matrices and 
the Sobolev norm $ \|  \ \|_s $ is equal to the $ s $-norm $ \norma \ \norma_s $, i.e.
\be\label{vettori=}
\forall w \in {\mathcal H}_B \, , \ \ \| w \|_s  =  \nors{w} \, , \quad  \forall s \geq 0 \, .
\ee
Then $ M w \in {\mathcal H}_C $ and the next lemma is a particular case of Lemma \ref{prodest}.

\begin{lemma}\label{sobonorm} {\bf (Sobolev norm)}
$ \forall s \geq s_0 $ there is $ C(s) \geq 1 $ such that, for any finite subset $ B, C \subset \Z^b \times {\frak I}$, for any  $ M \in \matr^B_C $, $ w \in {\mathcal H}_B  $, 
\be\label{opernormA}
\| Mw \|_s \leq C(s_0) \norma M\norma_{s_0} \|w\|_s + 
C(s) \norma M \norma_s \|w\|_{s_0} 
\, .
\ee
\end{lemma}

The following lemma is the analogue of the smoothing properties 
of the projection operators.

\begin{lemma} \label{norcompA}
{\bf (Smoothing)} 
Let $ M  \in \matr^B_C $. Then, $ \forall s' \geq  s \geq 0 $,
\be \label{Sm1A}
M_i^{i'}=0 \, , \ \forall |i - i' | < N \quad \Longrightarrow  \quad  \nors{M} \leq N^{-(s'-s)} \norma M \norma_{s'} \, , 
\ee
and,  for $ N \geq N_0 $, 
\be \label{Sm2} 
M_i^{i'}=0 \, , \ \forall |i - i' | > N 
\quad  \Longrightarrow  \quad 
\begin{cases} 
\norma M \norma_{s'}   \leq N^{s'-s} \nors{M} \, \\
\nors{M} \leq N^{s+b}  \| M \|_{0} \, . 
\end{cases} \ee
\end{lemma} 

\begin{pf}
Estimate (\ref{Sm1A}) and the first bound of (\ref{Sm2})
follow from the definition of the norms $\norma \  \norma_s $. The second bound 
of  (\ref{Sm2}) follows by the first bound in (\ref{Sm2}), noting that $ | M^{i'}_i | \leq \| M \|_0 $, $ \forall i, i' $, 
$$
\nors{M}  \leq N^{s}  \norma{M} \norma_0 \leq 
 N^s \sqrt{(2N+1)^{b}} \| M \|_{0} \leq   N^{s+b} \| M \|_{0}
$$
for $ N \geq N_0 $.
\end{pf}

In the next lemma 
we bound the $ s $-norm of a matrix in terms of the $ (s + b)$-norms of its lines. 

\begin{lemma} \label{norextracted}
{\bf (Decay along lines)}
Let $ M \in \matr^B_C$. Then, $ \forall s\geq 0$, 
\be\label{Mdecay} 
\nors{M} \lesssim \max_{i \in {\rm proj}_{\Z^b} C} \norma M_{\{i\}} \norma_{s+b} 
\ee
(we could replace the index $ b $ with any $ \a > b / 2 $).
\end{lemma}

\begin{pf}
For all $ i \in \ov{C} := {\rm proj}_{\Z^b} C $, $ i' \in \ov{B} := {\rm proj}_{\Z^b} B $, $ \forall s \geq 0 $, 
$$
|M_i^{i'}|  \leq \frac{ \norsb{M_{\{i\}}} }{\langle i - i' \rangle^{s+b}} \leq   \frac{m(s+b)}{\langle i - i' \rangle^{s+b}} 
$$
where $ m(s+b) := \max_{i \in \ov{C}} \norsb{M_{\{i\}}} $. 
As a consequence 
$$
\nors{M} = \Big( \sum_{n \in \ov{C} - \ov{B}} (M[n])^2 \langle n \rangle^{2s} \Big)^{1/2} \leq  
m(s+b) \Big( \sum_{n \in \Z^b} \langle n \rangle^{-2b}   \Big)^{1/2} 
$$
implying \eqref{Mdecay}.
\end{pf}

The $ L^2 $-norm  and $ s_0 $-norm of a matrix are related.

\begin{lemma}\label{Lem:L2s0} 
Let $ M \in {\cal M}_B^C $. Then, for $ s_0 > (d+\nu) / 2 $,
\be\label{schur}
\| M \|_0 \lesssim_{s_0} \norso{M} \, .  
\ee
\end{lemma}

\begin{pf}
Let $ m \in {\cal H}_+ $ be  such that $ M \prec m $ and $ \nors{M} = \| m \|_s $ for 
all $ s \geq 0 $, see (\ref{normdom}). Also for $H\in {\cal M}^{\{ 0\}}_C$, there is $ h \in {\cal H}_+ $
such that $ H \prec h $ and $ \nors{H} = \|h\|_s $, $ \forall s \geq 0 $. Lemma \ref{domprod} implies that 
$ M H \prec m h $ and so
$$
\norma MH \norma_0 \leq \|mh\|_0 
\leq |m|_{L^\infty} \|h\|_0 \lesssim_{s_0} \|m\|_{s_0}  \|h\|_0 = \norso{M} \norma H \norma_0 \, , \quad
\forall H\in {\cal M}^{\{ 0\}}_C \, .
$$
Then \eqref{schur} follows (recall \eqref{vettori=}).
\end{pf}

In the sequel we use the notion of left invertible operators.

\begin{definition} {\bf (Left Inverse)}\index{Left inverse}
A matrix $ M \in \matr_C^B $ is left invertible if there exists  
$ N \in \matr^C_B $  such that $ NM = {\rm Id}_B $. Then $ N $ is called a left inverse of  $ M $. 
\end{definition}

Note that $ M $ is left invertible if and only if $ M $ (considered as a linear map) is injective
(then ${\rm dim} \, {\mathcal H}_C \geq {\rm dim} \, {\mathcal H}_B $). The left inverses  of $ M $  are not unique
if ${\rm dim} \, {\mathcal H}_C > {\rm dim} \, {\mathcal H}_B $: 
they are uniquely defined only on the range of $ M $.

\smallskip

We shall often use the following perturbation lemma for left invertible operators. Note that 
the bound  (\ref{NR12}) for the perturbation in $ s_0 $-norm only, allows to estimate
the inverse (\ref{inv2}) also in $ s \geq s_0 $ norm. 

\begin{lemma} \label{leftinvA} {\bf (Perturbation of left invertible matrices)}
If $ M \in \matr_C^B $ has a left inverse $ N \in \matr^C_B $ , 
then there exists $ \d (s_0) > 0 $ such that,  
\be\label{NR12}
\forall P \in \matr_C^B \qquad { with} \qquad 
\norma N \norma_{s_0}  \norma P \norma_{s_0} \leq \d(s_0) \, , 
\ee 
the matrix $ M + P $ has a left inverse $ N_P $ that satisfies
\be \label{inv1A}
\norma N_P \norma_{s_0} \leq 2\norma N \norma_{s_0} \,,
\ee
and, $ \forall s \geq s_0 $,
\begin{align}  \label{inv12}
\norma N_P \norma_{s} & \leq   \big(1+ C(s)\norma N \norma_{s_0} \norma P \norma_{s_0} \big) \norma N \norma_s + 
C(s) \norma N \norma^2_{s_0} \norma P \norma_s \\
& \lesssim_s   \norma N \norma_s + \norma N \norma^2_{s_0} \norma P \norma_s \, . \label{inv2}
\end{align}
Moreover, 
\be\label{NR0}
\forall P \in \matr_C^B \qquad { with} \qquad \| N \|_0  \| P \|_0 \leq 1/2 \, , 
\ee 
the matrix $ M + P $ has a left inverse $ N_P $ that satisfies
\be \label{inv10}
\| N_P \|_0 \leq 2 \| N \|_0 \, .
\ee
\end{lemma}

\begin{pf}

\noindent
{\it Proof of \eqref{inv1A}.} The matrix $ N_P = A N $  with $ A  \in \matr^B_B $ is a left inverse of $ M + P $ if and only if
$$
I_B = A N (M+P) = A (I_B + NP) \, , 
$$
i.e. if and only if $ A $ is the inverse of $ I_B + NP \in \matr^B_B $. 
By  \eqref{norsoest} and \eqref{NR12} we have, taking $ \d(s_0) > 0 $ small enough,  
\be\label{NP1/2}
\norso{NP} \leq C(s_0) \norso{N} \norso{P} \leq C(s_0) \d(s_0) \leq 1/2 \, .
\ee
Hence  the matrix $ I_B + NP $ is invertible and
\be \label{devnr}
N_P = A N = (I_B + NP)^{-1} N = \sum_{p=0}^\infty (-1)^p (NP)^p N 
\ee
is a left inverse of $ M + P $.  
Estimate \eqref{inv1A} follows by 
\eqref{devnr} and \eqref{NP1/2}. 
\\[1mm]
{\it Proof of (\ref{inv12}).}  For all $ s \geq s_0 $, $\forall p \geq 1 $,  
\begin{align}
 \nors{(NP)^p N} & \stackrel{(\ref{interpm})} { \lesssim_s}  
 \norso{N} \nors{(NP)^p}  +   \nors{N}  \norso{(NP)^p} \nonumber \\
& \stackrel{ (\ref{Mnab}) } {\lesssim_s}  
\norso{N} (C(s_0)  \norso{NP})^{p-1} \nors{NP}  +  \nors{N} 
(C(s_0) \norso{NP})^{p}  \nonumber \\
& \stackrel{\eqref{NP1/2}, \eqref{interpm}} {\lesssim_s}  
 2^{-p} ( \norso{N} \norso{P}  \nors{N} +  \norso{N}^2 \nors{P}) \, . \label{stimahigh}
\end{align}
We derive \eqref{inv12} by
$$
\begin{aligned}
\nors{N_P} 
& \stackrel{(\ref{devnr})}   \leq  \nors{N} + \sum_{p\geq1} \nors{(NP)^pN} \\
& \stackrel{(\ref{stimahigh})} \leq  \nors{N} +
C(s) (  \norso{N} \norso{P}  \nors{N} +  \norso{N}^2 \nors{P}) \, .
\end{aligned}
$$
Finally  \eqref{inv10}
follows from \eqref{NR0}
as \eqref{inv1A} because the operatorial $ L^2 $-norm (see (\ref{L2norm}))
satisfies the algebra property $ \|  N P \|_0 \leq \| N \|_0 \| P \|_0 $. 
\end{pf}

\section{Multiscale step proposition}\label{App:MU}

This section is devoted to prove the multiscale step proposition \ref{propinvA}.

In the whole section  $\varsigma \in (0,1)$  is fixed and $\tau'>0$, $\Theta \geq 1$
are  real  parameters, on which we shall impose conditions in  Proposition \ref{propinvA}. 

Given $ \Omega, \Omega' \subset  E \subset \Z^b \times {\frak I} $ we define 
$$
{\rm diam}(E) := \sup_{k,k' \in E} |k-k'| \, , \qquad 
 {\rm d}(\Omega, \Omega') := \inf_{k \in \Omega, k' \in \Omega'} |k-k'|  \, , 
$$
where, for $ k= (i, \mathfrak a) $, $ k' := (i', \mathfrak a')$ we set 
$$ 
|k - k' | := \begin{cases} 
1  \quad \qquad \qquad \, \text{if} \ i = i' \, , \mathfrak a \neq \mathfrak a' \, ,   \cr 
0  \quad \qquad \qquad \, \text{if} \ i = i' \, , \mathfrak a = \mathfrak a' \, ,   \cr 
|i-i'|  \qquad \quad \text{if} \ i \neq i'  \, . 
\end{cases}
$$ 

\begin{definition}\label{goodmatrixA}
{\bf ($N$-good/bad matrix \cite{BBo10})} The matrix $ A \in {\cal M}_E^E $, with $ E \subset \Z^b \times {\frak I}$, 
$ {\rm diam}(E) \leq 4 N $,  is $ N $-good if $ A $ is invertible and
\be\label{NgoodmatA}
\forall s \in [s_0, s_1] \   , \ \ |A^{-1}|_s \leq N^{\tau'+\d s}.
\ee
Otherwise $ A $ is $ N $-bad.
\end{definition}

\begin{definition} \label{regularsA} {\bf (Regular/Singular sites \cite{BBo10}) } 
The index $ k := (i,\mathfrak a) \in \Z^{b} \times {\frak I}$  is  {\sc regular} for $A $ 
if $ |A_k^k| \geq \Theta $. Otherwise $ k $ is  {\sc singular}.
\end{definition}

\begin{definition}\label{ANregA}
{\bf ($(A,N)$-good/bad site \cite{BBo10})}
For $ A \in \matr^E_E $, we say that $ k \in E \subset \Z^b \times {\frak I}  $ is
\begin{itemize}
\item $(A,N)$-{\sc regular} if there is $ F \subset E$ such that
${\rm diam}(F) \leq 4N$,  ${\rm d}(k, E\backslash F) \geq N$ and
$A_F^F$ is $N$-good. 
\item $(A,N)$-{\sc good}  if it is regular for $A$ or $(A,N)$-regular. Otherwise we say that $ k $ is $(A,N)$-{\sc bad}.
\end{itemize}
\end{definition}

Let us consider  the new larger scale 
\be\label{newscaleA}
N' = N^\chi 
\ee
with $ \chi > 1 $.  

\smallskip

For a matrix $ A \in \matr_E^E  $ we define $ {\rm Diag}(A) := ( \d_{kk'} A_k^{k'})_{k, k' \in E} $.

\begin{proposition} {\bf (Multiscale step \cite{BBo10})} \label{propinvA}
Assume
\be\label{dtCA}
\varsigma \in (0,1/2) \, ,  \ \tau' > 2 \tau + b + 1 \, , \ C_1 \geq 2 \, , 
\ee
and, setting $ \kappa := \tau' + b + s_0 $,
\begin{align} \label{chi1A}
& \chi (\t' - 2 \t - b) >  3 (\kappa + (s_0+ b) C_1 ) \, , \ 
\chi \varsigma > C_1 \, , \\
& \label{s1A}
s_1 > 3 \kappa + \chi (\tau + b) + C_1 s_0 \, . 
\end{align}
For any given $ \Upsilon > 0 $, there exist $\Theta := \Theta (\Upsilon, s_1) > 0 $ large enough 
(appearing in  Definition \ref{regularsA}),  and  $  N_0 (\Upsilon, \Theta ,  s_1) \in \N $  
 such that:
\\[1mm]
$ \forall N \geq N_0(\Upsilon, \Theta  , s_1) $, 
$ \forall E \subset \Z^b \times {\frak I}$ with 
${\rm diam}(E) \leq 4N'=4N^\chi $ (see \eqref{newscaleA}), if $ A \in \matr_E^E $ satisfies 
\begin{itemize}
\item 
{\bf (H1)} $\norsone{A- {\rm Diag}(A)} \leq \Upsilon $ 
\item 
{\bf (H2)} $ \| A^{-1} \|_0 \leq (N')^{\tau}$
\item
{\bf (H3)} 
There is a partition of the  $(A,N)$-bad sites $ B = \cup_{\alpha} \Omega_\alpha$ with
\be\label{sepabadA}
{\rm diam}(\Omega_\alpha) \leq N^{C_1} \, , \quad {\rm d}(\Omega_\alpha , \Omega_\beta) \geq N^2 \ , \ \forall \alpha \neq \beta \, ,
\ee
\end{itemize}
then $ A $ is $ N' $-good. More precisely
\be\label{A-1altaA}
\forall s \in [s_0,s_1]  \ , \  \  \nors{A^{-1}} \leq \frac{1}{4} ({N'})^{\tau' } \Big( ({N'})^{\varsigma s}+ \nors{A- {\rm Diag}(A)} \Big) \, , 
\ee
and, for all $ s \geq s_1 $, 
\be\label{multi:s_APP}
|A^{-1} |_s \leq C(s) ({N'})^{\tau' } \big( ({N'})^{\varsigma s}+ | A- {\rm Diag}(A) |_{s} \big) \, .
\ee
\end{proposition}

The above proposition says, roughly, the following. 
If $ A $ has a sufficient off-diagonal decay (assumption (H1) and \eqref{s1A}), 
and if the sites that can not be inserted in  good ``small''  submatrices (of size $ O(N) $) along the diagonal of $ A $  are
sufficiently separated (assumption (H3)), then the $ L^2 $-bound $ (H2) $ for $ A^{-1} $ implies that the ``large'' matrix $ A $ 
(of size $ N' = N^\chi $ with $ \chi $ as in (\ref{chi1A})) is good, and 
$ A^{-1} $ satisfies also  the bounds \eqref{A-1altaA} in $ s $-norm for $ s > s_1 $.
Notice that the bounds for $ s > s_1 $
follow only by informations  on the $N$-good submatrices  in $ s_1 $-norm 
(see Definition \ref{goodmatrixA}) plus  the $s$-decay of $ A $. 
The link between the various constants is the following: 
\begin{itemize}
\item
According to (\ref{chi1A}) the exponent $ \chi $, which measures the new scale $ N'  \gg  N $, 
is large with respect to the size of the bad clusters $ \O_\a $, i.e. with respect to $ C_1 $. 
The intuitive meaning is that, for $ \chi $ large enough,
the ``resonance effects" due to the bad clusters are ``negligible" at the new larger scale. 
\item
The constant $ \Theta \geq 1 $ which defines the regular sites (see Definition \ref{regularsA}) 
must be large enough with respect to $ \Upsilon $, 
i.e. with respect to  the off diagonal part $ {\cal T} := A - {\rm Diag}(A) $, 
see (H1) and Lemma \ref{defmatrMN}.
\item
Note that $ \chi $ in (\ref{chi1A}) can be taken  
large independently of $ \t $, choosing, for example, $ \t' := 3 \t + 2b $.
\item
The Sobolev  
index $ s_1 $ has to be  large with respect to $ \chi $ and $ \t $, according to (\ref{s1A}).
This is also natural:  if the decay is sufficiently strong, 
then the ``interaction" between different clusters of $ N$-bad sites is weak enough.
\item
In (\ref{sepabadA})  we have fixed the separation $ N^2 $ between the bad clusters 
just for definiteness:
any separation $ N^\mu $, $ \mu > 0 $, would be sufficient. Of course,
the smaller $ \mu > 0 $ is, the larger the Sobolev exponent $ s_1 $ has to be.  
\end{itemize}

The proof of Proposition \ref{propinvA} is divided in several lemmas. In each of them we shall assume that 
the hypotheses of Proposition \ref{propinvA} are satisfied.  We set 
\be\label{AD+T}
{\cal T}  := A - {\rm Diag}(A) \, , 
\qquad \norsone{{\cal T}} \stackrel{(H1)} \leq \Upsilon \, . 
\ee
Call $G$ (resp. $B$) the set of the $(A,N)$-good (resp. bad) sites.
The partition
$$
 E = B \cup G 
 $$ 
 induces the orthogonal decomposition 
$$ 
{\mathcal H}_E = {\mathcal H}_B \oplus {\mathcal H}_G 
$$ 
and we write 
$$ 
u = u_B + u_G \qquad {\rm where}  \qquad u_B := \Pi_B u \, ,  \ u_G := \Pi_G u \, .
$$
We shall denote by $ I_G $, resp. $ I_B $, the restriction of the identity matrix to 
$ {\mathcal H}_G $, resp. 
$ {\mathcal H}_B $,  according to \eqref{defMBC}. 

The next Lemmas \ref{defmatrMN} and \ref{defAprime} 
say that the 
system $ A u = h $
can be nicely reduced along the good sites $ G $, giving rise to a (non-square) system $ A' u_B = Z h $, 
with a good control of the $ s $-norms of the matrices $ A' $ and $ Z $.
Moreover $ A^{-1} $ is a left inverse of $ A' $.

\begin{lemma} \label{defmatrMN}
{\bf (Semi-reduction on the good sites)} 
Let $ \Theta^{-1} \Upsilon \leq c_0 (s_1) $ be small enough. 
There exist $ {\cal M} \in \matr^E_G $, $ {\cal N} \in \matr^B_G $ satisfying, if $ N \geq N_1( \Upsilon)$ is large enough, 
\be  \label{Nm}
\norso{{\cal M}} \leq  c N^\kappa \, , \quad \norso{{\cal N}} \leq c \, \Theta^{-1} \Upsilon \, , 
\ee
for some $ c := c(s_1) > 0 $, and, $ \forall s \geq s_0 $, 
\be\label{Nmalta}
\begin{aligned}
& \nors{{\cal M}} \leq C(s) N^{2\kappa } 
\big( N^{s-s_0}+ N^{-b}\norma {\cal T} \norma_{s+b} \big) \, , \\
& \nors{{\cal N}}  \leq C(s) N^{\kappa } \big( 
N^{s-s_0}+ N^{-b} \norma {\cal T} \norma_{s+b} \big) \, , 
\end{aligned}
\ee
such that 
$$
Au = h \quad \Longrightarrow \quad u_G = {\cal N} u_B + {\cal M} h \, . 
$$
Moreover
\be \label{redex}
u_G = {\cal N} u_B + {\cal M} h \quad \Longrightarrow \quad
\forall k \ {\rm regular} \, , \ (Au)_k = h_k \, .
\ee
\end{lemma}

\begin{pf} It is based on ``resolvent identity" arguments. 
\\[1mm]
{\bf Step I.}  {\it There exist $ \Gamma , L  \in \matr^E_G $ satisfying  
\be\label{Gammas0} 
\norso{\Gamma} \leq  C_0 (s_1) \Theta^{-1} \Upsilon  \, , \quad \norso{L} \leq  N^{\kappa} \, , 
\ee 
and, $ \forall s \geq s_0 $, 
\be \label{Gammam}
\nors{\Gamma}  \leq C(s) N^\kappa 
\big( N^{s-s_0} + N^{-b} \norma {\cal T} \norma_{s+b} \big) \, , 
\quad \nors{L}  \leq  C(s) N^{\kappa +  s -s_0} \, , 
\ee
such that} 
\be\label{Au=h}
Au = h \quad  \Longrightarrow \quad u_G+ \Gamma u = Lh \, .
\ee
Fix any $ k \in G $ (see Definition \ref{ANregA}). If $ k $  is regular, let $ F := \{k \} $,  
and,  if $  k $ is not regular but $(A,N)$-regular, 
let $ F \subset E  $ such  that  $ {\rm d}(k, E\backslash F) \geq N $, $ {\rm diam}(F) \leq 4N $,
$ A_F^F  $ is  $N$-good.
We have
\be
\begin{aligned} \label{eqF}
A u = h  \quad  & \Longrightarrow \quad
A^F_F u_F + A^{E\backslash F}_F u_{E\backslash F} = h_F \\
&   \Longrightarrow \quad
u_F + Q u_{E\backslash F} = (A_F^F)^{-1} h_F 
\end{aligned}
\ee
where 
\be\label{defB}
Q := (A_F^F)^{-1} A^{E\backslash F}_F=(A_F^F)^{-1} {\cal T}^{E\backslash F}_F \in \matr_F^{E\backslash F} \, . 
\ee
The matrix $ Q $ satisfies 
\be \label{estB}
\norsone{Q} \stackrel{(\ref{algebra})} \leq C(s_1) \norsone{(A_F^F)^{-1}} \norsone{{\cal T}} 
\stackrel{(\ref{NgoodmatA}), (\ref{AD+T})}  \leq C(s_1) N^{\tau' + \varsigma s_1} \Upsilon
\ee
(the matrix  $ A_F^F  $ is  $N$-good).
Moreover, $ \forall s \geq s_0 $, using (\ref{interpm}) 
and  $ {\rm diam}(F) \leq 4N $, 
\begin{align}
\norma Q \norma_{s+b} & \lesssim_s 
 \norma (A_F^F)^{-1} \norma_{s+b} \norso{{\cal T}} + \norso{(A_F^F)^{-1}} \norma {\cal T} \norma_{s+b}  \nonumber \\
&  \stackrel{(\ref{Sm2})} {\lesssim_s}  
 N^{s+b -s_0} \norma (A_F^F)^{-1} \norma_{s_0} \norso{{\cal T}} + \norso{(A_F^F)^{-1}} \norma {\cal T} \norma_{s+b}  \nonumber \\ 
&  \stackrel{(\ref{NgoodmatA}), (\ref{AD+T})} {\lesssim_s}  
N^{(\varsigma -1) s_0 } 
\big( N^{s+b+ \tau' } \Upsilon + N^{\tau'+   s_0} \norma {\cal T} \norma_{s+b} \big)  \,. \label{Balta}
\end{align}
Applying the projector $ \Pi_{\{k\}} $ in (\ref{eqF}), we obtain
\be\label{primos}
Au=h \quad  \Longrightarrow \quad 
u_k + \sum_{k' \in E} \Gamma_k^{k'} u_{k'}= \sum_{k'\in E} L_k^{k'} h_{k'}  
\ee
that is  (\ref{Au=h})  with
\be\label{defGL}
\begin{aligned}
& \Gamma_k^{k'}  := 
\begin{cases} 0 \	 \  \ \quad{\rm if} \ \ k' \in F \\
Q_k^{k'}  \ \  \ \, {\rm if} \ \ k' \in E \setminus F\, , 
\end{cases}
\\
&  
L_k^{k'} := \begin{cases}
[(A_F^F)^{-1}]_k^{k'} \  \ \ \, {\rm if} \ \ k' \in F \\
0 \quad  \qquad \quad \ \,  \ \ {\rm if} \ \ k' \in E \setminus F. 
\end{cases}  
\end{aligned}
\ee
If $ k $ is regular then $ F = \{ k \} $, and, by Definition \ref{regularsA}, 
\be\label{invAi}
| A_k^k | \geq \Theta \, .
\ee
Therefore, by (\ref{defGL}) and (\ref{defB}), the $ k $-line of $ \Gamma $ satisfies 
\be\label{lin1}
\norsob{\Gamma_k} \leq \norsob{(A_k^k)^{-1} {\cal T}_k} \stackrel{ (\ref{invAi}), (\ref{AD+T})} {\lesssim_{s_0}} 
 \Theta^{-1} \Upsilon \, .
\ee
If $ k $ is not regular but $(A,N)$-regular, since  ${\rm d}(k, E \backslash F) \geq N $ 
we have, by (\ref{defGL}), that $\Gamma_k^{k'}=0$ for 
$ | k - k' | \leq N $. Hence, by Lemma \ref{norcompA}, 
\begin{align}
\norsob{\Gamma_k}  \stackrel{(\ref{Sm1A})} 
\leq  N^{-(s_1-s_0-b)} \norsone{\Gamma_k} 
& \stackrel{(\ref{defGL})} \leq   N^{-(s_1-s_0-b)}  \norsone{Q} \nonumber  \\
& \stackrel{(\ref{estB})} {\lesssim_{s_1}}    \Upsilon N^{\tau' + s_0 + b  - (1-\varsigma)s_1} \nonumber \\ 
& \lesssim_{s_1}    \Theta^{-1} \Upsilon  \label{lin2}
\end{align}
for $ N \geq N_0 (\Theta) $ large enough. Indeed
the exponent  $ \tau' + s_0 + b  - (1-\varsigma)s_1 <  0 $ because $ s_1 $ is large enough  
according to (\ref{s1A}) and $ \varsigma \in (0, 1/2) $ (recall $ \kappa := \tau' + s_0 + b $).
In both cases (\ref{lin1})-(\ref{lin2}) imply that each line $ \Gamma_k $ decays like
$$
\norsob{\Gamma_k}  \lesssim_{s_1} \Theta^{-1} \Upsilon \, ,\quad  \forall k \in G \, .
$$ 
Hence, by Lemma \ref{norextracted}, we get 
$$
\norso{\Gamma} \leq C'(s_1) \Theta^{-1} \Upsilon \, , 
$$
which is 
the first inequality in \eqref{Gammas0}. 
Likewise we prove the second estimate in \eqref{Gammas0}.
Moreover, $ \forall s \geq s_0 $, still by Lemma \ref{norextracted}, 
$$
\nors{\Gamma}  \lesssim \sup_{k \in G} \norma \Gamma_k \norma_{s+b}  
\stackrel{\eqref{defGL}} \lesssim \norsb{Q} \stackrel{\eqref{Balta}} {\lesssim_s} N^\kappa 
\big( N^{s-s_0} + N^{-b} \norsb{{\cal T}} \big)  
$$
where  $ \kappa := \tau' + s_0 + b $ and for $ N \geq N_0(\Upsilon)$.

The second estimate in (\ref{Gammam}) follows by $ \norso{L} \leq  N^{\kappa} $ (see (\ref{Gammas0})) and (\ref{Sm2})
(note that  by (\ref{defGL}), since diam$F \leq 4N$, we have $ L_k^{k'} = 0 $ for all $ |k - k' | > 4N $).
\\[2mm]
{\bf Step II.}  By (\ref{Au=h}) we have
\be\label{Auin}
Au = h \quad \Longrightarrow  \quad (I_G  + \Gamma^G)u_G  = Lh - \Gamma^B u_B \, .
\ee
By (\ref{Gammas0}), if $ \Theta $ is large enough (depending on $ \Upsilon$, namely on the potential $ V_0 $), 
we have $\norso{\Gamma^G} \leq 1/2 $. Hence,
by  Lemma \ref{leftinvA}, $ I_G  + \Gamma^G $ is invertible and 
\be\label{Linver0}
\norso{(I_G  + \Gamma^G)^{-1}} \stackrel{(\ref{inv1A})} \leq 2 \, , 
\ee  
and, $ \forall s \geq s_0 $ 
\be\label{Linver}
 \nors{(I_G  + \Gamma^G)^{-1}} 
\stackrel{(\ref{inv1A})} {\lesssim_s}  1+ \nors{\Gamma^G} 
\stackrel{(\ref{Gammam})} {\lesssim_s}  N^{\kappa }
\big(N^{s-s_0}  + N^{-b} \norma {\cal T} \norma_{s+b} \big) \, .
\ee
By (\ref{Auin}), we have 
$$ 
Au=h  \quad \Longrightarrow \quad u_G = {\cal M} h + {\cal N} u_B 
$$ 
with
\be\label{poi}
\begin{aligned} 
& {\cal M} := (I_G  + \Gamma^G)^{-1}L  \, , \\
& {\cal N} := - (I_G  + \Gamma^G)^{-1} \Gamma^B  \, , 
\end{aligned}
\ee
and  estimates (\ref{Nm})-(\ref{Nmalta}) follow by Lemma 
\ref{prodest}, (\ref{Linver0})-(\ref{Linver}) and (\ref{Gammas0})-(\ref{Gammam}).

Note that  
\be\label{conclus}
u_G + \Gamma u = Lh  \quad \iff \quad u_G = {\cal M} h + {\cal N} u_B \, . 
\ee
As a consequence, if $u_G = {\cal M} h + {\cal N} u_B$ then, by
(\ref{defGL}), for $ k $ regular,
$$
u_k + (A_k^k)^{-1} \sum_{k'\neq k} A_k^{k'} u_{k'}=(A_k^k)^{-1} h_k \, ,
$$
hence $ (Au)_k = h_k $, proving \eqref{redex}.
\end{pf}

\begin{lemma} \label{defAprime}
{\bf (Reduction on the bad sites)} We have
$$
Au = h \quad \Longrightarrow  \quad A' u_B = Zh  
$$
where
\be\label{A'Z}
\begin{aligned}
& A' := A^B + A^G {\cal N} \ \in  {\cal M}^B_E \, ,  \\
& Z := I_E - A^G {\cal M} \  \in  {\cal M}^E_E \, ,
\end{aligned}
\ee
satisfy
\be
\begin{aligned}
& \norso{A'} \leq c(\Theta)  \, , \\
&  \nors{A'} \leq  C(s, \Theta) N^{\kappa } (N^{s-s_0}
+ N^{-b}\norma {\cal T} \norma_{s+b})\, , \label{A'a} 
\end{aligned}
\ee
and
\be
\begin{aligned}
& \label{Za}
\norso{Z} \leq c N^\kappa \, , \\
&  \nors{Z} \leq C(s, \Theta) N^{2\kappa } (N^{s-s_0} 
+ N^{-b} \norma {\cal T} \norma_{s+b}) \, .
\end{aligned}
\ee
Moreover $(A^{-1})_B $ is a left inverse\index{Left inverse} of $ A' $. 
\end{lemma}

\begin{pf}
By Lemma \ref{defmatrMN},
$$
\begin{aligned}
Au = h \quad & \Longrightarrow \quad
\begin{cases}
A^G u_G + A^B u_B = h \\
u_G = {\cal N} u_B + {\cal M} h 
\end{cases} \\
&  \Longrightarrow \quad 
(A^G {\cal N} + A^B ) u_B = h-A^G {\cal M} h \, , 
\end{aligned}
$$
{\it i.e.} $ A' u_B = Z h $. Let us prove  estimates (\ref{A'a})-(\ref{Za}) for $ A' $ and $ Z $.
\\[1mm]
{\bf Step I.} {\it $ \forall \, k $ regular we have $ A'_k = 0 $, $ Z_k = 0 $.}
\\[1mm]
By (\ref{redex}), for all $ k $ regular, 
$$ 
\begin{aligned}
 \forall h \, , \  \forall u_B \in {\mathcal H}_B \, , \ & 
\Big(  A^G ({\cal N}u_B + {\cal M}h)+ A^B u_B \Big)_k = h_k \  , \\
& i.e. \quad  (A'u_B)_k = (Zh)_k \,  ,
\end{aligned}
$$
which implies $ A'_k = 0 $ and $ Z_k = 0 $. 
\\[1mm]
{\bf Step II.} {\it Proof of (\ref{A'a})-(\ref{Za})}.
\\[1mm] 
Call $ R \subset E$  the regular sites in $ E $. For all $ k \in E \backslash R $, 
we have $ |A_k^k| < \Theta $ (see  Definition \ref{regularsA}). Then 
(\ref{AD+T})  implies 
\be\label{passagge}
\begin{aligned}
& \norso{A_{E \backslash R}} \leq \Theta + \norso{{\cal T}}  \leq c(\Theta) \, , \\
&  \nors{A_{E \backslash R}} \leq \Theta + \nors{{\cal T}} \, , \quad  \forall s \geq s_0 \, .
\end{aligned}
\ee
By Step I and the definition of $ A' $ in (\ref{A'Z}) we get
$$
\nors{A'}  = \nors{A'_{E\backslash R}} \leq \nors{A^B_{E\backslash R}} + \nors{A^G_{E\backslash R} {\cal N}} \, .
$$
Therefore Lemma \ref{prodest}, (\ref{passagge}), (\ref{Nm}), (\ref{Nmalta}), imply  
$$
\nors{A'} \leq C(s, \Theta) N^{\kappa}  \big(N^{s-s_0} 
+ N^{-b} \norma {\cal T} \norma_{s+b} \big)  
\quad {\rm and} \quad \norso{A'} \leq c(\Theta) \, ,
$$
proving (\ref{A'a}). The bound (\ref{Za}) follows similarly.
\\[1mm]
{\bf Step III.} {\it $ (A^{-1})_B $ is a left inverse of $ A' $}.
\\[1mm] 
By 
$$ 
A^{-1} A' =  A^{-1} ( A^B + A^G {\cal N}) =  I^B_E + I^G_E {\cal N} 
$$ 
we get
$$
(A^{-1})_B A'=(A^{-1} A')_B=I^B_B - 0= I^B_B
$$ 
proving  that $(A^{-1})_B$ is a left inverse of $ A' $. 
\end{pf}

Now $ A' \in \matr^B_E $, and the set $B$ is partitioned in clusters $\Omega_\a$ of size
$O(N^{C_1})$, far enough one from another, see (H3).  Then, up to a remainder of very small $ s_0 $-norm 
(see (\ref{noroR0})), 
$ A' $ is defined by the submatrices $(A')^{\Omega_\a}_{\Omega'_\a}$ where $\Omega'_\a$ is some neighborhood
of $\Omega_\a$ (the distance between two distinct $\Omega'_\a$ and $\Omega'_\b$ remains large). Since
$ A' $ has a left inverse with $L^2$-norm $O({N'}^\tau)$, 
so have the submatrices $(A')^{\Omega_\a}_{\Omega'_\a}$. Since these
submatrices are of size $O(N^{C_1})$, the $s$-norms of their inverse will be estimated as
$O(N^{C_1s} {N'}^\tau)= O({N'}^{\tau + \chi^{-1} C_1 s})$, see (\ref{decaW0}). 
By Lemma \ref{leftinvA},  provided $ \chi $ is chosen large enough,
$ A' $ has a left inverse $ V $ with $s$-norms satisfying (\ref{LeftY}).
The details are given in the following lemma.

\begin{lemma} \label{defY} {\bf (Left inverse with decay)}
The matrix $ A' $ defined in Lemma \ref{defAprime} has a left inverse $ V $ which satisfies 
\be\label{LeftY}
\forall s \geq s_0 \ , \  \
\nors{V}\lesssim_s 
{N}^{2\chi \tau + \kappa + 2 (s_0+b) C_1} \big( N^{C_1 s} + \norsb{{\cal T}} \big)  \, . 
\ee
\end{lemma}

\begin{pf}
Define $ \dom \in \matr_E^B $ by
\be\label{grossodia}
\dom_{k'}^k :=\begin{cases}
(A')_{k'}^{k} \quad \, \rm{if}  \  \ ({\it k}, {\it k}') \in \cup_\a (\Omega_\alpha \times \Omega'_\alpha) \\
0 \qquad \quad \, \hbox{if}  \  \ ({\it k}, {\it k}') \notin \cup_\a (\Omega_\alpha \times \Omega'_\alpha) 
\end{cases}
\ee
where
\be \label{grossodia1}
\Omega'_\alpha := \big\{ k \in E \  : \ {\rm d}(k,\Omega_\alpha) \leq N^2/4 \big\} \, .
\ee
{\bf Step I.} {\it $ \dom $ has a left inverse $ W \in \matr^E_B $ with $ \| W \|_0 \leq 2 ({N'})^{\tau} $.}
\\[1mm]
We define $ \remain := A' - \dom $. By the definition (\ref{grossodia})-\eqref{grossodia1}, 
if  $ {\rm d}(k',k) < N^2 / 4 $ then $ \remain_{k'}^k = 0 $ and so
\begin{align} 
\norso{\remain} & 
\stackrel{\eqref{Sm1A}} \leq   
4^{s_1} N^{-2(s_1 - b - s_0)} \norma{\remain} \norma_{s_1 - b}  
\leq 4^{s_1} N^{-2(s_1- b - s_0)} \norma{A'} \norma_{s_1 - b} \nonumber  \\
&  \stackrel{(\ref{A'a}), (\ref{AD+T})}  {\lesssim_{s_1}}    N^{-2(s_1 - b - s_0)} N^{\kappa} ( N^{s_1-b-s_0}+N^{-b} 
\Upsilon) \nonumber \\
&  \lesssim_{s_1}   N^{2 \kappa - s_1}  \label{noroR0}
\end{align}
for $ N \geq N_0 (\Upsilon) $ large enough. 
Therefore 
\begin{eqnarray}\label{rsmall}
\| \remain \|_0 \| (A^{-1})_B \|_0 \stackrel{\eqref{schur}} {\lesssim_{s_0}}  \norso{\remain} \| A^{-1} \|_0 
& \stackrel{(\ref{noroR0}),(H2)}   {\lesssim_{s_1}}  &   N^{2\kappa - s_1} (N')^{\t} \nonumber \\
& \stackrel{(\ref{newscaleA})} = & C(s_1) N^{2\kappa  - s_1+ \chi \t}  \nonumber \\
&  \stackrel{(\ref{s1A})} \leq 1/2 
\end{eqnarray}
for $ N \geq N(s_1) $.
Since $ (A^{-1})_B  \in \matr^E_B $ is a left inverse of $ A' $ (see Lemma \ref{defAprime}), Lemma 
\ref{leftinvA} and (\ref{rsmall}) imply that  
$ \dom = A' - R $ has a left inverse $ W \in \matr^E_B $, and
\be\label{W00}
\| W \|_0 \stackrel{(\ref{inv10})} \leq 2  \| (A^{-1})_B \|_0 \leq  2  \| A^{-1} \|_0 \stackrel{(H2)}\leq
2 ({N'})^{\tau} \, . 
\ee
{\bf Step II.} {\it $ W_0 \in \matr^E_B $ defined by
\be\label{W0ii}
(W_0)^{k'}_k := \begin{cases}  W^{k'}_k \quad  \quad {\rm if} \ \ \, (k,k') \in \cup_\a
(\Omega_\alpha \times \Omega'_\alpha) \\ 
0 \, \qquad \quad \, \hbox{\rm if } \ \ (k,k') \not\in \cup_\a
(\Omega_\alpha \times \Omega'_\alpha) \end{cases}
\ee
is a left inverse of $ \dom $ and $ \nors{W_0} \leq C(s) {N}^{(s+b)C_1 + \chi \tau  } $, $ \forall s \geq s_0 $.}
\\[1mm]
Since $ W \dom = I_B $, we prove that $ W_0 $ is a left inverse of $ \dom $ showing that 
\be\label{vai}
(W-W_0)\dom = 0 \, .
\ee
Let us prove (\ref{vai}). 
For $ k \in B = \cup_\a \Om_\a $, there is $\a$ such that  $k \in \Omega_\alpha $,  and
\be\label{sum}
\forall k' \in B \ , \ \big( (W-W_0)\dom \big)_k^{k'}=
\sum_{q \notin \Omega'_\a} (W-W_0)_k^q \dom^{k'}_q
\ee
since $ (W - W_0)_k^q = 0 $ if $ q \in \Omega'_\a $, see the Definition (\ref{W0ii}). 
\\[1mm]
{\sc Case I:} $ k' \in \Om_\a $. Then $ \dom^{k'}_q = 0 $  in (\ref{sum}) and so $((W-W_0)\dom )_k^{k'} = 0 $.
\\[1mm]
{\sc Case II:} $ k' \in \Omega_\beta $ for some $\beta \neq \a $. Then,
since $ \dom^{k'}_q = 0 $ if $ q \notin \Omega'_\beta$, we obtain by (\ref{sum}) that
$$
\begin{aligned}
 ((W-W_0)\dom )_k^{k'} & =\sum_{q \in \Omega'_\beta} (W-W_0)_k^q \dom^{k'}_q
\stackrel{(\ref{W0ii})} = \sum_{q \in \Omega'_\beta} W_k^q \dom^{k'}_q \\
&  \stackrel{(\ref{grossodia})}  =
\sum_{k \in E} W_k^q \dom^{k'}_q
= (W \dom)_k^{k'} = (I_B)_k^{k'} = 0 \, .
\end{aligned}
$$
Since $ {\rm diam}(\Omega'_\a) \leq 2N^{C_1} $,  definition (\ref{W0ii}) implies 
$ (W_0)_k^{k'} = 0 $ for all $ | k - k' | \geq 2N^{C_1} $.  Hence,  $ \forall s \geq 0 $,
\be\label{decaW0}
\nors{W_0} \stackrel{(\ref{Sm2})} {\lesssim_s} N^{(s + b) C_1 } \| W_0 \|_0
\stackrel{(\ref{W00})} {\lesssim_s} {N}^{(s + b)C_1 +\chi \t } \, . 
\ee
{\bf Step III.} {\it $ A' $ has a left inverse $ V $ satisfying (\ref{LeftY}).}
\\[1mm]
Now $ A' = \dom + \remain $, $ W_0 $ is a left inverse of $ \dom $, and 
$$
\norso{W_0} \norso{\remain} 
\stackrel{(\ref{decaW0}), (\ref{noroR0})} 
\leq C(s_1) N^{(s_0 + b) C_1 + \chi \tau + 2 \kappa -  s_1} \stackrel{(\ref{s1A})} \leq 1/2
$$
(we use  also that $ \chi > C_1 $ by (\ref{chi1A}))
for $ N \geq N(s_1) $  large enough. Hence, by Lemma  \ref{leftinvA}, $A'$ has a left inverse $ V $ with
\be\label{Vs0}
\norso{V} \stackrel{(\ref{inv1A})} \leq 2 \norso{W_0} \stackrel{(\ref{decaW0})} \leq C N^{ (s_0 + b)C_1 + \chi \tau }
\ee
and, $ \forall s \geq s_0 $, 
\begin{eqnarray*}
\nors{V} & \stackrel{(\ref{inv1A})} {\lesssim_s} &   \nors{W_0} + \norso{W_0}^2  \nors{\remain} 
\lesssim_s \nors{W_0} + \norso{W_0}^2  \nors{A'}  \\
& \stackrel{(\ref{decaW0}), (\ref{A'a})} {\lesssim_s} & 
{N}^{2\chi \tau + \kappa + 2 (s_0+b) C_1} \big( N^{C_1 s} + \norsb{{\cal T}} \big)
\end{eqnarray*}
proving (\ref{LeftY}).
\end{pf}

\noindent
{\sc Proof of Proposition \ref{propinvA} completed.} 
Lemmata \ref{defmatrMN}, \ref{defAprime}, \ref{defY} imply
$$
Au = h  \quad \Longrightarrow \quad 
\begin{cases}
u_G = {\cal M} h + {\cal N}u_B \\ 
u_B = V Z h 
\end{cases}  
$$
whence 
\be\label{eccoci}
\begin{aligned}
& (A^{-1})_B=  V Z \\  
& (A^{-1})_G= {\cal M} + {\cal N} V Z={\cal M} + {\cal N}(A^{-1})_B \, . 
\end{aligned}
\ee
Therefore, $ \forall s \geq s_0 $, 
\begin{eqnarray}
\nors{(A^{-1})_B} & \stackrel{(\ref{eccoci}), (\ref{interpm})} 
{\lesssim_s} &  \nors{V} \norso{Z} + \norso{V} \nors{Z}  \nonumber \\
& \stackrel{(\ref{LeftY}), (\ref{Za}), (\ref{AD+T}), \eqref{Vs0}} {\lesssim_s} &  N^{2\kappa + 2 \chi \tau + 2 (s_0+b)C_1 } 
\big( N^{C_1 s} + \norsb{{\cal T}} \big) \nonumber  \\
& {\lesssim_s} & {(N')}^{\a_1} \big( {(N')}^{\a_2 s} + \nors{{\cal T}} \big) \nonumber 
\end{eqnarray}
using $ \norsb{{\cal T}} \leq C(s) (N')^b \nors{{\cal T}} $ (by (\ref{Sm2})) and defining
\be\label{DEFAL1AL2}
\a_1 := \, 2 \tau + b + 2 \chi^{-1}( \kappa +  C_1 (s_0 + b))\, , \quad  \a_2 := \chi^{-1} C_1 \, .  
\ee
We obtain the same bound for $ \nors{(A^{-1})_G} $.
Notice that by \eqref{chi1A} and \eqref{dtCA}, the exponents 
$ \a_1 , \a_2 $ in \eqref{DEFAL1AL2} satisfy 
\be\label{pro:al1al2}
\a_1 < \tau' \, , \quad  \a_2 < \varsigma \, . 
\ee
 Hence, for all  $s \geq s_0 $, 
\begin{align}
\nors{A^{-1} }  \leq  \nors{(A^{-1})_B} + \nors{(A^{-1})_G} 
& \leq  C(s) {(N')}^{\a_1} \big( {(N')}^{\a_2 s}+\nors{{\cal T}} \big) \label{newA-1s} \\
& \stackrel{\eqref{pro:al1al2}, \eqref{AD+T}} \leq 
C(s) {(N')}^{ \tau'} \big( {(N')}^{\varsigma s}+\nors{A - {\rm Diag}(A)} \big) \nonumber 
\end{align}
which is \eqref{multi:s_APP}. Moreover, 
for $ N \geq N(s_1) $ large enough, we have 
$$
\forall s \in [s_0, s_1] \, , \quad  C(s) {(N')}^{\a_1} \leq  \frac14 {(N')}^{\t'} \, , 
$$
and by \eqref{newA-1s} we deduce \eqref{A-1altaA}.

\chapter{Normal form close to an isotropic torus}\label{sec:2}

In this appendix we report the results in \cite{BBField}, which are used in Chapter \ref{sezione almost approximate inverse}.
Theorem \ref{thm1} provides, in a neighborhood of an isotropic invariant torus for an Hamiltonian vector field $ X_K $, 
symplectic variables in which the Hamiltonian $ K $ assumes the normal form \eqref{normalformH}. 
It is a classical result  of Herman \cite{Herman}, \cite{HF} 
that an invariant torus, densily 
filled by a quasi-periodic solution, is isotropic, see 
Lemma \ref{lem:iso} and Lemma \ref{lem:approclo}  for a more quantitative version. 
In view of the Nash-Moser iteration 
we need to perform an analogous construction for an only 
``approximately invariant" torus. The key step is Lemma \ref{modified} which  constructs, 
near 
an ``approximately invariant" torus,  an 
isotropic torus. This  appendix is written with  a self-contained character. 

\section{Symplectic coordinates near an invariant torus}

We consider the toroidal phase space
$$ 
{\cal P} := \T^\nu \times \R^\nu \times E \qquad {\rm where} \qquad  \T^\nu := \R^\nu \slash (2 \pi \Z)^\nu 
$$
is the standard flat torus and 
$ E $ is a real Hilbert space with scalar product $ \langle \  , \ \rangle $. 
We denote by $  u := (\teta, y, z )  $ the variables of $ {\cal P} $. 
We call $ (\teta, y ) $ the ``action-angle" variables and
$ z $  the ``normal" variables. 
We assume that $ E $ is endowed with a constant exact symplectic $ 2 $-form
\be\label{2foA}
\Omega_E (z,w) =  \langle \bar J z , w \rangle  \, , \quad \forall z, w \in E \, , 
\ee
where $ \bar J : E \to E  $ is an antisymmetric bounded linear operator with trivial kernel.
Thus $ {\cal P} $ is endowed with the symplectic $ 2 $-form
\be\label{2formA}
{\cal W} := (d y \wedge d \teta) \, \oplus \,  \Omega_E  
\ee
which is exact, namely
\be\label{exact}
{\cal W} = d \form 
\ee
where $ d$ denotes the exterior derivative and $ \form $ is the  Liouville $ 1 $-form\index{Liouville $ 1$-form} on $ {\cal P} $ defined by 
\be\label{lambda}
\begin{aligned}
& \qquad \qquad \qquad \quad \form_{(\teta, y, z)} : \R^\nu \times \R^\nu \times E \to \R \, ,  \\
& \form_{(\teta, y, z)}[ \hat \teta, \hat y, \hat z] := y \cdot \hat \teta + \frac12 \langle \bar J z, \hat z \rangle \, ,
\quad \forall (\hat \teta, \hat y, \hat z)  \in \R^\nu \times \R^\nu \times E \, ,
\end{aligned}
\ee
and the dot $ ``\cdot" $ denotes the usual scalar product of $ \R^\nu $. 

Given a Hamiltonian  function $ K :   {\cal D} \subset
 {\cal P } \to \R $, 
we consider the Hamiltonian system 
\be\label{HS0A}
u_t = X_K (u)  \qquad {\rm where} \qquad d K (u) [ \cdot ] = - {\cal W} ( X_K(u) , \, \cdot ) 
\ee
formally defines the Hamiltonian  vector field  $  X_K $. 
For infinite dimensional systems (i.e. PDEs) the Hamiltonian  $ K $ is, in general, 
well defined  and smooth  only on a dense subset 
$  {\cal D} =  \T^\nu \times \R^\nu \times E_1 \subset {\cal P} $ where  $ E_1 \subset E $
is a dense subspace of $ E $.  We require that, for all $ (\theta, y ) \in \T^\nu \times \R^\nu $, $ \forall z \in E_1 $, 
the Hamiltonian $ K $ admits 
a gradient  $\nabla_z K $, defined by 
\be\label{gradient}
d_z K(\theta, y, z ) [h] = \langle \nabla_z K (\theta, y, z), h \rangle \, , \quad \forall  h  \in E_1  \, , 
\ee
and that $\nabla_z K (\theta, y, z) \in E$ is in the space of definition of the (possibly unbounded)
operator  $ J := - {\bar J}^{-1} $.   
Then by \eqref{HS0A}, \eqref{2foA}, \eqref{2formA}, \eqref{gradient} the Hamiltonian vector field  
$ X_K  : \T^\nu \times \R^\nu \times E_1  $ $ \mapsto \R^\nu \times \R^\nu \times E  $
writes
$$ 
X_K  = (  \partial_y K , - \partial_\theta K,  J \nabla_z K ) \, , \qquad J : = - {\bar J}^{-1} \, . 
$$
A continuous curve 
 $ [t_0, t_1] \ni t \mapsto u (t) \in \T^\nu \times \R^\nu \times E $ is a solution of the Hamiltonian system \eqref{HS0A}
 if it is $ C^1 $ as a map from $ [t_0, t_1]$ to $ \T^\nu \times \R^\nu \times E_1 $ and  
$  u_t (t) = X_K ( u(t) )$, $ \forall t \in [t_0, t_1 ] $.  
For PDEs, the flow map  $ \Phi^t_K $ may not be well-defined everywhere. 
The next arguments, however,  will not require to solve the initial value problem, but only
a functional equation in order to find quasi-periodic solutions, 
see \eqref{Boundary value}.

\smallskip

We suppose that 
 \eqref{HS0A} possesses an embedded invariant torus\index{Torus embedding} 
\begin{align}\label{embedded torus}
&  \qquad \vphi \mapsto  i ( \vphi ) := (  \uth (\vphi), \uy (\vphi), \uz (\vphi) )  \\
& \label{smooth torus}
i \in C^1 (\T^\nu,  {\cal P}) \cap C^0 (\T^\nu,   {\cal P} \cap \T^\nu \times \R^\nu \times E_1) \, , 
\end{align}
which supports a quasi-periodic\index{Quasi-periodic solution} solution with non-resonant   frequency vector $ \om  \in \R^\nu $. More precisely we have that 
\be\label{invariatorus}
 i \circ \Psi^t_\om = \Phi^t_K \circ i   \, , \quad  \forall t \in \R \, , 
\ee
where $ \Phi^t_K $ denotes the flow generated by $ X_K $  and 
\be\label{linear flow}
\Psi^t_\om : \T^\nu \to \T^\nu \, , \quad \Psi^t_\om (\vphi ) := \vphi + \om t  \, , 
\ee
is  the translation flow of vector $\om$ on $ \T^\nu $.
Since $ \om \in \R^\nu $ is non-resonant, namely 
$$ 
\om \cdot \ell  \neq 0 \, , \quad  \forall \ell  \in \Z^\nu \setminus \{0\} \, , 
$$ 
each orbit of $ (\Psi^t_\om) $ is {\it dense} in $ \T^\nu $. 
Note that \eqref{invariatorus} {\it only} requires that the flow $ \Phi^t_K $ is well defined and smooth on
the compact manifold  $ {\cal T} := i ( \T^\nu )  \subset {\cal P } $ and 
 $ (\Phi^t_K)_{|{\cal T}} = i \circ \Psi^t_\om \circ i^{-1} $. This remark is important because, for PDEs, 
 the flow  could be ill-posed in a neighborhood of $ {\cal T} $. 
From a functional point of view \eqref{invariatorus} is equivalent to the equation
\be\label{Boundary value}
\om \cdot \pa_\vphi  i ( \vphi) - X_K ( i( \vphi )) = 0 \, . 
\ee
\begin{remark}
In the sequel we will formally differentiate several times the torus embedding $i$,
so that we assume  more regularity than \eqref{smooth torus}. In the framework
of a Nash-Moser scheme, the approximate torus embedding solutions $ i $ are indeed regularized 
at each step.
\end{remark}

We require that $ \uth : \T^\nu \to \T^\nu $ is a diffeomorphism of $ \T^\nu $ isotopic to the identity.
Then the embedded torus $ {\cal T} := i(\T^\nu )$ is a smooth graph over $ \T^{\nu } $. Moreover
the lift on $ \R^\nu $ of $\uth$ is a map
\be\label{inverse teta0}
\uth : \R^\nu \to \R^\nu , \quad 
\uth (\vphi) = \vphi  + \uTh ( \vphi ) \, ,  
\ee
where $  \uTh (\vphi) $ is $ 2 \pi $-periodic 
in each component $ \vphi_i  $, $ i=1, \ldots, \nu $,  
with invertible Jacobian matrix 
$$ 
D \uth (\vphi ) = {\rm Id} + D \uTh (\vphi ) \, , \quad  \forall \vphi \in \T^\nu \,.
$$
In the usual applications  $ D \uTh  $ is small and $ \om  $ is a Diophantine\index{Diophantine vector} vector, namely 
$$
| \om \cdot \ell  | \geq \frac{\g}{|\ell |^{\tau}} \, , \quad  \forall \ell  \in \Z^{\nu} \setminus \{0\} \, .
$$
The torus $  {\cal T} $ is  the graph of the function  (see \eqref{embedded torus} and \eqref{inverse teta0})
\be\label{torus graph0}
j := i \circ \uth^{-1} \, , \quad
j : \T^\nu \to \T^\nu \times \R^\nu \times E \, , \quad 
j (\teta ) := ( \teta,    \tilde \uy (\teta) , {\tilde \uz} (\teta)) \, ,
\ee
namely 
\be\label{torus graph}
{\cal T} = \big\{ ( \theta,  \tilde \uy (\teta), {\tilde \uz}(\teta)) \ ; \ \theta \in \T^\nu \big\} \, ,
\quad {\rm where} \quad   \tilde \uy  := \uy \circ \uth^{-1}  \, ,  \
{\tilde \uz}  := \uz \circ \uth^{-1} \, .
\ee
We first prove the isotropy of an invariant torus as in  \cite{Herman},
\cite{HF},
i.e. that the $ 2 $-form  
$ {\cal W} $ vanishes\index{Isotropic torus} on the tangent space to $ i ( \T^\nu ) \subset {\cal P }  $, 
\be\label{i star closed}
0 = i^* {\cal W} =  i^* d \form =   d (i^* \form) \, , 
\ee
or equivalently  the $ 1$-form $ i^* \form $ on $ \T^\nu $  is closed. 

\begin{lemma}\label{lem:iso}
The invariant torus $ i ( \T^\nu ) $ is isotropic.
\end{lemma}

\begin{pf}  
By \eqref{invariatorus} the pullback
\be\label{pull1}
( i \circ \Psi^t_\om )^* {\cal W} = (\Phi^t_K \circ i )^*  {\cal W} = 
  i^*  {\cal W}  \, . 
\ee
For smooth Hamiltonian systems in finite dimension $(\ref{pull1})$ is true because the 
$2$-form $ {\cal W} $ is invariant under the flow map $ \Phi^t_K $ (i.e. $(\Phi^t_K)^* {\cal W} = {\cal W} $). In our setting, the flow
$(\Phi^t_K)$ may not be defined everywhere, but 
 $ \Phi^t_K  $ is well defined  on $  i ( \T^\nu ) $ by the assumption \eqref{invariatorus},
 and still preserves $ {\cal W} $ on the manifold $  i ( \T^\nu ) $, see the proof of Lemma 
 \ref{lem:approclo} for details.
Next, denoting  the $ 2 $-form 
$$ 
(i^*  {\cal W})(\vphi) = \sum_{i<j} A_{ij} (\vphi) d \vphi_i \wedge d \vphi_j  \, , 
$$ 
we have  
$$
( i \circ \Psi^t_\om )^* {\cal W} = ( \Psi^t_\om )^* \circ i^* {\cal W} 
=   { \sum}_{i<j} A_{ij} (\vphi + \om t ) d \vphi_i \wedge d \vphi_j 
$$
and so \eqref{pull1} implies that 
$$ 
A_{ij} (\vphi + \om t ) = A_{ij} (\vphi ) \, , \quad   \forall t \in \R \, .
$$ 
Since 
the orbit  $ \{ \vphi + \om t \}$ is dense on $ \T^\nu $ ($ \om $ is non-resonant) 
and each function $ A_{ij}$ is continuous,  
it implies that 
$$ 
A_{ij} (\vphi ) = c_{ij}  \, , \ \forall \vphi \in \T^\nu \, , \quad i.e. \ \  
i^*  {\cal W}  = { \sum}_{i<j} c_{ij} d \vphi_i \wedge d \vphi_j  
$$ 
is constant.
But, by \eqref{exact}, the $ 2$-form  $ i^*  {\cal W}   = i^*  d \form = d (i^* \form )  $ 
is also exact. Thus each
$ c_{ij} = 0 $ namely $ i^*  {\cal W} = 0 $. 
\end{pf}
  
 We now consider the diffeomorphism of the phase space
\be\label{cG} 
\left(
\begin{array}{c}
\teta  \\
y   \\
z
\end{array}
\right) = 
G
\left(
\begin{array}{c}
\phi \\
\zeta  \\
w
\end{array}
\right) 
:= 
\left(
\begin{array}{c}
\!\! \!\! \!\!\!\! \!\! \!\!\!\! \!\! \!\!\!\! \!\! \!\!\!\! \!\! \!\!\!\! \! 
\!\! \!\! \!\!\!\! \!\! \!\!\!\! \!\! \!\!\!\! \!\! \!\!\!\! \!\! \!\!  
\!\! \!\! \!\!\!\! \!\! \!\!\!\! \!\! \!\! \uth (\phi)   \\
\uy (\phi) +  {[D \uth (\phi)]}^{-\top} \zeta -  
\big[ D {\tilde \uz} ( \uth (\phi ) )\big]^\top {\bar J} w   \\
\!\! \!\! \!\!\!\! \!\! \!\!\!\! \!\! \!\!\!\! \!\! \!\!\!\! \!\! \!\!\!\! \!\! \!\!
\!\! \!\! \!\!\!\! \!\! \!\!\!\! \!\! \!\!\!\! \!\! \!\!\!\! \!\! \!\!\!\! \!\! \!\! 
\!   \uz ( \phi ) + w
\end{array}
\right)
\ee
where  ${\tilde \uz} (\teta) := (\uz \circ  \uth^{-1}) (\teta ) $, see \eqref{torus graph}. The transposed operator 
$ \big[ D {\tilde \uz} ( \teta )\big]^\top : E \to \R^\nu $ is defined by the duality relation 
$$ 
\big[ D {\tilde \uz} ( \teta )\big]^\top w \cdot \hat \theta =  \langle w,  D {\tilde \uz} ( \teta ) [\hat \theta] \rangle \, , \quad 
\forall w \in E \, , \ \hat \theta \in \R^\nu \, . 
$$ 

\begin{lemma}\label{lemma:symplectic}
Let $ i $ be an isotropic torus embedding. Then $ G $ is symplectic. 
\end{lemma}

\begin{proof}
We may see $ G $ as the composition 
$ G := G_2 \circ G_1  $ of the  diffeomorphisms
$$
\left(
\begin{array}{c}
  \teta  \\
  y \\
 z
\end{array}
\right) = 
G_1 
\left(
\begin{array}{c}
  \phi  \\
  \zeta \\
 w
\end{array}
\right) 
:= 
\left(
\begin{array}{c}
   \uth (\phi )  \\
{[D \uth (\phi)]}^{-\top} \zeta \\
 w
\end{array}
\right) 
$$
and 
\be\label{G2}
\left(
\begin{array}{c}
  \teta  \\
  y \\
 z
\end{array}
\right) 
\mapsto G_2
\left(
\begin{array}{c}
  \teta  \\
  y \\
 z
\end{array}
\right) 
:=
\left(
\begin{array}{c}
 \!\! \!\! \!\! \!\! \!\! \!\! 
 \!\! \!\! \!\! \!\! \!\! \!\!\!\! \!\! \!\! \!\! \!\! \!\! \!\! \!\! \!\! \!\! \!\! \!\! \!\! \!\! \!\!\!\teta  \\
  \tilde \uy (\teta) + y  -   \big[ D {\tilde \uz}( \teta )\big]^\top {\bar J} z
   \\
\!\! \!\!\!\! \!\! \!\! \!\!\!\! \!\!\!\! \!\!\!\! \!\! \!\! \!\! \!\! \!\! \!\! \!\! \!\! \! {\tilde \uz}  (\teta) +z  
\end{array}
\right) 
\ee
where
$ \tilde \uy := \uy \circ \uth^{-1} $,  $ {\tilde \uz}  := \uz \circ \uth^{-1} $, see \eqref{torus graph}.
\smallskip
We claim that both $ G_1 $, $ G_2 $ are symplectic, whence the lemma follows.
\\[1mm] 
{\sc $ G_1 $  is symplectic.} 
Since $ G_1 $ is the identity in the third component, it is sufficient to check that 
the map 
$$ 
(\phi, \zeta) \mapsto \big( \uth (\phi ),   [D \uth (\phi)]^{-\top} \zeta   \big) 
$$
is a symplectic diffeomorphism on $ \T^\nu \times \R^\nu $, which is a direct calculus. 
\\[1mm]
{\sc $ G_2 $  is symplectic.} We prove that $ G_2^* \form - \form $ is closed and so (see \eqref{exact})
$$ 
G_2^* {\cal W} = G_2^* d \form= d G_2^* \form = d \form = {\cal W} \, . 
$$
By \eqref{G2} and the definition of pullback we have
\begin{align*}
(G_2^* \form )_{(\teta, y, z)}[\hat \teta, \hat y, \hat z] & = 
\big(   \tilde \uy (\teta) + y - \big[ D {\tilde \uz}( \teta )\big]^\top {\bar J} z \big) 
\cdot \hat \teta \\ 
& \quad  + \frac12 \langle {\bar J} ( {\tilde \uz}( \teta ) + z ) , \hat z + D \tilde \uz(\teta)[\hat \teta] \rangle \, .
\end{align*}
Therefore (recall \eqref{lambda})
\begin{align}
\big( (G_2^* \form )_{(\teta, y, z)} - \form_{(\teta, y, z)}\big) [\hat \teta, \hat y, \hat z] & =  
\big(   \tilde \uy (\teta) - \big[ D {\tilde \uz}( \teta )\big]^\top {\bar J} z  \big) 
\cdot \hat \teta + \frac12 \langle {\bar J} {\tilde \uz}( \teta )   , \hat z  \rangle \nonumber \\
&  +
\frac12 \langle {\bar J} {\tilde \uz}( \teta )   , D \tilde \uz(\teta) [ \hat \teta]  \rangle +
\frac12 \langle {\bar J} z  ,  D \tilde \uz(\teta)[ \hat \teta] \rangle \nonumber  \\
& =  
  \tilde \uy (\teta) \cdot \hat \teta + \frac12 \langle {\bar J} {\tilde \uz}( \teta )   , D \tilde \uz(\teta) [ \hat \teta] \rangle  \nonumber \\
 &   +  
\frac12 \langle {\bar J} {\tilde \uz}( \teta )   , \hat z  \rangle  +
\frac12 \langle {\bar J}  D {\tilde \uz}( \teta ) [\hat \teta]  ,  z  \rangle \, ,   \label{r1}
\end{align}
having used that $ {\bar J}^\top = - {\bar J} $. We first note that the $ 1 $-form 
\be\label{r2}
( \hat \teta, \hat y, \hat z ) \mapsto  \langle {\bar J} {\tilde \uz}( \teta ), \hat z  \rangle + 
 \langle {\bar J}  D {\tilde \uz}( \teta ) [\hat \teta ] ,  z  \rangle 
 = d (  \langle {\bar J} {\tilde \uz}( \teta )   , z  \rangle  )[  \hat \teta, \hat y, \hat z ] 
\ee
is exact. Moreover 
\be\label{r3}
  \tilde \uy (\teta) \cdot \hat \teta + \frac12 \langle {\bar J} {\tilde \uz}( \teta ), 
D \tilde \uz(\teta) [\hat \teta] \rangle = ( j^* \form )_\teta [\hat \teta ]  
\ee
(recall \eqref{lambda}) where  $ j := i \circ \uth^{-1} $ see \eqref{torus graph0}. 
Hence \eqref{r1}, \eqref{r2}, \eqref{r3} imply
$$
(G_2^* \form )_{(\teta, y, z)} - \form_{(\teta, y, z)} = \pi^*(j^* \form)_{(\teta, y, z)}  + 
d 	\Big(  \frac12 \langle {\bar J} {\tilde \uz}( \teta )   , z  \rangle  \Big) 
$$
where $\pi : \T^\nu \times \R^\nu \times E \to \T^\nu$ is the canonical projection. 

Since the torus $ j (\T^\nu) = i(\T^\nu) $
is isotropic 
the $ 1 $-form $ j^* \form $ on $\T^\nu$ is closed (as $i^* \form $, see \eqref{i star closed}).
This concludes the proof.
\end{proof}

Since $ G $ is symplectic 
 the transformed Hamiltonian vector field 
$$ 
G^* X_K := (DG)^{-1} \circ X_K \circ G = X_{\mathtt K} \, , \quad {\mathtt K} := K \circ G \, , 
$$
is still Hamiltonian. By construction 
(see \eqref{cG}) the torus $ \{ \zeta = 0,  w = 0 \} $ is invariant and 
\eqref{Boundary value}  implies 
$$
X_{\mathtt K} (\phi , 0,0)=(\om , 0,0)
$$ (see also Lemma \ref{lem:sc}).
As a consequence, the Taylor expansion of the transformed Hamiltonian $ {\mathtt K} $ 
in these new coordinates assumes the normal form 
\be\label{normalformH}
{\mathtt K} = {const}
+  \om \cdot \zeta +   \frac12 A(\phi) \zeta \cdot \zeta +  \langle C (\phi) \zeta , w \rangle +  
\frac12 \langle B (\phi) w, w \rangle  + O_3 (\zeta, w) 
\ee
where $ A (\phi) \in {\rm Mat}(\nu \times \nu ) $ is a  real symmetric matrix, 
$ B (\phi)$ is a self-adjoint operator of $ E $,  $ C (\phi) \in {\cal L}(\R^\nu, E) $, and 
$ O_3 (\zeta, w) $ collects all the terms at least cubic in $ (\zeta, w) $.

We have proved the following theorem:

\begin{theorem}\label{thm1} {\bf \cite{BBField}}
{\bf (Normal form close to an invariant isotropic torus)}
Let $ {\cal T} = i (\T^\nu) $ be an embedded  torus, see \eqref{embedded torus}-\eqref{smooth torus},  
 which is a smooth graph over $ \T^\nu $, see \eqref{torus graph0}-\eqref{torus graph},  
 invariant  for the Hamiltonian vector field $ X_K $, and on which the flow is conjugate to the translation flow of vector $\om$,
see \eqref{invariatorus}-\eqref{linear flow}. Assume moreover that ${\cal T}$ is  {\sc isotropic},  a property 
which is automatically verified  if $\om $ is non-resonant.

Then there exist symplectic coordinates  $(\phi, \zeta, w) $ 
in which $ {\cal T } $ is described by 
$$ \T^\nu \times \{ 0 \} \times \{ 0 \} $$ 
and
the Hamiltonian assumes the normal form \eqref{normalformH}, i.e. 
the torus  
$$ {\cal T} = G (\T^\nu \times \{ 0 \} \times \{ 0 \}) 
$$ where $ G $ is 
the symplectic diffeomorphism defined  in  \eqref{cG},  and 
the Hamiltonian  $ K \circ G $ has the Taylor expansion 
 \eqref{normalformH}  in a neighborhood  of the invariant torus. 
\end{theorem}

The normal form \eqref{normalformH} is  relevant  in view of a Nash-Moser approach, 
because it  provides a control of the linearized equations in the normal bundle of the torus. 
The linearized Hamiltonian system associated to $ {\mathtt K} $ at the trivial  solution $ (\phi, \zeta, w) (t) =  (\om t, 0, 0 ) $ 
is 
$$
\begin{cases}
\dot  \phi  - A( \om t ) \zeta -  [C(\om t )]^\top w =0
\cr 
\dot  \zeta  = 0    \cr
\dot w - J \big( B (\om t ) w +  C (\om t ) \zeta   \big)=0  
\end{cases}
$$
and note that the second equation is decoupled from the others. 
Inserting its constant solution $ \zeta (t) =  \zeta_0 $
in the third equation we are reduced to solve the quasi-periodically forced 
Hamiltonian linear 
equation in  $ w $,
$$
w_t - J  B (\om t ) w 
 =  g(\omega t) \, , \quad g(\omega t ) := J   C (\om t ) \zeta_0 \, . 
$$
This linear system may be studied  with both  the reducibility and  the 
multiscale  techniques presented 
 in Sections  \ref{sec:red} and  \ref{sec:MULTI}.
In particular, if the reducibility approach outlined in subsection \ref{sec:PR} applies,
there is a symplectic change of variable which makes $ B(\phi ) $ constant.

\section{Symplectic coordinates near an approximately invariant torus} 

In this section we report 
a construction of suitable symplectic coordinates near a torus which is only approximately invariant,  analogous  to the one in the previous section. 
 
For that, we first report  a basic fact about $ 1 $-forms on a torus. We regard  
a $ 1 $-form $ a = \sum_{i=1}^\nu a_i(\vphi) d \vphi_i $  equivalently 
as the vector field $  {\bf a}( \vphi ) = ( a_1 (\vphi) , \ldots, a_\nu (\vphi) ) $. 

Given a function $ g : \T^\nu \to \R $ with zero average, we  denote by 
$$ 
u := \Delta^{-1} g 
$$ 
the unique solution of $ \Delta u = g $ with zero average. 

\begin{lemma}\label{lem:deco1} {{\bf (Helmotz decomposition)}}
A smooth vector field $ \bf a $ on $ \T^\nu $ may be decomposed as the sum of a conservative and a 
divergence-free vector field: 
\be\label{decom1}
{\bf a} =  \nabla U +  {\bf c} + {\bf \rho}  \, ,  \quad U 
 : \T^\nu \to \R \, , \ \  {\bf c} \in \R^\nu \, , \ \   {\rm div} {\bf \rho } = 0 \, , \ \ \int_{\T^\nu} \rho d \vphi = 0  \, .  
\ee
The above decomposition is unique if we impose that the mean value of $U$ vanishes.
We have that
$$ 
U =   \Delta^{-1} ( {\rm div} \, {\bf a}  ) \, , 
$$ 
the components of $\rho$ are
\be\label{Deltagf}
\rho_j (\vphi ) =  \Delta^{-1}   \sum_{k=1}^\nu \partial_{\vphi_k} A_{kj} (\vphi)  \, , 
\quad  A_{kj} := \partial_{\vphi_k} a_j - \partial_{\vphi_j} a_k  \, , 
\ee
and  
$$  
{\bf c} = (c_j)_{j=1, \ldots, \nu} \, , \quad 
c_j = (2\pi)^{-\nu} \int_{\T^\nu} a_j (\vphi) \, d \vphi \, . 
$$ 
\end{lemma}

\begin{proof} Notice that 
${\rm div} (\vec a-\nabla U)=0$ if and only if $ {\rm div} \, {\bf a} = \Delta U  $. 
This equation has the solution $ U :=  $ $ \Delta^{-1} ( {\rm div} \, {\bf  a}  )$
(note that $ {\rm div} \, {\bf  a } $ has zero average). 
Hence 
\eqref{decom1} is achieved with $  {\bf  \rho } := {\bf  a}  - \nabla U  - {\bf c } $. 
By taking the $\vphi$-average we get that each 
$$ 
c_j = (2\pi)^{-\nu} \int_{\T^\nu} a_j (\vphi) \, d \vphi \, .
$$
Let us now prove the expression  \eqref{Deltagf} of $ \rho_j $.
We have 
$$ 
\partial_{\vphi_k} \rho_j - \partial_{\vphi_j} \rho_k = 
 \partial_{\vphi_k} a_j - \partial_{\vphi_j} a_k =:  A_{kj}  
 $$ 
because 
$$ 
\pa_{\vphi_j} \pa_{\vphi_k} U - \pa_{\vphi_k} \pa_{\vphi_j} U = 0 \, . 
$$
For each $ j=1, \ldots, \nu $ we differentiate
$ \partial_{\vphi_k} \rho_j - \partial_{\vphi_j} \rho_k = A_{kj} $ 
with respect to $ \vphi_k $ and we sum in $ k $, obtaining
$$ 
\Delta \rho_j - \sum_{k=1}^\nu \partial_{\vphi_k \vphi_j} \rho_k
=   \sum_{k=1}^\nu \partial_{\vphi_k} A_{kj} \, .
$$ 
Since 
$$ 
\sum_{k=1}^\nu \partial_{\vphi_k \vphi_j} \rho_k  =  \partial_{\vphi_j} {\rm div} {\bf  \rho } =  0 
$$ 
then
$ \Delta \rho_j =  \sum_{k=1}^\nu \partial_{\vphi_k} A_{kj}  $ and 
\eqref{Deltagf}  follows.
\end{proof}

\begin{corollary}\label{closed}
Any closed $ 1$-form on $ \T^\nu $ has the form  $ a (\vphi)  = {\bf c} + d U $ for some $ {\bf c} \in \R^\nu $. 
\end{corollary}

\begin{corollary}\label{nearby}
Let $ a (\vphi)  $ be a $ 1$-form on $ \T^\nu $, and let $\rho$ be defined by $(\ref{Deltagf})$.
Then  $ a - {\sum}_{j=1}^\nu \rho_j (\vphi) d \vphi_j $ is closed. 
\end{corollary}

We quantify how an embedded torus  $ i (\T^\nu) $ is approximately invariant\index{Approximately invariant torus}
 for the  Hamiltonian vector field $ X_{K } $
in terms of the ``error function"\index{Error function} 
\be\label{Z error}
Z(\vphi) := 
{\cal F}( i ) = 
(\om \cdot \partial_\vphi  i) (\vphi) - X_{K} (i (\vphi) )  \, . 
\ee
Consider the pullback 
$ 1 $-form  on $ \T^\nu $  (see \eqref{lambda}) 
\begin{align}\label{pull1form}
(i^* \form ) (\vphi ) & = \sum_{k=1}^\nu a_k (\vphi) d \vphi_k  
\end{align}
where
\begin{eqnarray} \label{pull1form coefficients}
a_k (\vphi)  &:= &\Big[\big[ D  \uth (\vphi) \big]^\top \uy(\vphi) + \frac12 [D \uz (\vphi)]^\top {\bar J} {\uz}(\vphi) \Big]_k  
\nonumber \\
&=& \uy (\varphi) \cdot \frac{\partial \uth}{\partial \varphi_k} (\varphi) +
\frac{1}{2} \la \bar{J} \uz (\varphi) , \frac{\partial \uz}{\partial \varphi_k} (\varphi)\ra \, ,
\end{eqnarray}
and 
the $ 2 $-form (recall \eqref{exact})
\be\label{istarW}
\begin{aligned}
&  i^* {\cal W} =  d  (i^* \form ) = { \sum}_{k < j }  A_{kj} (\vphi) d \vphi_k \wedge  d \vphi_j \, ,\\
&  \qquad   \ A_{kj} (\vphi) = \partial_{\vphi_k} a_j (\vphi) - \partial_{\vphi_j} a_k (\vphi) \, . 
\end{aligned}
\ee
We call  the coefficients $ (A_{kj})  $    the ``lack of isotropy" of the torus embedding
 $\varphi \mapsto  i (\vphi) $. 
In Lemma \ref{lem:approclo} below we quantify their size 
in terms of the error function $ Z $ defined in \eqref{Z error}.
We first recall that  the Lie derivative\index{Lie derivative} of a $ k $-form $ \b $ with respect to the vector field 
 $ Y $ is 
 $$ 
 L_Y  \b := \frac{d}{dt} \big[(\Phi^t_Y)^* \b \big]_{|t = 0} 
 $$
where $ \Phi^t_Y $ denotes the flow generated by $ Y $.

Given a function $ g(\vphi) $ with zero average, 
we denote by  $ u := (\omega \cdot \pa_\vphi)^{-1} g  $ 
the unique solution of $ \om  \cdot \pa_\vphi u = g $ with zero average.  

\begin{lemma}\label{lem:approclo}
The ``lack of isotropy" coefficients $ A_{kj} $ satisfy, $ \forall \vphi \in \T^\nu $,  
\begin{align}\label{coefficientiAkj}
(\om \cdot \partial_\vphi) {A_{kj}(\vphi) } & =  {\cal W} \big( DZ(\varphi)\underline{e}_k , Di (\vphi) \underline{e}_j \big) 
+ {\cal W} \big(Di (\vphi)\underline{e}_k, DZ(\varphi)\underline{e}_j  \big) 
\end{align}
where $ (\underline{e}_1 , \ldots , \underline{e}_\nu) $ denotes the canonical basis  of $ \R^\nu $. 
Thus,  since each $A_{kj}$ has zero mean value,  if the frequency vector $ \om  \in \R^\nu $ is non-resonant, then 
\be\label{Akj}
 A_{kj} (\vphi) = 
(\om \cdot \pa_\vphi)^{-1} \big(  {\cal W} \big( DZ(\varphi)\underline{e}_k , Di (\vphi) \underline{e}_j \big)  
 + {\cal W} \big(Di (\vphi)\underline{e}_k, DZ(\varphi)\underline{e}_j  \big)  \big) \, . 
\ee
\end{lemma}

\begin{proof}
We use Cartan's formula\index{Cartan formula}
$$ 
 L_\om (i^* {\cal W} )=d \big( (i^* {\cal W})(\om , \cdot ) \big) + \big(d(i^* {\cal W}) \big) (\om ,  \cdot ) \, .
 $$ 
Since $d(i^* {\cal W})=i^* d {\cal W} =0 $ by \eqref{exact} we get
\be \label{Lieder}  
L_\om (i^* {\cal W})=d \big( (i^* {\cal W}) (\om ,  \cdot ) \big) \, . 
\ee
Now we compute, for $\hat \phi \in \R^\nu $   
\begin{align*}
(i^* {\cal W}) (\om , \hat \phi) & =  {\cal W} (Di (\vphi) \om , Di (\vphi) \hat \phi) 
= {\cal W} ( X_K (i(\vphi)) 
+ Z(\vphi) , Di (\vphi) \hat \phi) \nonumber \\
& = -dK(i(\vphi))[Di (\vphi) \hat \phi] 
+ {\cal W} ( Z(\vphi), Di (\vphi) \hat \phi) \, .
\end{align*}
We obtain 
$$ 
(i^* {\cal W}) (\om ,  \cdot )=   \sum_{j=1}^\nu b_j(\vphi) d\vphi_j
$$
$$
b_j(\varphi)= (i^* {\cal W}) (\om , \underline{e}_j)= 
-\frac{\pa (K \circ i)}{\pa \vphi_j} (\vphi) 
+ {\cal W} (Z(\vphi) , Di(\vphi) \underline{e}_j) \, .
$$
Hence, by \eqref{Lieder}, the Lie derivative
\be\label{primaLieder}
L_\om (i^* {\cal W}) =  \sum_{k<j} B_{kj} (\vphi) d\vphi_k \wedge d \vphi_j
\ee
with
\begin{align} \label{espr1}
B_{kj}(\vphi) &= \frac{\pa b_j}{\pa \vphi_k} (\vphi) - \frac{\pa b_k}{\pa \vphi_j} (\vphi) \nonumber \\
&=  \frac{\pa }{\pa \vphi_k} ( {\cal W} (Z(\vphi),Di (\vphi) \underline{e}_j   )) - 
\frac{\pa }{\pa \vphi_j}  ( {\cal W} (Z(\vphi),Di (\vphi) \underline{e}_k   )) \nonumber \\
&= {\cal W} (DZ (\vphi) \underline{e}_k , Di (\vphi) \underline{e}_j) + 
{\cal W} (Di (\vphi) \underline{e}_k ,  DZ (\vphi) \underline{e}_j ) \, .
\end{align}
Recalling \eqref{linear flow} and \eqref{istarW} we have,
$ \forall \vphi \in \T^\nu $, 
$$
( \psi^t_\om )^* (i^* {\cal W}) (\vphi)  = i^* {\cal W} (\vphi + \om t ) =  \sum_{k < j } 
A_{kj}(\vphi + \om t) d \vphi_k \wedge d \vphi_j \, . 
$$
Hence the Lie derivative
\be\label{espr2}
L_\om  (i^* {\cal W}) (\varphi) = \sum_{k < j } ({\om \cdot \pa_\vphi}  A_{kj})(\vphi) d \vphi_k \wedge d \vphi_j \, .
\ee
Comparing \eqref{primaLieder}-\eqref{espr1} and \eqref{espr2} we deduce  \eqref{coefficientiAkj}.  
\end{proof}

The previous lemma provides another proof of Lemma \ref{lem:iso}. For an invariant torus embedding $ i(\vphi) $ 
the ``error function" $ Z(\vphi) = 0 $ (see  \eqref{Z error}) and so each $ A_{kj}  = 0 $. 
We now prove that near an approximate isotropic torus there is an isotropic torus.

\begin{lemma}\label{modified}
{\bf (Isotropic torus)}   
The torus embedding $ i_\d (\vphi) = (  \uth (\vphi) , y_\d (\vphi), \uz (\vphi) )$ defined by 
\be\label{modified action}
y_\d (\vphi) = \uy (\vphi) - [D  \uth(\vphi)]^{-\top} \rho (\vphi) \, , 
\quad  
 \rho_j := \Delta^{-1} \Big(  \sum_{k=1}^\nu \partial_{\vphi_k} A_{kj} (\vphi) \Big) \, , 
\ee
is  isotropic.
\end{lemma}

\begin{proof}
By Corollary \ref{nearby} the $ 1 $-form $ i^* \form - \rho $ is closed
with  $ \rho_j  $  defined in \eqref{modified action}, see also \eqref{Deltagf}, \eqref{pull1form}.  
Actually  
$$ 
i^* \form - \rho = i_\d^* \form 
$$ 
is the pullback of the  $ 1$-form  $ \form $  under
the  modified torus embedding $ i_\d $  
defined  in  \eqref{modified action}, see \eqref{pull1form coefficients}. Thus the torus $ i_\d (\T^\nu)  $ is isotropic. 
\end{proof}

In analogy with Theorem  \ref{thm1} we now introduce a symplectic set of coordinates $ (\phi, \zeta, w ) $ 
near the  isotropic torus $ {\cal T}_\d := i_\d (\T^\nu) $  
via the symplectic diffeomorphism 
\be\label{changevarappro} 
\left(
\begin{array}{c}
\teta  \\
y   \\
z
\end{array}
\right) = 
G_\d
\left(
\begin{array}{c}
\phi \\
\zeta  \\
w
\end{array}
\right) 
:= 
\left(
\begin{array}{c}
 \!\!\!\!\!\!\!\!\!\!\!\!\!\!\!\!\!\!\!\!\!\!\!\! \!\!\!\!\!\!\!\!\!\!\!\!\!\!\!\!\!\!\!\!\!\!\!\!\!\!\!\!\!\!\!\!\!\!\!\!\!\!\!\!\!\!\!\!\!\!\!\! \!
 \!\!\!\!\!\!\!\!\!\!\!\! \uth (\phi)   \\
y_\d (\phi) +  {[D  \uth (\phi)]}^{-\top} \zeta -  
\big[ D {\tilde \uz } ( \uth(\phi) )\big]^\top {\bar J} w   \\
 \!\!\!\!\!\!\!\!\!\!\!\!\!\!\!\!\!\!\!\!\!\!\!\! \!\!\!\!\!\!\!\!\!\!\!\!\!\!\!\!\!\!\!\!\!\!\!\!\!\!\!\!\!\!\!\!\!\!\!\! \!\!\!\!\!\!\!\!\!\!\!\! \!\!  \uz ( \phi ) + w
\end{array}
\right)
\ee
where 
$ \tilde \uz   := \uz \circ   \uth^{-1}  $.
The map $ G_\d $ is symplectic by Lemma \ref{lemma:symplectic} because 
$ i_\d $ is isotropic (Lemma \ref{modified}). 
In the new coordinates  $ (\phi, \zeta, w) $ 
the isotropic torus embedding $ i_\d $ is  trivial, namely 
$$ 
i_\d (\phi) = G_\d (\phi, 0, 0  ) \, . 
$$ 
Under the symplectic change of variable \eqref{changevarappro}, the Hamiltonian 
vector field 
$ X_{K} $  changes into 
\be\label{Kmu0}
X_{{\mathtt K}} =G_\delta^* X_{K} =  (DG_\d)^{-1}   X_{K} \circ  G_\d \qquad {\rm where} \qquad
{\mathtt K} := K \circ G_\d  \, .
\ee
The Taylor expansion of   the new Hamiltonian $ {\mathtt K} : \R^\nu \times \R^\nu \times E \to \R$ 
at the trivial torus  $ (\phi, 0, 0 ) $  is 
\be
\begin{aligned}\label{Kmud}
{\mathtt K}  &  =   
{\mathtt K}_{00}(\phi ) + {\mathtt K}_{10} ( \phi ) \cdot \zeta +   \langle {\mathtt K}_{01} ( \phi ),  w \rangle    \\
& \ \ \ +   \frac12 {\mathtt K}_{20} ( \phi ) \zeta \cdot \zeta  
 +   \langle {\mathtt K}_{11} (  \phi ) \zeta, w \rangle +  \frac12 \langle {\mathtt K}_{02} (\phi ) w, w \rangle + {\mathtt K}_{\geq 3}  ( \phi, \zeta, w)  
\end{aligned}
\ee
where $ {\mathtt K}_{\geq 3} $ collects all the terms at least cubic in the variables $ (\zeta , w )$.
The Taylor coefficients of $ {\mathtt K} $  are ${\mathtt K}_{00}(\phi) \in \R $,  
${\mathtt K}_{10}(\phi) \in \R^\nu $,  
${\mathtt K}_{01}(\phi) \in E $, 
${\mathtt K}_{20}(\phi) \in {\rm Mat}(\nu \times \nu ) $ is a  real symmetric matrix, 
${\mathtt K}_{02}(\phi)$ is a self-adjoint operator of $ E $ and 
${\mathtt K}_{11}(\phi) \in {\cal L}(\R^\nu, E) $. 

As seen in Theorem  \ref{thm1}, if $ i_\d $ were an invariant torus embedding,   
the coefficient $ {\mathtt K}_{00} (\phi
) = {\rm const} $, $ {\mathtt K}_{10}  (\phi) = \om $ and $ {\mathtt K}_{01}  (\phi)  = 0 $. 
We now express these coefficients 
in terms of the  error function $ Z_\d  := {\cal F} ( i_\d) $.  

\begin{lemma}\label{lem:sc}
The vector field 
\begin{align}
X_{{\mathtt K}} (\phi, 0, 0) & = 
\left(
\begin{array}{c}
\!\! {\mathtt K}_{10}( \phi) \!\!   \\
\!\!  -  \partial_{\phi} {\mathtt K}_{00}(\phi )  \!\! \\
\!\! J {\mathtt K}_{01} ( \phi ) \!\!
 \end{array}
\right)  = 
\left(
\begin{array}{c}
\!\! \om  \!\! \\
\!\! 0 \!\! \\
\!\! 0 \!\!
 \end{array}
\right)
 - \big( DG_\d (\phi, 0, 0) \big)^{-1} Z_\d (\phi)  \, . \label{apprVFA}
\end{align}
\end{lemma}

\begin{proof}
By \eqref{Kmu0} and $ i_\d (\phi ) = G_\d (\phi, 0, 0)  $, we have 
$$
\begin{aligned}
 X_{{\mathtt K}} (\phi, 0, 0) & =  
DG_\d (\phi, 0, 0)^{-1} X_{K} (  i_\d (\phi ) ) \\ 
& =  
DG_\d (\phi, 0, 0)^{-1} \big( \om  \cdot \pa_\vphi i_\d (\phi) - Z_\d(\phi) \big)  
\end{aligned}
$$
and \eqref{apprVFA} follows 
because $ DG_\d (\phi, 0, 0)^{-1}  Di_\d (\phi)[\om]  = (\om, 0, 0 ) $.
\end{proof}

We finally write 
the expression of the 
coefficients $ {\mathtt K}_{11} (\phi) $, $ {\mathtt K}_{20}(\phi) $  in terms of $ K $, 
which is used in Chapter \ref{sezione almost approximate inverse}. 

\begin{lemma}\label{smalldepd} The coefficients
\begin{align}
{\mathtt K}_{11}(\phi ) & = 
D_y  \nabla_z K (i_\delta(\phi )) [ D\theta_0 (\phi )]^{-\top}  
+ \bar J (D_\theta {\tilde z}_0) (\teta_0(\phi )) (D_y^2 K) (i_\delta(\phi )) 
[D \theta_0 (\phi )]^{-\top} \label{matK1} \\
{\mathtt K}_{20} ( \phi )  & = [D \teta_0 (\phi)]^{-1} (D_y^2 K)( i_\d (\phi) ) )
 [D \teta_0 (\phi)]^{-\top} \, . \label{matK2}
\end{align}
\end{lemma}

\begin{proof}
Formulas \eqref{matK1}\eqref{matK2} follow 
differentiating $ {\mathtt K} = K \circ G_\d $. 
\end{proof}

\backmatter

\printindex

\end{document}